\documentclass[letterpaper, 12pt]{memoir}

\usepackage{etex}
\usepackage[UKenglish]{babel}
\usepackage{amsmath}
\usepackage{amssymb}
\usepackage{amsthm}
\usepackage{amsfonts}
\usepackage{mathtools} 
\usepackage{mathabx}
\usepackage{bbold}
\usepackage{thmtools}
\usepackage{tikz-cd} 
\usetikzlibrary{arrows,arrows.meta}
\usepackage{bussproofs} 
\usepackage{stmaryrd} 
\usepackage{microtype}
\usepackage{mdframed} 
\usepackage{url}
\usepackage{wrapfig}
\usepackage{xspace}
\usepackage{microtype}
\usepackage{rotating} 
\usepackage{pdflscape} 
\usepackage{array,booktabs}
\usepackage[refpage, intoc]{nomencl} 
\usepackage{enumitem}
\usepackage{epigraph} 
\usepackage{geometry}
\usepackage{custom-titlepage}
\usepackage[pdftitle={Cartesian closed bicategories: type theory and coherence},
              pdfauthor={Philip Saville},
              pdfkeywords={bicategory, cartesian closed, coherence, 
              	Lambek-Curry-Howard correspondence, simply-typed lambda 
              	calculus, 2-dimensional type theory, 
              	normalisation-by-evaluation}]{hyperref}

\makeatletter
\def\maketag@@@#1{\hbox{\m@th\normalfont\normalsize#1}}
\makeatother

\setlist{itemsep = 0mm,
		topsep = 1mm}

\setlrmarginsandblock{2cm}{*}{1}
\setulmarginsandblock{2.2cm}{*}{1}
\checkandfixthelayout[nearest]

\DeclareMathAlphabet{\mathsf}{OT1}{kurier}{m}{n}
\SetMathAlphabet{\mathsf}{bold}{OT1}{kurier}{bx}{n}

\allowdisplaybreaks

\setsecnumdepth{subsection} 

\newcommand{\hide}[1]{}

\newenvironment{manyfigcap}
{
\begin{minipage}{\textwidth}
\centering
\vspace{10pt}
}
{
\end{minipage}
}

\newcommand{\blackqed}{$\hfill \ensuremath{\blacktriangleleft}$}

\theoremstyle{plain}

\theoremstyle{definition}
\newtheorem*{nonumthm}{Theorem}

\declaretheoremstyle[
headfont=\normalfont\bfseries,
bodyfont=\normalfont,
spaceabove=1em plus 0.75em minus 0.25em,
spacebelow=1em plus 0.75em minus 0.25em,
qed={\qedsymbol},
]{prooflesslemmastyle}

\declaretheoremstyle[
headfont=\normalfont\bfseries,
bodyfont=\normalfont,
spaceabove=1em plus 0.75em minus 0.25em,
spacebelow=1em plus 0.75em minus 0.25em,
]{lemmastyle}

\declaretheoremstyle[
headfont=\normalfont\bfseries,
bodyfont=\normalfont,
spaceabove=1em plus 0.75em minus 0.25em,
spacebelow=1em plus 0.75em minus 0.25em,
qed={\blackqed},
]{exmpstyle}

\declaretheorem[
style=exmpstyle,
title=Principle,
sharenumber=theoremlocal,
refname={principle,principles},
Refname={Principle,Principles},
]{myprinciplelocal}

\declaretheorem[
style=exmpstyle,
title=Example,
sharenumber=theorem,
refname={example,examples},
Refname={Example,Examples},
]{myexmp}

\declaretheorem[
style=exmpstyle,
title=Conjecture,
sharenumber=theoremlocal,
refname={conjecture,conjectures},
Refname={Conjecture, Conjectures},
]{myconjecturelocal}

\declaretheorem[
style=exmpstyle,
title=Definition,
sharenumber=theorem,
refname={definition,definitions},
Refname={Definition, Definitions},
]{mydefn}

\declaretheorem[
style=exmpstyle,
title=Notation,
sharenumber=theorem,
refname={notation, notations},
Refname={Notation, Notations}
]{mynotation}

\declaretheorem[
style=exmpstyle,
sharenumber=theorem,
title=Remark,
refname={remark, remarks},
Refname={Remark, Remarks},
]{myremark}

\declaretheorem[
style=prooflesslemmastyle,
title=Lemma,
sharenumber=theorem,
refname={lemma, lemmas},
Refname={Lemma, Lemmas}
]{prooflesslemma}

\declaretheorem[
style=prooflesslemmastyle,
title=Lemma,
sharenumber=theoremlocal,
refname={lemma, lemmas},
Refname={Lemma, Lemmas}
]{prooflesslemmalocal}

\declaretheorem[
style=prooflesslemmastyle,
title=Theorem,
sharenumber=theorem,
refname={theorem, theorems},
Refname={Theorem, Theorems}
]{prooflessthm}

\declaretheorem[
style=lemmastyle,
title=Theorem,
sharenumber=theorem,
refname={theorem, theorems},
Refname={Theorem, Theorems}
]{mythm}

\declaretheorem[
style=lemmastyle,
title=Corollary,
sharenumber=theorem,
refname={corollary, corollaries},
Refname={Corollary, Corollaries}
]{mycor}

\declaretheorem[
style=prooflesslemmastyle,
title=Corollary,
sharenumber=theorem,
refname={corollary, corollaries},
Refname={Corollary, Corollaries}
]{prooflesscor}

\declaretheorem[
style=prooflesslemmastyle,
title=Corollary,
sharenumber=theoremlocal,
refname={corollary, corollaries},
Refname={Corollary, Corollaries}
]{prooflesscorlocal}

\declaretheorem[
style=prooflesslemmastyle,
title=Proposition,
sharenumber=theorem,
refname={proposition, propositions},
Refname={Proposition, Propositions}
]{prooflesspropn}

\declaretheorem[
style=lemmastyle,
title=Proposition,
sharenumber=theorem,
refname={proposition, propositions},
Refname={Proposition, Propositions}
]{mypropn}

\declaretheorem[
style=prooflesslemmastyle,
title=Definition,
sharenumber=theoremlocal,
refname={proposition, propositions},
Refname={Proposition, Propositions}
]{prooflesspropnlocal}

\declaretheorem[
style=lemmastyle,
title=Lemma,
sharenumber=theoremlocal,
refname={lemma, lemmas},
Refname={Lemma, Lemmas}
]{lemmalocal}

\declaretheorem[
style=exmpstyle,
title=Definition,
sharenumber=theoremlocal,
refname={definition, definitions},
Refname={Definition, Definitions}
]{defnlocal}

\declaretheorem[
style=exmpstyle,
title=Construction,
sharenumber=theorem,
refname={construction, constructions},
Refname={Construction, Constructions}
]{myconstr}

\declaretheorem[
style=lemmastyle,
title=Lemma,
sharenumber=theorem,
refname={lemma, lemmas},
Refname={Lemma, Lemmas}
]{mylemma}

\newcommand{\subproof}[1]{\noindent\fbox{{#1}}\:}

\newcommand{\cf}{\emph{c.f.}}
\newcommand{\Def}[1]{\emph{#1}}

\newcommand{\eg}{\emph{e.g.}}
\newcommand{\ie}{\emph{i.e.}}

\newcommand{\st}{\mid} 
\newcommand{\rulename}[1]{\texttt{#1}}
\newcommand{\faketext}{\phantom{\Gamma}}

\newcommand{\U}{\mathrm{U}}
\newcommand{\T}{\mathrm{T}}
\newcommand{\K}{\mathrm{K}}
\renewcommand{\L}{\mathrm{L}}

\newcommand{\To}{\ensuremath{\Rightarrow}} 
\newcommand{\xra}[1]{\ensuremath{\xrightarrow{#1}}} 
\newcommand{\xla}[1]{\ensuremath{\xleftarrow{#1}}} 
\newcommand{\XRA}[1]{\ensuremath{\xRightarrow{#1}}}
\newcommand{\id}{\ensuremath{\mathrm{id}}} 
\newcommand{\Id}{\mathrm{Id}} 
\newcommand{\eval}{\mathrm{eval}} 
\newcommand{\evBar}{\mathrm{m}}
\newcommand{\evalMod}{\mathrm{E}} 

\newcommand{\cl}{\mathit{cl}} 

\newcommand{\wholecell}{[\underline{w}]}
\newcommand{\termobj}{{\mathbf{1}}}
\newcommand{\altexp}[2]{{\left[#1, #2 \right]}} 
\newcommand{\altexpsmall}[2]{[#1, #2]} 
\renewcommand{\exp}[2]{\expobj{#1}{#2}}
\newcommand{\glexp}[2]{{#1} \supset {#2}}
\newcommand{\gluedCo}[1]{\glued{\mathsf{e}}_{#1}}
\newcommand{\compact}[1]{\widetilde{{#1}}} 
\newcommand{\glueFun}{\mathfrak{J}} 

\renewcommand{\epsilon}{\varepsilon}

\DeclareMathOperator{\atasop}{\text{\textrm{@}}}
\newcommand{\concat}{\atasop}	
\newcommand{\psinv}[1]{{#1}^{{\star}}} 
\newcommand{\scriptsizepsinv}[1]{{#1}^{{\scriptscriptstyle\star}}} 
\newcommand{\iso}{\cong}
\newcommand{\op}[1]{{{#1}^\mathrm{op}}} 

\DeclareMathOperator*{\bilim}{bilim} 
\newcommand{\Psh}[1]{\mathcal{P}(#1)} 
\newcommand{\yon}{\mathrm{y}} 

\newcommand{\lan}[2]{\ensuremath{\mathrm{lan}_{#1} #2}} 
\newcommand{\seq}[1]{\langle #1 \rangle}
\newcommand{\seqlr}[1]{\left\langle #1 \right\rangle} 
\newcommand{\ind}[1]{{#1}_{\bullet}} 
\newcommand{\len}[1]{|#1|}
\newcommand{\gl}[1]{\mathrm{gl} {\left( {#1} \right)}}

\newcommand{\semlr}[1]{{\left\llbracket #1 \right\rrbracket}} 
\newcommand{\sem}[1]{\llbracket #1 \rrbracket} 
\newcommand{\usem}[2]{\ext{#2}{#1}} 
\newcommand{\hsem}[1]{h\sem{#1}} 
\newcommand{\adjDown}{{\mathbin{\rotatebox[origin=c]{90}{$\dashv$}}}}
\newcommand{\adjUp}{{\mathbin{\rotatebox[origin=c]{270}{$\dashv$}}}}
\newcommand{\Nat}{\mathbb{N}}
\newcommand{\Yon}{\mathrm{Y}} 
\newcommand{\Hom}{\mathrm{Hom}}
\newcommand{\twoHom}[2]{[{#1}, {#2}]}

\newcommand*\slice[2]{{#1} \mathclose{}/\mathopen{} {#2}}
\newcommand{\tens}{\otimes}
\newcommand{\cellOf}[1]{\overline{#1}}
\newcommand{\cospanCat}{\mathrm{C}}
\newcommand{\cospan}{(1 \xra{h_1} 0 \xla{h_2} 2)}
\newcommand{\constPseudofunctor}{\mathrm{const}}

\DeclareMathOperator*{\dotop}{.}
\renewcommand{\dot}{\mathbf{\dotop}}

\newcommand{\bicatA}{\mathcal{A}}
\newcommand{\bicatX}{\mathcal{X}}
\newcommand{\catA}{\ensuremath{\mathbb{A}}} 
\newcommand{\catB}{\ensuremath{\mathbb{B}}} 
\newcommand{\catC}{\ensuremath{\mathbb{C}}} 
\newcommand{\catD}{\ensuremath{\mathbb{D}}} 
\newcommand{\catJ}{\mathcal{J}}
\newcommand{\catX}{\mathbb{X}}
\newcommand{\catW}{\mathcal{W}}
\newcommand{\catV}{\mathcal{V}}
\newcommand{\cat}{\mathbb{C}}
\newcommand{\cloneOf}[1]{\mathrm{Cl}(#1)}

\newcommand{\stlcCloneTimes}[1]{\clone_{\Lambda^{\!{\times}\!}(#1)}}
\newcommand{\stlcCloneTimesExp}[1]{\clone_{\Lambda^{\!{\times},{\to}}(#1)}}
\newcommand{\altCat}{\mathcal{C}}
\newcommand{\baseCat}{\mathcal{B}}
\newcommand{\altaltCat}{\mathcal{D}}
\newcommand{\fpBicat}[1]{\left({{#1}, \Pi_n(-)}\right)} 
\newcommand{\ccBicat}[1]{\left({{#1}, \Pi_n(-), {\altTo}}\right)} 
\newcommand{\catOne}{\ensuremath{\mathbb{1}}} 
\newcommand{\Set}{\mathrm{Set}}
\newcommand{\Cat}{\mathbf{Cat}}
\newcommand{\funCat}[2]{\mathrm{Fun}(#1, #2)}

\newcommand{\MonCat}{\mathbf{MonCat}}

\newcommand{\termCatSymbol}{\mathcal{T}_{\mathrm{ps}}}

\newcommand{\termCatContextExt}{\termCatSymbol^{\text{@}, {\times}\!}}

\newcommand{\urestrict}[1]{{#1}{\big|}_1} 
\newcommand{\increstrict}[1]{\nucleus{#1}} 

\newcommand{\Con}{\mathrm{Con}} 

\newcommand{\CatCat}{\mathrm{Cat}}
\newcommand{\CartCatCat}{\mathrm{CartCat}}
\newcommand{\CCCatCat}{\mathrm{CCCat}}
\newcommand{\Bicat}{\mathbf{Bicat}}
\newcommand{\Gray}{\mathbf{Gray}}
\newcommand{\BicatCat}{\mathrm{Bicat}}
\newcommand{\Biclone}{\mathbf{Biclone}}
\newcommand{\BicloneCat}{\mathrm{Biclone}}
\newcommand{\Clone}{\mathrm{Clone}}
\newcommand{\CartClone}{\mathrm{CartClone}}
\newcommand{\CCClone}{\mathrm{CCClone}}
\newcommand{\CCBicloneCat}{\mathrm{CCBiclone}}
\newcommand{\CartBicloneCat}{\mathrm{CartBiclone}}
\newcommand{\CartBicatCat}{\mathrm{fp}\text{-}\mathrm{Bicat}}
\newcommand{\CCBicatCat}{\mathrm{cc}\text{-}\mathrm{Bicat}}
\newcommand{\Graph}{\mathrm{Grph}}
\newcommand{\TwoMultiGraph}{\mathrm{2}\text{-}\mathrm{MGrph}}
\newcommand{\BlankGraph}[1]{{{#1}}\text{-}\mathrm{Grph}}

\newcommand{\TwoGraph}{\BlankGraph{2}}
\newcommand{\MultiGraph}{\mathrm{MGrph}}
\newcommand{\Multicat}{\mathrm{MultiCat}}
\newcommand{\mCat}{\mathbb{L}}
\newcommand{\nCat}{\mathbb{M}}
\newcommand{\tensor}{\mathrm{T}} 
\newcommand{\mCatOf}[1]{\mathrm{M}{#1}}
\newcommand{\mbCat}{\mathcal{M}}

\newcommand{\altToHead}{\raisebox{.28mm}{\scriptsize$\rhd$}}
\DeclareMathOperator{\altTo}{\mathtt{=\!}\altToHead}
\DeclareMathOperator{\scriptsizealtTo}{\mathtt{=\!}\scriptsizealtToHead\!}
\newcommand{\scriptsizealtToHead}{{\raisebox{.21mm}{\tiny$\rhd$}}}
\DeclareMathOperator{\scriptsizescriptsizealtTo}{\mathtt{=\!}\scriptsizescriptsizealtToHead\!}
\newcommand{\scriptsizescriptsizealtToHead}{\raisebox{.05mm}{\tiny$\rhd$}}
\newcommand{\expobj}[2]{#1 \altTo #2} 
\newcommand{\scriptsizeexpobj}[2]{#1 \scriptsizealtTo #2} 
\newcommand{\scriptsizescriptsizeexpobj}[2]{#1 \scriptsizescriptsizealtTo #2} 

\newcommand{\subCat}{\mathcal{S}} 
\newcommand{\subFun}{E}

\renewcommand{\a}{\mathsf{a}}
\renewcommand{\r}{\mathsf{r}}
\renewcommand{\l}{\mathsf{l}}
\newcommand{\p}[2]{\mathsf{p}^{(#1)}_{#2}}
\newcommand{\pPicked}{\overline{\mathsf{p}}} 

\newcommand{\subbic}[1]{[#1]}
\newcommand{\natTrans}{\mathsf{k}}
\newcommand{\natCell}{\overline{\natTrans}}
\newcommand{\natTransAux}{\mathsf{n}}
\newcommand{\natCellAux}{\overline{\natTransAux}}
\newcommand{\altNat}{\mathsf{j}}
\newcommand{\altCell}{\overline{\altNat}}
\newcommand{\altaltNat}{\mathsf{m}}
\newcommand{\altaltCell}{\overline{\altaltNat}}

\newcommand{\phiTimes}{\Phi} 

\newcommand{\lin}{\mathcal{L}}
\newcommand{\prom}{\mathcal{P}}
\newcommand{\clone}{\mathbb{C}}
\newcommand{\restrictedPseudofunctor}[1]{\overline{#1}}
\newcommand{\restrictedPseudofunctorToProduct}[1]{\overline{#1}}
\newcommand{\nucleus}[1]{\overline{#1}}
\newcommand{\bicloneFromProducts}[1]{\mathrm{Bicl}({#1})}
\newcommand{\cloneFromProducts}[1]{\cloneOf{#1}}
\newcommand{\timesInSuper}{{\times\!}}
\newcommand{\timesArrowInSuper}{{\times}, {\to}}
\newcommand{\freeClone}[1]{\mathbb{FC}\mathrm{l}{\left(#1\right)}}
\newcommand{\freeCartClone}[1]{\mathbb{FC}\mathrm{l}^{\timesInSuper}{\left(#1\right)}}
\newcommand{\freeCartClosedClone}[1]{\mathbb{FC}\mathrm{l}^{\timesArrowInSuper}{\left(#1\right)}}
\newcommand{\altClone}{\mathbb{D}}
\newcommand{\cartClone}[2]{\left(#1,#2, {\Pi}_n(-) \right)}
\newcommand{\ccClone}[2]{\left(#1,#2, {\Pi}_n(-), {\exptype{}{}} \right)}
\newcommand{\biclone}{\mathcal{C}}
\newcommand{\freeBiclone}[1]{\mathcal{FC}\mathit{l}{\left(#1\right)}}
\newcommand{\freeBicat}[1]{\mathcal{FB}\mathit{ct}{(#1)}}
\newcommand{\freeCartBicat}[1]{\mathcal{FB}\mathit{ct}^{\timesInSuper}{(#1)}}
\newcommand{\freeCartClosedBicat}[1]{\mathcal{FB}\mathit{ct}^{\timesArrowInSuper}{(#1)}}
\newcommand{\freeCartBiclone}[1]
	{\mathcal{FC}\mathit{l}^{\timesInSuper}{\left(#1\right)}}
\newcommand{\freeCartClosedBiclone}[1]
	{\mathcal{FC}\mathit{l}^{\timesArrowInSuper}{\left(#1\right)}}
\newcommand{\altBiclone}{\mathcal{D}}
\newcommand{\synclonesymbol}{\mathrm{Syn}}
\newcommand{\synclone}[1]{\synclonesymbol(#1)}
\newcommand{\syncloneCart}[1]{\synclonesymbol^{\!{\times}\!}(#1)}
\newcommand{\syncloneCartClosed}[1]
	{\synclonesymbol^{\!{\times}, {\to}\!}(#1)}
\newcommand{\syncloneAtClosed}[1]
	{\termCatSymbol^{\text{@}, {\times}, {\to}\!}(#1)}

\newcommand{\hirclone}{\mathcal{H}}
\newcommand{\graph}{\mathcal{G}}

\newcommand{\nodes}[1]{{#1}_0}

\renewcommand{\cs}[2]{{#1} [ {#2} ]} 
\newcommand{\csbig}[2]{{#1} \big[ {#2} \big]} 
\newcommand{\uncs}[1]{\cslr{}{#1}} 
\newcommand{\uncssmall}[1]{\cs{}{#1}}
\newcommand{\cslr}[2]{{#1} {\left[ {#2} \right]}} 
\newcommand{\csthree}[3]{{#1}{\left[{#2}\right]\left[{#3}\right]}}
\newcommand{\csthreesmall}[3]{{#1}[{#2}][{#3}]}
\newcommand{\ms}[2]{{#1} \circ \seq{#2}}
\newcommand{\mslr}[2]{{#1} \circ \seqlr{#2}}
\newcommand{\msthree}[3]{({#1} \circ \seq{#2}) \circ \seq{#3}}
\newcommand{\msthreelr}[3]{\left({#1} \circ \seqlr{#2}\right) \circ \seqlr{#3}}
\newcommand{\clonetimes}{\boxtimes}
\newcommand{\bigclonetimes}{\bigboxtimes}
\newcommand{\pstrans}[1]{\widetilde{\mathsf{r}}(#1)}

\newcommand{\BE}[1]{\mathsf{BE}(#1)}

\renewcommand{\quote}{\mathsf{quote}}
\newcommand{\unquote}{\mathsf{unquote}}
\newcommand{\restr}[2]{ {{#1}{\big |}_{#2}} } 
\newcommand{\prodext}[2]{ \underline{#1} {#2} } 
\newcommand{\coerce}[1]{{\left[{#1}\right]}} 
\newcommand{\dterms}{\mathcal{K}}
\newcommand{\terms}{K} 
\newcommand{\langTerms}{L} 
\newcommand{\dlangTerms}{\mathcal{L}} 
\newcommand{\dneutTerms}{\mathcal{M}}
\newcommand{\dnormTerms}{\mathcal{N}}
\newcommand{\dvarTerms}{\mathcal{V}}
\newcommand{\neutTerms}{M}
\newcommand{\normTerms}{N}
\newcommand{\varTerms}{V}
\newcommand{\twoDisc}[1]{\mathrm{d}{#1}} 
\newcommand{\oneDisc}[1]{{\partial}{#1}} 

\newcommand{\context}{x_1 : A_1, \dots, x_n : A_n}

\newcommand{\hpi}[2]{\hcomp{\pi_{{#1}}}{#2}}
\newcommand{\varEnum}{\mathit{Var}}
\newcommand{\constrewr}{\kappa}
\newcommand{\prodop}{\textstyle{\prod}}
\newcommand{\altprod}{\mathrm{Prod}}

\newcommand{\evalterm}{\mathsf{eval}}
\newcommand{\heval}[1]{\hcomp{\evalterm}{#1}}
\newcommand{\genevalterm}[1]{\hcomp{\evalterm}{{#1}}}
\newcommand{\pair}[1]{\mathsf{tup}({#1})}
\newcommand{\pairlr}[1]{\mathsf{tup}{\left({#1}\right)}}
\newcommand{\pairbig}[1]{\mathsf{tup}{\big({#1}\big)}}
\newcommand{\pairName}{\mathsf{tup}}

\newcommand{\proj}[1]{{\varrho_{{#1}}}}
\newcommand{\indproj}[2]{{\varrho^{(#1)}_{{#2}}}}
\newcommand{\subid}[1]{\iota_{{#1}}}
\newcommand{\assoc}[1]{\mathsf{assoc}_{{#1}}}
\newcommand{\reduce}[1]{\mathsf{reduce}{\left(#1\right)}}
\newcommand{\reduceName}{\mathsf{reduce}}

\newcommand{\aeq}{=_\alpha}

\newcommand{\into}[1]{\llparenthesis \, #1  \, \rrparenthesis} 
\newcommand{\out}[1]{\overline{#1}}

\DeclareMathOperator{\binder}{ . }
\newcommand{\bind}{\binder}
\newcommand{\smallbind}{{\binder}}

\newcommand{\stlcTimes}{\Lambda^{\!{\times}\!}}
\newcommand{\stlc}{\Lambda^{\!{\times}, \to}}

\newcommand{\hir}{\mathrm{H}^{\mathrm{cl}}}

\newcommand{\lang}{\Lambda_{\mathrm{ps}}}
\newcommand{\langBiclone}{\lang^{\mathrm{bicl}}}
\newcommand{\langBicat}{\lang^{\mathrm{bicat}}} 
\newcommand{\langCartClosed}{\lang^{{\times}, {\to}}}
\newcommand{\langCart}{\lang^{{\times}}}

\newcommand{\dom}{\mathrm{dom}}

\newcommand{\fv}{\mathrm{fv}}

\newcommand{\sig}{\mathcal{S}}
\newcommand{\sigCat}[1]{{#1}\text{-}\mathrm{sig}}
\newcommand{\stlcTimesSigCat}{\sigCat{\stlcTimes}}
\newcommand{\stlcTimesUnSigCat}{\urestrict{\stlcTimesSigCat}}
\newcommand{\stlcSigCat}{\sigCat{\stlc}}
\newcommand{\stlcUnSigCat}{\urestrict{\stlcSigCat}}
\newcommand{\cartCatSigCat}{\sigCat{\langCart}}
\newcommand{\cartCatUnSigCat}{\urestrict{\cartCatSigCat}}
\newcommand{\langSigCat}{\sigCat{\langCartClosed}}
\newcommand{\langUnSigCat}{\urestrict{\langSigCat}}

\newcommand{\baseTypes}{\mathfrak{B}} 
\newcommand{\allTypes}[1]{\widetilde{#1}}

\newcommand{\un}{\mathsf{v}}
\newcommand{\co}{\mathsf{w}}
\newcommand{\unTimes}{\mathsf{u}^{{\times}}}
\newcommand{\coTimes}{\mathsf{c}^{{\times}}}
\newcommand{\unExp}{\mathsf{u}^{{\scriptsizeexpobj{}{}}}}
\newcommand{\coExp}{\mathsf{c}^{{\scriptsizeexpobj{}{}}}}

\newcommand{\exptype}[2]{\ensuremath{#1 \altTo #2}}

\newcommand{\rewrite}[2]{\ensuremath{#1 \To #2}}

\newcommand{\birewrite}[2]{\ensuremath{#1 \To #2}}
\newcommand{\app}[2]{\mathsf{app}{\left({#1},{#2}\right)}}
\newcommand{\lam}[2]{\lambda {#1} \dot {#2}}

\newcommand{\trans}[1]{{#1}^{\dagger}}
\newcommand{\altTrans}[1]{{#1}^\sharp}
\newcommand{\transExp}[1]{\mathsf{e}^{\dagger\!}({#1})}
\newcommand{\transExplr}[1]{\mathsf{e}^{\dagger\!}{\left({#1}\right)}}
\newcommand{\transTimes}[1]{\mathsf{p}^{\dagger\!}({#1})}
\newcommand{\transTimeslr}[1]{\mathsf{p}^{\dagger\!}{\left({#1}\right)}}

\newcommand{\etaTimes}[1]{\varsigma_{#1}}
\newcommand{\epsilonTimes}[1]{{\varpi_{{#1}}}}
\newcommand{\epsilonTimesInd}[2]{{\varpi^{(#1)}_{{#2}}}}
\newcommand{\counitCell}[2]{\mu^{(#1)}_{#2}}

\newcommand{\etaExp}[1]{\eta_{#1}}
\newcommand{\epsilonExpRewr}[1]{\epsilon_{#1}}
\newcommand{\genEpsilonExp}[1]{\beta_{#1}}

\newcommand{\epsilonExp}{\epsilon}
\newcommand{\wkn}[2]{\hcomp{{#1}}{\mathrm{inc}_{{#2}}}}

\newcommand{\horizCompSmall}[2]{ {#1}{\lbrace{#2}\rbrace} }
\newcommand{\horizComp}[2]{ {#1}{\left\lbrace{#2}\right\rbrace} }
\newcommand{\hcomp}[2]{\horizComp{#1}{#2}}
\newcommand{\hcompsmall}[2]{\horizCompSmall{#1}{#2}}
\newcommand{\hcompbig}[2]{ {#1}\big\{ {#2}\big\} }
\newcommand{\hcompthree}[3]{\horizComp{\horizComp{#1}{#2}}{#3}}
\newcommand{\hcompthreebigfst}[3]{\hcomp{{\hcompbig{#1}{#2}}}{{#3}}}

\DeclareMathOperator{\bulletop}{\bullet}
\renewcommand{\vert}{\bulletop}
\newcommand{\vertsub}[1]{\vert}

\newcommand{\nat}{\mathsf{nat}}

\newcommand{\phiCell}[2]{\phiTimes_{#2}}

\newcommand{\sub}[2]{\mathsf{sub}(#1; #2)} 
\newcommand{\subName}{\mathsf{sub}}
\newcommand{\cont}[2]{\mathsf{cont}(#1; #2)}
\newcommand{\contName}{\mathsf{cont}}

\newcommand{\post}[1]{\mathsf{post}(#1)}

\newcommand{\postName}{\mathsf{post}}
\newcommand{\fuse}{\mathsf{fuse}}
\newcommand{\swap}{\mathsf{swap}}
\newcommand{\push}[2]{\mathsf{push}(#1; #2)}
\newcommand{\pushName}{\mathsf{push}}

\newcommand{\unpack}{\mathsf{unpack}}

\newcommand{\plh}{%
  {\ooalign{$\phantom{0}$\cr\hidewidth$\scriptstyle\times$\cr}}%
}

\newcommand{\timessuper}{{\plh}}
\newcommand{\expsuper}{{\hspace{-.1mm}\scriptsizealtTo}}
\newcommand{\expsuperSmall}{{\hspace{-.1mm}\scriptsizescriptsizealtTo}}

\newcommand{\CCCTrans}{\alpha}
\newcommand{\CCCTransNat}{\overline{\CCCTrans}}
\newcommand{\CCCTransExp}{\CCCTrans^\expsuper}
\newcommand{\CCCTransExpSmall}{\CCCTrans^\expsuperSmall}

\newcommand{\CCCTransProd}{\CCCTrans^{\timessuper}}

\newcommand{\modif}{\Xi}
\newcommand{\altModif}{\Psi}

\newcommand{\ext}[1]{{{#1}^{{\#}}}}
\newcommand{\presBase}{\mathrm{q}}
\newcommand{\altPresBase}{\mathrm{u}}
\newcommand{\prodPres}{\presBase^{\timessuper}}
\newcommand{\altProdPres}{\altPresBase^{\timessuper}}
\newcommand{\expPres}{\presBase^\expsuper}
\newcommand{\altExpPres}{\altPresBase^\expsuper}

\renewcommand{\c}{c}

\newcommand{\inc}{\iota}

\newcommand{\squaresCat}[3]{\slice{#1}{#2}} 
\newcommand{\glued}[1]{\underline{#1}}
\newcommand{\glTup}[1]{ \{ #1 \}}
\newcommand{\glLam}{\glued{\mathsf{lam}}}
\newcommand{\lanext}[1]{{\seq{#1}}}  
\newcommand{\lanNat}{l}
\newcommand{\lanCell}{\overline{l}}

\newenvironment{td}[0]{\begin{center}\vspace{-0.0em}\begin{tikzcd}}{\end{tikzcd}\end{center}}
\newcommand{\equals}[1]{\overset{\small{\text{#1}}}{=}}
\newcommand{\twocell}[1]{\overset{\small{#1}}{\Leftarrow}}
\newcommand{\twocellRight}[1]{\overset{\small{#1}}{\Rightarrow}}
\newcommand{\twocellDown}[1]{\small{\Downarrow\!{#1}}}
\newcommand{\twocellUp}[1]{\small{\Uparrow\!{#1}}}
\newcommand{\twocellIso}[1]{\overset{\small{#1}}{\iso}}

\newcommand*{\corner}{\mbox{\LARGE{$\ulcorner$}}}
\newcommand*{\sidecorner}{\begin{rotate}{180}\corner\end{rotate}}

\newcommand*{\rightcorner}{\begin{rotate}{130}\corner\end{rotate}}

\newenvironment{rules}
{\begin{figure}[!h]
		\begin{mdframed}\centering
		}
		{\end{mdframed}
\end{figure}
}

\newcommand{\treeskip}{1em} 

\newenvironment{bprooftree}
  {\leavevmode\hbox\bgroup}
  {\DisplayProof\egroup}

\newcommand{\unaryRule}[3]{\begin{bprooftree}
\AxiomC{\ensuremath{#1}}
\RightLabel{\scriptsize #3}
\UnaryInfC{\ensuremath{#2}}
\end{bprooftree}\vspace{\treeskip}}

\newcommand{\binaryRule}[4]{\begin{bprooftree}
\AxiomC{\ensuremath{#1}}
\AxiomC{\ensuremath{#2}}
\RightLabel{\scriptsize #4}
\BinaryInfC{\ensuremath{#3}}
\end{bprooftree}\vspace{\treeskip}}

\newcommand{\trinaryRule}[5]{\begin{bprooftree}
\AxiomC{\ensuremath{#1}}
\AxiomC{\ensuremath{#2}}
\AxiomC{\ensuremath{#3}}
\RightLabel{\scriptsize #5}
\TrinaryInfC{\ensuremath{#4}}
\end{bprooftree}\vspace{\treeskip}}

\makeatletter
\renewcommand{\@marginparreset}{%
  \reset@font\footnotesize
  \footnote
  \raggedright
  \@setminipage
}
\makeatother

\usepackage{chngcntr}
\counterwithin{figure}{chapter}

\maxtocdepth{subsection}

\newcommand*{\nom}[2]{\nomenclature{\ensuremath{#1}}{#2}}

\renewcommand{\nompreamble}{\vspace{-1cm} With 
typing signature and page of first definition \vspace{2cm}}
\makeindex
\makenomenclature
\OnehalfSpacing

\author{Philip James Saville}
\title{Cartesian closed bicategories:\\ type theory and coherence}
\date{October 2019}

\begin{document}

\maketitle

\newgeometry{left=1.5in, right=1.5in}

\frontmatter

\thispagestyle{empty}
\section*{Abstract}
In this thesis I lift the Curry--Howard--Lambek correspondence between the 
simply-typed lambda calculus and cartesian closed categories
to the 
bicategorical setting, then use the resulting type theory to prove a coherence 
result for cartesian closed bicategories. 
Cartesian closed bicategories---\mbox{2-categories} `up to isomorphism' 
equipped with similarly weak products and exponentials---arise in 
logic, categorical algebra, and game semantics. However, calculations in such 
bicategories quickly fall into a quagmire of coherence data.  I show that 
there is at most one 2-cell between any parallel pair of 
1-cells in the free cartesian closed 
bicategory on a set and hence---in terms of the difficulty of 
calculating---bring the data of 
cartesian closed bicategories down to the familiar level of cartesian 
closed categories.

In fact, I prove this result in two ways. The first argument is closely related to Power's coherence theorem for bicategories with flexible bilimits. 
For the second, which is the central preoccupation of this thesis, 
the proof strategy has two parts: the construction of a type theory, and the 
proof that it satisfies a form of normalisation I call \Def{local coherence}.
I synthesise the type theory from algebraic principles using a novel 
generalisation of 
the (multisorted) abstract clones of universal algebra, called \Def{biclones}. 
The result brings together two extensions of the simply-typed lambda 
calculus: a 2-dimensional type theory in the style of Hilken, which 
encodes the 2-dimensional nature of a bicategory, and a 
version of explicit substitution, which encodes a composition 
operation that is 
only associative and unital up to isomorphism. 
For products and exponentials I develop the 
theory of cartesian and cartesian closed biclones and pursue a connection with 
the representable multicategories of Hermida.
Unlike preceding 2-categorical 
type theories, in which products and exponentials are encoded by postulating a 
unit and counit satisfying the triangle laws, the universal 
properties for products and exponentials are encoded using 
T. Fiore's biuniversal arrows. 

\enlargethispage{2\baselineskip}
Because the type theory is 
extracted from the construction of a free biclone, its syntactic model
satisfies a suitable 2-dimensional freeness universal property generalising the 
classical Curry--Howard--Lambek correspondence. One may therefore describe the 
type theory as an `internal language'. The relationship with the classical 
situation is made precise by a result establishing that the type theory I 
construct is the simply-typed lambda calculus up to isomorphism.

This relationship is exploited for the proof of local coherence. It is has been 
known for some time that one may 
use the normalisation-by-evaluation strategy to prove the simply-typed 
lambda calculus is strongly normalising. Using a bicategorical treatment of 
M. Fiore's categorical analysis of normalisation-by-evaluation, I prove a 
normalisation result which entails the coherence theorem for cartesian closed 
bicategories. In contrast to previous  coherence results for bicategories, the argument does not rely on 
the theory of rewriting or strictify using the Yoneda embedding. 
I prove bicategorical generalisations of a series 
of well-established category-theoretic results, present a notion of 
\Def{glueing of bicategories}, and bicategorify the folklore result providing 
sufficient conditions for a glueing category to be cartesian closed. Once these 
prerequisites have been met, the argument is remarkably similar to that in the 
categorical setting. 

\paragraph*{} 
This version of the thesis has been optimised for on-screen viewing: larger font, smaller margins, smaller page size. It works best when viewed with two pages side-by-side (`book view'). The content has not been changed.

\newpage
\thispagestyle{empty}
\section*{Declaration}
This dissertation is the result of my own work and includes nothing which is the outcome
of work done in collaboration except as declared in the Preface and specified in the text.
It is not substantially the same as any that I have submitted, or am concurrently
submitting, for a degree or diploma or other qualification at the University of Cambridge
or any other University or similar institution except as declared in the Preface and
specified in the text. I further state that no substantial part of my dissertation has
already been submitted, or is being concurrently submitted, for any such degree, diploma
or other qualification at the University of Cambridge or any other University or similar
institution except as declared in the Preface and specified in the text.
This dissertation does not exceed the prescribed limit of 60\,000 words.
\newpage

\thispagestyle{empty}
\section*{Acknowledgements}
First I have to thank my supervisor Marcelo Fiore. I have benefited from 
many hours of stimulating technical discussions, and owe an intellectual debt to his 
precise and thoughtful approach to problems. I am particularly grateful 
for the patient way he dealt with my (sometimes egregious) errors. I also owe thanks to Martin Hyland and Steve Awodey for examining this thesis, and to Andr{\'e} Joyal for suggesting Power's coherence proof for bicategories with finite bilimits could be adapted to cc-bicategories.
 
Thank you to the fellow occupants of FE14 and to the 
students of the PLS group, with whom I enjoyed many lunches and pub trips. Thank you 
especially to Ian Orton, Dylan McDermott, Hugo Paquet, Matthew Daggitt, 
and Michael Schaarschmidt, who 
all put up with me for three years, and to Alex Hickey for many chats over 
afternoon tea. Ian discovered a bug in an early version of the type theory that forms the first part of this thesis which, despite by panic at the time, greatly improved the end result. As with so many of my technical problems over the last few years, I was lucky to have Ian, Hugo and Dylan patiently spending their time to sanity-check my ideas and explain the basic concepts I was missing. 
Thank you also to Ohad Kammar for 
his forbearance and financial support as my intended submission date slipped back and back. 

I seem to have spent most of my time in Cambridge either doing mathematics or 
rowing. Sidney Sussex Boat Club was my pressure release valve and a source of 
great friendships, and I am incredibly grateful to everyone who made the club 
special. A special shout-out goes to `my' Lents 2019 crew. 

I cannot do justice to the friends and family who have helped me 
over the last four years: to all of you, thank you.

Finally, and most importantly, thank you to Aijing Wang for her love and care 
over the last four years. I wouldn't have done it without you.

\restoregeometry

\clearpage

\chapter*{Lay introduction}

This introduction is for the friends and family who have occasionally asked 
\emph{what it is I actually do}, and to whom I don't think I've ever managed a 
satisfactory answer. I hope this goes some way to explaining what the next 
200-odd pages are about.

Here's the three-sentence explanation. This thesis is about using 
\emph{category theory} and \emph{type theory} together to prove a 
\emph{coherence theorem}. I construct a 
type theory---a kind of mathematical language---to describe a 
category-theoretic structure which turns up in algebra and logic. Then, by 
proving a property of the type theory, I 
deduce the category-theoretic structure has a property called 
\emph{coherence}.  

Let's flesh that out a bit more.
Part I of the thesis is about 
\emph{syntax}, while Part II is about \emph{semantics}. The distinction 
between the two is one we are used to in our day-to-day lives. If you 
read a message from me and judge me for spelling `life' as `liffe', you are 
judging  the syntax: the string of symbols that make up the 
message. 
If you nonetheless grasped what I meant by the whole phrase `what have I been 
doing with my liffe', you understood the semantics: the meaning I was trying to 
convey. When a translator translates a sentence from English to Mandarin, they 
change the syntax (from Roman letters to Chinese characters), but maintain the 
semantics: a Chinese reader should finish the Chinese 
sentence understanding the same thing as an English reader who has just read 
the English sentence. 

The syntactic-semantic distinction is central to the study of programs and 
programming languages. On the syntactic side, there is the literal string of 
characters making up a program. If I write  
\texttt{print}(`hello world'), 
the computer has to break this up into the command (\texttt{print}) and the 
string 
that I'm telling it to print (\textit{hello world}), and act accordingly. If I 
write
$((3 + 6) \times 7)^2$, it has to break it up into the series of instructions
\begin{enumerate}
\item Add 3 to 6, \emph{then}
\item Multiply the result by 7, \emph{then}
\item Multiply this result by itself.
\end{enumerate}
Anyone who has sat down to write a program will know that a fair amount of time 
is spent chasing down the little syntactic mistakes (such as missing a crucial 
`;') that, as far as the computer is concerned, make what you have written 
unreadable.

Comparing programs only by their syntax is not very helpful, however. Here are 
three different programs 
that take in a number $x$ and give back another number:
\begin{equation} \label{eq:programs}
\frac{(\frac{x}{2} + 5) \times 6}{3}
\qquad\qquad 
(\frac{x}{2} + 5) \times 2
\qquad\qquad 
x + 10
\end{equation}
The string of symbols in each case is different, so syntactically they are 
different programs. But, as we learn in secondary 
school algebra, these all mean the same thing: they evaluate to the same 
answer. Intuitively, we can think of all these programs as the same. From the 
programmer's perspective, writing any one of these is as good as the other. So 
if the computer transforms between them (for example, because one of them is 
quicker to run), then the programmer doesn't care. But if the computer 
transforms one of these programs into $x + 1$, then they most certainly will. 

This suggests that we should study programming languages not just by thinking 
about the syntax, but by making precise our intuitive idea of what 
a program `says'. First we provide a mathematical description of what each 
part of a program means. For example, the 
command \texttt{add(2)(3)} `means' $2 + 3$. Then we say that 
two programs are the same if they have the same mathematical description. The 
idea is that the mathematics captures the meaning of the program (its 
semantics), and allows us to abstract away from its syntax. We can then prove 
all kinds of useful guarantees. For example, we can show that every 
syntactically correct program  
will eventually stop, and that the answer it will give is 
the one you would expect.

What does this have to do with category theory, type theory, or coherence? It 
turns out that type theory can be thought of as the logic of programs, and 
that category theory is one of the best ways of describing what these programs 
mean. 

Type theory grew up in the early 20th century in response to problems in logic, 
most famously Russell's paradox. One formulation of the paradox is this. 
Imagine you are a very organised person, and are constantly making lists: 
to-do lists, shopping lists, and so on. But one day you worry that you might be 
missing something, so you sit down to enumerate all the things that do not 
appear 
on any of your lists. Do you add \emph{this list} to this new list? If you do, 
it appears on a list, so shouldn't be on the list. If you don't, it doesn't 
appear on any list, so should be on the list. It seems neither choice is correct! The solution suggested by Russell is to 
stratify objects: at the first level are things that may appear in a list 
(things you need to do, food you need to buy), at the second level 
are lists of things in the first level, at the third level are lists of things 
at the second level, and so on. Every list has a level, and a list can 
only contain things at lower levels, so you never encounter the question of 
whether a list must contain the entry \emph{this list}. 

This kind of logic is governed by the principle that everything has a type, and a thing's 
type determines how it can behave. So you have a type of \emph{things that go 
in lists}, a type of \emph{lists of things that go in lists}, a type of lists 
of these lists, and so on. Similarly, you might have a type \texttt{nat} of natural (counting) numbers, and
the numbers $0, 1, \dots$ all have type \texttt{nat}. 
From this point of view, the expression $0 = 1$ is false, but expressions like 
$\frac{2}{0}$ or $\mathtt{print} + 2$ are ruled to be nonsense: the 
language of type theory simply doesn't allow you to form such expressions. With 
enough types and enough ways of forming new types, one 
can go a long way to formulating all of mathematics in a type theory.

This way of thinking has been absorbed into computer 
science as a way of structuring programs. When a programmer sits down to write 
a program, they have in mind some kind of input (say, a list of numbers) and an 
output (say, the highest number in the list). One can therefore think of a 
program as something that takes in something of some type, and gives out 
something of another type. For 
example, I can tell the computer that I want it to treat \texttt{add(2)(3)} as 
something of type \texttt{int}---as a whole number, obtained by adding 2 to 
3---or 
as something of type \texttt{string}---as a list of nine characters that happen 
to 
look like a command to add two numbers. If I declare \texttt{add(2)(3)} to be 
of 
type \texttt{string}, I can't treat it as a number: I can ask for its length 
(9), 
but can't multiply it by two. The more types you have, and the more 
constructions for new types you allow, the more precise you can make these 
restrictions.

Type theory, then, can be viewed in two ways. As a kind of logic, 
in which every true or false statement is attached to a type. Or as a 
programming language, in which the statements I can write down correspond to 
programs with a set input type and a set output type.

Thinking of programs as processes which take an input and return an output 
helps clarify the connection with category theory. Category theorists are 
mathematicians who truly believe that \emph{it's not about the destination, 
it's about the journey}. Instead of asking about particular objects, category 
theorists study the way things are 
related. The diagrams that you'll see if you flick through this thesis say 
exactly this: if you walk around the diagram following the arrows in one 
direction, and then walk around the diagram following the arrows in the other 
direction, the two walks will be equal. The 
fundamental idea is that, if I know all the ways to 
get into an object, and all the ways to get out of it, then I can discover 
everything I 
need to know. More than this: I can discover other, seemingly unrelated, 
objects that are related to the things around them in the same way. For 
example, the `if $\dots$ then' construction of logic, the collection of ways to 
assign an object of a set $B$ to every object of a set $A$, and the notion of 
group from algebra---which axiomatises the ways of rotating and reflecting 
shapes like triangles, squares, and cubes---are all examples of the same 
categorical construction. 

The categorical perspective has unearthed unexpected relationships between 
geometry, algebra, and logic, but it also plays an important role as a 
mathematical description for programming languages: category theory is the 
\emph{semantics} for the \emph{syntax} of type theory. For a type theorist, a 
program is a particular way of constructing objects of a certain type. For the 
category theorist, this is exactly a way of getting from one object (the input 
type) to another (the output type). Type theory and category theory are 
intertwined: by carefully choosing our categories, we can provide constructions 
that correspond exactly to the allowed type-theoretic expressions. By studying 
these categories, we can learn about type theory; by studying type theories, we 
can learn about their corresponding categories. Broadly speaking, this is the 
what I do in this thesis: I construct a type theory, show it corresponds to a 
special class of categories, and then---by proving something about the type 
theory---solve a problem about the class of categories.

The problem is called \emph{coherence}. The special categories I work 
with---the `cartesian closed bicategories' of the title---have 
uses in other areas of category theory, as well as in algebra and in the study 
of programming languages, but they are intricate. As well as the 
ways of getting from $A$ to $B$, they include the routes between these routes. Imagine 
$A$ and $B$ are Cambridge and Oxford. Then the routes between them might be 
walking directions for the various routes, and the routes-between-routes might 
be the ways you can change one set of directions into the other: change `left' 
for `right' at this junction, replace `100 yards' with `2 miles', and so on. 
Or you can imagine studying programs, and the ways of transforming them 
stage-by-stage into something that you can run in 0s and 1s on your hardware. 
In 
this example, you might have two programs with the same input type and the same 
output type---such as those in~(\ref{eq:programs}) above---and think about the 
ways of transforming 
one into another: replacing $\frac{y \times 6}{3}$ by $y \times 2$, and 
$\frac{x}{2} \times 2$ by just $x$, and so on. 

Precisely describing these two levels, and the ways they must interact, 
requires many 
axioms and many checks at every 
stage of a calculation. This quickly becomes tedious, and leads to 
proofs that are so 
long it is hard to check they are correct, let alone fit them onto a page so 
that they can be verified by the community. In this thesis I show 
that cartesian closed bicategories have the property that any equation you can 
write down for any cartesian closed bicategory (not relying on any special properties 
of a specific one) must hold. This means that those long tedious 
calculations are dramatically simplified: all those things that you had 
to check before are now guaranteed to hold by the theorem.

In Part I, then, I construct a type theory for describing cartesian 
closed bicategories. If a type theory is a logic for programs, this is a logic 
for programs \emph{and ways of transforming programs into one another}. I show 
that expressions in this type theory correspond exactly to data in any 
cartesian closed bicategory, so that a proof about the type theory is a proof 
about every cartesian bicategory. Then, in Part II, I prove a property 
of the type theory that guarantees that every cartesian closed bicategory is 
coherent. If you want to see what it all looks like, the type theory is in 
Appendix~\ref{chap:full-language}, and the big theorem is 
Theorem~\ref{thm:main-result-on-free-bicat}.

\clearpage

\tableofcontents

\mainmatter

\chapter{Introduction}


\section*{The Curry--Howard--Lambek correspondence and beyond}

The simply-typed lambda calculus lives a remarkable double life.
It can be seen 
as a term calculus for intuitionistic logic, or as the 
syntax of cartesian closed 
categories---a class of algebraic structures encompassing many important 
examples. This two-fold relationship, known as the 
\Def{Curry--Howard--Lambek correspondence}, is fundamental to the study of 
logic, type theory, and 
programming language theory. 

In this thesis we are largely concerned with the relationship between type 
theory and category theory. In the context of the simply-typed lambda calculus 
the crucial observation is due to Lambek~\cite{Lambek1980, Lambek1985}, who 
showed that the simply-typed lambda calculus may be interpreted in any 
cartesian closed category, that any cartesian closed category gives rise to a 
simply-typed lambda calculus, and moreover that these two operations are---in a 
suitable sense---mutually inverse. For a computer scientist, this
says that cartesian closed categories capture the meaning, or 
\Def{semantics}, of the simply-typed lambda calculus: to give a model of the 
simply-typed lambda calculus is to give a cartesian closed category. For a 
category theorist, this says that one 
may use the simply-typed lambda calculus as a convenient syntax or 
\Def{internal language} for constructing proofs in cartesian closed categories. 

The simply-typed lambda calculus is just the starting point. Internal languages 
are a key tool in topos theory~\cite{Makkai1977, Johnstone2002}, and there are 
well-known versions of Lambek's correspondence for linear 
logic~\cite{Benton1993}~(see~\eg~\cite{Mellies2009} for an overview) and 
Martin-L\"of type 
theory~\cite{Seely1984,Clairambault2014}. Meanwhile, categorical constructions 
such as monads have become standard for semantic descriptions of so-called 
`effectful programs', which display behaviours beyond merely 
computing some result~\cite{Moggi1989, Moggi1991}. 

Latent within each of these developments is the notion of \Def{reduction} or 
\Def{rewriting}.
In a Lambek-style semantics one begins with a type theory together with rules 
specifying how terms 
reduce to one another. These reduction rules generate an equational theory, and 
one identifies terms modulo this theory with morphisms 
in a suitable category. This is generally 
sufficient for type-theoretic applications, despite the loss of intensional 
information. To study the behaviour of reductions, however, this information 
must be retained.

One way to retain this information is 
through \Def{2-categories}. A 2-category consists of objects, 
morphisms, and \Def{2-cells} relating morphisms, subject to the usual unit and 
associativity laws. 
In the late 1980s multiple authors suggested 2-categories as a semantics for 
rewriting~(\eg~\cite{Rydeheard1987, 
Power1989}). In particular, Seely~\cite{Seely1987} sketched a connection 
between 2-categories equipped with a (lax) cartesian closed structure and the 
$\beta\eta$-rewriting rules of the simply-typed lambda calculus. In this model, 
$\eta$-expansion and $\beta$-reduction form the unit and counit of the 
adjunction defining 2-categorical cartesian closed structure. 
Hilken~\cite{Hilken1996} then 
took the identification between cartesian closed 2-categories and the rewriting 
theory of the simply-typed lambda calculus a step further by introducing a 
`2$\lambda$-calculus' consisting of types, terms, and \Def{rewrites} between 
terms. Syntactically, rewrites model reduction rules---for example, the 
$\beta\eta$-rules of the simply-typed lambda calculus---while semantically they 
play the role of 2-cells. 

Since Hilken's work, \Def{2-dimensional type theories} consisting of types, 
terms and rewrites have been employed for a range of applications, from 
rewriting 
theory~\cite{Hirschowitz2013} to 
the study of Martin-L\"of type theory and its connections to homotopy theory 
and higher category theory~(\eg~\cite{Garner2009,Licata2011, Licata2012}).
In this thesis I also connect 2-dimensional type theory to higher category 
theory, 
but with different aims. Here, the focus is on a class of higher categories 
of recent importance for applications in logic~\cite{FioreSpecies,Gambino2017,Olimpieri2020}, 
the 
semantics of programming 
languages~\cite{PaquetThesis}, and the study of category theory 
itself~\cite{Fiore2015,FioreOpetopicBonn} known as 
\Def{cartesian closed bicategories}. The copious data required to define a 
cartesian closed bicategory makes calculations within them a demanding 
undertaking: the aim of this thesis is to drastically reduce those demands. 

\section*{`The technical nightmares of bicategories'}

Suppose given a pair of spans 
$(A \leftarrow B \to C)$ and $(C \leftarrow D \to E)$ in a category with finite 
limits. By analogy with the category of sets, these could be thought of as 
`relations'
$A \rightsquigarrow C$ and $C \rightsquigarrow E$. How should the composite
$A \rightsquigarrow E$ be defined?
A natural suggestion is to take the pullback of $(B \to C \leftarrow D)$ and 
use the associated projection maps, thus:
\begin{td}[column sep = small, row sep = small]
\: & \: & B \times_C D
\arrow[phantom]{dd}[very near start, xshift=3mm, yshift=1mm]{\rightcorner}
\arrow{dl}
\arrow{dr} & \: & \: \\
\: & B \arrow{dl}\arrow{dr} & \: & D \arrow{dl}\arrow{dr} & \: \\
A & \: & C & \: & E 
\end{td} 
Because limits are only unique up to unique isomorphism, this 
definition does not satisfy the unit and associativity laws of a 2-category. 
However, such laws do hold up to specified isomorphism, and these isomorphisms 
satisfy coherence axioms. The resulting structure is called 
a \Def{bicategory}. Bicategories are rife in mathematics and theoretical computer science, arising
for instance in
algebra~\cite{Benabou1967,Street1995}, 
semantics of computation~\cite{CattaniFioreWinskel,Castellan2017}, 
datatype models~\cite{Abbott2003,Dagand2013}, 
categorical logic~\cite{FioreSpecies,Gambino2013}, and
categorical algebra~\cite{Fiore2015,Gambino2017,FGHW2017}.
More generally, one may (loosely) consider 
\Def{weak $n$-categories} to have
$k$-cells relating $(k-1)$-cells 
for $k=1, \dots, n$, such that the coherence axioms for $k$-cells are themselves 
witnessed by a specified $(k+1)$-cell. 

Weak higher category theory entails layers of complexity that do not exist at 
the 1-categorical level. 
Morphisms (more generally, $k$-cells) satisfying axioms up to some 
higher cell may exist in new relationships; specifying their behaviour leads to 
intimidating lists of axioms, for which the intuitive content is not 
immediately obvious.
Proofs become  
purgatorial exercises in drawing pasting diagram after pasting diagram, or 
diagram chases in which an intuitively-clear kernel is dominated by
endless structural isomorphisms shifting data back and forth. 	
Even at the level $k=2$, Lack---certainly a member of the higher-categorical 
\emph{cognoscenti}---refers to (strict) 2-category theory as a ``middle 
way'', avoiding ``some of the 
technical nightmares of bicategories''~\cite{Lack2010}. 

A small example highlights how the step from categories to bicategories 
blows up the length of a proof. Consider the following  
lemma, which is an elementary exercise in working with cartesian closed 
categories.

\begin{prooflesslemmalocal} \quad
\begin{enumerate}
\item Every object $X$ in a category with finite products $(\catC, \times, 1)$ 
has a canonical structure as a commutative comonoid, namely 
$\left( 1 \xla{!} X \xra{\Delta} X \times X\right)$.

\item Every endo-exponential $[\exp{X}{X}]$ in a cartesian closed category
$(\catC, \times, 1, {\exp{}{}})$ has a canonical structure as a monoid, namely 
\[
1 \xra{\Id_X} [\exp{X}{X}] \xla{\circ} [\exp{X}{X}] \times [\exp{X}{X}]
\]
\qedhere
\end{enumerate}
\end{prooflesslemmalocal}

Following the principle that higher 
categories behave in roughly the same manner as 1-categories so long as care 
is taken to specify the behaviour of the higher cells, one expects a 
version of this result to hold for cartesian closed bicategories. 
The bicategorical notion of monoid is called a 
\Def{pseudomonoid}~\cite{Day1997}. In a bicategory $\baseCat$ with finite 
products $(\times, 1)$, this is a structure
$(1 \xra{e} M \xla{m} M \times M)$
equipped with invertible 2-cells
$\alpha, \lambda$ and $\rho$ witnessing the categorical unit and 
associativity laws:
\begin{center}
\begin{tikzcd}
1 \times M 
\arrow[phantom]{dr}[description, xshift=2mm, yshift=1mm]{\twocellIso{\lambda}}
\arrow[bend right]{dr}[swap]{\simeq}
\arrow{r}{e \times M} &
M \times M 
\arrow{d}[description]{m} &
M \times 1
\arrow[phantom]{dl}[description, xshift=-2mm, yshift=1mm]{\twocellIso{\rho}}
\arrow[swap]{l}{M \times e} 
\arrow[bend left]{dl}{\simeq} \\
\: &
M &
\:
\end{tikzcd}
\:
\begin{tikzcd}
(M \times M) \times M 
\arrow[phantom]{drr}[description]{\twocellIso{\alpha}}
\arrow[swap]{d}{m \times M}
\arrow{r}{\simeq} &
M \times (M \times M) 
\arrow{r}{M \times m} &
M \times M 
\arrow{d}{m} \\
M \times M
\arrow[swap]{rr}{m} &
\: &
M
\end{tikzcd}
\end{center}
These 2-cells are required to satisfy two coherence laws, 
corresponding to the 
triangle and pentagon axioms for a monoidal category. Indeed, the prototypical 
example---obtained by instantiating the definition in $\Cat$---is of monoidal 
categories. Comparing with our categorical lemma suggests 
the following. 

\begin{myconjecturelocal} \quad \label{conj:nsr-objects}
\begin{enumerate}
\item Every object $X$ in a bicategory with finite products 
$(\baseCat, \times, 1)$ has a 
canonical structure as a 
commutative pseudocomonoid, with 1-dimensional 
structure  
$\left( 1 \xla{!} X \xra{\Delta} X \times X\right)$.

\item Every endo-exponential $[\exp{X}{X}]$ in a cartesian closed bicategory
$(\baseCat, \times, 1, {\exp{}{}})$ has a canonical structure as a 
pseudomonoid, with 1-dimensional structure 
\[
1 \xra{\Id_X} [\exp{X}{X}] \xla{\circ} [\exp{X}{X}] \times [\exp{X}{X}]
\]
\end{enumerate}
Moreover, in each case the 2-cells witnessing the 1-categorical axioms are 
canonical choices arising from the cartesian (closed) structure of $\baseCat$. 
\end{myconjecturelocal}

Constructing the witnessing 2-cells $\alpha, \lambda$ and $\rho$ is relatively 
straightforward: roughly 
speaking, one can translate each equality used in the categorical 
proof into a 2-cell, and then compose these together. The 
difficulty arises in checking the coherence laws, 
which entails a series of long diagram chases unfolding the properties of these 
composites. It is this extra work that makes bicategorical 
calculations more burdensome than their strict counterparts: it is not enough 
to merely witness the 
axioms---which corresponds to checking them in a strict setting---one 
must also check the 
witnesses are themselves \emph{coherent}. 

Not only do these checks entail extra work, they are often extremely tedious.
Generally one does not have to apply clever tricks or techniques, 
only plough through diagram chases until the result falls out. This is the 
case, for 
example, when one sits down to verify the coherence laws for 
Conjecture~\ref{conj:nsr-objects}. This leads to a false sense of security: it 
is tempting to believe 
that the coherence axioms `must' work out as expected, and that these extra 
checks may be omitted. As Power put 
it as long ago as 1989~\cite{Power1989bilimit}:
\begin{quotation}
The verification is almost always routine, and one's intuition is almost 
always 
vindicated; but to check the detail is often a very tedious job. Of course, 
one 
should still do it\dots[ignoring such details] can be dangerous, as 
illustrated in~\cite{Benabou1985}, because on rare occasions, one's intuition 
fails\dots
\end{quotation}

Despite these difficulties, higher categories---either as $\infty$-categories 
or as bicategories and tricategories---present an intuitively 
appealing 
and technically rich setting for studying phenomena arising throughout 
mathematics and theoretical computer science. Examples arise in 
topology~\cite{Leinster2004}, categorical logic~\cite{FioreSpecies}, 
categorical algebra~\cite{Benabou1967}, semantics 
of computation~\cite{Cattani1998}, and datatype semantics~\cite{Abbott2003}, to 
name but a few. The success of the `Australian school' of the 1970s and 1980s 
highlights especially the fruitfulness of studying
categorical constructions in the bicategorical setting 
(\eg~\cite{Street1972, 
Street1980, Blackwell1981}). 

One is, therefore, caught between interest and difficulty: one 
wants to be able to work in higher categories, but the technicalities of doing 
so are formidable. And the squeeze only becomes tighter as the structure 
becomes richer. The question then becomes: how can one construct a way out?

\subsection*{Coherence laws and coherence theorems}


One solution to the difficulties of working in a higher category is to develop 
a formal calculus that provides a 
pragmatic language for constructing and presenting proofs. 
In recent years there 
has been a great deal of work along these lines~(\eg~\cite{Riehl2017, 
Curien2019, Shulman2019a}), generally motivated 
by applications to $\infty$-categories~(although not always, 
see~\eg~\cite{Frey2019}). Much of the impetus stems from the connections 
between type theory, homotopy theory, and $\infty$-categories 
(\eg~\cite{Garner2009, Licata2011}), particularly the 
versions of Martin-L\"of type theory known as 
\Def{homotopy type theory} or \Def{univalent type 
theory}~(\eg~\cite{hottbook}). The 
type theory is generally strict---allowing for simpler reasoning---but  
satisfies an up-to-equivalence universal property interpreting it in the weak 
structure in question; this is analogous to the relationship between 
Martin-L\"of type theory with extensional identity types and locally cartesian 
closed categories~\cite{Clairambault2014}. 
A related strand of research is the development of 
\emph{computer-aided} systems such as Globular~\cite{globular}, which aim to 
provide interactive theorem-proving tools for certain weak $n$-categories. 

An alternative approach is to show that the weak structure in question is 
(weakly) equivalent to a strict structure: the so-called \emph{coherence} 
property. To paraphrase Jane Austen:
\begin{quotation}
It is a truth universally acknowledged, that a higher category in 
possession of a good structure, must be in want of a coherence theorem.
\end{quotation}
So long as equivalences are injective-on-cells in the appropriate 
sense, one 
can then parley this into a result proving that classes of diagrams always 
commute. 
Since Mac~Lane's first coherence theorem for monoidal 
categories, together with its pithy slogan \emph{all diagrams 
commute}~\cite{MacLane1963}, a cottage industry has sprung up proving 
coherence results in various forms
(notable examples include~\eg~\cite
{MacLane1985,
Power1989bilimit, 
Power1989monad, 
Joyal1993, 
Gordon1995}).
Coherence proofs often rely on the Yoneda embedding, which allows one 
to embed a weak structure (such as a bicategory) into a strict structure 
(such as the 2-category of $\Cat$-valued pseudofunctors), or on the 
sophisticated 
machinery of 2-dimensional universal algebra. Rewriting theory 
provides an alternative, syntactic, approach~(\eg~\cite{Houston2007, 
Forest2018}).

However, coherence turns out to be a subtle property. Certainly, one can not 
always show that \emph{all} diagrams commute: 
consider, for instance, the case of braided 
monoidal categories. In general, the dividing 
line between `coherent' and 
`non-coherent' definitions may not be where one would na{\"i}vely hope it to 
be, and the exact line can be surprising. Tricategories are not generally 
triequivalent to strict 
3-categories~\cite{Gordon1995}, and the tricategory $\Bicat$ is not 
triequivalent to the tricategory $\Gray$ of 2-categories, 2-functors, 
\emph{pseudo}natural transformations and modifications~\cite{Lack2007}. 

The 
difficulty, therefore, is twofold: first, to identify the boundaries between 
commutativity and its failure, and second, to prove that all diagrams within 
a conjectured boundary do in fact commute. 

\section*{Coherence for cartesian closed bicategories}

In this thesis I prove a coherence theorem for bicategories equipped with 
products and exponentials in an `up to 
equivalence' fashion. As far as I am aware, these were first studied 
in~\cite{Makkai1996}, and the coherence result I prove was first conjectured by 
Ouaknine~\cite{Ouaknine1997}. It is an unfortunate accident of terminology that 
there is no connection to the 
`cartesian bicategories' of 
Carboni \& Walters~\cite{Carboni1987, Carboni2008}, nor to the `closed 
cartesian bicategories' of Frey~\cite{Frey2019}. Precisely, the theorem is 
the following.

\begin{nonumthm}
The free cartesian closed bicategory on a set of 0-cells has at most one 2-cell 
between 
any parallel pair of 1-cells.
\hfill \qed
\end{nonumthm}

Note that this is a particularly concrete statement of coherence. In terms 
of Conjecture~\ref{conj:nsr-objects}, 
it goes further than showing that, once one has constructed 
witnessing 2-cells such as $\alpha, \lambda$ and $\rho$ using only the axioms 
of a cartesian closed 
bicategory, then the coherence laws will hold. The theorem  
also guarantees that there is a unique choice of witnessing 2-cells. 
Using this in tandem with a precise connection between the 2-cells of the free
cartesian closed bicategory and equality in the free cartesian closed category
(Section~\ref{sec:STLC-vs-pseudoSTLC}), 
we shall be able to show 
further that it suffices to calculate completely 1-categorically. 

This work was initially motivated by the difficulty of proving 
statements such as Conjecture~\ref{conj:nsr-objects} and the corresponding 
obstruction to the development of a theory of 
$\infty$-categories~\cite{FioreOpetopicBonn} in the cartesian closed 
bicategories of generalised species~\cite{FioreSpecies} and cartesian 
distributors~\cite{Fiore2015}. However, cartesian closed 
bicategories appear more widely, for example in categorical 
algebra~\cite{Gambino2017} and game semantics~\cite{Yamada2018, PaquetThesis}.

The strategy has two parts. First, I develop a type theory $\langCartClosed$ 
for cartesian closed 
bicategories and show that it satisfies a suitable 2-dimensional freeness 
property. This extends the classical Curry--Howard--Lambek correspondence 
to the 
bicategorical setting. 
The shape of the type theory follows the tradition of 
2-dimensional type theory instigated by Seely~\cite{Seely1987} and 
Hilken~\cite{Hilken1996}. The 
up-to-isomorphism nature of 
bicategorical composition is captured through an \emph{explicit substitution} 
operation (\cf~\cite{Abadi1989}).  
Second, I adapt the 
\emph{normalisation-by-evaluation} technique introduced by 
Berger \& Schwichtenberg~\cite{Berger1991} for proving normalisation of the 
simply-typed lambda calculus to extract the theorem above. Here I closely 
follow Fiore's categorical treatment of the proof~\cite{Fiore2002}.

Of course, for $\langCartClosed$ to be a type theory for cartesian closed 
bicategories, one must impose some constraints. I stipulate the following three
desiderata. 

\begin{quotation}
\textbf{Internal language.} The syntactic model of the type theory must 
be free, in an 
appropriately bicategorical sense. From a logical perspective, this corresponds 
to a soundness and completeness property. We shall not go so 
far as, say, constructing a triadjunction between a tricategory of signatures 
and the tricategory of cartesian closed bicategories. Instead, we prove strict 
universal properties~(\cf~\cite{Gurski2006}) wherever possible. As 
well as being readily verifiable, these properties are often easier to work 
with.
\end{quotation}

\begin{quotation}
\textbf{Relationship to STLC.} The type theory we construct must have the 
`flavour' of type theory. In particular, one should be able to recover the 
simply-typed lambda calculus (STLC) as some kind of fragment: following the 
intuition 
that cartesian closed bicategories are cartesian closed categories 
up-to-isomorphism, a corresponding property should relate 
the simply-typed lambda calculus to $\langCartClosed$. This also imposes 
restrictions on the form of judgements and derivations: they should be 
presented in a style recognisable as type theory. 
\end{quotation}

\begin{quotation}
\textbf{Usability.} This is connected to the preceding point. There is no gain 
in constructing a syntactic calculus that merely re-phrases the 
axioms of a cartesian closed bicategory. Instead, 
the type theory ought to 
be a reasonable tool for constructing proofs. Its equational 
theory ought to be kept small, and express requirements that are natural from 
the semantic perspective. 
\end{quotation}

These desiderata are not merely stylistic: they will play a key part in our 
eventual proof of coherence. The precise correspondence 
with the simply-typed lambda calculus, for example, will allow us to 
leverage the categorical arguments of~\cite{Fiore2002} in a particularly direct 
way. Moreover, they should also make the type theory amenable to deep embedding 
in proof assistants such as Agda~\cite{agda}, and to extension with further 
structure in future work.

\section*{Outline}

The thesis is in two parts. Part~\ref{part:lang} is devoted to the construction 
of $\langCartClosed$ and a proof of its free property. Part~\ref{part:nbe} 
covers the normalisation-by-evaluation proof.

In Chapter~\ref{chap:background} I present an overview of the basic theory 
of bicategories. Much of the theory is well-known, but I take the opportunity 
to develop it with a focus on T. Fiore's 
\Def{biuniversal arrows}~\cite[Chapter 9]{TFiore2006}. This bicategorification 
of universal arrows encompasses both biadjunctions and bilimits, and is 
particularly amenable to being translated into type theory. 

Chapter~\ref{chap:biclone-lang} constructs the core part of $\langCartClosed$, 
namely a type theory for mere bicategories. This type theory is synthesised 
from an algebraic description of bicategorical 
substitution, called a \Def{biclone}, which generalises 
the \Def{abstract clones} of universal algebra~(\eg~\cite{Cohn1981, CloneBookRef}). 
We also establish a coherence theorem for this fragment of the type theory, 
generalising the Mac~Lane-Par{\'e} coherence theorem for 
bicategories~\cite{MacLane1985}.

In Chapter~\ref{chap:fp-lang} we extend the type theory with finite products. 
We pursue a connection between the \Def{representable multicategories} of 
Hermida~\cite{Hermida2000}, introducing the notion of 
\Def{representable (bi)clone} and showing that it coincides with a 
notion of (bi)clone with cartesian structure. Thereafter we synthesise a type 
theory from the free such biclone, and show that its syntactic model is free. 

Chapter~\ref{chap:ccc-lang} follows a similar pattern: we define
cartesian closed biclones and extract a type theory from the construction 
of the free such. Establishing the free property for cc-bicategories throws up 
more complications than the preceding two chapters, so we spend some time over 
this. Thereafter we establish that the simply-typed lambda calculus embeds into 
$\langCartClosed$ and that, modulo the existence of invertible 
rewrites (2-cells), this restricts to a bijection on $\beta\eta$-equivalence 
classes of terms. We also observe that Power's coherence theorem for bicategories 
with flexible bilimits~\cite{Power1989bilimit} may be adapted to the case of 
cc-bicategories (Proposition~\ref{prop:yoneda-coherence-for-cc-bicats}).

In each of Chapters~\ref{chap:biclone-lang}--\ref{chap:ccc-lang}, 
the development is motivated by the 
construction of a version of the following diagram.  This provides a
technical statement of the intuitive fact that, in order to construct a 
type theory for cartesian or cartesian closed (bi)categories, it suffices to 
construct a type theory for the corresponding (bi)clones. As a slogan:
\emph{(bi)clones are the right intermediary between syntax and semantics}.
\begin{figure}[!h]
\centering
\begin{tikzpicture}
\node[align=center] (biclones) at (0,2) 
	{structured (bi)clones \\ {\scriptsize many-in one-out morphisms}};
\node[align=center] (signatures) at (-4,0) 
	{signatures};
\node[align=center] (unary-signatures) at (0,-2) 
	{unary signatures};
\node[align=center] (bicategories) at (4,0) 
	{structured (bi)categories \\ {\scriptsize one-in one-out morphisms}};

\draw[->, bend left=13] (signatures) to["{\scriptsize restriction}"] (unary-signatures);
\draw[draw=none] (signatures) -- (unary-signatures) 
	node[draw=none,fill=none,font=\scriptsize,midway] 
	{$\adjUp$};
\draw[->, bend left=13] (unary-signatures) to["{\scriptsize inclusion}"] (signatures);

\draw[->, bend left=13] (unary-signatures) to["{\scriptsize free}"] (bicategories);
\draw[draw=none] (bicategories) -- (unary-signatures) 
	node[draw=none,fill=none,font=\scriptsize,midway] 
	{$\adjUp$};
\draw[->, bend left=13] (bicategories) to (unary-signatures);

\draw[->, bend left=13] (biclones) to (signatures);
\draw[draw=none] (biclones) -- (signatures) 
	node[draw=none,fill=none,font=\scriptsize,midway] 
	{$\adjDown$};
\draw[->, bend left=13] (signatures) to["{\scriptsize free}"] (biclones);

\draw[->, bend left=13] (biclones) to["{\scriptsize restriction}"] (bicategories);
\draw[draw=none] (biclones) -- (bicategories) 
	node[draw=none,fill=none,font=\scriptsize,midway] 
	{$\adjDown$};
\draw[->, bend left=13] (bicategories) to (biclones);
\end{tikzpicture}
\end{figure}

We then move to the normalisation-by-evaluation proof. 
In Chapter~\ref{chap:calculations} we prove bicategorical correlates of three 
well-known facts about presheaf categories, namely:
\begin{enumerate}
\item Every presheaf category is complete,
\item Every presheaf category is cartesian closed,
\item For any presheaf $P$ and representable presheaf $\yon (X)$
on a small category with binary products, the 
exponential
$\altexp{\yon X}{P}$ is, up to isomorphism, the presheaf $P(- \times X)$.
\end{enumerate}
The reader willing to believe versions of these results for every 2-category 
$\Hom(\baseCat, \Cat)$ of $\Cat$-valued pseudofunctors may safely skip this 
chapter.

Chapter~\ref{chap:glueing} introduces the notion of 
\Def{glueing of bicategories} and establishes mild conditions for the glueing 
bicategory to be cartesian closed. In the 1-categorical setting, this implies 
the so-called \Def{fundamental lemma} of logical 
relations~\cite{Plotkin1973, Statman1985}. 

In Chapter~\ref{chap:nbe} we complete the proof of the main result 
via a bicategorical adaptation of Fiore's~\cite{Fiore2002}. 
Much of the apparatus required is contained in the preceding two chapters. 
Finally, Chapter~\ref{chap:conclusions} briefly lays out some 
applications and suggestions for further work.

Appendices~\ref{chap:summary:free-structures}--\ref{chap:full-language} 
contain an index of the 
bicategorical free 
constructions and 
syntactic models throughout this thesis, 
an overview of the cartesian 
closed structures we construct, and the complete set of rules for 
$\langCartClosed$ together with their semantic interpretation. 

\paragraph*{Previous publication.} The type theory $\langCartClosed$ was 
presented in the paper 
\emph{A type theory for cartesian closed bicategories}~\cite{LICS2019}. This is 
available online at
\url{https://ieeexplore.ieee.org/document/8785708}.

\section*{Contributions}

The most obvious contribution is the coherence theorem for cartesian closed 
bicategories. In fact, we prove this in three different ways: two closely-related arguments 
using the Yoneda lemma 
(Proposition~\ref{prop:yoneda-coherence-for-cc-bicats} 
and Theorem~\ref{thm:main-result-via-strictifying-first})
and the third by normalisation-by-evaluation
(Theorem~\ref{thm:main-result-on-free-bicat}).  
In each case the strategy is of interest in its own right. 
The arguments from the Yoneda argument extend Power's coherence argument for bicategories with 
flexible bilimits~\cite{Power1989bilimit} to closed structure for the first time. 
On the other hand, the normalisation-by-evaluation argument shows potential for 
further development.
First, it is plausible that, 
by further refining the normalisation-by-evaluation one would be able to extract a 
normalisation algorithm 
computing the canonical 2-cell between any given 1-cells in the free cartesian 
closed bicategory. Second,  the combination of syntactic and semantic methods 
employed here is a novel approach to 
proving higher-categorical coherence theorems (although Licata \& Harper have 
gone some way in this direction, using a groupoidal model to prove canonicity 
for their 2-dimensional type theory~\cite{Licata2012}). This approach may extend 
to situations where other proofs of coherence---employing either syntactic approaches 
or the apparatus of 2-dimensional universal algebra---are less successful.

From the type-theoretic perspective, I believe the view taken 
here---namely, that the appropriate mediator between syntax and semantics is 
some version of abstract clones---is a fruitful one. Indeed, the definition of 
the type theory $\langCartClosed$ follows automatically from the 
definition of cartesian closed biclones. As far as I am aware, this is the 
first attempt to construct a type theory describing higher categories from such 
universal-algebraic grounds, and the first to exploit the machinery of 
explicit substitution (although Curien's diagrammatic calculus 
for locally cartesian closed categories shows similar ideas~\cite{Curien1993}).


The theoretical development required for the normalisation 
proof---such as the work on bicategorical glueing in 
Chapter~\ref{chap:glueing}---lays important foundations for further work. For 
instance, the machinery of Part~\ref{part:nbe} is the groundwork for proving a 
conservative extension result for cartesian closed 
bicategories over 
bicategories with finite products in the style of~\cite{Lafont1987,Fiore2002sums}.

Finally, this thesis contains moderately detailed proofs of results that one 
would certainly \emph{expect} but I have not seen \emph{proved} in the 
literature, such 
as the cartesian closure of the 2-category $\Hom(\baseCat, \Cat)$ of 
$\Cat$-valued pseudofunctors, pseudonatural transformations and modifications. 
At the very least, I hope this saves others the work of reproducing the 
extensive calculations required.

\section*{Notation and prerequisites}
I have tried to keep the presentation self-contained and 
accessible to type theorists with a categorical bent, as well as to (higher) 
category theorists with less experience in type theory. I recap the bicategory 
theory we shall need, and do not employ 
any heavyweight results without proof. 
Similarly, I take the simply-typed lambda calculus and its semantics (as 
in~\eg~\cite{LambekAndScott, Crole1994}) as known, but do not assume
familiarity with strategies such as glueing or normalisation-by-evaluation.  
This occasionally requires recapitulating folklore or standard results, but I 
hope in these cases the presentation is original enough 
to be of interest in itself. 

I have attempted to generally (but not universally) maintain the following 
typographical conventions:
\begin{itemize}
\item Named 1-categories are written in $\mathrm{Roman}$ font (\eg~$\Set$); 
named higher categories are in $\mathbf{bold}$ font (\eg~$\Cat$). Arbitrary 
categories are written 
in blackboard bold $(\catC, \catD, \dots)$ and arbitrary bicategories in 
calligraphic font $(\baseCat, \altCat, \dots)$.  
\item 2-cells are denoted either by lower-case Greek letters 
$(\alpha, \beta, \tau, \sigma, \dots)$ or given suggestive names in 
$\mathsf{sans\text{-}serif}$~(\eg~$\pushName$). 
\end{itemize}
An index of notation covering most of the recurring 1- and 2-cells is on 
page~\pageref{nomenclature}.

I have also borrowed the convention of Troelstra \& 
Schwichtenberg~\cite{Troelstra2000} for denoting the 
end of environments. The end of a proof is marked by a 
white square $(\qedsymbol)$ and the end of a remark, definition or example by a 
black triangle $(\blacktriangleleft)$. 

%

%

\chapter{Bicategories, bilimits and biadjunctions} \label{chap:background}

This chapter introduces the basic theory of bicategories, bilimits and 
biadjoints. Much of the content is well-known, and many excellent overviews of 
the material are available~(\eg~\cite{Benabou1967, Street1980, Borceux1994, 
Street1995, Leinster2004}). The intention behind recapitulating it here is 
two-fold. 
Firstly, to fix notation. Second, to introduce concepts in a style  
that is convenient for later chapters. There are many equivalent ways of 
formulating basic notions such as adjunction, adjoint equivalence and 
universal arrow. In the categorical setting, translating between the various 
formulations is generally straightforward. Bicategorically, however, such 
translations can require extensive checking of coherence data. We avoid this by 
taking the most convenient definition for our purposes as primitive, and by focussing on the~\Def{biuniversal arrows}~of~\cite[Chapter 
9]{TFiore2006}. 
These 
capture both 
bicategorical limits and adjunctions---and thereby cartesian closed 
structure---in a uniform way. We therefore spend some time developing the 
theory of biuniversal arrows before showing how it specialises to standard results about bilimits and biadjunctions.   

\section{Bicategories}

The fundamental notion is that of a \Def{bicategory}, due to 
B{\'e}nabou~\cite{Benabou1967}. These structures often arise when one 
defines composition by a universal property. Such an operation will generally not be associative and unital up to equality, only up to some mediating isomorphisms. A classical example is the bicategory of spans over a category $\catC$ with pullbacks. The objects are those of $\catC$, the morphisms 
$A \rightsquigarrow B$ are spans $A \xla{f} X \xra{g} B$, and composition is given by pullback. 

\newpage
\begin{mydefn} 
A \Def{bicategory} $\baseCat$ consists of 
\begin{itemize} 
\item A class of objects $ob(\baseCat)$, 
\item For every $X, Y \in ob(\baseCat)$ a \Def{hom-category} $\big(\baseCat(X, Y), \vert, \id\big)$ with objects \Def{1-cells} $f : X \to Y$ and morphisms \Def{2-cells}~\mbox{$\alpha : f \To f' : X \to Y$}; composition of 2-cells is called \Def{vertical composition}, 
\item For every $X, Y, Z \in ob(\baseCat)$ an \Def{identity} functor 
\mbox{$\Id_X : \catOne \to \baseCat(X,X)$} 
(for $\catOne$ the terminal category) 
and a \Def{horizontal composition} functor 
\mbox{$\circ_{X,Y,Z} : \baseCat(Y, Z) \times \baseCat(X, Y) \to \baseCat(X, 
Z)$}, 
\item Invertible 2-cells
\vspace{-\parskip}
\begin{align*} 
\a_{h,g,f} : (h \circ g) \circ f &\To h \circ (g \circ f) : W \to Z \\ 
\l_f : \Id_X \circ f &\To f : W \to X \\ 
\r_g : g \circ \Id_X &\To g : X \to Y
\end{align*} 
for every $f : W \to X$, $g : X \to Y$ and $h : Y \to Z$, natural in each of 
their arguments and satisfying a \Def{triangle law} and a \Def{pentagon law} 
analogous to those for monoidal categories:
\begin{center}
\begin{tikzcd}[column sep = -1em]
\big((k \circ h) \circ g\big) \circ f 
\arrow[swap]{d}{\a_{k \circ h, g, f}}
\arrow{rr}{\a_{k,h,g} \circ f} &
\: &
\big( k \circ (h \circ g) \big) \circ f
\arrow{d}{\a_{k, h \circ g, f}} &
\: \\
(k \circ h) \circ (g \circ f)
\arrow[swap]{dr}{\a_{k, h, g \circ f}} &
\: &
k \circ \big( (h \circ g) \circ f\big) 
\arrow{dl}{k \circ \a_{h,g,f}} \\
\: &
k \circ \big( h \circ (g \circ f)\big)  &
\:
\end{tikzcd} 
\quad
\begin{tikzcd}[column sep = -0em]
(g \circ \Id_X) \circ f 
\arrow[swap]{dr}{\r_g \circ f}
\arrow{rr}{\a_{g, \Id, f}} &
\: &
g \circ (\Id_X \circ f) 
\arrow{dl}{g \circ \l_f} \\
\: &
g \circ f &
\:
\end{tikzcd}
\end{center}
\end{itemize}
The functorality of horizontal composition gives rise to the so-called 
\Def{interchange law}: for suitable 2-cells $\tau, \tau', \sigma, \sigma'$ we 
have \mbox{$(\tau' \vert \tau) \circ (\sigma' \vert \sigma) = (\tau' \circ 
\sigma') \vert (\tau \circ \sigma)$}. 
\end{mydefn}

\begin{mynotation} \label{not:whiskering}
In the preceding we employ the standard notation for the 
\Def{whiskering} operations.
For a 1-cell $f : X \to Y$ and 2-cells $\sigma : h \To h' : W \to X$ and $\tau 
: g \To g' : Y \to Z$ we write 
$f \circ \sigma$
and
$\tau \circ f$
for $\id_f \circ \sigma : f \circ h \To f \circ h'$ and 
$\tau \circ \id_f : g \circ f \To g' \circ f$, respectively. 
\end{mynotation} 

The category $\mathrm{Rel}$ of sets and relations may be viewed as a 
\Def{locally posetal} bicategory---\ie~a bicategory in which each hom-category is a poset---by stipulating that 
$R \leq S : A \to B$ if and only if $a R b$ implies $a S b$ for all $a \in A$ and $b \in B$. A relation $R : A \to B$  is equivalently a map 
$A \times B \to \{ 0, 1\}$. Replacing sets by categories, one obtains the bicategory 
$\mathbf{Prof}$: this has objects categories, 1-cells 
$\catC \nrightarrow \catD$ the functors $\op\catD \times \catC \to \Set$, and 2-cells natural transformations. The identity on $\catC$ is the hom-functor 
$\Hom(-, =)$, and composition is given using the universal property of a 
presheaf category~(see~\eg~\cite{Benabou2000}).

\begin{myremark} \label{rem:bicats-coherent} 
The coherence theorem for monoidal categories~\cite[Chapter VII]{cfwm} 
generalises to bicategories: any bicategory is biequivalent to a 
2-category~\cite{MacLane1985}. Loosely speaking, then, any diagram constructed 
from only the identity and the structural constraints $\a, \l, \r$ with the 
operations of horizontal and vertical composition must commute 
(see~\cite{Leinster2004} for a readable summary of the argument). We are 
therefore justified in treating $\a, \l$ and $\r$ as though they were the 
identity, and we will sometimes denote such 2-cells merely by $\iso$. 
\end{myremark} 
Every bicategory $\baseCat$ has three duals. Following the notation 
of~\cite[\S 1.6]{Lack2010}, these are
\begin{itemize}
\item $\op{\baseCat}$, obtained by reversing the 1-cells, 
\item $\baseCat^{\mathrm{co}}$, obtained by reversing the 2-cells,
\item $\baseCat^\mathrm{coop}$, obtained by reversing both.
\end{itemize}
We call the first option the \Def{opposite bicategory}. This is the only form 
of dual we shall employ in this thesis. 

A morphism of bicategories is called a \Def{pseudofunctor} (or 
\Def{homomorphism})~\cite{Benabou1967}. It is a mapping on objects, 1-cells and 
\mbox{2-cells} that preserves horizontal composition up to isomorphism. 
Vertical composition is preserved strictly.

\begin{mydefn} 
A \Def{pseudofunctor} $F: \baseCat \to \altCat$ between bicategories $\baseCat$ and $\altCat$ consists of 
\begin{itemize} 
\item A mapping $F : ob(\baseCat) \to ob(\altCat)$, 
\item A functor \mbox{$F_{X,Y} : \baseCat(X,Y) \to \altCat(FX, FY)$} for every $X,Y \in ob(\baseCat)$, 
\item An invertible 2-cell \mbox{$\psi_X : \Id_{FX} \To F(\Id_X)$} for every \mbox{$X \in ob(\baseCat)$}, 
\item An invertible 2-cell $\phi_{f,g} : F(f) \circ F(g) \To F(f \circ g)$ for every $g : X \to Y$ and $f : Y \to Z$, natural in $f$ and $g$,
\end{itemize} 
subject to two unit laws and an associativity law: 
\begin{center}
\begin{tikzcd}[column sep = 3.5em]
\Id_{FX'} \circ Ff \arrow{r}{\psi_{X'} \circ Ff} 
\arrow[swap]{d}{\l_{Ff}} &
F(\Id_{X'}) \circ F(f) \arrow{d}{\phi_{\Id_{X'}, f}} \\
Ff &
F(\Id_{X'} \circ f) \arrow{l}{F\l_f}
\end{tikzcd}
\quad\quad
\begin{tikzcd}[column sep = 3.5em]
Ff \circ \Id_{FX} \arrow{r}{F(f) \circ \psi_{X}} 
\arrow[swap]{d}{\r_{Ff}} &
F(f) \circ F(\Id_{X}) \arrow{d}{\phi_{f, \Id_{X}}} \\
Ff &
F(f \circ \Id_X) \arrow{l}{F\r_f}
\end{tikzcd}
\end{center}
\begin{td}[column sep = 4.5em, row sep = 1em]
\big(Fh \circ Fg\big) \circ Ff 
\arrow{r}{\a_{Fh, Fg, Ff}} 
\arrow[swap]{d}{\phi_{h,g} \circ Ff} &
Fh \circ \big( Fg \circ Ff \big) 
\arrow{r}{F(h) \circ \phi_{g,h}} &
F(h) \circ F(g \circ f)
\arrow{d}{\phi_{h, g \circ f}} \\
F(h \circ g) \circ Ff 
\arrow[swap]{r}{\phi_{h \circ g, f}} &
F\big( (h \circ g) \circ f\big) 
\arrow[swap]{r}{F\a_{h,g,f}} &
F\big( h \circ (g \circ f) \big) 
\end{td}
A pseudofunctor for which $\psi$ and $\phi$ are both the identity is 
called \Def{strict}. 
\end{mydefn} 

We often abuse notation by leaving $\psi$ and $\phi$ 
implicit when denoting a pseudofunctor.

\begin{myexmp}  \quad
\begin{enumerate}
\item A monoidal category is equivalently a \mbox{one-object} bicategory; a monoidal functor is equivalently a pseudofunctor between one-object bicategories,
\item A 2-category is equivalently a bicategory in which $\a, \l$ and $\r$ are all the identity. A strict pseudofunctor $F : \baseCat \to \altCat$ between 2-categories $\baseCat$ and $\altCat$ is equivalently a 2-functor.
\item For every \emph{locally small}  
bicategory $\baseCat$ (see~Notation~\ref{not:size-issues}) 
and $X \in \baseCat$ there exists the 
\Def{Yoneda pseudofunctor} $\Yon X : \baseCat \to \Cat$, defined by
$\Yon X := \baseCat(X, -)$. The 2-cells $\phi$ and $\psi$ are structural isomorphisms.
\qedhere
\end{enumerate} 
\end{myexmp} 

Morphisms of pseudofunctors are called \Def{pseudonatural 
transformations}~\cite{Gray1974}. These are \mbox{2-natural} transformations 
($\CatCat$-enriched natural transformations) in which every naturality square 
commutes up to a specified 2-cell. Morphisms of pseudonatural transformations 
are called~\Def{modifications}~\cite{Benabou1967, Street1980}. 

\begin{mydefn}
A \Def{pseudonatural transformation} $(\natTrans, \natCell) : F \To G : 
\baseCat \to \altCat$ between pseudofunctors 
$(F, \psi^F, \phi^F)$ and 
$(G, \psi^G, \phi^G)$ consists of the 
following data:
\begin{enumerate} 
\item A 1-cell $\natTrans_X : FX \to GX$ for every $X \in \baseCat$, 
\item An invertible 2-cell $\natCell_f : \natTrans_Y \circ Ff \To Gf \circ 
\natTrans_X : FX \to GY$ for every $f : X \to Y$ in 
$\baseCat$, natural in $f$ and satisfying the following unit and associativity 
laws for every $X \in \baseCat$, $f : X' \to X''$ and $g : X \to X'$ in 
$\baseCat$. : 
\begin{center}
\begin{tikzcd}[column sep = -0.5em, row sep = 1.8em]
\: &
(Gf \circ \natTrans_{X'}) \circ Fg 
\arrow{dr}{\a_{Gf, \natTrans, Fg}}  &
\: \\
(\natTrans_{X''} \circ Ff) \circ Fg
\arrow[swap]{d}{\a_{\natTrans, Ff, Fg}}
\arrow{ur}{\natCell_f \circ Fg} &
\:  &
Gf \circ (\natTrans_{X'} \circ Fg) 
\arrow{d}{G(f) \circ \natCell_g} \\
\natTrans_{X''} \circ (Ff \circ Fg) 
\arrow[swap]{d}{\natTrans_{X''} \circ \phi^F_{f,g}} &
\: &
Gf \circ (Gg \circ \natTrans_X) 
\arrow{d}{\a^{-1}_{Gf, Gg, \natTrans}} \\
\natTrans_{X''} \circ F(f \circ g) 
\arrow[swap]{dr}{\natCell_{fg}} &
\: &
(Gf \circ Gg) \circ \natTrans_X 
\arrow{dl}{\phi^G_{f,g} \circ \natTrans_X} \\
\: &
G(f \circ g) \circ \natTrans_X
\end{tikzcd}
\:\:
\begin{tikzcd}[column sep = 0.8em, row sep = 2em]
\: &
\natTrans_X \arrow{dr}{\l^{-1}_{\natTrans}} &
\: \\
\natTrans_X \circ \Id_{FX} 
\arrow[swap]{d}{\natTrans_X \circ \psi^F_X}
\arrow{ur}{\r_{\natTrans}} &
\: &
\Id_{GX} \circ \natTrans_X 
\arrow{d}{\psi^G_X \circ \natTrans_X} \\
\natTrans_X \circ F\Id_X 
\arrow[swap]{rr}{\natCell_{\Id_X}}  &
\: &
G\Id_X \circ \natTrans_X  
\end{tikzcd}
\end{center} 
\end{enumerate}
A pseudonatural transformation for which every $\natCell_f$ is the identity is 
called \Def{strict} or \Def{2-natural}.  
\end{mydefn}

\begin{myremark} \label{rem:2-cell-orientation}
Note that we orient the 2-cells of a pseudonatural transformation as in the 
following diagram:
\begin{td}
FX 
\arrow[phantom]{dr}[description]{\twocell{\natCell_f}}
\arrow[swap]{d}{\natTrans_X}
\arrow{r}{Ff} &
FY \arrow{d}{\natTrans_Y} \\

GX \arrow[swap]{r}{Gf} &
GY
\end{td}
This is the reverse of~\cite{Leinster1998} but follows the direction 
of~\cite{Benabou1967, Street1980}. Of course, since we require each 
$\natCell_f$ to be invertible, the two choices are equivalent.
\end{myremark}

\begin{mydefn}
A \Def{modification} $\modif : (\natTrans, \natCell) \to (\altNat, \altCell)$ 
between pseudonatural transformations $(\natTrans, \natCell), (\altNat, 
\altCell) : F \To G : \baseCat \to \altCat$ is a family of 2-cells $\modif_X 
: \natTrans_X \To \altNat_X$, such that the following commutes for every 
$f : X \to X'$ in 
$\baseCat$:\footnote{Leinster~\cite{Leinster2004} requires both the above 
coherence law and 
that the family of 2-cells $\modif_X$ be natural in $X$; this appears to be an 
oversight, as neither Leinster's own~\cite{Leinster1998} nor 
Street's~\cite{Street1995} mention naturality.}
\begin{td}
\natTrans_{X'} \circ Ff \arrow{r}{\natCell_f}
\arrow[swap]{d}{\modif_{X'} \circ Ff} &
Gf \circ \natTrans_{X} \arrow{d}{Gf \circ \modif_X} \\
\altNat_{X'}\circ Ff \arrow[swap]{r}{\altCell_f} &
Gf \circ \altNat_X
\end{td}
\end{mydefn} 

\begin{myexmp} \label{e:ProductBicat} 
For every pair of bicategories $\baseCat$ and $\altCat$ there exists a 
bicategory $\Hom(\baseCat, \altCat)$ of pseudofunctors, pseudonatural 
transformations and modifications. If $\altCat$ is a 2-category, so is 
$\Hom(\baseCat, \altCat)$. In particular, for every bicategory $\baseCat$ there 
exists a 2-category $\Hom(\baseCat, \Cat)$, which one might view as a 
bicategorical version of the covariant presheaf category $\Set^{\mathbb{C}}$. 
Where $\catC$ is a mere category, pseudofunctors $\catC \to \Cat$ are called 
\Def{indexed categories}~\cite{MacLane1985}. 
\end{myexmp} 

Bicategories, pseudofunctors, pseudonatural transformations and modifications organise themselves into a \emph{tricategory}~(weak 3-category, see~\cite{Gordon1995, Gurski2006, Gurski2013}) we denote $\Bicat$~\cite{Gordon1995}. 

\begin{mynotation} \label{not:size-issues}
A bicategory $\baseCat$ (resp. pseudofunctor $F$) is said to be \Def{locally $P$} if the property $P$ holds for each hom-category 
$\baseCat(X, Y)$ (resp. functor $F_{X,Y}$). 
In particular, a bicategory is \Def{locally small} if 
every hom-category is a set, and \Def{small} if it is locally small and its 
class of objects is a set. We shall use $\Cat$ to denote the 
2-category of small categories and stipulate that, whenever we write
$\Hom(\baseCat, \Cat)$, then it is assumed that $\baseCat$ is small.
As usual, such issues can be avoided using technical devices such as 
Groethendieck universes (see~\eg~\cite{Shulman2008}).
\end{mynotation}

The \emph{bicategorical Yoneda Lemma} takes the following form, due to 
Street~\cite{Street1980}.\footnote{The bicategorical Yoneda Lemma is an example 
of a result that one would certainly expect to hold---and is generally only 
ever stated in the literature---but for which the proof actually requires a 
significant amount of work: see \cite{Bakovic} for the gory details.}

\begin{prooflesslemma}
For any bicategory $\baseCat$ and pseudofunctor $F : \baseCat \to \Cat$, 
evaluating at the identity for each $B \in \baseCat$ provides the components 
$\Hom(\baseCat, \Cat)\big(\baseCat(B, -), F\big) \xra{\simeq} FB$ 
of an equivalence in $\Hom(\baseCat, \Cat)$. 
Hence, the Yoneda pseudofunctor 
$\Yon : \baseCat \to \Hom(\baseCat, \Cat) : X \mapsto \baseCat(X, -)$ 
is locally an equivalence.
\end{prooflesslemma}

Bicategories provide a convenient setting for abstractly describing many 
categorical concepts~(\eg~\cite{Lawvere2017}); this perspective that has been 
used to particular effect by the Australian school~(see for 
instance~\mbox{\cite{Lack2012, Lack2014}}). The following definition is a small 
example of this general phenomenon. 

\begin{mydefn} \label{def:AdjunctionInBicat} Let $\baseCat$ be a bicategory. 
\begin{enumerate} 
\item An \Def{adjunction} $(A,B, f, g, \un, \co)$ in $\baseCat$ is a pair of objects $(A, B)$ with arrows $f : A \leftrightarrows B : g$ and 2-cells \mbox{$\un : \Id_A \To g \circ f$} and $\co : f \circ g \To \Id_B$ such that the bicategorical triangle laws hold (\eg~\cite{Gurski2012}): 
\begin{center} 
\begin{tikzcd} 
f \arrow{r}{\r_f^{-1}} \arrow[equals]{d} & 
f \circ \Id_X \arrow{r}{f \circ \un} & 
f \circ ( g \circ f) \arrow{d}{\a^{-1}_{f,g,f}} \\ 
f & 
\Id_Y \circ f \arrow{l}{\l_f} & 
(f \circ g) \circ f \arrow{l}{\co \circ f} 
\end{tikzcd} 
\qquad\qquad
\begin{tikzcd} 
g \arrow{r}{\l_g^{-1}} \arrow[equals]{d} & 
\Id_Y \circ g \arrow{r}{\un \circ g} & 
(g \circ f ) \circ g \arrow{d}{\a_{g,f,g}} \\ 
g & 
g \circ \Id_X \arrow{l}{\r_g} & 
g \circ (f \circ g) \arrow{l}{g \circ \co} 
\end{tikzcd} 
\end{center} 
\item An \Def{equivalence} $(A, B, f, g, \un, \co)$ in $\baseCat$ is a 
pair of objects $(A,B)$ with arrows \mbox{$f : A \leftrightarrows B : g$} and 
invertible 2-cells $\un : \Id_A \XRA{\iso} g \circ f$ and $\co : f \circ 
g \XRA{\iso} \Id_B$,

\item An \Def{adjoint equivalence} is an adjunction that is also an 
equivalence.  
\end{enumerate}
If 1-cells $f$ and $g$ are part of an equivalence, we refer to $g$ as the 
\Def{pseudoinverse} of $f$. Pseudoinverses are unique up to invertible 2-cell. 
\end{mydefn} 

In $\Cat$, an (adjoint) equivalence is exactly an (adjoint) equivalence of 
categories. Moreover, just as in $\Cat$, every equivalence induces an adjoint 
equivalence with the same 1-cells (see~\eg~\cite{Leinster1998}).
The appropriate notion of equivalence between bicategories is called 
\Def{biequivalence}~\cite{Street1980}. 

\begin{mydefn} 
A \Def{biequivalence} between bicategories $\baseCat$ and $\altCat$ consists of 
pseudofunctors $F : \baseCat \leftrightarrows \altCat : G$ and chosen 
equivalences $G \circ F \simeq \id_\baseCat$ and $F \circ G \simeq \id_\altCat$ 
in the bicategories $\Hom(\baseCat, \baseCat)$ and $\Hom(\altCat, \altCat)$, 
respectively. 
\end{mydefn} 

By a result of Gurski~\cite{Gurski2012}, one may assume without loss of 
generality that a biequivalence is an \Def{adjoint biequivalence}, in which 
$F$ and $G$ also form a \Def{biadjunction} (see 
Definition~\ref{def:biadjunction-universal-arrow}). 

\begin{mynotation}
Following standard practice from $\Cat$,
we shall sometimes refer to a pair of arrows $f : A \leftrightarrows B : g$ as
an \emph{(adjoint) equivalence}, leaving the 2-cells implicit. When we wish 
to emphasise that these 2-cells are given as data, we refer to a 
\emph{chosen} or \emph{specified} equivalence.

Similarly, we may sometimes leave most of 
the data implicit and refer to the pseudofunctor 
$F$ on its own as a biequivalence. Unlike the 1-categorical case, however, we 
shall always assume this biequivalence to be chosen.
\end{mynotation}


\begin{myexmp} \quad
\begin{enumerate}
\item A biequivalence between one-object bicategories is exactly an equivalence 
of monoidal categories (that is, an equivalence in the 2-category $\MonCat$ of 
monoidal categories, monoidal functors and monoidal natural transformations). 
\item $\mathbf{Prof}$ is biequivalent to its opposite bicategory~\cite[Section 
7]{Day1997} 
(\cf~the fact that the category $\mathrm{Rel}$ is isomorphic to its opposite).
\qedhere
\end{enumerate}
\end{myexmp} 

Loosely speaking, an equivalence of categories relates objects that are the 
same up to isomorphism, and a biequivalence of bicategories relates objects 
that 
are the same up to equivalence. Indeed, since every pseudofunctor preserves 
(adjoint) equivalences, an 
(adjoint) 
equivalence $A \simeq B$ in a bicategory $\baseCat$ induces an (adjoint) 
equivalence $\baseCat(A, -) \simeq \baseCat(B, -)$ in 
$\Hom(\op\baseCat, \Cat)$ and hence an (adjoint) equivalence
$\baseCat(A, X) \simeq \baseCat(B, X)$ for every $X \in \baseCat$. 
One consequence is that, if the pseudofunctor $F : \baseCat \to 
\altCat$ is a biequivalence, then 
\begin{enumerate}
\item For every $C \in \altCat$ there exists an object $B \in \baseCat$ and an equivalence $C \simeq FB$,
\item $F$ is \Def{locally an equivalence}: for every $B, B' \in \baseCat$ the functor $F_{B,B'}$ is part of an equivalence of categories $\baseCat(B,B') \simeq \altCat(FB, FB')$; in particular, every $F_{B, B'}$ is fully faithful and essentially surjective.
\end{enumerate}
In the presence of the Axiom of Choice, this formulation is equivalent to the 
definition given above~(\eg~\cite[Proposition 1.5.13]{Leinster2004}).

In the categorical setting it is elementary to check that a natural 
isomorphism---as an iso in a functor category---is exactly a natural 
transformation for which every component is invertible. The bicategorical 
version of this result is the following.

\newpage
\begin{mylemma}  \label{lem:equivalence-of-pseudofunctors-pointwise}
Let $F, G : \baseCat \to \altCat$ be pseudofunctors and 
suppose $(\natTrans, \natCell) : F \To G$ is a pseudonatural transformation 
such that every $\natTrans_X : FX \to GX$ is part of a specified adjoint equivalence 
$(\natTrans_X, \psinv{\natTrans}_X, 
\co_X : \psinv{\natTrans}_X \circ \natTrans_X \To \Id_{FX},
\un_X : \Id_{FX} \To \natTrans_X \circ \psinv{\natTrans}_X)$. Then:
\begin{enumerate}
\item The family of 1-cells 
$\psinv{\natTrans}_X : GX \to FX$ are the components of a pseudonatural 
transformation $(\psinv{\natTrans}, \psinv{\natCell}) : G \To F$, where for 
$f : X \to Y$ the 2-cell
$\psinv{\natCell}_f$ is defined by commutativity of the following diagram:
\begin{td}[column sep = 6em]
\psinv{\natTrans}_{Y} \circ Gf
\arrow[swap]{d}{\iso}
\arrow{r}{\psinv{\natCell}_f} &
Ff \circ {\psinv{\natTrans}_X} \\

\psinv{\natTrans}_{Y} \circ (Gf \circ \Id_{GX})
\arrow[swap]{d}{\psinv{\natTrans}_{Y} \circ Gf \circ \un_X} &
\Id_{FY} \circ (Ff \circ {\psinv{\natTrans}_X}) 
\arrow[swap]{u}{\iso} \\

\psinv{\natTrans}_{Y} 
	\circ \left(Gf 
					\circ (\natTrans_X \circ \psinv{\natTrans}_X)\right)
\arrow[swap]{d}{\iso}  &
(\psinv{\natTrans}_{Y} \circ \natTrans_{Y}) 
	\circ (Ff \circ {\psinv{\natTrans}_X})
\arrow[swap]{u}{\co_{Y} \circ Ff \circ \psinv{\natTrans}_X} \\

\psinv{\natTrans}_{Y} 
	\circ \left((Gf \circ \natTrans_X) \circ \psinv{\natTrans}_X\right)
\arrow[swap]{r}[yshift=-0mm]
	{\psinv{\natTrans}_{Y} \circ \natCell_f^{-1} \circ \psinv{\natTrans}_X} &
\psinv{\natTrans}_{Y} 
	\circ \left((\natTrans_Y \circ Ff) \circ \psinv{\natTrans}_X\right)
\arrow[swap]{u}{\iso}
\end{td}

\item The pseudonatural transformations 
$(\natTrans, \natCell) : 
F \leftrightarrows G : 
(\psinv{\natTrans}, \psinv{\natCell})$
are the 1-cells of an equivalence $F \simeq G$ in $\Hom(\baseCat, \altCat)$. 
\end{enumerate}
\begin{proof}
To see that $(\psinv{\natTrans}, \psinv{\natCell})$ is a 
pseudonatural transformation, the naturality and the unit laws
follow from the corresponding laws for $\natCell_f$. For the associativity law 
the process is similar, except one also applies the triangle law 
relating $\un$ and $\co$.

For the second claim we construct invertible modifications 
$(\psinv{\natTrans}, \psinv{\natCell}) \circ (\natTrans, \natCell) \iso 
\Id_{F}$ and 
$\Id_{G} \iso 
(\natTrans, \natCell) \circ (\psinv{\natTrans}, \psinv{\natCell})$. The obvious 
choices  for the components are 
$\co_X : \psinv{\natTrans}_X \circ \natTrans_X \To \Id_{FX}$ and 
$\un_X : \Id_{GX} \To \natTrans_X \circ \psinv{\natTrans}_X$.
It remains to check the modification axiom. To this end, observe that for 
every $f : X \to Y$ in $\baseCat$,
is the composite
\[
(\psinv{\natTrans}_Y \circ \natTrans_Y) \circ Ff 
\XRA{\co_Y \circ Ff} 
\Id_{FY} \circ Ff 
\XRA{\iso}
Ff \circ \Id_{FX} 
\XRA{Ff \circ \co_X^{-1}}
Ff \circ (\psinv{\natTrans}_X \circ \natTrans_X)
\] 
Similarly, $\overline{(\natTrans \circ \psinv{\natTrans})}_f$ is the composite
\[
(\natTrans_Y \circ \psinv{\natTrans}_Y) \circ Gf 
\XRA{\un_Y^{-1} \circ Gf}
\Id_{GY} \circ Gf 
\XRA{\iso} 
Gf \circ \Id_{GX} 
\XRA{Gf \circ \un_X}
Gf \circ (\natTrans_Y \circ \psinv{\natTrans}_Y)
\]
One then sees that
\begin{td}[column sep = 3em]
(\psinv{\natTrans}_Y \circ \natTrans_Y) \circ Ff 
\arrow[
swap,
rounded corners,
to path=
{ -- ([xshift=1ex]\tikztostart.west)
-| ([xshift=-0.8cm]\tikztotarget.west)
-- ([xshift=-0.8cm]\tikztotarget.west)
-- (\tikztotarget.west)}, 
]{ddd}
\arrow[phantom]{ddd}[font=\scriptsize, xshift=-2.3cm, description]
	{\overline{(\psinv{\natTrans}_Y \circ \natTrans_Y)}_f}
\arrow{r}{\co_Y \circ Ff}
\arrow[swap]{d}{\co_Y \circ Ff} &
\Id_{FY} \circ Ff 
\arrow[equals]{dl}{} 
\arrow{ddd}{\iso} \\

\Id_{FY} \circ Ff 
\arrow[swap]{d}{\iso} &
\: \\

Ff \circ \Id_{FX}
\arrow[swap]{d}{Ff \circ \co_X^{-1}}
\arrow[equals]{dr}{} &
\: \\

Ff \circ (\psinv{\natTrans}_X \circ \natTrans_X) 
\arrow[swap]{r}{Ff \circ \co_X} &
Ff \circ \Id_{FX}
\end{td}
so that $(\co_X)_{X \in \baseCat}$ does indeed form a modification. The proof 
for $\un$ is similar.
\end{proof}
\end{mylemma}

This lemma is particularly useful when it comes to constructing a 
biequivalence: to construct an equivalence $F \circ G \simeq \id$ it suffices 
to construct a pseudonatural transformation for which each component is an 
equivalence. 

The lemma also justifies the following terminology. We call a pseudonatural 
transformation $(\natTrans, \natCell)$ a \Def{pseudonatural equivalence} if 
every component 
$\natTrans_X$ is an equivalence, and a 
\Def{pseudonatural isomorphism} if every $\natTrans_X$ is invertible.

\section{Biuniversal arrows} \label{sec:biuniversal-arrows}

In his famous textbook~\cite{cfwm}, Mac Lane makes precise the 
notion of universal property by introducing \Def{universal arrows}. 
The Yoneda Lemma, 
limits and adjunctions are then all 
characterised in these terms. We adopt a similar approach, focussing on T. 
Fiore's \emph{biuniversal arrows}~\cite{TFiore2006}.  As well as providing a 
uniform way to describe bicategorical limits and bicategorical 
adjunctions, this perspective is 
particularly amenable to syntactic description. Biuniversal 
arrows are 
fundamental to the type theoretic description of bicategorical products and 
exponentials we shall see in Chapters~\ref{chap:fp-lang} 
and~\ref{chap:ccc-lang}.

A detailed development of the relationship between biuniversal arrows and 
biadjoints, 
complete with proofs, is available in~\cite[Chapter 
9]{TFiore2006}. The other results in what follows are implicit in much 
historical work on bicategory theory (\eg~\cite{Street1980}), but---as far as I 
am aware---have not previously been collected together in this form.

We begin by recapitulating the notion of universal arrow and its bicategorical counterpart.

\begin{mydefn}
Let $F : \catB \to \catC$ be a functor and $C \in \catC$. A 
\Def{universal arrow from $F$ to $C$} 
is a pair $(R \in \catB, u : FR \to C)$ such that, for 
any $B \in \catB$ and $f : FB \to C$, there exists a unique $\trans{f} : B \to 
R$ such that $u \circ F\trans{f} = f$. 
\end{mydefn}

It is an exercise to show that every universal arrow  $(R, u)$ from $F$ 
to $C$ is equivalently a chosen family of natural isomorphisms 
$\catB(-, R) \iso \catC(F(-), C)$, 
or---equivalently again---a terminal object in the 
comma category $(F \downarrow C)$. It follows that a right adjoint to 
$F : \catB \to \catC$ is specified by a choice of universal arrow 
$\epsilon_C : FUC \to C$ for 
every $C \in \catC$. The mapping $U$ extends to a functor with $Uf := 
\trans{(f \circ \epsilon_{C})}$ for $f : C \to C'$. The counit is then 
$\epsilon$ and the unit $\eta$ arises by applying the universal property to the 
identity: $\eta_B := \trans{(\id_{FB})} : B \to UFB$. If both $\epsilon$ and 
$\eta$ are invertible, the result is an adjoint equivalence.

To define biuniversal arrows, one weakens the isomorphisms defining a 
universal arrow to equivalences. We take particular care in 
choosing how we spell these out. It is generally 
convenient to require adjoint equivalences; by the 
well-known lifting theorem (\eg~\cite[Proposition 1.5.7]{Leinster2004}) 
this entails no loss of generality, while providing a more structured object to 
work with. We also go beyond T.~Fiore's definition by requiring that each 
adjoint equivalence is determined by a choice of universal 
arrow. 

\begin{mydefn}[{\cf~\cite{TFiore2006}}] \label{def:biuniversal-arrow}
Let $F : \baseCat \to \altCat$ be a pseudofunctor and $C \in \altCat$. A 
\Def{biuniversal arrow from $F$ to $C$} consists of a pair 
$(R \in \baseCat, u : FR \to C)$ 
and, for every $B \in \baseCat$, a chosen adjoint equivalence of categories
\begin{align*} 
\baseCat(B, R) &\xra{\simeq} \altCat(FB, C) \\
(B \xra{h} R) &\mapsto (FB \xra{Fh} FR \xra{u} C) 
\end{align*}
specified by choosing a family of invertible 
universal 2-cells as the counit. 

Explicitly, a biuniversal arrow from $F$ to $C$ 
consists of the following data:
\begin{itemize}
\item A pair $(R \in \baseCat, u : FR \to C)$, 
\item For every $B \in \baseCat$ and $h : FB \to C$, a map 
$\psi_B(h) : B \to R$ 
and an invertible 2-cell $\epsilon_{B,h} : u \circ F\psi_B(h) \To h$, universal 
in the sense that for any map $f : B \to R$ and 2-cell 
$\tau : u \circ Ff \To h$ there exists a 2-cell 
$\trans{\tau} : f \To \psi_B(h)$, unique such that
\begin{equation} \label{eq:universal-arrow-UMP}
\begin{tikzcd}
FB \arrow[bend right=60, swap]{drr}{h} 
\arrow[phantom]{r}[description, font=\scriptsize]{\twocellDown{F\trans{\tau}}} 
\arrow[bend left]{r}{Ff} \arrow[bend right, swap]{r}{F\psi_B(h)} &
FR \arrow[phantom]{d}[description, font=\small, near end, yshift=-2mm]{{\Downarrow\!\epsilon_{B,h}}} \arrow{dr}{u} \\
\: &
\: &
C
\end{tikzcd}
\quad = \quad
\begin{tikzcd}
\: &
FR \arrow[phantom]{d}[description, font=\small]{{\Downarrow\!\tau}} \arrow{dr}{u} &
\: \\
FB \arrow[swap]{rr}{h} \arrow{ur}{Ff} &
\: &
C
\end{tikzcd}
\end{equation}
\end{itemize}
such that the 2-cell 
$\trans{(\id_{u \circ Ff})} : f \To \psi_B(u \circ Ff)$
is invertible for every $f : B \to R$. \qedhere
\end{mydefn}

\begin{myremark}
Pictorial representations such as~(\ref{eq:universal-arrow-UMP})
are known as
\Def{pasting diagrams}. It is a consequence of the coherence theorem for 
bicategories that, once a choice of bracketing is made for the source and target 1-cells, 
a pasting diagram identifies a unique 2-cell~(\mbox{\cf~\cite[Remark~3.1.16]{Gurski2006}}; for a detailed exposition, 
see~\cite[Appendix~A]{Verity1992}).
\end{myremark}

On the face of it, a biuniversal arrow is only local structure: the data 
imposes a requirement on each hom-category, but no global
constraints. This property will be particularly useful for our later work 
synthesising a type theory, where we shall encode bicategorical products and 
exponentials as biuniversal arrows. Global structure arises in 
the following 
way~\mbox{(\cf~\cite[III.2]{cfwm})}.

\begin{mylemma} \label{lem:biuniversal-arrow-equivalence}
Let $F : \baseCat \to \altCat$ be a pseudofunctor and $C \in \altCat$. There 
exists a biuniversal arrow $(R, u)$ from $F$ to $C$ if and only if there exists 
an equivalence of pseudofunctors $\baseCat(-, R) \simeq \altCat(F(-), C)$ in 
$\Hom(\op{\baseCat}, \Cat)$, 
\begin{proof}
For every equivalence of pseudofunctors $\baseCat(-, R) \xra{\gamma} 
\altCat(F(-), C)$ one obtains from the Yoneda Lemma an arrow 
$\gamma_R(\Id_R) : FR \to C$. 
This arrow is biuniversal: indeed, the image of 
$\gamma_{R}(\Id_R)$ under the pseudofunctor 
$\altCat(FR,C) \to 
	\Hom(\op{\baseCat}, \Cat)\big(\baseCat(-, R), \altCat(F(-), C)\big)$ given 
by the Yoneda Lemma is 
isomorphic to $\gamma$, and hence an equivalence. The converse 
is~\cite[Theorem~9.5]{TFiore2006}.
\end{proof}
\end{mylemma}

\begin{myremark}
In Chapter~\ref{chap:glueing} we shall see that a biuniversal arrow 
from $F : \baseCat \to \altCat$ to $C \in \altCat$ is equivalently a 
terminal object in the bicategorical comma category 
$(F \downarrow \constPseudofunctor_C)$, for 
$\constPseudofunctor_C$ the constant pseudofunctor at $C$. 
\end{myremark}

\paragraph{Elementary properties of biuniversal arrows.} 
Many standard properties of universal arrows---such as those 
in~\cite{cfwm}---extend to biuniversal 
arrows.
Biuniversal arrows are 
unique up to equivalence, and the $\trans{(-)}$ operation preserves both 
invertibility and naturality. 

\begin{mynotation}
In the next lemma, and throughout, we shall abuse notation by writing just $\iso$ 
for the invertible 2-cell filling a square. Unless marked otherwise, it is assumed
this 2-cell is oriented right-to-left (\cf~Remark~\ref{rem:2-cell-orientation}).
\end{mynotation}

\begin{prooflesslemma}[{\cite[Lemma 9.7]{TFiore2006}}] \label{lem:biuniversal-arrows-unique}
Let $F : \baseCat \to \altCat$ be a pseudofunctor and $C \in \altCat$. For any two biuniversal arrows $(R, u)$ and $(R', u')$ from $F$ to $C$ there exists an equivalence $e : R \to R'$ and an invertible 2-cell $\kappa$ filling
\begin{td}
FR \arrow{r}{u} \arrow[swap]{d}{Fe} \arrow[phantom]{dr}[description]{\twocellIso{\kappa}} &
C \arrow[equals]{d} \\
FR' \arrow[swap]{r}{u'} &
C
\end{td}
Moreover, for any other pair 
$(f : R \to R', \lambda : u' \circ Fe \XRA{\iso} u)$ 
filling the above diagram, $e$ and $f$ are isomorphic via $\trans{\lambda}$.
\end{prooflesslemma}

It follows from the essential uniqueness of equivalences that, 
if $u : FR \to C$ is a biuniversal arrow from $F$ to $C$ and 
$u' \iso u$, then $u'$ is also a biuniversal arrow from $F$ to $C$. The next lemma 
follows from further standard facts about adjoint equivalences of categories.

\begin{prooflesslemma} \label{lem:trans-natural-and-invertible}
Let $F : \baseCat \to \altCat$ be a pseudofunctor and $(R, u)$ a biuniversal 
arrow from $F$ to $C \in \altCat$. For every object $B \in \baseCat$, 
\begin{enumerate}
\item If $f : B \to R$ is any morphism and $\alpha : u \circ Ff \To h$ is 
invertible, then so is $\trans{\alpha}$.
\item If the 1-cells $h, h' : FB \to C$ and $f, f' : B \to R$ and 2-cells $\alpha : u \circ Ff \To h$ and $\beta : u \circ Ff' \To h'$ are related by a commutative diagram of 2-cells as on the left below
\begin{center}
\begin{tikzcd}
u \circ Ff \arrow[swap]{d}{u \circ F\sigma} \arrow{r}{\alpha_f} & 
h \arrow{d}{\tau} \\
u \circ Ff' \arrow[swap]{r}{\alpha_{f'}} &
h'
\end{tikzcd}
\qquad\qquad\qquad
\begin{tikzcd}[column sep=2.5em]
f  \arrow{r}{\trans{(\alpha_f)}} \arrow[swap]{d}{\sigma} & 
\psi_B(h) \arrow{d}{\psi_B(\tau)} \\
f' \arrow[swap]{r}{\trans{(\alpha_{f'})}} &
\psi_B(h')
\end{tikzcd}
\end{center}
then the diagram on the right above commutes. In particular, if
$\alpha : u \circ F(-) \To \id_{\altCat(FB,C)}$ is a natural transformation, 
then so is $\trans\alpha : \id_{\baseCat(B, R)} \To \psi_B(-)$. \qedhere
\end{enumerate}
\end{prooflesslemma}

It is sometimes convenient, for example when working with bilimits, to work with the notion of \emph{birepresentable pseudofunctor}. 

\begin{mydefn}[{\cite{Street1980}}]
Let $F : \baseCat \to \Cat$ be a pseudofunctor. A \Def{birepresentation} $(R, \rho)$ for $F$ consists of an object $R \in \baseCat$ and an equivalence $\rho : \baseCat(R, -) \xra{\simeq} H$ in $\Hom(\baseCat, \Cat)$. 
\end{mydefn}

Representable functors $F : \baseCat \to \Set$ correspond to universal arrows from the terminal object to $F$. Similarly, to relate biuniversal arrows to birepresentable functors we employ the dual notion of a biuniversal arrow from an object to a pseudofunctor.

\begin{mylemma}[{\cf~\cite[Proposition III.2.2]{cfwm}}]
A pseudofunctor $F : \baseCat \to \Cat$ is birepresentable if and only if there exists a biuniversal arrow from the terminal category $\catOne$ to $F$.
\begin{proof}
It is certainly the case that $\Cat(\catOne, F(-)) \simeq F$ in $\Hom(\baseCat, 
\Cat)$. From birepresentability and the closure of equivalences under 
composition one obtains
$\Cat(\catOne, F(-)) \simeq F \simeq \baseCat(R, -)$, so the result follows 
from Lemma~\ref{lem:biuniversal-arrow-equivalence}.
\end{proof}
\end{mylemma}

\subsection{Preservation of biuniversal arrows}  \label{sec:preservation-of-biuniversal-arrows} 

Preservation of biuniversal arrows will provide a 
systematic way to define preservation of bilimits and preservation of 
biadjoints. We begin by 
examining preservation of universal arrows. Using the fact that a right adjoint 
to $F : \catB \to \catC$ is completely specified by a choice of 
universal arrow
$\left( UC, F(UC) \to C \right)$
for each $C \in \catC$---namely, the counit---it is 
reasonable to define morphisms of 
universal arrows analogously to morphisms of 
adjunctions~\cite[Chapter IV]{cfwm}.

\begin{mydefn} \label{def:morphism-of-universal-arrows}
Let $F : \catB \to \catC$ and $F' : \catB' \to \catC'$ be functors 
and suppose $(R, u)$ is a universal arrow from $F$ to $C \in \catC$. A pair 
of functors $(K, L)$ \Def{preserves the universal arrow $(R, u)$} if the 
following 
diagram commutes
\begin{td}
\catB \arrow{r}{F} \arrow[swap]{d}{L} & 
\catC \arrow{d}{K} \\
\catB' \arrow[swap]{r}{F'} &
\catC'
\end{td}
and $F'LR = KFR \xra{Ku} KC$ is a universal arrow from $F'$ to $KR$.
\end{mydefn}

Equivalently, we ask that the functor $(F \downarrow C) \to (F' \downarrow KC)$ 
defined by $(B, h : FB \to C) \mapsto (LB, F'LB = KFB \xra{Kh} KC)$ preserves 
the terminal object. This is a slight weakening of the definition of 
transformation of adjunctions given in~\cite{cfwm}: Mac Lane asks that the 
unit (or counit) be \emph{strictly} preserved. 

The bicategorical translation is as one would expect.

\begin{mydefn} \label{def:preservation-of-biuniversal-arrows}
Let $F : \baseCat \to \altCat$ and $F' : \baseCat' \to \altCat'$ be 
pseudofunctors and suppose $(R, u)$ is a biuniversal arrow from $F$ to $C \in 
\altCat$. A triple of pseudofunctors and pseudonatural transformations 
$(K, L, \rho)$ as in the diagram 
\begin{equation} \label{eq:pres-of-biuniversal-arrows}
\begin{tikzcd}
\baseCat \arrow{r}{F} \arrow[swap]{d}{L} \arrow[phantom]{dr}[description]{\twocellRight{\rho}} & 
\altCat \arrow{d}{K} \\
\baseCat' \arrow[swap]{r}{F'} &
\altCat'
\end{tikzcd}
\end{equation}
\Def{preserves the biuniversal arrow $(R, u)$} if $F'LR \xra{\rho_R} KFR 
\xra{Ku} 
KC$ is a biuniversal arrow from $F'$ to $KC$.
\end{mydefn}

By Lemma~\ref{lem:biuniversal-arrow-equivalence}, if $(K, L, \rho)$ preserves 
the universal arrow $(R, u)$ as in~(\ref{eq:pres-of-biuniversal-arrows}) then 
one obtains a pseudonatural family of equivalences 
$\baseCat'(B', LR) \simeq \altCat'(F'B', KC)$. 

Just as an equivalence of categories preserves all `categorical' properties, so 
a biequivalence preserves all `bicategorical' properties. In particular, a 
biequivalence preserves all biuniversal arrows.

\begin{mylemma} \label{lem:biequivalences-preserve-biuniversal-arrows}
Let $H : \altCat \to \altaltCat$ be a biequivalence and 
$F : \baseCat \to \altCat$ be a pseudofunctor. If $(R, u)$ is a biuniversal 
arrow from $F$ to 
$C \in \altCat$, then 
$Hu$ is a biuniversal arrow from $HF$ to $HX$. Hence, the triple
$(H, \id_{\baseCat}, \id)$ preserves the biuniversal arrow.
\begin{proof}
Since $H$ is locally an equivalence, for every $B \in \baseCat$ 
there exists a composite adjoint 
equivalence of categories 
$\baseCat(B, R) \simeq \altCat(FB, C) 
				\overset{H_{FB,C}}{\simeq} \altaltCat(HFB, HC)$ 
taking $h : B \to R$ to $H(u \circ Fh)$. 
Since ${H(u) \circ HF(-)}$ is naturally isomorphic to 
this adjoint equivalence, it is an adjoint equivalence itself.
\end{proof}
\end{mylemma}

There are two ways of formulating that a functor $F$ 
preserves limits: one can either ask that the image of the terminal cone is 
also a terminal cone, or that the canonical map  $F(\lim H) \to \lim(FH)$ 
is an isomorphism. Similar considerations apply to preservation of biuniversal 
arrows.

\begin{mylemma} \label{lem:preserves-cones-equivalent-to-equivalence}
Consider a square of pseudofunctors $K, L, F, F'$ related by a pseudonatural 
transformation 
$(\rho, \overline{\rho}) : KF \To F'L$ 
as in~(\ref{eq:pres-of-biuniversal-arrows}), thus: 
\begin{td}
\baseCat \arrow{r}{F} \arrow[swap]{d}{L} 
\arrow[phantom]{dr}[description]{\twocellRight{\rho}} & 
\altCat \arrow{d}{K} \\
\baseCat' \arrow[swap]{r}{F'} &
\altCat'
\end{td}
For every 
pair of biuniversal arrows $(R, u)$ and $(R',u')$ from $F$ to $C \in \altCat$ 
and $F'$ to $KC \in \altCat'$, respectively, the following are equivalent:
\begin{enumerate}
\item $(K, L, \rho)$ preserves the biuniversal arrow $(R, u)$, 
\item The canonical map $\psi'_{LR}(Ku \circ \rho_R) : LR \to R'$ is an 
equivalence, where we write $\psi'_{LR}$ for the chosen pseudo-inverse to 
$u' \circ F'(-) : \baseCat'(LR, R') \to \altCat'(F'LR, KC)$.
\end{enumerate}
\begin{proof}
Suppose first that $\psi'_{LR}(Ku \circ \rho_R)$ is an equivalence. Since pseudofunctors preserve equivalences, the composite $\baseCat'(B', LR) \xra{\psi'_{LR}(Ku \circ \rho_R) \circ (-)} \baseCat'(B', R') \xra{u' \circ F'(-)} \altCat'(F'C', KC)$ is an equivalence. Hence $u' \circ F'(\psi'_{LR}(Ku \circ \rho_R))$ is a biuniversal arrow. But then the 2-cell $\epsilon'_{LR}(Ku \circ \rho_R)$ provides a natural isomorphism $u' \circ F'(\psi'_{LR}(Ku \circ \rho_R)) \XRA{\iso} Ku \circ \rho_R$, so $Ku \circ \rho_R$ is also a biuniversal arrow. 

The converse is a straightforward application of universality (\cf~also~Lemma~\ref{lem:biuniversal-arrows-unique}): if $(LR, Ku \circ \rho_R)$ and $(R', u')$ are both biuniversal arrows from $F'$ to $KC$, then the canonical arrows $LR \to R'$ and $R' \to LR$ obtained from the universal property must form an equivalence. 
\end{proof}
\end{mylemma}

It will be useful to define \Def{strict preservation} of 
biuniversal arrows. This strictness will play an important role in later 
chapters, where we will ask that the syntactic models of 
our type theories satisfy a strict freeness property. The aim of this 
definition is to ensure 
that the chosen structure witnessed by a biuniversal arrow (\eg~a bilimit) is 
taken to exactly the chosen structure in the target. 

\begin{mydefn} \label{def:strict-preservation-of-biuniversal-arrows}
Let $F : \baseCat \to \altCat$ and $F' : \baseCat' \to \altCat'$ be 
pseudofunctors and suppose $(R, u)$ and $(R', u')$ are biuniversal arrows from 
$F$ to $C \in \altCat$ and from $F'$ to $C' \in \altCat'$, respectively. A pair 
of pseudofunctors $(K, L)$ is a \Def{strict morphism of biuniversal arrows} 
from $(R,u)$ to $(R',u')$ if 
\begin{enumerate}
\item $K$ and $L$ are strict pseudofunctors such that $KF = F'L$, 
\item The data of the biuniversal arrow is preserved: $LR =  R'$, $KC = C'$ and $Ku = u'$, 
\item The mappings $\psi_B : \altCat(FB, C) \to \baseCat(B,R)$ and $\psi'_{B'} : 
\altCat'(F'B', C') \to \baseCat'(B',R')$ are preserved, so that $L\psi_B(f) = 
\psi_{LB}'K(f)$ for every $f : FB \to C$, 
\item For every $B \in \baseCat$ and equivalence 
$u \circ F(-) : \baseCat(B, R) \leftrightarrows \altCat(FB, C) : \psi_B$ 
the universal arrow 
\mbox{$\epsilon_{B,h}  : u \circ F\psi_B(h) \To h$} is strictly preserved, in 
the sense that $K_{FB,C}(\epsilon_{B,h}) = \epsilon'_{LB, Kh}$. \qedhere
\end{enumerate}
\end{mydefn}

In bicategory theory it is usually good practice to specify data up to 
equivalence, as pseudofunctors preserve equivalences but may not preserve 
isomorphisms or equalities. The preceding definition abuses this convention, 
and so is not `bicategorical' in style. A consequence is that an arbitrary
biequivalence may not
strictly preserve 
biuniversal arrows (\cf~the proof of 
Lemma~\ref{lem:biequivalences-preserve-biuniversal-arrows}). This level of 
strictness does, however, 
provide a way to talk about free bicategories-with-structure using the 
language of 1-category theory (\cf~\cite[Proposition 2.10]{Gurski2006}).

\begin{myremark}
We distinguish between \emph{preservation} of biuniversal arrows in 
the sense 
of Definition~\ref{def:preservation-of-biuniversal-arrows} and a 
\emph{morphism} of biuniversal arrows as in the preceding definition on the 
following
basis. In Definition~\ref{def:preservation-of-biuniversal-arrows} we 
require that the image of the given biuniversal arrow is a biuniversal 
arrow, but do not specify its exact nature. In 
the preceding definition, by contrast, we require that the pair $(K,L)$ takes 
the biuniversal 
arrow specified in the source to exactly the biuniversal arrow specified in the 
target.
\end{myremark}

Strict preservation of a biuniversal arrow is sufficient to imply 
preservation of the corresponding universal property, in the following sense.

\begin{mylemma} \label{lem:strict-preservation-strict-pres-UMP}
Let $F : \baseCat \to \altCat$ and $F' : \baseCat' \to \altCat'$ be 
pseudofunctors and suppose $(R, u)$ and $(R', u')$ are biuniversal arrows from 
$F$ to $C \in \altCat$ and from $F'$ to $C' \in \altCat'$, respectively. If 
$(K,L)$ is a strict morphism from $(R, u)$ to $(R', u')$, then for every $B \in 
\baseCat$, $h : B \to R$ and $\tau : u \circ Fh \To f$, $L\trans{\tau} = 
\trans{(K\tau)}$.
\begin{proof}
%
It suffices to show that $L\trans{\tau}$ satisfies the universal property of 
$\trans{(K\tau)}$. For this one observes that 
\begin{align*} 
\epsilon'_{LB, Kf} \vert F'L\trans{\tau} &= K(\epsilon_{B, f}) \vert 
KF(\trans{\tau}) \qquad \text{ by strict preservation } \\
	&= K(\epsilon_{B, f} \vert F\trans{\tau}) \\
	&= K\tau
\end{align*}
as required.
\end{proof}
\end{mylemma}

A strict morphism of 
biuniversal arrows $(K, L)$ defines a morphism of adjunctions (in the sense of 
Mac~Lane) at every hom-category. 
Indeed, it follows directly from the definition that for every $B \in \baseCat$ 
the following square commutes:
\begin{equation*} \label{eq:morphism-of-adjunctions}
\begin{tikzcd}
\baseCat(B, R) \arrow{rrr}{u_C \circ F(-)} \arrow[swap]{d}{L_{B,R}}  & 
\: &
\: &
\altCat(FB, C) \arrow{d}{K_{FB,C}} \\
\baseCat'(LB, LR) \arrow[equals]{r} &
\baseCat'(LB, R') \arrow[swap]{r}[yshift=-2mm]{u'_{LB} \circ F'(-)} & 
\altCat'(F'LB, C') \arrow[equals]{r} &
\altCat'(KFB, KC)
\end{tikzcd}
\end{equation*}
and $K_{FB,C}$ preserves the counit by assumption.

\section{Bilimits} \label{sec:bilimits}

We are now in a position to introduce bilimits and preservation of bilimits. 
The formulation in terms of biuniversal arrows is pleasingly concise. For every 
pair of bicategories $\catJ, \baseCat$ one has a 
\Def{diagonal pseudofunctor} 
$\Delta : \baseCat \to \Hom(\catJ, \baseCat)$ 
taking $B \in 
\baseCat$ to the constant pseudofunctor at $B$. Explicitly, $\Delta B : \catJ 
\to \baseCat$ takes a 2-cell $\tau : h \To h' : j \to j'$ to the identity 
2-cell $\id_B : \Id_B \To \Id_B : B \to B$. The 2-cell $\psi_j : \Id_{(\Delta 
B)(j)} \To (\Delta B)(\Id_j)$ is the identity and for a composite $j \xra{g} j' 
\xra{f} j''$ in $\catJ$ the 2-cell $\phi_{f,g} : (\Delta B)(f) \circ (\Delta 
B)(g) \To (\Delta B)(f \circ g)$ is $\l_{\Id_B} : \Id_B \circ \Id_B \To \Id_B$. 
A bilimit is then a biuniversal arrow. 

\begin{mydefn} \label{def:bilimit}
A \Def{bilimit} for $F : \catJ \to \baseCat$ is a biuniversal arrow from the diagonal pseudofunctor $\Delta : \baseCat \to \Hom(\catJ, \baseCat)$ to $F$. 
\end{mydefn}

Unwrapping the definition, we require a pair 
$(\bilim F, \lambda : \Delta(\bilim F) \To F)$ such that for every object $B 
\in \baseCat$ and cone 
(pseudonatural transformation) 
$\kappa 
: \Delta B \To F$ there exists a map $u_\kappa : B \to \bilim F$ and an 
invertible modification $\epsilon_{B, \kappa}$ filling
\begin{td}
\Delta B \arrow{rr}{\Delta (u_\kappa)} \arrow[swap]{dr}{\kappa} &
\: \arrow[phantom]{d}[description]{\twocell{\epsilon_{B, \kappa}}} &
\Delta(\bilim F) \arrow{dl}{\lambda} \\
\: &
F &
\end{td}
This modification is required to be universal in the sense that, for any 
1-cell
$v : B \to \bilim F$ 
and 2-cell
$\beta : \lambda \circ \Delta v \To \kappa$, 
there exists a unique 
$\trans{\beta} : v \To 
u_\kappa$ such that 
\begin{equation*} 
\begin{tikzcd}[column sep = 3em]
\Delta B \arrow[bend right=60, swap]{drr}{\kappa} 
\arrow[phantom]{r}[xshift=1.2mm, 
font=\scriptsize]{{\Downarrow\!\Delta\trans{\beta}}} 
\arrow[bend left]{r}{\Delta v} 
\arrow[bend right, swap]{r}[xshift=2mm]{\Delta u_\kappa} &
\Delta(\bilim F) \arrow[phantom]{d}[description, font=\small, near end, 
yshift=-2mm]{{\Downarrow\!\epsilon_{B,\kappa}}} \arrow{dr}{\lambda} \\
\: &
\: &
F
\end{tikzcd}
\quad = \quad
\begin{tikzcd}
\: &
\Delta(\bilim F) \arrow[phantom]{d}[description, 
font=\small]{{\Downarrow\!\beta}} \arrow{dr}{\lambda} &
\: \\
\Delta B \arrow[swap]{rr}{\kappa} \arrow{ur}{\Delta v} &
\: &
F
\end{tikzcd}
\vspace{-2mm}
\end{equation*}
Finally, we require that for every $w : B \to \bilim F$ the 2-cell 
$\trans{(\id_{\lambda \circ \Delta w})} :   
	w \To u_{\lambda \circ \Delta w}$
is invertible.

By Lemma~\ref{lem:biuniversal-arrow-equivalence} this definition can be 
rephrased as a pseudonatural family of adjoint equivalences 
$\baseCat(B, \bilim F) \simeq \Hom(\catJ, \baseCat)(\Delta B, F)$. It 
therefore coincides with that of Street~\cite{Street1980} in terms of 
birepresentations. We say that a bicategory $\baseCat$ is \Def{bicomplete} or 
\Def{admits all bilimits} if for every small bicategory
$\catJ$ 
and pseudofunctor $F : \catJ \to \baseCat$ the bilimit $\bilim F$ exists 
in $\baseCat$.

\paragraph*{Preservation of bilimits.} We define preservation of bilimits as 
preservation of the corresponding biuniversal arrows, via the following lemma. 

\begin{mylemma}
For any 
bicategory $\catJ$ and pseudofunctor $H : \baseCat \to \altCat$ the following 
diagram commutes up to canonical isomorphism:
\begin{equation} \label{eq:bilimit-preservation-diagram}
\begin{tikzcd}
\baseCat 
\arrow[phantom]{dr}{\twocell{\iso}}
\arrow{r}{\Delta^{\baseCat}} 
\arrow[swap]{d}{H} &
\Hom(\catJ, \baseCat) \arrow{d}{H \circ (-)} \\
\altCat \arrow[swap]{r}{\Delta^{\altCat}} &
\Hom(\catJ, \altCat)
\end{tikzcd}
\end{equation}
\begin{proof}
Let us write $H_\ast := H \circ (-)$.
Unwinding the respective definitions,  
$(H_\ast \circ \Delta^{\baseCat})B : \catJ \to \altCat$ 
is the pseudofunctor sending every $j \in \catJ$ to $HB$, every 
$p : j \to j'$ to $H\Id_B$ and every 2-cell $\sigma : p \To p'$ to
the identity. This coincides with 
$(\Delta^{\altCat} \circ H)B$ everywhere except that 
$(\Delta^{\altCat} \circ H)(B)(j \xra{p} j') = \Id_{HB}$. 
So for every $B \in \baseCat$ there exists a pseudonatural isomorphism 
$\alpha_B := 
	(H_\ast \circ \Delta^{\baseCat})B \To (\Delta^{\altCat} \circ H)B$ 
with components $\alpha_B(j) := \Id_{HB}$ for all $j \in \catJ$. 
The witnessing 2-cell is the evident composite of $\psi^H$ with structural isomorphisms. 
Thus one obtains an invertible 1-cell $\alpha_B$ in $\Hom(\catJ, \altCat)$ for every
$B \in \baseCat$. To extend this to a pseudonatural isomorphism, one takes
$\overline{\alpha}_f : 
	\alpha_{B'} \circ H_\ast(\Delta^{\baseCat}f) 
	\To 
	\Delta^{\altCat}(Hf) \circ \alpha_B$ 
(for $f : B \to B'$)
to be the invertible modification with components given by the structural isomorphism
$\Id_{HB'} \circ Hf \XRA{\iso} Hf \circ \Id_{HB}$. 
Then $(\alpha, \overline{\alpha})$ is the required isomorphism. 
\end{proof}
\end{mylemma}

Thus, assuming the bilimit exists in $\altCat$, 
we say that $H$ \Def{preserves the bilimit} of $F : \catJ \to \baseCat$ if 
$\left(H_\ast, H, (\alpha, \overline{\alpha})\right)$ preserves the biuniversal arrow 
$(\bilim F,\lambda)$. 
By 
Lemma~\ref{lem:preserves-cones-equivalent-to-equivalence}, this condition is 
equivalent to requiring that the canonical map $H(\bilim F) \to \bilim (HF)$ is 
an equivalence. 

The general perspective of biuniversal arrows leads to a straightforward proof 
that biequivalences preserve all bilimits. 

\begin{mycor} \label{cor:biequivalences-preserve-bilimits} 
For any biequivalence  $H :  \baseCat \leftrightarrows \baseCat' : G$,
\begin{enumerate}
\item \label{c:preservation-of-bilimits} $H$ preserves all bilimits that exist 
in $\baseCat$,
\item If $\baseCat$ has all $\catJ$-bilimits then $\baseCat'$ has all 
$\catJ$-bilimits.
\end{enumerate}
\begin{proof}
For~(1), suppose $F : \catJ \to \baseCat$ has a bilimit. By 
Lemma~\ref{lem:biequivalences-preserve-biuniversal-arrows} one obtains a 
biuniversal arrow from $H_\ast \circ \Delta$ to $H_\ast(F)$, which 
by~(\ref{eq:bilimit-preservation-diagram}) is biuniversal from 
$\Delta^{\baseCat'}H$ to $HF$. So the bilimit is preserved. 

For (2), suppose $F : \catJ \to \baseCat'$. Then $GF : \catJ \to \baseCat$ has a bilimit and hence, by the previous part, so does $HGF : \catJ \to \baseCat'$. Since $HG \simeq \id_{\baseCat'}$, it follows that $F$ has a bilimit.
\end{proof}
\end{mycor}

Two other classes of pseudofunctors that one would certainly expect to preserve 
bilimits are right biadjoints (see 
Definition~\ref{def:biadjunction-universal-arrow}) and birepresentables. This 
is indeed the case.

\newpage
\begin{mylemma} \label{lem:representables-and-adjoints-preserve-bilimits}
\quad
\begin{enumerate}
	\item If the pseudofunctor $F : \baseCat \to \altCat$ has a left biadjoint, 
	then $F$ preserves all bilimits that exist in $\baseCat$.
	\item If $F : \baseCat \to \Cat$ is a birepresentable pseudofunctor, then $F$ preserves all bilimits that exist in $\baseCat$. 
\end{enumerate}
\begin{proof}
These are~\cite[\S 1.32]{Street1980} and \cite[\S 1.20]{Street1980}, respectively.
\end{proof}
\end{mylemma}

\section{Biadjunctions}

Recalling that an adjunction is specified by a choice of universal arrows, we 
define a biadjunction by a choice of biuniversal arrows (\cf~\cite{Power1998}). 

\begin{mydefn} \label{def:biadjunction-universal-arrow}
Let $F : \baseCat \to \altCat$ be a pseudofunctor. To specify a \Def{right biadjoint} 
to $F$ is to specify a biuniversal arrow 
$(UC, u_C : FUC \to C)$ from $F$ to $C$ for every $C \in \altCat$.
\end{mydefn}

Spelling out the definition, to give a right 
biadjoint $U : \altCat \to \baseCat$ to $F$ is to give:
\begin{itemize}
\item A mapping $U : ob(\altCat) \to ob(\baseCat)$, 
\item A family of 1-cells $(u_C : FUC \to C)_{C \in \altCat}$, 
\item For every $B \in B$ and $h : FB \to C$ a 1-cell $\psi_B(h) : B \to UC$ 
and an invertible 2-cell $\epsilon_{B,h} : u_C \circ F\psi_B(h) \To h$ that is 
universal in the sense of~(\ref{eq:universal-arrow-UMP})~(p.~\pageref{eq:universal-arrow-UMP}), 
such that the unit 
$\eta_h := \trans{(\id_{u_C \circ Fh})} : h \To \psi_B(u_C \circ Fh)$ 
is invertible for every $h$.
\end{itemize}
One thereby obtains a pseudofunctor $U : \altCat \to \baseCat$ by setting 
$U(C) := UC$ on objects, 
$U({C \xra{g} C'}) := \psi_{UC}(g \circ u_C)$ and 
$U(g \XRA{\sigma} g') := \trans{\left((\sigma \circ u_C) \vert \epsilon_{UC, g}\right)}$.
By  Lemma~\ref{lem:biuniversal-arrow-equivalence}, this definition is 
equivalent to asking for a pair of pseudofunctors 
\mbox{$F : \baseCat \leftrightarrows \altCat : U$} together with a pseudonatural family of 
equivalences $\baseCat(B, UC) \simeq \altCat(FB, C)$. For detailed proofs of 
this and related results, see~\cite[Chapter 9]{TFiore2006}. 

The biuniversal arrow formulation of biadjoints, relying as it does on universal properties 
at each level, is perhaps easiest to work with 
when it comes to calculations (\cf~\cite{FioreSpecies}). As we shall see in 
Chapters~\ref{chap:fp-lang} and~\ref{chap:ccc-lang}, it is also particularly 
amenable to being expressed syntactically.

\begin{myremark} \label{rem:bilimits-as-biadjoints}
The definition of bilimit can now be rephrased in the 
following fashion: the pseudofunctor $\mathrm{bilim} : \Hom(\catJ, \baseCat) 
\to \baseCat$, when it exists, is right biadjoint to the diagonal 
pseudofunctor~(\cf~\cite[Remark 9.2.1]{TFiore2006}).
\end{myremark}

We have chosen to place bilimits and biadjoints on a similar 
footing by presenting them both as instances of biuniversal arrows. The 
preceding remark indicates that the theory of bilimits could alternatively be 
phrased using biadjoints. For example, one may use the fact that a 
right biadjoint preserves all bilimits, together with the 
observation that every biequivalence can be `upgraded' to an adjoint 
biequivalence~\cite{Gurski2012}, to obtain an 
alternative proof of 
Corollary~\ref{cor:biequivalences-preserve-bilimits}(\ref{c:preservation-of-bilimits}).

\paragraph*{Preservation of biadjunctions.} 
We shall use the notion of
preservation of biadjunctions to define preservation 
of exponentials. 

\begin{mydefn}
For any biadjoint pair 
${F : \baseCat \leftrightarrows \altCat : U}$ and pseudofunctor
$F' : \baseCat' \to \altCat'$, 
we say that the triple $(K, L, \rho)$ as below 
\begin{equation} \label{eq:preservation-of-biadjoints}
\begin{tikzcd}
\baseCat \arrow{r}{F} \arrow[swap]{d}{L} \arrow[phantom]{dr}[description]{\twocellRight{\rho}} & 
\altCat \arrow{d}{K} \\
\baseCat' \arrow[swap]{r}{F'} &
\altCat'
\end{tikzcd}
\end{equation}
\Def{preserves the biadjunction} if $(K, L, \rho)$ preserves each biuniversal 
arrow $u_C : FUC \to C$. 
\end{mydefn}

A triple $(K, L, \rho)$ preserving a biadjunction preserves the corresponding 
counits up to isomorphism. By definition, 
whenever 
$(K, L, \rho)$ preserves the 
biadjunction $F \dashv U$ as in~(\ref{eq:preservation-of-biadjoints}), then 
$
F'LUC \xra{\rho_{UC}} KFUC \xra{Ku_C} KC
$ 
is a biuniversal arrow from $F'L$ to $KC$. The next lemma entails that,
if $F'$ has a right adjoint $U'$, then
\[
F'U'KC \xra{\simeq} F'LUC \xra{\rho_{UC}} KFUC \xra{Ku_C} KC
\]
is another such biuniversal arrow. By 
Lemma~\ref{lem:biuniversal-arrows-unique}, this must be
canonically isomorphic to the biuniversal arrow $u'_{KC}$ 
witnessing the biadjunction $F' \dashv U'$. 

\begin{mylemma} \label{lem:preservation-of-biadjoints-to-equivalence}
Let $(K, L, \rho)$ preserve the biadjunction $F \dashv U$ as 
in~(\ref{eq:preservation-of-biadjoints}) and suppose $F'$ has a right 
biadjoint $U'$. Then $U'K \simeq LU$. 
\begin{proof}
The definition of preservation of a biuniversal arrow, together with the 
definition of a biadjunction, entails that for any $B \in \baseCat$ and $C \in 
\altCat$:
\[
\baseCat'(B, LUC) \simeq \altCat'(F'B, KC) \simeq \baseCat'(B, U'KC)
\]
By Lemma~\ref{lem:biuniversal-arrow-equivalence} these equivalences may equally 
be expressed as equivalences of pseudofunctors. Hence, $\Yon \circ (LU) \simeq 
\Yon \circ (U'K)$, for 
$\Yon : \baseCat' \to \Hom\left(\op{(\baseCat')}, \Cat\right)$ the 
Yoneda embedding. The Yoneda Lemma then entails that $LU \simeq U'K$, as 
claimed. 
\end{proof}
\end{mylemma}

We end this chapter by instantiating  
Lemma~\ref{lem:biequivalences-preserve-biuniversal-arrows} in the particular 
case of biadjunctions.

\newpage
\begin{mylemma} \label{lem:biequivalences-preserves-biadjoints}
Suppose that $F : \baseCat \to \altCat$ has a right biadjoint $U$ and that 
$H : \altCat \leftrightarrows \altCat' : G$ is a biequivalence. Then 
$HF : \baseCat \leftrightarrows \altCat' :  UG$ is a biadjunction.
\begin{proof}
By Lemma~\ref{lem:biequivalences-preserve-biuniversal-arrows}, each biuniversal 
arrow 
$u_C : FUC \to C$ 
defining the biadjunction $F \dashv U$ is preserved. In 
particular, taking $C' \in \altCat'$ such that $GC' \simeq C$ and the 
biuniversal arrow $u_{GC'} : FUGC' \to GC'$, one obtains a 
biuniversal arrow $HFUGC' \to HGC'$ from $HF$ to $HGC'$. But from the 
biequivalence one has an adjoint
equivalence $HG \simeq \id_{\altCat'}$ for which the component at $C'$ is an 
adjoint equivalence 
$HGC' \simeq C'$. Composing, 
there exists a biuniversal arrow
$(HF)(UG)C' \to C'$ 
from $HF$ 
to 
$C'$, as required. 
\end{proof}
\end{mylemma}

\part{A type theory for cartesian closed bicategories}
\label{part:lang}


\chapter{A type theory for biclones} \label{chap:biclone-lang}

In this chapter we begin our construction of the type theory $\langCartClosed$ 
for cartesian 
closed bicategories. We focus on the bicategorical fragment: we construct a 
type theory $\langBicat$ for bicategories and use it to recover a version of 
the Mac~Lane-Par{\'e} coherence theorem for bicategories~\cite{MacLane1985}.

The work is driven by the theory of \emph{biclones}, a 
bicategorification of the \emph{abstract clones} of universal 
algebra~\cite{Cohn1981}. Abstract clones axiomatise the notion of 
equational 
theory with variables and a substitution operation, and provide a natural 
intermediary 
between syntax (in the form of the set of terms generated from
operators over a set of variables) and semantics (in the form of categorical 
algebraic theories) (see~\eg~\cite[p.129]{CloneBookRef}). 
Biclones will play the same role in our construction, axiomatising syntax with 
an up-to-isomorphism substitution operation. We shall then synthesise 
the rules of our type theory $\langBiclone$ from biclone structure. 

The resulting type theory varies from classical type theories such as the 
simply-typed lambda calculus in two important respects. First, we make 
use of a form of~\emph{explicit substitution}~\cite{Abadi1989}; second, it is 
\emph{2-dimensional} in the sense that judgements relate types, terms and 
\emph{rewrites} between terms. 

These two developments both arise in the study of 
rewriting in the lambda calculus, but have previously only been studied 
independently. 
Explicit substitution calculi were first studied as versions of the 
lambda calculus closer to machine implementation~\cite{Abadi1989} 
and have found applications in proof theory~\cite{Ritter2000} and programming 
language theory~\cite{Levy1999}. Much recent research~(\eg~\cite{DiCosmo1997, 
Ritter1999}) has focussed on 
Melli{\`e}s' observation that, contrary to what one might expect 
from the lambda calculus, such calculi may not be strongly 
normalising~\cite{Mellies1995}~(see~\eg~\cite{Rose2011} for an overview).

Two-dimensional type theories, on the other hand, first arose from Seely's 
observation~\cite{Seely1987} that 
$\eta$-expansion  and $\beta$-reduction form the unit and counit of a 
\emph{lax} 
(directed) cartesian closed structure, a perspective advocated further by 
Jay \& Ghani~\cite{Ghani1995,Jay1995} and put to use by 
Hilken~\cite{Hilken1996} for a proof-theoretic account of rewriting. In the 
strict setting, 
Hirschowitz~\cite{Hirschowitz2013} and 
Tabereau~\cite{Tabareau2011} have constructed 2-dimensional type theories to 
describe 2-categorical structures in rewriting theory and programming language 
design, respectively.
The connection with intensional equality, meanwhile, has recently sparked
significant interest in type theories with a notion of `rewrite' or `equality' 
motivated by the connection between higher category theory, topology and type 
theory. 
Examples include Licata \& Harper's 2-dimensional directed type 
theory~\cite{Licata2011, Licata2012}, Riehl \& Shulman's type theory for 
synthetic 
$\infty$-categories~\cite{Riehl2017}, and Garner's 2-dimensional type 
theory~\cite{Garner2009}.  

The type theory we shall construct brings together a novel combination of 
explicit substitution and 2-dimensional judgements. Following Hilken, we relate 
terms by separate syntactic entities called \Def{rewrites}, and interpret these 
as 2-cells. This 
contrasts with many type theories motivated by connections with homotopy type 
theory (\eg~the 
Riehl-Shulman and Garner type theories), which capture 2-cells using 
Martin-L\"of style identity types. The relationship between the 
two approaches remains to be explored. 

\paragraph*{Outline.} The chapter breaks up into three parts. In 
Section~\ref{sec:bicategorical-type-theory} we 
consider the appropriate form of signature for a 2-dimensional type theory 
and construct the free biclone over such a signature. This drives the second 
part (Section~\ref{sec:type-theory-for-biclones}), where we synthesise the 
type theory $\langBiclone$ and show that it is 
the internal language of biclones; as a corollary, we obtain an internal 
language for bicategories. Finally, in Section~\ref{sec:coherence-for-biclones} 
we use $\langBiclone$ to prove a coherence 
result for biclones, amounting to a form of normalisation for the 
corresponding type theory.

\section{Bicategorical type theory}
\label{sec:bicategorical-type-theory}

\subsection{Signatures for 2-dimensional type theory}

A signature for the simply-typed lambda calculus is specified by a 
choice of base types and constants (sometimes called a 
\emph{$\lambda\!{\times}$-signature}~\cite{Crole1994}).  A natural way of 
packaging such data, exemplified by Lambek \& Scott~\cite{LambekAndScott}, is 
as a graph. Taking inspiration from Lambek's notion of multicategories as 
models of \emph{deductive systems}~\cite{Lambek1969, LambekAndScott}, one may 
extend this using a \emph{multigraph} (\cf~\cite{Lambek1989, Hermida2000, 
Leinster2004}). 
Here, one thinks of a judgement 
$({x_1 : A_1, \,\dots\, , x_n : A_n \vdash t : B})$ as 
corresponding to an edge with source $(A_1, \,\dots\, , A_n)$ and target $B$.%
\footnote{This should not be confused with the terminology in 
graph theory, where a multigraph sometimes refers to a graph in which there are 
allowed to be multiple edges between nodes~(\eg~\cite[p.10]{graphtheorybook}).}  

\begin{mydefn} 
A \Def{multigraph} $\graph$ consists of a set $\nodes\graph$ of \Def{nodes} 
together with a set $\graph(A_1, \,\dots\, , A_n; B)$ of \Def{edges} from 
$(A_1, 
\dots, A_n)$ to $B$ for every $A_1, \,\dots\, , A_n, B \in \nodes{\graph}$ (we 
allow 
$n=0$). A \Def{homomorphism} of multigraphs 
$h = (h, h_{A_1, \,\dots\, , A_n; B}) : \graph \to \graph'$ consists 
of a function $h : \nodes{\graph} \to \nodes{\graph}'$ together with 
functions 
$h_{A_1, \,\dots\, , A_n; B} : 
	\graph(A_1, \,\dots\, , A_n; B) \to \graph'(hA_1, \,\dots\, , hA_n; hB)$ 
for every 
$A_1, \,\dots\, , A_n, B \in \nodes{\graph} \:\: (n \in \Nat)$. We denote the 
category of multigraphs and multigraph homomorphisms by $\MultiGraph$. The full 
subcategory $\Graph$ of \Def{graphs} has objects those multigraphs $\graph$ 
such that $\graph(A_1, \,\dots\, , A_n; B) = \emptyset$ whenever $n \neq 1$. 
\end{mydefn} 

\begin{myexmp} \label{ex:lambda-calc-multigraph}
Every graph freely generates a \Def{typed 
$\lambda$-calculus}~\cite{LambekAndScott} with types the nodes and a unary 
operator for each edge. Conversely, the simply-typed lambda calculus over a 
fixed set of base types determines a multigraph with nodes the types and edges 
$(A_1, \,\dots\, , A_n) \to B$ the derivable terms \mbox{$x_1 : A_1, \,\dots\, 
, x_n : 
A_n \vdash t : B$} up to $\alpha$-equivalence (we assume a fixed enumeration of 
variables $x_1, x_2, \dots$ determining the name of the $i$th variable in the 
context). 
\end{myexmp}

%
In this vein, the appropriate notion of signature for a 2-dimensional type theory is a form of `2-multigraph' (\cf~\cite[Chapter 2]{Gurski2013}). 

\begin{mynotation}  \label{not:ind}
In the following definition, and throughout, we write $\ind{A}$ for a finite 
sequence $\seq{A_1, \,\dots\, , A_n}$.\footnote{This notation is  adopted from 
homological algebra, where one writes $\ind{X}$ for a chain complex \mbox{$X_1 
\to X_2 \to \cdots$}~(\eg~\cite{weibel}).} Following 
Example~\ref{ex:lambda-calc-multigraph}, we use Greek letters 
$\Gamma, \Delta, \dots$ to denote sequences $\seq{A_1, \,\dots\, , A_n}$ in 
which 
the names of the terms $A_i$ are not of importance. We use 
$\Gamma_1, \Gamma_2$ or $\Gamma_1 \concat \Gamma_2$ to denote the concatenation 
of $\Gamma_1$ and $\Gamma_2$, and write $\len{\Gamma}$ for the length of 
$\Gamma$.
\end{mynotation} 

\begin{mydefn}
A \Def{2-multigraph} $\graph$ is a set of \Def{nodes} $\nodes{\graph}$ equipped 
with a multigraph $\graph(\ind{A}; B)$ of \Def{edges} and \Def{surfaces} for every 
$A_1, \,\dots\, , A_n, B \in \nodes{\graph}$ (we allow $n=0$). A 
\Def{homomorphism} 
of 2-multigraphs $h = (h, h_{\ind{A}, B}, h_{f,g}) : \graph \to 
\graph'$ is a map $h : \nodes{\graph} 
\to \nodes{\graph}'$ together with functions
\begin{align*}
h_{A_1, \,\dots\, , A_n; B} : \graph(\ind{A}; B) &\to \graph'(hA_1, \,\dots\, , 
hA_n; hB) \\
h_{f,g} : \graph(\ind{A}; B)(f,g) &\to \graph'(hA_1, \,\dots\, , hA_n; 
hB)(hf,hg)
\end{align*}
for every $A_1, \,\dots\, , A_n, B \in \nodes{\graph} \:\: (n \in \Nat)$ and 
$f, g \in \graph(\ind{A}; B)$. We denote the category of \mbox{2-multigraphs} 
by 
$\TwoMultiGraph$. The full subcategory $\TwoGraph$ of \Def{2-graphs} is formed 
by restricting to 2-multigraphs $\graph$ such that $\graph(A_1, \,\dots\, , 
A_n; B) 
= \emptyset$ whenever $n \neq 1$.
\end{mydefn}

\begin{myexmp} \label{exmp:multigraphs} \quad
\begin{enumerate}
\item Every category determines a graph; every bicategory determines a 2-graph.
\item \label{c:monoidal-cat} Every monoidal category $(\cat, \tens, I)$ 
determines a multigraph $\graph_\cat$ with nodes 
$\nodes{(\graph_\cat)} := ob(\catC)$
and 
$\graph_\altCat(X_1, \,\dots\, , X_n; Y) := \cat(X_1 \tens \dots \tens X_n, Y)$
(for some chosen bracketing of the $n$-ary tensor product).
\item More generally, every \Def{multicategory}~\cite{Lambek1969} determines a multigraph.
\qedhere
\end{enumerate}
\end{myexmp}

We shall see in Chapter~\ref{chap:fp-lang} that every bicategory with finite 
products determines a \mbox{\Def{bi-multicategory}} and every bi-multicategory 
determines a 2-multigraph.

\subsection{Biclones}

We turn to constructing 
bicategorical substitution structure over a 2-multigraph.  As indicated above, 
our approach is to 
bicategorify the notion of \Def{abstract clone}~\cite{Cohn1981}. 

\paragraph*{Abstract clones.} Abstract 
clones  provide a presentation-independent 
description of (algebraic) equational theories with variables and substitution. 
A leading example is the \emph{clone of operations} given by 
the set of terms over a fixed signature, subject to the substitution operation. 
We shall recall only the basic properties we require: for an introduction to 
the theory of clones from the perspective of universal 
algebra, see~\eg~\cite{CloneBookRef, Taylor1999}.

\begin{mydefn} A \Def{(sorted) abstract clone} $(S, \clone)$ consists of a set $S$ of \Def{sorts} with
\begin{itemize} 
\item A set $\clone(X_1, \,\dots\, , X_n; Y)$ of \Def{operations}
$t : X_1, \,\dots\, , X_n \to Y$ for each 
${X_1, \,\dots\, , X_n, Y \in S} \:~({n \in \Nat})$,
\item Distinguished \Def{projections} 
$\p{i}{\ind{X}} \in \clone(X_1, \,\dots\, , X_n; X_i) \: (i = 1, \,\dots\, , 
n)$ for 
each ${X_1, \,\dots\, , X_n \in S} \:~ {(n \in \Nat)}$, 
\item For all sequences of sorts $\Gamma$ and sorts 
$Y_1, \,\dots\, , Y_n, Z \:\: (n \in \Nat)$ a 
\Def{substitution} function
\begin{gather*}
\subName_{\Gamma, \ind{Y}, Z} : \clone(\ind{Y}; Z) \times \prodop_{i=1}^n 
\clone(\Gamma; Y_i) \to \clone(\Gamma; Z)
\end{gather*}
we denote by $\subName\big(f, (g_1, \,\dots\, , g_n)\big) := \cslr{f}{g_1, 
\,\dots\, , 
g_n}$,
\end{itemize} 
such that 
\begin{enumerate} 
\item $\cslr{t}{\p{1}{\ind{X}}, \,\dots\, , \p{n}{\ind{X}}} = t$ for all $t \in 
\clone(\ind{X}; Y)$, 
\item $\cslr{\p{k}{\ind{Y}}}{t_1, \,\dots\, , t_n} = t_k \:\: (k=1, \,\dots\, , 
n)$ for 
all $\left(t_i \in \clone(\Gamma; Y_i)\right)_{i=1,\dots,n}$, 
\item $\cslr{\cslr{t}{\ind{u}}}{\ind{v}} = \cslr{t}{\cslr{\ind{u}}{\ind{v}}}$ 
for all $v_j \in \clone(\ind{W}; X_j)$, $u_i \in \clone(\ind{X}; Y_i)$ and $t 
\in \clone(\ind{Y}; Z)$ ($ i = 1, \,\dots\, , n$ and $j = 1, \,\dots\, , m$). 
\end{enumerate} 
We write $\cslr{(\cslr{t}{\ind{u}})}{\ind{v}}$ for the iterated substitution 
$\cslr{\cslr{t}{u_1, \,\dots\, , u_n}}{v_1, \,\dots\, , v_m}$; by 
default, we 
bracket substitution to the left. An operation of form $t : X \to Y$ is called 
\Def{unary}. 

A morphism $h = (h, h_{\ind{X}; Y}) : (S, \clone) \to (S', \clone')$ of 
abstract clones is a map 
$h : {S \to S'}$ together with functions 
$h_{\ind{X}; Y} : 
	\clone(X_1, \,\dots\, , X_n; Y) \to \clone'(hX_1, \,\dots\, , hX_n; hY)$ 
for each ${X_1, \,\dots\, , X_n, Y \in S}$, such that 
the projections and substitution operation are preserved. We denote the 
category of clones by $\Clone$. 
\end{mydefn} 

Following the terminology for multicategories, we occasionally refer to the 
operations $t : X_1, \,\dots\, , X_n \to Y$ of a clone as \Def{multimaps} or 
\Def{arrows}. Where the 
context is unambiguous, we refer to a sorted clone $(S, \clone)$ simply as an 
\Def{$S$-clone} and denote it by $\clone$; a clone with a single 
sort is called \Def{mono-sorted}. 

\newpage
\begin{myexmp} \label{ex:CartesianCategoryIsAClone} \quad 
\begin{enumerate}
\item Every clone $(S, \clone)$ defines a 
category $\overline\clone$ by restricting to the unary operations.  We call this the \Def{nucleus} of $(S, \clone)$. 
Composition is given by substitution in $(S, \clone)$ and the identity on $X \in S$ is
$\p{1}{X}$.

\item \label{c:cart-cat-to-clone}  Any small category $\cat$ 
with finite products defines an $ob(\cat)$-clone $\cloneOf\cat$ with 
\[
\cloneOf{\cat}(X_1, \,\dots\, , X_n; Y) 
	:= \cat(X_1 \times \dots \times X_n, Y)
\]
The projections are the projections in $\cat$; the substitution 
$\cslr{t}{u_1, \,\dots\, , u_n}$ is the composite 
$t \circ \seq{u_1, \,\dots\, , u_n}$. 
\qedhere
\end{enumerate}
\end{myexmp} 

One may read the two cases just presented as follows: every Lawvere theory 
defines a mono-sorted clone, and every mono-sorted clone defines a Lawvere 
theory. In fact, the full subcategory of $\Clone$ consisting 
of just the mono-sorted clones is 
equivalent to the category of Lawvere theories (see~\eg~\cite{CloneBookRef}).  
This makes precise the sense in which clones capture a notion of algebraic 
theory. In the next chapter we shall explore 
the relationship between multi-sorted clones and cartesian categories more 
generally.


\paragraph*{Clones and type-theoretic syntax.}  The definition of abstract 
clone isolates 
three axioms sufficient to describe substitution.  The next example 
shows how a clone augments a graph with a notion of substitution 
(\cf~Example~\ref{ex:lambda-calc-multigraph}).

\begin{myexmp} \label{ex:stlc-a-clone}
For a chosen set of base types $\baseTypes$ and multigraph $\graph$ with nodes 
generated 
by the grammar
\[
X, Y ::= B \st X \times Y \st \exptype{X}{Y} \qquad (B \in \baseTypes)
\]
the corresponding 
lambda calculus may be 
equipped with a simultaneous substitution operation 
$\left( t, (u_1, \,\dots\, , u_n) \right) \mapsto t[u_1 / x_1, \,\dots\, , u_n 
/ x_n]$
which respects the typing 
in the sense that the following rule is admissible:
\begin{prooftree}
\AxiomC{$x_1 : A_1, \,\dots\, , x_n : A_n \vdash t : B$}
\AxiomC{$(\Delta \vdash u_i : A_i)_{i=1, \,\dots\, , n}$}
\BinaryInfC{$\Delta \vdash t[u_1 / x_1, \,\dots\, , u_n / x_n]$}
\end{prooftree}
One therefore obtains a clone with sorts the types and 
multimaps $X_1, \,\dots\, , X_n \to Y$ the $\alpha$-equivalence classes of 
derivable 
terms $x_1 : X_1, \,\dots\, , x_n : X_n \vdash {t : Y}$. 
The three axioms 
encapsulate 
the following standard properties of simultaneous substitution (\cf~the 
\emph{syntactic substitution lemma}~\cite[p.27]{Barendregt1985}):
\begin{center}
$x_k[u_1 / x_1, \,\dots\, , u_n / x_n] = u_k$ \qquad $t[x_1 / x_1, \,\dots\, , 
x_n / x_n] = t$ \qquad $t[u_i/x_i][v_j/y_j] = t\big[u_i[v_j/y_j] / x_i\big]$
\end{center}
One still obtains a clone if one takes $\alpha\beta\eta$-equivalence classes of 
terms; we denote this by $\stlcCloneTimesExp{\graph}$.
\end{myexmp}

Example~\ref{ex:stlc-a-clone} exemplifies the way in which clones provide an 
algebraic description of (type-theoretic) syntax. The tradition of categorical 
algebra, on the other hand, describes such syntax through the construction of a 
\Def{syntactic category}, for which one aims to prove a freeness 
universal property. Generally some massage is required to account for the fact 
that categorical morphisms take a single object as their domain, but terms may 
exist in contexts 
of arbitrary length. For example, one may take contexts as objects and 
morphisms as lists of terms~(\eg~\cite{Pitts2000}), or restrict to unary 
contexts and take morphisms to be single terms~(\eg~\cite{Crole1994}). It turns 
out that, if one 
employs the latter strategy, the  
relationship between the clone-theoretic and category-theoretic perspectives is 
particularly tight.

\begin{mylemma} \quad \label{lem:graph-multi-graph-adj-and-free-clone}
\begin{enumerate}
\item \label{c:graph-in-multigraph} The inclusion
$\Graph \hookrightarrow \MultiGraph$ 
has a right adjoint given by restricting to edges of the form $X \to Y$. 
\item \label{c:clone-to-multigraph} The forgetful functor 
$\Clone \to \MultiGraph$ 
taking a clone to 
its underlying multigraph has a left adjoint.
\item \label{c:cat-to-clone} The functor 
$\overline{(-)} : \Clone \to \CatCat$ 
restricting a clone to its nucleus has a left adjoint.
\end{enumerate}
\begin{proof}
For~(\ref{c:graph-in-multigraph}) define a functor 
$\lin : \MultiGraph \to \Graph$ 
by taking $\lin\graph$ to be the graph with nodes 
exactly the nodes of $\graph$ and edges $(\lin\graph)(X, Y) := \graph(X, Y)$. 
The action on homomorphisms is similar: for 
$h : \graph \to \graph'$ one obtains $\lin(h)$ by restricting to edges of the 
form $X \to Y$. Then, where 
$\inc : \Graph \hookrightarrow \MultiGraph$
denotes the obvious embedding, a multigraph homomorphism
$h : \inc(\graph) \to \graph'$
is a map on nodes 
$h : \nodes{(\inc\graph)} \to \nodes{\graph}'$
together with maps
$h_{\ind{X}; Y} : (\inc\graph)(X_1, \,\dots\, , X_n; Y) \to
		\graph'(hX_1, \,\dots\, , hX_n; hY)$
for each $X_1, \,\dots\, , X_n, Y \in \nodes{(\inc\graph)} \:\: (n \in \Nat)$.
Since $(\inc\graph)(X_1, \,\dots\, , X_n; Y)$ is empty except when $n=0$, this 
is 
equivalently a graph homomorphism 
$\graph \to \lin\graph'$. 

For~(2) we construct the free clone $\freeClone{\graph}$ on a multigraph 
$\graph$. The construction is similar to that for the free 
multicategory on a 
multigraph~(\cf~\cite[Chapter 2]{Leinster2004}).
The sorts are the nodes of $\graph$, and the operations are given by 
the following 
deductive system:
\begin{center}
\unaryRule
	{c \in \graph(X_1, \,\dots\, , X_n; Y)}
	{c \in \freeClone{\graph}(X_1, \,\dots\, , X_n; Y)}
	{}  \qquad
\unaryRule
	{X_i \in \{ X_1, \,\dots\, , X_n \}}
	{\p{i}{X_1, \,\dots\, , X_n} \in \freeClone{\graph}(X_1, \,\dots\, , X_n; 
	X_i)}
	{} \\

\binaryRule
	{f \in \freeClone{\graph}(X_1, \,\dots\, , X_n; Y)}
	{\big(g_i \in \freeClone{\graph}(\Gamma; X_i)\big)_{i=1,\dots,n}}
	{\cslr{f}{g_1, \,\dots\, , g_n} \in \freeClone{\graph}(\Gamma; Y)}
	{}
\end{center}
subject to the equational theory requiring the three axioms of a clone. To see 
this is free,  observe that for any clone $(S, \clone)$ and multigraph 
homomorphism
$h : \graph \to \clone$ 
from $\graph$ to the multigraph underlying $(S, \clone)$, 
the unique clone homomorphism 
$\ext{h} : \freeClone{\graph} \to \clone$ 
extending $h$ must be defined by
\[
\ext{h}(c) := h(c) \qquad 
\ext{h}{(\p{i}{\ind{A}})} := \p{i}{\ext{h}{\ind{A}}} \qquad
\ext{h}{\left(\cslr{f}{g_1, \,\dots\, , g_n} \right)} :=
	\cslr{(\ext{h}f)}{(\ext{h}g_1), \,\dots\, , (\ext{h}g_n)}
\]
For~(\ref{c:cat-to-clone}), let $\catC$ be a category. Define a clone
$\prom \catC$ with sorts the objects of $\catC$ and hom-sets constructed as 
follows:
\begin{center}
\unaryRule
	{f \in \catC(X, Y)}
	{f \in (\prom \catC)(X; Y)}
	{} \qquad
\unaryRule
	{X_i \in \{ X_1, \,\dots\, , X_n \}}
	{\p{i}{X_1, \,\dots\, , X_n} \in (\prom \catC)(X_1, \,\dots\, , X_n; X_i)}
	{} \vspace{0.5\treeskip} \\

\binaryRule
	{f \in (\prom \catC)(X_1, \,\dots\, , X_n; Y)}
	{\big(g_i \in (\prom \catC)(\Gamma; X_i)\big)_{i=1,\dots,n}}
	{\cslr{f}{g_1, \,\dots\, , g_n} \in (\prom \catC)(\Gamma; Y)}
	{}
\end{center}
The equational theory $\equiv$ is the three laws of a clone, augmented by
\begin{center}
\unaryRule{\phantom{\p{1}{}}}{\p{1}{X} \equiv \id_X \in (\prom \catC)(X;X)}{}
\qquad
\binaryRule
	{f \in \catC(Y, Z)}
	{g \in \catC(X, Y)}
	{f \circ g \equiv \cs{f}{g} \in (\prom\catC)(X;Z)}
	{} \vspace{-\treeskip}
\end{center}
For any clone $(T, \altClone)$, a 
clone homomorphism
$h : \prom \catC \to \altClone$ 
consists of a map of objects $ob(\catC) \to T$
together with substitution-preserving mappings
$(\prom\catC)(X_1, \,\dots\, , X_n; Y) \to \altClone(X_1, \,\dots\, , X_n; Y)$
for each $X_1, \,\dots\, , X_n, Y \in ob(\catC) \:\: (n \in \Nat)$.  
Restricting to 
unary operations, this is exactly a functor $\catC \to \overline\altClone$. 
Conversely, since any clone homomorphism is fixed on the projections, a functor 
$\catC \to \overline\altClone$ corresponds uniquely to a clone homomorphism
$\prom\catC \to \altClone$.
\end{proof}
\end{mylemma}

In the light of the preceding lemma one obtains the diagram below. The 
adjunction between the 1-category $\CatCat$ 
and $\Graph$ 
is the usual 
free-forgetful adjunction, and the functor
$\overline{(-)} : \Clone \to \CatCat$ 
restricts a clone $(S, \clone)$ to its unary operations (\ie~its nucleus).
The outer square commutes on the nose and hence the inner square 
commutes up to natural isomorphism.
\begin{equation} \label{eq:clone-cat-multigraph-diagram}
\begin{tikzcd}[column sep = 5em, row sep=2.5em]
\: &
\Clone 
\arrow[bend left]{dr}{\overline{(-)}}
\arrow[bend right]{dl}[swap]{\text{forget}} &
\: \\
\MultiGraph
\arrow[phantom]{ur}[description]{\adjDown}
\arrow[bend right = 23, yshift=.5mm]{ur}[swap]{\freeClone{-}}
\arrow[bend right]{dr}[swap]{\lin} &
\: &
\CatCat
\arrow[phantom]{ul}[description]{\adjDown}
\arrow[bend left = 23, yshift=.5mm]{ul}[near end]{\prom}
\arrow[bend left]{dl}{\text{forget}} \\
\: &
\Graph 
\arrow[phantom]{ur}[description]{\adjUp}
\arrow[bend left, yshift=-.5mm]{ur}[near start]{\mathbb{FC}\mathrm{at}}
\arrow[phantom]{ul}[description]{\adjUp}
\arrow[bend right = 23, hookrightarrow, yshift=-.5mm]{ul} &
\:
\end{tikzcd} 
\end{equation}
Indeed, examining the constructions one sees that 
$\nucleus{(-)} \circ \prom \iso \id_{\CatCat}$ and hence that 
\begin{equation} \label{eq:nucleus-gives-free-cat}
\CatCat(\mathbb{FC}\mathrm{at}(\graph), \catC)
	\iso 
\CatCat{\left(\nucleus{\prom(\mathbb{FC}\mathrm{at}(\graph))}, \catC\right)}
	\iso 
\CatCat(\nucleus{\freeClone{\inc\graph}}, \catC)
\end{equation}
For our purposes, the moral is the following: to provide a 
type-theoretic description of the free category on a graph, 
it is sufficient to describe the free clone on a multigraph. One thereby 
obtains a more natural type theory---one does not need to restrict the rules to 
unary contexts---and the commutativity of this diagram guarantees that, when 
one does perform such a restriction, the result is
(up to isomorphism) as intended.

Our aim in what follows is to lift this story to the bicategorical setting, and 
use it to extract a type theory for bicategories. We begin by defining a 
bicategorified notion of clone.

\paragraph*{Biclones.}
Abstract clones may be defined in any cartesian category (and much 
more generally, see~\cite{Staton2013, Fiore2017}). The bicategorified version 
arises by instantiating this definition in $\Cat$ and weakening 
the axioms to natural isomorphisms. 

\begin{mydefn} 
A \Def{(sorted) biclone} $(S, \biclone)$ is a set $S$ of \Def{sorts} equipped with the following data:
\begin{itemize} 
\item For all $X_1, \,\dots\, , X_n, Y \in S \:\: (n \in \Nat)$ a category 
$\biclone(X_1, \,\dots\, , X_n;Y)$ 
with objects \Def{multimaps} ${f : \ind{X} \to Y}$ and morphisms 
\Def{2-cells} 
\mbox{$\alpha : f \To g : \ind{X} \to Y$}, subject to a 
\Def{vertical composition} operation,
\item Distinguished \Def{projection} functors $\p{i}{\ind{X}} : \catOne \to 
\biclone(X_1, \,\dots\, , X_n; X_i) \: (i = 1, \,\dots\, , n)$ for all 
$X_1, \,\dots\, , X_n \in S \:\: (n \in \Nat)$, 
\item For all sequences of sorts $\Gamma$ and sorts 
$Y_1, \,\dots\, , Y_n, Z \:\: (n \in \Nat)$ a 
\Def{substitution} functor
\begin{gather*}
\subName_{\Gamma, \ind{Y}, Z} : \biclone(\ind{Y}; Z) \times \prodop_{i=1}^n 
\biclone(\Gamma; Y_i) \to \biclone(\Gamma; Z)
\end{gather*}
we denote by $\subName\big(f, (g_1, \,\dots\, , g_n)\big) := \cslr{f}{g_1, 
\,\dots\, , 
g_n}$,
\item Natural families of invertible \Def{structural isomorphisms} 
\begin{align*}
\assoc{t, \ind{u}, \ind{v}} : \cslr{\cslr{t}{u_1, \,\dots\, , u_n}}{\ind{v}} 
&\To 
\cslr{t}{\cslr{u_1}{\ind{v}}, \,\dots\, , \cslr{u_n}{\ind{v}}}  \\
\subid{u} : u &\To \cslr{u}{\p{1}{\ind{X}}, \,\dots\, , \p{n}{\ind{X}}} \\
\indproj{k}{u_1, \,\dots\, , u_n} : \cslr{\p{k}{\ind{Y}}}{u_1, \,\dots\, , u_n} 
&\To u_k 
\quad (k = 1,\dots, n)
\end{align*}
for every $t \in \biclone(\ind{Y}, Z)$, \mbox{$u_j \in \biclone(\ind{X}, 
Y_j)$}, \mbox{$v_i \in \biclone(\ind{W}, X_i)$} and $u \in \biclone(\ind{X}, 
Y)$ ($i = 1, \,\dots\, , n$ and $j = 1, \,\dots\, , m$),
\end{itemize}
This data is subject to coherence laws corresponding to the triangle and 
pentagon laws of a bicategory:
\begin{td}[ampersand replacement = \&, column sep = 7em] 
\cslr{t}{\ind{v}} 
\arrow{r}[yshift=0mm]{\cslr{\subid{t}}{\ind{v}}} 
\arrow[swap, equals]{d}{} \& 
\cslr{\cslr{t}{\p{1}{}, \,\dots\, , \p{n}{}}}{\ind{v}} 
\arrow{d}{\assoc{t; \p{\bullet}{}; \ind{v}} } \\ 
\cslr{t}{\ind{v}} \& 
\cslr{t}{\cslr{\p{1}{}}{\ind{v}}, \,\dots\, , \cslr{\p{n}{}}{ \ind{v}}} 
\arrow{l}[yshift=0mm]
	{\cslr{t}{\indproj{1}{\ind{v}}, \,\dots\, , \indproj{n}{\ind{v}}}} 
\end{td} 
\begin{td}[column sep = 7em]
\cslr{\cslr{\cslr{t}{\ind{u}}}{\ind{v}}}{\ind{w}} 
\arrow[swap]{d}{\assoc{\cs{t}{\ind{u}}; \ind{v}; \ind{w}}} 
\arrow{r}[yshift=0mm]{\cslr{\assoc{t; \ind{u}; \ind{v}}}{\ind{w}}} &
\cslr{\cslr{t}{\cslr{\ind{u}}{\ind{v}}}}{\ind{w}} 
\arrow[]{r}[yshift=0mm]{\assoc{t; \cs{\ind{u}}{\ind{v}}; \ind{w}}} &
\cslr{t}{\cslr{\cslr{\ind{u}}{\ind{v}}}{\ind{w}}} 
\arrow{d}
	{\cslr{t}
		{\assoc{\ind{u}; \ind{v}; \ind{w}}}} \\
\cslr{\cslr{t}{\ind{u}}}{\cslr{\ind{v}}{\ind{w}}} 
\arrow[swap]{rr}[yshift=-0mm]
	{\assoc{t; \ind{u}; \cs{\ind{v}}{\ind{w}}}} &
\: & 
\cslr{t}{\cslr{\ind{u}}{\cslr{\ind{v}}{\ind{w}}}}
\end{td}	
\end{mydefn} 

\begin{myremark} \label{rem:biclones-and-invertibility}
Note that an invertible 2-cell is simply an iso in the relevant hom-category, 
but the definition of invertible multimap is more subtle (see 
Definition~\ref{def:clone-1-cell-invetibility}).
\end{myremark}

We direct the 2-cells to match the definition of a \Def{skew monoidal 
category}~\cite{Szlachanyi2012}; the definition should therefore generalise to 
the lax setting.  When we wish to emphasise the set of sorts, we call a biclone 
$(S, \biclone)$ an \Def{$S$-biclone}; where the set of sorts is clear from 
context, we refer to a biclone $(S, \biclone)$ simply by $\biclone$. One 
obtains a \Def{2-clone}---a clone enriched over $\CatCat$---when all 
the structural isomorphisms $\assoc{}, \subid{}, \indproj{i}{} \: 
(i=1,\dots,n)$ are the identity. The second half of this chapter will be 
devoted to a coherence theorem showing that every freely-generated biclone is 
suitably equivalent to a 2-clone. 

\begin{myexmp}[{\cf~Example~\ref{ex:CartesianCategoryIsAClone}}] 
\label{ex:CartesianBicatABiclone} 
\quad
\begin{enumerate}
\item Every clone defines a \Def{locally discrete} biclone, in which each 
hom-category is discrete.
\item Every bicategory $\baseCat$ with finite products defines a biclone; if 
$\baseCat$ is a 2-category with strict (2-categorical) products, this is a 
2-clone.
\item \label{c:bicat-as-biclone} Every biclone $(S, \biclone)$ gives rise to a 
bicategory $\overline{\biclone}$ by taking the \Def{unary} 
hom-categories,~\ie~by taking $\overline{\biclone}(X,Y) := \biclone(X;Y)$. We call this
the \Def{nucleus} of $(S, \biclone)$.  
\qedhere
\end{enumerate}
\end{myexmp} 

One may think of a biclone as a generalised deductive system in 
which the multimaps $f : A_1, \,\dots\, , A_n \to B$ are judgements 
$A_1, \,\dots\, , A_n \vdash f : B$, related by proof transformations
$\tau : f \To f'$ (\cf~\cite{Seely1987}). Conversely, 
Example~\ref{ex:CartesianBicatABiclone}(\ref{c:bicat-as-biclone}) shows that a 
type theory for biclones would encompass bicategories as a special case. 
In Lemma~\ref{lem:free-bicategory-on-a-multigraph} we shall see 
that the type theory 
describing the free biclone on a 2-graph restricts to a type theory for the 
free bicategory on a 2-graph~(\cf~diagram~(\ref{eq:clone-cat-multigraph-diagram})).

\begin{myremark}
Biclones are objects worthy of further study in their own right. Thinking of 
them as `bicategorified clones' suggests a connection---to be fleshed 
out---with some 
notion of `bicategorical Lawvere theory', and with pseudomonads. On the other 
hand, biclones provide a categorical description of certain kinds of explicit 
substitution; possible connections with the categorical semantics of the 
simply-typed lambda calculus with explicit substitution (\eg~\cite{Ghani1999}) 
remain to be explored. 
\end{myremark}

\newpage
\paragraph{Free biclones and free bicategories.}
Defining a free biclone requires an appropriate notion 
of morphism.
The definitions are natural extensions of those for bicategories.

\begin{mydefn} \label{def:biclone-pseudofunctor}
A \Def{pseudofunctor} $F: (S, \biclone) \to (S', \biclone')$ between biclones 
consists of a mapping $F : ob(\biclone) \to ob(\biclone')$ equipped with:
\begin{itemize} 
\item A functor \mbox{$F_{\ind{X};Y} : \biclone(X_1, \,\dots\, , X_n; Y) \to 
\biclone'(FX_1, \,\dots\, , FX_n; FY)$} for all 
${X_1, \,\dots\, , X_n, Y \in S} \:\: {(n \in \Nat)}$,
\item Invertible 2-cells \mbox{$\psi^{(i)}_{\ind{X}} : \p{i}{F\ind{X}} \To 
F(\p{i}{\ind{X}})$} $(i = 1,\dots, n)$ for each $X \in S$, 
\item An invertible 2-cell $\phi_{t,\ind{u}} : \cslr{(Ft)}{Fu_1, \,\dots\, , 
Fu_n} 
\To 
F(\cslr{t}{u_1, \,\dots\, , u_n})$ for every \mbox{$(u_j : \ind{X} \to Y_i)_{j 
= 1, 
\dots, n}$} and $t : \ind{Y} \to Z$, natural in $t$ and $u_1, \,\dots\, , u_n$,
\end{itemize} 
subject to the following three coherence laws for 
$i=1, \,\dots\, , n$:
\begin{equation} \label{eq:pseudofunctor:left-unit}
\begin{tikzcd}[ampersand replacement = \&]
\cslr{\p{i}{F\ind{X}}}{Fu_1, \,\dots\, , Fu_n} 
\arrow{r}{\indproj{i}{F\ind{u}}} 
\arrow[swap]{d}{\cslr{\psi^{(i)}_{\ind{X}}}{F\ind{u}}} \&
Fu_i \\
\cslr{(F\p{i}{\ind{X}})}{F\ind{u}} 
\arrow[swap]{r}{\phi_{\p{i}{}, \ind{u}}} \&
F(\cslr{\p{i}{\ind{X}}}{\ind{u}}) 
\arrow[swap]{u}{F\indproj{i}{\ind{u}}}
\end{tikzcd}
\end{equation} 
\begin{equation} \label{eq:pseudofunctor:right-unit}
\begin{tikzcd}[column sep = 8em, ampersand replacement = \&]
F(t)\arrow[]{r}{F\subid{t}} 
\arrow[swap]{d}{\subid{Ft}} \&
F{\left(\cs{t}{\p{1}{\ind{X}}, \,\dots\, , \p{n}{\ind{X}}}\right)} \\
\cs{(Ft)}{\p{1}{F\ind{X}}, \,\dots\, , \p{n}{F\ind{X}}} 
\arrow[swap]{r}[yshift=-0mm]{\cs{(Ft)}{\psi^{(1)}, \,\dots\, , \psi^{(1)}}} \&
\cs{(Ft)}{F\p{1}{\ind{X}}, \,\dots\, , F\p{n}{\ind{X}}} 
\arrow[swap]{u}{\phi_{t; \p{\bullet}{}}}
\end{tikzcd}
\end{equation} 
\vspace{3mm}
\begin{equation} \label{eq:pseudofunctor:assoc}
\begin{tikzcd}[column sep = 8em]
\csthree{F(t)}{F\ind{u}}{F\ind{v}} 
\arrow{r}{\assoc{Ft; F\ind{u}; F\ind{v}}} 
\arrow[swap]{d}{\cslr{\phi_{t; \ind{u}}}{F\ind{v}}}  &
\cslr{F(t)}{\cslr{F\ind{u}}{F\ind{v}}} 
\arrow{d}{\cslr{F(t)}{\phi_{\ind{u}; \ind{v}}}} \\
\cslr{F(\cslr{t}{\ind{u}})}{F\ind{v}} 
\arrow[swap]{d}{\phi_{\cs{t}{\ind{u}}; \ind{v}}}  & 
\cslr{F(t)}{F(\cslr{\ind{u}}{\ind{v}})} 
\arrow{d}{\phi_{t; \cs{\ind{u}}{\ind{v}}}} \\
F(\csthree{t}{\ind{u}}{\ind{v}}) 
\arrow[swap]{r}{F\assoc{t; \ind{u}; \ind{v}}} &
F(\cslr{t}{\cslr{\ind{u}}{\ind{v}}})
\end{tikzcd}
\end{equation}
A pseudofunctor for which $\phi$ and every $\psi^{(1)}, \dots, \psi^{(n)}$ is 
the identity is called \Def{strict}.
\end{mydefn}

\begin{myexmp} \label{exmp:pseudofunctors-restrict}
Every pseudofunctor of biclones $F : (S, \biclone) \to (T, \altBiclone)$ restricts
to a pseudofunctor of bicategories $\overline{F} : \overline\biclone \to \overline\altBiclone$
between the nucleus of $(S, \biclone)$ and the nucleus of $(T, \altBiclone)$ 
(recall Example~\ref{ex:CartesianBicatABiclone}(\ref{c:bicat-as-biclone})).  
\end{myexmp}

The construction of the free biclone on a 2-multigraph follows the pattern of 
its 1-categorical counterpart.

\begin{myconstr}[Free biclone on a 2-multigraph] \label{constr:free-biclone}
Let $\graph$ be a 
\mbox{2-multigraph}. Define a biclone $\freeBiclone{\graph}$ as follows. The 
sorts are 
nodes of $\graph$ and the hom-categories are defined by the following 
deductive system:
\begin{center}
\unaryRule{c \in \graph(A_1, \,\dots\, , A_n; B)}{c \in 
\freeBiclone{\graph}(A_1, \dots, A_n; B)}{} \vspace{-0.5\treeskip}
\qquad
\unaryRule
	{\constrewr \in \graph(A_1, \,\dots\, , A_n; B)(c,c')}
	{\constrewr \in \freeBiclone{\graph}(A_1, \,\dots\, , A_n; B)}
	{} \\

\unaryRule
	{}
	{\p{i}{A_1, \,\dots\, , A_n} \in \freeBiclone{\graph}(A_1, \,\dots\, , A_n; 
	A_i)}
	{$(1 \leq i \leq n)$} \vspace{0.5\treeskip} \\

\binaryRule
	{f \in \freeBiclone{\graph}(A_1, \,\dots\, , A_n; B)}
	{\big(g_i \in \freeBiclone{\graph}(\ind{X};A_i)\big)_{i=1,\dots,n}}
	{\cslr{f}{g_1, \,\dots\, , g_n} \in \freeBiclone{\graph}(\ind{X}; B)}
	{} \vspace{0.5\treeskip} 

\binaryRule
	{\tau \in \freeBiclone{\graph}(A_1, \,\dots\, , A_n; B)(f,f')}
	{\big(\sigma_i \in 
		\freeBiclone{\graph}(\ind{X}; A_i)(g_i,g_i')\big)_{i=1,\dots,n}}
	{\cslr{\tau}{\sigma_1, \,\dots\, , \sigma_n} \in 
		\freeBiclone{\graph}(\ind{X}; B)
			(\cslr{f}{g_1, \,\dots\, , g_n}, \cslr{f'}{g_1', \,\dots\, , g_n'})}
	{} \vspace{\treeskip}

\unaryRule
	{f \in \freeBiclone{\graph}(\ind{A}; B)}
	{\id_f \in \freeBiclone{\graph}(\ind{A}; B)(f,f)}
	{} 
\quad
\binaryRule
	{\tau \in \freeBiclone{\graph}(\ind{A}; B)(f', f'')}
	{\sigma \in \freeBiclone{\graph}(\ind{A}; B)(f, f')}
	{\tau \vert \sigma \in \freeBiclone{\graph}(\ind{A}; B)(f,f'')}
	{} 
	
\trinaryRule
	{f \in \freeBiclone{\graph}(\ind{B}; C)}
	{\big( g_i \in \freeBiclone{\graph}(\ind{A}; B_i) \big)_{i=1,\dots,n}}
	{\big( h_j \in \freeBiclone{\graph}(\ind{X}; B_j) \big)_{j=1,\dots,m}}
	{\assoc{f, \ind{g}, \ind{h}} 
		\in \freeBiclone{\graph}(\ind{X}; C)
			(\csthree{f}{\ind{g}}{\ind{h}}, \cs{f}{\cs{\ind{g}}{\ind{h}}})}
	{} \vspace{0.5\treeskip}

\unaryRule
	{f \in \freeBiclone{\graph}(A_1, \,\dots\, , A_n; B)}
	{\subid{f} \in \freeBiclone{\graph}(\ind{A}; B)
		{\left(f, \cs{f}{\p{1}{\ind{A}}, \,\dots\, , \p{n}{\ind{A}}}\right)}}
	{} \vspace{0.5\treeskip}

\unaryRule
	{\big(g_i \in \freeBiclone{\graph}(\ind{X}; A_i)\big)_{i=1,\dots,n}}
	{\indproj{i}{A_1, \,\dots\, , A_n} \in 
		\freeBiclone{\graph}(\ind{X}; A_i)
			(\cslr{\p{i}{A_1, \,\dots\, , A_n}}{g_1, \,\dots\, , g_n}, g_i)}
	{$(1 \leq i \leq n)$}
\end{center}
The equational theory $\equiv$ requires that
\begin{itemize}
\item Every $\freeBiclone{\graph}(A_1, \,\dots\, , A_n; B)$ forms a category 
with composition the $\vert$ operation and identity on 
$f \in \freeBiclone{\graph}(A_1, \,\dots\, , A_n; B)$ given by $\id_f$,
\item The operation 
$\big(f, (g_1, \,\dots\, , g_n)\big) \mapsto \cslr{f}{g_1, \dots, g_n}$ is 
functorial with respect to this category structure, 
\item The families of 2-cells 
$\assoc{}, \subid{}$ and $\indproj{i}{} \: (i=1,\dots,n)$ are invertible, 
natural and satisfy the triangle and pentagon 
laws of a biclone.  \qedhere
\end{itemize}
\end{myconstr}

It is clear that this construction yields a 
biclone. Indeed, Lambek's definition of 
the internal language of a multicategory~\cite{Lambek1989} transfers readily to 
clones, and the preceding construction 
may be used to extend this definition to biclones. The only 
adjustment is that the operation symbols 
$f: A_1, \,\dots\, , A_n \to B$ 
are now related by transformations $\tau  : f \To f'$. The judgements in 
our type theory $\langBiclone$ will match these sequents precisely. 

We shall, so far as possible, phrase the free properties we prove in terms of a 
unique strict pseudofunctor of biclones~(\cf~\cite[Proposition 2.10]{Gurski2013}): this obviates the need to work with 
uniqueness up to 2-cell, in which the 2-cells may themselves only be unique up 
to a unique 3-cell. In particular, we 
bicategorify diagram~(\ref{eq:clone-cat-multigraph-diagram}) by using 
1-categories of bicategorical objects (biclones and bicategories) in which the 
morphisms are \emph{strict} pseudofunctors. Write $\BicloneCat$ and $\BicatCat$ 
for these two categories.
The relevant freeness universal property of 
Construction~\ref{constr:free-biclone} is therefore the following.

\begin{mylemma} \label{lem:free-biclone}
The forgetful functor $\BicloneCat \to \TwoMultiGraph$ taking a biclone to its 
underlying 2-multigraph has a left adjoint.
\begin{proof}
Let $\graph$ be a 2-multigraph and $(T, \altClone)$ be a biclone. We show that 
for every 2-multigraph morphism $h : \graph \to \altClone$ there exists a 
unique strict pseudofunctor of biclones 
$h^\sharp : \freeBiclone{\graph} \to \graph$ such 
that 
$h^\sharp \circ \inc = h$, for $\inc : \graph \to \freeBiclone{\graph}$ the 
inclusion.

Define $\ext{h}$ by induction as follows:
\begin{align*}
\ext{h}(c) &:= h_{\ind{A}; B}(c) \qquad\quad \text{ for } c 
\in \graph(A_1, \,\dots\, , A_n; B) \\
\ext{h}(\constrewr) &:= h_{\ind{A}; B}(\constrewr) 
\qquad\quad \text{ for } 
\constrewr \in \graph(A_1, \,\dots\, , A_n; B)(c,c') \\ 
\ext{h}(\id_f) &:= \id_{\ext{h}(f)} \\
\ext{h}(\tau \vert \sigma) &:= \ext{h}(\tau) \vert \ext{h}(\sigma)
\end{align*}
We then require that $\ext{h}$ strictly preserves the projections, the 
substitution operations and the structural isomorphisms. This is a strict 
pseudofunctor $\freeBiclone{\graph} \to \altClone$ extending $h$. 
Uniqueness follows 
because any strict pseudofunctor must strictly preserve projections and the 
substitution operations, and so also strictly preserve the structural 
isomorphisms.
\end{proof}
\end{mylemma}

The proof of 
Lemma~\ref{lem:graph-multi-graph-adj-and-free-clone} extends straightforwardly 
to an adjunction between $\TwoGraph$ and $\TwoMultiGraph$. The following lemma 
therefore completes our bicategorical adaptation 
of diagram~(\ref{eq:clone-cat-multigraph-diagram}). 

\newpage
\begin{mylemma} \quad \label{lem:free-bicategory-on-a-multigraph}
\begin{enumerate}
\item \label{c:free-bicat} The forgetful functor $\BicatCat \to \TwoGraph$ 
taking a bicategory to its underlying 2-graph has a left adjoint
(\cf~\cite[Proposition 2.10]{Gurski2013}).
\item \label{c:free-biclone}  The functor 
$\overline{(-)}: \BicloneCat \to \BicatCat$ 
restricting a biclone to its nucleus (recall Example~\ref{ex:CartesianBicatABiclone}) has a left 
adjoint.
\end{enumerate}
\begin{proof}
For~(\ref{c:free-bicat}) we define the free bicategory $\freeBicat{\graph}$ on 
a  
2-graph $\graph$ as the 
following deductive system (\cf~the description of bicategories as a 
generalised algebraic theory~\cite{Ouaknine1997}):
\begin{figure}[!h]
\centerfloat
\unaryRule{c \in \graph(A, B)}{c \in \freeBicat{\graph}(A, B)}{} 
\vspace{-0.5\treeskip}
\quad
\unaryRule
	{\constrewr \in \graph(A, B)(c,c')}
	{\constrewr \in \freeBicat{\graph}(A, B)}
	{} 
\quad
\unaryRule
		{\phantom{\constrewr \in \freeBicat{\graph}(A, B)}}
		{\Id_A \in \freeBicat{\graph}(A, A)}
		{} \vspace{-0.5\treeskip}

\binaryRule
	{f \in \freeBicat{\graph}(A, B)}
	{g \in \freeBicat{\graph}(X; A)}
	{f \circ g \in \freeBicat{\graph}(X; B)}
	{} 
\binaryRule
	{\tau \in \freeBicat{\graph}(A, B)(f,f')}
	{\sigma \in \freeBicat{\graph}(X, A)(g,g')}
	{\tau \circ \sigma \in \freeBicat{\graph}(X; B)(f \circ g, f' \circ g')}
	{} \vspace{0.5\treeskip}

\unaryRule
	{f \in \freeBicat{\graph}(A, B)}
	{\id_f \in \freeBicat{\graph}(A, B)(f,f)}
	{}
\quad
\binaryRule
	{\tau \in \freeBicat{\graph}(A, B)(f', f'')}
	{\sigma \in \freeBicat{\graph}(A, B)(f, f')}
	{\tau \vert \sigma \in \freeBicat{\graph}(A, B)(f,f'')}
	{}

\trinaryRule
	{f \in \freeBicat{\graph}(B, C)}
	{g \in \freeBicat{\graph}(A, B)}
	{h \in \freeBicat{\graph}(X, B)}
	{\a_{f;g;h}
		\in \freeBiclone{\graph}(X; C)
			(\csthree{f}{g}{h}, \cs{f}{\cs{g}{h}})}
	{} \vspace{0.5\treeskip}

\unaryRule
	{f \in \baseCat(A, B)}
	{\l_f \in \freeBicat{\graph}(A,B)(\Id_B \circ f, f)}
	{} \vspace{-\treeskip}
\quad
\unaryRule
	{f \in \freeBicat{\graph}(A, B)}
	{\r_f \in \freeBicat{\graph}(A, B)\left(f \circ \Id_A, f\right)}
	{} \vspace{-\treeskip}
\end{figure}

subject to an equational theory requiring
\begin{itemize}
\item Every $\freeBicat{\graph}(A, B)$ forms a category with 
composition the $\vert$ operation and identity on 
$f \in \freeBicat{\graph}(A, B)$ given by $\id_f$,
\item The operation 
$(f,g) \mapsto f \circ g$ 
is functorial with respect to this category 
structure, 
\item The families of 2-cells $\a, \l$ and $\r$ are invertible, natural and 
satisfy 
the triangle and pentagon laws of a bicategory.  
\end{itemize}
Since strict pseudofunctors are determined on all the structural data, any 
2-graph homomorphism 
$h : \graph \to \altCat$ to the 2-graph 
underlying a 
bicategory $\altCat$ determines a unique strict pseudofunctor 
$\ext{h} : \freeBiclone{\graph} \to \altCat$ restricting to $h$ on $\graph$.

For~(\ref{c:free-biclone}), let $\baseCat$ be any bicategory. Define a biclone
$\prom \baseCat$ as follows. The sorts are objects of $\baseCat$ and the 
hom-categories 
$(\prom\baseCat)(X_1, \,\dots\, , X_n; Y)$ are those given by the deductive 
system of
Construction~\ref{constr:free-biclone}, adapted by replacing the first two 
rules by 
\begin{center}
\unaryRule{f \in \baseCat(X, Y)}{f \in (\prom\baseCat)(X;Y)}{} 
\vspace{-\treeskip}
\qquad\qquad
\unaryRule
	{\constrewr \in \baseCat(X,Y)(f,f')}
	{\constrewr \in (\prom\baseCat)(X;Y)(f,f')}
	{}  \vspace{-\treeskip}
\end{center}
and augmenting the equational theory with rules ensuring the biclone and 
bicategory structures coincide wherever possible:
\begin{center}
\unaryRule
		{\phantom{f \in \baseCat(Y, Z)}}
		{\p{1}{X} \equiv \Id_X \in (\prom \baseCat)(X;X)}{}
\binaryRule
	{f \in \baseCat(Y, Z)}
	{g \in \baseCat(X, Y)}
	{f \circ g \equiv \cs{f}{g} \in (\prom\baseCat)(X;Z)}
	{} 
\unaryRule
	{f \in \baseCat(X, Y)}
	{(\id_f)_{\baseCat} \equiv (\id_f)_{\prom\baseCat} \in (\prom\baseCat)(X; Y)}
	{} \vspace{0.5\treeskip}

\binaryRule
	{\tau \in \baseCat(Y, Z)(f,f')}
	{\sigma \in \baseCat(X, Y)(g,g')}
	{\tau \circ \sigma \equiv \cs{\tau}{\sigma} 
		\in (\prom\baseCat)(X;Z)(\cs{f}{g}, \cs{f'}{g'})}
	{} \quad
\binaryRule
	{\tau \in \baseCat(X, Y)(f,f')}
	{\sigma \in \baseCat(X, Y)(f',f'')}
	{\tau \vert_\baseCat \sigma \equiv \tau \vert_{\prom\baseCat} \sigma 
		\in (\prom\baseCat)(X;Y)(f, f'')}
	{} \vspace{0.5\treeskip}
	
\trinaryRule
	{f \in \freeBicat{\graph}(B, C)}
	{g \in \freeBicat{\graph}(A, B)}
	{h \in \freeBicat{\graph}(X, B)}
	{\assoc{f,g,h} \equiv \a_{f, g, h} \in \freeBicat{\graph}(X, C)}
	{} \vspace{0.5\treeskip}

\unaryRule
	{f \in \baseCat(X, Y)}
	{\subid{f} \equiv \r_f^{-1} : (\prom\baseCat)(X, Y)(f, \cs{f}{\p{1}{X})}}
	{}
\unaryRule
	{f \in \baseCat(X, Y)}
	{\indproj{1}{f} \equiv \l_f : 
		(\prom\baseCat)(X, Y)(\cs{\p{1}{Y}}{f}, f)}
	{} \vspace{-\treeskip}
\end{center}
The free property is a simple extension of that for clones  
(Lemma~\ref{lem:graph-multi-graph-adj-and-free-clone}(\ref{c:cat-to-clone})).
\end{proof}
\end{mylemma} 
One therefore obtains the following diagram of 
adjunctions, 
generalising diagram~(\ref{eq:clone-cat-multigraph-diagram}). 
As for~(\ref{eq:clone-cat-multigraph-diagram}), the outer diagram commutes on the 
nose so the inner diagram commutes up to isomorphism.
\begin{equation} \label{eq:biclone-bicat-multigraph-diagram}
\begin{tikzcd}[column sep = 5.5em, row sep=3.5em]
\: &
\BicloneCat 
\arrow[bend left]{dr}{\overline{(-)}}
\arrow[bend right]{dl}[swap]{\text{forget}} &
\: \\
\TwoMultiGraph
\arrow[phantom]{ur}[description]{\adjDown}
\arrow[bend right = 20, yshift=.5mm]{ur}[swap]{\freeBiclone{-}}
\arrow[bend right]{dr}[swap]{\lin} &
\: &
\BicatCat
\arrow[phantom]{ul}[description]{\adjDown}
\arrow[bend left = 20, yshift=.5mm]{ul}[near end]{\prom}
\arrow[bend left]{dl}{\text{forget}} \\
\: &
\TwoGraph 
\arrow[phantom]{ur}[description]{\adjUp}
\arrow[bend left, yshift=-.5mm]{ur}{\freeBicat{-}}
\arrow[phantom]{ul}[description]{\adjUp}
\arrow[bend right = 20, yshift=-.5mm, hookrightarrow]{ul} &
\:
\end{tikzcd} 
\end{equation}
It follows that, modulo a natural isomorphism, 
the free bicategory on a 2-graph $\graph$ is obtained as the nucleus of
the free biclone on $\graph$ (regarded as a 2-multigraph). 
Indeed, examining the constructions one sees that 
$\nucleus{(-)} \circ \prom \iso \id_{\BicatCat}$, yielding the following 
chain of natural isomorphisms (\cf~equation~(\ref{eq:nucleus-gives-free-cat})):
\begin{equation} \label{eq:nucleus-gives-free-bicat}
\BicatCat(\freeBicat{\graph}, \baseCat)
	\iso 
\BicatCat{\left(\nucleus{\prom(\freeBicat{\graph})}, \baseCat\right)}
	\iso 
\BicatCat(\nucleus{\freeBiclone{\inc\graph}}, \baseCat)
\end{equation}
For us, the moral is the following:
Construction~\ref{constr:free-biclone} gives precisely the rules required to 
freely define bicategorical substitution structure. In Section~\ref{sec:type-theory-for-biclones}, we shall use this to construct a type theory for bicategories. Before that, we finish giving 
the definitions required to specify an equivalence of biclones. These will be a 
key part of the coherence result at the end of this chapter. 

\paragraph{Relating biclone pseudofunctors.} The definition of transformation 
between biclone 
homomorphisms is 
rather involved. There is a well-known notion of transformation between maps of 
multicategories (\eg~\cite[Definition 2.3.5]{Leinster2004}), but the cartesian 
nature of biclone substitution means the definition is not directly applicable. 
However, every clone canonically gives rise to a multicategory---we discuss 
this in some detail in 
Section~\ref{sec:product-structure-from-rep}---and this suggests the 
definition of transformation should be a bicategorical adaptation of that 
for multicategory maps. 
The definition of modification is then fixed.

The following notation is intended to be reminiscent of the notation 
$f \times g$ for the action of the categorical 
cartesian product on morphisms.  

\begin{mynotation} \label{not:clone-product-notation}
For multimaps $(f_i : \Gamma_i \to Y_i)_{i=1,\dots, n}$ and in a (bi)clone, 
one obtains the composite
\[
\Gamma_1, \,\dots\, , \Gamma_n \xra{[ \p{1+\sum_{i=1}^{k-1} \len{\Gamma_i}}{}, 
\,\dots\, , \p{\len{\Gamma_k} + \sum_{i=1}^{k-1} \len{\Gamma_i}}{}]} \Gamma_k 
\xra{f_k} Y_k
\]
for $k=1,\dots,n$. 
For $h : Y_1, \dots Y_n \to Z$ we therefore define 
$\cslr{h}{\bigclonetimes_{i=1}^n f_i} = 
\cslr{h}{f_1 \clonetimes \,\cdots\, \clonetimes f_n} : 
\Gamma_1, \,\dots\, , \Gamma_n \to 
Z$ 
to be 
the composite 
\[
\cslr{h}{\cslr{f_1}{\p{1}{}, \,\dots\, , \p{\len{\Gamma_1}}{}}, 
\dots, \cslr{f_n}{\p{1 + \sum_{i=1}^{n-1} \len{\Gamma_i}}{}, \,\dots\, , 
\p{\len{\Gamma_n} + \sum_{i=1}^{n-1} \len{\Gamma_i}}{}}}
\]
\end{mynotation}

In particular, for $(g_j : \Gamma \to X_j)_{j=1,\dots,n}$, $(f_i : X_i \to 
Y_i)_{i=1,\dots,n}$ and $h : Y_1, \,\dots\, , Y_n \to Z$ there exists a 
canonical isomorphism 
\[
\mathsf{f}_{h; \ind{f}; \ind{g}} : 
	\csthree{h}
			{f_1 \clonetimes \,\cdots\, \clonetimes f_n}
			{g_1, \,\dots\, , g_n} \To 
	\cslr	{h}
		{\cslr{f_1}{g_1}, \,\dots\, , \cslr{f_n}{g_n}}
\]
given by applying $\assoc{}$ twice and then the projections $\indproj{i}{}$.
\nom{\mathsf{f}_{h; \ind{f}; \ind{g}}}
	{The canonical 2-cell 
	$\mathsf{f}_{h; \ind{f}; \ind{g}} : 
	\csthree{h}
			{f_1 \times \dots \times f_n}
			{g_1, \,\dots\, , g_n} \To 
	\cslr	{h}
		{\cslr{f_1}{g_1}, \,\dots\, , \cslr{f_n}{g_n}}$ in a biclone}

\begin{mydefn} \label{def:biclone-transformation} 
Let $F, G : (\biclone, S) \to (\biclone', S')$ be pseudofunctors of biclones. A 
\Def{transformation} $(\alpha, \overline{\alpha}) : F \To G$ consists of the 
following data:
\begin{enumerate} 
\item A multimap $\alpha_X : FX \to GX$ for every $X \in S$, 
\item An invertible 2-cell 
\begin{equation} \label{eq:def-of-pseudonatural-transformation}
\overline{\alpha}_t : \cslr{\alpha_Y}{Ft} \To 
	\cslr{G(t)}{\alpha_{X_1} \clonetimes \,\cdots\, \clonetimes \alpha_{X_n}} : 
	FX_1, \,\dots\, , FX_n \to GY
\end{equation}
for every $t : X_1, \,\dots\, , X_n \to Y$ in $\biclone$, natural in $t$ and 
satisfying the following two laws for $k = 1, \,\dots\, , n$: 
\vspace{-2mm}
\begin{equation*}
\makebox[\textwidth]{
\begin{tikzcd}[column sep = small, ampersand replacement = \&]
\csthree{\alpha_Y}{F(t)}{F\ind{u}} 
\arrow[swap]{d}{\cslr{\overline{\alpha}_t}{F\ind{u}}} 
\arrow{r}[yshift=0mm]{\assoc{\alpha; Ft; F\ind{u}}} \& 
\cslr{\alpha_Y}{\cslr{F(t)}{F\ind{u}}} 
\arrow{r}[yshift=2mm]{\cslr{\alpha_Y}{\phi_{t; \ind{u}}}} \&
\cslr{\alpha_Y}{F(\cslr{t}{\ind{u}})}  
\arrow{r}[yshift=0mm]{\overline{\alpha}_{\cslr{t}{\ind{u}}}} \&
\cslr{G(\cslr{t}{\ind{u}})}{\bigclonetimes_{i=1}^n \alpha_{X_i}} \\
\csthree{G(t)}{\bigclonetimes_{i=1}^n \alpha_{X_i}}{F\ind{u}}  
\arrow[swap]{d}{\mathsf{f}_{Gt; \ind{\alpha}; F\ind{u}}} \&
\: \&
\: \&
\: \\
\cslr{G(t)}{\cslr{\alpha_{X_1}}{Fu_1}, \,\dots\, , \cslr{\alpha_{X_n}}{Fu_n}} 
\arrow[swap]{d}{\cslr{G(t)}{\overline{\alpha}_{u_1}, \,\dots\, , 
\overline{\alpha}_{u_n}}} \&
\: \&
\: \&
\: \\
\cslr{G(t)}{\cslr{G(\ind{u})}{\bigclonetimes_{i=1}^n \alpha_{X_i}}} 
\arrow[swap]{rrr}{\assoc{Gt; G\ind{u}; \bigclonetimes_i \alpha_{X_i}}^{-1}} \&
\: \&
\: \&
\csthree{G(t)}{G(\ind{u})}{\bigclonetimes_{i=1}^n \alpha_{X_i}} 
\arrow[swap]{uuu}{\cslr{\phi_{t; \ind{u}}}{\bigclonetimes_{i=1}^n \alpha_{X_i}}}
\end{tikzcd}
}
\end{equation*}
\begin{td}[column sep = 9em, row sep = 3em]
\cslr{\p{k}{G\ind{X}}}{\alpha_{X_1} \clonetimes \,\cdots\, \clonetimes 
\alpha_{X_n}} 
\arrow[swap]{d}{\indproj{k}{(\bigclonetimes_i \alpha_{X_i})}} 
\arrow{r}[]
	{\cslr{\psi^{(k)}_{\ind{X}}}
	{\alpha_{X_1} \clonetimes \,\cdots\, \clonetimes \alpha_{X_n}}} 
	&

\cslr{G(\p{k}{\ind{X}})}{\alpha_{X_1} \clonetimes \,\cdots\, \clonetimes 
\alpha_{X_n}} \\

\cslr{\alpha_{X_k}}{\p{k}{F\ind{X}}} 
\arrow[swap]{r}{\cslr{\alpha_{X_k}}{\psi^{(k)}_{\ind{X}}}} &

\cslr{\alpha_{X_k}}{F\p{k}{\ind{X}}} 
\arrow[swap]{u}{\overline{\alpha}_{(\p{k}{\ind{X}})}}
\end{td}
\qedhere
\end{enumerate} 
\end{mydefn}

\begin{mydefn}
Let $(\alpha, \overline{\alpha}), (\beta, \overline{\beta}) : F \To G$ be 
transformations of pseudofunctors \mbox{$(S, \biclone) \to (S', \biclone')$}. A 
\Def{modification} $\modif : (\alpha, \overline{\alpha}) \to (\beta, 
\overline{\beta})$ consists of a 2-cell $\modif_X: \alpha_X \To \beta_X$ for 
every $X \in S$, such that the following diagram commutes 
for every $t : X_1, \,\dots\, , X_n \to Y$:
\begin{td}[column sep = huge, row sep = 1.5em]
\cslr{\alpha_Y}{Ft} \arrow{r}{\cslr{\modif_Y}{Ft}} 
\arrow[swap]{d}{\overline{\alpha}_t} &
\cslr{\beta_Y}{Ft} \arrow{d}{\overline{\beta}_t} \\
\cslr{G(t)}{\alpha_{X_1} \clonetimes \,\cdots\, \clonetimes \alpha_{X_n}} 
\arrow[swap]{r}[yshift=-2mm]
	{\cslr{G(t)}{\modif_{X_1} \clonetimes \,\cdots\, \clonetimes \modif_{X_n}}} 
	&
\cslr{G(t)}{\beta_{X_1} \clonetimes \,\cdots\, \clonetimes \beta_{X_n}}
\vspace{-2mm}
\end{td} 
\end{mydefn}


It is natural to conjecture that biclones together with their pseudofunctors, 
transformations and modifications form a tricategory $\Biclone$ into which 
$\Bicat$ embeds as a sub-tricategory. We do not pursue such considerations 
here, but we do give the definition of equivalence they would suggest.

\begin{mydefn}
A \Def{biequivalence} between biclones $(S, \biclone)$ and $(S', \biclone')$ consists of
\begin{itemize}
\item Pseudofunctors $F : \biclone \leftrightarrows \biclone' : G$, 
\item Pairs of transformations 
$(\alpha, \overline{\alpha}) 
	: F \circ G \leftrightarrows \id_{\biclone'} : 
	(\alpha', \overline{\alpha'})$ 
and 
\mbox{$(\beta, \overline{\beta}) 
	:  G \circ F \leftrightarrows \id_{\biclone} : 
		(\beta', \overline{\beta'})$},
\item Invertible modifications 
$\modif : \alpha \circ \alpha' \to \id_{\id_{\biclone'}}$, 
$\modif' : \id_{FG} \to \alpha' \circ \alpha$, 
$\altModif : \beta \circ \beta' \to \id_{\id_\biclone}$ and 
$\altModif' : \id_{GF} \to \beta' \circ \beta$.  \qedhere
\end{itemize}
\end{mydefn}

\begin{mylemma} \label{lem:biclone-equivalence-fully-faithful}
For any biequivalence $F : (S, \biclone) \leftrightarrows (S',\biclone') : G$ of biclones, 
\begin{enumerate}
\item The pseudofunctor $F$ is a \Def{local equivalence}, \ie~every $F_{X_1, 
\,\dots\, , X_n; Y} : \biclone(X_1, \,\dots\, , X_n; Y) \to \biclone'(FX_1, 
\,\dots\, , FX_n; FY)$ is full, faithful and essentially surjective, 
\item For every $X' \in S'$ there exists $X \in S$ such that $FX \simeq X'$ in $\biclone'$. 
\end{enumerate}
\begin{proof}
Just as for categories and for bicategories, \cf~\cite[p. 173]{awodey}.
\end{proof}
\end{mylemma}

\section{The type theory \texorpdfstring{$\langBiclone$}{for biclones}}
\label{sec:type-theory-for-biclones}

We now turn to constructing the type theory $\langBiclone$ that will be the 
internal language of biclones. Following the general philosophy of Lambek's 
internal language for multicategories~\cite{Lambek1989}, our 
approach is to 
define a 
term calculus for the rules of 
Construction~\ref{constr:free-biclone}. Thus, for 
every rule in the construction we postulate an introduction 
rule in the type theory. These rules are collected in 
Figures~\ref{r:basic-terms}--\ref{r:basic-rewrites}. Note that we slightly 
abuse notation by simultaneously introducing the structural 
isomorphisms (corresponding to $\assoc{}, \subid{}$ and $\indproj{k}{}$) and 
their inverses.  

The equational theory $\equiv$ is derived directly from the axioms 
of a biclone; the rules are collected together in 
\mbox{Figures~\ref{r:hom-categories-axioms}--\ref{r:congruence-laws}}. The 
typing rules respect this equational theory in the following sense. 

\begin{prooflesslemma}
For any 2-multigraph $\graph$ and derivable judgements $\Gamma \vdash \tau 
\equiv \tau' : \rewrite{t}{t'} : B$ in $\langBiclone(\graph)$, the judgements 
$\Gamma \vdash \tau : \rewrite{t}{t'} : B$ and \mbox{$\Gamma \vdash \tau' : 
\rewrite{t}{t'} : B$} are derivable.
\end{prooflesslemma}

We denote the type theory over a fixed 2-multigraph $\graph$ by 
$\langBiclone(\graph)$; when we do not wish to specify a particular choice of 
signature, we simply write $\langBiclone$.

In what follows we provide a more leisurely introduction to $\langBiclone$ and 
establish some basic meta-theoretic properties.

\paragraph*{Judgements.} We must capture the fact that a biclone 
has both 1-cells and 2-cells: for this we follow the tradition of 2-dimensional 
type theories consisting of types, terms and 
\Def{rewrites} \mbox{(\cf~\cite{Seely1987, Hilken1996, Hirschowitz2013})}. 
Accordingly, there are two 
forms of typing judgement. 
Alongside the usual $\Gamma \vdash t : A$ to indicate `term $t$ has type $A$ in 
context $\Gamma$', we write $\Gamma \vdash \tau : \rewrite{t}{t'} : A$ to 
indicate `$\tau$ is a rewrite from term $t$ of type $A$ to term $t'$ of type 
$A$, in context $\Gamma$'.

Contexts are finite lists of (variable, type) pairs 
in which variable names must not occur more than once: the relevant rules are 
given 
in Figure~\ref{r:biclone:contexts}. Writing $\varEnum$ for the set of 
variables, any context $\Gamma$ determines a finite partial function from 
variables to types; we write $\dom(\Gamma)$ for the domain of this function. 
The concatenation of contexts $\Gamma$ and $\Delta$ satisfying
$\dom(\Gamma) \cap \dom(\Delta) = \emptyset$ is denoted
$\Gamma \concat \Delta$. 

\begin{rules}
\unaryRule	{\faketext}
			{\diamond \mathrm{\: ctx}}
			{\qquad\qquad} 
\binaryRule	{\Gamma \mathrm{\: ctx}}
			{x \notin \dom(\Gamma)}
			{\Gamma, x : A \mathrm{\: ctx}} 
			{$\big(A \in \nodes{\graph}\big)$} \vspace{-\treeskip}
\caption{Context-formation rules for $\langBiclone(\graph)$.\label{r:biclone:contexts}}
\end{rules}

\paragraph*{Raw terms.} Following the template provided by clones, we may 
capture constants in a signature---that is, edges in a 2-multigraph---by 
constants in the type theory, and projections by 
variables. The outstanding question is how to model the substitution operation 
of a biclone. This cannot be 
the standard meta-operation of substitution: 
Construction~\ref{constr:free-biclone} requires that substitution is not 
associative on the nose, only up to the $\assoc{}$ 2-cell.
Our solution is to model the substitution operation of 
the free 
biclone by a form of \emph{explicit substitution}~\cite{Abadi1989}. 
For every family of terms $u_1, 
\dots, u_n$ and term $t$ with free variables among $x_1, \,\dots\, , x_n$ we 
postulate a term $\hcomp{t}{x_1 \mapsto u_1, \,\dots\, , x_n \mapsto u_n}$; 
this is 
the formal analogue of the term $t[u_1/x_1, \,\dots\, , u_n/x_n]$ defined by 
the 
meta-operation of capture-avoiding substitution~(\cf~\cite{Abadi1989, 
Ritter1997}).  
The variables $x_1, \,\dots\, , x_n$ are bound by this 
operation. For a fixed 
2-multigraph 
$\graph$ the \emph{raw terms} are therefore variables, constant terms and 
explicit 
substitutions, as in the grammar
\[
t, u_1, \,\dots\, , u_n ::= x \st c(x_1, \,\dots\, , x_n) \st \hcomp{t}{x_1 
\mapsto u_1, \,\dots\, , x_n \mapsto u_n} 
\qquad\:\: (c \in \graph(A_1, \,\dots\, , A_n; B))
\] 
One may
think of constants $c(x_1, \,\dots\, , x_n)$ as $n$-ary operators: indeed, for 
every sequence of $n$ terms $(u_1, \,\dots\, , u_n)$ explicit substitution defines a 
mapping 
\begin{center}
$(u_1, \,\dots\, , u_n) \mapsto 
	\hcomp{c(x_1, \,\dots\, , x_n)}{x_1 \mapsto u_1,\dots, x_n \mapsto u_n}$
\end{center}
This is emphasised by the following notational 
convention. 

\begin{mynotation} \label{not:BicatConventions}
We adopt the following abuses of notation:
\begin{enumerate}
\item Writing $\hcomp{t}{x_i \mapsto u_i}$ or just $\hcomp{t}{u_i}$ for 
$\hcomp{t}{x_1 \mapsto u_1, \,\dots\, , x_n \mapsto u_n}$, 
\item \label{c:constants} Writing $\hcomp{c}{u_1, \,\dots\, , u_n}$ for the 
explicit substitution $\hcomp{c(x_1, \,\dots\, , x_n)}{x_i \mapsto u_i}$ 
whenever $c$ is a constant. \qedhere
\end{enumerate}
\end{mynotation}

\begin{myremark} 
Alternative notations for explicit substitution include $t\seq{x := u}$ and the 
let-binding operation {\ttfamily let $x = u$ in $t$} (\eg~\cite{Ritter1997, 
pcfbookref}). 
\end{myremark} 

\newpage
\paragraph*{$\alpha$-equivalence on terms.} We work with terms up to 
$\alpha$-equivalence defined in the standard way (\cf~\cite{Ritter1997}).

\begin{mydefn} \label{def:biclone-alpha-eq}
For any 2-multigraph $\graph$ we define the \Def{$\alpha$-equivalence relation} $\aeq$ on raw terms by the rules
\begin{center}
\unaryRule{\phantom{t \aeq t'}}{t \aeq t}{refl} \quad
\unaryRule{t \aeq t'}{t' \aeq t}{symm} \quad
\binaryRule{t \aeq t'}{t' \aeq t''}{t \aeq t''}{trans} \\[-\treeskip]
\trinaryRule{t[y_i /x_i] \aeq t'[y_i / x_i']}{(u_i \aeq u_i')_{i= 1, \,\dots\, 
, n}}{y_1, \,\dots\, , y_n \text{ fresh }}{\hcomp{t}{x_1 \mapsto u_1, \,\dots\, 
, x_n \mapsto u_n} \aeq \hcomp{t}{x_1' \mapsto u_1', \,\dots\, , x_1' \mapsto 
x_n'}}{}
\end{center}
The simultaneous substitution operation $t[u_i / x_i]$ is defined by
\begin{align*}
x_k[u_i / x_i] &:= u_k \\
c(x_1, \,\dots\, , x_n)[u_i / x_i] &:= \hcomp{c}{u_1, \,\dots\, , u_n} \\ 
(\hcomp{t}{z_j \mapsto u_j})[v_i/x_i] &:= \hcomp{t}{z_j \mapsto u_j[v_i/x_i]}  
\end{align*} 
where in the final rule we assume that each $z_j$ does not occur among the 
$x_i$ or freely in any of the $v_i$.
\end{mydefn}

\paragraph*{Raw rewrites.}
Following the pattern set for terms, we define the class of \emph{raw rewrites} 
between terms by the following grammar, where $t, \ind{u}$ and $\ind{v}$ are 
(families of) terms, $x_1, \,\dots\, , x_n$ are variables and $1 \leq i \leq n$:
\[
\tau, \sigma, \sigma_1, \,\dots\, , \sigma_n ::= \assoc{t; \ind{u}; \ind{v}} 
\st \subid{t} \st \indproj{i}{\ind{u}} \st \id_t \st \constrewr(x_1, \,\dots\, 
, x_n) \st \tau \vert \sigma \st \hcomp{\tau}{x_1 \mapsto \sigma_n, \,\dots\, , 
x_n \mapsto \sigma_n}
\] 
with a family of inverses (for $i=1, \dots, n$), as follows:
\[
\assoc{t; \ind{u}; \ind{v}}^{-1} \st \subid{t}^{-1} \st \indproj{-i}{\ind{u}}
\]
Taking the rewrites in turn, we have invertible \Def{structural rewrites} 
$\assoc{}, \subid{}$ and $\indproj{i}{} \: (i = 1, \,\dots\, , n)$ and an 
identity rewrite $\id_t$ for every term $t$. Next, for every constant 
$\constrewr \in \graph(A_1,\dots, A_n; B)$
we have a
constant rewrite
$\constrewr(x_1, \,\dots\, , x_n)$. 
Vertical composition is captured by a binary operation on 
rewrites~(\cf~\cite{Hilken1996, Hirschowitz2013, Licata2017}), while the 
explicit substitution 
operation mirrors that for terms. (Note that vertical composition follows 
function composition order, not diagrammatic order.) We adopt the standard 
category-theoretic convention of 
writing $t$ for $\id_t$ where no ambiguity may arise, as well as adapting the 
conventions of Notation~\ref{not:BicatConventions} to rewrites. In particular, 
one obtains \Def{whiskering} operations $\hcomp{t}{\sigma}$ and 
$\hcomp{\tau}{u}$ for terms $t, u$ and rewrites 
$\tau : t \To t', \sigma : u \To u'$. 

\paragraph*{$\alpha$-equivalence on rewrites.}
The $\alpha$-equivalence relation extends to rewrites in the way one would 
expect: as for terms, the substitution operation 
binds the variables being explicitly substituted for. The 
definition of the meta-operation of substitution on rewrites is analogous to 
that employed by Hilken~\cite{Hilken1996} and 
Hirschowitz~\cite{Hirschowitz2013}.  

\begin{mydefn}
For any 2-multigraph $\graph$ we define the \Def{$\alpha$-equivalence relation} $\aeq$ on rewrites by the rules
\begin{center}
\unaryRule{\phantom{\tau \aeq \tau'}}{\tau \aeq \tau}{refl}
\quad
\unaryRule{\tau \aeq \tau'}{\tau' \aeq \tau}{symm}
\vspace{-\treeskip}
\quad
\binaryRule{\tau \aeq \tau'}{\tau' \aeq \tau''}{\tau \aeq \tau''}{trans} 
\vspace{-\treeskip}
\\
\unaryRule{t \aeq t'}{\subid{t} \aeq \subid{t'}}{}
\qquad
\unaryRule{u_1 \aeq u_1' \quad \dots \quad u_n \aeq u_n'}{\indproj{k}{u_1, 
\dots, u_n} \aeq \indproj{k}{u_1', \,\dots\, , u_n'}}{$1 \leq k \leq n$}
\vspace{-\treeskip} \\

\trinaryRule
	{(u_j \aeq u_j')_{j = 1, \,\dots\, , m}}
	{(v_i \aeq v_i')_{i = 1, \,\dots\, , n}}
	{t \aeq t'}
	{\assoc{t, \ind{v}, \ind{u}} \aeq \assoc{t', \ind{v}', \ind{u}'}}
	{} 
	
\binaryRule
	{\tau \aeq \tau'}
	{\sigma \aeq \sigma'}
	{\tau \vert \sigma \aeq \tau' \vert \sigma'}{}

\trinaryRule
	{\tau[y_i /x_i] \aeq \tau'[y_i / x_i']}
	{(\sigma_i \aeq \sigma_i')_{i= 1, \,\dots\, , n}}
		{y_1, \,\dots\, , y_n \text{ fresh }}
	{\hcomp{\tau}{x_1 \mapsto \sigma_1, \,\dots\, , x_n \mapsto \sigma_n} 
		\aeq 
	\hcomp{\tau}{x_1' \mapsto \sigma_1', \,\dots\, , x_1' \mapsto \sigma_n'}}
	{} \vspace{-\treeskip}
\end{center}
The meta-operation of capture-avoiding substitution is extended to rewrites as follows:
\begin{align*}
\subid{u}[u_i/x_i] &:= \subid{u[u_i/x_i]} \\
\indproj{k}{t_1, \,\dots\, , t_n}[u_i/x_i] &:= \indproj{k}{\ind{t}[u_i/x_i]} \\
\assoc{t,\ind{u},\ind{v}}[u_i/x_i] &:= \assoc{t[u_i/x_i], \ind{u}[u_i/x_i], \ind{v}[u_i/x_i]} \\
\constrewr(x_1, \,\dots\, , x_n)[u_i / x_i] &:= \hcomp{\constrewr}{u_1, 
\,\dots\, , u_n} \\
(\tau' \vert \tau)[u_i/x_i] &:= \tau'[u_i/x_i] \vert \tau[u_i/x_i] \\
\id_t[u_i/x_i] &:= \id_{t[u_i/x_i]} \\
(\hcomp{\tau}{z_j \mapsto \sigma_j})[u_i/x_i] &:= 	\hcomp{\tau}{z_j \mapsto \sigma_j[u_i/x_i]}
\end{align*}
where in the final rule we assume that each $z_j$ does not occur among the 
$x_i$ or freely in any 
of the $u_i$. These rules extend to the inverses of rewrites in the obvious 
fashion.
\end{mydefn}

A structural induction shows the typing judgement respects $\alpha$-equivalence.

\begin{prooflesslemma} 
Let $\graph$ be a 2-multigraph. Then in $\langBiclone(\graph)$:
\begin{enumerate} 
\item If $\Gamma \vdash t : B$ and $t =_\alpha t'$ then $\Gamma \vdash t' :B$,
\item If $\Gamma \vdash \tau : \rewrite{t}{t'} : B$ and $\tau =_\alpha \tau'$ then $\Gamma \vdash \tau : \rewrite{t}{t'} : B$. \qedhere
\end{enumerate} 
\end{prooflesslemma} 

In an explicit substitution calculus the structural operations manifest 
themselves in a correspondingly explicit manner. Indeed, the fact that 
$\langBiclone$ admits arbitrary context renamings follows immediately from the 
\rulename{horiz-comp} rule.

\begin{mydefn}
Let $\Gamma := (x_i : A_i)_{i=1,\dots,n}$ and $\Delta := (y_j : 
B_j)_{j=1,\dots,m}$ be contexts. A \Def{context renaming} \mbox{$r : \Gamma \to 
\Delta$} is a mapping $r : \{x_1, \,\dots\, , x_n\} \to \{y_1, \,\dots\, , 
y_m\}$ on variables which respects typing in the sense that whenever $r(x_i) = 
y_j$ then $A_i = B_j$.
\end{mydefn}

The following rules are then derivable for any context renaming $r$.

\begin{rules} 
\begin{prooftree}
\AxiomC{$\Gamma \vdash t : A$}
\AxiomC{$r : \Gamma \to \Delta$}
\BinaryInfC{$\Delta \vdash \hcomp{t}{x_1 \mapsto r(x_1), \dots, x_n \mapsto r(x_n)} : A$}
\end{prooftree}
\begin{prooftree}
\AxiomC{$\Gamma \vdash \tau: \rewrite{t}{t'} : A$}
\AxiomC{$r : \Gamma \to \Delta$}
\BinaryInfC{$\Delta \vdash \hcomp{\tau}{x_i \mapsto r(x_i)} : \rewrite{\hcomp{t}{x_i \mapsto r(x_i)}}{\hcomp{t'}{x_i \mapsto r(x_i)}} : A$}
\end{prooftree}
\caption{Context renaming as a derived rule (for $\Gamma = (x_i : A_i)_{i=1, \dots, n}$)\label{fig:context-renaming}}
\end{rules}

Weakening arises as a special case: for a fresh variable 
$x \notin \dom(\Gamma)$, one takes the 
inclusion 
\mbox{$\mathrm{inc}_{x} : \Gamma \hookrightarrow \Gamma, x : A$}. 

\begin{mynotation} \label{not:context-morphism}
For a context renaming $r$ we write $\hcomp{t}{r}$ and $\hcomp{\tau}{r}$ for the terms and rewrites formed using the admissible rules of Figure~\ref{fig:context-renaming}.
\end{mynotation}

\begin{figure*}[!ht]
\centering
{\small
\begin{minipage}{\textwidth}
\begin{mdframed}
\centering
\input{rules/basic-terms}
\caption{Introduction rules on basic terms \label{r:basic-terms}}
\end{mdframed}
\end{minipage}

\begin{minipage}{\textwidth}
\begin{mdframed}
\centering
\input{rules/structural-rewrites}
\caption{Introduction rules on structural rewrites \label{r:structural-rewrites}}
\end{mdframed}
\end{minipage}

\begin{minipage}{\textwidth}
\begin{mdframed}
\centering
\input{rules/basic-rewrites}
\input{rules/rewrites-constructors} 
\caption{Introduction rules on basic rewrites \label{r:basic-rewrites}}
\end{mdframed}
\end{minipage}
}

\begin{manyfigcap}
Introduction rules for terms, structural rewrites and basic rewrites in $\langBiclone(\graph)$. 
\end{manyfigcap}

\label{fig:common-rules}
\end{figure*}

%
%

\begin{figure*}[!h]

\centering
{\small

\begin{minipage}{\textwidth}
\begin{mdframed}
\centering
\input{rules/categorical-vert-comp}
\caption{Categorical structure of vertical composition \label{r:hom-categories-axioms}}
\end{mdframed}
\end{minipage}

\begin{minipage}{\textwidth}
\begin{mdframed}
\centering
\input{rules/interchange}
\caption{Preservation rules \label{r:hcomp-functor}}
\end{mdframed}
\end{minipage}

\begin{minipage}{\textwidth}
\begin{mdframed}
\centering
\input{rules/structural-rewrites-nat}
\caption{Naturality rules on structural rewrites \label{r:structural-rewrites-nat}}
\end{mdframed}
\end{minipage}

\begin{minipage}{\textwidth}
\begin{mdframed}
\centering
\input{rules/biclone-laws}
\caption{Biclone laws \label{r:biclone-laws}}
\end{mdframed}
\end{minipage}
}

\begin{manyfigcap}
Equational theory for structural rewrites in $\langBiclone(\graph)$.
\end{manyfigcap}
\end{figure*}


\begin{figure*}[!h]

\centering
{\small

\begin{minipage}{\textwidth}
\begin{mdframed}
\centering
\input{rules/structural-rewrites-invertibility}
\caption{Invertibility of the structural rewrites}
\end{mdframed}
\end{minipage}

\begin{minipage}{\textwidth}
\begin{mdframed}
\centering
\input{rules/congruences}
\caption{Congruence laws \label{r:congruence-laws}}
\end{mdframed}
\end{minipage}
}

\begin{manyfigcap}
Equational theory for structural rewrites in $\langBiclone(\graph)$.
\end{manyfigcap}
\end{figure*}


\clearpage

\paragraph{Well-formedness properties of $\langBiclone$.} \label{sec:PropertiesOfLangBicat}

We finish this introduction to $\langBiclone$ by showing that it satisfies 
versions of the standard syntactic properties of, for example, the simply-typed 
lambda 
calculus~(\cf~\cite[Chapter 4]{Crole1994}). The intention is to justify the 
claim that the properties one would expect by analogy with the simply-typed 
lambda calculus do in fact 
hold. The proofs are all straightforward structural inductions. 

\begin{mydefn}  \label{def:biclone-free-vars}
Fix a 2-multigraph $\graph$. We define the \Def{free variables in a term} $t$ in $\langBiclone(\graph)$ as follows: 
\begin{align*}
\fv(x_i) &:= \{ x_i \} \qquad  \text{for } x_i \text{ a variable,} \\
\fv\big(c(x_1, \,\dots\, , x_n)\big) &:= \{x_1, \,\dots\, , x_n\} \qquad \text{ 
for } c \in \graph(A_1, \,\dots\, , A_n; B), \\
\fv(\hcomp{t}{x_1 \mapsto u_1, \,\dots\, , x_n \mapsto u_n}) &:= (\fv(t) - \{ 
x_1, \,\dots\, , x_n\}) \cup \textstyle{\bigcup}_{i=1}^n \fv(u_i)
\end{align*}
Similarly, define the \Def{free variables in a rewrite $\tau$} in $\langBiclone(\graph)$ as follows:
\begin{align*}
\fv\big(\subid{t}\big) &:= \fv(t) \\
\fv\big( \indproj{k}{u_1, \,\dots\, , u_n}\big) &:= \fv(u_k) \\
\fv\big(\assoc{t,\ind{v}, \ind{u}}\big) &:= \textstyle{\bigcup}_{i=1}^n \fv(u_i) \\
\fv(\id_t) &:= \fv(t)  \\
\fv(\tau' \vert \tau) &:= \fv(\tau') \cup \fv(\tau)  \\
\fv\big(\sigma(x_1, \,\dots\, , x_n)) &:= \{ x_1, \,\dots\, , x_n\} \text{ for 
} \sigma \in \graph(A_1, \,\dots\, , A_n; B)(c, c')  \\
\fv(\hcomp{\tau}{x_1 \mapsto \sigma_1, \,\dots\, , x_n \mapsto \sigma_n}) &:= 
(\fv(\tau) - \{ x_1, \,\dots\, , x_n\}) \cup \textstyle{\bigcup}_{i=1}^n 
\fv(\sigma_i) 
\end{align*}
We define the free variables of a specified inverse $\sigma^{-1}$ to be exactly the free variables of $\sigma$. An occurrence of a variable in a term (rewrite) is \Def{bound} if it is not free. 
\end{mydefn} 

\begin{prooflesslemma}
Let $\graph$ be a 2-multigraph. For any derivable judgements $\Gamma \vdash u : B$ and \mbox{$\Gamma \vdash \tau : \rewrite{t}{t'} : B$} in $\langBiclone(\graph)$, 
\begin{enumerate} 
\item $\fv(u) \subseteq \dom(\Gamma)$, 
\item $\fv(\tau) \subseteq \dom(\Gamma)$, 
\item The judgements $\Gamma \vdash t : B$ and $\Gamma \vdash t' : B$ are both derivable. 
\end{enumerate} 
Moreover, for any context $\Gamma := (x_i  : A_i)_{i = 1, \,\dots\, , n}$  and 
derivable terms $(\Delta \vdash u_i : A_i)_{i= 1, \,\dots\, , n}$, 
\begin{enumerate}
\item If $\Gamma \vdash t : B$, then $\Delta \vdash t[u_i / x_i] : B$, 
\item If $\Gamma \vdash \tau : \rewrite{t}{t'} : B$, then $\Delta \vdash \tau[u_i / x_i] : \rewrite{t[u_i / x_i]}{t[u_i / x_i]} : B$. \qedhere
\end{enumerate}
\end{prooflesslemma}

%

\subsection{The syntactic model} \label{sec:bicat:syntactic-model}

The rules of $\langBiclone$ are synthesised directly from the construction of 
the free biclone on a \mbox{2-multigraph}. It is not surprising, therefore, that its 
syntactic model satisfies the same free property, justifying our description of 
$\langBiclone$ as a type theory for biclones. In this section we spell 
out the construction and show that it restricts to bicategories.

Constructing the syntactic model is a matter of reversing the correspondence between the rules of $\langBiclone$ and Construction~\ref{constr:free-biclone}.

\begin{myconstr} \label{constr:langbiclone-syntactic-model}
For any 2-multigraph $\graph$ define the \Def{syntactic model} 
$\synclone{\graph}$ of $\langBiclone(\graph)$ as follows. The sorts are nodes 
$A, B, \dots$ of $\graph$. For $A_1, \,\dots\, , A_n, B \in \nodes{\graph}$ 
the hom-category $\synclone{\graph}(A_1, \,\dots\, , A_n; B)$ has objects 
$\alpha$-equivalence classes of terms \mbox{$(x_1 : A_1, \,\dots\, , x_n : A_n 
\vdash t : B)$} derivable in $\langBiclone(\graph)$. We assume a fixed 
enumeration $x_1, x_2, \dots$ of variables, and that the variable name in the 
$i$th position is determined by this enumeration. Morphisms in 
$\synclone{\graph}(A_1, \,\dots\, , A_n; B)$ are $\alpha{\equiv}$-equivalence 
classes of rewrites \mbox{$(x_1 : A_1, \,\dots\, , x_n : A_n \vdash \tau : 
\rewrite{t}{t'} : B)$}. Composition is vertical composition and the identity is 
$\id_t$. 

The substitution operation $\big( t, (u_1, \,\dots\, , u_n) \big) \mapsto 
\cslr{t}{u_1, \,\dots\, , u_n}$ is explicit substitution
\begin{align*} 
t, (u_1, \,\dots\, , u_m) \mapsto \hcomp{t}{x_1 \mapsto u_1, \,\dots\, , x_n 
\mapsto u_n} \\ 
\tau, (\sigma_1, \,\dots\, , \sigma_m) \mapsto \hcomp{\tau}{x_1 \mapsto 
\sigma_1, \,\dots\, , x_n \mapsto \sigma_n}
\end{align*} 
and the projections $(A_1, \,\dots\, , A_n) \to A_k$ are instances of the 
\rulename{var} rule \mbox{$x_1 : A_1, \,\dots\, , x_n : A_n \vdash x_k : A_k$} 
for $k = 1, \,\dots\, , n$. The 2-cells $\assoc{}, \subid{}$ and 
$\indproj{k}{}$ are the corresponding structural rewrites.
\end{myconstr}

\begin{mynotation} \label{not:multiple-variable-names}
We shall generally play fast and loose with the 
requirement that the variables 
in a context $(x_1 : A_1, \,\dots\, , x_n : A_n)$ are labelled in turn by the 
enumeration $x_1, \,\dots\, , x_n, \dots$. We will allow ourselves to pick 
more meaningful variable names as a simple form of syntactic sugar, and rely on 
the fact that the proper variable names can always be recovered when required.
\end{mynotation}

The equational theory guarantees that $\synclone{\graph}$ is a biclone. The 
proof of the free property mirrors Lemma~\ref{lem:free-biclone}.

\begin{mylemma} \label{lem:strict-pseudofunctor-biclone}
For any 2-multigraph $\graph$, biclone $(S, \biclone)$ and 2-multigraph 
homomorphism $h : \graph \to \biclone$ there exists a unique strict 
pseudofunctor $h\sem{-} : \synclone{\graph} \to \biclone$ such that  $h\sem{-} 
\circ \inc = h$, for $\inc : \graph \hookrightarrow \synclone{\graph}$ the 
inclusion.
\begin{proof}
Fix a context $\Gamma := (x_i : A_i)_{i = 1,\dots, n}$. We define $h\sem{-}$ by 
induction on the derivation of judgements in 
$\langBiclone$:
\newpage
\begin{align*} 
h\sem{B} &:= h(B) \qquad \text{on types} \\[5pt]
h\sem{\Gamma \vdash c(x_1, \,\dots\, , x_n) : B} &:= h(c) \qquad \text{ for } c 
\in \graph(\ind{A};B) \\ 
h\sem{\Delta \vdash \hcomp{t}{x_i \mapsto u_i} : B} &:= 
	\cs
		{\big(h\sem{\Gamma \vdash t : B}\big)}
		{h\sem{\Delta \vdash \ind{u} :\ind{A}}} \\[5pt]
h\sem{\Gamma \vdash \id_t : \rewrite{t}{t} : B} &:= \id_{h\sem{\Gamma \vdash t : B}} \\
h\sem{\Gamma \vdash \constrewr(\ind{x}) : \rewrite{c(\ind{x})}{c'(\ind{x})} : B} &:= h(\constrewr) \qquad \text{ for } \constrewr \in \graph(\ind{A}, B)(c,c')  \\
h\sem{\Gamma \vdash \tau' \vert \tau : \rewrite{t}{t''} : B} &:= h\sem{\Gamma \vdash \tau' : \rewrite{t'}{t''} : B} \vert h\sem{\Gamma \vdash \tau : \rewrite{t}{t'} : B} \\
h\sem{\hcomp{\tau}{x_i \mapsto \sigma_i}} &:= 
	\cs
		{\big(h\sem{\Gamma \vdash \tau : \rewrite{t}{t'} : B}\big)}
		{h\sem{\Delta \vdash \ind{\sigma} : 
			\rewrite{\ind{u}}{\ind{u}'} : \ind{A}}}
\end{align*} 
where we omit the full typing derivation 
$\Delta \vdash \hcomp{\tau}{x_i \mapsto \sigma_i} : 
	\rewrite{\hcomp{t}{x_i \mapsto u_i}}{\hcomp{t'}{x_i \mapsto u_i'}} : B$ in 
the final case for reasons 
of space. In order for $h\sem{-}$ to be strict we must require that it strictly 
preserves the $\assoc{}, \subid{}$ and $\indproj{k}{}$ 2-cells. Uniqueness 
holds just as in Lemma~\ref{lem:free-biclone}.
\end{proof}
\end{mylemma}

\begin{prooflessthm} \label{thm:bic:syntactic-model-is-free}
For any 2-multigraph $\graph$, the syntactic model $\synclone{\graph}$ of 
$\langBiclone(\graph)$ is the free biclone on $\graph$.
\end{prooflessthm}

A type theory satisfying a property of this form, and which is therefore sound and complete for reasoning in the freely constructed structure, is often referred to as the \Def{internal language} or \Def{internal logic} (\eg~\cite{Makkai1977, LambekAndScott, Crole1994, Gambino2013}). This terminology is used with varying degrees of precision, and generally not in the precise sense of Lambek%
~\cite[Definition~5.3]{Lambek1989}; nonetheless, we may now justifiably state that $\langBiclone$ is the internal language of biclones. 


%

By the theorem, we may identify $\synclone{\graph}$ with the free biclone
$\freeBiclone{\graph}$ on $\graph$. 
The diagram of adjunctions~(\ref{eq:biclone-bicat-multigraph-diagram})~(p.~\pageref{eq:biclone-bicat-multigraph-diagram}) then entails 
that for a 2-graph $\graph$ the nucleus of $\synclone{\graph}$---%
obtained by restricting the syntactic model of $\langBiclone$ to unary multimaps---%
is the free bicategory on $\graph$.
Equivalently, one may restrict the type theory $\langBiclone$ to unary contexts
and construct its syntactic model as in Construction~\ref{constr:langbiclone-syntactic-model}.
Let $\langBicat$ denote the type theory obtained by replacing the 
context-formation rules of Figure~\ref{r:biclone:contexts} with the single rule
of Figure~\ref{r:biclone:unary-contexts}.
\begin{rules}
\unaryRule	{\phantom{\Gamma}}
			{x : A \mathrm{\: ctx}} 
			{$\big(A \in \nodes{\graph}\big)$} \vspace{-\treeskip}
\vspace{4mm}
\caption{Context-formation rule for $\langBicat(\graph)$.\label{r:biclone:unary-contexts}}
\end{rules}

\newpage
\begin{myconstr} \label{constr:bicat-termcat} 
For any 2-graph $\graph$, define a bicategory $\urestrict{\synclone{\graph}}$ 
as follows. 
Objects are unary contexts $(x : A)$ for $x$ a \emph{fixed} variable name. The 
hom-category \mbox{$\urestrict{\synclone{\graph}}\big((x:A),(x:B)\big)$} has 
objects 
$\alpha$-equivalence classes of derivable terms \mbox{$(x : A \vdash t :B)$} in 
$\langBicat$
and morphisms $\alpha{\equiv}$-equivalence classes of rewrites 
\mbox{$(x : A \vdash \tau : \rewrite{t}{t'} : B)$} in $\langBicat$. Vertical 
composition is the $\vert$ operation. 
Horizontal 
composition is given by explicit substitution and the identity on $(x : A)$ by 
the \rulename{var} rule $( x: A \vdash x : A)$. The structural isomorphisms 
$\l, \r$ and $\a$ are $\proj{}$, $\subid{}^{-1}$ and $\assoc{}$, respectively. 
\end{myconstr} 


\begin{myremark}
The structural isomorphism $\r$ is given by $\subid{}^{-1}$ because we 
have 
directed the structural isomorphisms in a 
biclone to match that of a skew monoidal category, but followed B{\'e}nabou's 
convention~\cite{Benabou1967} directing the unitors in a bicategory to remove 
compositions with the identity.
\end{myremark}

The required theorem follows immediately from Theorem~\ref{thm:bic:syntactic-model-is-free}
and the chain of isomorphisms~(\ref{eq:nucleus-gives-free-bicat})~(p.~\pageref{eq:nucleus-gives-free-bicat}).


\begin{prooflessthm} \label{thm:unary-syntactic-model-of-lang-biclone-free-bicat}
For any 2-graph $\graph$, the syntactic model $\urestrict{\synclone{\graph}}$ of 
$\langBicat(\graph)$ is the free bicategory on $\graph$.
\end{prooflessthm}

The restriction to a fixed variable name is necessary for the free property to 
be strict. Without such a restriction there are countably many equivalent 
objects $(x_1 : A), (x_2 : A), \dots$ in $\urestrict{\synclone{\graph}}$, and 
the action of 
the pseudofunctor defined in Lemma~\ref{lem:strict-pseudofunctor-biclone} is 
unique only up to its action on each variable name. The next lemma shows 
that---up to biequivalence---this restriction is immaterial. 

\begin{mylemma} \label{lem:subbicategory-equivalence-biequivalence}
Let $\baseCat$ be a bicategory and $\subCat$ a sub-bicategory. Suppose that for 
every $X \in \baseCat$ there exists a chosen $\subbic{X} \in \subCat$ with a 
specified adjoint equivalence 
$f_X : X \leftrightarrows \subbic{X} : g_X$ in $\baseCat$ such that
\begin{enumerate}
\item For $X \in \subCat$ the equivalence $X \simeq \subbic{X}$ is the identity, and
\item If $h : X \to Y$ is a 1-cell in $\subCat$, then so is the composite
$(g_Y \circ h) \circ f_X :\subbic{X} \to \subbic{Y}$.
\end{enumerate}
Then $\baseCat$ and $\subCat$ are biequivalent.
\begin{proof}
Let us denote the 2-cells witnessing the equivalence $X \simeq \subbic{X}$
by
\begin{align*}
\un_X &: \Id_{\subbic{X}} \To g_X \circ f_X \\
\co_X &: f_X \circ g_X \To \Id_X
\end{align*}
There exists an evident pseudofunctor $\inc: \subCat \hookrightarrow \baseCat$ 
given by the inclusion. In the other direction, we define $\subFun : \baseCat 
\to \subCat$ by setting 
\begin{center}
$\subFun(X) :=  \subbic{X}$ \quad and \quad $\subFun(\tau : t \To t' : X \to Y) 
:= 
(g_Y \circ \tau) \circ f_X$
\end{center}
We then define $\psi_X := \Id_{\subbic{X}} \XRA{\un_X} g_X \circ f_X 
\XRA{\iso} (g_X \circ \Id_X) \circ h_X = \subFun(\Id_X)$. For a composable pair 
$X \xra{u} Y \xra{t} Z$ we define $\phi_{t,u}$ by commutativity of the 
following diagram:
\begin{td}[column sep = 7em]
\left(g_Z \circ (t \circ f_Y) \right)) 
	\circ \left(g_Y \circ (u \circ f_X)\right)
\arrow[swap]{d}{\iso}
\arrow{r}{\phi_{t, u}} &
g_Z \circ \left( (t \circ u) \circ f_X \right) \\

\left(g_Z \circ t\right) \circ \left((f_Y \circ g_Y) \circ (u \circ f_X) \right)
\arrow[swap]{r}[yshift = -1mm]
	{(g_Z \circ t) \circ (\co_Y \circ (u \circ f_X))} &
\left(g_Z \circ t\right) \circ \left(\Id_Y \circ (u \circ f_X) \right)
\arrow[swap]{u}{\iso}
\end{td}
The unit and associativity laws for a pseudofunctor follow from coherence and 
the triangle laws of an adjoint equivalence. We then need to construct 
pseudonatural transformations $(\alpha, \overline{\alpha}) : \id_{\baseCat} 
\leftrightarrows \inc \circ \subFun : (\beta, \overline{\beta})$ and $(\gamma, 
\overline{\gamma}) : \id_{\subCat} \leftrightarrows \subFun \circ \inc : 
(\delta, \overline{\delta})$. 

For $\alpha$, we take $\alpha_X := g_X$ and $\overline{\alpha}_t$ to be the 
composite
\begin{td}[column sep = 4em]
g_Y \circ t 
\arrow{r}{\overline{\alpha}_t}
\arrow[swap]{d}{\iso} &
\left( g_Y \circ \left(t \circ f_X\right) \right) \circ g_X \\
\left( g_Y \circ t \right) \circ \Id_X 
\arrow[swap]{r}{g_Y \circ t \circ \co_X^{-1}} &
\left( g_Y \circ t \right) \circ \left(f_X \circ g_X \right) 
\arrow[swap]{u}{\iso}
\end{td}
for $t : X \to Y$. For $\beta$ and $\overline{\beta}$ the idea is the same. We 
define $\beta_X := f_X$ and for $t : X \to Y$ we set
\begin{td}[column sep = 4em]
f_Y \circ \left( g_Y \circ \left( t \circ f_Y \right) \right) 
\arrow[swap]{d}{\iso}
\arrow{r}{\overline{\beta}_t} &
t \circ f_X \\
\left( f_Y \circ g_Y \right) \circ \left( t \circ f_X \right)
\arrow[swap]{r}{\co_Y \circ t \circ f_X} &
\Id_Y \circ \left( t \circ f_X \right)
\arrow[swap]{u}{\iso}
\end{td}
The definitions of $(\gamma, \overline{\gamma})$ and 
$(\delta, \overline{\delta})$ are identical. One then obtains modifications 
$\modif : \id \xra{\iso} \alpha \circ \beta$ and 
$\altModif : \beta \circ \alpha \xra{\iso} \id$ by taking 
$\modif_X := \Id_X \XRA{\un_X} g_X \circ f_X$ and 
$\altModif_X := f_X \circ g_X \XRA{\co_X} X$; similarly 
$\gamma \circ \delta \iso \id$ 
and $\delta \circ \gamma \iso \id$. 
\end{proof}
\end{mylemma}

Hence, $\langBicat$ is the internal language for bicategories. If one restricts 
to a single variable name the universal property is strict, else 
it is up to biequivalence. In the next section we 
show that the syntactic model of $\langBiclone$ is biequivalent as a biclone to 
the syntactic model of a strict type theory. From this we deduce a 
coherence result for biclones, which amounts to a form of normalisation for the 
rewrites of $\langBiclone$. All of this will restrict to unary contexts, and 
hence to $\langBicat$, recovering a version of the coherence theorem of Mac 
Lane \& Par{\'e}~\cite{MacLane1985}.

\section{Coherence for biclones} \label{sec:coherence-for-biclones} 

In practice, the coherence theorem for bicategories~\cite{MacLane1985} entails 
that one may treat any bicategory as though it were a 2-category: roughly, one 
may assume that the structural isomorphisms $\a, \l$ and $\r$ behave as though 
they were the identity (see~\eg~\cite[Chapter 1]{Leinster2004} for a detailed 
exposition). In terms of $\langBicat$, this amounts to treating  $\assoc{}, 
\indproj{i}{}$ and $\subid{}$ as though they were all identities. Our aim in 
this section is to extend this result to $\langBiclone$.

The motivation is three-fold. 
First, the coherence theorem will simplify 
the calculations we shall require in 
future chapters. Second, the proof involves some of the calculations we shall 
need to extend when it comes to defining a pseudofunctorial interpretation of 
the full type theory $\langCartClosed$ (see Section~\ref{sec:syntactic-model-free-property}). 
Finally, 
the proof strategy is of 
interest in itself.  The strategy may be 
regarded as a version of Mac Lane's classical strategy for monoidal 
categories~\cite[Chapter VII]{cfwm}, in which the syntax of the respective 
type theories provide \Def{structural induction} principles. 
It is 
reasonable to imagine that one may prove similar results for monoidal 
bicategories (via a linear calculus), tricategories (via a 3-dimensional 
calculus) or even higher-dimensional structures, by an analogous strategy. 

To foreshadow the coherence result we shall prove in later chapters, let us 
make precise the notion of normalisation we are interested in. We wish to lift 
the standard notion of normalisation for systems such as the (untyped) 
$\lambda$-calculus (\eg~\cite{Girard1989}) to a normalisation property on 
rewrites. More precisely, we wish to consider versions of 
\emph{abstract reduction systems}~\cite{Huet1980} in which one also tracks how 
a reduction might happen; that is, the possible	 
\emph{witnesses} of a reduction. 
Our notion of normalisation then becomes: there is 
at most one witness to any possible reduction. 
This suggests the following definitions. We use the term 
\Def{constructive} by analogy with constructive proofs, in which one requires 
an explicit witness to the truth of a statement, to emphasise 
that we are requiring an explicit witnesses to the existence of a 
reduction. 

\begin{mydefn} \quad
\begin{enumerate}
\item An \Def{abstract reduction system (ARS)} $(A, {\to})$ is a set $A$ 
equipped with a 
binary \Def{reduction relation} ${\to} \subseteq A \times A$.
\item A \Def{constructive abstract reduction system (CARS)} consists of a set 
$A$ together with a family of sets $W_A(a, b)$ of \Def{reduction witnesses} 
indexed by $a, b \in A$. A CARS is \Def{coherent} if for every
$a, b \in A$ and $u, v \in  W_A(a, b)$, 
one has $u = v$.  \qedhere
\end{enumerate}
\end{mydefn}

In a CARS we are not merely interested in the existence of a reduction: we are 
also interested in the equality relation on reductions. In particular, an ARS 
in the usual sense is a 
CARS in which every $W(a,a')$ is either empty or a singleton: either $a$ 
reduces to $a$, or it does not. 

The term `coherent' is motivated by the following example. 

\begin{myexmp} \quad
\begin{enumerate}
\item Every graph 
$\graph$ defines a CARS $A(\graph)$ with underlying set $\nodes{\graph}$ and 
reduction witnesses $W_{A(\graph)}(t, t') := \graph(t, t')$. 
\item 
Every category $\catC$ defines a CARS 
$\overline{\catC}$ on $ob(\catC)$ by taking 
$W_{\overline{\catC}}(A, B) := \catC(A, B)$. 
The coherence theorem for monoidal categories of~\cite[Chapter VII]{cfwm} 
then states that the CARS corresponding to the free monoidal category on one 
generator is coherent. \qedhere
\end{enumerate}
\end{myexmp}

In the bicategorical setting, we are interested in coherence in each hom-category.

\begin{mydefn} \quad
\begin{enumerate}
\item A 2-multigraph $\graph$ is \Def{locally coherent} if for every 
$A_1, \dots, A_n, B \in \nodes{\graph}$ the associated CARS 
$A\big( \graph(A_1, \,\dots\, ,A_n; B) \big)$ is coherent.

\item A biclone (bicategory) is \Def{locally coherent} if its underlying 2-multigraph is locally coherent. \qedhere
\end{enumerate}
\end{mydefn}

Spelling out the definitions, a 2-multigraph $\graph$ is locally coherent if 
for all 
edges $e$ and $e'$ in $\graph(A_1, \,\dots\, , A_n; B)$ there exists at most one 
surface $\constrewr : e \To e'$, and a biclone is 
locally coherent if there is at most one 2-cell between any parallel pair of 
terms. The coherence theorem for bicategories~\cite{MacLane1985} can therefore 
be rephrased as stating that the free bicategory on a 2-multigraph is locally 
coherent.

Now, every type theory consisting of types, terms and rewrites has an 
underlying 2-multigraph with nodes given by the types, edges 
$A_1, \,\dots\, , A_n \to B$ by the $\alpha$-equivalence classes of derivable 
terms 
\mbox{$x_1 : A_1, \,\dots\, , x_n : A_n \vdash t :B$} and surfaces by the
derivable rewrites 
modulo $\alpha$-equivalence and the equational theory. We call the type theory 
\Def{locally coherent} if this 2-multigraph is locally coherent. We spend the 
rest of this chapter proving that $\langBiclone$ is locally coherent.

Our strategy is the following. We shall adapt the calculi of 
Hilken~\cite{Hilken1996} and
Hirschowitz~\cite{Hirschowitz2013} to construct a type theory that matches 
$\langBiclone$ but has a strict substitution operation; the syntactic model 
will be the free 2-clone (\cf~Construction~\ref{constr:free-biclone}). We shall
then construct an equivalence between the two 
syntactic models by induction on the respective type theories. We finish by 
briefly commenting how the result restricts to bicategories.

\subsection{A strict type theory}

The first step is the construction of a strict type theory. Since we draw 
heavily on previous work, our presentation will be brief. 
Fix some 2-multigraph $\graph$. The type theory $\hir(\graph)$ (where $H$ stands for both \emph{Hilken} and \emph{Hirschowitz}) is constructed as follows. Contexts are as in $\langBiclone$. The \emph{raw terms} are either variables or constants, given by the following grammar: 
\[
u_1, \,\dots\, , u_n ::= x \st c(u_1, \,\dots\, , u_n)
\]
As for $\langBiclone$, we think of constants $c(x_1, \,\dots\, , x_n)$ as 
$n$-ary \emph{operators}. The \emph{raw rewrites} are vertical composites of 
identity maps and constant rewrites: 
\[
\sigma_1, \,\dots\, , \sigma_n, \tau, \sigma ::= \id_t \st \constrewr(u_1, 
\,\dots\, , u_n) \st c(\sigma_1, \,\dots\, , \sigma_n) \st \tau \vert \sigma
\qquad (u_1, \dots, u_n \text{ terms})
\]
Note that we require two forms of constant rewrite, corresponding to 
substitution of terms into rewrites and substitution of rewrites into terms: 
these form the right and left whiskering operations in the syntactic 
model.

The typing rules for $\hir(\graph)$ are collected in 
Figures~\ref{fig:hilken-intro-rules}--\ref{ref:hir:cong-rules}.
\vspace{\baselineskip}

\begin{rules}
\unaryRule{\faketext}{x_1 : A_1, \dots, x_n : A_n \vdash x_k : A_k}{var} \vspace{.5\treeskip}

\binaryRule{c \in \graph(A_1, \dots, A_n; B)}{(\Delta \vdash u_i : A_i)_{i=1, \dots, n}}{x_1 : A_1, \dots, x_n : A_n \vdash c(u_1, \dots, u_n) : B}{const}
\unaryRule{\Gamma \vdash t : B}{\Gamma \vdash \id_t : \rewrite{t}{t} : B}{id}

 \binaryRule	{\Gamma \vdash \tau' : \rewrite{t'}{t''} : B}
 {\Gamma \vdash \tau : \rewrite{t}{t'} : B}
 {\Gamma \vdash \tau' \vert \tau : \rewrite{t}{t''} : B}
 {vert-comp}

 \binaryRule		{\constrewr \in \graph(A_1, \dots, A_n; B)(c, c')}
 	{(\Delta \vdash u_i : A_i)_{i = 1, \dots, n}}
 	{\Delta \vdash \constrewr(u_1, \dots, u_n) : \rewrite{c(u_1, \dots, u_n)}{c'(u_1, \dots, u_n)} : B}
 	{right-whisker}
 
 \binaryRule		{c \in \graph(A_1, \dots, A_n; B)}
 	{(\Delta \vdash \sigma_i : \rewrite{u_i}{u_i'} : A_i)_{i = 1, \dots, n}}
 	{\Delta \vdash c(\sigma_1, \dots, \sigma_n) : \rewrite{c(u_1, \dots, u_n)}{c'(u_1, \dots, u_n)} : B}
 	{left-whisker}
\caption{Introduction rules for $\hir(\graph)$.\label{fig:hilken-intro-rules}}
\end{rules}

\begin{figure}[ht]
\begin{mdframed}
\centering
\input{rules/categorical-vert-comp}
\caption{Categorical rules for vertical composition \label{fig:hir:cat-rules}}
\end{mdframed}
\end{figure}

\begin{figure}[ht]
\begin{mdframed}
\centering
\begin{small}
\trinaryRule	{c \in \graph(A_1, \dots, A_n; B)}
	{(\Delta \vdash \sigma_i' : \rewrite{u_i'}{u_i''} : A_i)_{i = 1, \dots, n}}
	{(\Delta \vdash \sigma_i : \rewrite{u_i}{u_i'} : A_i)_{i=1, \dots, n}}
	{\Delta \vdash c(\tau_1', \dots, \tau_n') \vert c(\tau_1, \dots, \tau_n) \equiv c(\tau_1' \vert \tau_1, \dots, \tau_n' \vert \tau_n) : \rewrite{c(u_1, \dots, u_n)}{c(u_1'', \dots, u_n'')} : B}
	{}
\end{small}

\binaryRule	{c \in \graph(A_1, \dots, A_n; B}
{(\Delta \vdash u_i : A_i)_{i=1, \dots, n}}
{\Delta \vdash c(\id_{u_1}, \dots, \id_{u_n}) \equiv \id_{c(u_1, \dots, u_n)} : \rewrite{c(u_1, \dots, u_n)}{c(u_1, \dots, u_n)} : B}
{}

\begin{small}
\binaryRule	{\constrewr \in \graph(A_1, \dots, A_n; B)(c,c')}
{(\Delta \vdash \sigma_i : \rewrite{u_i}{u_i'} : A_i)_{i=1, \dots, n}}
{\Delta \vdash \constrewr(u_1', \dots, u_n') \vert c(\sigma_1, \dots, \sigma_n) \equiv c'(\sigma_1, \dots, \sigma_n) \vert \constrewr(u_1, \dots, u_n) : \rewrite{c(\ind{u})}{c'(\ind{u}')} : B}
{}
\end{small}

\caption{Compatibility laws for constants \label{ref:hir:compat-rules}}	
\end{mdframed}
\end{figure}

\begin{figure}[ht]
\begin{mdframed}
\centering
\unaryRule	{\Gamma \vdash \tau : \rewrite{t}{t'} : A}
{\Gamma \vdash \tau \equiv \tau : \rewrite{t}{t'} : A}
{refl}
\unaryRule	{\Gamma \vdash \tau \equiv \tau' : \rewrite{t}{t'} : A}
{\Gamma \vdash \tau' \equiv \tau : \rewrite{t}{t'} : A}
{symm}
\binaryRule	{\Gamma \vdash \tau' \equiv \tau'' : \rewrite{t}{t'} : A}
{\Gamma \vdash \tau \equiv \tau' : \rewrite{t}{t'} : A}
{\Gamma \vdash \tau \equiv \tau'' : \rewrite{t}{t'} : A}
{trans}

\binaryRule	{\Gamma \vdash \tau' \equiv \sigma' : \rewrite{t'}{t''} : A}
{\Gamma \vdash \tau \equiv \sigma : \rewrite{t}{t'} : A}
{\Gamma \vdash \tau' \vert \tau \equiv \sigma' \vert \sigma  : \rewrite{t}{t''} : A}
{}
\binaryRule	{c \in \graph(A_1, \dots, A_n; B)}
{(\Delta \vdash \sigma_i \equiv \sigma' : \rewrite{u_i}{u_i'} : A_i)_{i=1, \dots, n}}
{\Delta \vdash c(\sigma_1, \dots, \sigma_n) \equiv c(\sigma_1', \dots, \sigma_n') : \rewrite{c(u_1, \dots, u_n)}{c(u_1', \dots, u_n')}}
{}

\caption{Congruence rules \label{ref:hir:cong-rules}}	
\end{mdframed}
\end{figure}


\FloatBlock

For $\hir$ to be a strict biclone we require a strictly associative and unital substitution operation. Accordingly, we define substitution of terms into terms, of terms into rewrites, and 
of rewrites into terms as follows.
\begin{align*}
x_k[u_i / x_i] &:= u_k  \\
c(u_1, \,\dots\, , u_n)[v_j / y_j] &:= c\big(u_1[v_j/y_j], \,\dots\, , 
u_n[v_j/y_j]\big) 
\\[5pt]
\id_t[u_i/x_i] &:= \id_{t[u_i/x_i]} \\
(\tau' \vert \tau)[u_i/x_i] &:= \tau'[u_i/x_i] \vert \tau[u_i/x_i] \\
c(\sigma_1, \,\dots\, , \sigma_n)[u_i/x_i] &:= c\big(\sigma_1[u_i/x_i] \dots, 
\sigma_n[u_i/x_i]\big)\\
\sigma(u_1, \,\dots\, , u_n)[v_j/y_j] &:= \sigma\big(u_1[v_j/y_j], \,\dots\, , 
u_n[v_j/y_j]\big) \\[5pt]
x_k[\sigma_i / x_i] &:= \sigma_k \\
c(u_1, \,\dots\, , u_n)[\sigma_j / y_j] &:= c\big(u_1[\sigma_j / y_j], 
\,\dots\, , u_n[\sigma_j/y_j]\big)
\end{align*}
The Substitution Lemma holds for all three forms of substitution. 

\begin{prooflesslemma}
For any 2-multigraph $\graph$, the following rules are admissible in 
$\hir(\graph)$:
\begin{center}
\binaryRule	{x_1 : A_1, \,\dots\, , x_n : A_n \vdash t : B}
{(\Delta \vdash u_i : A_i)_{i = 1, \,\dots\, ,n}}
{\Delta \vdash t[u_i/x_i] : B}
{}
\binaryRule	{x_1 : A_1, \,\dots\, , x_n : A_n \vdash \tau : \rewrite{t}{t'} : B}
{(\Delta \vdash u_i : A_i)_{i = 1, \,\dots\, ,n}}
{\Delta \vdash \tau[u_i/x_i] : \rewrite{t[u_i/x_i]}{t'[u_i/x_i]} : B}
{}
\binaryRule	{x_1 : A_1, \,\dots\, , x_n : A_n \vdash t : B}
{(\Delta \vdash \sigma_i : \rewrite{u_i}{u_i'} : A_i)_{i = 1, \,\dots\, ,n}}
{\Delta \vdash t[\sigma_i/x_i] : \rewrite{t[u_i/x_i]}{t[u_i'/x_i]} : B}
{} 

\hfill\qedhere
\end{center}
\end{prooflesslemma}

As there are no operations that bind variables, the definition of 
$\alpha$-equivalence is trivial. The equational theory $\equiv$ is defined in 
Figures~\ref{fig:hir:cat-rules}--\ref{ref:hir:cong-rules}. The rules diverge 
from $\langBiclone$ most importantly in 
Figure~\ref{ref:hir:compat-rules}, which 
ensures the 
meta-operation of substitution is functorial,
and that the two different ways of composing with constant rewrites are equal. 
This guarantees that the composites 
$\tau[u_i' / x_i] \vert t[\sigma_i / x_i]$ and 
$t'[\sigma_i / x_i] \vert \tau[u_i/ x_i]$ coincide 
(\cf~the \emph{permutation equivalence} of~\cite{Hirschowitz2013}). 

Following the pattern of~\cite{Hilken1996, Hirschowitz2013}, we define a 
substitution operation making the following rule admissible, where
$\tau[\sigma_i / x_i] := t'[\sigma_i/ x_i] \vert \tau[u_i/x_i]$:
\begin{center}
\binaryRule	{x_1 : A_1, \,\dots\, , x_n : A_n \vdash \tau : \rewrite{t}{t'} : B}
{(\Delta \vdash \sigma_i : \rewrite{u_i}{u_i'} : A_i)_{i=1, \,\dots\, , n}}
{\Delta \vdash \tau[\sigma_i/x_i] : \rewrite{t[u_i/x_i]}{t'[u_i'/x_i]} : B}
{subst} \vspace{-\treeskip}
\end{center} 
We could have defined vertical composition by whiskering in the 
opposite order, thus:
$\tau[\sigma_i / x_i] := \tau[u_i' / x_i] \vert t[\sigma_i / x_i]$.
The next lemma guarantees that these two coincide. The proof is by structural 
induction, using 
Figure~\ref{ref:hir:compat-rules} for the constant cases.

\begin{prooflesslemma}
For any 2-multigraph $\graph$, the following rule is admissible in $\hir(\graph)$:
\begin{center}
\binaryRule	
	{x_1 : A_1, \,\dots\, , x_n : A_n \vdash \tau : \rewrite{t}{t'} : B}
	{(\Delta \vdash \sigma_i : \rewrite{u_i}{u_i'} : A_i)_{i=1, \,\dots\, , n}}
	{\Delta \vdash 
		t'[\sigma_i/x_i] \vert \tau[u_i/x_i] 
			\equiv 
		\tau[u_i'/x_i] \vert t[\sigma_i/x_i] : 
		\rewrite{t[u_i/x_i]}{t'[u_i'/x_i]} : B}
{}  \vspace{-\treeskip}
\end{center}\qedhere
\end{prooflesslemma}

Further structural inductions establish the key properties we shall be relying on. 

\begin{prooflesslemma} \label{lem:clone-structure-of-langBiclone}
For any 2-multigraph $\graph$ and terms $t, u_1, \,\dots\, , u_n$ in 
$\langBiclone(\graph)$:
\begin{enumerate}
\item $x_k[u_i/x_i] = u_k$,
\item $t[x_i/x_i] = t$, 
\item $t[u_i/x_i][v_j/y_j] = t\big[u_i[v_j/y_j]/x_i\big]$.
\end{enumerate}
Moreover, for any rewrites $\tau, \sigma_1, \,\dots\, , \sigma_n$,
\begin{enumerate}
\item $\id_{x_k}[\sigma_i/x_i] \equiv \sigma_k$,
\item $\tau[\id_{x_i}/x_i] \equiv \tau$, 
\item $\tau[\sigma_i/x_i][\mu_j/y_j] \equiv \tau\big[\sigma_i[\mu_j/y_j]/x_i\big]$. \qedhere
\end{enumerate}
\end{prooflesslemma}

Hence the three laws of an abstract clone hold on both terms and rewrites. It 
is similarly straightforward to establish that $t[\sigma_i' \vert \sigma_i / 
x_i] \equiv t[\sigma_i' / x_i] \vert t[\sigma_i /x_i]$ and hence deduce the 
\emph{interchange law} 
$(\tau' \vert \tau)[\sigma_i' \vert \sigma_i / x_i] 
	\equiv \tau'[\sigma_i'/x_i] \vert \tau[\sigma_i / x_i]$. 
Finally we observe that 
$\id_t[\id_{u_i}/ x_i] \equiv \id_{t[u_i/x_i]}$. 
Together these considerations establish the following does indeed define a strict biclone.

\begin{myconstr} \label{constr:free-2-clone}
For any 2-multigraph $\graph$, define a strict biclone $\hirclone(\graph)$ as 
follows. The sorts are nodes in $\graph$. The 1-cells are terms 
$({x_1 : A_1, \,\dots\, , x_n : A_n \vdash t : B})$ derivable in 
$\hir(\graph)$, for $x_1, x_2, \dots$ a chosen enumeration of variables, and 
the 2-cells are ${\equiv}$-classes of rewrites 
$
({x_1 : A_1, \,\dots\, , x_n : A_n} \vdash \tau : \rewrite{t}{t'} : B)
$.
Composition 
is the $\vert$ operation 
and the identity on a term-in-context $t$ is $\id_t$. 

Substitution is the meta-operation of 
substitution in $\hir(\graph)$:
\begin{align*}
t, (u_1, \,\dots\, , u_n) \mapsto t[u_1/x_1, \,\dots\, , u_n /x_n] \\ 
\tau, (\sigma_1, \,\dots\, , \sigma_n) \mapsto \tau[\sigma_1/x_1, \,\dots\, , 
\sigma_n / x_n] 
\end{align*} 
The projections $\p{i}{\ind{A}} : A_1, \,\dots\, , A_n \to A_i$ are given 
by  the \rulename{var} rule.
\end{myconstr}

It is not hard to see that $\hirclone(\graph)$ is the free 2-clone on $\graph$. 

\begin{mylemma} \label{lem:strict-biclone-free-property}
For any 2-multigraph $\graph$, strict biclone $(T, \altBiclone)$ and 
2-multigraph homomorphism $h : \graph \to \altBiclone$, there exists a unique 
strict pseudofunctor $h\sem{-} : \hirclone(\graph) \to \altBiclone$ such that  
$h\sem{-} \circ \inc = h$, for $\inc : \graph \hookrightarrow 
\hirclone(\graph)$ the inclusion.
\begin{proof}
A straightforward adaptation of the proof of Lemma~\ref{lem:strict-pseudofunctor-biclone}. The most significant work is showing that the pseudofunctor $h\sem{-}$ respects substitution, in the sense that 
\begin{align*}
h\sem{\Delta \vdash \tau[\sigma_i / x_i] : \rewrite{t[u_i / x_i]}{t'[u_i'/x_i]} : B} & \\
 &\hspace{-40mm}= \cslr{\big(h\sem{x_1 : A_1, \,\dots\, , x_n : A_n \vdash \tau 
 : 
 \rewrite{t}{t'} : B}\big)}{\Delta \vdash \ind{\sigma} : 
 \rewrite{\ind{u}}{\ind{u}'} : \ind{A}} 
\end{align*}
for all judgements $x_1 : A_1, \,\dots\, , x_n : A_n \vdash \tau : 
\rewrite{t}{t'} : B$ and \mbox{$(\Delta \vdash \sigma_i : \rewrite{u_i}{u_i'} : 
A_i)_{i=1,\dots,n}$}. This is proven by two structural inductions, one for each 
of the whiskering operations.
\end{proof}
\end{mylemma}

\subsection{Proving biequivalence}

The next stage of the proof is to construct a biequivalence of biclones 
$\hirclone(\graph) \simeq \synclone{\graph}$ over a fixed 2-multigraph 
$\graph$. We shall then see how this restricts to a biequivalence 
of bicategories when $\graph$ is a 2-graph and $\hir$ and 
$\langBiclone$ are restricted to unary contexts.

Fix a 2-multigraph $\graph$. We begin by constructing pseudofunctors 
\mbox{$\into{-} : \hirclone(\graph) \leftrightarrows \synclone{\graph} : 
\out{({-})}$}. The definition of $\out{({-})}$ is simpler, so we do this first. 
Intuitively, this mapping is a \Def{strictification} evaluating away explicit 
substitutions; for constants we exploit the fact the underlying signatures are 
the same. 

\begin{myconstr} \label{constr:biclone-out}
For any 2-multigraph $\graph$, we define a mapping from raw terms in 
$\langBiclone(\graph)$ to raw terms in $\hir(\graph)$ as follows: 
\begin{align*}
\out{x_k} &:= x_k \\
\out{c(x_1, \,\dots\, , x_n)} &:= c(x_1, \,\dots\, , x_n) \\
\out{\hcomp{t}{x_i \mapsto u_i}} &:= \out{t}[\out{u_i}/x_i]
\end{align*} 
This extends to a map on raw rewrites:
\begin{equation*}
\begin{aligned}[c]
\out{\assoc{t, \ind{u}, \ind{v}}} &:= \id_{\out{t}[\out{u_i}/x_i][\out{v_j}/y_j]} \\
\out{\subid{t}} &:= \id_{\out{t}} \\
\out{\indproj{k}{\ind{u}}} &:= \id_{\out{u_k}} 
\end{aligned}
\qquad\qquad\quad
\begin{aligned}[c]
\out{\id_t} &:= \id_{\out{t}} \\
\out{\constrewr(x_1, \,\dots\, , x_n)} &:= \constrewr(x_1, \,\dots\, , x_n) \\
\out{\tau \vert \sigma} &:= \out{\tau} \vert \out{\sigma} \\
\out{\hcomp{\tau}{x_i \mapsto \sigma_i}} &:= \out{\tau}[\out{\sigma_i} / x_i] 
\end{aligned} 
\end{equation*}
\qedhere
\end{myconstr}


This mapping respects typing and the equational theory. 

\begin{mylemma} \label{lem:out-well-defined}
For any 2-multigraph $\graph$,
\begin{enumerate}
\item For all derivable terms $t, t'$ in $\langBiclone(\graph)$, if $t \aeq t'$ 
then $\out{t} = \out{t'}$, 
\item For all derivable rewrites $\tau, \tau'$ in $\langBiclone(\graph)$, if 
$\tau \aeq \tau'$ then $\out{\tau} = \out{\tau'}$, 
\item If $\Gamma \vdash t : B$ in $\langBiclone(\graph)$ then $\Gamma \vdash \out{t} : B$ in $\hir(\graph)$, 
\item If $\Gamma \vdash \tau : \rewrite{t}{t'} : B$ in $\langBiclone(\graph)$ then $\Gamma \vdash \out{\tau} : \rewrite{\out{t}}{\out{t'}} : B$ in $\hir(\graph)$,
\item If $\Gamma \vdash \tau \equiv \tau' : \rewrite{t}{t'} : B$ in $\langBiclone(\graph)$ then $\Gamma \vdash \out{\tau} \equiv \out{\tau'} : \rewrite{\out{t}}{\out{t'}} : B$ in $\hir(\graph)$. \qedhere
\end{enumerate}
\begin{proof}
By structural induction.
\end{proof}
\end{mylemma}

\begin{mypropn} \label{prop:out-is-a-pseudofunctor}
For any 2-multigraph $\graph$ the mapping $\out{({-})}$ extends to a 
pseudofunctor $\synclone{\graph} \to \hirclone(\graph)$. 
\begin{proof}
By Lemma~\ref{lem:out-well-defined} and the definition of $\out{({-})}$ on 
identities and vertical compositions, the mapping $\out{({-})}$ defines a 
functor \mbox{$\synclone{\graph}(\ind{A}; B) \to \hirclone(\ind{A}; B)$} on 
each hom-category by 
\[
\out{(\Gamma \vdash \tau : \rewrite{t}{t'} : B)} := 
(\Gamma \vdash \out{\tau} : \rewrite{\out{t}}{\out{t'}} : B)
\] 
For preservation 
of projections and substitution, one notes that 
\[
\out{x_1 : A_1, \,\dots\, , x_n : A_n \vdash x_k : A_k} = (x_1 : A_1, \,\dots\, 
, x_n : A_n \vdash x_k : A_k)
\]%
and that, for $\Gamma = (x_i : A_i)_{i=1, \,\dots\, , n}$,
\begin{align*}
\cslr{\out{(\Gamma \vdash t : B)}}{\out{\Delta \vdash u_1 : A_1}, \,\dots\, , 
\out{\Delta \vdash u_n : A_n}} &= \cslr{(\Gamma \vdash \out{t} : B)}{\Delta 
\vdash \out{\ind{u}} : \ind{A}} \\
&= (\Delta \vdash \out{t}[\out{u_i}/x_i] : B) \\
&= \out{\Delta \vdash \hcomp{t}{x_i \mapsto u_i} : B}
\end{align*}
so $\out{({-})}$ is indeed a strict pseudofunctor.
\end{proof}
\end{mypropn}

Now we turn to defining the pseudofunctor $\into{-} : \hirclone(\graph) \to 
\synclone{\graph}$. The mapping we choose makes precise the sense in which 
$\hir$ is a fragment of $\langBiclone$. 

\begin{myconstr}
For any 2-multigraph $\graph$, define a mapping from raw terms in $\hir(\graph)$ to raw terms in $\langBiclone(\graph)$ as follows: 
\begin{align*}
\into{x_k} &:= x_k \\
\into{c(u_1, \,\dots\, , u_n)} &:= \hcomp{c}{ \into{u_1}, \,\dots\, , 
\into{u_n}}
\end{align*} 
Extend this to a map on raw rewrites as follows:
\begin{equation*}
\begin{aligned}[c]
\into{\id_t} &:= \id_{\into{t}} \\
\into{\tau \vert \sigma} &:= \into{\tau} \vert \into{\sigma} 
\end{aligned} 
\qquad\quad\qquad
\begin{aligned}[c]
\into{c(\sigma_1, \,\dots\, , \sigma_n)} &:= \hcomp{c}{x_i \mapsto 
\into{\sigma_i}} \\
\into{\constrewr(u_1, \,\dots\, , u_n)} &:= \hcomp{\constrewr}{x_i \mapsto 
\into{u_i}}
\end{aligned}
\end{equation*}
\qedhere
\end{myconstr}

Once again, the mapping respects typings and the equational theory. 

\begin{prooflesslemma} \label{lem:into-well-defined}
For any 2-multigraph $\graph$,
\begin{enumerate}
\item For all derivable terms $t, t'$ in $\hir(\graph)$, if $t = t'$ then 
$\into{t} \aeq \into{t'}$, 
\item For all derivable rewrites $\tau, \tau'$ in $\hir(\graph)$, if 
$\tau = \tau'$ then $\into{\tau} \aeq \into{\tau'}$, 
\item If $\Gamma \vdash t : B$ in $\hir(\graph)$ then $\Gamma \vdash \into{t} : B$ in $\langBiclone(\graph)$, 
\item If $\Gamma \vdash \tau : \rewrite{t}{t'} : B$ in $\hir(\graph)$ then $\Gamma \vdash \into{\tau} : \rewrite{\into{t}}{\into{t'}} : B$ in $\langBiclone(\graph)$,
\item If $\Gamma \vdash \tau \equiv \tau' : \rewrite{t}{t'} : B$ in $\hir(\graph)$ then $\Gamma \vdash \into{\tau} \equiv \into{\tau'} : \rewrite{\into{t}}{\into{t'}} : B$ in $\langBiclone(\graph)$. \qedhere
\end{enumerate}
\end{prooflesslemma}

It is immediate from the preceding lemma that $\into{-}$ defines a functor 
\mbox{$\hirclone(\graph)(\ind{A}; B) \to \synclone{\graph}(\ind{A}; B)$} on 
each 
hom-category, 
and that
$\into{-}$ strictly preserves 
identities. For preservation of substitution, however, we are 
required to construct a family of 2-cells \mbox{$\hcomp{\into{t}}{x_i \mapsto 
\into{u_i}} \To \into{t[u_i/x_i]}$}. This should be compared 
to~\cite{Ritter1997}, where a similar 
translation is constructed at the meta-level.

\newpage
\begin{myconstr} \label{constr:def-of-sub}
For any 2-multigraph $\graph$, define a family of rewrites $\subName$ in 
$\langBiclone(\graph)$ so that the rule
\begin{center}
\binaryRule	{x_1 : A_1, \,\dots\, , x_n : A_n \vdash \into{t} : B}
	{(\Delta \vdash \into{u_i} : A_i)_{i=1, \,\dots\, , n}}
	{\Delta \vdash \sub{t}{\ind{u}} : 
		\rewrite
			{\hcomp{\into{t}}{x_i \mapsto \into{u_i}}}
			{\into{t[u_i/x_i]}} 
		: B}
	{} \vspace{-\treeskip}
\end{center}
is admissible by setting
\begin{align*}
\sub{x_k}{\ind{u}} &:= \hcomp{x_k}{x_i \mapsto \into{u_i}} \XRA{\indproj{k}{\into{\ind{u}}}} \into{u_k} \\
\sub{c(\ind{u})}{\ind{v}} &:= \hcompthree{c}{u_i}{v_j} \XRA{\assoc{c(\ind{x}), \ind{u}, \ind{v}}} \hcomp{c}{\hcomp{u_i}{v_j}} \XRA{\hcomp{c}{\sub{u_i}{\ind{v}}}} \hcomp{c}{\into{u_i[v_j / y_j]}} \qedhere
\end{align*}
\end{myconstr}

We establish the various properties required of $\mathsf{sub}$ by induction. 
The naturality of structural rewrites implies the following.

\begin{prooflesslemma}
For any 2-multigraph $\graph$, the following judgements are derivable in 
$\synclone{\graph}$:
\begin{center}
\binaryRule	
	{\Gamma \vdash \into{t} : B}
	{(\Delta \vdash \into{\sigma_i} : 
		\rewrite{\into{u_i}}{\into{u_i'}} : A_i)_{i=1, \,\dots\, , n}}
	{\Delta \vdash 
		\sub{t}{\ind{u}'} \vert \hcomp{\into{t}}{\into{\sigma_i}} 
			\equiv 
		\into{t[\sigma_i/x_i]} \vert \sub{t}{\ind{u}} : 
	\rewrite{\hcomp{\into{t}}{\into{u_i}}}{\hcomp{\into{t'}}{\into{u_i}}} : B}
	{} \vspace{\treeskip}

\binaryRule	
	{\Gamma \vdash \into{\tau} : \rewrite{\into{t}}{\into{t'}} : B}
	{(\Delta \vdash \into{u_i} : A_i)_{i=1, \,\dots\, , n}}
	{\Delta \vdash 
		\sub{t'}{\ind{u}} \vert \hcomp{\into{\tau}}{\into{u_i}} 
			\equiv 
		\into{\tau[u_i/x_i]} \vert \sub{t}{\ind{u}} : 
		\rewrite
			{\hcomp{\into{t}}{\into{u_i}}}
			{\hcomp{\into{t}}{\into{u_i'}}} : B}
	{}
\end{center}
Hence the following judgement is derivable:
\begin{center}
\binaryRule	{\Gamma \vdash \into{\tau} : \rewrite{\into{t}}{\into{t'}} : B}
{(\Delta \vdash \into{\sigma_i} : \rewrite{\into{u_i}}{\into{u_i'}} : 
A_i)_{i=1, \,\dots\, , n}}
{\Delta \vdash \sub{t'}{\ind{u}'} \vert \hcomp{\into{\tau}}{\into{\sigma_i}} \equiv \into{\tau[\sigma_i/x_i]} \vert \sub{t}{\ind{u}} : \rewrite{\hcomp{\into{t}}{\into{u_i}}}{\hcomp{\into{t'}}{\into{u_i'}}} : B}
{}
\end{center}
and the $\subName$ rewrites are natural.
\end{prooflesslemma}

Next we want to prove the three coherence laws for a pseudofunctor. The law for 
$\indproj{i}{}$~(\ref{eq:pseudofunctor:left-unit}) holds by definition. We 
prove the other two laws using correlates of Mac Lane's 
original five axioms of a monoidal category~\cite{MacLane1963}. 

\begin{mylemma} \label{lem:biclone-maclane-laws-derivable}
For any biclone $(S, \biclone)$ the following diagrams commute:
\begin{center}
\begin{tikzcd}
\p{k}{} \arrow{r}{\subid{}} &
\cslr{\p{k}{}}{\p{1}{}, \,\dots\, , \p{n}{}} \\
\cslr{\p{k}{}}{\p{1}{}, \,\dots\, , \p{n}{}} \arrow{u}{\indproj{k}{}} 
\arrow[equals]{ur} &
\:
\end{tikzcd}
\qquad\quad
\begin{tikzcd}
\cslr{\p{k}{}}{\p{1}{}, \,\dots\, , \p{n}{}} \arrow{r}{\indproj{k}{}} &
\p{k}{} \\
\p{k}{}  \arrow{u}{\subid{}} \arrow[equals]{ur} &
\:
\end{tikzcd}
\end{center}
\begin{center}
\begin{tikzcd}
\csthree{t}{\ind{u}}{\p{1}{}, \,\dots\, , \p{n}{}} \arrow{r}{\assoc{}} &
\cslr{t}{\cslr{\ind{u}}{\p{1}{}, \,\dots\, , \p{n}{}}} \\
\cslr{t}{\ind{u}} \arrow{u}{\subid{}} \arrow[swap]{ur}{\cslr{t}{\subid{}, 
\dots, \subid{}}} &
\:
\end{tikzcd}
\qquad\quad
\begin{tikzcd}
\cslr{\cslr{\p{k}{}}{\ind{u}}}{\ind{v}} \arrow{r}{\indproj{k}{}} &
\cslr{u_k}{\ind{v}} \\
\csthree{\p{k}{}}{\ind{u}}{\ind{v}} \arrow{u}{\assoc{}} 
\arrow[swap]{ur}{\cslr{\indproj{k}{}}{\ind{v}}} &
\:
\end{tikzcd}
\end{center}
\begin{proof}
By adapting Kelly's arguments for monoidal categories~\cite{Kelly1964}.
\end{proof}
\end{mylemma}

\begin{mylemma} \label{lem:pseudofun-laws-for-sub}
For any 2-multigraph $\graph$ and derivable terms \mbox{$(x_1 : A_1, \,\dots\, 
, x_n  : A_n \vdash \into{t} : C)$}, \mbox{$(y_1 : B_1, \,\dots\, , y_m : B_m 
\vdash u_i : A_i)_{i=1,\dots, m}$} and $(\Delta \vdash v_j : B_j)_{j = 1, 
\,\dots\, , m}$ in $\langBiclone(\graph)$, the following diagrams commute in 
$\synclone{\graph}$:
\begin{center}
\begin{tikzcd}[column sep = large]
\hcomp{\into{t}}{x_i \mapsto x_i} \arrow{r}{\sub{t}{\ind{x}}} &
\into{t} \\
\into{t} \arrow{u}{\subid{}} \arrow[equals]{ur} &
\:
\end{tikzcd}
\qquad
\begin{tikzcd}[column sep = huge]
\hcompthree{\into{t}}{\into{u_i}}{\into{v_j}} \arrow{r}{\hcomp{\sub{t}{\ind{u}}}{v_j}}  \arrow[swap]{d}{\assoc{}} &
\hcomp{\into{t[u_i / x_i]}}{\into{v_j}} \arrow{dd}{\sub{t[u_i/x_i]}{\ind{v}}}  \\
\hcomp{\into{t}}{\hcomp{\into{u_i}}{\into{v_j}}} \arrow[swap]{d}{\hcomp{\into{t}}{\sub{u_i}{\ind{v}}}} &
\: \\
\hcomp{\into{t}}{\into{u_i[v_j / y_j]}} \arrow[swap]{r}{\sub{t}{\ind{u}[v_j/y_j]}} &
\into{t\big[u_i[v_j/y_j]/x_i\big]}
\end{tikzcd} 
\end{center}
\begin{proof}
Both claims are proven by induction using the laws of Lemma~\ref{lem:biclone-maclane-laws-derivable}. For the unit law one uses the two laws on $\subid{}$; for the associativity law one uses naturality and the law relating $\indproj{i}{}$ and $\assoc{}$. 
\end{proof}
\end{mylemma}

We have shown that $\mathsf{sub}$ is natural and satisfies the three laws of a pseudofunctor. 

\begin{prooflesscor}
For any 2-multigraph $\graph$ the mapping $\into{-}$ extends to a 
pseudofunctor \mbox{$\hirclone(\graph) \to \synclone{\graph}$}. 
\end{prooflesscor}

\paragraph{Relating the two composites.}
With the two pseudofunctors in hand, we next examine the composites 
$\into{-} \circ \out{({-})}$ and $\out{({-})} \circ \into{-}$. Our first 
observation is that the strictification of an already-strict term $\into{t}$ is 
simply $t$. 

\begin{mylemma} \label{lem:composite-to-identity}
For any 2-multigraph $\graph$, the composite $\out{({-})} \circ \into{-}$ is 
the identity on $\hirclone(\graph)$. 
\begin{proof}
On objects the claim is trivial. On multimaps one proceeds inductively:
\begin{gather*}
x_k \mapsto \into{x_k} = x_k \mapsto \out{x_k} = x_k \\
c(u_1, \,\dots\, , u_n) \mapsto \hcomp{c}{\into{u_1}, \,\dots\, , \into{u_n}} 
\mapsto 
c(x_1, \,\dots\, , x_n)\left[ \out{\into{u_i}} /x_i \right] = c(u_1, \,\dots\, 
, u_n)
\end{gather*}
The induction for 2-cells is similar:
\begin{gather*}
\id_t \mapsto \id_{\into{t}} \mapsto \id_{\out{\into{t}}} = \id_t 
\qquad\qquad\qquad 
\text{ by the preceding } \\
\tau' \vert \tau \mapsto \into{\tau'} \vert \into{\tau} \mapsto 
\out{\into{\tau'}} \vert \out{\into{\tau'}} = \tau' \vert \tau 
	\qquad \text{ by inductive hypothesis} \\[1em]
\constrewr(u_1, \,\dots\, , u_n) \mapsto \hcomp{\constrewr}{\into{u_1}, 
\,\dots\, , \into{u_n}} \mapsto \constrewr(x_1, \,\dots\, , 
x_n)[\out{\into{u_i}} / x_i] = \constrewr(u_1, \,\dots\, , u_n) \\
c(\sigma_1, \,\dots\, , \sigma_n) \mapsto \hcomp{c}{\into{\sigma_1}, \,\dots\, 
, \into{\sigma_n}} \mapsto c(x_1, \,\dots\, , x_n)[\out{\into{\sigma_i}} / x_i] 
= c(\sigma_1, \,\dots\, , \sigma_n)
\end{gather*}
\end{proof}
\end{mylemma}

We finish our construction of the biequivalence 
$\hirclone(\graph) \simeq \synclone{\graph}$ by defining an invertible 
pseudonatural transformation 
$\into{-} \circ \out{({-})} \iso \id_{\synclone{\graph}}$. This amounts to 
defining a reduction procedure 
within $\langBiclone(\graph)$ taking a term to one in which explicit 
substitutions occur as far to the left as possible. The $\subName$ rewrites of 
Construction~\ref{constr:def-of-sub} will play a crucial role.

\begin{myconstr} \label{constr:biclone-reduce}
For any 2-multigraph $\graph$, define a rewrite $\mathsf{reduce}$ typed by the 
rule 
\begin{center}
\unaryRule	{\Gamma \vdash t : B}
{\Gamma \vdash \reduce{t} : \rewrite{t}{\into{\out{t}}} : B}
{}
\end{center}
inductively as follows:
\begin{small}
\begin{align*}
\reduce{x_k} &:= x_k \XRA{\id_{x_k}} x_k \\
\reduce{c(x_1, \,\dots\, , x_n)} &:= c(x_1, \,\dots\, , x_n) \XRA{\subid{}} 
\hcomp{c}{x_1, \,\dots\, , x_n} = \out{c(x_1, \,\dots\, , x_n)} \\
\reduce{\hcomp{t}{x_i \mapsto u_i}} &:= \hcomp{t}{x_i \mapsto u_i} \XRA{\hcomp{\reduce{t}}{\reduce{u_i}}} \hcomp{\into{\out{t}}}{x_i \mapsto \into{\out{u_i}}} \XRA{\sub{\out{t}}{\ind{\out{u}}}} \into{\out{t}[\out{u_i} / x_i]}  
\end{align*}
\end{small}
\end{myconstr}

We think of $\mathsf{reduce}$ as a \Def{normalisation} procedure on terms. When 
such a procedure is defined as a meta-operation, it passes through the term
constructors; in $\langBiclone$, it is natural. 

\begin{mylemma}
For any 2-multigraph $\graph$, the following rule is admissible in 
$\langBiclone(\graph)$:
\begin{center}
\unaryRule	{\Gamma \vdash \tau : \rewrite{t}{t'} : B}
{\Gamma \vdash 
\into{\out{\tau}} \vert \reduce{t} \equiv \reduce{t'} \vert \tau : 
\rewrite{t}{\into{\out{t'}}} :  B}
{}
\end{center}
\begin{proof}
By induction on the derivation of $\tau$. For the structural maps one uses the 
fact the structural maps are all natural; for $\subid{}$ and $\assoc{}$ one 
also makes use of the unit and associativity laws of 
Lemma~\ref{lem:pseudofun-laws-for-sub}, respectively. The other cases are 
straightforward.
\end{proof}
\end{mylemma} 

Terms in which no substitutions occur do not reduce any further.

\begin{mylemma}
For any 2-multigraph $\graph$ and judgement 
$\Gamma \vdash t : B$ derivable in 
$\hir(\graph)$, the rule
\begin{center}
\unaryRule	{\Gamma \vdash \into{t} : B}
{\Gamma \vdash \reduce{\into{t}} \equiv \id_{\into{t}} : \rewrite{\into{t}}{\into{t}} : B}
{}
\end{center} 
is admissible in $\langBiclone(\graph)$.
\begin{proof}
The claim is well-typed because $\into{\out{\into{t}}} = \into{t}$ by Lemma~\ref{lem:composite-to-identity}. The result then follows by structural induction: the \rulename{var} case holds by definition, while the \rulename{const} case is just the triangle law of a biclone.
\end{proof}
\end{mylemma}

The $\reduceName$ rewrite is central to our definition of the invertible 
transformation \mbox{$\id_{\synclone{\graph}} \To \into{\out{({-})}}$}; 
the rest of the work is book-keeping. We define a transformation of pseudofunctors 
(Definition~\ref{def:biclone-transformation}) as follows. Take the 
identity $\indproj{1}{B} : B \to B$ on multimaps; as a term this is 
\mbox{$(x_1 : B  \vdash x_1 : B)$}. 
For each $\Gamma := (x_i : A_i)_{i=1, \,\dots\, , n}$ and 
derivable term $(\Gamma \vdash t : B)$ we are now required to give a 2-cell
\[
(\Gamma \vdash \hcomp{x_1}{x_1 \mapsto t} : B) \To (\Gamma \vdash \hcomp{\into{\out{t}}}{x_i \mapsto \hcomp{x_i}{x_i \mapsto x_i}} : B)
\]
For this, take the composite $\pstrans{t}$ defined by
\begin{equation} \label{eq:def-of-pstrans}
\begin{tikzcd}
\hcomp{x_1}{x_1 \mapsto t} 
\arrow[swap, Rightarrow]{d}{\indproj{1}{}} 
\arrow[dashed, Rightarrow]{rr}{\pstrans{t}} &
\: &
\hcomp{\into{\out{t}}}{x_i \mapsto \hcomp{x_i}{x_i \mapsto x_i}} \\
t \arrow[swap, Rightarrow]{r}{\reduce{t}} &
\into{\out{t}} 
\arrow[Rightarrow, swap]{r}{\subid{}} &
\hcomp{\into{\out{t}}}{x_i \mapsto x_i} 
\arrow[Rightarrow,swap]{u}{\hcomp{\into{\out{t}}}{x_i \mapsto \subid{}}}
\end{tikzcd}
\end{equation}
in context $\Gamma$. The composite is natural because $\mathsf{reduce}$ is.

\begin{mycor}
For any 2-multigraph $\graph$, the multimaps $\indproj{1}{B} : B \to B$ 
together 
with the 2-cells $\pstrans{t}$ defined 
in~(\ref{eq:def-of-pseudonatural-transformation}) form an invertible 
transformation $\id_{\synclone{\graph}} \XRA{\iso} \into{\out{({-})}}$. 
\begin{proof}
By induction, the 2-cell $\mathsf{reduce}$ is invertible, so $\pstrans{t}$ is 
invertible for every derivable term $t$. It remains to check the two axioms, 
for which one uses naturality and the laws of 
Lemma~\ref{lem:biclone-maclane-laws-derivable}.
\end{proof}
\end{mycor}

Let us summarise what we have seen in this section. We have a pair of 
pseudofunctors \mbox{$\into{-} : \hirclone(\graph) \leftrightarrows 
\synclone{\graph} : \out{({-})}$} related by invertible transformations 
$\into{-} \circ \out{({-})} \iso \id_{\synclone{\graph}}$ and \newline%
\mbox{$\out{({-})} \circ \into{-} = \id_{\hirclone(\graph)}$}. Together these 
form the claimed biequivalence.

\begin{prooflessthm} \label{thm:biclones-coherence-biequivalence}
For any 2-multigraph $\graph$, the pseudofunctors $\into{-} : 
\hirclone(\graph) \leftrightarrows \synclone{\graph} : \out{({-})}$ form a 
biequivalence of biclones.
\end{prooflessthm}

We restate the result as a statement of coherence in the style of~\cite{Joyal1993}.

\begin{prooflesscor} \label{cor:coherence-for-free-biclones}
For any 2-multigraph $\graph$, the free biclone on $\graph$ is biequivalent to the free strict biclone on $\graph$.
\end{prooflesscor}

We can use Lemma~\ref{lem:biclone-equivalence-fully-faithful} to parlay the 
preceding corollary into a normalisation result for $\langBiclone$. Since we 
have no control over the behaviour of constant rewrites, we restrict to 
2-multigraphs with no surfaces.

\newpage
\begin{mythm}
Let $\graph$ be a 2-multigraph such that for any nodes $A_1, \,\dots\, , A_n, B 
\in \nodes{\graph}$ and edges \mbox{$f, g : A_1, \,\dots\, , A_n \to B$} the 
set $\graph(\ind{A}; B)(f,g)$ of surfaces $f \To g$ is empty. Then 
$\langBiclone(\graph)$ is locally coherent.
\begin{proof}
The approach is standard (\cf~\cite[p. 16]{Leinster2004}). Suppose given a pair 
of rewrites in $\langBiclone(\graph)$ typed by $\Gamma \vdash \tau : 
\rewrite{t}{t'} : B$ and $\Gamma \vdash \sigma : \rewrite{t}{t'} : B$. Since 
there are no constant rewrites, the definition of $\out{(-)}$ entails that 
$\out{\tau} 
= \id_{\out{t}} = \out{\sigma}$ in $\hir(\graph)$. By 
Lemma~\ref{lem:biclone-equivalence-fully-faithful} the pseudofunctor 
$\out{(-)}$ is locally faithful, so $\tau \equiv \sigma$, as required.
\end{proof}
\end{mythm}


Loosely speaking, any diagram of rewrites in $\langBiclone$ formed from 
$\assoc{}, \subid{}, \indproj{i}{}$ and $\id$ using the operations of vertical 
composition and explicit substitution must commute. We shall freely make use of 
this property from now on. 

Adapting the preceding argument to apply to bicategories---and hence recover 
a version of the classic result of~\cite{MacLane1985}---is a minor adjustment. 
Fix a 2-graph $\graph$. Restricting the construction of $\hirclone(-)$ to unary 
contexts and a fixed variable name 
(\cf~Construction~\ref{constr:bicat-termcat}) yields a 2-category; this 
is free on $\graph$ by  
Lemma~\ref{lem:free-bicategory-on-a-multigraph}. Similarly, the biequivalence 
of 
biclones $\into{-} : \hirclone(\graph) \leftrightarrows \synclone{\graph} : 
\out{({-})}$ restricts to a biequivalence of bicategories. One therefore 
obtains the following.

\begin{prooflesscor}  \label{cor:coherence-for-bicats}
For any 2-graph $\graph$, the free bicategory on $\graph$ is biequivalent to the free 2-category on $\graph$.
\end{prooflesscor}

Alternatively, one may observe that since the internal language for bicategories $\langBicat$ is constructed by restricting the internal language $\langBiclone$ for biclones to unary contexts, any composite of the rewrites $\assoc{}, \subid{}$ and $\indproj{i}{}$ in $\langBicat$ must exist in $\langBiclone$. Hence the local coherence of $\langBiclone$ entails the local coherence of $\langBicat$. 

\begin{prooflesscor}
Let $\graph$ be a 2-graph such that for any nodes $A, B \in \nodes{\graph}$ and edges \mbox{$f, g : A \to B$} the set $\graph(A, B)(f,g)$ of surfaces $f \To g$ is empty. Then $\langBicat(\graph)$ is locally coherent. 
\end{prooflesscor}

%

\chapter{A type theory for fp-bicategories} \label{chap:fp-lang}

In this chapter we extend the type theory $\langBiclone$ with finite 
products. We 
develop a theory of product structures in biclones, and use this to synthesise 
our type theory $\langCart$. Along the way we pursue a 
connection with the \Def{representable multicategories} of 
Hermida~\cite{Hermida2000}. Hermida's work can be seen as bridging 
multicategories and monoidal categories; we show that similar connections 
hold between clones and cartesian categories, and also between biclones and 
bicategories with finite products. The resulting translation mediates between 
products presented by biuniversal arrows (in the style of Hermida's 
representability) and the presentation in terms of natural isomorphisms or 
pseudonatural equivalences. 

With this abstract framework in place, we examine its 
implications for the construction 
of an internal language for biclones with finite products and---by 
extension---for bicategories with finite products. The resulting type theory 
provides a calculus for the kind of 
universal-property reasoning commonly employed when dealing with (bi)limits, 
and contrasts with previous work on type-theoretic descriptions of 
2-dimensional 
cartesian (closed) structure, 
in which products are defined by an invertible unit and counit satisfying the 
triangle laws of an adjunction~(\eg~\cite{Seely1987, Hilken1996, 
Hirschowitz2013}).  

\section{fp-Bicategories} \label{sec:fp-bicategories}
Let us begin by recalling the notions of bicategory with finite products 
and product-preserving 
pseudofunctor. It will be convenient to directly consider
all finite products, so that the bicategory is equipped with $n$-ary products 
for each $n \in \Nat$. 
This reduces the need to deal with the equivalent objects given by 
re-bracketing binary products. 
To 
avoid confusion with the `cartesian bicategories' of 
Carboni and Walters~\cite{Carboni1987, Carboni2008}, we call a bicategory 
with all finite products an \Def{fp-bicategory}. (We will, however, freely make 
use of the term `cartesian' when defining finite products in (bi)clones and 
(bi)multicategories.) 

We define $n$-ary products in a bicategory as a bilimit over a discrete 
bicategory (set) with $n$ objects. As we saw in 
Remark~\ref{rem:bilimits-as-biadjoints}, 
this can be expressed equivalently as a 
right biadjoint.  For bicategories $\baseCat_1, \,\dots\, , \baseCat_n$ the 
\Def{product bicategory} 
$\prod_{i=1}^n \baseCat_i$ has objects 
$(B_1, \,\dots\, , B_n) \in \prod_{i=1}^n ob(\baseCat_i)$ and structure given 
pointwise. An 
fp-bicategory is a bicategory $\baseCat$ equipped with a right biadjoint to the 
diagonal pseudofunctor 
$\Delta^n : \baseCat \to \baseCat^{{\times} n} : B 
\mapsto (B, \,\dots\, , B)$ for every $n \in \Nat$. Applying 
Definition~\ref{def:biadjunction-universal-arrow} in this context, one may 
equivalently ask for a biuniversal arrow 
$(\pi_1, \,\dots\, , \pi_n) : 
\Delta^n\big(\prodop_n(A_1, \,\dots\, , A_n)\big) \to (A_1, \,\dots\, , A_n)$ for 
every 
$A_1, \,\dots\, , A_n \in \baseCat \:\: (n \in \Nat)$. 

\begin{mydefn}  \label{def:fp-bicat}
An \Def{fp-bicategory} $\fpBicat{\baseCat}$ is a bicategory $\baseCat$ equipped 
with the following data
for every $A_1, \,\dots\, , A_n \in \baseCat \:\: (n \in \Nat)$: 
\begin{enumerate} 
\item \label{c:fp-bicat:obj} A chosen object $\prodop_n(A_1, \,\dots\, , A_n)$, 
\item \label{c:fp-bicat:proj} Chosen arrows  
$\pi_k : \prodop_n(A_1, \,\dots\, , A_n) \to A_k \:\: (k=1,\dots,n)$, referred 
to as 
\Def{projections}, 
\item \label{c:fp-bicat:adj-equiv} For every $X \in \baseCat$ an adjoint 
equivalence
\begin{equation} \label{diag:fp-bicategory}
\begin{tikzcd}
\baseCat{\left(X, \prodop_n (A_1, \,\dots\, ,A_n)\right)} 
\arrow[bend left = 20]{r}{(\pi_1 \circ -, \,\dots\, , \pi_n \circ -)} 
\arrow[phantom]{r}[xshift=-0.3em]{\adjUp{\simeq}} &
\prod_{i=1}^n \baseCat(X, A_i) 
\arrow[bend left = 20]{l}{\seqlr{-, \,\dots\, , =}}
\end{tikzcd}
\end{equation}
defined by choosing a family of universal arrows we denote 
$\epsilonTimes{} 
= (\epsilonTimesInd{1}{}, \,\dots\, , \epsilonTimesInd{n}{})$.
\end{enumerate}
We call the right adjoint $\seqlr{-, \,\dots\, , =}$ the \Def{$n$-ary 
tupling}.   
\end{mydefn} 

\begin{myremark} \label{rem:strictness-of-products-defs}
The preceding definition admits two degrees of strictness. Requiring the 
equivalence~(\ref{diag:fp-bicategory}) to be an isomorphism, and $\baseCat$ to 
be a 2-category, yields the definition of \Def{2-categorical} ($\CatCat$-enriched) 
products. These products are not strict in the 1-categorical sense, however: as 
the example of $(\Cat, \times, \catOne)$ shows, it 
may not be the case that $(A \times B) \times C = A \times (B \times C)$. In 
this thesis, we shall generally write \emph{strict} to 
mean only that~(\ref{diag:fp-bicategory}) is an isomorphism, and 
specify explicitly when we mean the stronger sense.
\end{myremark}

Explicitly, the universal arrows of~(\ref{diag:fp-bicategory}) may be specified 
as follows. For any finite family of 1-cells 
\mbox{$(t_i : X \to A_i)_{i=1, \,\dots\, , n}$}, one requires a 1-cell 
$\seqlr{t_1, \,\dots\, , t_n} : X \to \prod_n(A_1, \,\dots\, , A_n)$ 
and a family of invertible 2-cells 
\mbox{$(\epsilonTimesInd{k}{t_1, \,\dots\, , t_n} : 
	\pi_k \circ \seqlr{\ind{t}} \To t_k)_{k = 1, \,\dots\, , n}$}. 
\nom{\epsilonTimesInd{k}{t_1, \,\dots\, , t_n}}
	{The $k$th component of the counit for product structure, of type
			$\pi_k \circ \seqlr{\ind{t}} \XRA\iso t_k$}
These 2-cells are universal in the sense that, for any family of 2-cells 
\mbox{$(\alpha_i : \pi_i \circ u \To t_i : \Gamma \to A_i)_{i = 1, \,\dots\, , 
n}$}, 
there exists a 2-cell 
$\transTimes{\alpha_1, \,\dots\, , \alpha_n} : 
	u \To \seqlr{t_1, \,\dots\, , t_n} : {\Gamma \to \prod_{i=1}^n A_i}$, 
	unique 
	such that 
\nom{\transTimes{\alpha_1, \,\dots\, , \alpha_n}}
	{The unique mediating 2-cell
			$u \To \seqlr{t_1, \,\dots\, , t_n}$ corresponding to 
			$\alpha_i : \pi_i \circ u \To t_i \:\: {(i=1, \,\dots\, , n)}$}
\begin{equation} \label{eq:p-trans-ump}
\epsilonTimesInd{k}{t_1, \,\dots\, , t_n} \vert \big(\pi_k \circ 
\transTimes{\alpha_1, \,\dots\, , \alpha_n}\big) = \alpha_k : \pi_k \circ u \To 
t_k
\end{equation}
for $k = 1, \,\dots\, , n$. One thereby obtains a functor $\seqlr{-, \,\dots\, 
, =}$ and an adjoint equivalence as in~(\ref{diag:fp-bicategory}) with counit 
$\epsilonTimes{}= (\epsilonTimesInd{1}{}, \,\dots\, , \epsilonTimesInd{n}{})$ 
and 
unit 
$\transTimes{\id_{\pi_1 \circ t}, \,\dots\, , \id_{\pi_n \circ t}} : 
	t \To \seqlr{\pi_1 \circ t, \,\dots\, , \pi_n \circ t}$. 
\nom{\etaTimes{t}}{The unit for product structure, of type
		$t \XRA\iso \seqlr{\pi_1 \circ t, \,\dots\, , \pi_n \circ t}$}
This defines a \Def{lax} $n$-ary product structure: one merely obtains an 
adjunction in~(\ref{diag:fp-bicategory}). One turns this into a bicategorical 
(\emph{pseudo}) product by further requiring the unit and counit to be 
invertible. The \Def{terminal object} $\termobj$ arises as 
$\prodop_0()$. 

\begin{myremark} \label{rem:unary-products}
Throughout we shall assume that the chosen unary product structure on an 
fp-bicategory is trivial, in the sense that
$\prodop_1(A) = A$, $\seq{t} = t$ and 
$\epsilonTimesInd{1}{A} = \l_A : \Id \circ t \To t$. 
\end{myremark}

\begin{mynotation}  \label{not:products} \quad
\begin{enumerate}
\item We denote the unit $
	\transTimes{\Id_{\pi_1 \circ t}, \,\dots\, , \Id_{\pi_n \circ t}} : 
	t \To \seqlr{\pi_1 \circ t, \,\dots\, , \pi_n \circ t}$
by $\etaTimes{t}$. (We reserve $\eta$ and~$\epsilon$ for the unit and counit of 
exponential structure.)
\item \label{c:types} We write $A_1 \times \dots \times A_n$ or $\prod_{i=1}^n 
A_i$ for $\prod_n (A_1, \,\dots\, , A_n)$,
\item \label{c:maps} We write $\seqlr{f_i}_{i=1, \,\dots\, , n}$ or simply 
$\seqlr{\ind{f}}$ for the $n$-ary tupling $\seqlr{f_1, \,\dots\, , f_n}$, 
\item Following the 1-categorical notation, for any family of 1-cells 
$f_i : A_i \to A_i' \: \: (i=1,\dots,n)$ we write 
$\prod_n(f_1, \,\dots\, , f_n)$ or $\prod_{i=1}^n f_i$ for the 
$n$-ary tupling 
$\seqlr{f_1 \circ \pi_1, \,\dots\, , f_n \circ \pi_n} : 
\prod_{i=1}^n A_i \to \prod_{i=1}^n A_i'$, 
and likewise on 2-cells.  \qedhere
\end{enumerate}
\end{mynotation} 

One must take treat the $\prod_i f_i$ notation with some care. In a 1-category, 
the morphism
$f \times A = f \times \id_A$ is equal to the pairing 
$\seqlr{f \circ \pi_1, \pi_2}$. In an fp-bicategory, this may not be the case:  
$f \times A = f \times \Id_A = \seqlr{f \circ \pi_1, \Id_A \circ \pi_2}$.

\begin{myremark} \label{rem:nAryProductsDeterminedUpToEquiv} 
Like any biuniversal arrow, products are unique up to 
equivalence~(\cf~Lemma~\ref{lem:biuniversal-arrows-unique}).  Explicitly, given
adjoint equivalences ${(g : C \leftrightarrows \prodop_{i=1}^n B_i : h)}$
and ${(e_i : B_i \leftrightarrows A_i : f_i)_{i=1, \,\dots\, , n}}$
in a bicategory $\baseCat$, the composite 
\begin{equation*} 
\begin{tikzcd} 
& 
\baseCat(X, \prod_{i=1}^n B_i) 
\arrow[bend left = 20]{r}{(\pi_1 \circ -, \,\dots\, , \pi_n \circ -)} 
\arrow[bend left = 20]{dl}{h\circ-}
\arrow[phantom]{r}[xshift=0em]{\adjUp{\simeq}} 
& 
\prod_{i=1}^n \baseCat(X, B_i)
\arrow[bend left = 20]{l}{\seqlr{-, \,\dots\, , =}} 
\arrow[bend left = 20]{dr}{\Pi_{i=1}^n(e_i \circ -)} 
\arrow[phantom]{dr}[description]{\adjUp{\simeq}} 
&
\\
\arrow[phantom]{ru}[description]{\adjUp{\simeq}} 
\baseCat(X, C)
\arrow[bend left = 20]{ur}{g\circ-}
& 
& 
& 
\prod_{i=1}^n \baseCat(X, A_i) 
\arrow[bend left = 20]{ul}{\Pi_{i=1}^n(f_i \circ -)} 
\end{tikzcd}
\end{equation*} 
yields an adjoint equivalence
\begin{equation*} 
\begin{tikzcd} 
\baseCat(X, C)
\arrow[bend left = 20]{rr}
{(\,
  ((e_1\circ\pi_1)\circ g)\circ- , \dots , ((e_n\circ \pi_n)\circ g)\circ-
  \,)} 
\arrow[phantom]{rr}[description, xshift=2mm]{\adjUp{\simeq}} 
&& 
\prod_{i=1}^n \baseCat(X, A_i) 
\arrow[bend left = 20]{ll}{h\circ\seqlr{f_1\circ- , \dots , f_n\circ=}} 
\end{tikzcd}
\end{equation*} 
presenting 
$C$ as the product of $A_1,\ldots,A_n$. 
\end{myremark} 

One may generally think of bicategorical product structure as an intensional 
version of the familiar categorical structure, except the usual equations 
(\eg~\cite{Gibbons1997}) are now witnessed by natural families of invertible 
2-cells. It will be useful to have 
explicit names for these 2-cells. 

\newpage
\begin{myconstr}  \label{constr:fuse-post-semantically}
Let $\fpBicat{\baseCat}$ be an fp-bicategory. We define the following families 
of invertible 2-cells:
\begin{enumerate}
\item For $(h_i : Y \to A_i)_{i=1,\dots, n}$ and $g : X \to Y$, we define
\nom{\post{\ind{h}; g}}{The canonical 2-cell 
	$\seqlr{h_1, \,\dots\, , h_n} \circ g \To \seqlr{h_1 \circ g, \,\dots\, , 
	h_n \circ 
	g}$}
\[
\post{\ind{h}; g} : \seqlr{h_1, \,\dots\, , h_n} \circ g \To \seqlr{h_1 \circ 
g, 
\dots, h_n \circ g}
\]
as $\transTimes{\alpha_1, \,\dots\, , \alpha_n}$, where $\alpha_k$
is the composite
\[
\pi_k \circ \left( \seq{h_1, \,\dots\, , h_n} \circ g\right)
\XRA\iso
\left( \pi_k \circ \seq{h_1, \,\dots\, , h_n} \right) 
	\circ g
\XRA{\epsilonTimesInd{k}{} \circ g}
h_k \circ g
\]
for $k=1,\dots,n$.
\item For $(h_i : A_i \to B_i)_{i=1,\dots,n}$ and 
$(g_i : X \to A_i)_{i=1,\dots,n}$, we define 
\nom{\fuse(\ind{h}; \ind{g})}{The canonical 2-cell 
	$ \left(\prod_{i=1}^n h_i\right) \circ \seqlr{g_1, \,\dots\, , g_n} \To 
	\seqlr{h_1 \circ g_1, \,\dots\, , h_n \circ g_n}$}
\[
\fuse(\ind{h}; \ind{g}) : \left(\prodop_{i=1}^n h_i\right) \circ \seqlr{g_1, 
\dots, g_n} \To \seqlr{h_1 \circ g_1, \,\dots\, , h_n \circ g_n}
\] as $\transTimes{\beta_1, \,\dots\, , \beta_n}$, where $\beta_k$ is defined 
by the 
diagram
\vspace{-.5cm}
\begin{td}
\pi_k 
	\circ \left( \left(\prodop_{i=1}^n h_i \right) 
	\circ \seq{g_1, \,\dots\, , g_n}\right)
\arrow[swap]{d}{\iso} 
\arrow{rr}{\beta_k} &
\: &
h_k \circ g_k \\ 
\left( \pi_k \circ \prodop_{i=1}^n h_i \right) \circ \seq{g_1, \,\dots\, , g_n}
\arrow[swap]{r}[yshift=-1mm]{\epsilonTimesInd{k}{} \circ \seq{g_1, \,\dots\, , 
g_n}} 
&
\left(h_k \circ \pi_k\right) \circ \seq{g_1, \,\dots\, , g_n} 
\arrow[swap]{r}[yshift=-1mm]{\iso} &
h_k \circ \left( \pi_k \circ \seq{g_1, \,\dots\, , g_n} \right)
\arrow[swap]{u}{h_k \circ \epsilonTimesInd{k}{}}
\end{td}
for $k=1,\dots,n$.  

\item \label{c:phiTimes-defined} For 
$(h_i : A_i \to B_i)_{i=1,\dots,n}$ and 
$(g_j : X_j \to A_j)_{j=1,\dots,n}$
we define 
\nom{\phiTimes_{\ind{h}, \ind{g}}}{The canonical 2-cell 
	$\big(\prod_{i=1}^n h_i \big) \circ \big( \prod_{i=1}^n g_i \big) \To
	\prod_{i=1}^n (h_i g_i)$ witnessing the pseudofunctorality of 
	$\prod_n(-, \,\dots\, , 	=)$}
\[
\phiTimes_{\ind{h}, \ind{g}} : 
\left(\prodop_{i=1}^n h_i \right) \circ \left( \prodop_{i=1}^n g_i \right) \To
\prodop_{i=1}^n (h_i g_i)
\]
to be the composite 
$\seq{\a^{-1}_{h_1, g_1, \pi_1}, \dots, \a^{-1}_{h_n, g_n, \pi_n}} 
	\vert \fuse{(\ind{h}; g_1 \circ \pi_1, \dots, g_n \circ \pi_n)}$. 
This 2-cell witnesses the pseudofunctoriality of $\prod_n{(-, \,\dots\, , =)}$. \qedhere
\end{enumerate}
\end{myconstr}

Informally, one can use the preceding construction to translate
a sequence of equalities relating the product structure of a cartesian 
category into a composite of 
invertible 2-cells---the difficulty, as outlined in the introduction to this 
thesis, is verifying such a composite satisfies the required coherence laws. As 
a further 
step to simplifying this effort, we observe that each of the 2-cells just 
constructed is 
natural and satisfies the expected equations. The many isomorphisms required to 
state these lemmas in their full bicategorical generality tend to obscure the 
`self-evident' nature of these results, so we state them for 2-categories with 
pseudo (bicategorical) products.

\newpage
\begin{prooflesslemma} \label{lem:PseudoproductCanonical2CellsLaws}
Let $\baseCat$ be a 2-category with finite pseudo-products. Then for all 
families of suitable
1-cells 
$f,g,h,f_i, g_i, h_i \: \: (i=1, \,\dots\, , n)$, the following diagrams 
commute 
whenever they are well-typed:
\begin{figure}[!h]
\centerfloat
\begin{minipage}{0.52\textwidth} 
\begin{equation} \label{c:PostAndId}
\begin{tikzcd} 
\seqlr{f_1, \,\dots\, , f_n} \arrow[swap, equals]{dr}{} \arrow[equals]{r}{} &
\seqlr{f_1, \,\dots\, , f_n} \circ \Id \arrow{d}{\postName} \\
\: &
\seqlr{f_1 \circ \Id,  \dots, f_n \circ \Id}
\end{tikzcd}
\end{equation} 
\end{minipage}
\:\:\quad
\begin{minipage}{0.52\textwidth}
\begin{equation}
\begin{tikzcd} \label{c:FuseAndEta}
\prod_{i=1}^n f_i 
\arrow{r}[yshift=2mm]{(\prod_i f_i) \circ \etaTimes{\Id}} 
\arrow[swap, equals]{dr}{} &
\left(\prod_{i=1}^n f_i\right) \circ \seqlr{\pi_1, \,\dots\, ,  \pi_n} 
\arrow{d}{\fuse} \\
\: &
\seqlr{f \circ \pi_1, \,\dots\, , f_n \circ \pi_n}
\end{tikzcd}
\end{equation}
\end{minipage}
\end{figure}

\phantom{some fake text to move the qed}
\begin{figure}[!h]
\begin{minipage}{0.52\textwidth}
\begin{equation} \label{c:PostAndEta}
\begin{tikzcd}[column sep = 3em]
f \circ g 
\arrow{r}{\etaTimes{f} \circ g} 
\arrow[swap]{dr}{\etaTimes{fg}} &
\seqlr{\pi_1 \circ f, \,\dots\, , \pi_n \circ f} \circ g \arrow{d}{\postName} \\
\: &
\seqlr{\pi_1 \circ f \circ g, \,\dots\, , \pi_n \circ f \circ g}
\end{tikzcd}
\end{equation}
\end{minipage}
\quad
\begin{minipage}{0.44\textwidth}
\begin{equation} \label{c:DoublePostLaw}
\begin{tikzcd}[column sep = 3em]
\seqlr{\ind{f}} \circ g \circ h \arrow{r}[yshift=0mm]{\postName \circ h} 
\arrow[swap]{dr}{\postName} &
\seqlr{\ind{f} \circ g} \circ h \arrow{d}{\postName} \\
\: &
\seqlr{\ind{f} \circ g \circ h}
\end{tikzcd}
\end{equation}
\end{minipage}
\end{figure}
\vspace{2mm}
\begin{figure}[!h]
\begin{equation} \label{c:FuseMonoidal}
\begin{tikzcd}
\big(\prod_{i=1}^n f_i\big)  \circ \big(\prod_{i=1}^n g_i\big)  
\circ \seqlr{h_1, \,\dots\, , h_n} 
\arrow{r}[yshift=2mm]{\phiTimes_{\ind{f}, \ind{g}} \circ \seqlr{h_1, \,\dots\, 
, 
h_n}} 
\arrow[swap]{d}{(\prod_{i} f_i) \circ \fuse} &
\prod_{i=1}^n (f_i \circ g_i) \circ \seqlr{h_1, \,\dots\, , h_n} 
\arrow{d}{\fuse} \\
\big(\prod_{i=1}^n f_i\big) \circ \seqlr{g_1 \circ h_1, \,\dots\, , g_n \circ 
h_n} 
\arrow[swap]{r}{\fuse} &
\seqlr{f_1 \circ g_1 \circ h_1, \,\dots\, , f_n \circ g_n \circ h_n}
\end{tikzcd}
\end{equation}
\vspace{2mm}
\begin{equation} \label{c:PostAndFuse}
\begin{tikzcd}
\big(\prod_{i=1}^n f_i\big)  \circ \seqlr{g_1, \,\dots\, , g_n} \circ h 
\arrow{r}[yshift=2mm]{(\prod_i f_i) \circ \postName} 
\arrow[swap]{d}{\fuse \circ h} &
\big(\prod_{i=1}^n f_i\big)  \circ \seqlr{g_1 \circ h, \,\dots\, , g_n \circ h} 
\arrow{d}{\fuse} \\
\seqlr{f_1 \circ g_1, \,\dots\, , f_n \circ g_n} \circ h 
\arrow[swap]{r}{\postName} &
\seqlr{f_1 \circ g_1 \circ h, \,\dots\, , f_n \circ g_n \circ h}
\end{tikzcd}
\end{equation}
\hfill
\qedhere
\end{figure}

\end{prooflesslemma}

In Lemma~\ref{lem:properties-of-post} we shall see that these laws hold equally 
within the syntax of the type theory $\langCartClosed$ for fp-bicategories.

The restriction to a base 2-category, rather than a bicategory, turns out 
to be of no great consequence. Indeed, Power's coherence 
result restricts as follows to 
fp-bicategories.

\begin{mypropn}[{\cite[Theorem 4.1]{Power1989bilimit}}] 
\label{prop:power-coherence}
Every fp-bicategory is biequivalent to a 2-category with 
strict (2-categorical) products. 
\begin{proof}
We present Power's proof, adapted to the special case of products. Let 
$\fpBicat{\baseCat}$ be an fp-bicategory. By the Mac Lane-Par{\'e} coherence 
theorem, $\baseCat$ is biequivalent to a 2-category; by 
Lemma~\ref{lem:biequivalences-preserve-biuniversal-arrows}, this is a 
2-category with bicategorical products. We may therefore assume 
without loss of generality that $\fpBicat\baseCat$ is a 2-category with 
bicategorical 
products. Now let 
$\Yon : \baseCat \to \Hom(\op{\baseCat}, \Cat)$ 
be the Yoneda embedding and $\overline{\baseCat}$ be the 
closure 
of 
$ob(\Yon \baseCat)$ in $\Hom(\op{\baseCat}, \Cat)$ under equivalences. The 
Yoneda embedding factors as a composite
$\baseCat \xra{i} \overline{\baseCat} \xra{j} \Hom(\op{\baseCat}, \Cat)$.
Since $\Yon$ is locally an equivalence, the inclusion $i : \baseCat \to 
\overline\baseCat$ is a biequivalence. Choose a pseudoinverse 
$k : \overline\baseCat \to \baseCat$. 

Now, for any 
$P_1, \,\dots\, , P_n \in \overline{\baseCat} \:\: (n \in \Nat)$
a 2-categorical product 
$\prodop_n(jP_1, \,\dots\, , jP_n)$
exists (pointwise) in the 2-category $\Hom(\op\baseCat, \Cat)$: one can show 
this by a direct calculation or by applying general theory 
as in~\cite[Proposition 3.6]{Power1989bilimit} 
(see also Chapter~\ref{chap:calculations}). We show this product also lies 
in
$\overline\baseCat$. Since an isomorphism of hom-categories is certainly an 
equivalence of hom-categories, $\prodop_n(jP_1, \,\dots\, , jP_n)$ is (up to 
equivalence) the bicategorical product of $jP_1, \,\dots\, , jP_n$ in 
$\Hom(\op\baseCat, \Cat)$. Moreover, since $i$ and $k$ form a biequivalence, 
$\Yon \circ k = (j \circ i) \circ k \simeq j \circ \id_{\overline\baseCat} = 
j$. 
So, 
applying the uniqueness of products up to equivalence and the fact that $\Yon$ 
preserves products 
(Lemma~\ref{lem:representables-and-adjoints-preserve-bilimits}): 
\[
\prodop_n(jP_1, \,\dots\, , jP_n) \simeq 
\prodop_n{\left((\Yon k) P_1, \,\dots\, , (\Yon k) P_n \right)} \simeq
\Yon{\left( \prod_n(kP_1, \,\dots\, , kP_n) \right)}
\]
Since 
$\Yon{\left( \prod_n(kP_1, \,\dots\, , kP_n) \right)}$ certainly lies in 
$\overline\baseCat$, 
it follows that 
$\prodop_n(jP_1, \,\dots\, , jP_n)$ also lies in $\overline\baseCat$, as 
claimed. 
\end{proof}
\end{mypropn}

This result obviates the need to deal with the 
various 2-cells of Construction~\ref{constr:fuse-post-semantically}. The reader 
may therefore simplify some of the longer 2-cells we shall construct (for 
example, in Chapter~\ref{chap:glueing}). However, we shall \emph{not} rely on 
it in what follows.

\subsection{Preservation of products}

\paragraph*{fp-Pseudofunctors.} Defining preservation of products is 
straightforward: it is 
just an instance of preservation of bilimits. We ask that for 
each $n \in \Nat$ the biuniversal arrow defining the $n$-ary product is 
preserved. 
Strict preservation of these biuniversal arrows amounts to requiring that the 
chosen product structure in the domain is taken to exactly the chosen product 
structure in the target.

\begin{mydefn}  \label{def:fp-pseudofunctor}
An \Def{fp-pseudofunctor} $(F, \prodPres)$
\nom{\prodPres_{\ind{A}}}
	{An equivalence 
		$\prodop_{i=1}^n (FA_i) \to 
			F\left( \prodop_{i=1}^n A_i \right)$
	forming part of the data of an fp-pseudofunctor}
between \mbox{fp-bicategories} 
$\fpBicat{\baseCat}$ and $\fpBicat{\altCat}$ is a pseudofunctor $F : \baseCat 
\to \altCat$ equipped with specified adjoint equivalences
\[
\seqlr{F\pi_1, \,\dots\, , F\pi_n}
:
F{\left(\prodop_{i=1}^n A_i \right)} 
\leftrightarrows
\prodop_{i=1}^n (FA_i)
: 
\prodPres_{\ind{A}}
\]
for every $A_1, \,\dots\, , A_n \in \baseCat \:\: (n \in \Nat)$. We denote 
the 2-cells witnessing these equivalences as follows:
\begin{align*}
\unTimes_{\ind{A}} : \Id_{(\prod_i FA_i)} \To 
	\seqlr{F\pi_1, \,\dots\, , F\pi_n} \circ \prodPres_{\ind{A}} \\
\coTimes_{\ind{A}} : \prodPres_{\ind{A}} \circ \seqlr{F\pi_1, \,\dots\, , 
F\pi_n} 
	\To \Id_{(F\Pi_i A_i)}
\end{align*}
\nom{\unTimes_{\ind{A}}}
	{A 2-cell $\Id_{(\prod_i FA_i)} \To
		\seqlr{F\pi_1, \,\dots\, , F\pi_n} \circ \prodPres_{\ind{A}}$, 
	part of the data of an \mbox{fp-pseudofunctor} $(F, \prodPres)$}
\nom{\coTimes_{\ind{A}}}
	{A 2-cell $\prodPres_{\ind{A}} \circ \seqlr{F\pi_1, \,\dots\, , F\pi_n} 
		\To \Id_{(F\prod_i A_i)}$, part 
	of the data of an fp-pseudofunctor $(F, \prodPres)$}
We call $(F, \prodPres)$ \Def{strict} if $F$ is strict and satisfies
\begin{align*}
F{\left(\prodop_n(A_1,\ldots,A_n)\right)} &= \prodop_n(FA_1,\ldots,FA_n) \\
F(\pi_i^{A_1,\ldots,A_n}) &= \pi_i^{FA_1,\ldots,FA_n} \\
F\seqlr{t_1,\ldots,t_n} &=\seqlr{Ft_1,\ldots,Ft_n} \\
F\varpi^{(i)}_{t_1,\ldots,t_n} &= \varpi^{(i)}_{Ft_1,\ldots,Ft_n} \\
\prodPres_{A_1,\ldots,A_n} &= \Id_{\Pi_n(FA_1,\ldots,FA_n)}
\end{align*}
with adjoint equivalences canonically induced by the 
$2$-cells 
$\transTimes{\r_{\pi_1},\ldots,\r_{\pi_n}}
 : \Id\XRA\iso\seqlr{\pi_1,\ldots,\pi_n}$.
\end{mydefn} 

By Lemma~\ref{lem:strict-preservation-strict-pres-UMP}, a strict 
fp-pseudofunctor commutes with the $\transTimes{-, \,\dots\, , =}$ 
operation on 2-cells: 
$F{\left(\transTimes{\alpha_1, \dots, \alpha_n}\right)} =
	\transTimes{F\alpha_1, \dots, F\alpha_n}$.

\begin{myremark} \label{rem:fp-pseudofunctor-for-different-prod}
The fact that products are unique up to equivalence has the following 
consequence 
for fp-pseudofunctors. If $\baseCat$ is a bicategory equipped with two product 
structures, say 
$\fpBicat{\baseCat}$ and 
$\big(\baseCat, \altprod_n(-) \big)$, 
then for any fp-pseudofunctor
$(F, \prodPres) : \fpBicat{\baseCat} \to \fpBicat\altCat$ 
there exists an (equivalent) fp-pseudofunctor 
$\big(\baseCat, \altprod_n(-) \big) \to \fpBicat{\altCat}$
with witnessing equivalence
\[
F{\left(\altprod_n(A_1, \,\dots\, , A_n) \right)} 
	\simeq 
F{\left(\prodop_n(A_1, \,\dots\, , A_n)\right)}
	\xra{\prodPres_{\ind{A}}} \prodop_n(FA_1, \,\dots\, , FA_n)
\]
arising from the tupling map 
$\seq{\pi_1, \,\dots\, , \pi_n} : 
\altprod_n(A_1, \,\dots\, , A_n) \to \prodop_n(A_1, \,\dots\, , A_n)$. 
\end{myremark}


We saw in Lemma~\ref{lem:preservation-of-biadjoints-to-equivalence} that, 
when a biadjunction is preserved, one obtains an equivalence of pseudofunctors relating the two biadjunctions. We 
shall make use of the following concrete instance of this fact. 

\begin{mylemma} \label{lem:fp-pseudofunctor-nat-2-cell}
For any fp-pseudofunctor 
$(F, \prodPres) : \fpBicat\baseCat \to \fpBicat{\altCat}$ the family of 
\mbox{1-cells} 
$\prodPres_{\ind{A}} : \prod_{i=1}^n FA_i \to F(\prod_{i=1}^n A_i)$ are the 
components of 
a pseudonatural transformation 
$\prod_{i=1}^n \left(F(-), \,\dots\, , F(=)\right) \To 
	{\left( F \circ \prod_{i=1}^n\right)}{(-, \,\dots\, , =)}$,
and hence an equivalence in
$\Hom{\left(\prod_{i=1}^n \baseCat, \altCat\right)}$. 
\begin{proof}
The witnessing 2-cells 
$\nat_{\ind{f}}$ filling
\nom{\nat_{\ind{f}}}{The 2-cells 
	$\prodPres_{\ind{A}} \circ \prod_{i=1}^n Ff_i \To
		F(\prodop_{i=1}^n f_i) \circ \prodPres_{\ind{A}}$ witnessing that
		$\prod_{i=1}^n \left(F(-), \,\dots\, , F(=)\right) \simeq 
		\left( F \circ \prod_{i=1}^n\right){(-, \,\dots\, , =)}$ for every
		fp-pseudofunctor $(F, \prodPres)$}
\begin{td}[column sep = 3em, row sep = 2em]
	\prod_{i=1}^n FA_i
	\arrow[swap]{d}{\prodPres_{\ind{A}}}
	\arrow{r}{\prod_i Ff_i} &
	\prod_{i=1}^n FA_i' 
	\arrow{d}{\prodPres_{\ind{A}'}}
	\arrow[phantom]{dl}[description]{\twocell{\nat_{\ind{f}}}} \\
	
	F(\prod_{i=1}^n A_i)
	\arrow[swap]{r}{F(\prod_i f_i)} &
	F(\prod_{i=1}^n A_i')
\end{td}
are defined as the following composite:
\begin{td}[column sep = 7em, row sep = 2.5em]
\prodPres_{\ind{A}'} \circ \prod_{i=1}^n Ff_i 
\arrow[swap]{d}{\iso}
\arrow{r}{\nat_{\ind{f}}} &
F{\left(\prodop_{i=1}^n f_i \right)} \circ \prodPres_{\ind{A}} \\

\big(\prodPres_{\ind{A}'} \circ \left(\prod_{i=1}^n Ff_i\right)\big)
	\circ \Id_{(\prod_n F\ind{A})}
\arrow[swap]{d}
	{\prodPres_{\ind{A}'} \circ \left(\prod_{i=1}^n Ff_i\right) 
		\circ \unTimes_{\ind{A}}} &
\Id_{F(\prod_n \ind{A}')} 
	\circ \left(F(\prodop_{i=1}^n f_i)  
	\circ \prodPres_{\ind{A}}\right)
\arrow[swap]{u}{\iso} \\

\big(\prodPres_{\ind{A}'} 
	\circ \prod_{i=1}^n F(f_i)\big) 
	\circ  \left( \seqlr{F(\ind{\pi})} 
 	\circ \prodPres_{\ind{A}} \right)
\arrow[swap]{d}{\iso} &
\left(\prodPres_{\ind{A}'} 
	\circ \seq{F\ind{\pi}}\right) 
	\circ \big(F(\prodop_{i=1}^n f_i)  
	\circ \prodPres_{\ind{A}}\big)
\arrow[swap]{u}
	{\coTimes_{\ind{A}'} \circ F(\prod_i f_i)  \circ \prodPres_{\ind{A}}} \\

\prodPres_{\ind{A}'} 
	\circ \big(\big(\prod_{i=1}^n F(f_i)
	\circ  \seqlr{F(\ind{\pi})} \big) 
 	\circ \prodPres_{\ind{A}} \big)
\arrow[swap]{d}{\prodPres_{\ind{A}'} \circ \fuse \circ \prodPres_{\ind{A}}} &
\prodPres_{\ind{A}'} 
	\circ \big(\big(\seq{F\ind{\pi}} 
	\circ F(\prodop_{i=1}^n f_i) \big) 
	\circ \prodPres_{\ind{A}}\big)
\arrow[swap]{u}{\iso}
 \\

\prodPres_{\ind{A}'} 
	\circ \left(\seqlr{F(\ind{f}) \circ  F(\ind{\pi})} 
 	\circ \prodPres_{\ind{A}}\right) 
\arrow[swap]{d}
	{\prodPres_{\ind{A}'} \circ \seqlr{\phi^F_{\ind{f}; \ind{\pi}}} 
	 	\circ \prodPres_{\ind{A}}} 	&
\prodPres_{\ind{A}'} 
	\circ \left(\seqlr{F(\ind{\pi}) 
 	\circ F{\left(\prodop_{i=1}^n f_i\right)}} 
  	\circ \prodPres_{\ind{A}}\right) 
\arrow[swap]{u}
	{\prodPres_{\ind{A}'} \circ \fuse^{-1} 
	  	\circ \prodPres_{\ind{A}} } \\

\prodPres_{\ind{A}'} 
	\circ \left(\seqlr{F(\ind{f} \circ \ind{\pi})} 
 	\circ \prodPres_{\ind{A}}\right)
\arrow[swap]{r}[yshift=-2mm]
	{\prodPres_{\ind{A}'} \circ 
		\seqlr{F(\epsilonTimesInd{-1}{}), \,\dots\, , 
		F(\epsilonTimesInd{-n}{})} 
	 	\circ \prodPres_{\ind{A}} } &
\prodPres_{\ind{A}'} 
	\circ \left(\seqlr{F(\ind{\pi} \circ \prodop_{i=1}^n f_i)} 
 	\circ \prodPres_{\ind{A}}\right)
 \arrow[swap]{u}
 	{\prodPres_{\ind{A}'} \circ \seqlr{(\phi^F_{\ind{\pi}; \prod_i f_i})^{-1}} 
  		\circ \prodPres_{\ind{A}}} 
\end{td}
\end{proof}
\end{mylemma}

In a cartesian category it is is often useful to `unpack'
an $n$-ary tupling from inside a cartesian functor in the following manner:
\begin{align*}
\seqlr{F\pi_1, \,\dots\, , F\pi_n} \circ F\seq{f_1, \,\dots\, , f_n}
&= \seqlr{F(\ind{\pi}) \circ F\seqlr{f_1, \,\dots\, , f_n}} \\
&= \seqlr{F(\ind{\pi} \circ \seqlr{f_1, \,\dots\, , f_n})} \\
&= \seqlr{Ff_1, \,\dots\, , Ff_n}
\end{align*}
In an fp-bicategory, one obtains a natural family of 2-cells we call 
$\unpack$.  

\enlargethispage*{2\baselineskip}
\begin{myconstr} \label{constr:def-of-unpack}
For any fp-pseudofunctor $F : \fpBicat{\baseCat} \to \fpBicat{\altCat}$ 
the invertible 2-cell
$\unpack_{\ind{f}} : \seqlr{F\pi_1, \,\dots\, , F\pi_n} \circ 
F\seq{f_1, \,\dots\, , f_n} 
\To \seqlr{Ff_1, \,\dots\, , Ff_n} : {FX \to \prod_{i=1}^n FB_i}$
\nom{\unpack_{\ind{f}}}{The 2-cell
	$\seqlr{F\pi_1, \,\dots\, , F\pi_n} \circ F\seqlr{f_1, \,\dots\, , f_n} 
	\To \seqlr{Ff_1, \,\dots\, , f_n}$ `unpacking' an $n$-ary tupling} 
is defined to be 
$\transTimes{\tau_1, \,\dots\, , \tau_n}$, where $\tau_k \:\: 
(k=1, \,\dots\, , n)$ is 
given by the following diagram:
\begin{td}
\pi_k
	\circ \left(\seqlr{F\pi_1, \,\dots\, , F\pi_n} 
	\circ F{\seqlr{f_1, \,\dots\, , f_n}}\right)
\arrow{r}{\tau_k}
\arrow[swap]{d}{\iso} &
Ff_k \\
\left( \pi_k 
	\circ \seqlr{F\pi_1, \,\dots\, , F\pi_n}\right)
	\circ F{\seqlr{f_1, \,\dots\, , f_n}}
\arrow[swap]{d}{\epsilonTimesInd{k}{} \circ F\seqlr{f_1, \,\dots\, , f_n}} & 
\: \\
F(\pi_k) \circ F{\seqlr{f_1, \,\dots\, , f_n}} 
\arrow[swap]{r}{\phi^F_{\pi_k, \seqlr{\ind{f}}}} &
F\left(\pi_i \circ \seqlr{f_1, \,\dots\, , f_n}\right) 
\arrow[swap]{uu}{F\epsilonTimesInd{k}{}}
\end{td}
\end{myconstr}

As with the 2-cells of Construction~\ref{constr:fuse-post-semantically}, it is 
useful to have certain coherence properties ready-made. For $\unpack$ one has 
the following.

\begin{prooflesslemma} \label{lem:unpack-prod-pres}
For any fp-pseudofunctor 
$(F, \prodPres) : \fpBicat{\baseCat} \to \fpBicat{\altCat}$ 
and family of 1-cells
$(f_i : X_i \to Y_i)_{i=1, \,\dots\, , n}$ in $\baseCat$, the following diagram commutes:
\begin{equation*}
\begin{tikzcd}[column sep = 5em]
\left(\seqlr{F\pi_1, \,\dots\, , F\pi_n} 
	\circ F(\prod_{i=1}^n f_i) \right)
	\circ \prodPres_{\ind{X}}  
\arrow[swap]{d}{\iso} 
\arrow{r}{\unpack \circ \prodPres_{\ind{X}}} &
\seqlr{F(f_1 \circ \pi_1), \,\dots\, , F(f_n \circ \pi_n)} \circ 
\prodPres_{\ind{X}} 
\\
\seqlr{F\pi_1, \,\dots\, , F\pi_n} 
	\circ \left(  F(\prod_{i=1}^n f_i)
	\circ  \prodPres_{\ind{X}} \right)
\arrow[swap]{d}{\seqlr{F\pi_1, \,\dots\, , F\pi_n} \circ \nat_{\ind{f}}}  &
\seqlr{Ff_1 \circ F\pi_1, \,\dots\, , Ff_n \circ F\pi_n} \circ 
\prodPres_{\ind{X}} 
\arrow[swap]{u}
	{\seqlr{\phi^F_{f_1, \pi_1},\dots, \phi^F_{f_n, \pi_n}} \circ 
	\prodPres_{\ind{X}}}  \\
\seqlr{F\pi_1, \,\dots\, , F\pi_n} \circ \left(\prodPres_{\ind{Y}} 
	\circ \left(\prod_{i=1}^n Ff_i\right)\right)  
\arrow[swap]{d}{\iso} &
\left((\prod_{i=1}^n Ff_i) 
	\circ \seqlr{F\pi_1, \,\dots\, , F\pi_n}\right) 
	\circ \prodPres_{\ind{X}}
\arrow[swap]{u}{\fuse \circ \prodPres_{\ind{X}}} \\
\left(\seqlr{F\pi_1, \,\dots\, , F\pi_n} \circ \prodPres_{\ind{Y}}\right) 
	\circ \left(\prod_{i=1}^n Ff_i\right)  
\arrow[swap]{d}{(\unTimes_{\ind{Y}})^{-1} \circ (\prod_i Ff_i)} &
(\prod_{i=1}^n Ff_i) 
	\circ \left(\seqlr{F\pi_1, \,\dots\, , F\pi_n} 
	\circ \prodPres_{\ind{X}}\right)
\arrow[swap]{u}{\iso} \\
\Id_{(\prod_i FY_i)} \circ \left(\prod_{i=1}^n Ff_i\right)
\arrow[swap]{r}{\iso} &
(\prod_{i=1}^n Ff_i) \circ \Id_{(\prod_i FX_i)}
\arrow[swap]{u}{(\prod_i Ff_i) \circ \unTimes_{\ind{X}}} \\
\end{tikzcd}
\end{equation*}
\end{prooflesslemma}

\paragraph*{Morphisms of fp-pseudofunctors.}
The tricategorical nature of $\Bicat$ leads naturally to a consideration of 2- 
and 3-cells relating fp-pseudofunctors. Experience from the 1-categorical 
setting, however, suggests that 
new definitions are not 
needed. For cartesian functors 
$F, G : \fpBicat{\catC} \to \fpBicat{\catD}$ it is elementary to check that 
every natural transformation 
$\alpha : F \To G$ satisfies
\begin{equation} \label{eq:all-nat-trans-are-cartesian}
\begin{tikzcd}
	F\big(\prod_{i=1}^n A_i\big) \arrow{rr}[yshift=0mm]{\seqlr{F\pi_1, 
	\,\dots\, , 
	F\pi_n}} \arrow[swap]{d}{\CCCTrans_{(\prod_n \ind{A})}} & 
	\: & 
	\prod_{i=1}^n F(A_i) \arrow{d}{\prod_{i=1}^n \CCCTrans_{A_i}} \\ 
	G\big(\prod_{i=1}^n A_i\big) \arrow[swap]{rr}[yshift=0mm]{\seqlr{G\pi_1, 
	\dots, G\pi_n}} & 
	\: & 
	\prod_{i=1}^n G(A_i) 
\end{tikzcd}
\end{equation}

The corresponding bicategorical fact is the following: every pseudonatural 
transformation extends canonically to an 
\Def{fp-transformation}~(\cf~the 
\Def{monoidal pseudonatural transformations} of~\cite[Chapter 3]{Houston2007}).

\newpage
\begin{mydefn} \label{def:fp-transformation}
Let $(F, \prodPres)$ and $(G, \altProdPres)$ be \mbox{fp-pseudofunctors} 
$\fpBicat{\baseCat} \to \fpBicat{\altCat}$. An \Def{fp-transformation}
$(\CCCTrans, \CCCTransNat, \CCCTransProd)$ is a pseudonatural transformation
$(\CCCTrans,\CCCTransNat) : F \To G$ equipped with a 2-cell
$\CCCTransProd_{A_1, \,\dots\, , A_n}$
as in the following diagram for every 
$A_1, \,\dots\, , A_n \in \baseCat \:\: (n \in \Nat)$:
\begin{td}
F{\big(\prod_{i=1}^n A_i\big)}
\arrow[phantom]{drr}[description]{\twocell{\CCCTransProd_{A_1, \,\dots\, , 
A_n}}}
\arrow{rr}[yshift=0mm]{\seqlr{F\pi_1, \,\dots\, , F\pi_n}} 
\arrow[swap]{d}{\CCCTrans_{(\prod_n \ind{A})}} & 
\: & 
\prod_{i=1}^n F(A_i) \arrow{d}{\prod_{i=1}^n \CCCTrans_{A_i}} \\ 
G{\big(\prod_{i=1}^n A_i\big)} 
\arrow[swap]{rr}[yshift=0mm]{\seqlr{G\pi_1, \dots, G\pi_n}} & 
\: & 
\prod_{i=1}^n G(A_i) 
\end{td}
These 2-cells are required to satisfy 
\begin{equation*} 
\begin{tikzcd}[column sep = 2em]
\pi_k 
	\circ \left( \left(\prod_{i=1}^n \CCCTrans_{A_i}\right) 
	\circ \seqlr{F\pi_1, \,\dots\, , F\pi_n}\right) 
\arrow[swap]{d}{\iso} 
\arrow{rr}{\pi_k 
\circ \CCCTransProd_{A_1, \,\dots\, , A_n}} &
\: &
\pi_k 
	\circ \left(\seqlr{G\pi_1, \,\dots\, , G\pi_n} 
	\circ \CCCTrans_{(\prod_n \ind{A})}\right) 
\arrow{d}{\iso}  \\
\left( \pi_k \circ \prod_{i=1}^n \CCCTrans_{A_i} \right)
	\circ \seqlr{F\pi_1, \,\dots\, , F\pi_n}
\arrow[swap]{d}{\epsilonTimesInd{k}{} \circ \seqlr{F\ind{\pi}}} &
\: &
\left(\pi_k 
	\circ \seqlr{G\pi_1, \,\dots\, , G\pi_n}\right) 
	\circ \CCCTrans_{(\prod_n \ind{A})} 
\arrow{dd}{\epsilonTimesInd{k}{}\circ \CCCTrans_{(\prod_n \ind{A})}} \\
\left(\CCCTrans_{A_k} \circ \pi_k\right) \circ \seqlr{F\pi_1, \,\dots\, , 
F\pi_n} 
\arrow[swap]{d}{\iso} &
\: &
\: \\
\CCCTrans_{A_k} 
	\circ \left(\pi_k 
	\circ \seqlr{F\pi_1, \,\dots\, , F\pi_n}\right)
\arrow[swap]{r}[yshift=0mm]{\CCCTrans_{A_k} \circ \epsilonTimesInd{k}{}} &
\CCCTrans_{A_k} \circ F\pi_k 
\arrow[swap]{r}[yshift=0mm]{\CCCTransNat_{\pi_k}} &
G\pi_k \circ \CCCTrans_{(\prod_n \ind{A})}
\end{tikzcd}
\end{equation*}
\end{mydefn}

\begin{mylemma} \label{lem:fp-trans-from-pseudo-trans}
Let $(F, \prodPres)$ and $(G, \altProdPres)$ be \mbox{fp-pseudofunctors} 
$\fpBicat{\baseCat} \to \fpBicat{\altCat}$ and $(\CCCTrans, 
\CCCTransNat) : F \To G$ a pseudonatural transformation. Then, where
$\CCCTransProd_{A_1, \,\dots\, , A_n}$ is defined to be the composite
\begin{td}[column sep = 4.5em]
\left(\prod_{i=1}^n \CCCTrans_{A_i}\right)
	\circ \seqlr{F\pi_1, \,\dots\, , F\pi_n} 
\arrow[swap]{d}{\fuse} \arrow[]{r}{\CCCTransProd_{A_1, \,\dots\, , A_n}} &
\seqlr{G\pi_1, \,\dots\, , G\pi_n} \circ \CCCTrans_{A_1 \times \dots \times 
A_n} 
\\

\seqlr{\CCCTrans_{A_1} \circ F\pi_1, \,\dots\, , \CCCTrans_{A_n} \circ F\pi_n} 
\arrow[swap, yshift=-2mm]{r}{\seqlr{\CCCTransNat_{\pi_1}, \,\dots\, , 
\CCCTransNat_{\pi_n}}} &
\seqlr{G\pi_1 \circ \CCCTrans_{(\prod_n \ind{A})}, 
		\dots, 
		G\pi_n \circ \CCCTrans_{(\prod_n \ind{A})}} 
\arrow[swap]{u}{\postName^{-1}} 
\end{td}
the triple $(\CCCTrans, \CCCTransNat, \CCCTransProd)$ is an 
fp-transformation.
\begin{proof}
A straightforward diagram chase unwinding the definitions of $\fuse$ and 
$\postName$.
\end{proof}
\end{mylemma}

In a similar vein, one might define an \Def{fp-biequivalence} of 
fp-bicategories 
to consist of a pair of fp-pseudofunctors $(F, \prodPres)$ and 
$(G,\altProdPres)$, with fp-transformations
$FG \leftrightarrows \id$ and $GF \leftrightarrows \id$ and invertible 
modifications forming equivalences $FG \simeq \id$ and $GF \simeq \id$. The 
composition of \mbox{fp-transformations} is the usual composition of 
pseudonatural 
transformations, with the composite witnessing 2-cell 
for~(\ref{eq:all-nat-trans-are-cartesian}) given by the evident pasting 
diagram. 
However, this apparently more-structured notion of biequivalence may always be constructed 
from a biequivalence of the underlying bicategories. 

\begin{mylemma} \label{lem:fp-biequivalence-from-biequivalence} 
For any fp-bicategories $\fpBicat{\baseCat}$ and $\fpBicat{\altCat}$, there 
exists an \mbox{fp-biequivalence} 
$\fpBicat{\baseCat} \simeq \fpBicat{\altCat}$ 
if and only if there exists a biequivalence of the 
underlying bicategories.
\begin{proof}
The reverse direction is 
immediate. The forward direction follows from
Lemma~\ref{lem:biequivalences-preserve-biuniversal-arrows} and 
Lemma~\ref{lem:fp-trans-from-pseudo-trans}.
\end{proof}
\end{mylemma}

In this thesis we will only ever be concerned with the existence of a 
biequivalence between fp-bicategories, not its particular structure. It will 
therefore suffice to work with biequivalences throughout.

\section{Product structure from representability}
\label{sec:product-structure-from-rep}

In Chapter~\ref{chap:biclone-lang} we saw that a type theory for 
biclones---and, by restriction to unary contexts, bicategories---could be 
extracted directly from the construction of the free biclone on a signature. 
In order to take a similar 
approach in the case of fp-bicategories, we develop the 
theory of product structures in biclones. 

What does it mean to define products in a biclone? As usual, the categorical 
case is informative. Thinking of (sorted) clones as cartesian versions of 
multicategories suggests that products in a clone ought to arise in a way 
paralleling tensor products in a multicategory. Translating the work of 
Hermida~\cite{Hermida2000} to clones in the most 
na{\" i}ve way possible, one might require 
a family of arrows 
$\rho_{\ind{X}} : X_1, \,\dots\, , X_n \to \prodop_n(X_1, \,\dots\, , X_n)$
in a clone $\clone$ 
inducing isomorphisms
$\clone(X_1, \,\dots\, , X_n; A) \iso 
\clone{\left( \prodop_n(X_1, \,\dots\, , X_n) ; A \right)}$ 
by precomposition. On the other hand, Lambek~\cite{Lambek1989} defines
products in a multicategory $\mCat$ by requiring isomorphisms of the form 
$\mCat{\left(\Gamma; \prodop_n(X_1, \,\dots\, , X_n) \right)} \iso
\prod_{i=1}^n \mCat(\Gamma; A_i)$. Connecting these two approaches to product 
structure will be the focus of the 
next section. 

Taking multicategories as our starting point, we shall study two forms of 
universal property, corresponding to Hermida's and Lambek's definitions, 
respectively. 
We shall show how these notions may be applied to clones and, 
moreover, demonstrate that for clones they actually 
coincide~(Theorem~\ref{thm:representable-clone-cartesian-equivalence}). 
 
Thereafter, in Section~\ref{sec:products-in-stlc}, we shall see how one can 
extract the usual product structure of the 
simply-typed lambda calculus from the theory of such \emph{cartesian clones}. 
This will provide the template for lifting this work to the bicategorical 
setting, and hence for the product structure of the type theory $\langCart$.

\subsection{Cartesian clones and representability} 
\label{sec:cartesian-clones-and-representability}

We start by recalling a little of the theory of (representable) 
multicategories and their relationship to monoidal categories. Extensive 
overviews are available in~\cite{Leinster2004, Yau2016}. 

\paragraph*{Representable multicategories.}
The notion of \Def{multicategory} is a crucial part of Lambek's extended study 
of deductive systems~\cite{Lambek1969, Lambek1980, Lambek1985, Lambek1989}. The 
motivating example takes objects to be types in some sequent calculus and 
multimaps $X_1, \,\dots\, , X_n \vdash Y$ to be derivable sequents; composition 
is 
given by a cut rule. Lambek defines tensor products and (left and right) 
internal homs in a multicategory by the existence of certain natural 
isomorphisms. More recent work by Hermida~\cite{Hermida2000} connects these 
ideas to the categorical setting by making precise the correspondence between 
monoidal categories and so-called \Def{representable} multicategories. 

\begin{mydefn}[{\cite{Lambek1969, Lambek1989}}] 
A \Def{multicategory} $\mCat$ consists of the following data:
\begin{itemize}
\item A set $ob(\mCat)$ of \Def{objects}, 
\item For every sequence $X_1, \,\dots\, , X_n \:\: (n \in \Nat)$ of objects 
and 
object $Y$ a 
\Def{hom-set} $\mCat(X_1, \,\dots\, , X_n; Y)$ consisting of \Def{multimaps} or 
\Def{arrows} 
$f : X_1, \,\dots\, , X_n \to Y$ (here $n$ may be zero). As with (bi)clones, we 
sometimes 
denote sequences $X_1, \,\dots\, , X_n$ by Greek letters $\Gamma, \Delta, 
\dots$ to 
emphasise the connection with contexts,
\item For every $X \in ob(\mCat)$ an \Def{identity multimap} $\id_X : X \to X$,
\item For every set of sequences $\Gamma_1, \,\dots\, , \Gamma_n$ and objects 
$Y_1, \,\dots\, , Y_n, Z$, a \Def{composition} operation
\[ 
\circ_{\ind{\Gamma}; \ind{Y};Z} : 
\mCat(Y_1, \,\dots\, , Y_n; Z) \times \prodop_{i=1}^n \mCat(\Gamma_i; Y_i) 
\to 
\mCat(\Gamma_1, \,\dots\, , \Gamma_n; Z) 
\] 
we denote by 
$\circ_{\ind{\Gamma}; \ind{Y};Z}\big( f, (g_1, \,\dots\, , g_n)\big) 
	:= \ms{f}{g_1, \,\dots\, , g_n}$.
\end{itemize}
This is subject to three axioms requiring that composition is associative and 
unital.
We call multimaps of the form $X \to Y$ \Def{linear}.
\end{mydefn} 

\begin{mynotation} \label{not:multicat-vs-clone-sub}
Note that we write composition in a multicategory as 
$\ms{f}{g_1, \,\dots\, , g_n}$
and substitution in a clone as
$\cs{f}{g_1, \,\dots\, , g_n}$. 
\end{mynotation}

Multicategories are also known as \Def{coloured (planar) operads} 
(\eg~\cite{Yau2016}). Multicategories form a category $\Multicat$ of 
multicategories and their functors, and also a 2-category of multicategories, 
multicategory functors, and transformations (\eg~\cite[Chapter 
2]{Leinster2004}). 

\begin{mydefn} \label{def:functor-and-nat-trans-for-multicat}
\quad 
\begin{enumerate}
\item A \Def{functor} $F : \mCat \to \nCat$ between multicategories $\mCat$ and 
$\nCat$ consists of:
\begin{itemize}
\item A mapping $F : ob(\mCat) \to ob(\nCat)$ on objects,
\item For every \mbox{$X_1, \,\dots\, , X_n, Y \in \mCat \:\: (n \in \Nat)$} a 
mapping on hom-sets 
$F_{\ind{X}; Y} : \mCat(X_1, \,\dots\, , X_n; Y) \to \nCat(FX_1, \,\dots\, , 
FX_n; FY)$,
\end{itemize}
such that composition and the identity are preserved.

\item \label{c:multi-nat-trans} A \Def{transformation} $\alpha : F \To G$ 
between multicategory functors $F, G : \mCat \to \nCat$ is a family of 
multimaps $(\alpha_X : FX \to GX)_{X \in ob(\mCat)}$ such that for every $f : 
X_1, \,\dots\, , X_n \to Y$ the equation $Ff \circ (\alpha_{X_1}, \,\dots\, , 
\alpha_{X_n}) = \alpha_Y \circ (Gf)$ holds. \qedhere
\end{enumerate}
\end{mydefn}

From the perspective of deductive systems, moving from multicategories 
to clones amounts to changing the composition operation from a cut rule to a 
substitution operation. The  
composition operation of a multicategory is \emph{linear}: given maps 
$(h_i : \Gamma \to Y_i)_{i=1, \,\dots\, , m}$ and 
$f : Y_1, \,\dots\, , Y_m \to Z$ in a multicategory, the composite 
$\ms{f}{h_1, \,\dots\, , h_m}$ 
has 
type $\Gamma, \,\dots\, , \Gamma \to Z$.  By contrast, the substitution 
operation in a clone is 
\emph{cartesian}: given maps 
$h_i$ and $f$ as above, the substitution 
$\cslr{f}{h_1, \,\dots\, , h_m}$ has type $\Gamma \to Z$. 

Every multicategory $\mCat$ defines a category $\overline{\mCat}$ by 
restricting 
to linear morphisms. Conversely, every monoidal category 
$(\catC, \otimes, I)$ canonically defines a multicategory with objects those 
of $\catC$ and multimaps $X_1, \,\dots\, , X_n \to Y$ given by morphisms 
$X_1 \otimes \cdots \otimes X_n \to Y$ (for a specified bracketing of the 
$n$-ary tensor product).  A 
natural question is therefore the 
following: under what conditions is the category 
$\overline\mCat$ corresponding to a multicategory monoidal? Hermida answers 
this 
by showing that 
there exists a 2-equivalence between the 2-category $\MonCat$ of monoidal 
categories and the 2-category of representable multicategories.

\begin{mydefn} \label{def:representable-multicategory}
A \Def{representable multicategory} $\mCat$ is a multicategory equipped with a 
chosen object 
$\tensor_n(X_1, \,\dots\, , X_n) \in \mCat$ and a chosen multimap 
$\rho_{X_1, \,\dots\, , X_n} : X_1, \,\dots\, , X_n \to 
\tensor_n(X_1, \,\dots\, , X_n)$ 
for every $X_1, \,\dots\, , X_n \in \mCat \:\: (n \in \Nat)$
such that 
\begin{enumerate} 
\item \label{c:universal-arrow-multimap} Each chosen 
$\rho_{X_1, \,\dots\, , X_n}$ 
is \Def{representable}: for every $Y \in \mCat$, 
precomposition 
with $\rho_{X_1, \,\dots\, , X_n}$ induces an isomorphism 
\mbox{$\mCat(X_1, \,\dots\, , X_n; Y) 
	\iso 
	\mCat{\left(\tensor_n(X_1, \,\dots\, , X_n), Y\right)}$} of hom-sets, and
\item \label{c:universal-arrow-closed-composition}  The representable arrows are closed under composition.  \qedhere
\end{enumerate} 
\end{mydefn} 

Thus, a multimap $\rho_{\ind{X}}$ is representable if and only if for every 
$h : X_1, \,\dots\, , X_n \to Y$ there exists a unique 
multimap 
$\altTrans{h} : \prodop_n(X_1, \,\dots\, , X_n) \to Y$ such that 
$\altTrans{h} \circ 
\rho_{X_1, \,\dots\, , X_n} = h$. 

\begin{myremark}
It is common to refer to the arrows $\rho_{\ind{X}}$ of the preceding 
definition as \Def{universal}; we change the terminology slightly because we 
will imminently define a multicategorical version of universal arrows in the 
sense of Chapter~\ref{chap:background}. The two concepts are related: 
the representability condition~(\ref{c:universal-arrow-multimap}) above is 
equivalent 
to requiring that 
each $\mCat(X_1, \,\dots\, , X_n; -) : \mCat \to \Set$ is representable, which 
is in 
turn equivalent to specifying a universal arrow from the terminal set to this 
functor~(\cf~\cite[Chapter III]{cfwm}). 
\end{myremark}

We briefly recapitulate Hermida's construction.

\newpage
\begin{mylemma}[{\cite[Definition 9.6]{Hermida2000}}] 
\label{lem:multicat-to-monoidal-cat}
For every representable multicategory $\mCat$, the associated category 
$\overline{\mCat}$ is monoidal.
\begin{proof}
The tensor product $X \otimes Y$ is $\tensor_2(X,Y)$ and the unit $I$ arises 
from 
the empty sequence, as 
$\tensor_0()$. The map $f \otimes g$ is defined by the universal property, as 
the unique linear map filling the following diagram:
\begin{td}[column sep = 3em]
\tensor_2(X, Y) \arrow[dashed]{r}{f \otimes g} &
\tensor_2(X', Y') \\
X, Y \arrow{u}{\rho_{X,Y}} \arrow[swap]{r}{(f,g)} &
X', Y' \arrow[swap]{u}{\rho_{X',Y'}}
\end{td}
\end{proof}
\end{mylemma}

The second condition~(\ref{c:universal-arrow-closed-composition}) is necessary: 
it allows one to use the universal property to check the axioms of a monoidal 
category involving iterated tensors $(A \otimes B) \otimes C$ (\cf~the 
preservation conditions for lifting monoidal structure to a category of 
algebras~\cite{Seal}, in particular the \Def{left-linear classifiers} 
of~\cite{ct2018talk}). 

\paragraph*{Cartesian multicategories.}
Representability is a universal property that allows us to construct 
\Def{monoidal} structure. To construct \Def{cartesian} structure, 
however, one requires more. In particular, one ought to obtain 
Lambek's definition of \Def{cartesian multicategory}~\cite[\S 4]{Lambek1989}, 
requiring multimaps 
$\pi_i : \prodop_n(A_1, \,\dots\, , A_n) \to A_i \:\: (i=1, \,\dots\, , n)$
inducing natural isomorphisms 
$\mCat{\left(\Gamma; \prodop_n(X_1, \,\dots\, , X_n) \right)} \iso
\prod_{i=1}^n \mCat(\Gamma; A_i)$.
Next we shall see how to obtain a definition equivalent to Lambek's, but phrased in terms of universal arrows. This will be the starting point for our comparison between product structure and representability.

\begin{mydefn} \label{def:multicategory-universal-arrow}
Let $F : \mCat \to \nCat$ be a functor of multicategories and $X \in \nCat$. A 
\Def{universal arrow from $F$ to $X$} is a pair $(R, u : FR \to X)$ such that 
for every $h : FA_1, \,\dots\, , FA_n \to X$ there exists a unique multimap 
$\trans{h} : A_1, \,\dots\, , A_n \to R$ such that $u \circ (F\trans{h}) = h$. 
\end{mydefn}

\begin{myremark}
One could define universal arrows slightly more generally, by taking a 
universal arrow from $F$ to $X$ to be a \emph{sequence} of objects $R_1, 
\,\dots\, , 
R_n$ with a universal multimap $FR_1, \,\dots\, , FR_n \to X$. The definition 
given seems sufficient for our purposes, so we do not seek this extra 
generality.
\end{myremark}

As in the categorical case, we can rephrase the definition of universal arrow 
as a natural isomorphism. 

\newpage
\begin{mylemma} \label{lem:multicat-universal-arrow-to-natural-isomorphism}
Let $F  : \mCat \to \nCat$ be a functor of multicategories and $X \in \nCat$. 
The following are equivalent:
\begin{enumerate}
\item A specified universal arrow $(R,u)$ from $F$ to $X$,
\item A choice of object $R \in \mCat$ and an isomorphism 
$\mCat(A_1, \,\dots\, , A_n; R) \iso \nCat(FA_1, \,\dots\, , FA_n; X)$, 
multinatural in 
the sense that for 
any $f : A_1, \,\dots\, , A_n \to B$ the following diagram commutes:
\begin{td}
\mCat(B; R) 
\arrow{r}{\iso} 
\arrow[swap]{d}{\ms{(-)}{f}} & 
\nCat(FB; X) 
\arrow{d}{\ms{(-)}{Ff}} \\
\mCat(A_1, \,\dots\, , A_n; R) 
\arrow[swap]{r}{\iso} & 
\nCat(FA_1, \,\dots\, , FA_n; X) 
\end{td}
\end{enumerate}
\begin{proof}
The direction \mbox{(1)\To(2)} is clear. For the reverse, denote the 
isomorphism by $\phi_{\ind{A}} : \mCat(A_1, \,\dots\, , A_n; R) \to \nCat(FA_1, 
\dots, FA_n; X)$ and its inverse by $\psi_{\ind{A}}$. We show that $u := 
\phi_{R}(\id_R) : FR \to X$ is a universal arrow by showing that
that $\psi_{\ind{A}}(-)$ is 
inverse to 
$\ms{\phi_R(\id_R)}{F(-)}$.

First, for any 
$h : FA_1, \,\dots\, , FA_n \to X$, naturality of $\phi$ with respect to the 
multimap 
$\psi_{\ind{A}}(h) : A_1, \,\dots\, , A_n \to R$ 
gives the equation
$\ms{\phi_R(\id_R)}{F\psi_{\ind{A}}(h)} = 
	\phi_{\ind{A}}\psi_{\ind{A}}(h) = h$.
Second, let 
$g : A_1, \,\dots\, , A_n \to R$.  The naturality of $\psi$ with respect to 
$g$ entails that 
$\psi_{\ind{A}}{\left( \ms{\phi_R(\id_R)}{Fg} \right)}
	= \ms{\psi_R\phi_R(\id_R)}{g} 
	= g$, 
as required.
\end{proof}
\end{mylemma}

The category of multicategories $\Multicat$ has products given 
as follows. For multicategories $\mCat$ and $\nCat$ the product 
$\mCat \times \nCat$ has objects pairs 
$(M, N) \in ob(\mCat) \times ob(\nCat)$
and hom-sets
\begin{center}
${(\mCat \times \nCat)}{\left( (A_1, B_1), \,\dots\, , (A_n, B_n); (X, Y) 
\right)} :=
	\mCat(A_1,\dots, A_n; X) \times \nCat(B_1, \,\dots\, , B_n; Y)$
\end{center}
Composition is defined pointwise:
\begin{equation} \label{eq:product-of-multicats-composition}
\small
\begin{tikzcd}[column sep = -12em]
\mCat(\ind{A}; X) \times \nCat(\ind{B}; Y) \times 
\prod_{i=1}^n \left( \mCat(\Gamma_i, A_i) \times \nCat(\Delta_i, B_i) \right)
\arrow[swap]{dr}{\iso}
\arrow{rr}{\circ_{\mCat \times \nCat}} &
\: &
\mCat(\ind{\Gamma}; X) \times \nCat(\ind{\Delta}; Y) \\
\: &
\left( 
	\mCat(\ind{A}; X) \times \prod_{i=1}^n ( \mCat(\Gamma_i, A_i)
\right) \times 
\left( 
	\nCat(\ind{B}; Y) \times \prod_{i=1}^n \nCat(\Delta_i, B_i) 
\right)
\arrow[swap]{ur}{\circ_\mCat \times \circ_\nCat} &
\: 
\end{tikzcd}
\normalsize
\end{equation}
The product structure is then almost identical to that in $\CatCat$. Then for every 
multicategory $\mCat$ and $n \in \Nat$ there exists a diagonal functor 
$\Delta^n : \mCat \to \mCat^{\!{\times} n} : X \mapsto (X, \,\dots\, , X)$, 
and  
Definition~\ref{def:multicategory-universal-arrow} provides a natural 
notion of multicategory with finite products. 

\begin{mydefn} \label{def:fp-multicategory}
A \Def{cartesian multicategory} is a multicategory $\mCat$ equipped with a 
choice of universal arrow 
$\Delta^n\prodop_n(X_1, \,\dots\, , X_n) \to (X_1, \,\dots\, , X_n)$ from
$\Delta^n$ to $(X_1, \dots, X_n)$ for every 
$X_1, \,\dots\, , X_n \in \mCat \:\: (n \in \Nat)$. 
\end{mydefn}

Applying Lemma~\ref{lem:multicat-universal-arrow-to-natural-isomorphism}, 
asking for a multicategory to have finite products is equivalent to asking for 
a chosen sequence of linear multimaps 
$(\pi_i : \prodop_n(X_1, \,\dots\, , X_n) \to X_i)_{i=1, \,\dots\, ,n}$, 
inducing a 
multinatural family of isomorphisms 
\begin{equation} \label{eq:fp-multicategory}
\mCat{\left(\Gamma; \prodop_n(X_1, \,\dots\, , X_n)\right)} 
	\iso 
{\mCat^{\!{\times} n}}
	{\big((\Gamma, \,\dots\, , \Gamma); (X_1, \,\dots\, , X_n)\big)} = 
\prodop_{i=1}^n \mCat(\Gamma; X_i)
\end{equation}
for every $X_1, \,\dots\, , X_n \in \mCat \:\: (n \in \Nat)$.
One thereby recovers Lambek's definition of cartesian products in a 
multicategory~\cite[\S 4]{Lambek1989}.

\paragraph*{Cartesian clones.}
We wish to extend the two definitions we have just seen from multicategories to 
clones. Thinking of (sorted) clones as cartesian versions of multicategories suggests 
the following construction, in which we re-use the notation of 
Notation~\ref{not:clone-product-notation}~(p.~\pageref{not:clone-product-notation}).

\begin{myconstr}  \label{constr:multicategory-from-clone}
Every clone $(S, \clone)$ canonically defines a multicategory $\mCatOf{\clone}$ 
with 
\begin{itemize} 
\item $ob(\mCatOf C) := S$, 
\item $(\mCatOf{\clone})(X_1, \,\dots\, , X_n; Y) := \clone(X_1, \,\dots\, , 
X_n; Y)$ 
\end{itemize} 
Composition is defined as follows. For every 
family of multimaps 
$g_i : \Gamma_i \to Y_i 
	\:\: (i=1, \,\dots\, , n)$ 
and multimap 
$f : Y_1, \,\dots\, , Y_n \to Z$ we define the composite 
$\ms{f}{g_1, \,\dots\, , g_n}$ in $\mCatOf{\clone}$ to be the substitution 
$\cslr{f}{g_1 \clonetimes \cdots \clonetimes g_n}$ in $\clone$. 
The identity $\id_{X,X} \in (\mCatOf C)(X; X)$ is the unary projection 
$\p{1}{} \in \clone(X,X)$, and the axioms follow directly from the three laws of a 
clone. 
\end{myconstr} 

\begin{mynotation}
Motivated by the preceding construction, we shall sometimes write 
$\id_A$ for the projection $\p{1}{1} : A \to A$ in a clone, and refer to it as 
the \Def{identity} on $A$. 
\end{mynotation}

It is clear that this construction extends to a faithful functor
$\mCatOf{(-)} : \Clone \to \Multicat$, yielding a commutative diagram
\begin{equation} \label{eq:clone-multicat-restriction}
\begin{tikzcd}
\Clone 
\arrow[swap]{dr}{\overline{(-)}}
\arrow{rr}{\mCatOf{(-)}} &
\: &
\Multicat
\arrow{dl}{\overline{(-)}} \\
\: &
\CatCat &
\: 
\end{tikzcd}
\end{equation}
in which the downward arrows restrict to unary/linear arrows.
We define  
representability and products in $\Clone$ by applying the definition to the 
image of $\mCatOf{(-)}$.

\begin{mydefn}  \label{def:representable-and-cartesian-clones} \quad
\begin{enumerate}
\item A \Def{representable clone} is a clone $(S, \clone)$ equipped with a choice of representable structure on $\mCatOf\clone$. 
\item A \Def{cartesian clone} is a clone $(S, \clone)$ equipped with a choice of cartesian structure on $\mCatOf\clone$.
\qedhere
\end{enumerate}
\end{mydefn}

\begin{myexmp} \label{ex:fp-cat-to-cartesian-clone}
Every category with finite products $\fpBicat\catC$ defines a clone $\cloneOf\catC$ 
(recall Example~\ref{ex:CartesianCategoryIsAClone}(\ref{c:cart-cat-to-clone}) on page~\pageref{ex:CartesianCategoryIsAClone}). This clone is cartesian, with product structure 
exactly as in $\catC$. 
\end{myexmp}

A clone may therefore be equipped with two kinds of tensor. In the 
representability case, one asks for representable arrows 
$X_1, \,\dots\, , X_n \to \tensor_n(X_1, \,\dots\, , X_n)$. In the cartesian 
case, 
one asks for universal arrows
$\prodop_n(X_1, \,\dots\, , X_n) \to X_i$ for $i=1, \,\dots\, , n$. In terms of 
the 
internal language, these may be thought of as \Def{tupling} and 
\Def{projection} operations, respectively. Identifying representable arrows with a tupling operation (an identification we shall make precise in Corollary~\ref{cor:representable-to-cartesian-clone}), the question then becomes: how does 
one construct a tupling operation given only projections, and how does one 
construct projections given only a tupling operation?

In the light of 
Lemma~\ref{lem:multicat-universal-arrow-to-natural-isomorphism},
we can already construct a tupling operation from projections, and so from cartesian structure. 
If $\mCatOf{\clone}$ has finite products 
witnessed by a universal arrow 
$\pi = (\pi_1, \,\dots\, , \pi_n) : 
	\prodop_n(X_1, \,\dots\, , X_n) \to (X_1, \,\dots\, , X_n)$ 
for each 
$X_1, \,\dots\, , X_n \in S \:\: (n \in \Nat)$, then for every sequence of 
objects $\Gamma$
one obtains a mapping 
$\psi_\Gamma : 
\prod_{i=1}^n(\mCatOf{\clone})(\Gamma; X_i) 
	\to 
	(\mCatOf{\clone}){\big(\Gamma; \prodop_n(X_1, \,\dots\, , X_n)\big)}$ 
such that the following equations hold for every multimap 
$h : \Gamma \to \prodop_n(X_1, \,\dots\, , X_n)$ and 
sequence of multimaps
$(f_i : \Gamma \to X_i)_{i=1, \,\dots\, , n}$:
\begin{equation} \label{eq:equations-from-universal-arrow}
\psi_\Gamma(\cslr{\pi_1}{h}, \,\dots\, , \cslr{\pi_n}{h}) = h \quad \text{ and 
} 
\quad 
\cslr{\pi_i}{\psi_\Gamma(f_1, \,\dots\, , f_n)} = f_i
\:\:\: (i=1, \,\dots\, , n)
\end{equation} 
Thus, $\psi_\Gamma(-, \dots, =)$ provides a `tupling' operation.
This is substantiated by the next lemma.

\begin{mydefn} \label{def:clone-1-cell-invetibility}
Let $(S, \clone)$ be a clone. A multimap $f : X_1, \,\dots\, , X_n \to Y$ in 
$\clone$
is \Def{invertible} or an \Def{iso} if there exists a family of unary multimaps
$(g_i : Y \to X_i)_{i=1, \,\dots\, , n}$ in $\clone$
such that $\cs{f}{g_1, \,\dots\, , g_n} = \id_Y$ and
$\cs{g_i}{f} = \p{i}{\ind{X}}$ for $i = 1, \,\dots\, , n$. 
If there exists an invertible multimap 
$f :  X_1, \,\dots\, , X_n \to Y$ we say $ X_1, \,\dots\, , X_n$ and 
$Y$ are \Def{isomorphic}, and write 
$ X_1, \,\dots\, , X_n \iso Y$.
\end{mydefn}

A small adaptation of the usual categorical proof shows that inverses in a 
clone are unique, in the sense that if $f$ has inverses $(g_1, \dots, g_n)$
and $(g_1', \dots, g_n')$ then $g_i = g_i'$ for $i=1,\dots,n$. 

\begin{mylemma} \label{lem:tupling-from-cartesian-clone-gives-identity}
Let $(S, \clone)$ be a cartesian clone. Then, where the $n$-ary product 
of $X_1, \dots, X_n \in S \:\: (n \in \Nat)$ is
witnessed by the universal arrow
$(\pi_1, \dots, \pi_n) : \prodop_n(X_1, \dots, X_n) \to (X_1, \dots, X_n)$,
\[
\cslr
	{\psi_{\ind{X}}(\p{1}{\ind{X}}, \,\dots\, , \p{n}{\ind{X}})}
	{\pi_1, \,\dots\, , \pi_n} 
= 
\id_{\prod_n(X_1, \dots, X_n)}
\]
Hence $X_1, \dots, X_n \iso \prodop_n(X_1, \dots, X_n)$.
\begin{proof}
For the first part one uses the two equations 
of~(\ref{eq:equations-from-universal-arrow}): 
\begin{align*}
\cslr
	{\psi_{\ind{X}}(\p{1}{\ind{X}}, \,\dots\, , \p{n}{\ind{X}})}
	{\pi_1, \,\dots\, , \pi_n} 
&=
\psi_{(\prod_n \ind{X})}
	\left( 
	\cslr{\ind{\pi}}
		{\cslr
			{\psi_{\ind{X}}(\p{1}{\ind{X}}, \,\dots\, , \p{n}{\ind{X}})}
			{\pi_1, \,\dots\, , \pi_n}}\right) 
	&\text{ by } 
	(\ref{eq:equations-from-universal-arrow}) \\
&= \psi_{(\prod_n \ind{X})}
	\left( 
		\csthree{\ind{\pi}}
				{\psi_{\ind{X}}(\p{1}{\ind{X}}, \,\dots\, , \p{n}{\ind{X}})}
				{\pi_1, \,\dots\, , \pi_n} 
	\right) \\
&= \psi_{(\prod_n \ind{X})}
	\left( 
		\cslr	{\p{\bullet}{\ind{X}}}
			{\pi_1, \,\dots\, , \pi_n} 
	\right)
	&\text{ by  } 
	(\ref{eq:equations-from-universal-arrow}) \\
&= \psi_{(\prod_n \ind{X})}\left( \pi_1, \,\dots\, , \pi_n \right) \\
&= \psi_{(\prod_n \ind{X})}\left( 
	\cslr{\pi_1}{\id_{(\prod_n \ind{X})}}, \,\dots\, , 
	\cslr{\pi_n}{\id_{(\prod_n \ind{X})}} \right) \\
&= \id_{(\prod_n \ind{X})} 
	&\text{ by  } 
	(\ref{eq:equations-from-universal-arrow})
\end{align*}
Then $(\pi_i : \prodop_n(X_1, \dots, X_n) \to X_i)_{i=1, \dots, n}$ and
$\psi_{\ind{X}}(\p{1}{\ind{X}}, \dots, \p{n}{\ind{X}})$ 
form the claimed isomorphism.
\end{proof}
\end{mylemma}

We now turn to examinining how representability (thought of as `tupling')
gives rise to `projections'. The next lemma is the key construction.

\begin{mylemma} \label{lem:product-structure-from-representable-arrow}
For any representable clone $(S, \clone)$ and 
$X_1, \,\dots\, , X_n \in S \:\: (n \in \Nat)$ 
there exist multimaps
$\pi_i : \tensor_n(X_1, \,\dots\, , X_n) \to X_i \:\: (i=1, \,\dots\, , n)$
such that
\begin{center}
$\pi_i \circ \rho_{\ind{X}} = \p{i}{\ind{X}}$
\quad and \quad 
$\cslr{\rho_{\ind{X}}}{\pi_1, \,\dots\, , \pi_n} = \id_{\prod \ind{X}}$ 
\end{center}
where $\rho_{\ind{X}}$ is the representable arrow.
\begin{proof}
By representability, we may define 
$\pi_i := \altTrans{(\p{i}{\ind{X}})}$. The first 
claim then holds by 
assumption. For the second, observing that 
$\altTrans{(\rho_{\ind{X}})} = \id_{\prod \ind{X}}$, it suffices to show that 
$\csthree{\rho_{\ind{X}}}{\pi_1, \,\dots\, , \pi_n}{\rho_{\ind{X}}} = 
\rho_{\ind{X}}$. But this is straightforward: 
\[
\csthree{\rho_{\ind{X}}}{\pi_1, \,\dots\, , \pi_n}{\rho_{\ind{X}}} = 
\cslr{\rho_{\ind{X}}}{\cslr{\ind{\pi}}{\rho_{\ind{X}}}} = 
\cslr{\rho_{\ind{X}}}{\p{1}{}, \,\dots\, , \p{n}{}} = \rho_{\ind{X}}
\]
\end{proof}
\end{mylemma} 

Another important consequence of 
Lemma~\ref{lem:product-structure-from-representable-arrow} is that, in the case of clones, representable arrows are always closed under composition.

\begin{mylemma} \label{lem:multicat-from-clone-universal-arrows-compose}
For any clone $(S, \clone)$, the multicategory $\mCatOf{\clone}$ is 
representable if and only if for every 
$X_1, \,\dots\, , X_n \in S \:\: (n \in \Nat)$ there exists a 
chosen object $\tensor_n(X_1, \,\dots\, , X_n)$ and a representable multimap 
$\rho_{\ind{X}} : X_1, \,\dots\, , X_n \to \tensor_n(X_1, \,\dots\, , X_n)$. 
\begin{proof}
It suffices to show that, for any clone $(S, \clone)$, the representable 
multimaps in $\mCatOf{\clone}$ are closed under composition. Suppose given 
representable  multimaps
\begin{gather*}
\rho_{\ind{X}} : X_1, \,\dots\, , X_n \to \tensor_n(X_1, \,\dots\, , X_n) \\ 
\rho_{\ind{Y}} : Y_1, \,\dots\, , Y_m \to \tensor_m(Y_1, \,\dots\, , Y_m) \\
\rho_{(\tensor \ind{X}, \tensor \ind{Y})} : \tensor_n \ind{X}, \tensor_m 
\ind{Y} 
\to \tensor_2(\tensor_n \ind{X}, \tensor_m \ind{Y})
\end{gather*}
We want to show that the composite 
$\ms
	{\rho_{(\tensor \ind{X}, \tensor \ind{Y})}}
	{\rho_{\ind{X}}, \rho_{\ind{Y}}}$ 
in $\mCatOf{\clone}$, which is the 
composite 
$\cslr
	{\rho_{(\tensor \ind{X}, \tensor \ind{Y})}}
	{\rho_{\ind{X}} \clonetimes \rho_{\ind{Y}}} 
= 
\cslr
	{\rho_{(\tensor \ind{X}, \tensor \ind{Y})}}
	{\cslr{\rho_{\ind{X}}}
		{\p{1}{}, \,\dots\, , \p{n}{}}, 
		\cslr{\rho_{\ind{Y}}}{\p{n+1}{}, \,\dots\, , \p{n+m}{}}
	}$ in $\clone$, is 
representable.

\newpage
By Lemma~\ref{lem:product-structure-from-representable-arrow}, we may define 
multimaps
\begin{align*}
\pi_i^X : \tensor_n (X_1, \,\dots\, , X_n) &\to X_i \qquad \text{ for } i=1, 
\,\dots\, , 
n 
\\
\pi_j^Y : \tensor_m(Y_1, \,\dots\, , Y_m) &\to Y_j \qquad \text{ for } j=1, 
\,\dots\, ,m 
\\
\pi_1^{X,Y} : \tensor_2(\tensor_n \ind{X}, \tensor_m \ind{Y}) 
	&\to \tensor_n \ind{X} \\
\pi_2^{X,Y} : \tensor_2(\tensor_n \ind{X}, \tensor_m \ind{Y}) 
	&\to \tensor_m \ind{Y}
\end{align*}
Then, setting
\begin{align*}
Z_i := 	\begin{cases}
			X_i \qquad &\text{ for } i= 1, \,\dots\, , n \\
			Y_{i-n} \qquad &\text { for } i = n+1, \,\dots\, , n+m
		\end{cases}
\end{align*}
we define
$\overline{\pi}_i : 
	\tensor_2 (\tensor_n \ind{X}, \tensor_m \ind {Y}) \to Z_i$
by iterated applications of $\pi_i$:
\begin{equation} \label{eq:def-of-ri-maps}
\overline{\pi}_i := 	\begin{cases}
			\cslr{\pi_i^X}{\pi_1^{X,Y}} \qquad \text{ for } 1 \leq i \leq n \\
			\cslr{\pi_{i-n}^Y}{\pi_2^{X,Y}} \qquad \text{ for } n + 1 \leq i 
			\leq 
			n + m
		\end{cases}
\end{equation}
The rest of the proof revolves around proving the following two equalities in 
$\clone$:
\begin{equation} \label{eq:adjoint-equivalence-following-rep-arrow}
\begin{tikzcd}
X_1, \,\dots\, , X_n, Y_1, \,\dots\, , Y_m \arrow[]{r}{\p{i}{}} 
\arrow[swap]{d}{[\rho_{\ind{X}} \clonetimes \rho_{\ind{Y}}]} &
Z_i \\
\tensor_n \ind{X}, \tensor_m \ind{Y} 
\arrow[swap]{r}{\rho_{(\tensor \ind{X}, \tensor \ind{Y})}} &
\tensor_2(\tensor_n \ind{X}, \tensor_m \ind{Y}) 
\arrow[swap]{u}{\overline{\pi}_i}
\end{tikzcd}
\end{equation}
\begin{equation} \label{eq:adjoint-equivalence-following-ris}
\begin{tikzcd}
\tensor_2(\tensor_n \ind{X}, \tensor_m \ind{Y}) 
\arrow[swap]{d}{[\overline{\pi}_1, \,\dots\, , 
\overline{\pi}_{n+m}]} \arrow[equals]{r} &
\tensor_2(\tensor_n \ind{X}, \tensor_m \ind{Y}) \\
X_1, \,\dots\, , X_n, Y_1, \,\dots\, , Y_m 
\arrow[swap]{r}[yshift=-2mm]{[\rho_{\ind{X}} \clonetimes \rho_{\ind{Y}}]} &
\tensor_n \ind{X}, \tensor_m \ind{Y} 
\arrow[swap]{u}{\rho_{(\tensor \ind{X}, \tensor \ind{Y})}}
\end{tikzcd}
\end{equation}
Indeed, if these two diagrams commute, then for any 
$g : X_1, \,\dots\, , X_n, Y_1, \,\dots\, , Y_m \to A$ one may define 
$\altTrans{g} : \tensor_2 (\tensor_n \ind{X}, \tensor_m \ind{Y}) \to A$ to be 
the composite $\cslr{g}{\overline{\pi}_1, \,\dots\, , \overline{\pi}_{n+m}}$. 
It then follows that
that $\altTrans{(-)}$ is 
the inverse to precomposing with $\overline{\rho} := \cslr{\rho_{(\tensor \ind 
{X}, 
\tensor \ind {Y}})}{\rho_{\ind{X}} \clonetimes \rho_{\ind{Y}}}$:
\[
\csthree{g}{\overline{\pi}_1, \,\dots\, , 
\overline{\pi}_{n+m}}{\overline{\rho}} = 
\cslr{g}{\cslr{\overline{\pi}_1}{\overline{\rho}}, \,\dots\, , 
\cslr{\overline{\pi}_{n+m}}{\overline{\rho}}} 
\overset{(\ref{eq:adjoint-equivalence-following-rep-arrow})}{=} 
\cslr{g}{\p{1}{}, 
\dots, \p{n+m}{}} = g 
\]
while, for any $h : \tensor_2(\tensor_n \ind{X}, \tensor_m \ind{Y}) \to A$, 
\[
\csthree{h}{\overline{\rho}}{\overline{\pi}_1, \,\dots\, , 
\overline{\pi}_{n+m}} 
\overset{(\ref{eq:adjoint-equivalence-following-ris})}{=} 
\cslr{h}{\p{1}{\tensor 
(\tensor \ind {X}, \tensor \ind{Y})}} = h 
\]
as required.

It therefore remains to establish the commutativity of the two diagrams above.
We compute~(\ref{eq:adjoint-equivalence-following-rep-arrow}) directly. For 
example, for $1 \leq i \leq n$, unfolding the universal property of each of the 
projections gives
\begin{align*}
\csthree{\overline{\pi}_i}
	{\rho_{(\tensor \ind{X}, \tensor \ind{Y})}}
	{\rho_{\ind{X}} \clonetimes \rho_{\ind{Y}}} 
&= 
\csthree{\cslr{\pi_i^X}{\pi_1^{X,Y}}}
	{\rho_{(\tensor \ind{X}, \tensor \ind{Y})}}
	{\rho_{\ind{X}} \clonetimes \rho_{\ind{Y}}} \\
&=
\csthree
	{\pi_i^X}
	{\cslr{\pi_1^{X,Y}}
		{\rho_{(\tensor \ind{X}, \tensor\ind{Y})}}}
	{\rho_{\ind{X}} \clonetimes \rho_{\ind{Y}}} \\
&= \csthree
		{\pi_i^X}
		{\p{1}{(\tensor\ind{X}, \tensor\ind{Y})}}
		{\rho_{\ind{X}} \clonetimes \rho_{\ind{Y}}} \\
&= \cslr
	{\pi_i^X}
	{\cslr
		{\p{1}{(\tensor\ind{X}, \tensor\ind{Y})}}
		{\rho_{\ind{X}} \clonetimes \rho_{\ind{Y}}}} \\
&= \cslr
	{\pi_i^X}
	{\cslr
		{\rho_{\ind{X}}}
		{\p{1}{}, \,\dots\, , \p{n}{}}} 	 \\
&= \cslr{\cslr{\pi_i^X}{\rho_{\ind{X}}}}{\p{1}{}, \,\dots\, , \p{n}{}} \\
&= \cslr{\p{i}{}}{\p{1}{}, \,\dots\, , \p{n}{}} \\
&= \p{i}{}
\end{align*}
as required. For~(\ref{eq:adjoint-equivalence-following-ris}), 
Lemma~\ref{lem:product-structure-from-representable-arrow} entails that
\[
\cslr{\rho_{\ind{X}}}{\overline{\pi}_1, \,\dots\, , \overline{\pi}_{n}} = 
\cslr{\rho_{\ind{X}}}{\cslr{\pi_1^X}{\pi_1^{X,Y}}, \,\dots\, , 
\cslr{\pi_n^X}{\pi_1^{X,Y}}} = 
\csthree{\rho_{\ind{X}}}{\ind{\pi}^X}{\pi_1^{X,Y}} 
= \pi_1^{X,Y}
\]
and hence that
\begin{align*}
\csthree
	{\rho_{(\tensor \ind{X}, \tensor \ind{Y})}}
	{\cslr{\rho_{\ind{X}}}{\p{\bullet}{}},\cslr{\rho_{\ind{Y}}}{\p{\bullet}{}}}
	{\ind{\overline{\pi}}}
&= 
\cslr
	{\rho_{(\tensor \ind{X}, \tensor \ind{Y})}}
	{\cslr
		{\rho_{\ind{X}}}{\ind{\overline{\pi}}},
	\cslr
		{\rho_{\ind{Y}}}{\ind{\overline{\pi}}}} \\
	&= \cslr{\rho_{(\tensor \ind{X}, \tensor \ind{Y})}}{\pi_1^{X,Y}, 
	\pi_2^{X,Y}} \\
	&= \id_{\tensor(\tensor \ind{X}, \tensor \ind{Y})}
\end{align*}
as required.
\end{proof}
\end{mylemma} 

We now make precise the sense in which the inverse to 
precomposing with a representable arrow provides a tupling operation. The product structure on a representable clone is, as expected, 
given by the
1-cells constructed in 
Lemma~\ref{lem:product-structure-from-representable-arrow}.

\begin{mylemma} \label{lem:clone-representable-iff-cartesian}
For any clone $(S, \clone)$, the following are equivalent:
\begin{enumerate}
\item $(S, \clone)$ is representable, 
\item $(S, \clone)$ is cartesian.
\end{enumerate}
\begin{proof}
\subproof{$\To$} We prove the forward direction first. Suppose $\rho_{\ind{X}} 
: X_1, \,\dots\, , X_n \to \tensor_n(X_1, \,\dots\, , X_n)$ is representable; 
we claim 
the required universal arrow is given by the sequence of multimaps 
$(\pi_1, \,\dots\, , \pi_n) : 
\Delta\tensor_n(X_1, \,\dots\, , X_n) \to (X_1, \,\dots\, , X_n)$
defined in Lemma~\ref{lem:product-structure-from-representable-arrow}. 
To this end, let $(f_i : \Gamma \to X_i)_{i=1, \,\dots\, , n}$ in $\clone$. We 
set 
$\psi_{\Gamma}(f_1, \,\dots\, , f_n) : \Gamma 
\to \tensor_n(X_1, \,\dots\, , X_n)$ to be the composite 
$\cslr{{\rho_{\ind{X}}}}{{f_1, \,\dots\, , f_n}}$.  
By Lemma~\ref{lem:product-structure-from-representable-arrow}, 
\[
\pi_i \circ \big( \psi_{\Gamma}(f_1, \,\dots\, , f_n) \big) = 
\cslr{\pi_i}{\cslr{\rho_{\ind{X}}}{f_1, \,\dots\, , f_n}} = 
\cslr{\p{i}{\ind{X}}}{f_1, \,\dots\, , f_n} = f_i
\]
for $i=1, \,\dots\, , n$, so it remains to show that
$\psi_{\Gamma}(\cslr{\pi_1}{h}, \,\dots\, , \cslr{\pi_n}{h}) = h$
for every
$h : {\Gamma \to \tensor_n(X_1, \,\dots\, , X_n)}$. Applying the lemma again, 
\[
\psi_{\Gamma}(\cslr{\pi_1}{h}, \,\dots\, , \cslr{\pi_n}{h}) = 
\cslr{\rho_{\ind{X}}}{\cslr{\pi_1}{h}, \,\dots\, , \cslr{\pi_n}{h}} =
\csthree{\rho_{\ind{X}}}{\pi_1, \,\dots\, , \pi_n}{h} = h
\] 
as required.

\subproof{$\Leftarrow$} 
We claim that 
$\rho_{\ind{X}} := 
	\psi_{\ind{X}}(\p{1}{\ind{X}}, \,\dots\, , \p{n}{\ind{X}}) : 
	X_1, \,\dots\, , X_n \to \prodop_n(X_1, \,\dots\, , X_n)$ is 
representable. 

To this end, suppose $h : X_1, \,\dots\, , X_n \to A$. We define $\trans{h} : 
\prodop_n(X_1, \,\dots\, , X_n) \to A$ to be the composite $\cslr{h}{\pi_1, 
\,\dots\, , 
\pi_n}$. Then 
\begin{align*}
\cslr{\trans{h}}{\rho_{\ind{X}}} &= \csthree{h}{\pi_1, \,\dots\, , 
\pi_n}{\psi_\Gamma(\p{1}{\ind{X}}, \,\dots\, , \p{n}{\ind{X}})} \\ 
	&= \cslr{h}{\cslr{\ind{\pi}}{\psi_\Gamma(\p{1}{\ind{X}}, \,\dots\, , 
	\p{n}{\ind{X}})}} \\ &= \cslr{h}{\p{1}{\ind{X}}, \,\dots\, , 
	\p{n}{\ind{X}}} \\ 
	&= h
\end{align*}
so the existence part of the claim holds. It remains to check the 
equality $\trans{(\cslr{f}{\rho_{\ind{X}}})} = f$ 
for an arbitrary $f : \prodop_n(X_1, \,\dots\, , X_n) \to A$. Examining the 
equality
\[
\trans{(\cslr{f}{\rho_{\ind{X}}})} = 
\csthree{f}{\rho_{\ind{X}}}{\pi_1, \,\dots\, , \pi_n} = 
\cslr{f}
	{\cslr{\psi_{\ind{X}}(\p{1}{\ind{X}}, \,\dots\, , \p{n}{\ind{X}})}
		{\pi_1, \,\dots\, , \pi_n}}
\] 
it suffices to show that 
$\cslr
	{\psi_{\ind{X}}(\p{1}{\ind{X}}, \,\dots\, , \p{n}{\ind{X}})}
	{\pi_1, \,\dots\, , \pi_n}$ 
is the identity. This is Lemma~\ref{lem:tupling-from-cartesian-clone-gives-identity}.
\end{proof}
\end{mylemma}

We summarise the last two results in the following theorem. The final case 
is Lemma~\ref{lem:multicat-universal-arrow-to-natural-isomorphism}.

\enlargethispage*{4\baselineskip}
\begin{prooflessthm} \label{thm:representable-clone-cartesian-equivalence}
For any clone $(S, \clone)$, the following are equivalent:
\begin{enumerate}
\item \label{c:representable} $(S, \clone)$ is representable, 
\item For every $X_1, \,\dots\, , X_n \in S \:\: (n \in \Nat)$ there exists a 
choice 
of object $\prodop_n(X_1, \,\dots\, , X_n) \in S$ together with a
representable multimap 
$\rho_{\ind{X}} : X_1, \,\dots\, , X_n \to \prodop_n(X_1, \,\dots\, , X_n)$,
\item \label{c:fp-clone} $(S, \clone)$ is cartesian, 
\item \label{c:nat-isoms-multicat} For any 
$X_1, \,\dots\, , X_n \in S \:\: (n \in \Nat)$ 
there exists a chosen object $\prodop_n(X_1, \,\dots\, , X_n) \in S$ and an 
isomorphism 
$(\mCatOf{\clone})\big(\Gamma; \prodop_n(X_1, \,\dots\, , X_n)\big) \iso 
\prod_{i=1}^n(\mCatOf{\clone})(\Gamma; X_i)$, 
multinatural in the sense that for 
any $f : \Gamma \to A$ the following diagram commutes:
\begin{td}
(\mCatOf{\clone})\big(\Gamma; \prodop_n(X_1, \,\dots\, , X_n)\big) 
\arrow{r}{\iso} & 
\prod_{i=1}^n(\mCatOf{\clone})(\Gamma; X_i) \\
(\mCatOf{\clone})\big(A; \prodop_n(X_1, \,\dots\, , X_n)\big) 
\arrow[swap]{r}{\iso} 
\arrow{u}{\ms{(-)}{f}} & \prod_{i=1}^n(\mCatOf{\clone})(A; X_i) 
\arrow[swap]{u}{\ms{(-)}{f}}
\end{td}
\end{enumerate}
\end{prooflessthm}

In the case of clones, therefore, the two approaches to 
defining 
product structure---Hermida's representability or Lambek's natural 
isomorphisms---actually coincide.
We tie this back to Hermida's equivalence between monoidal categories and 
representable multicategories with the following observation.

\begin{mycor} \label{cor:representable-to-cartesian-clone}
For any representable clone $(S, \clone)$, the monoidal structure on the 
category $\overline{\mCatOf{\clone}}$ associated to $\mCatOf{\clone}$ is 
cartesian. 
\begin{proof}
The required natural isomorphism follows by restricting the 
isomorphism~(\ref{eq:fp-multicategory}) to linear multimaps. Explicitly, 
the $n$-ary product of $X_1, \,\dots\, , X_n$ is $\prodop_n(X_1, \,\dots\, , 
X_n)$, and the 
projections are \mbox{$\pi_i : \prodop_n(X_1, \,\dots\, , X_n) \to X_i$}. The 
$n$-ary tupling of maps $(f_i : A \to X_i)_{i=1, \,\dots\, , n}$ is given via 
the representable arrow $\rho_{\ind{X}}$ for $X_1, \,\dots\, , X_n$, as 
$\cslr{\rho_{\ind{X}}}{f_1, \,\dots\, , f_n}$.
\end{proof}
\end{mycor}

%
%
%

It is reasonable to suggest that one 
could refine Hermida's 2-equivalence between monoidal categories and 
representable multicategories to a 2-equivalence between cartesian categories 
and representable clones; the calculations required would take us beyond the 
theory we shall actually need, so we do not pursue the point here. Instead we 
turn to the syntactic implications of 
the theory just developed.


\subsection{From cartesian clones to type theory}
\label{sec:products-in-stlc}

%

\paragraph*{From cartesian clones to cartesian categories.}
In Chapter~\ref{chap:biclone-lang} we saw that the free category on a graph 
could be constructed by restricting the free clone on that graph to its unary 
operations. This fact extends to cartesian clones and cartesian categories.
To show this, we need to enrich our notion of signature to include 
product structure. The definition was already hinted at in 
Example~\ref{ex:stlc-a-clone}.

\begin{mydefn}
A \Def{$\stlcTimes$-signature} $\sig = (\baseTypes, \graph)$ consists of 
\begin{enumerate}
\item A set of base types $\baseTypes$, 
\item A multigraph $\graph$ with nodes
generated by the grammar
\begin{equation} \label{eq:stlc-times-sig-grammar}
A_1, \,\dots\, , A_n ::= 
	B \st 
	\prodop_n(A_1, \,\dots\, , A_n) \qquad (B \in \baseTypes, n \in \Nat)
\end{equation}
\end{enumerate}
If the graph $\graph$ is a 2-graph we call the signature \Def{unary}. A 
\Def{homomorphism} of $\stlcTimes$-signatures
$h : {\sig \to \sig'}$ 
is a multigraph homomorphism 
$h : {\graph \to \graph'}$
which respects the product structure in the sense that
$h{\left( \prodop_n(A_1, \,\dots\, , A_n) \right)} 
= \prodop_n\left( hA_1, \,\dots\, , hA_n \right)$. 
We denote the category of $\stlcTimes$-signatures and their homomorphisms by
$\stlcTimesSigCat$, and the full subcategory of unary $\stlcTimes$-signatures by
$\stlcTimesUnSigCat$. 
\end{mydefn}

\begin{mynotation} \label{not:stlc-times-alltypes-not}
For any $\stlcTimes$-signature $\sig = (\baseTypes, \graph)$ we write 
$\allTypes\baseTypes$ for the set generated from $\baseTypes$ by the 
grammar~(\ref{eq:stlc-times-sig-grammar}) (equivalently, the set $\nodes\graph$ 
of nodes in $\graph$).
In particular, when the signature is just a set (\ie~the graph $\graph$ has no edges) we denote the 
signature $\sig = (\baseTypes, \sig)$ simply by $\allTypes\baseTypes$. 
\end{mynotation}

The following lemma mirrors the situation for graphs and 2-multigraphs.

\begin{mylemma} \label{lem:stlcTimesSig-reflective-subcat}
The embedding $\inc : \stlcTimesUnSigCat \hookrightarrow \stlcTimesSigCat$ 
has a right adjoint. 
\begin{proof}
Define the functor 
$\widetilde{\lin} : \stlcTimesSigCat \to \stlcTimesUnSigCat$ 
to be the 
restriction of the corresponding functor 
$\lin : \MultiGraph \to \Graph$. Thus, 
$\widetilde{\lin}$ restricts a signature $(\baseTypes, \graph)$ to the 
signature with base types 
$\baseTypes$ and multigraph $\lin\graph$ containing only edges of the form 
$X \to Y$. This is a right adjoint to the given 
inclusion because $\lin$ is right adjoint to the inclusion 
$\Graph \hookrightarrow \MultiGraph$.
\end{proof}
\end{mylemma}

Every cartesian category 
$\fpBicat{\catC}$
has an underlying unary 
$\stlcTimes$-signature with edges $X \to Y$ given by morphisms $X \to Y$ in 
$\catC$ (\cf~\cite[Theorem~4.9.2]{Crole1994}).  Similarly, every cartesian 
clone 
$\cartClone{S}{\clone}$ has an underlying 
$\stlcTimes$-signature with the edges given by
multimaps. We wish to construct the free cartesian 
clone over such a signature. 
Theorem~\ref{thm:representable-clone-cartesian-equivalence} 
guarantees that it is sufficient to add a representable arrow
$A_1, \,\dots\, , A_n \to \prodop_n(A_1, \,\dots\, , A_n)$
for every sequence of types
$A_1, \,\dots\, , A_n \:\: (n \in \Nat)$.
For the construction we follow the forward direction of the proof of 
Lemma~\ref{lem:clone-representable-iff-cartesian}.

\begin{myconstr} \label{constr:free-cartesian-clone}
For any $\stlcTimes$-signature $\sig = (\baseTypes, \graph)$, define a clone 
$(\nodes{\graph}, \freeCartClone{\sig})$ with sorts generated from 
$\baseTypes$ by the rules
\[
A_1, \,\dots\, , A_n ::= 
	B \st 
	\prodop_n(A_1, \,\dots\, , A_n) \qquad (B \in \baseTypes, n \in \Nat)
\] 
as the 
following deductive system:
\begin{center}
\unaryRule
		{c \in \graph(A_1, \,\dots\, , A_n; B)}
		{c \in \freeCartClone{\sig}(A_1, \,\dots\, , A_n; B)}
		{} \\

\unaryRule
		{\faketext}
		{\p{i}{A_1, \,\dots\, , A_n} \in 
			\freeCartClone{\sig}(A_1, \,\dots\, , A_n; A_i)}
		{$(1 \leq i \leq n)$} \vspace{0.5\treeskip} \\

\binaryRule
	{f \in \freeCartClone{\sig}(A_1, \,\dots\, , A_n; B)}
	{\big(g_i \in \freeCartClone{\sig}(\ind{X}; A_i)\big)_{i=1,\dots,n}}
	{\cslr{f}{g_1, \,\dots\, , g_n} \in \freeCartClone{\sig}(\ind{X}; B)}{}

\unaryRule
	{\faketext}
	{\pairName_{\ind{A}} \in 
		\freeCartClone{\sig}\left(A_1, \,\dots\, , A_n; 
			\prodop_n(A_1, \,\dots\, , A_n) \right)}
	{}

\unaryRule 
	{\faketext}
	{\mathsf{proj}^{(i)}_{\ind{A}} \in \freeCartClone{\sig}
		\left(\prodop_n(A_1, \,\dots\, , A_n); A_i\right)}
	{$(1 \leq i \leq n)$}
\end{center}
subject to an equational theory requiring
\begin{itemize}
\item The clone laws hold with projection $\p{i}{\ind{A}}$ and substitution 
		$\cslr{f}{g_1, \,\dots\, , g_n}$,  
\item 
	$\cslr{\mathsf{proj}^{(i)}_{\ind{A}}}{\pairName_{\ind{A}}}  
		\equiv \p{i}{\ind{A}}$ for $i=1, \,\dots\, , n$,
\item 
	$\cslr{\pairName_{\ind{A}}}
		{\mathsf{proj}^{(n)}_{\ind{A}}, \,\dots\, , 
		\mathsf{proj}^{(n)}_{\ind{A}}}
		\equiv \p{1}{(\prod_n \ind{A})}$.
\qedhere
\end{itemize}
\end{myconstr}

The clone $\freeCartClone{\sig}$ is cartesian because it is representable. 
Indeed, for any $A_1, \,\dots\, , A_n, B \in \nodes\graph$, the equational laws 
ensure that the map 
$(-) \circ \pairName_{\ind{A}}$ 
has inverse
$\cslr{(-)}{\mathsf{proj}^{(n)}_{\ind{A}}, \,\dots\, , 
\mathsf{proj}^{(n)}_{\ind{A}}}$,
giving rise to the required natural isomorphism
$\freeCartClone{\sig}{\left( \prodop_n(A_1, \,\dots\, , A_n); B\right)}
	\iso \freeCartClone{\sig}(A_1, \,\dots\, , A_n; B)$.

In order to state that this construction yields the free cartesian clone, we 
need to define a notion of product-preserving clone homomorphism.  This is the 
clone-theoretic 
translation of
Definition~\ref{def:morphism-of-universal-arrows}, requiring that the 
universal arrow is preserved.  

\begin{mydefn}
A \Def{cartesian clone homomorphism} 
$h : \cartClone{S}{\clone} \to \cartClone{T}{\altClone}$ 
is a clone homomorphism $h : (S, \clone) \to (T, \altClone)$ such that
the canonical map
$
\psi_{\prod \ind{A}}(h\pi_1, \,\dots\, , h\pi_n) : 
h{\left( \prodop_n(A_1, \,\dots\, , A_n) \right)} \to 
\prodop_n \left( hA_1, \,\dots\, , A_n \right)
$ 
is invertible
for every 
$A_1, \,\dots\, , A_n \in S \:\: (n \in \Nat)$.

We call $h$ \Def{strict} if 
\begin{align*}
h{\left(\prodop_n(A_1, \,\dots\, , A_n)\right)} &=
	\prodop_n\left(hA_1, \,\dots\, , hA_n\right) \\
h(\pi_i^{\ind{A}}) &= 
		{\left( \prodop_n(hA_1, \,\dots\, , hA_n) \xra{\pi_i} 
			h(A_i) \right)} \qquad (i=1, \,\dots\, , n)
\end{align*}
for every $A_1, \,\dots\, , A_n \in S \:\: (n \in \Nat)$.
\end{mydefn}

\begin{mylemma} \label{lem:free-cartesian-clone}
For any cartesian clone $\cartClone{T}{\altClone}$, $\stlcTimes$-signature 
$\sig$ 
and
$\stlcTimes$-signature homomorphism $h : \sig \to \altClone$, there exists a 
unique strict cartesian
clone homomorphism $\ext{h} : \freeCartClone{\sig} \to \altClone$ such that
$\ext{h} \circ \inc = h$, for 
$\inc : \sig \hookrightarrow \freeCartClone{\sig}$ 
the inclusion.
\begin{proof}
We define $\ext{h}$ by induction. The requirement that 
$\ext{h} \circ \inc = h$
completely determines the action of $\ext{h}$ on objects, and also entails that
$\ext{h}(c) = h(c)$ on constants.
On multimaps, the clone homomorphism axioms require that we set
\begin{align*}
\ext{h}(\p{i}{\ind{A}}) &:= \p{i}{\ext{h}{\ind{A}}} \\
\ext{h}(\cslr{f}{g_1, \,\dots\, , g_n}) &:= 
	\cslr{\ext{h}(f)}{\ext{h}(g_1), \,\dots\, , \ext{h}(g_n)}
\end{align*}
The definition on 
$\mathsf{proj}^{(i)}$ is determined by the hypothesis. Finally, on $\pairName$ 
we
set
$\ext{h}\left(\pairName_{\ind{A}}\right) := \rho_{\ext{h}(\ind{A})}$, 
so that $\ext{h}$ sends $\pairName_{\ind{A}}$ to the representable arrow on 
$A_1, \,\dots\, , A_n$ (which exists by 
Lemma~\ref{lem:clone-representable-iff-cartesian}).
For uniqueness, it remains to show that the action of $\ext{h}$ on 
$\pairName$ is determined by the hypotheses. For this, consider
\newpage
\begin{align*}
\rho_{(\ext{h}\ind{A})} 
&= 
\cslr{\rho_{(\ext{h}\ind{A})}}
	{\p{1}{\ext{h}(\ind{A})}, \,\dots\, , \p{n}{\ext{h}(\ind{A})}} \\
&= 
\cslr{\rho_{(\ext{h}\ind{A})}}
	{\ext{h}(\p{1}{\ind{A}}), \,\dots\, , \ext{h}(\p{n}{\ind{A}})} \\
&= 
\cslr{\rho_{(\ext{h}\ind{A})}}
	{\ext{h}(\cslr{\mathsf{proj}^{(1)}}{\rho_{\ind{A}}}), \,\dots\, , 
	\ext{h}(\cslr{\mathsf{proj}^{(n)}}{\rho_{\ind{A}}})} \qquad 
		&\text{by Lemma~\ref{lem:product-structure-from-representable-arrow}} \\
&= 
\cslr{\rho_{(\ext{h}\ind{A})}}
	{\cslr{\ext{h}(\mathsf{proj}^{(1)})}{\ext{h}(\rho_{\ind{A}})}, \,\dots\, , 
	\cslr{\ext{h}(\mathsf{proj}^{(n)})}{\ext{h}(\rho_{\ind{A}})}} \\
&=
\cslr{\rho_{(\ext{h}\ind{A})}}
	{\cslr{\pi_1}{\ext{h}(\rho_{\ind{A}})}, \,\dots\, , 
	\cslr{\pi_n}{\ext{h}(\rho_{\ind{A}})}}  
	&\text{by~cartesian}\\
&=
\csthree
	{\rho_{(\ext{h}\ind{A})}}
	{\pi_1, \,\dots\, , \pi_n}
	{\ext{h}(\rho_{\ind{A}})} 
	&\text{by Lemma~\ref{lem:product-structure-from-representable-arrow}} \\
&= 
\cslr{\p{1}{(\prod_n \ind{A})}}
	{\ext{h}(\rho_{\ind{A}})} \\
&= 
\ext{h}(\rho_{\ind{A}})
\end{align*}
Hence, the action of any clone homomorphism satisfying the two hypotheses is 
completely determined, and $\ext{h}$ is unique.
\end{proof}
\end{mylemma}

%

The term calculus corresponding to the deductive system of
Construction~\ref{constr:free-cartesian-clone} 
is specified by the following rules:
\begin{enumerate}
\item For every sequence of types $A_1, \,\dots\, , A_n \:\: (n \in \Nat)$, 
there 
exists a type $\prodop_n(A_1, \,\dots\, , A_n)$, 
\item For every context $x_1 : A_1, \,\dots\, , x_n : A_n$ there exists a 
multimap 
with components $A_1, \,\dots\, , A_n \to \prodop_n(A_1, \,\dots\, , A_n)$; 
that is, a 
rule
\begin{equation} \label{eq:single-var-pairing-rule}
\unaryRule
	{\faketext}
	{x_1 : A_1, \,\dots\, , x_n : A_n \vdash 
		\seqlr{x_1, \,\dots\, , x_n} : \prodop_n(A_1, \,\dots\, , A_n)}
	{} \vspace{-\treeskip}
\end{equation}
\item \label{c:inverse-to-pi-operations} An inverse to precomposing with
$\seq{x_1, \,\dots\, , x_n}$; following the proof of the forward direction of 
Lemma~\ref{lem:clone-representable-iff-cartesian}, we require multimaps
\begin{center}
\unaryRule{\faketext}
		{p : \prodop_n(A_1, \,\dots\, , A_n) \vdash \pi_i(p) : A_i}
		{$(1 \leq i \leq n)$} \vspace{-\treeskip}
\end{center}
such that the equations of 
Lemma~\ref{lem:product-structure-from-representable-arrow} hold,~\ie~that the 
equations
\begin{center}
$\pi_i(\seqlr{x_1, \,\dots\, , x_n}) \equiv x_i \:\: (i = 1, \,\dots\, , n)$ 
\quad and \quad 
$p \equiv \seqlr{\pi_1(p), \,\dots\, , \pi_n(p)}$
\end{center}
obtained by substitution both hold for any 
$x_1 : A_1, \,\dots\, , x_n : A_n$ and 
$p : \prodop_n(A_1, \,\dots\, , A_n)$.
\end{enumerate}

Thus, we recover the laws for products in the 
simply-typed lambda calculus, restricted to variables, from purely 
clone-theoretic reasoning. The usual 
rules, defined on all 
terms, also arise from our abstract considerations. Inspecting the 
proof of Lemma~\ref{lem:clone-representable-iff-cartesian}, one 
sees that for every $(t_i : \Gamma \to X_i)_{i=1, \,\dots\, , n}$ the 
corresponding 
multimap $\Gamma \to \prodop_n(X_1, \,\dots\, , X_n)$ is given by the composite 
$\cslr{\rho_{\ind{X}}}{t_1, \,\dots\, , t_n}$. Translating this into the syntax 
and 
using the standard equality 
$\seqlr{x_1, \,\dots\, , x_n}[t_i/x_i] = \seqlr{t_1, \,\dots\, , t_n}$ defining 
the 
meta-operation of substitution, one arrives at the rule
\begin{center}
\unaryRule
	{(\Gamma \vdash t_i : A_i)_{i=1,\dots,n}}
	{\Gamma \vdash \seqlr{t_1, \,\dots\, , t_n} : \prodop_n(A_1, \,\dots\, , 
	A_n)}
	{} \vspace{-\treeskip}
\end{center}
which, in the presence of substitution, is 
equivalent modulo admissibility to~(\ref{eq:single-var-pairing-rule}).
This is subject to the two equations
$\pi_i\left(\seqlr{t_1, \,\dots\, , t_n}\right) \equiv t_i \:\: (i = 1, 
\,\dots\, , n)$ 
and  
$t \equiv \seqlr{\pi_1(t), \,\dots\, , \pi_n(t)}$.

We therefore recover a presentation of products---modulo 
$\beta\eta$---in the 
simply-typed lambda calculus. More precisely, it is straightforward to see that 
for any $\stlcTimes$-signature $\sig$ the 
clone $\freeCartClone{\sig}$ of Construction~\ref{constr:free-cartesian-clone}
is canonically isomorphic to the syntactic clone 
$\stlcCloneTimes{\sig}$ 
of the simply-typed lambda calculus with products but not 
exponentials (recall Example~\ref{ex:stlc-a-clone} on page~\pageref{ex:stlc-a-clone}).  
Lemma~\ref{lem:free-cartesian-clone} then implies that $\stlcTimes(\sig)$ 
is the internal language of the free cartesian clone on $\sig$. 

We are ultimately interested in the internal language of the free cartesian 
\emph{category} on a (unary) signature. For this we need to show that the cartesian category $\overline{\stlcCloneTimes{\sig}}$, obtained by restricting $\stlcCloneTimes{\sig}$ to unary morphisms, is
the free cartesian category on $\sig$. This is the content of the next 
lemma, in which we call a cartesian functor \Def{strict} if it strictly 
preserves the product-forming operation and each projection. We write $\CartClone$ and 
$\CartCatCat$ for the categories of cartesian clones and cartesian 
categories with their strict morphisms.

As a technical convenience---in 
order to obtain a strict universal property---we shall assume that all the 
cartesian categories (resp. cartesian clones) under consideration have unary 
products given in the canonical way: for every object $A$ the unary product 
$\prodop_1(A)$ is exactly $A$ (recall from Remark~\ref{rem:unary-products} that 
this is a standing assumption for fp-bicategories). 

\begin{mylemma} \label{lem:free-cart-clone-on-cart-cat}
The functor $\overline{(-)} : \CartClone \to \CartCatCat$ restricting a 
cartesian clone to its nucleus has a left adjoint.
\begin{proof}
We show that for any cartesian category $\left(\catC, \Pi_n(-)\right)$, 
cartesian clone $\cartClone{T}{\altClone}$ and strict cartesian functor
$F : \catC \to \overline\altClone$
there exists a cartesian clone $\prom \catC$ and a strict cartesian clone 
homomorphism $\ext{F}: \prom\catC \to \altClone$, unique such that
$\overline{\ext{F}} = F$.

Define $\prom \catC$ as follows. The sorts are the objects of $\catC$ 
and for hom-sets we take 
\begin{center}
$(\prom\catC)(X_1, \,\dots\, , X_n; Y) := 
\catC({X_1 \times \cdots \times X_n}; Y)$
\end{center}
The substitution $\cs{t}{u_1, \,\dots\, , u_n}$ is defined to be the composite
$t \circ \seq{u_1, \,\dots\, , u_n}$ and the projections $\p{i}{\ind{X}}$ are 
the 
projections $\pi_i : \prodop_n(X_1, \,\dots\, , X_n) \to X_i$ for $i=1, 
\,\dots\, , n$.
Since we assume the unary product structure on $\catC$ is the identity, its 
cartesian structure immediately defines a 
cartesian structure on $\prom\catC$. Note in particular that $\prom\catC$ has 
the property that 
$(\prom\catC)(X_1, \,\dots\, , X_n; Y) = 
	(\prom\catC){\left( \prodop_n(X_1, \,\dots\, , X_n); Y\right)}$. 

Now, $\overline{\prom\catC}$ is the cartesian category with objects those of 
$\catC$ and hom-sets of form $\catC{\left(\prodop_1(X), Y \right)}$. So 
$\overline{\prom\catC} = \catC$. We therefore take the unit to be 
$\eta_\catC := \id_{\catC}$.

Next suppose that 
$F : \catC \to \overline{\altClone}$
is a strict cartesian functor. The functor $\ext{F}$ is exactly $F$ on objects, 
while for a multimap $t : X_1, \,\dots\, , X_n \to Y$ in $\prom\catC$ we define
\[
\ext{F}(t) := 
\big( FX_1, \,\dots\, , FX_n 
		\xra{\psi_{F\ind{X}}(\p{1}{}, \,\dots\, , \p{n}{})} 
	\prodop_{i=1}^n FX_i =
	F{\left( \prodop_{i=1}^n X_i \right)} 
		\xra{Ft}
	FY \big)
\]
By the assumption that unary products are the identity, $\ext{F}(u) = F(u)$ for every unary morphism 
$u : X \to Y$. In particular, this holds for the projections $\pi_i$, so 
$\ext{F}$ is a strict cartesian clone homomorphism. 

Finally, suppose that $G : \prom\catC \to \altClone$ is any strict cartesian 
clone homomorphism satisfying $\overline{G} = F$. Since $ob{\prom\catC} = ob\catC$
we must have $FX = GX$ on objects. On arrows, note first that $G$
preserves the tupling operation:
\begin{align*}
G(\psi_{\ind{X}}(\p{1}{}, &\dots, \p{n}{})) \\
	&= \cs{\Id_{\prod_n G\ind{X}}}{G(\psi_{\ind{X}}(\p{1}{}, \dots, \p{n}{}))} \\
	&= \csthree
			{\psi_{G\ind{X}}(\p{1}{}, \dots, \p{n}{})}
			{\pi_1, \dots, \pi_n}
			{G(\psi_{\ind{X}}(\p{1}{}, \dots, \p{n}{}))} 
			&\text{ by Lemma~(\ref{lem:tupling-from-cartesian-clone-gives-identity})} \\
	&= \csthree
			{\psi_{G\ind{X}}(\p{1}{}, \dots, \p{n}{})}
			{G\pi_1, \dots, G\pi_n}
			{G(\psi_{\ind{X}}(\p{1}{}, \dots, \p{n}{}))} 
			&\text{ by strict preservation} \\
	&= 	\cslr
			{\psi_{G\ind{X}}(\p{1}{}, \dots, \p{n}{})}
			{G(\cs{\ind{\pi}}{\psi_{\ind{X}}(\p{1}{}, \dots, \p{n}{})})} \\
	&= \cslr
			{\psi_{G\ind{X}}(\p{1}{}, \dots, \p{n}{})}
			{G(\p{1}{}), \dots, G(\p{n}{})} 
			&\text{ by equation~(\ref{eq:equations-from-universal-arrow})} \\
	&= \psi_{G\ind{X}}(\p{1}{}, \dots, \p{n}{})
\end{align*}
It follows that, for any $t : X_1, \dots, X_n \to Y$ in $\prom\catC$, 
\begin{align*}
\ext{F}(t) &= \cs{(Ft)}{\psi_{F\ind{X}}(\p{1}{}, \dots, \p{n}{})} \\
			&= \cs{(\overline{G}t)}{\psi_{G\ind{X}}(\p{1}{}, \dots, \p{n}{})} \\
			&= \cs{(Gt)}{\psi_{G\ind{X}}(\p{1}{}, \dots, \p{n}{})} \\
			&= G(\cs{t}{\psi_{\ind{X}}(\p{1}{}, \dots, \p{n}){}}) \\
			&= G(t \circ \seq{\pi_1, \dots, \pi_n}) \\
			&= Gt
\end{align*}
where the penultimate equality uses the fact that the cartesian structure of
the clone $\prom\catC$ is inherited from that of the category $\catC$.
Hence $G = \ext{F}$, as required.
\end{proof}
\end{mylemma}

\newpage
With this lemma in hand, one obtains a diagram 
restricting~(\ref{eq:clone-cat-multigraph-diagram})~(p.~\pageref{eq:clone-cat-multigraph-diagram}) 
to the cartesian setting;
the construction of the free cartesian category $\mathbb{FC}\mathrm{at}^{{\times}}(\sig)$
on a unary $\langCart$-signature $\sig$ is standard
(\cf~the construction of the free cartesian closed category in~\cite[Chapter~4]{Crole1994}):
\begin{equation} 
\label{eq:cart-clone-cat-multigraph-diagram}
\begin{tikzcd}[column sep = 5.5em, row sep = 3em]
\: &
\CartClone 
\arrow[bend left = 18]{dr}{\overline{(-)}}
\arrow[bend right = 18]{dl}[swap]{\text{forget}} &
\: \\
\stlcTimesSigCat
\arrow[phantom]{ur}[description]{\adjDown}
\arrow[bend right = 18]{ur}[swap]{\freeCartClone{-}}
\arrow[bend right = 18]{dr}[swap]{\widetilde{\lin}} &
\: &
\CartCatCat
\arrow[phantom]{ul}[description]{\adjDown}
\arrow[bend left = 18]{ul}[near end]{\prom}
\arrow[bend left = 18]{dl}{\text{forget}} \\
\: &
\stlcTimesUnSigCat 
\arrow[phantom]{ur}[description]{\adjUp{}}
\arrow[bend left = 18]{ur}{\mathbb{FC}\mathrm{at}^{{\times}}(-)}
\arrow[phantom]{ul}[description]{\adjUp{}}
\arrow[bend right = 18, hookrightarrow]{ul} &
\:
\end{tikzcd} 
\end{equation}
Moreover, the outer diagram commutes and, as we observed in the proof of the preceding 
lemma, $\nucleus{(-)} \circ \prom = \id_{\CartCatCat}$. One thereby obtains the 
following chain of natural isomorphisms~(\cf~equation~(\ref{eq:nucleus-gives-free-cat})):
\begin{equation} \label{eq:free-cartcat-suffices-to-restrict}
\CartCatCat(\mathbb{FC}\mathrm{at}^{{\times}}(\sig), \catC)
	= 
\CartCatCat{\left(\nucleus{\prom(\mathbb{FC}\mathrm{at}^{{\times}}(\sig))}, \catC\right)}
	\iso 
\CartCatCat{\left(\nucleus{\freeCartClone{\inc\sig})}, \catC\right)}
\end{equation}
Hence, 
just as it was sufficient to construct an 
internal language 
for (bi)clones to describe (bi)categories, so it is sufficient to construct an 
internal language for cartesian clones---namely the simply-typed lambda 
calculus with just products---to describe cartesian categories.

Our aim in the next section
is to reverse this process: we shall lift the theory just presented to the 
bicategorical setting, and use it to extract a principled construction of the 
type theory $\langCart$ with finite products. 

\subsection{Cartesian biclones and representability}

\paragraph*{Representable bi-multicategories.}

Our first step is to bicategorify the 
definition of multicategory. Multicategories can be defined in any monoidal 
category (\eg~\cite[Definition~11.2.1]{Yau2016}); taking the definition in 
$\Cat$ with the product monoidal structure and weakening the equalities to 
isomorphisms suggests the following definition (\cf~also the definition of 
\Def{cartesian 2-multicategory}~\cite{Licata2017}).

\begin{mydefn}
A \Def{bi-multicategory} $\mbCat$ consists of the following data:
\begin{itemize}
\item A class $ob(\mbCat)$ of objects, 
\item For every $X_1, \,\dots\, , X_n, Y \in ob(\mbCat) \:\: (n \in \Nat)$ a 
\Def{hom-category} 
$\left(\mbCat(X_1, \,\dots\, , X_n; Y), \vert, \id\right)$ 
consisting of \Def{multimaps} or \Def{1-cells} 
$f : X_1, \,\dots\, , X_n \to Y$ and \Def{2-cells} $\tau : f \To f'$, subject 
to a 
\Def{vertical composition} operation,
\item For every $X \in ob(\mbCat)$ an \Def{identity} functor $\Id_X : \catOne 
\to \mbCat(X; X)$,
\item For every family of sequences $\Gamma_1, \,\dots\, , \Gamma_n$ and 
objects 
$Y_1, \,\dots\, , Y_n, Z \:\: (n \in \Nat)$ a 
\Def{horizontal composition functor}:
\[
\circ_{\ind{\Gamma}; \ind{Y}; Z} : \mbCat(Y_1, \,\dots\, , Y_n; Z) \times 
\prodop_{i=1}^n \mbCat(\Gamma_i; Y_i) \to \mbCat(\Gamma_1, \,\dots\, , \Gamma_n; 
Z)
\]
We denote the composition $\circ_{\ind{\Gamma}; \ind{Y}; Z}\big( f, (g_1, 
\dots, g_n)\big)$ by $\ms{f}{g_1, \,\dots\, , g_n}$,
\item Natural families of invertible 2-cells
\begin{gather*}
\a_{f; \ind{g}; \ind{h}} : \msthree{f}{\ind{g}}{h_1^{(1)}, \,\dots\, , 
h_{m_1}^{(1)}, \,\dots\, , h_1^{(n)}, \,\dots\, , h_{m_n}^{(n)}} \To 
\ms{f}{\ms{g_1}{\ind{h}^{(1)}}, \,\dots\, , \ms{g_n}{\ind{h}^{(n)}}} \\
\r_f : f \To \ms{f}{\Id_{Y_1}, \,\dots\, , \Id_{Y_n}} \\
\l_f : \ms{\Id_Z}{f} \To f
\end{gather*}
for all $f : Y_1, \,\dots\, , Y_n \to Z$, 
$(g_i : X_1^{(i)}, \,\dots\, , X_{m_n}^{(i)} \to Y_i)_{i=1, \,\dots\, , n}$ and 
${(h_j^{(i)} : \Delta_{j}^{(i)} \to X_j^{(i)})_{\substack{j=1, \,\dots\, , m_i \\ i=1, \,\dots\, , n}}}$.
\end{itemize}
\enlargethispage*{3\baselineskip} 
This data is subject to a triangle law and a pentagon law:
\begin{td}[column sep = 6em]
\ms{f}{g_1, \,\dots\, , g_n} 
\arrow{r}{\ms{\r_{f}}{g_1, \,\dots\, , g_n}} 
\arrow[equals]{d} &
\msthree{f}{\Id, \,\dots\, , \Id}{g_1, \,\dots\, , g_n} 
\arrow{d}{\a_{(f; \Id_{\ind{Y}}; \ind{g})}} \\
\ms{f}{g_1, \,\dots\, , g_n} &
\ms{f}{\ms{\Id}{g_1, \,\dots\, , g_n}, \,\dots\, , \ms{\Id}{g_1, \,\dots\, , 
g_n}} 
\arrow{l}{\ms{f}{\l_{g_1}, \,\dots\, , \l_{g_n}}}
\end{td}
\begin{td}[column sep = huge]
\ms{\big( \msthree{f}{\ind{g}}{\ind{h}} \big)}{\ind{i}} 
\arrow{rr}{\a_{(\ms{f}{\ind{g}}; \ind{h}; \ind{i})}} 
\arrow[swap]{d}{\ms{\a_{(f; \ind{g}; \ind{i})}}{\ind{i}}} &
\: &
\ms{(\ms{f}{\ind{g}})}{\ms{\ind{h}}{\ind{i}}} 
\arrow{d}{\a_{(f; \ind{g}; \ms{\ind{h}}{\ind{i}})}} \\
\ms{ \big(\ms{f}{\ms{\ind{g}}{\ind{h}}}\big) }{\ind{i}} 
\arrow[swap]{r}{\a_{(f; \ms{\ind{g}}{\ind{h}}; \ind{i})}} &
\ms{f}{\msthree{\ind{g}}{\ind{h}}{\ind{i}}} 
\arrow[swap]{r}[yshift=-2mm]
	{\ms{f}
		{\a_{(g_1; \ind{h}; \ind{i})}, \,\dots\, , \a_{(g_n; \ind{h}; 
		\ind{i})}}} &
\ms{f}{\ms{\ind{g}}{\ms{\ind{h}}{\ind{i}}}}
\end{td}
A multimap (resp. 2-cell) of form $f : X \to Y$ 
(resp. $\tau : f \To f' : X \to Y$) is called \Def{linear}.\qedhere
\end{mydefn}

\begin{mynotation}
Note that, just as for clones and multicategories, we use square brackets to 
denote biclone substitution and angle brackets to denote bi-multicategory 
composition~(\cf~Notation~\ref{not:multicat-vs-clone-sub}).
\end{mynotation}

\begin{myremark}
It is natural to conjecture that a construction similar to 
Construction~\ref{constr:free-biclone} would enable one to construct the free 
bi-multicategory on a 2-multigraph and hence a \Def{linear} version of 
$\langBiclone$. Then the argument of Section~\ref{sec:coherence-for-biclones} 
should readily extend to a coherence theorem for bi-multicategories.
\end{myremark}

Examples of bi-multicategories arise naturally, mirroring the 1-categorical 
situation. Every bi-multicategory $\mbCat$ gives rise to a bicategory 
$\overline{\mbCat}$ by restricting to the linear multimaps and their 2-cells 
\big(\cf~Example~\ref{ex:CartesianBicatABiclone}(\ref{c:bicat-as-biclone})\big),
 and---by the following lemma---every monoidal bicategory gives rise to a 
 bi-multicategory (\cf~\cite[Definition 9.2]{Hermida2000}).

\newpage
\begin{mylemma}
Every monoidal bicategory $(\baseCat, \tens, I)$ induces 
a bi-multicategory. 
\begin{proof}
By the coherence theorem for tricategories~\cite{Gordon1995}, we may assume 
without loss of generality that the monoidal bicategory is a \Def{Gray monoid}, 
\ie~a monoid in the monoidal category $\mathrm{Gray}$ 
(see~\eg~\cite[Chapter 3]{Gurski2013} and~\cite[Definition 3.8]{Houston2007}). 
Since Gray monoids also satisfy a 
coherence theorem, we may assume that the underlying 
bicategory $\baseCat$ is a 2-category, and that any pair of composites of the 
structural equivalences 
$a_{A,B,C} : (A \tens B) \tens C \to A \tens (B \tens C)$, 
$l_A : I \tens A \to A$ and
$r_A : A \tens I \to A$
are related by a unique isomorphism (see~\cite[Theorem~10.4]{Gurski2006} 
and~\cite[Theorem~4.1]{Houston2007}). 

The bi-multicategory 
$\int \baseCat$ 
has objects those of $\baseCat$ and hom-categories
$(\int \baseCat)(X_1, \,\dots\, , X_n; Y) := \baseCat(X_1 \tens \cdots \tens 
X_n, 
Y)$, 
where we specify the left-most bracketing 
$\big( ((X_1 \tens X_2) \tens X_3) \tens \cdots \big) \tens X_n$. 

For sequences of objects 
$\Gamma_i := (A_j^{(i)})_{j=1, \,\dots\, , m_i} 
\:\: (i=1, \,\dots\, , n)$
and multimaps
${(g_i : \Gamma_i \to X_i)_{i=1, \,\dots\, ,n}}$ 
and
$f : X_1 \tens \cdots \tens X_n \to Y$, 
the composite 
$\ms{f}{g_1, \,\dots\, , g_n}$ is defined to be 
\[
A_1^{(1)}
	\tens \cdots \tens 
	A_1^{(i)} 
	\tens \cdots \tens
	A_{m_i}^{(i)} 
	\tens \cdots \tens 
	A_1^{(n)} 
	\tens \cdots \tens 
	A_{m_n}^{(n)}
\xra{\simeq} 
\bigotimes_{i=1}^n \Gamma_i 
\xra{\bigotimes_{i=1}^n g_i}
X_1 \tens \cdots \tens X_n 
\xra{f}
Y
\]
where the equivalence is the canonical such. By the coherence theorem for Gray 
monoids, there is a unique choice of isomorphism for each of the structural 
2-cells, and these must satisfy the triangle and pentagon laws.
\end{proof}
\end{mylemma}

For morphisms of bi-multicategories we borrow the terminology from $\Bicat$. 
Thus, bi-multicategories are related by \Def{pseudofunctors}, 
\Def{transformations} 
and \Def{modifications}. 

\begin{mydefn} \quad 
\begin{enumerate}
\item A \Def{pseudofunctor} $F : \mbCat \to \mbCat'$ of bi-multicategories 
consists 
of:
\begin{enumerate}
\item A map $F : ob(\mbCat) \to ob(\mbCat')$ on objects,
\item A functor $F_{\ind{X}; Y} : 
\mbCat(X_1, \,\dots\, , X_n; Y) \to \mbCat'(FX_1, \,\dots\, , FX_n; FY)$ for 
every 
sequence of objects 
$X_1, \,\dots\, ,X_n, Y \in ob(\mbCat) \:\: (n \in \Nat)$, 
\item An invertible 2-cell 
$\psi_X : \Id_{FX} \To F\Id_X$ for every $X \in ob(\mbCat)$, 
\item An invertible 2-cell 
$\phi_{f; \ind{g}} : 
	\ms{F(f)}{Fg_1, \,\dots\, , Fg_n} \To 
	F\left( \ms{f}{g_1, \,\dots\, , g_n} \right)$ for every 
${f : X_1, \,\dots\, , X_n \to Y} \:\: {(n \in \Nat)}$ and 
$(g_i : \Gamma_i \to X_i)_{i=1, \,\dots\, , n}$ in $\mbCat$,  
natural in the sense of 
Definition~\ref{def:functor-and-nat-trans-for-multicat}(\ref{c:multi-nat-trans}).
\end{enumerate}
This data is subject to the following three coherence laws:
\small
\begin{center}
\begin{tikzcd}[column sep = 1.2em]
\ms{\Id_{FZ}}{Ff} 
\arrow{r}{\l_{Ff}} 
\arrow[swap]{d}{\ms{\psi_Z}{Ff}} &
Ff \\
\ms{F(\Id_Z)}{Ff} 
\arrow[swap]{r}[yshift=-2mm]{\phi_{(\Id_Z; f)}} &
F\big( \ms{\Id_Z}{f} \big) 
\arrow[swap]{u}{F\l_f}
\end{tikzcd}
\:\:\:
\begin{tikzcd}[column sep = 1.2em]
Ff 
\arrow{r}{F\r_f} 
\arrow[swap]{d}{\r_{Ff}} &
F\left( \ms{f}{\Id_{Y_1}, \,\dots\, , \Id_{Y_n}} \right) \\
\ms{F(f)}{\Id_{FY_1},\dots, \Id_{FY_n}}  
\arrow[swap]{r}[yshift=-2mm]{\ms{F(f)}{\psi_{Y_1}, \,\dots\, , \psi_{Y_n}}} &
\ms{F(f)}{F\Id_{Y_1},\dots, F\Id_{Y_n}} 
\arrow[swap]{u}{\phi_{(f; \Id_{F\ind{Y}})}}
\end{tikzcd}
\end{center}
\normalsize
\begin{td}[column sep = 7em]
\msthree{Ff}{F\ind{g}}{F\ind{h}} 
\arrow{r}[yshift=0mm]{\a_{(Ff; F\ind{g}; F\ind{h})}} 
\arrow[swap]{d}{\ms{\phi_{(f; \ind{g})}}{F\ind{h}}} &
\mslr{F(f)}{\ms{Fg_1}{F\ind{h}^{(1)}}, \,\dots\, , \ms{Fg_n}{F\ind{h}^{(n)}}} 
\arrow{d}{\ms{F(f)}{\phi_{(g_1; \ind{h})}, \,\dots\, , \phi_{(g_n; \ind{h})}}} 
\\

\mslr{F\left(\ms{f}{\ind{g}}\right)}
	{F\ind{h}} 
\arrow[swap]{d}{\phi_{(\ms{f}{\ind{g}}; \ind{h})}} &
\mslr{Ff}
	{F\big(\ms{g_1}{\ind{h}^{(1)}}\big), \,\dots\, , 
		F\big(\ms{g_n}{\ind{h}^{(n)}}\big)} 
\arrow{d}{\phi_{(f; \ms{\ind{g}}{\ind{h}^{(\bullet)}})}} \\
F\left( \msthree{f}{\ind{g}}{\ind{h}} \right) 
\arrow[swap]{r}{F\a_{(f; \ind{g}; \ind{h})}} &
F\left( \mslr{f}
		{\ms{g_1}{\ind{h}^{(1)}}, \,\dots\, , \ms{g_n}{\ind{h}^{(n)}}}
	\right)
\end{td}

\item A \Def{transformation} $(\alpha, \overline\alpha) : F \To F'$ between 
pseudofunctors 
$F, F' : \mbCat \to \mbCat$ of bi-multicategories consists of
\begin{enumerate}
\item A linear multimap $\alpha_X : FX \to F'X$ for every 
$X \in \mbCat$, 
\item A 2-cell
$\overline{\alpha}_f : \ms{\alpha_Z}{Ff} \To 
	\ms{Gf}{\alpha_{Y_1}, \,\dots\, , \alpha_{Y_n}}$ 
for every 
$f : Y_1, \,\dots\, , Y_n \to Z$ in $\mbCat$, natural in $f$ in the sense of 
Definition~\ref{def:functor-and-nat-trans-for-multicat}(\ref{c:multi-nat-trans}).
\end{enumerate}

\enlargethispage{4\parskip}
This data is subject to the following associativity and unit laws
for every 
$f : Y_1, \,\dots\, , Y_n \to Z$ 
and 
$(g_i : \Gamma_i \to Y_i)_{i=1, \,\dots\, , n}$
in $\mbCat$:
\begin{td}[column sep = 3.5em]
\ms{\Id_{GY}}{\alpha_Y} 
\arrow[swap]{d}{\l_{\alpha_Y}} 
\arrow[]{rr}{\ms{\psi_Y}{\alpha_Y}} &
\: &
\ms{G\Id_Y}{\alpha_Y} \\

\alpha_Y 
\arrow[swap]{r}{\r_{\alpha_Y}} &
\ms{\alpha_Y}{\Id_{FY}} \arrow[swap]{r}{\ms{\alpha_Y}{\psi_Y}} &
\ms{\alpha_Y}{F\Id_Y} \arrow[swap]{u}{\overline{\alpha}_{\Id_Y}} 
\end{td}
\begin{small}
\vspace{-2mm}
\begin{td}[column sep = 2em]
\msthreelr{\alpha_Y}{Ff}{F\ind{g}} 
\arrow{r}[yshift=0mm]{\a_{(\alpha_Y; Ff; F\ind{g})}} 
\arrow[swap]{d}{\ms{\overline{\alpha}_{f}}{F\ind{g}}} &
\mslr{\alpha_Y}{\left(\ms{F(f)}{F\ind{g}}\right)} 
\arrow{r}[yshift=2mm]{\ms{\alpha_Y}{\phi_{(f; \ind{g}})}} &
\mslr{\alpha_Y}{F\left(\ms{f}{\ind{g}}\right)} 
\arrow{dddd}{\overline{\alpha}_{\ms{f}{\ind{g}}}} \\
\msthreelr{G(f)}{\alpha_{Y_1}, \,\dots\, , \alpha_{Y_n}}{F\ind{g}} 
\arrow[swap]{d}{\a_{(Gf; \alpha_{\ind{Y}}; F\ind{g})}} 
&
\: &
\: \\
\mslr{G(f)}{\ms{\alpha_{Y_1}}{Fg_1}, \,\dots\, , \ms{\alpha_{Y_n}}{Fg_n}} 
\arrow[swap]{d}
	{\ms{G(f)}{\overline{\alpha}_{g_1}, \,\dots\, , \overline{\alpha}_{g_n}}} 
&
\: &
\: \\
\mslr{G(f)}{\ms{Gg_1}{\alpha_{\Gamma_1}}, \,\dots\, , 
\ms{Gg_n}{\alpha_{\Gamma_n}}} 
\arrow[swap]{d}{\a^{-1}_{(Gf; G\ind{g}; \ind{\alpha})}} &
\: &
\: \\
\msthreelr{G(f)}{Gg_1, \,\dots\, , Gg_n}{\ind{\alpha}} 
\arrow[swap]{rr}{\ms{\phi_{(f; \ind{g})}}{\ind{\alpha}}} &
\: &
\ms{G\big(\ms{f}{\ind{g}}\big)}{\ind{\alpha}} 
\end{td}
\end{small}
Note that, where 
$\Gamma_i := A_1^{(i)}, \,\dots\, , A_{m_i}^{(i)}$, 
we write $\alpha_{\Gamma_i}$ for the sequence 
$\alpha_{A_1^{(i)}}, \,\dots\, , \alpha_{A_{m_i}^{(i)}}$.

\item A \Def{modification} $\modif : (\alpha, \overline{\alpha}) \to (\beta, 
\overline{\beta})$ between transformations $(\alpha, \overline{\alpha}), 
(\beta, 
\overline\beta) : F \To F'$ is a family of 2-cells $\modif_X : \alpha_X \To 
\beta_X$ such that the following diagram commutes for every 
$f : Y_1, \,\dots\, , Y_n \to Z$:
\vspace{-2mm}
\begin{td}[column sep = 8em]
\ms{\alpha_Z}{Ff} \arrow{r}{\ms{\modif_Z}{Ff}} 
\arrow[swap]{d}{\overline{\alpha}_f} &
\ms{\beta_Z}{Ff} \arrow{d}{\overline{\beta}_f} \\
\ms{G(f)}{\alpha_{Y_1}, \,\dots\, , \alpha_{Y_n}} 
\arrow[swap]{r}[yshift=-2mm]{\ms{G(f)}{\modif_{Y_1}, \,\dots\, , \modif_{Y_n}}} 
&
\ms{G(f)}{\beta_{Y_1}, \,\dots\, , \beta_{Y_n}}
\end{td} 
\end{enumerate}
\vspace{-3\parskip}
\qedhere
\end{mydefn}

One would expect that bi-multicategories, pseudofunctors, transformations and 
modifications organise themselves into a tricategory; we do not pursue such 
considerations here. Instead, we lift Hermida's notion of representability to  
bi-multicategories.
As usual, it is convenient to require as much as possible of the definition to 
be data.

\begin{mydefn}
A \Def{representable bi-multicategory} $(\mbCat, \tensor_n)$ consists of the 
following data:
\begin{enumerate}
\item For every $X_1, \,\dots\, , X_n \in \mbCat \:\: (n \in \Nat)$, a chosen  
object
$\tensor_n(X_1, \,\dots\, , X_n) \in \mbCat$ and chosen
\Def{birepresentable}
multimap
$\rho_{\ind{X}} : X_1, \,\dots\, , X_n \to \tensor_n(X_1, \,\dots\, , X_n)$, 
such that the birepresentable multimaps are closed under composition,

\item For every $A, X_1, \dots, X_n \in \mbCat \:\: (n \in \Nat)$, an adjoint 
equivalence
\begin{td}
\mbCat{\left(\tensor_n(X_1, \,\dots\, , X_n); A \right)}
\arrow[bend left = 20]{r}{\ms{(-)}{\rho_{\ind{X}}}}
\arrow[phantom]{r}[description]{\adjUp{\simeq}} &
\mbCat(X_1, \,\dots\, , X_n; A)
\arrow[bend left= 20]{l}{\psi_{\ind{X}}}
\end{td}
specified by a choice of universal arrow
$\epsilon_{\ind{X}}$. \qedhere
\end{enumerate}
\end{mydefn}

The birepresentability of $\rho_{\ind{X}}$ entails the following. For 
every 
$f : X_1, \,\dots\, , X_n \to A$ we require a choice of multimap 
$\psi_{\ind{X}}(f) : \tensor_n(X_1, \,\dots\, , X_n) \to A$ and 2-cell 
$\epsilon_{\ind{X}; f} : \ms{\psi_{\ind{X}}(f)}{\rho_{\ind{X}}} \To f$. This 
2-cell is universal in the sense that for any 
$g : \tensor_n(X_1, \,\dots\, , X_n)\to A$ and 
$\sigma : \ms{g}{\rho_{\ind{X}}} \To f$ there exists a unique 2-cell 
$\trans{\sigma} : g \To \psi_{\ind{X}}(f)$ such that
\begin{equation} \label{eq:birepresentability-ump}
\begin{tikzcd}
\ms{g}{\rho_{\ind{X}}} \arrow{rr}{\ms{\trans{\sigma}}{\rho_{\ind{X}}}} 
\arrow[swap]{dr}{\sigma} &
\: &
\ms{\psi_{\ind{X}}(f)}{\rho_{\ind{X}}} \arrow{dl}{\epsilon_{\ind{X}; f}} \\
\: &
f &
\:
\end{tikzcd} 
\end{equation}

\begin{myremark}
Hermida's construction suggests that every representable bi-multicategory 
ought to induce a monoidal bicategory, and indeed that there exists a  
triequivalence between representable bi-multicategories and monoidal 
bicategories. 
Here we shall restrict ourselves to proving that every 
representable biclone induces an fp-bicategory: a considerably easier 
task, as one only needs to check a universal property, rather 
than many coherence axioms. 
\end{myremark}

Following the 1-categorical template of 
Section~\ref{sec:cartesian-clones-and-representability}, we next examine the 
construction of finite products in a 
bi-multicategory. To avoid the double prefix in `fp-bi-multicategories' we 
refer to such objects as `cartesian bi-multicategories'.

\paragraph*{Cartesian bi-multicategories.}

Once again, we translate between the categorical and 
bicategorical settings by replacing 
universal arrows with biuniversal arrows. 

\begin{mydefn} \label{def:universal-arrow-bi-multicat} Let $F : \mbCat \to 
\mbCat'$ be a pseudofunctor of 
bi-multicategories and 
$X 
\in \mbCat'$. A \Def{biuniversal arrow $(R, u)$ from $F$ to $X$} 
consists of
\begin{enumerate}
\item An object $R \in \mbCat$, 
\item A linear multimap $u : FR \to X$, 
\item For every $A \in \mbCat$, a chosen adjoint equivalence 
\begin{td}
\mbCat(A_1, \,\dots\, , A_n; R) 
\arrow[bend left = 20]{r}{\ms{u}{F(-)}}
\arrow[phantom]{r}[description]{\adjUp{\simeq}} &
\mbCat'(FA_1, \,\dots\, , FA_n; X)
\arrow[bend left= 20]{l}{\psi_{\ind{A}}}
\end{td}
\vspace{-1\parskip}
specified by a choice of universal arrow 
$\epsilon_h : \ms{u}{F\psi_{\ind{A}}(h)} \To h : FA_1, \,\dots\, , FA_n \to X$
(\cf~Definition~\ref{def:biuniversal-arrow}).
\qedhere
\end{enumerate}
\end{mydefn}

We translate this into a `global' definition in the by-now-familiar way. 

\begin{mylemma} \label{lem:biuniversal-arrow-iff-family-of-equivalences}
For any pseudofunctor of bi-multicategories 
$F : \mbCat \to \mbCat'$ and 
$X \in \mbCat'$, the following are equivalent:
\begin{enumerate}
\item A choice of biuniversal arrow from $F$ to $X$, 
\item Chosen adjoint equivalences 
	$\kappa_{\ind{A}} : 
		\mbCat(A_1, \,\dots\, , A_n; R) \leftrightarrows 
		\mbCat'(FA_1, \,\dots\, , FA_n; X) : 
	\delta_{\ind{A}}$ 
	for $A_1, \,\dots\, , A_n \in \mbCat \:\: (n \in \Nat)$,
	specified by a choice of universal arrow and
	pseudonatural in the sense that for every 
$f : A_1, \,\dots\, , A_n \to R$ and 
$(g_i : \Gamma_i \to A_i)_{i=1, \,\dots\, , n}$ 
there exists an 
invertible 2-cell 
$\nu_{f ; \ind{g}} : 
	\ms{\kappa_{\ind{A}}(f)}{Fg_1, \,\dots\, , Fg_n} \To 
	\kappa_{\ind{A}}\left(\ms{f}{g_1, \,\dots\, , g_n}\right)$, 
multinatural in 
$f, g_1, \,\dots\, , g_n$ and satisfying
\begin{align}
&\qquad\quad\begin{tikzcd}[ampersand replacement = \&] 
\label{eq:biuniversal-arrow-nat-unit}
\kappa_{\ind{A}}(f) 
\arrow{r}{\kappa_{\ind{A}}(\r_f)} 
\arrow[swap]{d}{\r_{\kappa_{\ind{A}}(f)}} \&
\kappa_{\ind{A}}\left(\ms{f}{\Id_{\ind{A}}}\right) \\
\ms{\kappa_{\ind{A}}(f)}{\Id_{\ind{A}}} 
\arrow[swap]{r}[yshift=-2mm]{\ms{\kappa_{\ind{A}}(f)}{\ind{\psi}}}  \&
\ms{\kappa_{\ind{A}}(f)}{F\Id_{\ind{A}}} 
\arrow[swap]{u}{(\nu_{f; \Id_{\ind{A}}})} \\
%
%
\faketext \&
\: 
\end{tikzcd} \\
\label{eq:biuniversal-arrow-nat-assoc}
&\begin{tikzcd}[ampersand replacement = \&, column sep = 3em]
\msthree{\kappa_{\ind{A}}(f)}{F\ind{g}}{F\ind{h}} 
\arrow[swap]{d}{\ms{\nu_{(g; \ind{f})}}{F\ind{h}}} 
\arrow{r}[yshift=2mm]{\a_{(\kappa_{\ind{A}}(f); F\ind{g}; F\ind{h})}} \&
\mslr{\kappa_{\ind{A}}(f)}
	{\ms{Fg_1}{F\ind{h}^{(1)}}, \,\dots\, , \ms{Fg_n}{F\ind{h}^{(n)}}} 
\arrow{d}{\ms{\kappa_{\ind{A}}(f)}
	{\phi_{(g_1; \ind{h})}, \,\dots\, , \phi_{(g_n; \ind{h})}}} \\
\mslr{\kappa_{\ind{\Gamma}}{\left( \ms{f}{\ind{g}} \right)}}{F\ind{h}} 
\arrow[swap]{d}{\nu_{(\ms{f}{\ind{g}}; h)}} \&
\mslr{\kappa_{\ind{A}}(f)}
	{F{\big(\ms{g_1}{\ind{h}^{(1)}}\big)}, \,\dots\, , 
		F{\big(\ms{g_n}{\ind{h}^{(n)}}\big)}} 
\arrow{d}{\nu_{(f; \ms{\ind{g}}{\ind{h}})}} \\
\kappa_{\ind{\Delta}}{\left( \msthree{f}{\ind{g}}{\ind{h}} \right)} 
\arrow[swap]{r}{\kappa_{\ind{\Delta}}(\a_{(f; \ind{g}; \ind{h})})} \&
\kappa_{\ind{\Delta}}{\left( 
		\mslr{f}{\ms{g_1}{\ind{h}^{(1)}}, \,\dots\, , \ms{g_n}{\ind{h}^{(n)}}} 
	\right)}
\end{tikzcd}
\end{align}
for $\Gamma_i := X_1^{(i)}, \,\dots\, , X_{m_i}^{(i)}$ and $(h_j^{(i)} : 
\Delta_j^{(i)} \to X_j^{(i)})_{\substack{j=1, \,\dots\, , m_i \\ i=1,\dots,n}}$.
\end{enumerate}
\begin{proof}
\subproof{(1)$\To$(2)} By biuniversality, $\ms{u}{F(-)}$ is part of an adjoint 
equivalence for every $A_1, \,\dots\, , A_n \in \mbCat \:\: (n \in \Nat)$, so 
it remains to check 
pseudonaturality. Taking $\kappa_{\ind{A}}$ to be $\ms{u}{F(-)}$, we are 
required to provide 2-cells $\nu_{f; \ind{g}}$ of type 
$\msthree{u}{Ff}{Fg_1, \,\dots\, , Fg_n} \To 
	\ms{u}{F\big( \ms{f}{g_1, \,\dots\, , g_n} \big)}$, 
for which we take 
$\big(\ms{u}{\phi_{f; \ind{g}}}\big) \vert \a_{u; Ff; F\ind{g}}$. 
The naturality 
condition and two axioms~(\ref{eq:biuniversal-arrow-nat-unit}) 
and~(\ref{eq:biuniversal-arrow-nat-assoc}) then follow directly from the 
coherence laws of a pseudofunctor. 

\subproof{(2)$\To$(1)} This direction is a little more delicate, but we can 
follow the template provided by 
Lemma~\ref{lem:multicat-universal-arrow-to-natural-isomorphism}. Let us first 
make explicit the content of the adjoint equivalence 
\[
{\kappa_{\ind{A}} : 
\mbCat(A_1, \,\dots\, , A_n; R)} \leftrightarrows \mbCat'(FA_1, \,\dots\, , 
FA_n; X) : 
\delta_{\ind{A}}
\]
Choosing a universal arrow entails that for every 
$f : FA_1, \,\dots\, , FA_n \to X$ 
there exists a multimap
$\delta_{\ind{A}}(f) : A_1, \,\dots\, , A_n \to R$ and a 2-cell 
${\overline{\delta}_f : \kappa_{\ind{A}}\delta_{\ind{A}}(f)} \To f$, universal 
in 
the sense that for any 
$g : A_1, \,\dots\, , A_n \to R$ and 
$\sigma : \kappa_{\ind{A}}(g) \To f$ there exists a unique 2-cell 
$\altTrans{\sigma} : g \To \delta_{\ind{A}}(f)$ such that
\begin{equation} \label{eq:multinaturality-adjoint-eq-universality}
\begin{tikzcd}
\kappa_{\ind{A}}(g) 
\arrow[swap]{dr}{\sigma}
\arrow{rr}{\kappa_{\ind{A}}(\altTrans{\sigma})}  &
\: &
\kappa_{\ind{A}}\delta_{\ind{A}}(f) 
\arrow{dl}{\overline{\delta}_f} \\ 
\: &
f &
\:
\end{tikzcd}
\end{equation}
We claim that $u := \kappa_R(\Id_R) : FR \to X$ is biuniversal. Thus, for 
every $f : FA_1, \,\dots\, , FA_n \to X$ we need to provide 
an arrow $\overline{f} : A_1, \dots A_n \to R$ and a universal 2-cell 
$\epsilon_{\ind{A}; f} : \ms{u}{F\overline{f}} \To f$. 

For the arrow we take $\overline{f} := \delta_{\ind{A}}(f)$. For the 2-cell we 
make use of the naturality condition to define $\epsilon_{\ind{A}; f}$ as the 
invertible composite
\begin{td}[column sep = 5em]
\ms{u}{F\delta_{\ind{A}}(f)}
\arrow[equals]{d}
\arrow{rr}{\epsilon_{\ind{A}; f}} &
\: &
f \\

\ms{\kappa_R(\Id_R)}{F\delta_{\ind{A}}(f)}
\arrow[swap]{r}
	{\nu_{(\Id_R; \delta_{\ind{A}}(f))}} &
\kappa_{\ind{A}}\left( \ms{\Id_R}{\delta_{\ind{A}}(f)} \right) 
\arrow[swap]{r}{\kappa_{\ind{A}}(\l_{\delta_{\ind{A}}(f)})} &
\kappa_{\ind{A}}\delta_{\ind{A}}(f)
\arrow[swap]{u}{\overline{\delta}_f}
\end{td}
To establish universality, let
$g : A_1, \,\dots\, , A_n \to R$ 
be a multimap and 
$\gamma : \ms{u}{Fg} \To f$ 
be any 2-cell. We need to show there 
exists a unique 2-cell 
$\trans{\gamma} : g \To \overline{f}$ such that 
\begin{equation} \label{eq:biuniversal-arrow-nat-equivalence-diagram}
\begin{tikzcd}
\ms{u}{Fg} \arrow{rr}{\ms{u}{F\trans{\gamma}}} \arrow[swap]{dr}{\gamma} &
\: &
\ms{u}{F\overline{f}} \arrow{dl}{\epsilon_{\ind{A}; f}} \\
\: &
f
\end{tikzcd}
\end{equation}
By the universal property~(\ref{eq:multinaturality-adjoint-eq-universality}), 
to define 
$\trans{\gamma} : g \To \overline{f} = \delta_{\ind{A}}(f)$ it 
suffices to define a 2-cell $\kappa_{\ind{A}}(g) \To f$, for which we 
take
\[
\alpha_{\gamma, f, g} := \kappa_{\ind{A}}(g) \XRA{\kappa_{\ind{A}}(\l_g^{-1})} 
\kappa_{\ind{A}}(\ms{\Id_R}{g}) \XRA{\nu_{\Id_R; g}^{-1}} 
\ms{\kappa_{\ind{A}}(\Id_R)}{Fg}\XRA{\gamma} f
\]
We define $\trans{\gamma} := \altTrans{(\alpha_{\gamma, f, g})}$. That this 
fills~(\ref{eq:biuniversal-arrow-nat-equivalence-diagram}) is an easy check 
using the definition and naturality of $\nu$. For uniqueness, suppose $\sigma : 
g \To \overline{f} = \delta_{\ind{A}}(f)$ also 
fills~(\ref{eq:biuniversal-arrow-nat-equivalence-diagram}). By the universal 
property defining $\trans\gamma$ it suffices to show 
that $\sigma$ is the unique 2-cell corresponding to $\alpha_{\gamma, f, g}$ 
via~(\ref{eq:multinaturality-adjoint-eq-universality}). 
This follows from the naturality of $\nu$ and $\l$ and the 
definition of $\alpha_{\gamma, f, g}$. 

This completes the construction of 
an \emph{adjunction} 
$\mbCat(A_1, \,\dots\, , A_n; R) \leftrightarrows \mbCat'(FA_1, \,\dots\, , 
FA_n; X)$; to 
show this is an adjoint equivalence, we need to show the 
unit is also invertible. But the unit is given by applying the $\trans{(-)}$ 
operation to the identity, \ie~by applying the $\altTrans{(-)}$ operation to an 
invertible 2-cell. This is invertible by 
Lemma~\ref{lem:trans-natural-and-invertible}. 
\end{proof}
\end{mylemma}

The definition of product of multicategories lifts straightforwardly to 
bi-multicategories. For bi-multicategories $\mbCat$ and $\mbCat'$, the bi-multicategory
$\mbCat \times \mbCat'$ has objects pairs 
$(X, X') \in ob(\mbCat) \times ob(\mbCat')$ 
and composition as in~(\ref{eq:product-of-multicats-composition}) 
on page~\pageref{eq:product-of-multicats-composition}. The 
structural isomorphisms are given pointwise. Then there exists a canonical 
\Def{diagonal pseudofunctor} 
$\Delta^n : \mbCat \to \mbCat^{\times n}$ for every
bi-multicategory $\mbCat$ and
$n \in \Nat$.  
This suggests the following definition. 

\begin{mydefn}
A \Def{cartesian bi-multicategory} $\fpBicat\mbCat$ consists of a
bi-multicategory $\mbCat$ equipped with the following data
for every 
$X_1, \,\dots\, , X_n \in \mbCat \:\: (n \in \Nat)$: 
\begin{enumerate}
\item A chosen object $\prodop_n(X_1, \,\dots\, , X_n)$, 
\item A choice of biuniversal arrow 
$\pi = (\pi_1, \,\dots\, , \pi_n) : 
	\Delta^n {\left(\prodop_n(X_1, \,\dots\, ,X_n)\right)} \to (X_1, \,\dots\, 
	, X_n)$
from $\Delta^n$ to $(X_1, \,\dots\, , X_n) \in \mbCat^{\times n}$. 
\qedhere
\end{enumerate}
\end{mydefn}

By the preceding lemma, a bi-multicategory is cartesian if and only if there 
exists a pseudonatural 
 family of adjoint equivalences
\begin{equation*}
\mbCat\big(\Gamma; \prodop_n(X_1, \,\dots\, , X_n)\big) \simeq 
\mbCat^{\times n}{\left(\Delta^n(\Gamma); (X_1, \,\dots\, , X_n) \right)} = 
\prod_{i=1}^n \mbCat(\Gamma ; X_i)
\end{equation*}
The universal property therefore manifests itself as follows. For every 
sequence of multimaps 
$(t_i : \Gamma \to X_i)_{i=1, \,\dots\, , n}$
there exists a multimap 
$\pair{t_1, \,\dots\, , t_n} : \Gamma \to \prodop_n(X_1, \,\dots\, , X_n)$
and a 2-cell
$\epsilonTimes{}$ 
with components
$\epsilonTimesInd{i}{\ind{t}} :
	\ms{\pi_i}{\pair{t_1, \,\dots\, , t_n}} \To t_i$
for $i=1, \,\dots\, , n$. This 2-cell is universal in the sense that, if 
$u : \Gamma \to \prodop_n(X_1, \,\dots\, , X_n)$
and
$\alpha_i : \ms{\pi_i}{u} \To t_i$
for $i=1, \,\dots\, , n$, then there exists a unique 2-cell
$\transTimes{\alpha_1, \,\dots\, , \alpha_n} : 
	u \To \pair{t_1, \,\dots\, , t_n}$
filling the following diagram for $i=1, \,\dots\, , n$:
\begin{equation} \label{eq:cartesian-bimulticat-ump}
\begin{tikzcd}
\ms{\pi_i}{u} 
\arrow[swap]{dr}{\alpha_i}
\arrow{rr}{\ms{\pi_i}{\alpha}} &
\: &
\ms{\pi_i}{\pair{t_1, \,\dots\, , t_n}}
\arrow{dl}{\epsilonTimesInd{i}{\ind{t}}} \\
\: &
t_i &
\:
\end{tikzcd}
\end{equation}
Finally, the unit 
$\eta_u := 
	\transTimes{\id_{\ms{\pi_1}{u}}, \,\dots\, , \id_{\ms{\pi_n}{u}}} : 
	u \To \pair{\ms{\pi_1}{u}, \,\dots\, , \ms{\pi_n}{u}}$ 
is required to be invertible for 
every $u : \Gamma \to \prodop_n(X_1, \dots X_n)$. 

Our next task is to extend the theory of representable and cartesian 
bi-multicategories to biclones.

\paragraph*{Cartesian biclones.}

As we did for clones, we define products in a biclone by first defining a 
bi-multicategory structure on each biclone
(\cf~Construction~\ref{constr:multicategory-from-clone}). 

\begin{myconstr} \label{constr:bimulticat-from-biclone}
Every biclone $(S, \biclone)$ canonically defines a bi-multicategory 
$\mCatOf{\biclone}$ as follows:
\begin{itemize}
\item $ob(\mCatOf\biclone) := S$, 
\item $(\mCatOf{\biclone})(X_1, \,\dots\, , X_n; Y) := \biclone(X_1, \,\dots\, 
, X_n; 
Y)$, 
\item $\Id_X := \p{1}{1} : \catOne \to (\mCatOf{\biclone})(X; X)$, 
\item The composition functor $(\mCatOf{\biclone})(Y_1, \,\dots\, , Y_n; Z) 
\times 
\prod_{i=1}^n (\mCatOf{\biclone})(\Gamma_i; Y_i) \to 
(\mCatOf{\biclone})(\Gamma_1, \,\dots\, , \Gamma_n; Z)$ is defined by
\[
\ms{f}{g_1, \,\dots\, , g_n} := \cslr{f}{g_1 \clonetimes \cdots \clonetimes g_n}
\]
using the notation of Notation~\ref{not:clone-product-notation},
%
\item The unitor structural isomorphisms are defined as follows, for 
$f : X_1, \,\dots\, , X_n \to Y$:
\begin{gather*}
\r_f := f \XRA{\subid{}} \cslr{f}{\p{1}{\ind{X}}, \,\dots\, , \p{n}{\ind{X}}} 
\XRA{\cslr{f}{\indproj{-1}{}, \,\dots\, , \indproj{-1}{}}} 
\cslr{f}{\cslr{\p{1}{X_1}}{\p{1}{\ind{X}}}, \,\dots\, , 
\cslr{\p{1}{X_n}}{\p{n}{\ind{X}}}} 
\\
\l_f := \cslr{\p{1}{Y}}{\cslr{f}{\p{1}{\ind{X}}, \,\dots\, , \p{n}{\ind{X}}}} 
\XRA{\indproj{1}{}} \cslr{f}{\p{1}{\ind{X}}, \,\dots\, , \p{n}{\ind{X}}} 
\XRA{\subid{}^{-1}} f
\end{gather*}
The associativity structural isomorphism is a little complex. Suppose given 
sequences of objects 
$\Gamma_i := B_1^{(i)}, \,\dots\, , B_{m_i}^{(i)} \:\: (i=1, \dots, n)$
and multimaps 
$(g_i : \Gamma_i \to Y_i)_{i=1, \dots,n}$ 
and 
$f :   {Y_1, \,\dots\, , Y_n \to Z}$.
Moreover suppose that 
$\Delta_j^{(i)} := A^{(i,j)}_1, \,\dots\, , A^{(i,j)}_{k(i,j)}$, and that
$h_{j}^{(i)} : \Delta_j^{(i)} \to B_{j}^{(i)}$ 
for 
${j = 1, \,\dots\, , m_i}$ and 
$i = 1, \,\dots\, , n$. 

Now, writing
$\pPicked(R)$ for the 
projection picking out the element $R$ in the codomain, there exists a map
\begin{equation} \label{eq:first-morphism}
\cslr
	{h_j^{(i)}}
	{\pPicked(A^{(i,j)}_1), \,\dots\, , \pPicked(A^{(i,j)}_{k(i,j)})}
	:
	\Delta_1^{(1)}, \,\dots\, , \Delta_{m_1}^{(1)}, 
	\dots, 
	\Delta_1^{(n)}, \,\dots\, , \Delta_{m_n}^{(n)} 
	\to
	B_j^{(i)} 
\end{equation}
for every $i=1, \,\dots\, , n$ and $j= 1, \,\dots\, , m_i$.
On the other hand, one may first project out from the full sequence
$\Delta_1^{(1)}, \,\dots\, , \Delta_{m_1}^{(1)}, 
	\dots, 
	\Delta_1^{(n)}, \,\dots\, , \Delta_{m_n}^{(n)}$
to the subsequence $\Delta_1^{(i)}, \,\dots\, , \Delta_{m_i}^{(i)}$ and then 
project 
again before applying $h_j^{(i)}$.  Abusively writing 
$\left[\pPicked(\Delta_1^{(i)}), \,\dots\, , 
\pPicked(\Delta_{m_i}^{(i)})\right]$
for the sequence
$\left[\pPicked(A_1^{(i,1)}), \,\dots\, , \pPicked(A_{k(i, 
m_i)}^{(i,m_i)})\right]$, 
one thereby obtains
\begin{equation} \label{eq:second-morphism}
\csthree
	{h_j^{(i)}}
	{\pPicked(A^{(i,j)}_1), \,\dots\, , \pPicked(A^{(i,j)}_{k(i,j)})}
	{\pPicked(\Delta_1^{(i)}), \,\dots\, , \pPicked(\Delta_{m_i}^{(i)})}
\end{equation}
The pair of parallel multimaps~(\ref{eq:first-morphism}) 
and~(\ref{eq:second-morphism}) are related by a canonical composite of 
structural isomorphisms:
\begin{equation} \label{eq:multicat-from-biclone-canonical-iso}
\begin{aligned}
h_j^{(i)} 
	&\uncs{\pPicked(A^{(i,j)}_1), \,\dots\, , 
	\pPicked(A^{(i,j)}_{k(i,j)})}
	\uncs{\pPicked(\Delta_1^{(i)}), \,\dots\, , \pPicked(\Delta_{m_i}^{(i)})}  
	\\
	&\qquad\qquad\iso 
\cslr
	{h_j^{(i)}}
	{\dots, 
		\cs
			{\pPicked(A^{(i,j)}_l)}
			{\pPicked(\Delta_1^{(i)}), \,\dots\, , 
			\pPicked(\Delta_{m_i}^{(i)})}, 
	\dots}	\\
	&\qquad\qquad\iso 
\cslr
	{h_j^{(i)}}
	{\pPicked(A^{(i,j)}_1), \,\dots\, , \pPicked(A^{(i,j)}_{k(i,j)})}
\end{aligned}
\end{equation}
Making use of the same notation, 
$\msthree{f}{g_1, \,\dots\, , g_n}{h_1^{(1)},\dots, h_{m_1}^{(1)}, \,\dots\, , 
h_1^{(n)}, 
\dots, h_{m_n}^{(n)}}$ is
\begin{center}
$\csthree{f}
	{\dots, \cslr{g_i}{\pPicked(B_1^{(i)}), \,\dots\, , 
		\pPicked(B_{m_i}^{(i)})}, \dots}
	{\dots, 
		\cslr{h_j^{(i)}}{\pPicked(\Delta_1^{(i)}), \,\dots\, , 
		\pPicked(\Delta_{m_j}^{(j)})}, \dots}$
\end{center}
and 
$\mslr{f}
	{\ms{g_1}{h_1^{(1)},\dots, h_{m_1}^{(1)}}, \,\dots\, , 
		\ms{g_n}{h_1^{(n)}, \,\dots\, , h_{m_n}^{(n)}}}$ 
is
\begin{center}
$\cslr{f}
	{\dots, 
		\csthree{g_i}
				{\dots, \cslr{h_j^{(i)}}{\pPicked(A_1^{(i,j)}), \,\dots\, ,
					\pPicked(A_{k(i,j)}^{(i,j)})},\dots}
				{\pPicked(\Delta_1^{(i)}), \,\dots\, , 
					\pPicked(\Delta_{m_i}^{(i)})}, \dots
	}$
\end{center}
so $\a_{f; \ind{g}; \ind{h}}$ is the composite
\small
\begin{td}
\csthree{f}
	{g_1 \clonetimes \cdots \clonetimes g_n}
	{h_1^{(1)} \clonetimes \cdots \clonetimes h_j^{(i)} \clonetimes \cdots 
	\clonetimes h_{m_n}^{(n)}} 
\arrow[swap]{d}{\mathsf{f}_{f; \ind{g}; \ind{h}}} \\

%

\cslr{f}{\cslr{g_1}{h_1^{(1)} \clonetimes \cdots \clonetimes h_{m_1}^{(1)}}, 
\,\dots\, , 
\cslr{g_n}{h_1^{(n)} \clonetimes \cdots \clonetimes h_{m_n}^{(n)}}} 
\arrow[swap]{d}{\iso} 
\arrow{d}{(\ref{eq:multicat-from-biclone-canonical-iso})} \\

\cslr{f}
	{\dots, 
		\cslr{g_i}
			{\dots, 
				\csthree
					{h_j^{(i)}}
					{\pPicked(A^{(i,j)}_1), \,\dots\, , 
						\pPicked(A^{(i,j)}_{k(i,j)})}
					{\pPicked(\Delta_1^{(i)}), \,\dots\, , 
					\pPicked(\Delta_{m_i}^{(i)})}
					,\dots}, 
\dots}
\arrow[swap]{d}{\iso} \\

\cslr{f}{\dots, 
	\csthree{g_i}
		{\dots, 
			\cslr
				{h_j^{(i)}}
					{\pPicked(A_1^{(i,j)}), \,\dots\, , 
						\pPicked(A_{k(i,j)}^{(i,j)})},\dots}
		{\pPicked(\Delta_1^{(i)}), \,\dots\, , 
					\pPicked(\Delta_{m_i}^{(i)})}, \dots
	}
\end{td}
%
%
%
%
\normalsize
where the final isomorphism is the evident composite of structural isomorphisms 
in $(S, \biclone)$ and $\mathsf{f}_{f; \ind{g}; \ind{h}}$ is defined after
Notation~\ref{not:clone-product-notation} (page~\pageref{not:clone-product-notation}). 
\end{itemize}
The two coherence laws hold by the coherence of biclones.
\end{myconstr}

We now see where the awkwardness in the definition of pseudofunctors and 
transformations of biclones arises (Definitions~\ref{def:biclone-pseudofunctor} 
and~\ref{def:biclone-transformation}): the more natural definitions are for 
bi-multicategories, and the versions for biclones arise via 
Construction~\ref{constr:bimulticat-from-biclone}.

\begin{mynotation}
Following the preceding construction, we sometimes write 
$\Id_A$ for the projection $\p{1}{A}$ in a biclone, and refer to it
as \Def{the identity on $A$}.
\end{mynotation}

\begin{myremark} \label{rem:biclone-to-unary}
For a biclone $(S, \biclone)$, the bicategory 
$\overline\biclone$ 
obtained by restricting to unary hom-categories is biequivalent to
the restriction $\overline{\mCatOf\biclone}$ of the corresponding 
bi-multicategory to linear hom-categories
\big(\cf~(\ref{eq:clone-multicat-restriction})\big). Indeed, the objects and 
hom-categories are equal: the only difference is that for
$f : X \to Y$ and $g : Y \to Z$ in $(S, \biclone)$ the corresponding composite 
in $\overline\biclone$ is $\cslr{f}{g}$ while in 
$\overline{\mCatOf\biclone}$ it is
$\cslr{f}{\cslr{g}{\p{1}{Y}}}$.
\end{myremark}

%
%

The definitions of representable and cartesian biclones are now induced from 
their bi-multicategorical 
counterparts~(\cf~Definition~\ref{def:representable-and-cartesian-clones}).

\begin{mydefn} \quad 
\begin{enumerate}
\item A \Def{representable biclone} is a biclone $(S, \biclone)$ equipped with 
a choice of representable structure $\tensor_n(-)$ on $\mCatOf\biclone$.
\item A \Def{cartesian biclone} is a biclone $(S, \biclone)$
equipped with a choice of cartesian structure $\prodop_n(-)$ on 
$\mCatOf\biclone$.
\qedhere
\end{enumerate}
\end{mydefn}

\begin{myremark} \label{rem:biclone-unary-products-are-trivial}
As for fp-bicategories, we stipulate that the unary product structure in a 
cartesian biclone is the identity~(\cf~Remark~\ref{rem:unary-products}).
\end{myremark}

For a clone $(S, \clone)$, the mapping $\cslr{(-)}{h}$ composing with a single 
multimap $h : {X_1, \,\dots\, , X_n \to R}$ is equal to the mapping 
$\ms{(-)}{h}$ 
performing the same composition in $\mCatOf\clone$, since for any 
$g : R \to A$ one has
$\ms{g}{h} 	\overset{\text{def.}}{=} \cslr{g}{\cslr{h}{\p{1}{\ind{X}}, 
\,\dots\, , 
\p{n}{\ind{X}}}} 
			= \cslr{g}{h}$. 
In the world of biclones, however, the functors $\cslr{(-)}{h}$ and 
$\ms{(-)}{h}$ 
are related by a structural 
isomorphism~(\cf~Remark~\ref{rem:biclone-to-unary}). Since  
$(\mCatOf{\biclone})(\Gamma; A) = \biclone(\Gamma;A)$
for every $\Gamma$ and $A$, a choice 
of adjoint equivalence
$ \psi_{\ind{X}} : 
	(\mCatOf{\biclone})(X_1, \,\dots\, , X_n; A) \leftrightarrows 
	(\mCatOf{\biclone})(R; A) 
	: \ms{(-)}{h}$
is equivalently a choice of adjoint equivalence
$ \psi_{\ind{X}}' : 
	\biclone(X_1, \,\dots\, , X_n; A) \leftrightarrows 
	\biclone(R; A) 
	: \cslr{(-)}{h}$. 
(To see this, apply the fact that for 
any morphisms 
$f : X \to Y$ and $g, g' : Y \to X$ in a 2-category, if 
$g \iso g'$ then
$f$ and $g$ are the 1-cells of an equivalence $X \simeq Y$ if and only if 
$f$ and $g'$ are the 1-cells of such an equivalence.) 

It follows that a representable biclone
$(S, \biclone, \tensor_n)$
is equivalently a biclone 
$(S, \biclone)$ 
equipped with a choice of object $\tensor_n(X_1, \,\dots\, , X_n)$ 
and multimap 
$\rho_{\ind{X}} : X_1, \,\dots\, , X_n \to \tensor_n(X_1, \,\dots\, , X_n)$ 
for every
${X_1, \,\dots\, , X_n \in S} \:\: (n \in \Nat)$, 
together with a choice of adjoint equivalence
$\biclone(X_1, \,\dots\, , X_n; A) \simeq 
	\biclone\left(\tensor_n(X_1, \,\dots\, , X_n); A \right)$
induced by pre-composing with $\rho_{\ind{X}}$ for every $A \in S$. Explicitly, 
this entails that for every 
$t : X_1, \,\dots\, , X_n \to A$ 
there exists a chosen multimap
$\psi_{\ind{X}}(t) : \tensor_n(X_1, \,\dots\, , X_n) \to A$
and a 2-cell 
$\epsilon_{\ind{X}; f} : \cslr{\psi_{\ind{X}}(f)}{\rho_{\ind{X}}} \To f$, 
universal in the sense that for any 
$g : \tensor_n(X_1, \,\dots\, , X_n)\to A$ and 
$\sigma : \cslr{g}{\rho_{\ind{X}}} \To f$ there exists a unique 2-cell 
$\trans{\sigma} : g \To \psi_{\ind{X}}(f)$ such that
\begin{equation} \label{eq:birepresentability-biclone-ump}
\begin{tikzcd}
\cslr{g}{\rho_{\ind{X}}} \arrow{rr}{\cslr{\trans{\sigma}}{\rho_{\ind{X}}}} 
\arrow[swap]{dr}{\sigma} &
\: &
\cslr{\psi_{\ind{X}}(f)}{\rho_{\ind{X}}} \arrow{dl}{\epsilon_{\ind{X}; f}} \\
\: &
f &
\:
\end{tikzcd} 
\end{equation}

A similar story holds for cartesian biclones. For a sequence of multimaps
${(\pi_i : R \to X_i)_{i=1, \dot, n}}$ 
and $u : \Gamma \to A_i$ 
in the bi-multicategory $\mCatOf\biclone$ associated to a cartesian biclone
$\cartClone{S}{\biclone}$,
there exists the following composite of 
structural isomorphisms:
\[
\ms{\pi_i}{u} = 
	\cslr{\pi_i}{\cslr{u}{\p{1}{\Gamma}, \,\dots\, , \p{\len{\Gamma}}{\Gamma}}}
	\iso \csthree{\pi_i}{u}{\p{1}{\Gamma}, \,\dots\, , \p{\len{\Gamma}}{\Gamma}}
	\iso \cslr{\pi_i}{u}
\]
It follows that the functor
$\left( \ms{\pi_1}{-}, \,\dots\, , \ms{\pi_n}{-} \right) 
	: (\mCatOf{\biclone})(\Gamma; R) \to 
		\prodop_{i=1}^n (\mCatOf{\biclone})(\Gamma; X_i)$ 
is naturally isomorphic to the functor
$\left( \cslr{\pi_1}{-}, \,\dots\, , \cslr{\pi_n}{-} \right) 
	: \biclone(\Gamma; R) \to 
		\prodop_{i=1}^n \biclone(\Gamma; X_i)$.
A cartesian biclone 
$\cartClone{S}{\biclone}$ is therefore
equivalently a biclone equipped with a choice of object
$\prodop_n(X_1, \,\dots\, , X_n)$ 
and multimaps
$\left(\pi_i : \prodop_n(X_1, \,\dots\, , X_n) \to X_i\right)_{i=1, \,\dots\, 
,n}$
for every sequence 
${X_1, \,\dots\, , X_n \in S} \:\: {(n \in \Nat)}$, 
together with a choice of adjoint equivalence
$\biclone{\left(\Gamma; \prodop_n(X_1,\dots, X_n) \right)}
	\simeq \prodop_{i=1}^n \biclone(\Gamma; X_i)$. 
The counit of this adjoint equivalence is then characterised by the following universal property. For 
every sequence of multimaps 
${(t_i : \Gamma \to X_i)_{i=1, \,\dots\, , n}}$
there exists a multimap 
$\pair{t_1, \,\dots\, , t_n} : \Gamma \to \prodop_n(X_1, \,\dots\, , X_n)$
and a 2-cell
$\epsilonTimes{}$ 
with components
$\epsilonTimesInd{i}{\ind{t}} :
	\cslr{\pi_i}{\pair{t_1, \,\dots\, , t_n}} \To t_i$
for $i=1, \,\dots\, , n$. This 2-cell is universal in the sense that, if 
$u : \Gamma \to \prodop_n(X_1, \,\dots\, , X_n)$
and
$\alpha_i : \cslr{\pi_i}{u} \To t_i$
for $i=1, \,\dots\, , n$, then there exists a unique 2-cell
$\transTimes{\alpha_1, \,\dots\, , \alpha_n} : 
	u \To \pair{t_1, \,\dots\, , t_n}$
filling the following diagram for $i=1, \,\dots\, , n$:
\begin{equation} \label{eq:cartesian-biclone-ump}
\begin{tikzcd}
\cslr{\pi_i}{u} 
\arrow[swap]{dr}{\alpha_i}
\arrow{rr}{\cslr{\pi_i}{\alpha}} &
\: &
\cslr{\pi_i}{\pair{t_1, \,\dots\, , t_n}}
\arrow{dl}{\epsilonTimesInd{i}{\ind{t}}} \\
\: &
t_i &
\:
\end{tikzcd}
\end{equation} 
Rather than translating between compositions $\ms{f}{\ind{g}}$ and 
$\cslr{f}{\ind{g}}$ throughout, in what follows we employ the biclone version 
of the universal property.

\begin{myremark}
\label{rem:universal-arrow-biclone}
We have just shown that a \Def{biuniversal arrow 
in a biclone}---defined exactly as in 
Definition~\ref{def:universal-arrow-bi-multicat}---exists if and only if there 
exists a biuniversal arrow in the corresponding bi-multicategory. 
\end{myremark}

\begin{myexmp} \label{ex:fp-bicategory-defines-a-cartesian-biclone}
Every fp-bicategory $\fpBicat\baseCat$ 
defines a biclone 
$\bicloneFromProducts{\baseCat}$ with sorts $ob(\baseCat)$ and hom-categories
$\bicloneFromProducts{\baseCat}(X_1, \dots, X_n; Y)
	:= \baseCat{\left(\prodop_n(X_1, \dots, X_n), Y\right)}$~(\cf~Example~\ref{ex:fp-cat-to-cartesian-clone} on 
	page~\pageref{ex:fp-cat-to-cartesian-clone}).
The substitution $\cs{f}{g_1, \dots, g_n}$ is
$f \circ \seq{g_1, \dots, g_n}$.
This biclone is 
cartesian: for the adjoint 
equivalence~(\ref{eq:cartesian-biclone-ump}) one takes the adjoint 
equivalence defining finite products in $\baseCat$.  
\end{myexmp}

\paragraph*{The equivalence between representability and cartesian structure.}
Our aim now is to prove a version of 
Theorem~\ref{thm:representable-clone-cartesian-equivalence} for 
biclones, establishing that a biclone admits a representable structure
(embodied by~(\ref{eq:birepresentability-biclone-ump})) if and only if it 
admits a cartesian structure (embodied by~(\ref{eq:cartesian-biclone-ump})). In 
the 1-categorical case the key to this equivalence is the 
construction of a 
sequence of multimaps
$\pi_i : \tensor_n(X_1, \,\dots\, , X_n) \to X_i$
satisfying two equations
for $i=1, \,\dots\, , n$. 
The corresponding bicategorical construction is 
up-to-isomorphism.

\begin{mylemma} \label{lem:bicat-product-structure-from-birepresentability}
For any representable biclone $(S, \biclone, \tensor_n)$ and 
$X_1, \,\dots\, , X_n \in S \:\: (n \in \Nat)$ 
there exist multimaps
$\pi_i : \tensor_n(X_1, \,\dots\, , X_n) \to X_i$
and invertible 2-cells
$\counitCell{i}{\ind{X}} : \cs{\pi_i}{\rho_{\ind{X}}} \To \p{i}{\ind{X}}$ 
and
$\etaTimes{\ind{X}} : 
	\Id_{\tensor_n(X_1, \dots, X_n)}
		\To 
	\cs{\rho_{\ind{X}}}{\pi_1, \dots, \pi_n}$
(for $i=1, \,\dots\, , n$), as in the diagrams below:
\begin{center}
\begin{small}
\begin{tikzcd}[column sep = 1em]
\: &
\tensor_n(X_1, \,\dots\, , X_n) \arrow{dr}{\pi_i} 
\arrow[phantom]{d}[description, font=\scriptsize]
	{\twocellDown{\counitCell{i}{\ind{X}}}} &
\: \\
X_1, \,\dots\, , X_n 
\arrow[swap]{rr}{\p{i}{\ind{X}}} \arrow{ur}{\rho_{\ind{X}}} &
\: &
X_i
\end{tikzcd}
\end{small}
\:\:
\begin{small}
\begin{tikzcd}[column sep = 1em]
\: &
X_1, \,\dots\, , X_n \arrow{dr}{\rho_{\ind{X}}} &
\: \\
\tensor_n(X_1, \,\dots\, , X_n) \arrow[swap]{rr}{\Id} \arrow{ur}{[\pi_1, 
\,\dots\, , 
\pi_n]} &
\: \arrow[phantom]{u}[description, font=\scriptsize]
	{\twocellUp{\etaTimes{\ind{X}}}} &
\tensor_n(X_1, \,\dots\, , X_n)
\end{tikzcd}
\end{small}
\end{center}
\begin{proof}
Define $\pi_i := \psi_{\ind{X}}(\p{i}{\ind{X}})$. 
For $\counitCell{i}{\ind{X}}$, 
we may immediately take the universal 2-cell  
$\epsilon_{\ind{X}; \p{i}{}}$ of~(\ref{eq:birepresentability-biclone-ump}). 
%
For $\etaTimes{\ind{X}}$ we apply the universal 
property~(\ref{eq:birepresentability-biclone-ump}) to the structural isomorphism
$\indproj{1}{(\tensor_n\ind{X})}$
to obtain an invertible 2-cell 
$\trans{(\indproj{1}{\ind{X}})} : 
	\Id_{(\tensor_n \ind{X})} \To 
	\psi_{\ind{X}}(\rho_{\ind{X}})$. 
We complete the construction by defining a 2-cell
$\cslr{\rho_{\ind{X}}}{\pi_1, \,\dots\, , \pi_n} \To 
\psi_{\ind{X}}(\rho_{\ind{X}})$. 
Define 
$\alpha_{\ind{X}}$
to be the composite 
\[
\csthree{\rho_{\ind{X}}}
		{\pi_1, \,\dots\, , \pi_n}
		{\rho_{\ind{X}}}
\XRA{\iso}
\cslr
	{\rho_{\ind{X}}}
	{\cslr{\ind{\pi}}{\rho_{\ind{X}}}} 
\XRA{\cslr{\rho_{\ind{X}}}{\counitCell{\bullet}{\ind{X}}}} 
\cslr
	{\rho_{\ind{X}}}
	{\p{1}{\ind{X}}, \,\dots\, , \p{n}{\ind{X}}} 
\XRA{\subid{}^{-1}} 
\rho_{\ind{X}}
\] 
Since this composite is invertible, by the universal 
property~(\ref{eq:birepresentability-biclone-ump}) there exists an invertible 
2-cell 
$\trans{(\alpha_{\ind{X}})} : 
	\cslr{\rho_{\ind{X}}}{\pi_1, \,\dots\, , \pi_n} 
	\To 
	\psi_{\ind{X}}(\rho_{\ind{X}})$. 
We therefore define $\etaTimes{\ind{X}}$ to be the composite
\[
\Id_{(\tensor \ind{X})} 
	\XRA{\trans{\indproj{1}{\ind{X}}}} 
\psi_{\ind{X}}(\rho_{\ind{X}}) 
	\XRA{(\trans{\alpha}_{\ind{X}})^{-1}} 
\cslr{\rho_{\ind{X}}}{\pi_1, \,\dots\, , \pi_n}
\]
\end{proof}
\end{mylemma}

To bicategorify Lemma~\ref{lem:clone-representable-iff-cartesian} we shall also 
employ a kind of `mirror image' of the preceding lemma, capturing the crucial 
construction available in the presence of cartesian structure; this should be 
compared to the discussion preceding 
Definition~\ref{def:clone-1-cell-invetibility} (page~\pageref{def:clone-1-cell-invetibility}).
Just as we had to generalise the notion of isomorphism for the clone case, so we need to generalise the notion of (adjoint) equivalence for the biclone case.

\newpage
\begin{mydefn} \label{def:multicat-adjoint-equivalence}
Let $(S, \biclone)$ be a biclone. 
\begin{enumerate}
\item An \Def{adjunction} 
$X_1 \dots, X_n \leftrightarrows Y$ in 
$(S, \biclone)$ 
consists of 1-cells $e : X_1, \,\dots\, , X_n \to Y$ and 
${f_i : Y \to X_i} \:\: (i = 1, \,\dots\, , n)$ with 2-cells 
\begin{align*}
\eta : \p{1}{Y} &\To \cslr{e}{f_1, \,\dots\, , f_n} : Y \to Y \\
\epsilon_i : \cslr{f_i}{e} &\To \p{i}{X_1, \,\dots\, , X_n} : X_1, \,\dots\, , 
X_n \to 
X_i 
\quad (i = 1, \,\dots\, , n)
\end{align*}
such that the following diagrams commute for $i=1, \,\dots\, ,n$:
\begin{center}
\vspace{-2\parskip}
\begin{minipage}{0.45\textwidth}
\begin{equation} \label{eq:triangle-law-one}
\begin{tikzcd}[column sep = small] 
\cslr{\p{1}{Y}}{e} 
\arrow[swap]{d}{\indproj{1}{e}} 
\arrow{r}{\cslr{\eta}{e}} & 
\cslr{\cslr{e}{\ind{f}}}{e} 
\arrow{r}{\assoc{e; \ind{f};e}} & 
\cslr{e}{\cslr{\ind{f}}{e}} 
\arrow{d}{\cslr{e}{\epsilon_1, \,\dots\, , \epsilon_n}} 
\\ 
e \arrow[swap]{rr}{\subid{e}} & 
\: & 
\cslr{e}{\p{1}{\ind{X}}, \,\dots\, , \p{n}{\ind{X}}} & 
\end{tikzcd} 
\vspace{2mm}
\end{equation}
\end{minipage}
\quad
\begin{minipage}{0.45\textwidth}
\begin{equation} \label{eq:triangle-law-two}
\begin{tikzcd}[column sep = small] 
f_i 
\arrow{r}{\subid{{f_i}}} 
\arrow[swap, equals]{d}{} & 
\cslr{f_i}{\p{1}{Y}} 
\arrow{r}{\cslr{f_i}{\eta}} & 
\cslr{f_i}{\cslr{e}{f_1, \,\dots\, , f_n}}
\arrow{d}{\assoc{f_i;e;\ind{f}}^{-1}} \\ 
f_i & 
\cslr{\p{i}{}}{f_1, \,\dots\, , f_n} 
\arrow{l}{\indproj{i}{\ind{f}}} & 
\cslr{\cslr{f_i}{e}}{f_1, \,\dots\, , f_n} 
\arrow{l}[yshift=-2mm]
	{\cslr{\epsilon_i}{f_1, \,\dots\, , f_n}} 
\end{tikzcd} 
\end{equation}
\end{minipage} 
\vspace{-\parskip}
\end{center}
\item An \Def{equivalence} in $(S, \biclone)$ consists of 1-cells 
$e : X_1, \,\dots\, , X_n \to Y$ and 
$f_i : Y \to X_i \: (i = 1, \,\dots\, , n)$ with 
invertible 2-cells 
\begin{align*}
\eta : \p{1}{Y} &\XRA{\iso} \cslr{e}{f_1, \,\dots\, , f_n} : Y \to Y \\ 
\epsilon_i : \cslr{f_i}{e} &\XRA{\iso} \p{i}{X_1, \,\dots\, , X_n} : 
	X_1, \,\dots\, , X_n \to X_i \:\: (i = 1, \,\dots\, , n)
\end{align*}
\item A \Def{adjoint equivalence} in $(S, \biclone)$ is an adjunction for which 
$\eta$ and $\epsilon_i$ are invertible for $i = 1, \,\dots\, , n$.  \qedhere
\end{enumerate}
\end{mydefn}  

In particular, a unary (adjoint) equivalence $X \simeq Y$ is just an (adjoint) 
equivalence in the usual, bicategorical sense. 

\begin{mylemma} \label{lem:biclone-with-products-adjoint-eq}
For any sequence of objects 
$X_1, \,\dots\, , X_n \:\: (n \in \Nat)$ in a cartesian biclone 
$\cartClone{S}{\biclone}$, there exists an adjoint equivalence 
between $X_1, \,\dots\, , X_n \simeq \prodop_n(X_1, \,\dots\, , X_n)$.
\begin{proof}
We employ the notation of~(\ref{eq:cartesian-biclone-ump}) for cartesian 
structure.
For the 2-cell
$\cs{\pi_i}{\pairName(\p{1}{\ind{X}}, \,\dots\, , \p{n}{\ind{X}})} 
	\To \p{i}{\ind{X}}$ 
we can immediately take $\epsilonTimesInd{i}{\ind{X}}$. The real work is in 
providing a 2-cell  
$\gamma : \Id_{(\prod\ind{X})} \To 
	\cslr{\pairName(\p{1}{}, \,\dots\, , \p{n}{})}{\pi_1, \,\dots\, , \pi_n}$. 
By the 	universality of the counit 
$\epsilonTimes{} = (\epsilonTimesInd{1}{}, \,\dots\, , \epsilonTimesInd{n}{})$ 
it 
suffices to define a 
family of 
invertible 2-cells 
$\zeta_i : 
	\cslr	{\pi_i}
		{\cslr{\pairName(\p{1}{}, \,\dots\, , \p{n}{})}
			{\pi_1, \,\dots\, ,\pi_n}} 
	\To \pi_i$
for $i=1, \,\dots\, , n$. 
We may then define $\gamma$ to be the composite
\[
\Id_{(\prod \ind{X})} \XRA{\etaTimes{\Id_{(\prod \ind{X})}}}  
\pairName(\cslr{\ind{\pi}}{\Id_{(\prod \ind{X})}})  
\XRA{\pairName(\subid{}^{-1}, 
\dots, \subid{}^{-1})} \pairName(\ind{\pi}) \XRA{(\transTimes{\zeta_1, 
\,\dots\, , 
\zeta_n})^{-1}} \cslr{\pairName(\p{\bullet}{})}{\ind{\pi}} 
\]
where $\etaTimes{}$ is the unit of the adjoint equivalence witnessing $(\pi_1, 
\dots,\pi_n)$ as a biuniversal arrow. The 2-cells $\zeta_i$ are defined as 
follows:
\[
\cslr{\pi_i}{\cslr{\pairName(\p{1}{},\dots, \p{n}{})}{\ind{\pi}}} 
\XRA{\assoc{}^{-1}} \csthree{\pi_i}{\pairName(\p{1}{}, \,\dots\, , 
\p{n}{})}{\ind{\pi}} \XRA{\cslr{\epsilonTimesInd{i}{\ind{X}}}{\ind{\pi}}} 
\cslr{\p{i}{}}{\ind{\pi}} \XRA{\indproj{i}{}} \pi_i
\]
Since each $\zeta_i$ is 
invertible, $\transTimes{\zeta_1, \,\dots\, , \zeta_n}$ is also invertible. 
Checking 
that 
diagram~(\ref{eq:triangle-law-two}) commutes is 
straightforward; 
for~(\ref{eq:triangle-law-one}) one must use the universal property, checking 
that both routes around the diagram are the unique 2-cell corresponding to the 
composite
\[
\cslr
	{\pi_i}
	{\csthreesmall
		{\pairName(\p{1}{}, \,\dots\, , \p{n}{})}
		{\ind{\pi}}
		{\pairName(\p{1}{}, \,\dots\, , \p{n}{})}} 
\XRA{\cslr{\pi_i}{\ind{\beta}}} \cslr{\pi_i}{\pairName(\p{1}{}, \,\dots\, , 
\p{n}{})} 
\XRA{\epsilonTimesInd{i}{\ind{X}}} \p{i}{}
\]
where $\beta_i$ is defined to be
\small
\[
\csthree{\pairName(\p{\bullet}{})}{\ind{\pi}}{\pairName(\p{\bullet}{})} 
\XRA{\assoc{}} 
\cslr{\pairName(\p{\bullet}{})}{\cslr{\ind{\pi}}{\pairName(\p{\bullet}{})}} 
\XRA{\cslr{\pairName(\p{\bullet}{})}{\epsilonTimesInd{\bullet}{\ind{X}}}} 
\cslr{\pairName(\p{\bullet}{})}{\p{\bullet}{}} \XRA{\subid{}^{-1}} 
\pairName(\p{\bullet}{})
\] 
\normalsize
for $i=1,\dots,n$. 
\end{proof}
\end{mylemma}

As for clones, the extra structure of a biclone entails that birepresentable arrows are closed under composition. The strategy for the proof is 
familiar from Lemma~\ref{lem:multicat-from-clone-universal-arrows-compose}.

\begin{mylemma}
A biclone $(S, \biclone)$ admits a representable structure if and only if for 
every 
$X_1, \,\dots\, , X_n \in \mbCat \:\: ({n \in \Nat})$ there exists a chosen 
object 
$\tensor_n(X_1, \,\dots\, , X_n) \in \mbCat$ 
and a birepresentable multimap 
$\rho_{\ind{X}} : X_1, \,\dots\, , X_n \to \tensor_n(X_1, \,\dots\, , X_n)$.
\begin{proof}
It suffices to show that birepresentable multimaps are closed under 
composition. Mirroring the proof of 
Lemma~\ref{lem:multicat-from-clone-universal-arrows-compose}, suppose given 
birepresentable multimaps
\begin{gather*}
\rho_{\ind{X}} : X_1, \,\dots\, , X_n \to \tensor_n(X_1, \,\dots\, , X_n) \\ 
\rho_{\ind{Y}} : Y_1, \,\dots\, , Y_m \to \tensor_m(Y_1, \,\dots\, , Y_m) \\
\rho_{(\prod \ind{X}, \prod \ind{Y})} : \tensor_n \ind{X}, \tensor_m \ind{Y} 
\to \tensor_2(\tensor_n \ind{X}, \tensor_m \ind{Y})
\end{gather*}
We want to show that the composite $\rho_{(\prod \ind{X}, \prod \ind{Y})} \circ 
(\rho_{\ind{X}}, \rho_{\ind{Y}})$ in $\mCatOf{\biclone}$, which is the 
composite $\overline{\rho} := \cslr{\rho_{(\prod \ind{X}, \prod 
\ind{Y})}}{\cslr{\rho_{\ind{X}}}{\p{1}{}, \,\dots\, , \p{n}{}}, 
\cslr{\rho_{\ind{Y}}}{\p{n+1}{}, \,\dots\, , \p{n+m}{}}}$ in $\biclone$, is 
birepresentable. Define projections $\pi_i^X : \tensor_n(X_1, \,\dots\, , X_n) 
\to 
X_i$, $\pi_j^Y : \tensor_m(Y_1, \,\dots\, , Y_m) \to Y_j$ and $\pi^{X,Y}$ as in 
the 
proof of Lemma~\ref{lem:multicat-from-clone-universal-arrows-compose}, and 
likewise define a family of multimaps $\overline{\pi}_i : \tensor_{2}(\tensor_n 
\ind{X}, 
\tensor_m \ind{Y}) \to Z_i$ for $i=1, \,\dots\, , n+m$ (where $Z_i$ is $X_i$ 
for $1 
\leq i \leq n$ and $Y_{i-n}$ for $n+1 \leq i \leq n+m$) as 
in~(\ref{eq:def-of-ri-maps}).  Finally, for $1 \leq i \leq n$ define an 
invertible 2-cell $\beta^{(1)} : \cslr{\rho_{\ind{X}}}{\overline{\pi}_1, 
\,\dots\, , 
\overline{\pi}_n} \To 
\pi_1^{X, Y} : \tensor_2( \tensor_n \ind{X}, \tensor_m \ind{X} ) \to \tensor_n 
\ind{X}$ by
\begin{td}
\cslr{\rho_{\ind{X}}}{\overline{\pi}_1, \,\dots\, , \overline{\pi}_n} 
\arrow[equals]{d} 
\arrow{r}{\beta^{(1)}} &
\pi_1^{X,Y} \\

\cslr{\rho_{\ind{X}}}{\cslr{\pi_1^X}{\pi_1^{X,Y}}, \,\dots\, , 
\cslr{\pi_n^X}{\pi_1^{X,Y}}} 
\arrow[swap]{d}{\assoc{\rho_{\ind{X}}; \ind{\pi}; \pi_1}^{-1}} &
\: \\

\csthree{\rho_{\ind{X}}}{\pi_1^X, \,\dots\, , \pi_n^X}{\pi_1^{X,Y}} 
\arrow[swap]{r}{\cslr{\etaTimes{\ind{X}}^{-1}}{\pi_1^{X,Y}}} &
\cslr{\Id_{(\tensor \ind{X}}}{\pi_1^{X,Y})} 
\arrow[swap]{uu}{\indproj{1}{\pi_1}}
\end{td}
We define $\beta^{(2)} : \cslr{\rho_{\ind{Y}}}{r_{n+1}, \,\dots\, , 
\overline{\pi}_{n+m}} \To \pi_2^{X, Y} : \tensor_2( \tensor_n \ind{X}, 
\tensor_m 
\ind{X} ) \to \tensor_m \ind{Y}$ similarly. 

We are now in a position to define the pseudo-inverse to 
$\ms{(-)}{\overline{\rho}} : \mbCat\big( \tensor_2(\tensor_n \ind{X}, \tensor_m 
\ind{Y}); A \big) \to \mbCat(X_1, \,\dots\, , X_n, Y_1, \,\dots\, , Y_m; A)$. 
For $h : 
X_1, \,\dots\, , X_n, Y_1, \,\dots\, , Y_m \to A$ we define 
$\overline{\psi}(h)$ to be 
the 
composite 
\[\tensor_2(\tensor_n \ind{X}, \tensor_m \ind{Y}) 
\xra{[\overline{\pi}_1, \,\dots\, , 
\overline{\pi}_{n+m}]} X_1, \,\dots\, , X_n, Y_1, \,\dots\, , Y_m \xra{h} A
\] in 
$\biclone$; this 
mapping is clearly functorial. It therefore suffices to construct natural 
isomorphisms $\id_{\mbCat( \tensor(\tensor \ind{X}, \tensor \ind{Y}); A )} \iso 
\overline{\psi}\big( \ms{(-)}{\overline{\rho}} \big)$ and $\id_{\mbCat(X_1, 
\dots, X_n, Y_1, \,\dots\, , Y_m; A)} \iso 
\ms{\big(\overline{\psi}(-)\big)}{\overline{\rho}}$; this lifts to an adjoint 
equivalence between the same 1-cells by the usual well-known argument 
(\eg~\cite[IV.3]{cfwm}). 

To this end, let us define invertible 2-cells $\tau$ and 
$\sigma_i \: (i=1, \,\dots\, , n+m)$ that will make up the bulk of the required 
isomorphisms. 
The 2-cell $\tau$ is defined as follows:
\begin{small}
\begin{td}[column sep = 1.5em]
\csthree{\rho_{(\tensor \ind{X}, \tensor 
\ind{Y})}}{\cslr{\rho_{\ind{X}}}{\p{1}{}, 
\dots, \p{n}{}}, \cslr{\rho_{\ind{Y}}}{\p{n+1}{}, \,\dots\, , 
\p{n+m}{}}}{\overline{\pi}_1, \,\dots\, , \overline{\pi}_{n+m}} 
\arrow[swap]{d}{\iso} 
\arrow{r}{\tau} &
\Id_{\tensor (\tensor \ind{X}, \tensor \ind{Y})} \\

\cslr{\rho_{(\tensor \ind{X}, \tensor 
\ind{Y})}}{\cslr{\rho_{\ind{X}}}{\cslr{\p{\bullet}{}}{\ind{\overline{\pi}}}}, 
\cslr{\rho_{\ind{Y}}}{\cslr{\p{\bullet}{}}{\ind{\overline{\pi}}}}} 
\arrow[swap]{d}{\cslr{\rho_{(\tensor \ind{X}, \tensor 
\ind{Y})}}{\cslr{\rho_{\ind{X}}}{\indproj{\bullet}{}}, 
\cslr{\rho_{\ind{Y}}}{\indproj{\bullet}{}}}} &
\: \\

\cslr{\rho_{(\tensor \ind{X}, \tensor 
\ind{Y})}}{\cslr{\rho_{\ind{X}}}{\overline{\pi}_1, 
\dots, \overline{\pi}_n}, \cslr{\rho_{\ind{Y}}}{\overline{\pi}_{n+1}, \,\dots\, 
, 
\overline{\pi}_{n+m}}} \arrow[swap]{r}[yshift=-2mm]{\cslr{\rho_{(\tensor 
\ind{X}, 
\tensor \ind{Y})}}{ \beta^{(1)}, \beta^{(2)}}} &
\cslr{\rho_{(\tensor \ind{X}, \tensor \ind{Y})}}{\pi_1^{X,Y}, \pi_2^{X,Y}} 
\arrow[swap]{uu}{\etaTimes{(\tensor \ind{X}, \tensor \ind{Y})}^{-1}}
\end{td}
\end{small}
The 2-cells $\sigma_1, \,\dots\, , \sigma_n$, on the other hand, are defined by 
the 
following 
diagram; the definitions of $\sigma_{n+1}, \,\dots\, , \sigma_{n+m}$  are the same, modulo the obvious adjustments.
\begin{small}
\begin{td}
\cslr{\overline{\pi}_i}{\cslr{\rho_{(\tensor \ind{X}, \tensor 
\ind{Y})}}{\cslr{\rho_{\ind{X}}}{\p{1}{}, \,\dots\, , \p{n}{}}, 
\cslr{\rho_{\ind{Y}}}{\p{n+1}{}, \,\dots\, , \p{n+m}{}}}} \arrow[equals]{d} 
\arrow[]{r}{\sigma_i} &
\p{i}{X_1, \,\dots\, , X_n, Y_1, \,\dots\, , Y_m} \\

\cslr{\cslr{\pi_i^X}{\pi_1^{X,Y}}}{\cslr{\rho_{(\tensor \ind{X}, \tensor 
\ind{Y})}}{\cslr{\rho_{\ind{X}}}{\p{1}{}, \,\dots\, , \p{n}{}}, 
\cslr{\rho_{\ind{Y}}}{\p{n+1}{}, \,\dots\, , \p{n+m}{}}}} \arrow[swap]{d}{\iso} 
&
\: \\

\csthree{\pi_i^X}{\cslr{\pi_1^X}{\rho_{(\tensor \ind{X}, \tensor 
\ind{Y})}}}{\cslr{\rho_{\ind{X}}}{\p{1}{}, \,\dots\, , \p{n}{}}, 
\cslr{\rho_{\ind{Y}}}{\p{n+1}{}, \,\dots\, , \p{n+m}{}}} 
\arrow[swap]{d}{\csthree{\pi_i^X}{\counitCell{1}{\tensor\ind{X}, 
\tensor\ind{Y}}}{\cslr{\rho_{\ind{X}}}{\p{\bullet}{}}, 
\cslr{\rho_{\ind{Y}}}{\p{\bullet}{}}}}  &
\: \\

\csthree{\pi_i^X}{\p{1}{\ind{X}}}{\cslr{\rho_{\ind{X}}}{\p{1}{}, \,\dots\, , 
\p{n}{}}, 
\cslr{\rho_{\ind{Y}}}{\p{n+1}{}, \,\dots\, , \p{n+m}{}}} \arrow[swap]{d}{\iso} &
\: \\

\csthree{\pi_i^X}{\rho_{\ind{X}}}{\p{1}{}, \,\dots\, , \p{n}{}} 
\arrow[swap]{r}{\cslr{\counitCell{i}{\ind{X}}}{\p{1}{}, \,\dots\, , \p{n}{}}} &
\cslr{\p{i}{}}{\p{1}{}, \,\dots\, , \p{n}{}} 
\arrow[swap]{uuuu}{\indproj{i}{\p{\bullet}{}}}
\end{td}
\end{small}
The required natural isomorphisms are then defined to be the composites
\begin{gather*}
\ms{\overline{\psi}(g)}{\overline{\rho}} = \csthree{g}{\overline{\pi}_1, 
\,\dots\, , 
\overline{\pi}_{n+m}}{\overline{\rho}} \XRA{\assoc{}} 
\cslr{g}{\cslr{\ind{r}}{\overline{\rho}}} \XRA{\cslr{g}{\ind{\sigma}}} 
\cslr{g}{\p{1}{},\dots, \p{n+m}{}} \XRA{\subid{}^{-1}} g \\
\overline{\psi}(\ms{h}{\overline{\rho}}) = 
\csthree{h}{\overline{\rho}}{\overline{\pi}_1, \,\dots\, , 
\overline{\pi}_{n+m}} 
\XRA{\assoc{}} \cslr{h}{\cslr{\overline{\rho}}{\overline{\pi}_1, \,\dots\, , 
\overline{\pi}_{n+m}}} \XRA{\cslr{h}{\tau}} \cslr{h}{\Id_{\tensor 
(\tensor\ind{X}, 
\tensor\ind{Y})}} \XRA{\subid{}^{-1}} h 
\end{gather*}
for $g : \tensor_2(\tensor_n\ind{X}, \tensor_m\ind{Y}) \to A$ and $h : X_1, 
\dots, X_n, Y_1, \,\dots\, , Y_m \to A$. 
\end{proof}
\end{mylemma}

We now prove the central result of this section.

\begin{mylemma}
A biclone $(S, \biclone)$ admits a choice of representable structure if and 
only if it admits a choice of cartesian structure. 
\begin{proof}
\subproof{$\To$} Let $\rho_{\ind{X}} : X_1, \,\dots\, , 
X_n \to \tensor_n(X_1, \,\dots\, , X_n)$ be a birepresentable multimap. We 
claim the 
sequence of multimaps 
$(\pi_i : \tensor_n(X_1,\dots, X_n) \to X_i)_{i=1,\dots,n}$ 
defined in Lemma~\ref{lem:bicat-product-structure-from-birepresentability}
form a biuniversal multimap. We are therefore required to provide a mapping 
$\pairName : 
\prod_{i=1}^n \mbCat(\Gamma; X_i) \to 
	\mbCat(\Gamma; \tensor_n(X_1, \,\dots\, , X_n)\big)$ 
and a universal 2-cell with	components 
$\epsilonTimesInd{i}{\ind{X}} 
: \cslr{\pi_i}{\pairName(f_1, \,\dots\, , f_n)} \To f_i$ 
for $i=1,\dots,n$. We define 
$\pairName(f_1, \,\dots\, , f_n) := \cslr{\rho_{\ind{X}}}{f_1, \,\dots\, , 
f_n}$ and set 
$\epsilonTimesInd{i}{\ind{X}}$ to be the composite
\[
\cslr{\pi_i}{\cslr{\rho_{\ind{X}}}{f_1, \,\dots\, , f_n}} \XRA{\assoc{}^{-1}} 
\csthree{\pi_i}{\rho_{\ind{X}}}{f_1, \,\dots\, , f_n} 
\XRA{\cslr{\counitCell{i}{\ind{X}}}{\ind{f}}} \cslr{\p{i}{}}{f_1, \,\dots\, , 
f_n} 
\XRA{\indproj{i}{}} f_i 
\]
For universality, suppose $g : \Gamma \to \tensor_n(X_1, \,\dots\, , X_n)$ and 
$\alpha_i : \cslr{\pi_i}{g} \To f_i$ for $i = 1,\dots, n$. We define 2-cell 
$\transTimes{\alpha_1, \,\dots\, , \alpha_n} : g \To \pairName(f_1, \,\dots\, , 
f_n)$ 
by the commutativity of the following diagram:
\begin{equation} \label{eq:translating-universal-arrows}
\begin{tikzcd}[column sep = 6	em]
g 
\arrow[swap]{d}{\indproj{-1}{g}}
\arrow{rr}{\transTimes{\alpha_1, \,\dots\, , \alpha_n}} &
\: &
\cslr{\rho_{\ind{X}}}{f_1, \,\dots\, , f_n}  \\
\cslr{\Id_{(\tensor \ind{X}})}{g} 
\arrow[swap]{r}{\cslr{\etaTimes{\ind{X}}}{g}} &
\csthree{\rho_{\ind{X}}}{\pi_1, \,\dots\, , \pi_n}{g} 
\arrow[swap]{r}{\assoc{\rho_{\ind{X}}; \ind{\pi}; g}} &
\cslr{\rho_{\ind{X}}}{\cslr{\pi_1}{g},\dots, \cslr{\pi_n}{g}} 
\arrow[swap]{u}{\cslr{\rho_{\ind{X}}}{\ind{\alpha}}}
\end{tikzcd}
\end{equation}
where we employ the 2-cell $\etaTimes{\ind{X}}$ defined in 
Lemma~\ref{lem:bicat-product-structure-from-birepresentability}. For the 
existence part of the claim, we need to check that the composite 
\[
\cslr{\pi_i}{g} \XRA{\cslr{\pi_i}{\transTimes{\alpha_1, \,\dots\, , \alpha_n}}} 
\cslr{\pi_i}{\pairName(f_1, \,\dots\, , f_n)} 
\XRA{\epsilonTimesInd{i}{\ind{X}}} f_i
\] 
is equal to $\alpha_i$ for $i=1, \,\dots\, , n$. Most of the calculation is 
straightforward; the key lemma is that the following diagram 
commutes for $i=1,\dots, n$:
\begin{td}[column sep = 4em]
\pi_i 
\arrow[swap]{d}{\subid{\pi_i}} 
\arrow[equals]{r} &
\pi_i \\

\cslr{\pi_i}{\Id_{(\tensor \ind{X})}} 
\arrow[swap]{d}{\cslr{\pi_i}{\etaTimes{\ind{X}}}} &
\: \\

\cslr{\pi_i}{\cslr{\rho_{\ind{X}}}{\pi_1, \,\dots\, , \pi_n}} 
\arrow[swap]{d}{\assoc{\pi_i; \rho_{\ind{X}}; \ind{\pi}}^{-1}} &
\: \\

\csthree{\pi_i}{\rho_{\ind{X}}}{\pi_1, \,\dots\, , \pi_n} 
\arrow[swap]{r}{\cslr{\counitCell{i}{\ind{X}}}{\ind{\pi}}} &
\cslr{\p{i}{}}{\pi_1, \,\dots\, , \pi_n} 
\arrow[swap]{uuu}{\indproj{i}{\ind{\pi}}}
\end{td}
For uniqueness, let 
$g : \Gamma \to \tensor_n(X_1, \,\dots\, , X_n)$
be any multimap and suppose that
$\sigma : g \To \pairName(f_1, \,\dots\, , f_n)$ satisfies 
$\epsilonTimesInd{i}{\ind{X}} \vert \cslr{\pi_i}{\sigma} = \alpha_i$ for 
$i=1,\dots,n$. Substituting this equation into the definition of 
$\transTimes{\alpha_1, \,\dots\, , \alpha_n}$ and using the above diagram, one 
sees 
that $\sigma = \transTimes{\alpha_1, \,\dots\, , \alpha_n}$ as required. 

Finally, it remains to 
check that the unit and counit of the adjunction we have just constructed are 
invertible. The counit is the universal 2-cell, which is certainly invertible. 
The unit is constructed by applying $\transTimes{-,\dots, =}$ to the identity, 
which is invertible since it is a composite of invertible 2-cells.

\subproof{$\Leftarrow$} For the converse, we claim that 
$\rho_{\ind{X}} := \pairName(\p{1}{\ind{X}}, \,\dots\, , \p{n}{\ind{X}}) : 
X_1, \,\dots\, , X_n 
\to \prodop_n(X_1, \,\dots\, , X_n)$ is birepresentable. We therefore need to 
supply 
a mapping $\psi_{\ind{X}} : (\mCatOf\biclone)(X_1, \,\dots\, , X_n; A) \to 
(\mCatOf\biclone)\big( \prodop_n(X_1, \,\dots\, , X_n); A\big)$ and a universal 
2-cell $\epsilon_{A, g} : \cslr{\psi_{\ind{X}}(g)}{\rho_{\ind{X}}} \To g$. We 
define $\psi_{\ind{X}}(g) := \cslr{g}{\pi_1, \,\dots\, , \pi_n}$ and set 
$\epsilon_{A, 
g}$ to be the invertible composite
\begin{td}[row sep = 2.5em, column sep = 4em]
\csthree
	{g}
	{\pi_1, \,\dots\, , \pi_n}
	{\pairName(\p{1}{\ind{X}}, \,\dots\, ,  \p{n}{\ind{X}})} 
\arrow[swap]{d}
	{\assoc{g; \ind{\pi}; \pairName(\p{\bullet}{})}^{-1}} 
\arrow{r}{\epsilon_{A,g}} &
g \\

\cslr{g}{\cslr{\ind{\pi}}{\pairName(\p{1}{\ind{X}}, \p{n}{\ind{X}})}} 
\arrow[swap]{r}{\cslr{g}{\epsilonTimesInd{\bullet}{\ind{X}}}} &
\cslr{g}{\p{1}{\ind{X}},\dots, \p{n}{\ind{X}}} 
\arrow[swap]{u}{\subid{g}^{-1}}
\end{td}
For universality, let $f : \prodop_n(X_1, \,\dots\, , X_n) \to A$ by any 
multimap and $\delta : \cslr{f}{\pairName(\p{1}{}, \,\dots\, , \p{n}{})} \To g$ 
be any 2-cell. We 
define $\trans{\delta}$ as the following invertible composite, using the 2-cell 
$\gamma$ from the adjoint equivalence of 
Lemma~\ref{lem:biclone-with-products-adjoint-eq}:
\begin{equation*} 
f \XRA{\subid{}} 
\cslr{f}{\p{1}{(\prod \ind{X}})} 
\XRA{\cs{f}{\gamma^{-1}}}
\cslr{f}{\cslr{\pairName(\p{\bullet}{\ind{X}})}{\pi_1, \,\dots\, , \pi_n}} 
\XRA{\assoc{}^{-1}} 
\csthree{f}{\pairName(\p{\bullet}{\ind{X}})}{\ind{\pi}} 
\XRA{\cslr{\delta}{\ind{\pi}}} 
\cslr{g}{\ind{\pi}}
\end{equation*}
The rest of the proof is a diagram chase. To check the existence part of the 
universal property one uses law~(\ref{eq:triangle-law-one}) of an adjoint 
equivalence; for uniqueness one uses~(\ref{eq:triangle-law-two}). Since 
$\trans{\delta}$ is invertible whenever $\delta$ is, the unit is invertible and 
one obtains the required adjoint equivalence.
\end{proof}
\end{mylemma}

We collect these results together to obtain a bicategorical version of
Theorem~\ref{thm:representable-clone-cartesian-equivalence}. The final 
case is Lemma~\ref{lem:biuniversal-arrow-iff-family-of-equivalences}.

\begin{prooflessthm}  \label{thm:representable-iff-cartesian-biclone}
Let $(S, \biclone)$ be a biclone. Then the following are equivalent:
\begin{enumerate}
\item \label{c:biclone-representable} $(S, \biclone)$ admits a representable 
structure,
\item For every 
$X_1, \,\dots\, , X_n \in S \:\: (n \in \Nat)$ there exists a choice of object
$\prodop_n(X_1, \,\dots\, , X_n)$ and a
birepresentable multimap 
$\rho_{\ind{X}} : X_1, \,\dots\, , X_n \to \prodop_n(X_1, \,\dots\, , X_n)$,
\item \label{c:biclone-is-fp} $(S, \biclone)$ admits a cartesian structure, 
\item \label{c:pseudonatural-equivalences} For every 
$X_1, \,\dots\, , X_n \in S \:\: (n \in \Nat)$ 
there exists a choice of object $\prodop_n(X_1, \,\dots\, , X_n)$ together with 
a 
chosen family of
adjoint equivalences $(\mCatOf\biclone)\big(\Gamma; \prodop_n(X_1, \,\dots\, , 
X_n)\big) \simeq \prod_{i=1}^n (\mCatOf\biclone)(\Gamma; X_i)$, pseudonatural 
in the sense of 
Lemma~\ref{lem:biuniversal-arrow-iff-family-of-equivalences}(2).  \qedhere
\end{enumerate} 
\end{prooflessthm}

Restricting to unary hom-categories, case~(\ref{c:pseudonatural-equivalences}) 
of the theorem entails the following. 

\begin{prooflesscor} \label{cor:representable-biclone-to-fp-bicat}
For any representable biclone $(S, \biclone, \tensor_n)$, the nucleus 
$\overline\biclone$ is an fp-bicategory with product structure defined as in $\biclone$.	
\end{prooflesscor}

\subsection{Synthesising a type theory for fp-bicategories} 
\label{sec:synthesising-langcart}

\paragraph*{fp-Bicategories from cartesian biclones.} 
On page~\pageref{eq:cart-clone-cat-multigraph-diagram} we 
used diagram~(\ref{eq:cart-clone-cat-multigraph-diagram}) and the isomorphisms following 
to argue that, in order to 
construct a type theory describing cartesian categories, it is sufficient to 
construct a type theory for cartesian clones. Moreover, we showed how such a 
type theory could be synthesised from the construction of the free cartesian 
clone on a $\stlcTimes$-signature.

We repeat this process to synthesise the type theory $\langCart$. The starting 
point is an appropriate notion of signature. 
To extend from clones to biclones we extended from multigraphs to 
\mbox{2-multigraphs}; to extend from cartesian clones to cartesian biclones we extend 
$\stlcTimes$-signatures in the same way.

\begin{mydefn}
A \Def{$\langCart$-signature} $\sig = (\baseTypes, \graph)$ consists of 
\begin{enumerate}
\item A set of base types $\baseTypes$, 
\item  A 2-multigraph $\graph$ for which the set of nodes $\nodes{\graph}$ is 
generated by the grammar
\begin{equation} \label{eq:langcart-sig-grammar}
A_1, \,\dots\, , A_n ::= 
	B \st 
	\prodop_n(A_1, \,\dots\, , A_n) \qquad (B \in \baseTypes, n \in \Nat)
\end{equation}
\end{enumerate}
A \Def{homomorphism} $h : \sig \to \sig'$ of $\langCart$-signatures
is a 2-multigraph homomorphism $h : {\graph \to \graph'}$ that respects 
products, 
in the sense that 
$h_0{\left( \prodop_n(A_1, \,\dots\, , A_n) \right)} = 
	\prodop_n\left( h_0A_1, \,\dots\, , h_0A_n\right)$
for all ${A_1, \,\dots\, , A_n \in \nodes{\graph}} \:\: (n \in \Nat)$. 

We denote the category of $\langCart$-signatures by $\cartCatSigCat$ and 
the full sub-category of \Def{unary} $\langCart$-signatures---in which the 
2-multigraph $\graph$ is a 2-graph---by $\cartCatUnSigCat$.
\end{mydefn}

Every cartesian bi-multicategory (resp. cartesian biclone) determines an 
$\langCart$-signature, and every fp-bicategory 
determines a unary $\langCart$-signature. 

\begin{mynotation}[{\cf~Notation~\ref{not:stlc-times-alltypes-not}}] 
\label{not:langcart-times-alltypes-not}
For any $\langCart$-signature $\sig = (\baseTypes, \graph)$ we write 
$\allTypes\baseTypes$ for the set generated from $\baseTypes$ by the 
grammar~(\ref{eq:langcart-sig-grammar}).
In particular, when the signature is just a set (\ie~the graph $\graph$ has no edges) we denote the 
signature $\sig = (\baseTypes, \sig)$ simply by $\allTypes\baseTypes$. 
\end{mynotation}

The following result is proven in exactly the same way as 
Lemma~\ref{lem:stlcTimesSig-reflective-subcat}.

\begin{prooflesslemma} \label{lem:inc-adjoint-cartesian-structure-signatures}
The inclusion $\inc : \cartCatUnSigCat \hookrightarrow \cartCatSigCat$ has a 
right adjoint. 
\end{prooflesslemma}

The construction of the free cartesian clone on a cartesian category 
(Lemma~\ref{lem:free-cart-clone-on-cart-cat}) relies crucially on the identity 
$\seq{\pi_1, \,\dots\, , \pi_n} = \id_{(\prod_{i=1}^n X_i)}$ in a cartesian 
category so we cannot directly import this into the bicategorical setting. 
In place of diagram~(\ref{eq:cart-clone-cat-multigraph-diagram}), therefore, one obtains
a slightly restricted result. We will construct the following diagram of adjunctions, 
in which $\CartBicloneCat$ denotes the category 
of cartesian biclones and strict pseudofunctors strictly preserving the product 
structure, and $\CartBicatCat$ denotes the category of fp-bicategories and 
strict fp-pseudofunctors:
\begin{equation} 
\label{eq:cart-biclone-cat-multigraph-diagram}
\begin{tikzcd}[column sep = 3em]
\: &
\CartBicloneCat 
\arrow[bend right = 20]{dl} &
\: \\
\cartCatSigCat
\arrow[phantom]{ur}[description]{\adjDown}
\arrow[bend right = 20]{ur}
\arrow[bend right=20]{dr} &
\: &
\CartBicatCat
\arrow[bend left=20]{dl} \\
\: &
\cartCatUnSigCat 
\arrow[phantom]{ur}[description]{\adjUp{}}
\arrow[bend left=20]{ur}
\arrow[phantom]{ul}[description]{\adjUp{}}
\arrow[bend right = 20]{ul} &
\:
\end{tikzcd} 
\end{equation}
We shall then show that the free fp-bicategory on a unary $\langCart$-signature $\sig$ is obtained by
restricting the construction of the free cartesian biclone on $\sig$ to unary multimaps. 
Thus, the internal language of the free 
fp-bicategory on $\sig$ is the internal language of the free cartesian biclone on $\sig$, in which
every rule is restricted to unary multimaps. 
Here some care is required:
as we shall see, this is not the same as taking the nucleus of the free cartesian biclone.

Let us begin by making precise the notion of a (strict) morphism of cartesian 
biclones. The notion of biuniversal arrow for biclones is defined exactly as 
for bi-multicategories (Definition~\ref{def:universal-arrow-bi-multicat}); the 
corresponding notion of preservation extends that for bicategories
(Definition~\ref{def:strict-preservation-of-biuniversal-arrows}).

\begin{mydefn} 
Let $F : (S, \biclone) \to (T, \altBiclone)$ and $F' : (S', \biclone') \to (T', 
\altBiclone')$ be pseudofunctors of biclones and suppose $(R, u)$ and $(R', 
u')$ are biuniversal arrows from $F$ to $C \in T$ and from $F'$ to $C' \in T'$, 
respectively. A pair of pseudofunctors 
$(K : \altBiclone \to \altBiclone', L : \biclone \to \biclone')$ is a 
\Def{strict morphism of biuniversal arrows} from $(R,u)$ to $(R',u')$ if 
\begin{enumerate}
\item $K$ and $L$ are strict pseudofunctors satisfying $KF = F'L$, 
\item $LR =  R'$, $KC = C'$ and $Ku = u'$, 
\item The mappings $\psi_B : \altBiclone(FB, C) \to \biclone(B,R)$ and 
$\psi'_{B'} 
: \altBiclone'(F'B', C') \to \biclone'(B',R')$ are preserved, so that 
$L\psi_B(f) = 
\psi_{LB}'K(f)$ for every $f : FB \to C$, 
\item For every $B \in S$ and equivalence $\cs{u}{F(-)} : \baseCat(B, R) 
\leftrightarrows \altCat(FB, C) : \psi_B$ the universal arrow 
\mbox{$\epsilon_{B,h}  : \cs{u}{F\psi_B(h)} \To h$} is strictly preserved, in 
the sense that $K_{FB,C}(\epsilon_{B,h}) = \epsilon_{LB, Kh}$. \qedhere
\end{enumerate}
\end{mydefn}

We instantiate this in the case of cartesian biclones using the notation
of~(\ref{eq:cartesian-biclone-ump}) (page~\pageref{eq:cartesian-biclone-ump}).

\newpage
\begin{mydefn}
A \Def{cartesian pseudofunctor} 
$(F, \prodPres) : \cartClone{S}{\biclone} \to \cartClone{S'}{\biclone'}$ 
of cartesian biclones is a pseudofunctor $F : \biclone \to \biclone'$ equipped 
with a choice of equivalences 
$\pair{F\pi_1, \,\dots\, , F\pi_n} : 
	F{\left(\prodop_n(A_1, \dots A_n)\right)} \leftrightarrows
	\prodop_n \left(FA_1, \,\dots\, , FA_n \right) : \prodPres_{\ind{A}}$
for each ${A_1, \,\dots\, , A_n \in S} \:\: {(n \in \Nat)}$.

We call $(F, \prodPres)$ \Def{strict} if $F$ is a strict pseudofunctor and 
satisfies
\begin{align*}
F{\left(\prodop_n(A_1,\ldots,A_n)\right)} &= \prodop_n(FA_1,\ldots,FA_n) \\
F(\pi_i^{A_1,\ldots,A_n}) &= \pi_i^{FA_1,\ldots,FA_n} \\
F{\left(\pair{t_1,\ldots,t_n}\right)} &=\pair{Ft_1,\ldots,Ft_n} \\
F\varpi^{(i)}_{t_1,\ldots,t_n} &= \varpi^{(i)}_{Ft_1,\ldots,Ft_n} \\
\prodPres_{A_1,\ldots,A_n} &= \Id_{\Pi_n(FA_1,\ldots,FA_n)}
\end{align*}
and the equivalences are canonically induced by the 
$2$-cells 
$	\Id 	
	\XRA\iso 
	\pair{\cs{\pi_1}{\Id}, \,\dots\, , \cs{\pi_n}{\Id}} 
	\XRA\iso 
	\pair{\pi_1,\ldots,\pi_n}$. 
\end{mydefn}

If $(F, \prodPres) : \cartClone{S}{\biclone} \to \cartClone{S'}{\biclone'}$ is 
a cartesian pseudofunctor of biclones, one obtains an fp-pseudofunctor
between the associated fp-bicategories by restriction. 
To complete our diagram of adjunctions~(\ref{eq:cart-biclone-cat-multigraph-diagram}) 
it remains to construct 
free cartesian biclones and free fp-bicategories. We begin 
with the former. 

Theorem~\ref{thm:representable-clone-cartesian-equivalence} presents us 
with a choice. We can encode either representability (via the universal 
property~(\ref{eq:birepresentability-biclone-ump})) or cartesian structure (via 
the universal 
property~(\ref{eq:cartesian-biclone-ump})). In type-theoretic terms, this 
amounts to defining the universal property with respect to a pairing operation 
$x_1 : X_1, \,\dots\, , x_n : X_n \vdash 
\seq{x_1, \,\dots\, , x_n} : \prodop_n(X_1, \,\dots\, , X_n)$
or, alternatively, to defining the universal property with respect to 
projections 
$\left( p : \prodop_n(X_1, \,\dots\, , X_n) 
	\vdash \pi_i(p) : X_i \right)_{i=1, \,\dots\, , n}$.
We choose the latter because it more 
closely matches our definition of fp-bicategory.

\begin{myconstr} \label{constr:free-cart-biclone}
For any $\langCart$-signature $\sig$, define a cartesian biclone 
$\freeCartBiclone{\sig}$ 
with sorts 
\[
A_1, \,\dots\, , A_n ::= 
	B \st 
	\prodop_n(A_1, \,\dots\, , A_n) \qquad (B \in \baseTypes, n \in \Nat)
\] 
by extending 
the construction of the free biclone (Construction~\ref{constr:free-biclone}) 
with the following rules:
\begin{center}
\unaryRule
	{\faketext}
	{\pi_i^{\ind{A}} \in \freeCartBiclone{\sig}\left(\prodop_n(A_1, \,\dots\, , 
	A_n); 
	A_i\right)}
	{$(1 \leq i \leq n)$}
	
\unaryRule
	{\left(  t_i \in \freeCartBiclone{\sig}(\Gamma; A_i) \right)_{i=1, \,\dots\, , 
	n}}
	{\pair{t_1, \,\dots\, , t_n} \in 
		\freeCartBiclone{\sig}\left(\Gamma; \prodop_n(A_1, \,\dots\, , A_n)\right)}
	{}

\unaryRule
	{\left(  t_i \in \freeCartBiclone{\sig}(\Gamma; A_i) \right)_{i=1, \,\dots\, , 
	n}}
	{\epsilonTimesInd{i}{\ind{t}} \in 
		\freeCartBiclone{\sig}\left(\Gamma; A_i\right)
			\left(\pair{t_1, \,\dots\, , t_n}, t_i \right)}
	{$(1 \leq i \leq n)$}

\unaryRule
	{\left( 
		\alpha_i \in 
			\freeCartBiclone{\sig}(\Gamma; A_i)(\cs{\pi_i^{\ind{A}}}{u}, t_i) 
	\right)_{i=1, \,\dots\, , n}}
	{\transTimes{\alpha_1, \,\dots\, , \alpha_n} \in 
		\freeCartBiclone{\sig}\left(\Gamma; \prodop_n(A_1, \,\dots\, , A_n) \right)
			\left(u, \pair{t_1, \,\dots\, , t_n} \right)}
	{}
\end{center}
Moreover, extend the equational theory $\equiv$ of 
Construction~\ref{constr:free-biclone} with the following rules 
encoding the universal property~(\ref{eq:cartesian-biclone-ump}):
\begin{itemize}
\item If $\alpha_i : u \To t_i : \Gamma \to A_i$ for $i=1, \,\dots\, , n$, then 
$\alpha_i \equiv \epsilonTimesInd{i}{\ind{t}} \vert \transTimes{\alpha_1, 
\dots, \alpha_n}$ 
for $i=1, \,\dots\, , n$, 

\item If $\gamma : u \To \pair{t_1, \,\dots\, , t_n} :
			\Gamma \to \prodop_n(A_1, \,\dots\, , A_n)$, then
		$\gamma \equiv 
			\transTimes{\epsilonTimesInd{1}{\ind{t}} \vert 
					\cs{\Id_{\pi_1}}{\gamma}, 
				\dots, \epsilonTimesInd{n}{\ind{t}} \vert 
					\cs{\Id_{\pi_n}}{\gamma}}$,

\item If $\alpha_i \equiv \alpha_i'$ for $\alpha_i, \alpha_i'$ 2-cells
of type $\cs{\pi_i^{\ind{A}}}{u} \To t_i$ for $i=1, \,\dots\, , n$, then
$\transTimes{\alpha_1, \,\dots\, , \alpha_n} \equiv 
	\transTimes{\alpha_1', \,\dots\, , \alpha_n'}$. 
\end{itemize}
Finally, we require that every $\epsilonTimesInd{i}{\ind{t}}$ and
$\etaTimes{t} 
	:= \transTimes{\Id_{\cs{\pi_1}{t}}, \,\dots\, , \Id_{\cs{\pi_n}{t}}}$ 
is invertible.
\end{myconstr}

\begin{mylemma}  \label{lem:free-cart-biclone-cartesians-structure}
For any $\langCart$-signature $\sig$ and any finite family of 2-cells  
$(\alpha_i : \rewrite{\hcomp{\pi_i}{u}}{t_i} : 
	\Gamma \to A_i)_{i=1, \dots, n}$ 
in $\freeCartBiclone{\sig}$, then 
$\transTimes{\alpha_1, \,\dots\, , \alpha_n}$ is the unique 
2-cell $\gamma$ (modulo ${\equiv}$) such that 
$\alpha_i \equiv \epsilonTimesInd{i}{\ind{t}} \vert \gamma$ 
for $i=1, \,\dots\, , n$.
\begin{proof}
The existence part of the claim is immediate. For uniqueness, if
$\gamma$ satisfies the given 
equation then 
$\gamma \equiv 
\transTimes{\epsilonTimesInd{1}{\ind{t}} \vert 
					\cs{\Id_{\pi_1}}{\gamma}, 
				\dots, \epsilonTimesInd{n}{\ind{t}} \vert 
					\cs{\Id_{\pi_n}}{\gamma}} 
\equiv \transTimes{\alpha_1, \,\dots\, , \alpha_n}$, as 
claimed. 
\end{proof}
\end{mylemma}

It follows that $\freeCartBiclone{\sig}$ is cartesian. The associated free property 
is 
then straightforward. 

\begin{mylemma} \label{lem:free-cart-biclone-proved}
For any $\langCart$-signature $\sig$, cartesian biclone 
$\cartClone{T}{\altBiclone}$ 
and 
$\langCart$-signature homomorphism $h : \sig \to \altBiclone$  from $\sig$ to 
the 
$\langCart$-signature underlying $\cartClone{T}{\altBiclone}$ there exists a 
strict cartesian pseudofunctor 
$\ext{h} : \freeCartBiclone{\sig} \to \altBiclone$, 
unique such that $\ext{h} \circ \inc = h$, for 
$\inc : \sig \hookrightarrow \freeCartBiclone{\sig}$
the inclusion. 
\begin{proof}
We extend the pseudofunctor $\ext{h}$ defined in Lemma~\ref{lem:free-biclone} 
by setting
\begin{align*}
\ext{h}{\left( \prodop_n(A_1, \,\dots\, , A_n) \right)} &:=
	\prodop_n {\left( \ext{h}(A_1), \,\dots\, , \ext{h}(A_n) \right)} \\[5pt]
\ext{h}(\pi_i^{\ind{A}}) &:= \pi_i^{\ext{h}(\ind{A})} \\
\ext{h}(\pair{t_1, \,\dots\, , t_n}) &:= 
	\pair{\ext{h}(t_1), \,\dots\, , \ext{h}(t_n)} \\[5pt]
\ext{h}(\epsilonTimesInd{i}{\ind{t}}) &:=
	\epsilonTimesInd{i}{\ext{h}(\ind{t})} \\
\ext{h}{\left( \transTimes{\alpha_1, \,\dots\, , \alpha_n} \right)} &:=
	\transTimes{\ext{h}(\alpha_1), \,\dots\, , \ext{h}(\alpha_n)}
\end{align*}
It is clear this defines a strict cartesian pseudofunctor. For uniqueness, all 
the cases apart from $\transTimes{\alpha_1, \,\dots\, , \alpha_n}$ are 
determined by 
the definition of strict cartesian pseudofunctor. To complete the proof, we 
adapt the argument of Lemma~\ref{lem:strict-preservation-strict-pres-UMP}. 
For any strict cartesian pseudofunctor 
$F : \freeCartBiclone{\sig} \to \altBiclone$ and 2-cells
$(\alpha_i : \cs{\pi_i^{\ind{A}}}{u} \To t_i : \Gamma \to A_i)_{i=1, \,\dots\, 
, n}$,
\enlargethispage*{3\baselineskip}
\begin{align*}
\epsilonTimesInd{i}{F\ind{t}} \vert 
	F{\left(\transTimes{\alpha_1, \,\dots\, , \alpha_n}\right)} 
	&= F(\epsilonTimesInd{i}{\ind{t}}) \vert 
		F{\left(\transTimes{\alpha_1, \,\dots\, , \alpha_n}\right)} \\
	&= F{\left(\epsilonTimesInd{i}{F\ind{t}} \vert 
		\transTimes{\alpha_1, \,\dots\, , \alpha_n}\right)} \\
	&= F\alpha_i 
\end{align*}
for $i=1, \,\dots\, , n$. Hence, by the universal 
property~(\ref{eq:cartesian-biclone-ump}) of a cartesian biclone, 
$F{\left(\transTimes{\alpha_1, \,\dots\, , \alpha_n}\right)} = 
	\transTimes{F\alpha_1, \,\dots\, , F\alpha_n}$.
\end{proof}
\end{mylemma}

\begin{myremark}
The preceding proof should be compared to that for the free cartesian 
clone on a $\stlcTimes$-signature 
(Lemma~\ref{lem:free-cart-clone-on-cart-cat}). The argument for uniqueness lifts to 2-cells by virtue of the fact that 
pseudofunctors strictly preserve vertical composition.
\end{myremark}

It remains to 
construct the free fp-bicategory on a unary $\stlcTimes$-signature and relate it to the free 
cartesian biclone over the same signature. The proof is straightforward: one 
restricts Lemma~\ref{lem:free-cart-biclone-proved} to unary multimaps 
and observes the same 
universal property holds. Example~\ref{exmp:nucleus-of-fp-not-free} shows that 
it is important to restrict every rule 
to unary multimaps---\ie~require that $\len{\Gamma}=1$
for every rule in Construction~\ref{constr:free-cart-biclone}---rather than simply taking the nucleus
of $\freeCartBiclone\sig$.

\begin{prooflesslemma} \label{lem:free-fp-bicat-from-free-biclone} 
For any unary $\langCart$-signature $\sig$, let $\freeCartBicat{\sig}$ denote the 
fp-bicategory obtained by restricting every rule of 
Construction~\ref{constr:free-cart-biclone} to unary multimaps and 2-cells 
between them, and let $h : \sig \to \altCat$ be a 
$\langCart$-signature homomorphism from $\sig$ to 
the 
$\langCart$-signature underlying an fp-bicategory $\fpBicat{\altCat}$. Then 
there exists a 
strict fp-pseudofunctor 
$\ext{h} : \freeCartBicat{\sig} \to \altCat$, 
unique such that $\ext{h} \circ \inc = h$, for 
$\inc : \sig \hookrightarrow \freeCartBicat{\sig}$
the inclusion. 
\end{prooflesslemma}

\begin{myexmp} \label{exmp:nucleus-of-fp-not-free}
Fix a $\langCart$-signature $\sig = (\baseTypes, \graph)$. Then the nucleus 
$\nucleus{\freeCartBiclone{\sig}}$
of 
$\freeCartBiclone{\sig}$ is not isomorphic to 
$\freeCartBicat\sig$. Roughly speaking, the composite
$\cs{\p{1}{A, B}}{\pi_1, \pi_2} : A \times B \to A$ 
exists in the free cartesian biclone on a signature $\sig$, but not in the free 
fp-bicategory on $\sig$. Let us make this precise. 

Since the freeness universal property of $\freeCartBicat{\sig}$ is strict we may exploit the following principle, which 
restates the fact that free objects are unique up to \emph{canonical} 
isomorphism: if 
$\baseCat$ and $\baseCat'$ are both the free fp-bicategory on $\sig$, then the 
canonical map $\baseCat \to \baseCat'$ extending the unit is an isomorphism. We 
claim that the canonical map 
$\ext\inc : \freeCartBicat{\sig} \to \nucleus{\freeCartBiclone\sig}$ 
extending the inclusion $\inc : \sig \hookrightarrow 
\nucleus{\freeCartBiclone\sig}$
is not an 
isomorphism. Since an isomorphism is necessarily a bijection on hom-sets, it suffices to find 
a morphism in $\nucleus{\freeCartBiclone{\sig}}$ that is not in the image of 
$\ext\inc$. We claim that, where $X, Y \in \allTypes\baseTypes$, then
$\cs{\p{1}{X,Y}}{\pi_1, \pi_2} : X \times Y \to X$ is not in the image of
$\ext\inc$. To see this is the case, observe that a morphism $h$ is in the 
image of $\ext\inc$ if and only if it falls into one of the following 
(disjoint) sets:
\begin{enumerate}
\item The \Def{basic maps} $\pi_i, \eval$ and $\Id$, 
\item Maps in the image of an \Def{operator}: $\lambda f$ or $\seq{f_1, \dots, f_n}$
for $f, f_1, \dots, f_n$ in the image of $\ext\inc$, 
\item The \Def{composites} $f \circ g$ where $f$ and $g$ are both in the image 
of 
$\ext\inc$. 
\end{enumerate}
It is clear that $\cs{\p{1}{X,Y}}{\pi_1, \pi_2}$ is not of any of these types, 
and so is not in the image of $\ext\inc$. It follows that 
$\ext\inc$ is not an isomorphism, and hence that
$\nucleus{\freeCartBiclone{\sig}}$ is not the free fp-bicategory on $\sig$. 
\end{myexmp}

Lemma~\ref{lem:free-fp-bicat-from-free-biclone} guarantees that the free fp-bicategory 
on a $\langCart$-signature $\sig$ arises by restricting every rule of the type theory 
for cartesian biclones to unary contexts 
and constructing the syntactic model. Hence, it suffices to construct 
a type theory for cartesian biclones. We do this by extending the type 
theory $\langBiclone$ for biclones with rules corresponding to those of  
Construction~\ref{constr:free-cart-biclone}.

\section{The type theory \texorpdfstring{$\langCart$}{for bicategories with 
finite products}}

For a $\langCart$-signature $\sig = (\baseTypes, \graph)$ we denote the 
associated type theory by $\langCart(\sig)$. The types of $\langCart(\sig)$
are the nodes of $\graph$. The rules are all those of $\langBiclone$ together 
with those of Figures~\ref{r:fp:products-terms}--\ref{r:fp:inverses}. Note that 
we specify the invertibility of the unit and counit by introducing explicit 
inverses for these rewrites (Figure~\ref{r:fp:inverses}). 

The tupling operation is functorial with respect to vertical composition and 
the unit of the 
adjunction is 
obtained by applying the universal property to the identity (see also 
Lemma~\ref{lem:langcart-admissible-rules}).

\begin{mydefn}  \label{def:fp:unit} \quad
\begin{enumerate}
\item For any family of derivable rewrites \mbox{$(\Gamma \vdash \tau_i : 
\rewrite{t_i}{t_i'} : A_i)_{i=1,\dots, n}$} we define $\pair{\tau_1, \,\dots\, 
, 
\tau_n} : \pair{t_1, \,\dots\, , t_n} \To \pair{t_1', \,\dots\, , t_n'}$ to be 
the 
rewrite $\transTimes{\tau_1 \vert \epsilonTimesInd{1}{t_1, \,\dots\, , t_n}, 
\,\dots\, , 
\tau_n \vert \epsilonTimesInd{n}{t_1, \,\dots\, , t_n}}$ in context $\Gamma$. 
\item For any derivable term $\Gamma \vdash t : \prodop_n(A_1, \,\dots\, , 
A_n)$ we 
define the unit $\etaTimes{t} : t \To \pair{\hcomp{\pi_1}{t}, \,\dots\, , 
\hcomp{\pi_n}{t}}$ to be the rewrite $\transTimes{\id_{\hcomp{\pi_1}{t}}, 
\dots, \id_{\hcomp{\pi_n}{t}}}$ in context $\Gamma$.   \qedhere
\end{enumerate}
\end{mydefn}

The rules of $\langCart$ provide a relatively compact way to construct the 
structure required for cartesian clones. In particular, the focus on (global) 
biuniversal arrows and (local) universal arrows---and the corresponding fact 
that one does not need to specify a triangle law relating the unit and 
counit---contrasts 
with all previous work on type theories for cartesian closed 2-categories~\cite{
Seely1987, 
Hilken1996, 
Tabareau2011, 
Hirschowitz2013}, 
which encode the pairing and projection 
operations on rewrites directly. Reproducing 
the triangle-law approach in the context of fp-bicategories would require:

\begin{enumerate}
\item For every sequence of types $A_1, \,\dots\, , A_n$ a product type 
$\prodop_n(A_1, \,\dots\, , A_n)$, 
\item Projection and tupling operations on terms as in the usual simply-typed 
lambda calculus, 
\item Tupling and projection operations on rewrites,
\item An invertible unit $\etaTimes{u} : u \To \seqlr{\pi_1(u), \,\dots\, , 
\pi_n(u)}$ 
in context $\Gamma$ for every $\Gamma \vdash u : \prodop_n(A_1, \,\dots\, , 
A_n)$ 
and an invertible counit $\epsilonTimesInd{i}{\ind{t}} : \hpi{i}{\seqlr{t_1, 
\dots, t_n}} \To t_i \:\: (i=1, \,\dots\, ,n)$ in context $\Gamma$ for every 
$(\Gamma \vdash t_i : A_i)_{i=1, \,\dots\, , n}$.
\end{enumerate} 
This data must be subject to an equational theory requiring naturality of each 
$\etaTimes{u}$ and $\epsilonTimesInd{i}{\ind{t}}$, the two triangle laws, 
functorality of the tupling and projection operations on rewrites, and that the 
equational theory is a congruence with respect to these operations. Such an 
approach, therefore, requires many more rules. Moreover, the calculus of 
(bi)universal arrows provided by $\langCart$ captures a categorical style of 
reasoning, because the syntax allows one to manipulate the 
universal property through primitives in the type theory.

\setlength{\floatsep}{5pt plus 1.0pt minus 2.0pt} 

\begin{figure*}[!h]
{\small
\begin{minipage}{\textwidth}
\begin{mdframed}
\centering
\input{rules/cart/products-terms}
\vspace{-\treeskip}
\caption{\label{r:fp:products-terms} Terms for product structure}
\end{mdframed}
\end{minipage}

\begin{minipage}{\textwidth}
\begin{mdframed}
\centering
\input{rules/cart/products-rewrites}
\caption{Rewrites for product structure \label{r:fp:ProductsRewrites}}
\end{mdframed}
\end{minipage}

\begin{minipage}{\textwidth}
\begin{mdframed}
\centering
\input{rules/cart/trans-ump}
\input{rules/cart/extra-congruence}
\caption{Universal property and congruence laws for  $\transTimes{\alpha_1, \dots,\alpha_n}$\label{r:fp:umpTrans}}
\end{mdframed}
\end{minipage}

\begin{minipage}{\textwidth}
\begin{mdframed}
\centering
\input{rules/cart/invertibility-intro}
\input{rules/cart/invertibility-congruences}
\caption{Inverses for the unit and counit\label{r:fp:inverses}}
\end{mdframed}
\end{minipage}
}
\begin{manyfigcap}
Rules for $\langCart(\graph)$.
\end{manyfigcap}
\end{figure*}

\paragraph{$\alpha$-equivalence and free variables.}

The well-formedness properties of $\langBiclone$ extend to $\langCart$; we 
briefly note them here. As we have not introduced any binding constructs, the 
definition of $\alpha$-equivalence extends straightforwardly from that for 
$\langBiclone$. 

\begin{mydefn} \label{def:fp:alpha-equivalence}
For any $\langCart$-signature $\sig$ we extend 
Definition~\ref{def:biclone-alpha-eq} 
to define the \Def{$\alpha$-equivalence relation} 
$\aeq$ for $\langCart(\sig)$. For terms we take the same 
set of 
rules; the substitution operation $t[u_i / x_i]$ is extended by the rules
\begin{center}
$\pi_k(p)[u/p] := \hpi{k}{u}$ \:\: and \:\: $\pair{t_1, \,\dots\, , t_n}[u_i / 
x_i] 
:= \pair{t_1[u_i / x_i], \,\dots\, , t_n[u_i / x_i]}$
\end{center}
For rewrites, we add the rules
\begin{center}
\unaryRule{(t_i \aeq t_i')_{i= 1, \,\dots\, , n}}{\epsilonTimesInd{k}{t_1, 
\,\dots\, , 
t_n} \aeq \epsilonTimesInd{k}{t_1', \,\dots\, , t_n'}}{$(1 \leq k \leq n)$}
\qquad
\unaryRule{\sigma_1 \aeq \sigma_1' \quad \dots \quad \sigma_n \aeq 
\sigma_n'}{\transTimes{\sigma_1, \,\dots\, , \sigma_n} \aeq 
\transTimes{\sigma_1', 
\dots, \sigma_n'}}{} \vspace{-\treeskip}
\end{center}
where the meta-operation of capture-avoiding substitution is extended by the 
rules
\begin{center}
$\epsilonTimesInd{k}{t_1, \,\dots\, , t_n}[u_i/x_i] := 
\epsilonTimesInd{k}{t_1[u_i/x_i], \,\dots\, , t_n[u_i/x_i]}$ \:\: and \:\: 
$\transTimes{\ind{\alpha}}[u_i/x_i] := \transTimes{\ind{\alpha}[u_i/x_i]}$
\end{center}
Finally, we define $\fv(\sigma^{-1}) := \fv(\sigma)$. 
\end{mydefn}

As for $\langBicat$, we work up to $\alpha$-equivalence of terms and rewrites, 
silently identifying terms and rewrites with their $\alpha$-equivalence 
classes. 

Extending the definition of free variables is similarly straightforward.

\begin{mydefn} \label{def:fp:free-vars}
Fix a $\langCart$-signature $\sig$. We define the \Def{free variables in a 
term} $t$ 
in $\langCart(\sig)$ by extending Definition~\ref{def:biclone-free-vars} as 
follows: 
\begin{center}
$\fv\big(\pair{t_1, \,\dots\, , t_n}\big) := \bigcup_{i=1}^n \fv(t_i)$ \:\: and 
\: 
\: $\fv\big(\pi_k(p)\big) := \{ p \}$
\end{center}
Define the \Def{free variables in a rewrite $\tau$} in 
$\langCart(\sig)$ by extending Definition~\ref{def:biclone-free-vars} as 
follows:
\begin{center}$\fv(\epsilonTimesInd{k}{t_1, \,\dots\, , t_n}) := \fv(t_k)$ 
\quad
and
\quad
$\fv\big(\transTimes{\alpha_1, \,\dots\, , \alpha_n}\big) := \bigcup_{i=1}^n 
\fv(\alpha_i)$
\end{center}
We define the free variables of a specified inverse $\sigma^{-1}$ to be exactly 
the free variables of $\sigma$. An occurrence of a variable in a term (resp. 
rewrite) is \Def{bound} if it is not free. 
\end{mydefn} 

The next two lemmas---both of which are proven by structural induction---show 
that the preceding definitions behave in the way one would expect.

\begin{prooflesslemma} 
Let $\sig$ be a $\langCart$-signature. Then in $\langCart(\sig)$:
\begin{enumerate} 
\item If $\Gamma \vdash t : B$ and $t =_\alpha t'$ then $\Gamma \vdash t' :B$,
\item If $\Gamma \vdash \tau : \rewrite{t}{t'} : B$ and $\tau =_\alpha \tau'$ 
then $\Gamma \vdash \tau' : \rewrite{t}{t'} : B$,
\item If $\tau_i \aeq \tau_i'$ for $i=1, \,\dots\, , n$, then 
$\pair{\tau_1, \,\dots\, , \tau_n} \aeq \pair{\tau_1', \,\dots\, , \tau_n'}$, 
\item If $u \aeq u'$ then $\etaTimes{u} \aeq \etaTimes{u'}$. \qedhere
\end{enumerate} 
\end{prooflesslemma} 

\begin{prooflesslemma} 
Let $\sig$ be a $\langCart$-signature. For any derivable judgements $\Gamma 
\vdash u : 
B$ and \mbox{$\Gamma \vdash \tau : \rewrite{t}{t'} : B$} in 
$\langCart(\sig)$, 
\begin{enumerate} 
\item $\fv(u) \subseteq \dom(\Gamma)$, 
\item $\fv(\tau) \subseteq \dom(\Gamma)$, 
\item The judgements $\Gamma \vdash t : B$ and $\Gamma \vdash t' : B$ are both 
derivable. 
\end{enumerate} 
Moreover, whenever $(\Delta \vdash u_i : A_i)_{i= 1, \,\dots\, , n}$ and 
$\Gamma := 
(x_i  : A_i)_{i = 1, \,\dots\, , n}$, then
\begin{enumerate}
\item If $\Gamma \vdash t : B$, then $\Delta \vdash t[u_i / x_i] : B$, 
\item If $\Gamma \vdash \tau : \rewrite{t}{t'} : B$, then $\Delta \vdash 
\tau[u_i / x_i] : \rewrite{t[u_i / x_i]}{t'[u_i / x_i]} : B$. \qedhere
\end{enumerate}
\end{prooflesslemma}

\subsection{\texorpdfstring{The syntactic model for $\langCart$}{The syntactic model}}  
\label{sec:fp:syntactic-model}

Lemma~\ref{lem:free-fp-bicat-from-free-biclone}
guarantees that, in order to construct a type theory for fp-bicategories, it suffices to construct a type theory for cartesian biclones. To verify that 
$\langCart$ is such a type theory, furthermore, it suffices to show that its 
syntactic model is canonically isomorphic to the free cartesian biclone 
$\freeCartBiclone{\sig}$ over the same signature in the category 
$\CartBicloneCat$.

The syntactic model is constructed by 
extending~Construction~\ref{constr:langbiclone-syntactic-model}.

\begin{myconstr} \label{constr:langcart-syntactic-model}
For any $\langCart$-signature $\sig$ define the \Def{syntactic model} 
$\syncloneCart{\sig}$ of $\langCart(\sig)$ as follows. The sorts are 
nodes $A, B, \dots$ of $\graph$. For 
$A_1, \,\dots\, , A_n, B \in \baseTypes \:\: (n \in \Nat)$ the hom-category 
$\syncloneCart{\sig}(A_1, \,\dots\, , A_n; B)$ 
has objects $\alpha$-equivalence classes of terms 
\mbox{$(x_1 : A_1, \,\dots\, , x_n : A_n \vdash t : B)$} derivable in 
$\langCart(\sig)$. We assume a fixed 
enumeration $x_1, x_2, \dots$ of variables, and that the variable name in the 
$i$th position is determined by this enumeration. Morphisms in 
$\syncloneCart{\sig}(A_1, \,\dots\, , A_n; B)$ are $\alpha{\equiv}$-equivalence 
classes of rewrites \mbox{$(x_1 : A_1, \,\dots\, , x_n : A_n \vdash \tau : 
\rewrite{t}{t'} : B)$}. Composition is vertical composition with identity  
$\id_t$; the substitution operation is explicit substitution and the structural 
rewrites are $\assoc{}, \subid{}$ and $\indproj{i}{}$. 
\end{myconstr}

Inspecting each rule in turn, one sees that 
$\syncloneCart{\sig}$ is merely $\freeCartBiclone{\sig}$, 
presented with the notation $x_1 : X_1, \,\dots\, , x_n : X_n \vdash t : B$ 
instead 
of $t : X_1, \,\dots\, , X_n \to B$. We make this statement precise by 
establishing 
it satisfies the same universal property. 

Lemma~\ref{lem:free-cart-biclone-cartesians-structure}, restated in 
type-theoretic notation, becomes the following.

\begin{mylemma} \label{lem:fp:trans-ump}
For any $\langCart$-signature $\sig$, if the judgements $(\Gamma \vdash 
\alpha_i : 
\rewrite{\hcomp{\pi_i}{u}}{t_i} : A_i)_{i = 1, \,\dots\, , n}$ are derivable in 
$\langCart(\sig)$ then $\transTimes{\alpha_1, \,\dots\, , \alpha_n}$ is the 
unique 
rewrite $\gamma$ (modulo $\alpha{\equiv}$) such that the equality
\begin{equation} \label{eq:fp:counit-ump}
\Gamma \vdash \epsilonTimesInd{k}{t_1, \,\dots\, , t_n} \vert 
\hcomp{\pi_k}{\gamma} 
\equiv \alpha_k : \rewrite{\hcomp{\pi_i}{u}}{t_k} : A_k
\end{equation}
is derivable for $k = 1, \,\dots\, , n$.
\begin{proof}
By \rulename{U1} (Figure~\ref{r:fp:umpTrans}) the rewrite 
$\transTimes{\alpha_1, \,\dots\, , \alpha_n}$ certainly 
satisfies~(\ref{eq:fp:counit-ump}). For any other $\gamma$ satisfying the 
equation, $\gamma \overset{\rulename{U2}}{\equiv} 
\transTimes{\epsilonTimesInd{1}{\ind{t}} \vert \hcomp{\pi_1}{\gamma}, \,\dots\, 
, 
\epsilonTimesInd{n}{\ind{t}} \vert \hcomp{\pi_n}{\gamma}} 
\overset{\text{cong}}{\equiv} \transTimes{\alpha_1, \,\dots\, , \alpha_n}$, as 
claimed. 
\end{proof}
\end{mylemma}

\begin{myremark} \label{rem:fp:biuniversal-universal-interleaving}
In the light of the preceding lemma, for any $\langCart$-signature $\sig$ the 
mappings
\begin{align*}
(\alpha_1, \,\dots\, , \alpha_n) &\mapsto \transTimes{\alpha_1, \,\dots\, , 
\alpha_n} \\
(\epsilonTimesInd{1}{\ind{t}} \vert \hcomp{\pi_1}{\tau}, \,\dots\, , 
\epsilonTimesInd{n}{\ind{t}} \vert \hcomp{\pi_n}{\tau}) &\mapsfrom \tau
\end{align*} 
define the following bijective correspondence of rewrites, derivable in 
$\langCart(\sig)$:
\begin{center}
\begin{bprooftree} 
\AxiomC{$\hcomp{\pi_k}{u} \To t_k \quad (k = 1, \,\dots\, , n)$} 
\doubleLine 
\UnaryInfC{$u \To \pair{t_1, \,\dots\, , t_n}$} 
\end{bprooftree} 
\end{center}
It is natural to conjecture that a calculus for 
fp-\emph{tri}categories (resp. fp-$\infty$-categories) would have three (resp. a countably 
infinite tower of) such correspondences. Similar considerations will apply to 
exponentials.
\end{myremark}

It also follows from the preceding lemma that $\syncloneCart{\sig}$ is 
cartesian: the adjoint equivalence is exactly
\begin{align*}
\syncloneCart{\sig}\big(\Gamma, \prodop_n(A_1, \,\dots\, , A_n)\big) 
&\xra{\simeq} 
\prodop_{i=1}^n \syncloneCart{\sig}(\Gamma; A_i) \\
\big(\Gamma \vdash u : \prodop_n(A_1, \,\dots\, , A_n)\big) 
&\mapsto 
(\Gamma \vdash \hpi{i}{u} : A_i)_{i=1,\dots,n}
\end{align*} 
where the pseudoinverse 
$
\prod_{i=1}^n \syncloneCart{\sig}(\Gamma; A_i) \to 
	\syncloneCart{\sig}\big(\Gamma, \prodop_n(A_1, \,\dots\, , A_n)\big)
$
is the $\pairName$ operation. The universal property of 
$\syncloneCart{\sig}$ interprets each term as its corresponding construct.


\begin{mypropn} \label{propn:synclonecart-free-property}
For any $\langCart$-signature $\sig = (\baseTypes, \graph)$, cartesian biclone 
$\cartClone{T}{\altClone}$ and $\langCart$-signature homomorphism 
$h : \sig \to \biclone$, 
there exists a unique strict cartesian pseudofunctor 
$h\sem{-} : \syncloneCart{\sig} \to \biclone$ 
such that 
$h\sem{-} \circ \iota = h$, 
for 
$\iota : \sig \hookrightarrow \syncloneCart{\sig}$ the inclusion.
\begin{proof}
The pseudofunctor is constructed by induction on the syntax of 
$\langCart(\sig)$ as follows:
\begin{align*}
h\sem{B} &:= h(B) \qquad \text{\normalsize on base types} \\
h\sem{\prodop_m(B_1, \,\dots\, , B_m)} &:= 
	\prodop_m \left(h\sem{B_1}, \,\dots\, , h\sem{B_m}\right)
\\[6pt]
h\sem{\Gamma \vdash x_k :A_i} &:= \p{k}{h\sem{A_1}, \,\dots\, , h\sem{A_n}} \\ 
h\sem{\Gamma \vdash c(x_1, \,\dots\, , x_n) : B} &:= h(c) \qquad 
\text{\normalsize 
for } c \in 
\graph(\ind{A};B) \\ 
h\sem{\Delta \vdash \hcomp{t}{x_i\mapsto u_i} : B} &:= 
	\cslr
		{\big(h\sem{\Gamma \vdash t : B}\big)}
		{h\sem{\Delta \vdash \ind{u}: \ind{A}}} \\ 
h\sem{\Gamma \vdash \pair{t_1, \,\dots\, , t_m} : \prodop_m(B_1, \,\dots\, , 
B_m)} &:= 
	\pair{h\sem{\Gamma \vdash t_1 : B_1}, \,\dots\, , 
			h\sem{\Gamma \vdash t_m : B_m}} \\
h\sem{p : \prodop_m(B_1, \,\dots\, , B_m) \vdash \pi_k(p) : B_k} &:= 
	\pi_k^{h\sem{B_1}, \,\dots\, , h\sem{B_m}} \\[6pt]
h\sem{\Gamma \vdash \id_t : \rewrite{t}{t} : B} &:= 
	\id_{h\sem{\Gamma \vdash t : B}} \\
h\sem{\Gamma \vdash \constrewr(\ind{x}) : \rewrite{c(\ind{x})}{c'(\ind{x})} : 
B} 
	&:= h(\constrewr) \qquad \text{\normalsize for } \constrewr \in 
	\graph(\ind{A}, B)(c,c')  \\
h\sem{\Gamma \vdash 
	\epsilonTimesInd{k}{t_1, \,\dots\, , t_m} : 
	\rewrite{\hcomp{\pi_k}{\pair{t_1, \,\dots\, , t_m}}}{t_k} : B_k} &:= 
	\epsilonTimesInd{k}{h\sem{t_1}, \,\dots\, , h\sem{t_m}} \\
h\sem{\Gamma \vdash 
		\transTimes{\alpha_1, \,\dots\, , \alpha_m} : 
			\rewrite{u}{\pair{\ind{t}}} : \prodop_m \ind{B}} &:= 
		\transTimes{h\sem{\Gamma \vdash 
			\ind{\alpha} : \rewrite{\hcomp{\ind{\pi}}{u}}{\ind{t}} :\ind{B}}}
			\\
h\sem{\Gamma \vdash \tau' \vert \tau : \rewrite{t}{t''} : B} &:= 
	h\sem{\Gamma \vdash \tau' : \rewrite{t'}{t''} : B} \vert 
	h\sem{\Gamma \vdash \tau : \rewrite{t}{t'} : B} \\
h\sem{\Delta \vdash \hcomp{\tau}{\sigma_i} 
	: \rewrite{\hcomp{t}{u_i}}{\hcompsmall{t'}{u'_i}} : B} 
	&:= 
		\cslr
			{\big(h\sem{\Gamma \vdash \tau : 
				\rewrite{t}{t'} : B}\big)}
			{h\sem{\sigma_1}, \,\dots\, , h\sem{\sigma_n}}
\end{align*} 
where $\Gamma := (x_i : A_i)_{i=1,\dots,n}$ and we
abbreviate 
$h\sem{\Delta \vdash \sigma_i : 
				\rewrite{u_i}{u_i'} : A_i}$
by $h\sem{\sigma_i}$ in the final rule. It is clear that this defines a 
strict pseudofunctor; 
the $\transTimes{\alpha_1, \,\dots\, , \alpha_m}$ case is required by the 
strict preservation of universal and biuniversal arrows 
(\cf~Lemma~\ref{lem:free-cart-biclone-proved}). 
\end{proof}
\end{mypropn}

%

Lemma~\ref{lem:free-fp-bicat-from-free-biclone}, together with the preceding 
proposition, entail that the free fp-bicategory on a unary 
$\langCart$-signature is obtained as follows. First, one restricts
$\langCart$ to unary contexts. Then one constructs the syntactic model 
in the same manner as Construction~\ref{constr:langcart-syntactic-model}, except 
morphisms and 2-cells are equivalence classes of terms and rewrites in  
this restricted type theory.
Thus, define $\urestrict{\langCart}$
to be the type theory obtained by restricting  $\langCart$ to contexts of 
the form 
$x : A$ (defined by Figure~\ref{r:biclone:unary-contexts} on page~\pageref{r:biclone:unary-contexts}. 
The resulting free property is the 
following.


\begin{mythm} \label{thm:unary-contexts-fp-bicat}
For any unary $\langCart$-signature $\sig$, the bicategory 
$\urestrict{\syncloneCart{\sig}}$ constructed by restricting 
Construction~\ref{constr:langcart-syntactic-model} to the type theory 
$\urestrict{\langCart}$
is the free fp-bicategory on $\sig$, in the 
sense of Lemma~\ref{lem:free-fp-bicat-from-free-biclone}.
\begin{proof}
For any fp-bicategory $\fpBicat{\altCat}$ and $\langCart$-signature homomorphism
$h : \sig \to \altCat$ the extension fp-pseudofunctor
$\ext{h} : \urestrict{\syncloneCart{\sig}} \to \altCat$ is defined inductively 
as in Proposition~\ref{propn:synclonecart-free-property}, with the following 
adjustments:
\begin{align*}
%
h\sem{x : A \vdash x : A} &:= \Id_{h\sem{A}} \\ 
h\sem{z : Z \vdash \hcomp{t}{x \mapsto u} : B} &:= 	
		{h\sem{x : A \vdash t : B}} \circ
		{h\sem{z : Z \vdash u: A}} \\ 
h\sem{x : A \vdash \pair{\ind{t}} : \prodop_m(B_1, \,\dots\, , B_m)} &:= 
	\seq{h\sem{x : A \vdash t_1 : B_1}, \,\dots\, , 
			h\sem{x : A \vdash t_m : B_m}} \\
%
h\sem{z : Z \vdash \hcomp{\tau}{\sigma} 
	: \rewrite{\hcomp{t}{u}}{\hcompsmall{t'}{u'}} : B} &:= 
	{h\sem{x : A \vdash \tau : \rewrite{t}{t'} : B}} \circ 
		{h\sem{z : Z \vdash \sigma : \rewrite{u}{u'} : A}}
\end{align*} 
\end{proof}
\end{mythm}



\begin{myremark} \label{rem:fp:order-matter-for-type-theory}
As with the construction of $\freeCartBicat{\sig}$, it is important that we first 
restrict $\langCart$ to unary contexts, then 
construct the syntactic model (recall 
Example~\ref{exmp:nucleus-of-fp-not-free}).
\end{myremark}

In the semantics of the simply-typed lambda calculus it is common to restrict
the syntactic model to unary contexts in order to achieve the desired universal property~(see~\eg~\cite[Chapter~4]{Crole1994}).
Hence, we are still justified in calling $\langCart$ 	the internal language of fp-bicategories. 

\subsection{\texorpdfstring{Reasoning within $\langCart$}{Reasoning with 
type-theoretic products}}

In later chapters we shall reason within $\langCart$---and its extension 
$\langCartClosed$ for cartesian closed bicategories---to prove various 
properties of the syntactic models and their semantic interpretation. We 
collect together some results to simplify such calculations.

All the rules of the triangle-law approach to defining products are derivable. 
For example, from Lemma~\ref{lem:fp:trans-ump} one recovers the functoriality 
of the tupling operation and the unit-counit presentation of products 
(see Figure~\ref{fig:fp:admissible-rules}). These derived rules should be 
compared to the primitive rules of~\cite{Seely1987, Hilken1996}.

\setlength{\floatsep}{5pt plus 1.0pt minus 2.0pt}
 
\begin{figure*}[!h]
\begin{minipage}{\textwidth}

\centering
{
\begin{minipage}{\textwidth}
\begin{mdframed}
\begin{minipage}{\textwidth}
\centering 

{\small
\unaryRule{(\Gamma \vdash \id_{t_i} : \rewrite{t_i}{t_i} : A_i)_{i = 1, \dots, n}}
{\Gamma \vdash \pair{\id_{t_1}, \dots, \id_{t_n}} \equiv \id_{\pair{t_1, \dots, t_n}} : \rewrite{\pair{t_1, \dots, t_n}}{\pair{t_1, \dots, t_n}} : \prodop_n(A_1, \dots, A_n)}
{}
}

{\small
\binaryRule{(\Gamma \vdash \tau_i' : \rewrite{t_i'}{t_i''} : A_i)_{i = 1, \dots, n}}
{(\Gamma \vdash \tau_i : \rewrite{t_i}{t_i'} : A_i)_{i = 1, \dots, n}}
{\Gamma \vdash \pair{\tau_1', \dots, \tau_n'} \vert \pair{\tau_1, \dots, \tau_n} \equiv \pair{\tau_1' \vert \tau_1, \dots, \tau_n' \vert \tau_n} : \rewrite{\pair{\ind{t}}}{\pair{\ind{t}''}}: \prodop_n(\ind{A})}
{}}
\end{minipage}


\begin{minipage}{\textwidth}
\centering 

\begin{small}
\unaryRule	{\Gamma \vdash \sigma : \rewrite{u}{u'} : \prodop_n(A_1, \dots, A_n)}
{\Gamma \vdash \etaTimes{u'} \vert \sigma \equiv \pair{\hcomp{\pi_1}{\sigma}, \dots, \hcomp{\pi_n}{\sigma}} \vert \etaTimes{u} : \rewrite{u}{\pair{\hcomp{\ind{\pi}}{u'}}} : \prodop_n(A_1, \dots, A_n)}
{$\etaTimes{}$-nat}
\end{small}

{\small
\unaryRule	{(\Gamma \vdash \tau_i : \rewrite{t_i}{t_i'} : A_i)_{i = 1, \dots, n}}
{\Gamma \vdash \epsilonTimesInd{k}{t_1', \dots, t_n'} \vert \hcomp{\pi_k}{\pair{\tau_1, \dots, \tau_n}} \equiv \tau_k \vert \epsilonTimesInd{k}{t_1, \dots, t_n}  : \rewrite{\hcomp{\pi_k}{\pair{\ind{t}}}}{t_k} : A_k }
{$\epsilonTimesInd{k}{}$-nat $(1 \leq k \leq n)$}
}
\end{minipage}

\begin{minipage}{\textwidth}
\centering 
{\small
\unaryRule{\Gamma \vdash \pair{t_1, \dots, t_n} : \prodop_n(A_1, \dots, A_n)}
{\Gamma \vdash \pair{\epsilonTimesInd{1}{\ind{t}}, \dots, \epsilonTimesInd{n}{\ind{t}}} \vert \etaTimes{\pair{\ind{t}}} \equiv \id_{\pair{\ind{t}}} : \rewrite{\pair{\ind{t}}}{\pair{\ind{t}}} : \prodop_n(\ind{A})}
{triangle-law-1}
}
\unaryRule{\Gamma \vdash \hcomp{\pi_k}{u} : A_k}
{\Gamma \vdash \epsilonTimesInd{k}{t_1, \dots, t_n} \vert \hcomp{\pi_k}{\etaTimes{u}} \equiv \id_{\hcomp{\pi_k}{u}} : \rewrite{\hcomp{\pi_k}{u}}{\hcomp{\pi_k}{u}} : A_k}
{triangle-law-2 $(1 \leq k \leq n)$}
\end{minipage}

{
\caption{Admissible rules for $\langCart(\graph)$\label{fig:fp:admissible-rules}}
}
\end{mdframed}
\vspace{1.5em}
\end{minipage}
}
\end{minipage}
\end{figure*}

\begin{mylemma} \label{lem:langcart-admissible-rules}
For any $\langCart$-signature $\sig$, the rules of 
Figure~\ref{fig:fp:admissible-rules} 
are all admissible.
\begin{proof}
The proofs are all similar; we prove naturality of $\etaTimes{}$ as an example 
of equational reasoning in $\langCart(\sig)$. One can either use the universal 
property (Lemma~\ref{lem:fp:trans-ump}) or reason directly using both the 
equational rules~\rulename{U1} and~\rulename{U2}. We opt for the former. Let 
$\Gamma \vdash \sigma : \rewrite{u}{u'} : 
	\prodop_n(A_1, \,\dots\, ,A_n)$
be any rewrite. Then for $k=1,\dots, n$:
\begin{align*}
\epsilonTimesInd{k}{\ind{\pi}{u'}} \vert \hcomp{\pi_k}{\etaTimes{u'} \vert 
\sigma} &\equiv \epsilonTimesInd{k}{\ind{\pi}{u'}} \vert 
\hcomp{\pi_k}{\etaTimes{u'}} \vert \hcomp{\pi_k}{\sigma} \\
		&\overset{\rulename{U1}}{\equiv} \id_{\hcomp{\pi_k}{u}} \vert 
		\hcomp{\pi_k}{\sigma} \\
		&\equiv \hcomp{\pi_k}{\sigma}
\end{align*}
\begin{align*}
\epsilonTimesInd{k}{\ind{\pi}{u'}} \vert 
\hcomp{\pi_k}{\pair{\hcomp{\pi_1}{\sigma}, \,\dots\, , \hcomp{\pi_n}{\sigma}} 
\vert 
\etaTimes{u}} &\equiv \epsilonTimesInd{k}{\ind{\pi}{u'}} \vert 
\hcomp{\pi_k}{\pair{\hcomp{\pi_1}{\sigma}, \,\dots\, , \hcomp{\pi_n}{\sigma}}} 
\vert 
\hcomp{\pi_k}{\etaTimes{u}} \\ 
	&\overset{\rulename{U1}}{\equiv} \hcomp{\pi_k}{\sigma} \vert 
	\epsilonTimesInd{k}{\hcomp{\ind{\pi}}{u}} \vert \hcomp{\pi_k}{\etaTimes{u}} 
	\\
	&\equiv \hcomp{\pi_k}{\sigma}  
\end{align*}
Applying the universal property of $\transTimes{\hcomp{\pi_1}{\sigma}, 
\,\dots\, , 
\hcomp{\pi_n}{\sigma}}$, one sees that 
$\etaTimes{u'} \vert \sigma \equiv 
\pair{\hcomp{\pi_1}{\sigma}, \,\dots\, , \hcomp{\pi_n}{\sigma}}$, as required.
\end{proof}
\end{mylemma}

\FloatBlock

We also give the 
syntactic 
constructions of the 2-cells 
$\postName$ and $\fuse$ 
(recall~Construction~\ref{constr:fuse-post-semantically} on page~\pageref{constr:fuse-post-semantically}). 
Intuitively, $\postName$ witnesses the identity 
$\seqlr{t_1, \,\dots\, , t_n}[u_i / x_i] = 
	\seqlr{ t_1[u_i / x_i], \,\dots\, , t_n[u_i / x_i]}$ for 
capture-avoiding substitution in the simply-typed lambda calculus.

\begin{myconstr} \label{constr:post}
Let $\sig$ be a $\langCart$-signature. Define a 
2-cell 
$\postName$ in $\langCart(\sig)$ with typing
\begin{center}
\binaryRule{x_1 : A_1, \,\dots\, , x_n \vdash \pair{t_1, \,\dots\, , t_m} : 
\prodop_m(B_1, \,\dots\, , B_m)}
			{(\Delta \vdash u_i : A_i)_{i=1,\dots,n}}
			{\Delta \vdash \post{\ind{t}; \ind{u}} : \rewrite{\hcomp{\pair{t_1, 
			\dots, t_m}}{u_i}}{\pair{\hcomp{t_1}{u_i}, \,\dots\, , 
			\hcomp{t_m}{u_i}}} : \prodop_m(B_1, \,\dots\, , B_m)  }
			{} \vspace{-.5\treeskip}
\end{center}
by setting $\post{\ind{t}; \ind{u}} := \transTimes{\alpha_1, \,\dots\, , 
\alpha_m}$ 
where
\[
\alpha_k := \hcomp{\pi_k}{\hcomp{\pair{t_1, \,\dots\, , t_m}}{u_i}} 
\XRA{\assoc{}^{-1}} \hcompthree{\pi_k}{\pair{t_1, \,\dots\, , t_m}}{u_i} 
\XRA{\hcomp{\epsilonTimesInd{k}{}}{u_i}} \hcomp{t_k}{u_i} 
\] 
Also define a 2-cell $\fuse$ with signature
\begin{small}
\begin{center}
\binaryRule{(x_i : A_i \vdash t_i : A_i)_{i=1, \,\dots\, , n}}{(\Delta \vdash 
u_i : 
A_i)_{i=1, \,\dots\, , n}}{\Delta \vdash \fuse(\ind{t}; \ind{u}) : 
\rewrite{\hcomp{\pair{\hcomp{\ind{t}}{\ind{\pi}(p)}}}{\pair{u_1, \,\dots\, , 
u_n}}}{\pair{\hcomp{t_1}{u_1}, \,\dots\, , \hcomp{t_n}{u_n}}} : \prodop_n(B_1, 
\dots, B_n)}{} \vspace{-.5\treeskip}
\end{center}
\end{small}
by setting $\fuse(\ind{t}; \ind{u}) := \transTimes{\beta_1, \,\dots\, , 
\beta_n}$ 
for $\beta_k$ the composite
\begin{td}
\hpi{k}{\hcomp{\pair{\hcomp{\ind{t}}{\ind{\pi}(p)}}}{\pair{u_1, \,\dots\, , 
u_n}}} 
\arrow[swap, Rightarrow]{d}{\assoc{}^{-1}} \arrow[Rightarrow]{r}{\beta_k} &
\hcomp{t_k}{u_k} \\

\hcompthree{\pi_k}{\pair{\hcomp{\ind{t}}{\ind{\pi}(p)}}}{\pair{u_1, \,\dots\, , 
u_n}} \arrow[swap, Rightarrow]{d}{\hcomp{\epsilonTimesInd{k}{}}{\pair{u_1, 
\dots, u_n}}} & 
\: \\

\hcomp{\hcomp{t_k}{\pi_k(p)}}{\pair{u_1, \,\dots\, , u_n}} \arrow[swap, 
Rightarrow]{r}{\assoc{}} &
\hcomp{t_k}{\hpi{k}{\pair{u_1, \,\dots\, , u_n}}} \arrow[swap, 
Rightarrow]{uu}{\hcomp{t_k}{\epsilonTimesInd{k}{}}}
\end{td}
\end{myconstr}

Since they are defined by applying the universal property to
rewrites that are both natural and invertible, it follows that $\postName$ and 
$\fuse$ are also invertible, as well as being natural in 
the sense that the following rules are admissible:
\vspace{-2mm}
\begin{center}
\binaryRule{(x_1 : A_1, \,\dots\, , x_n : A_n \vdash \tau_j : 
\rewrite{t_j}{t_j'} : B_j)_{j=1,\dots, m}}
		{(\Delta \vdash \sigma_i : \rewrite{u_i}{u_i'}  : A_i)_{i=1,\dots,n}}
		{\Delta \vdash \post{\ind{t}'; \ind{u}'} \vert 
		\hcomp{\pair{\ind{\tau}}}{\sigma_i} \equiv 
		\pair{\hcomp{\ind{\tau}}{\sigma_i}} \vert \post{\ind{t};\ind{u}} : 
		\rewrite{\hcomp{\pair{\ind{t}}}{u_i}}{\pair{\hcomp{\ind{t}'}{u_i'}}}: 
		\prodop \ind{B}}{} \vspace{-1.5\treeskip}

\begin{prooftree}
\AxiomC{$(x_i : A_i \vdash \tau_i : \rewrite{t_i}{t_i'} : A_i)_{i=1, \,\dots\, 
, n}$}
\AxiomC{$(\Delta \vdash \sigma_i : \rewrite{u_i}{u_i'} : A_i)_{i=1, \,\dots\, , 
n}$}
\BinaryInfC{$\Delta \vdash \fuse(\ind{t}'; \ind{u}') \vert 
\hcomp{\pair{\hcomp{\ind{\tau}}{\ind{\pi}(p)}}}{\pair{\ind{\sigma}}} \equiv 
\pair{\hcomp{\ind{\tau}}{\ind{\sigma}}} \vert \fuse(\ind{t}; \ind{u}) : 
$\hspace{35mm}}
\noLine
\UnaryInfC{\hspace{40mm}$: 
\rewrite{\hcomp{\pair{\hcomp{\ind{t}}{\ind{\pi}(p)}}}{\pair{u_1, \,\dots\, , 
u_n}}}{\pair{\hcomp{t_1'}{u_1'}, \,\dots\, , \hcomp{t_n'}{u_n'}}} : 
\prodop_n\ind{B}$}
\end{prooftree}
\end{center}
\normalsize
Moreover, 
the proofs of Lemma~\ref{lem:PseudoproductCanonical2CellsLaws} translate 
readily to the type theory. 

\newpage
\begin{mylemma}  \label{lem:properties-of-post}
Let $\Gamma := (x_i : A_i)_{i=1,\dots,n}$ and $\Delta := (y_l : B_l)_{l = 
1,\dots, k}$ be contexts and suppose $(\Delta \vdash \sigma_i : 
\rewrite{u_i}{u_i'} : A_i)_{i=1, \,\dots\, , n}$. Then
\begin{enumerate}
\item (Naturality). If  $(\Gamma \vdash \tau_j : \rewrite{t_j}{t_j'} : 
B_j)_{j=1,\dots,m}$, then 
%
{\tikzcdset{arrow style=tikz, arrows={Rightarrow}}
\begin{td}[column sep = 3em]
\hcomp{\pair{t_1, \,\dots\, , t_m}}{\ind{u}} \arrow{r}{\postName} 
\arrow[swap]{d}{\hcomp{\pair{\tau_1, \,\dots\, , \tau_m}}{\ind{\sigma}}} &
\pair{\hcomp{t_1}{\ind{u}}, \,\dots\, , \hcomp{t_m}{\ind{u}}} 
\arrow{d}{\pair{\hcomp{\tau_1}{\ind{\sigma}}, \,\dots\, , 
\hcomp{\tau_m}{\ind{\sigma}}}} \\

\hcomp{\pair{t_1', \,\dots\, , t_m'}}{\ind{u}'} \arrow[swap]{r}{\postName} &
\pair{\hcomp{t_1'}{\ind{u}'}, \,\dots\, , \hcomp{t_m'}{\ind{u}'}} 
\end{td}
\item (Compatibility with $\subid{}$). If $(\Gamma \vdash t_m : 
B_m)_{j=1,\dots,m}$ then 
\begin{td}
\pair{t_1, \,\dots\, , t_m} \arrow[swap]{dr}{\pair{\subid{}, \,\dots\, , 
\subid{}}} 
\arrow{r}{\subid{}} &
\hcomp{\pair{t_1, \,\dots\, , t_m}}{\ind{x}} \arrow{d}{\postName} \\
\: &
\pair{\hcomp{t_1}{\ind{x}}, \,\dots\, , \hcomp{t_m}{\ind{x}}}
\end{td} 

\item (Compatibility with $\assoc{}$). \label{c:post-and-assoc} For terms 
$(\Gamma \vdash t_m : C_m)_{j=1,\dots,m}$ and $(\Sigma \vdash v_l : 
B_l)_{l=1,\dots, k}$ then
\begin{td}[column sep = 4em]
\hcompthree{\pair{t_1, \,\dots\, , t_m}}{\ind{u}}{\ind{v}} 
\arrow[swap]{dd}{\assoc{}} \arrow{r}{\hcomp{\postName}{\ind{v}}} &
\hcomp{\pair{\hcomp{t_1}{\ind{u}}, \,\dots\, , \hcomp{t_m}{\ind{u}}}}{\ind{v}}  
\arrow{d}{\postName} \\
\: &
\pair{\hcompthree{t_1}{\ind{u}}{\ind{v}}, \,\dots\, , 
\hcompthree{t_m}{\ind{u}}{\ind{v}}} \arrow{d}{\pair{\assoc{}, \,\dots\, , 
\assoc{}}} 
\\
\hcomp{\pair{t_1, \,\dots\, , t_m}}{\hcomp{\ind{u}}{\ind{v}}} 
\arrow[swap]{r}{\postName} &
\pair{\hcomp{t_1}{\hcomp{\ind{u}}{\ind{v}}}, \,\dots\, , 
\hcomp{t_m}{\hcomp{\ind{u}}{\ind{v}}}} 
\end{td}

\item (Compatibility with $\etaTimes{}$). If $\Gamma \vdash t : \prodop_m(B_1, 
\dots, B_m)$ then
\begin{td}
\hcomp{t}{\ind{u}} \arrow{r}{\hcomp{\etaTimes{}}{\ind{u}}} 
\arrow[swap]{d}{\etaTimes{}} &
\hcomp{\pair{\hpi{1}{t}, \,\dots\, , \hpi{m}{t}}}{\ind{u}} \arrow{d}{\postName} 
\\
\pair{\hpi{1}{\hcomp{t}{\ind{u}}}, \,\dots\, , \hpi{m}{\hcomp{t}{\ind{u}}}} &
\pair{\hcompthree{\pi_1}{t}{\ind{u}}, \,\dots\, , 
\hcompthree{\pi_m}{t}{\ind{u}}} 
\arrow{l}[yshift=-2mm]{\pair{\assoc{}, \,\dots\, , \assoc{}}}
\end{td}
}
\end{enumerate}
\begin{proof}
The proofs are straightforward calculations using the universal property of 
Lemma~\ref{lem:fp:trans-ump}. For example, for naturality we simply observe 
that 
\begin{align*}
\epsilonTimesInd{k}{\hcomp{t_1'}{\ind{u}'}, \,\dots\, , \hcomp{t_m'}{\ind{u}'}} 
&\vert \hpi{k}{\pair{\hcomp{\tau_1}{\ind{\sigma}}, \,\dots\, , 
\hcomp{\tau_m}{\ind{\sigma}}} \vert \post{\ind{t}; \ind{u}}} \\
		&=  \epsilonTimesInd{k}{\hcomp{t_1'}{\ind{u}'}, \,\dots\, , 
		\hcomp{t_m'}{\ind{u}'}} \vert 
		\hpi{k}{\pair{\hcomp{\tau_1}{\ind{\sigma}}, \,\dots\, , 
		\hcomp{\tau_m}{\ind{\sigma}}}} \vert \hpi{k}{\post{\ind{t}; \ind{u}}} 
		\\ 
		&= \hcomp{\tau_k}{\ind{\sigma}} \vert \epsilonTimesInd{k}{t_1, 
		\,\dots\, , 
		t_m} \vert \hpi{k}{\post{\ind{t}; \ind{u}}} \\
		&= \hcomp{\tau_k}{\ind{\sigma}} \vert \hcomp{\epsilonTimesInd{k}{t_1, 
		\dots, t_m}}{\ind{u}} \vert \assoc{\pi_k(p); \pair{t_1, \,\dots\, , 
		t_m}; 
		\ind{u}}^{-1} 
\end{align*}
and that
\begin{align*}
\epsilonTimesInd{k}{\hcomp{t_1'}{\ind{u}'}, \,\dots\, , \hcomp{t_m'}{\ind{u}'}} 
&\vert \hpi{k}{\post{\ind{t}'; \ind{u}'} \vert \hcomp{\pair{\tau_1, \,\dots\, , 
\tau_m}}{\ind{\sigma}}} \\
		&= \epsilonTimesInd{k}{\hcomp{t_1'}{\ind{u}'}, \,\dots\, , 
		\hcomp{t_m'}{\ind{u}'}} \vert \hpi{k}{\post{\ind{t}'; \ind{u}'}} \vert 
		\hpi{k}{\hcomp{\pair{\tau_1, \,\dots\, , \tau_m}}{\ind{\sigma}}}  \\
		&= \hcomp{\epsilonTimesInd{k}{t_1',\dots, t_m'}}{\ind{u}'} \vert 
		\assoc{\pi_k(p); \pair{t_1', \,\dots\, , t_m'}; \ind{u}'}^{-1} \vert 
		\hpi{k}{\hcomp{\pair{\tau_1, \,\dots\, , \tau_m}}{\ind{\sigma}}} \\
		&=  \hcomp{\epsilonTimesInd{k}{t_1',\dots, t_m'}}{\ind{u}'} \vert 
		\hcompthree{\pi_k}{\pair{\tau_1, \,\dots\, , \tau_m}}{\ind{\sigma}} 
		\vert 
		\assoc{\pi_k(p); \pair{t_1, \,\dots\, , t_m}; \ind{u}}^{-1}  \\
		&= \hcomp{\tau_k}{\ind{\sigma}} \vert \hcomp{\epsilonTimesInd{k}{t_1, 
		\dots, t_m}}{\ind{u}} \vert \assoc{\pi_k(p); \pair{t_1, \,\dots\, , 
		t_m}; 
		\ind{u}}^{-1} 
\end{align*}
Hence, by the universal property of Lemma~\ref{lem:fp:trans-ump}, the required 
equality holds. The other cases are similar.
\end{proof}
\end{mylemma}

\subsection{Products from context extension} 
\label{sec:products-from-context-extension}

We end this chapter by noting a `degenerate' or `implicit' way for a deductive 
system to 
exhibit product structure. The construction gives rise to a syntactic model 
that is an fp-bicategory, but does not arise via a cartesian biclone or 
provide a type-theoretic description of bicategorical products. 
While this structure is not in the vein of those we have discussed above, 
it will play an important role: exponentials in the simply-typed lambda 
calculus are 
defined with respect to these products. The product structure is given by 
context concatenation.  

\begin{myconstr} \label{constr:context-cart-termcat} 
For any $\langCart$-signature $\sig$, define a bicategory 
$\termCatContextExt(\sig)$ as 
follows. Fix an enumeration of variables $x_1, \,\dots\, , x_n, \dots$. The 
objects 
are then contexts $\Gamma, \Delta, \dots$ in which the $i$th entry has variable 
name $x_i$. The 1-cells $\Gamma \to (y_j : B_j)_{j = 1, \,\dots\, , m}$ are 
$m$-tuples of \mbox{$\alpha$-equivalence} classes of terms $(\Gamma \vdash t_j 
:B_j)_{j= 1, \,\dots\, , m}$ derivable in $\langCart(\sig)$; the 
\mbox{2-cells}\hide{\mbox{$(\Gamma \vdash t_j :B_j)_{j= 1, \,\dots\, , m} \To 
(\Gamma \vdash t_j' :B_j)_{j= 1, \,\dots\, , m}$}} are $m$-tuples of 
\mbox{$\alpha{\equiv}$-equivalence} classes of rewrites \mbox{$(\Gamma \vdash 
\tau : \rewrite{t_j}{t'_j} : B_j)_{j= 1, \,\dots\, , m}$}. 

Vertical composition is given pointwise by the $\vert$ operation, and 
horizontal composition by explicit substitution: 
\begin{align*} 
(t_1, \,\dots\, , t_l), (u_1, \,\dots\, , u_m) 
	&\mapsto (\hcomp{t_1}{x_i \mapsto u_i}, \dots, \hcomp{t_m}{x_i \mapsto u_i}) \\ 
(\tau_1, \,\dots\, , \tau_l), (\sigma_1, \,\dots\, , \sigma_m) 
	&\mapsto 
(\hcomp{\tau_1}{x_i \mapsto \sigma_i}, \,\dots\, , \hcomp{\tau_m}{x_i \mapsto 
\sigma_i})
\end{align*} 
The identity on $\Delta = (y_j : B_j)_{j = 1, \,\dots\, , m}$ is the 
\rulename{var} 
rule \mbox{$(\Delta \vdash y_j : B_j)_{j = 1, \,\dots\, , m}$}, and the structural 
isomorphisms $\l, \r$ and $\a$ are given pointwise by $\proj{}$, 
$\subid{}^{-1}$ and $\assoc{}$, respectively. 
\end{myconstr} 

Since $\langCart$ comes equipped with a product structure, this bicategory has 
two product structures: 
one given by the product 
structure in the type theory, and the other by context extension. We emphasise 
this with the notation.

The type-theoretic  
product structure is induced from that on the full sub-bicategory of unary 
contexts via the 
following lemma, which 
can be seen as the type-theoretic translation of 
Lemma~\ref{lem:biclone-with-products-adjoint-eq} on page~\pageref{lem:biclone-with-products-adjoint-eq}.

\begin{mylemma} \label{lem:fp-bicat:equivalences} For any 
$\langCartClosed$-signature 
$\sig$ and context 
$\Gamma = (x_i : A_i)_{i=1, \,\dots\, ,n}$, there exists an adjoint equivalence 
$\Gamma \leftrightarrows \big(p : \prodop_n(A_1, \,\dots\, , A_n)\big)$
in $\termCatContextExt(\sig)$.
\begin{proof}
Take the 1-cells
\begingroup
\addtolength{\jot}{.5em}
\begin{align*}
\left(\Gamma \vdash \pair{x_1, \,\dots\, , x_n} : \prodop_n(A_1, \,\dots\, , 
A_n)\right) 
: 
\Gamma 
\to (p : \prodop_n(A_1, \,\dots\, , A_n)) \\
\left(p : \prodop_n(A_1, \,\dots\, , A_n) 
	\vdash \pi_i(p) : A_i \right)_{i=1,\dots,n} 
:  (p : \prodop_n(A_1, \,\dots\, , A_n)) \to \Gamma
\end{align*}
\endgroup
For the unit and counit of the 
required adjoint equivalence we take 
\[
\left( \Gamma \vdash \epsilonTimesInd{i}{\ind{x}}
	\rewrite{\hpi{i}{\pair{x_1, \dots, x_n}}}{x_i} : A_i\right)_{i=1, \dots, n}
\]
and the composite
\begin{td}
p
\arrow[]{rr}{} 
\arrow[swap]{d}{\etaTimes{p}} &
\: &
\hcomp{\pair{x_1, \dots, x_n}}{\pi_i(p)} \\
\pair{\hcomp{\pi_1}{p}, \,\dots\, , \hcomp{\pi_n}{p}}  
\arrow[swap]{r}[yshift=-2mm]
	{\pair{\indproj{-1}{}, \,\dots\, , \indproj{-n}{}}} &
\pair
	{\hcomp{x_1}{\hcomp{\ind{\pi}}{p}}, 
		\dots, 
	 \hcomp{x_n}{\hcomp{\ind{\pi}}{p}}
	} 
\arrow[swap]{r}[yshift=-2mm]{\post{\ind{x}; \hcomp{\ind{\pi}}{p}}^{-1}} & 
\hcomp{\pair{x_1, \,\dots\, , x_n}}{\hcomp{\ind{\pi}}{p}}
\arrow[swap]{u}
	{\hcomp{\pair{x_1, \,\dots\, , x_n}}{\subid{\ind{\pi}(p)}^{-1}}}
\end{td}
The proof then amounts to making use of naturality to the 
point where one can apply the triangle laws of 
Figure~\ref{fig:fp:admissible-rules}.
\end{proof}
\end{mylemma} 

\begin{myremark} \label{rem:syntactic-cc-bicategories-biequivalent}
The preceding lemma, together with 
Lemma~\ref{lem:subbicategory-equivalence-biequivalence}
on page~\pageref{lem:subbicategory-equivalence-biequivalence}, in fact entails that 
$\termCatContextExt(\sig) \simeq \urestrict{\syncloneCart{\sig}}$ for every 
unary
$\langCart$-signature $\sig$. 
\end{myremark}

We define the product 
$(x_i^{(1)} : 
	A_i^{(1)})_{i = 1, \,\dots\, , m_1} \times \dots \times (x_i^{(n)} : 
		A_i^{(n)})_{i = 1, \,\dots\, , m_n}$ 
of arbitrary contexts
to be the product 
$(p_1 : \prodop_{i=1}^{m_1} A_i^{(1)}) 
	\times 
	\dots 
	\times 
(p_n : \prodop_{i=1}^{m_n} A_i^{(n)})$ 
of the corresponding unary contexts. The $i$th projection is the 
$\len{\Gamma^{(i)}}$-tuple 
\begin{equation} \label{eq:context-ext-type-theoretic-product-proj}
\left(p : \prodop_n \big( \prodop_{\len{\Gamma^{(1)}}} \ind{A}^{(1)}, 
\,\dots\, , 
		\prodop_{\len{\Gamma^{(n)}}} \ind{A}^{(n)} \big) \vdash 
		\hpi{j}{\pi_i(p)} : A^{(i)}_j\right)_{j=1, \dots, \len{\Gamma^{(i)}}} 
\end{equation}
and the tupling 
of $n$ maps $(\Delta \to \Gamma^{(i)})_{i=1, \dots, n}$, that is, of 
$\len{\Gamma^{(i)}}$-tuples 
$(\Delta \vdash t_j^{(i)} : A_j^{(i)})_{
	\substack{j=1, \dots, \len{\Gamma^{(i)}} \\ 
				i = 1, \dots, n}}$, 
is 
\[
\Delta \vdash \pairlr{ \pair{\ind{t}^{(1)}}, \dots, \pair{\ind{t}^{(n)}} } : 
	\prodop_n 
		\big( \prodop_{\len{\Gamma^{(1)}}} \ind{A}^{(1)}, 
		\,\dots\, , 
		\prodop_{\len{\Gamma^{(n)}}} \ind{A}^{(n)} \big)
\]
The counit $\epsilonTimesInd{i}{}$ is the composite indicated by the pasting 
diagram
\begin{td}[column sep = 8em, row sep = 3em]
\prodop_n 
		\big( \prodop_{\len{\Gamma^{(1)}}} \ind{A}^{(1)}, 
		\,\dots\, , 
		\prodop_{\len{\Gamma^{(n)}}} \ind{A}^{(n)} \big)
\arrow{r}{\pi_i(p)} 
\arrow[phantom]{dr}[description, near start, xshift=-3mm]
	{\twocellIso{\epsilonTimesInd{i}{}}} &
\prodop_{\len{\Gamma^{(i)}}} \ind{A}^{(i)}
\arrow{r}{(\pi_1(p), \dots, \pi_{\len{\Gamma^{(i)}}}(p))}
\arrow[phantom, near end]{d}{\iso} &
\Gamma^{(i)} \\

\Delta
\arrow{u}{ \pairlr{ \pair{\ind{t}^{(1)}}, \dots, \pair{\ind{t}^{(n)}} }}
\arrow[swap, bend right = 12]{ur}{\pair{\ind{t}^{(i)}}}
\arrow[bend right, swap]{urr}{t_1^{(i)}, \dots, t_{\len{\Gamma^{(i)}}}^{(i)}} &
\: &
\: 
\end{td}
That is, the $\len{\Gamma^{(i)}}$-tuple with $j$th component the composite 
rewrite
\begin{equation*} 
\begin{tikzcd}[column sep = 6em]
\hcompthree
	{\pi_j}
	{\pi_i(p)}
	{\pairlr{ \pair{\ind{t}^{(1)}}, \dots, \pair{\ind{t}^{(n)}} }}
\arrow{r}
\arrow[swap]{d}{\iso} &
t_j^{(i)} \\
\hpi{j}{\hpi{i}{\pairlr{ \pair{\ind{t}^{(1)}}, \dots, \pair{\ind{t}^{(n)}} }}}
\arrow[swap]{r}
	{\hpi{j}
		{\epsilonTimesInd{i}{}}} &
\hpi{j}{\pair{t_1^{(i)}, \dots, t_{\len{\Gamma^{(i)}}}^{(i)} }}
\arrow[swap]{u}
	{\epsilonTimesInd{j}{}}
\end{tikzcd}
\end{equation*}
The next lemma encapsulates the required universal property. 

\begin{mylemma} \label{lem:context-extension-bicat-type-product-structure}
For any unary $\langCart$-signature $\sig$, the 1-cell
\[
\left(p : \prodop_n \big( \prodop_{\len{\Gamma^{(1)}}} \ind{A}^{(1)}, 
\,\dots\, , 
		\prodop_{\len{\Gamma^{(n)}}} \ind{A}^{(n)} \big) \vdash 
		\hpi{j}{\pi_i(p)} : A^{(i)}_j\right)_{j=1, \dots, \len{\Gamma^{(i)}}} 
\]
of~(\ref{eq:context-ext-type-theoretic-product-proj}) is a biuniversal arrow 
defining an fp-structure on $\termCatContextExt(\sig)$.
\begin{proof}
Taking the structure described above, it remains to check the universal 
property of the counit. Suppose that 
$\Delta \vdash u : \big( \prodop_{\len{\Gamma^{(1)}}} \ind{A}^{(1)}, 
		\,\dots\, , 
		\prodop_{\len{\Gamma^{(n)}}} \ind{A}^{(n)} \big)$
and that 
$(\Delta \vdash t_j^{(i)} : A_j^{(i)})_{j=1, \dots, \len{\Gamma^(i)}}$
for $i=1, \dots, n$, and
consider a family of rewrites
\[\left(\Delta \vdash \alpha_j^{(i)} : 
	\rewrite{\hcompthree{\pi_j}{\pi_i(p)}{u}}
			{t_j^{(i)}} : A_j^{(i)} \right)_{
	\substack{j=1, \dots, \len{\Gamma^{(i)}} \\ i= 1, \dots, n }}
\]
One thereby obtains composites
$
\widetilde{\alpha}_j^{(i)} :=
	\hcomp{\pi_j}{\hpi{i}{u}} 
		\XRA\iso
	\hcompthree{\pi_j}{\pi_i(p)}{u}
		\XRA{\alpha_j^{(i)}}
	t_j^{(i)}
$
for $j=1, \dots, \len{\Gamma^{(i)}}$ and $i=1, \dots, n$. Applying the 
universal property of $\epsilonTimes{}$ (Lemma~\ref{lem:fp:trans-ump}) for each 
$i$, one obtains
$\transTimes
	{\widetilde{\alpha}_1^{(i)}, \dots, 
	\widetilde{\alpha}_{\len{\Gamma^{(i)}}}^{(i)}} 
	: \hpi{k}{u} \To 
	\pair{t_1^{(i)}, \dots, t_{\len{\Gamma^{(i)}}}^{(i)}}$
for $i=1, \dots, n$. Finally applying the universal property to this family of 
rewrites, one obtains
\[
\transTimeslr
	{\transTimes
	{\widetilde{\alpha}_1^{(1)}, \dots, 
	\widetilde{\alpha}_{\len{\Gamma^{(1)}}}^{(1)}}, \dots, 
	\transTimes
	{\widetilde{\alpha}_1^{(n)}, \dots, 
	\widetilde{\alpha}_{\len{\Gamma^{(n)}}}^{(n)}}} :
	u \To 
	\pairlr{ \pair{\ind{t}^{(1)}}, \dots, \pair{\ind{t}^{(n)}} }
\]
To see that this 2-cell satisfies the required universal property, apply 
the corresponding property from
Lemma~\ref{lem:fp:trans-ump} twice. 
\end{proof}
\end{mylemma}

We now turn to the second, strict, product structure. This arises from context extension. 
Constructing products in this way is a standard method in the categorical 
setting (\eg~\cite{Pitts2000}) and is also employed by Hilken~\cite{Hilken1996} 
in the 2-categorical case to obtain a strict product. Taken on its own, 
however, it does not enable one to reason about products within the type 
theory. 

\newpage
\begin{mylemma} \label{lem:context-product-structure}
For any $\langCart$-signature $\sig$ the syntactic model 
$\termCatContextExt(\sig)$ of 
$\langCart(\sig)$ is an 
\mbox{fp-bicategory} with product structure given by context extension.
\begin{proof}
We claim first that every context 
$\Gamma := (x_i : A_i)_{i=1, \,\dots\, , n}$ is the $n$-ary product 
$\prod_{i=1}^n (x_i : A_i)$ of unary contexts 
$(x_1 : A_1), \,\dots\, , (x_n : A_n)$.
Define projections 
$\pi_k : \Gamma \to A_k$ for $k =1, \,\dots\, , n$ by $\Gamma \vdash x_k : 
A_k$. 
Then, given 1-cells 
$\Delta \vdash t_i : A_i$ for $i = 1, \,\dots\, , n$, define the $n$-ary 
tupling to 
be the $n$-tuple $(\Delta \vdash t_i : A_i)_{i=1, \,\dots\, , n}$. The unit and 
counit are the 2-cells with components $\indproj{-i}{}$ and $\indproj{i}{}$, 
respectively.

We extend this to all contexts in the obvious way. For contexts $\Gamma_i$ 
$(i = 1, \,\dots\, , n)$ such that $\Gamma_i := (x_j : A_j^{(i)})_{j=1, 
\,\dots\, , 
\len{\Gamma_i}}$ the product $\prod_{i=1}^n \Gamma_i$ is the concatenated 
context $\Gamma_1, \,\dots\, , \Gamma_n$ (the enumeration of variables ensures 
no 
variable names are duplicated). The $k$th projection is the 
$\len{\Gamma_k}$-tuple 
$(\Gamma_1, \,\dots\, , \Gamma_n \vdash x_j : A_j^{(k)})_{1 + \sum_{l=1}^{k-1} 
\len{\Gamma_l} 
\leq j \leq \len{\Gamma_k} + \sum_{l=1}^{k-1} \len{\Gamma_l}}$ and the $n$-ary 
tupling of 
1-cells 
$(\bar{t}_i : \Delta \to \Gamma_i)_{i=1, \,\dots\, , n}$ with $\bar{t}_i := 
(\Delta 
\vdash t_j^{(i)} : A_j^{(i)})_{j=1, \,\dots\, , \len{\Gamma_i}}$ is just the 
unfolded 
$\sum_{i=1}^n \len{\Gamma_i}$-tuple $(\Delta \vdash t_j^{(i)} : 
A_j^{(i)})_{\substack{i 
=1, \,\dots\, , n \\ j=1, \,\dots\, , \len{\Gamma_i}}}$. The unit and counit 
are as in 
the unary case.
\end{proof}
\end{mylemma}

%

\chapter{A type theory for cartesian closed bicategories} \label{chap:ccc-lang}

We now build on the preceding chapters, and the type theory 
$\langCart$, to construct a type theory for cartesian closed 
bicategories. First we extend the theory of clones with finite 
products to include exponentials via a version of Lambek's
\Def{internal hom} of a multicategory~\cite{Lambek1989}. Next we extend this to (cartesian) 
biclones and use it to extract a type theory $\langCartClosed$ for 
which the syntactic model is free among cartesian closed biclones.
The proof of the corresponding bicategorical free property, however, 
throws up a subtlety: 
exponentials in the Lambek style are 
defined as a right (bi)adjoint to \emph{context extension} rather than the 
type-theoretic product. In terms of the syntactic models of the preceding 
chapter, exponentials appear with respect to the context extension product 
structure, rather than the type-theoretic product structure (recall 
Section~\ref{sec:products-from-context-extension}).
As we shall see, it follows that the restriction of $\langCartClosed$ to unary 
contexts cannot 
satisfy a strict free property mirroring that 
of $\langBicat$ and $\langCart$. We address this 
by showing that the syntactic model of $\langCartClosed$ is biequivalent to the 
cartesian closed bicategory enjoying such a strict free property. 
(Table~\ref{table:all-free-constructions} on 
page~\pageref{table:all-free-constructions} provides an index of 
the various free constructions and syntactic models we employ.)
We end the chapter by making precise the claim that $\langCartClosed$ is the 
simply-typed lambda calculus up to isomorphism.

\section{Cartesian closed bicategories}
\label{sec:cc-bicategories}

Let us start by recapitulating the definition of cartesian closed bicategory. To 
give a cartesian closed structure on an fp-bicategory $\fpBicat{\baseCat}$ is 
to specify a biadjunction $(-) \times A \dashv 
(\expobj{A}{-})$ for every $A \in \baseCat$. Following 
Definition~\ref{def:biadjunction-universal-arrow}, this amounts to choosing an 
object $(\expobj{A}{B})$ and a biuniversal arrow $\eval_{A,B} : (\expobj{A}{B}) 
\times A \to B$ for every $A, B \in \baseCat$. We unfold the definition as 
follows. 

\begin{mydefn} \label{def:cc-bicat}
A \Def{cartesian closed bicategory} or \mbox{\Def{cc-bicategory}} is an 
fp-bicategory $\fpBicat{\baseCat}$ equipped with the following data 
for every $A, B \in \baseCat$:
\begin{enumerate}
\item \label{c:cc-bicat-obj} A chosen object $(\expobj{A}{B})$, 
\item \label{c:cc-bicat-map} A specified 1-cell $\eval_{A,B} : (\expobj{A}{B}) 
\times A \to B$, 
\item \label{c:cc-bicat-adj-equiv} For every $X \in \baseCat$, an adjoint 
equivalence 
\begin{equation} \label{diag:cc-bicat}
\begin{tikzcd} 
\baseCat(X, \expobj{A}{B}) \arrow[bend left = 20]{r}{\eval_{A,B} \circ (- 
\times A)} \arrow[phantom]{r}[xshift=0em]{\adjUp{\simeq}} & 
\baseCat(X \times A, B) \arrow[bend left = 20]{l}{\lambda} 
\end{tikzcd} 
\end{equation}
specified by a family of universal arrows 
$\epsilon_f : \eval_{A,B} \circ (\lambda f \times A) \To f$. 
\end{enumerate}
We call the functor $\lambda(-)$ \Def{currying} and refer to $\lambda f$ as the 
\Def{currying of $f$}.
\end{mydefn}

\begin{myremark}
As for products, we shall call an exponential structure \Def{strict} if the 
equivalences~(\ref{diag:cc-bicat}) are isomorphisms. When the underlying 
bicategory $\baseCat$ is a 2-category, this yields the definition of cartesian 
closure in the $\CatCat$-enriched sense
(\cf~Remark~\ref{rem:strictness-of-products-defs}).
\end{myremark}


Explicitly, the equivalences~(\ref{diag:cc-bicat}) are given by the following 
universal property.  For every 
1-cell $t : X \times A \to B$ we require a 1-cell 
$\lambda t : X \to (\expobj{A}{B})$ and an invertible 
2-cell $\epsilonExp_t : \eval_{A,B} \circ (\lambda t \times A) \To t$,
\nom{\epsilonExp_t}
	{The counit for exponential structure, of type
		$\eval_{A,B} \circ (\lambda t \times A) \XRA\iso t$} 
universal in the sense that for any 2-cell 
$\alpha : \eval_{A,B} \circ (u \times A) \To t$ there exists a unique 2-cell 
$\transExp{\alpha} : u \To \lambda t$ 
\nom{\transExp{\alpha}}
	{The unique mediating 2-cell $u \To \lambda t$ corresponding to
		$\alpha : \eval_{A,B} \circ (u \times A) \To t$}
such that 
$\epsilon_t \vert \big(\eval_{A,B} \circ (\transExp{\alpha} \times A)\big) = 
\alpha$. Moreover, we require that the unit 
$\etaExp{t} := \transExp{\id_{\eval_{A,B} \circ (t \times A)}}$
\nom{\etaExp{t}}
	{The unit for exponential structure, of type
		$t \XRA\iso \lambda\left( \eval_{A,B} \circ (t \times A) \right)$}
is also invertible. 

\begin{mynotation}  \label{not:exponentials}
Following the categorical notation, for 1-cells 
$f : A' \to A$ and $g : B \to B'$ we write 
$(\expobj{f}{g}) : (\expobj{A}{B}) \to (\expobj{A'}{B'})$ for the exponential 
transpose of the composite
$(g \circ \eval_{A,B}) \circ (\Id_{\scriptsizeexpobj{A}{B}} \times f)$, thus:
\[
(\expobj{f}{g}) := 
\lambda\big( 
(\expobj{A}{B}) \times A' \xra{(\scriptsizescriptsizeexpobj{A}{B}) \times f} 
(\expobj{A}{B}) \times A \xra{\eval_{A,B}}
B \xra{g} B'
\big)\]
and likewise on 2-cells. 
\end{mynotation} 

As for products, 1-category theoretic notation can be misleading when the 
identity is referred to explicitly. Consider the identities
\begin{align*}
(\expobj{f}{\Id_B}) &= 
\lambda{\left( 
	(\Id_B \circ \eval_{A,B}) \circ (f \times \Id_A)
\right)} \\
(\expobj{\Id_A}{g}) &= 
\lambda{\left( 
	(g \circ \eval_{A,B}) \circ (\Id_{\scriptsizeexpobj{A}{B}}\times \Id_A)
\right)} 
\end{align*}
In a \emph{2-category}
with pseudo-products and pseudo-exponentials, one may safely write 
$(\expobj{f}{\Id_B})$ as simply $\lambda( \eval_{A,B} \circ (f \times A))$, 
but cannot simplify $(\expobj{\Id_A}{g})$ in 
a similar way to $\lambda(g \circ \eval_{A,B})$. Note, however, that this 
simplification 
is possible in the presence of strict products, when the unit is an identity.

\begin{myremark} \label{rem:exponentials-up-to-equiv}
The uniqueness of exponentials up to equivalence manifests itself in the same
way as for products.  For instance, given an adjoint 
equivalence~$e: E \simeq (\expobj{A}{B}): f$, the object $E$ inherits an
exponential structure by composition with $e$ and $f$
(\cf~Remark~\ref{rem:nAryProductsDeterminedUpToEquiv}).
\end{myremark}

In Construction~\ref{constr:fuse-post-semantically} we saw that standard 
properties of cartesian categories are witnessed by natural families of 2-cells 
in an fp-bicategory. The same principle holds for cc-bicategories.  

\begin{myconstr}  \label{constr:push-semantically}
Let $\ccBicat{\baseCat}$ be a cc-bicategory. 
For $g : X \to Y$ and $f : Y \times A \to B$ we define 
$\pushName(f,g) : 
\lambda(f) \circ g \To \lambda\big(f \circ (g \times A)\big)$ 
\nom{\pushName(f,g)}
	{The canonical 2-cell 
		$\lambda(f) \circ g \To \lambda\big(f \circ (g \times A)\big)$}
as 
$\transExp{\tau}$, for $\tau$ the composite
\begin{td}[column sep = 3em]
\eval_{A,B} \circ \left(( \lambda f \circ g ) \times A\right)
\arrow[swap]{d}{\eval \circ (\phiTimes_{f,g})^{-1}}
\arrow{r}{\tau} &
f \circ (g \times A) \\

\eval_{A,B} \circ \left(( \lambda f \times A) \circ (g \times A)\right)
\arrow[swap]{r}{\iso} &
\left( \eval_{A,B} \circ (\lambda f \times A) \right) 
	\circ (g \times A)
\arrow[swap]{u}{\epsilonExp_f \circ (g \times A)}
\end{td}
where $\phiTimes_{f,g} : (f \times A) \circ (g \times A) \To (fg \times A)$ 
witnesses $\prodop_2(-, =)$ as a pseudofunctor 
(recall~Construction~\ref{constr:fuse-post-semantically}(\ref{c:phiTimes-defined})).
 
\end{myconstr}

This family of 2-cells is natural in each of its arguments and satisfies the 
expected equations, some of which are collected in the following lemma. As for 
Lemma~\ref{lem:PseudoproductCanonical2CellsLaws}, we 
assume the underlying bicategory is strict for the sake of clarity.

\begin{prooflesslemma} \label{lem:PseudoCCCCanonical2CellsLaws}
Let $\ccBicat{\baseCat}$ be a 2-category with finite pseudo-products and 
pseudo-exponentials.  Then for all 1-cells 
$f,g$ and $h$, the following diagrams commute 
whenever they are well-typed:
\begin{equation}
\label{c:PushAndPsi}
\begin{tikzcd}
(\lambda f) \circ \Id \arrow[equals]{d} \arrow{r}{\pushName} &
\lambda \big(f \circ (\Id \times A)\big) \\
\lambda f \arrow[swap]{r}{\lambda(f \circ \etaTimes{f})} &
\lambda (f \circ \seq{\pi_1, \pi_2}) 
\arrow[equals]{u}{}
\end{tikzcd}
\end{equation}

\begin{equation} \label{c:PhiAndEtaExp}
\begin{tikzcd} 
f \circ g 
\arrow{r}{\etaExp{f \circ g}} 
\arrow[swap]{d}{\etaExp{f} \circ g} &
\lambda\big(\eval \circ (fg \times A) \big) \\
\lambda \big(\eval \circ (f \times A)\big) \circ g
\arrow[swap]{r}{\pushName} &
\lambda \big(\eval \circ (f \times A) \circ (g \times A) \big) 
\arrow[swap]{u}{\lambda(\eval \circ \phiTimes_{f,g; \Id})}
\end{tikzcd}
\end{equation}

\begin{equation} \label{c:PushAndPhi}
\begin{tikzcd}
(\expobj{f}{g}) \circ \Id \arrow{r}{\pushName} \arrow[equals]{d} &
\lambda\big( g \circ \eval \circ ((\expobj{A}{B}) \times f) \circ (\Id \times 
B)\big) 
\arrow{d}{\lambda(g \circ \eval \circ \phiTimes_{\Id; f, \Id})} \\
(\expobj{f}{g}) \arrow[equals]{r} &
\lambda\big(g \circ \eval \circ ((\expobj{A}{B}) \times f)\big)
\end{tikzcd}
\end{equation}

\begin{equation} \label{c:DoublePost}
\begin{tikzcd}
\lambda (f) \circ g \circ h 
\arrow{r}{\pushName \circ h} 
\arrow[swap]{d}{\pushName} &
\lambda\big(f \circ (g \times A)\big) \circ h 
\arrow{r}{\pushName} &
\lambda\big(f \circ (g \times A) \circ (h \times A)\big) 
\arrow[]{d}{\lambda(f \circ \phiTimes_{g,h; \Id})} \\
\lambda\big(f \circ ((g \circ h) \times A)\big) 
\arrow[equals]{rr} &
\: &
\lambda\big(f \circ (gh \times A)\big) 
\end{tikzcd}
\end{equation}
\end{prooflesslemma}

A pseudofunctor between cartesian closed bicategories is cartesian closed if it preserves both the biuniversal arrows defining products 
and the biuniversal arrows defining exponentials. 

\begin{mydefn}  
A \Def{cartesian closed pseudofunctor} or \Def{cc-pseudofunctor} between 
cc-bicategories $\ccBicat{\baseCat}$ and 
\mbox{$\ccBicat{\altCat}$} is an \mbox{fp-pseudofunctor} 
$(F, \prodPres)$ equipped with specified adjoint equivalences 
\[
\evBar_{A,B}
:
F(\expobj{A}{B})
\leftrightarrows 
(\expobj{FA}{FB}) 
: 
\expPres_{A,B}
\]
for every $A, B \in 
\baseCat$, where $\evBar_{A,B} :  F(\expobj{A}{B}) \to (\expobj{FA}{FB})$ is 
the exponential transpose of 
$F(\eval_{A,B}) \circ \prodPres_{\scriptsizeexpobj{A}{B}, A}$. 
\nom{\evBar_{A,B}}
	{The canonical map $F(\expobj{A}{B}) \to (\expobj{FA}{FB})$ for an 
	\mbox{fp-pseudofunctor} $(F, \prodPres)$, defined
		as the transpose of 
		$F(\eval_{A,B}) \circ \prodPres_{\scriptsizeexpobj{A}{B}, A}$}%
\nom{\expPres_{A,B}}
	{An equivalence 
		$(\expobj{FA}{FB}) \to 
			F\left(\expobj{A}{B} \right)$
	forming part of the data of a cc-pseudofunctor}%
We denote the 2-cells witnessing that $\expPres_{A,B}$ 
and $\evBar_{A,B}$ form an equivalence by
\begin{align*}
\unExp_{A,B} : 
		\Id_{(\scriptsizeexpobj{FA}{FB})} \To 
			\evBar_{A,B} \circ \expPres_{A,B} \\
\coExp_{A,B} : 
		\expPres_{A,B} \circ \evBar_{A,B} \To 
			\Id_{F(\scriptsizeexpobj{A}{B})}
\end{align*}%
\nom{\unExp_{A,B}}
	{A 2-cell $\Id_{(\expobj{FA}{FB})} \To 
				\evBar_{A,B} \circ \expPres_{A,B}$, 
	part of the data of a cc-pseudofunctor $(F, \prodPres, \expPres)$}%
\nom{\coExp_{A,B}}
	{A 2-cell $\expPres_{A,B} \circ \evBar_{A,B} \To 
				\Id_{F(\expobj{A}{B})}$, part 
	of the data of a cc-pseudofunctor $(F, \prodPres, \expPres)$}%
A cc-pseudofunctor $(F, \prodPres, \expPres)$ is \Def{strict} if 
$(F, \prodPres)$ is a strict fp-pseudofunctor such that
\begin{align*}
F(\expobj{A}{B}) &= (\expobj{FA}{FB}) \\
F(\eval_{A,B}) &=\eval_{FA,FB} \\
F(\lambda t) &=\lambda(Ft) \\
F(\epsilonExp_t) &=\epsilonExp_{Ft} \\
\expPres_{A,B} &=\Id_{\scriptsizeexpobj{FA}{FB}}
\end{align*}
with equivalences canonically induced by the \mbox{$2$-cells}
\begin{center}
$\transExp{\eval_{FA, FB} \circ \kappa}
 : \Id_{(\scriptsizeexpobj{FA}{FB})}
   \XRA\iso
   \lambda(\eval_{FA,FB} \circ \Id_{(\scriptsizeexpobj{FA}{FB}) \times FA})$
\end{center}
for $\kappa$ is the canonical isomorphism 
$
\Id_{\scriptsizeexpobj{FA}{FB}} \times FA
\iso
\Id_{(\scriptsizeexpobj{FA}{FB}) \times FA}
$. 
\end{mydefn}

\begin{myremark}[{\cf~Remark~\ref{rem:fp-pseudofunctor-for-different-prod}}] 
\label{rem:cc-pseudofunctor-for-different-struct}
If $\baseCat$ is a bicategory equipped with two cartesian closed 
structures, say 
$\ccBicat{\baseCat}$ and 
$\big(\baseCat, \altprod_n(-), \altexp{-}{-} \big)$, 
then for any cc-pseudofunctor
$(F, \prodPres, \expPres) : \ccBicat{\baseCat} \to \ccBicat\altCat$ 
there exists an (equivalent) cc-pseudofunctor 
\[
\big(\baseCat, \altprod_n(-), \altexp{-}{-} \big) \to \ccBicat{\altCat}
\]
with witnessing equivalences arising from the uniqueness of 
products and exponentials up to equivalence. 
\end{myremark}

\paragraph*{cc-Biequivalences from biequivalences.} In the preceding chapter 
(page~\pageref{def:fp-transformation}) we 
saw that, so far as we are concerned, it is unnecessary to 
distinguish between pseudonatural transformations and their product-respecting 
counterparts. A similar situation holds in the cartesian closed case. For 
cartesian closed pseudofunctors 
$(F, \prodPres, \expPres), 
(G, \altProdPres, \altExpPres) : \ccBicat{\baseCat} \to \ccBicat{\altCat}$, a 
\mbox{\Def{cc-transformation}} $F \To G$ is an 
fp-transformation 
$(\CCCTransNat, \CCCTrans, \CCCTransProd) : 
	(F, \prodPres) \To (G,\altProdPres)$ 
(recall Definition~\ref{def:fp-transformation})
equipped with a 2-cell $\CCCTransExp_{A,B} \:\: (A, B \in \baseCat)$ as 
in the diagram below
\begin{td}[column sep = 5em]
F(\expobj{A}{B}) \times FA 
\arrow[
swap,
rounded corners,
to path=
{ -- ([yshift=0ex]\tikztostart.north)
|- ([yshift=5.5ex]\tikztostart.east)
-| ([yshift=5.5ex]\tikztotarget.south)
-- (\tikztotarget.north)}, 
]{rr}
\arrow[swap, phantom]{rr}[yshift=7.2ex, font=\scriptsize]
	{\eval_{FA,FB} \circ (\evBar^F_{A,B} \times FA)}
\arrow{r}[yshift=0mm]{\evBar^F_{A,B} \times FA} 
\arrow[swap]{d}
	{\CCCTrans_{\scriptsizescriptsizeexpobj{A}{B}} \times \CCCTrans_A} & 
(\expobj{FA}{FB}) \times FA \arrow{r}[yshift=0mm]{\eval_{FA,FB}} 
\arrow[phantom, description]{d}{\overset{\CCCTransExpSmall_{A,B}}{\Leftarrow}} & 
FB \arrow{d}{\CCCTrans_B} \arrow{d}{\CCCTrans_B} \\ 
G(\expobj{A}{B}) \times GA 
\arrow[
swap,
rounded corners,
to path=
{ -- ([yshift=0ex]\tikztostart.south)
|- ([yshift=-5.5ex]\tikztostart.east)
-| ([yshift=-5.5ex]\tikztotarget.north)
-- (\tikztotarget.south)}, 
]{rr}
\arrow[swap, phantom]{rr}[yshift=-7.2ex, font=\scriptsize]
	{\eval_{GA,GB} \circ (\evBar^G_{A,B} \times GA)}
\arrow[swap]{r}[yshift=-0mm]{\evBar^G_{A,B} \times GA} & 
(\expobj{GA}{GB}) \times GA 
\arrow[swap]{r}[yshift=-0mm]{\eval_{GA,GB}} & 
GB 
\end{td} 
such that the following pasting diagram is equal to 
$\CCCTransNat_{\eval_{A,B}}$:
\begin{center}
\scalebox{1}{ 
\begin{tikzcd}[ampersand replacement = \&, column sep = 4em, row sep = 3em] 
\: \& 
\: \& 
F\big((\expobj{A}{B}) \times A\big) 
\arrow{dr}{F\eval_{A,B}} \& 
\: \\ 
F\big( (\expobj{A}{B}) \times A \big) 
\arrow[
swap,
rounded corners,
to path=
{ -- ([yshift=0ex]\tikztostart.north)
|- ([yshift=2.5cm]\tikztostart.east)
-| ([yshift=2.5cm]\tikztotarget.south)
-- (\tikztotarget.north)}, 
]{rrr}
\arrow[phantom]{rrr}[yshift=2.7cm, font = \scriptsize]{F\eval_{A,B}}
\arrow[phantom]{rrr}[yshift=2.1cm, xshift=-2mm, font=\small]{\iso}
\arrow{r}{\seq{F\pi_1, F\pi_2}} 
\arrow[swap]{d}[description]
	{\CCCTrans_{(\scriptsizeexpobj{A}{B}) \times A}} 
\arrow[bend left = 12]{urr}
	{\Id_{F((\scriptsizescriptsizeexpobj{A}{B}) \times A)}}
\arrow[phantom, description]{dr}
	{\overset{
		\CCCTransProd_{\scriptsizescriptsizeexpobj{A}{B}, A}}
		{\Leftarrow}
	} \& 
F(\expobj{A}{B}) \times FA 
\arrow{d}[description]
	{\CCCTrans_{\scriptsizescriptsizeexpobj{A}{B}} \times \CCCTrans_A} 
\arrow{rr}{\eval_{FA, FB} \circ (\evBar^F_{A,B} \times FA)} 
\arrow[phantom, description]{u}{\iso} 
\arrow{ur}[xshift=1mm]{\prodPres_{(\scriptsizescriptsizeexpobj{A}{B}, A)}} \& 
\arrow[phantom, description]{d}
	{\overset{\CCCTransExpSmall_{A,B}}{\Leftarrow}} 
\arrow[phantom, description]{u}[yshift=1mm]
	{\overset{\epsilonExp}{\iso}} \& 
FB 
\arrow{d}[description]{\CCCTrans_{B}} \\ 
G\big( (\expobj{A}{B}) \times A \big) 
\arrow[
swap,
rounded corners,
to path=
{ -- ([yshift=0ex]\tikztostart.south)
|- ([yshift=-2.5cm]\tikztostart.east)
-| ([yshift=-2.5cm]\tikztotarget.north)
-- (\tikztotarget.south)}, 
]{rrr}
\arrow[phantom]{rrr}[yshift=-2.7cm, font = \scriptsize]{G\eval_{A,B}}
\arrow[phantom]{rrr}[yshift=-2.1cm, xshift=-2mm, font=\small]{\iso}
\arrow[swap]{r}{\seq{G\pi_1, G\pi_2}} 
\arrow[bend right = 12, swap]{drr}
	{\Id_{G((\scriptsizescriptsizeexpobj{A}{B}) \times A)}} 
\& 
G(\expobj{A}{B}) \times GA 
\arrow[swap]{rr}{\eval_{GA, GB} \circ (\evBar_{A,B}^G \times GA)} 
\arrow[swap]{dr}[xshift=1mm]{\prodPres_{\scriptsizescriptsizeexpobj{A}{B},B}} 
\arrow[phantom, description]{d}{\iso} \& 
\arrow[phantom, description]{d}[yshift=-1mm]{\overset{\epsilonExp}{\iso}} \& 
GB \\ 
\: \& 
\: \& 
G\big((\expobj{A}{B}) \times A\big) 
\arrow[swap]{ur}{G\eval_{AB}} \& 
\: 
\end{tikzcd} 
}
\end{center}
We call the transformation \Def{strong} if every $\CCCTransNat_{f}$, 
$\CCCTransProd_{A_1, \,\dots\, , A_n}$ and $\CCCTransExp_{A,B}$ is invertible. 

In a 
cc-bicategory, every fp-transformation---and hence every pseudonatural 
transformation---lifts canonically to a cc-transformation: one simply inverts 
the coherence law to obtain a definition of $\CCCTransExp_{A,B}$. 
Moreover, by Lemma~\ref{lem:biequivalences-preserve-biuniversal-arrows}  
every biequivalence extends canonically to a cc-pseudofunctor. Thus, in order 
to construct a 
\Def{cc-biequivalence} between 
cc-bicategories---namely a 
biequivalence of the underlying bicategories in which the pseudofunctors are 
cc-pseudofunctors and the pseudonatural transformations are 
cc-transformations---it suffices to construct a biequivalence of the underlying 
bicategories (\cf~Lemma~\ref{lem:fp-biequivalence-from-biequivalence}).

\begin{prooflesslemma} \label{lem:bicat-to-cartclosedbicat-equivalence-lift}
Let $\ccBicat{\baseCat}$ and $\ccBicat{\altCat}$ be cc-bicategories.
Then there exists a 
biequivalence $\baseCat \simeq \altCat$ if and only if 
there exists a cc-biequivalence $\ccBicat{\baseCat} \simeq \ccBicat\altCat$. 
\end{prooflesslemma}

%

\subsection{Coherence via the Yoneda embedding.}
It turns out that one may refine the Yoneda-style proof of coherence for 
fp-bicategories given on page~\pageref{prop:power-coherence}~%
(Proposition~\ref{prop:power-coherence}) to encompass exponentials.%
\footnote{I am grateful to Andr{\'e} Joyal for suggesting this is possible, especially so because at the time I thought it was not.} 
 The 
proof does not go through verbatim, because the exponentials in 
$\Hom(\baseCat, \Cat)$ are not generally strict. The solution is to first 
strictify the bicategory $\baseCat$ to a 2-category $\altCat$, then pass to 
the 2-category $\twoHom{\altCat}{\Cat}$ of 
2-functors, 2-natural transformations, and modifications. This is cartesian 
closed as a 2-category---and hence as a bicategory---by general
enriched category theory~\cite[Example~5.2]{Day1970}. 

\begin{mypropn} \label{prop:yoneda-coherence-for-cc-bicats}
For any cc-bicategory $\ccBicat{\baseCat}$ there exists a 
strictly cartesian closed 2-category $\ccBicat{\altCat}$ such that
$\baseCat \simeq \altCat$.
\begin{proof}
By Proposition~\ref{prop:power-coherence} we may assume without loss of 
generality that $\baseCat$ is a 2-category with 2-categorical products and 
pseudo-exponentials. It therefore admits a \emph{2-categorical} Yoneda embedding
$\Yon : \baseCat \hookrightarrow \twoHom{\op\baseCat}{\Cat}$. Let 
$\overline\baseCat$ denote the closure of $\Yon(ob(\baseCat))$ under 
equivalences and factor the Yoneda embedding as 
$\baseCat \xra{i} \overline\baseCat \xra{j} \twoHom{\op\baseCat}{\Cat}$. By the 
2-categorical Yoneda lemma, $i$ is a biequivalence. 

The rest of the argument runs as for Proposition~\ref{prop:power-coherence}. 
For any $P, Q \in \overline\baseCat$ the strict exponential
$(\exp{jP}{jQ})$ exists in $\twoHom{\op\baseCat}{\Cat}$. But then
\[
(\exp{jP}{jQ}) 
= \left(\exp{(\Yon i^{-1})P}{(\Yon i^{-1})Q}\right)
\simeq \Yon{\left( \exp{i^{-1}P}{i^{-1}Q} \right)}
\]
so the exponential $(\exp{jP}{jQ}) \in \overline\baseCat$, as required.
\end{proof}
\end{mypropn}

In a sense, of course, this proposition solves the problem we set ourselves in 
the introduction to this thesis: cc-bicategories are coherent. However, the 
normalisation-by-evaluation proof is valuable in itself. First, it is a new approach to higher-categorical 
coherence; second, the speculation that it may be refinable to a normalisation 
algorithm on 2-cells; and third, it makes use of machinery that will 
play an important role in other, further developments. We therefore keep this 
result in mind, but do not let it deter us from our work in the rest of this 
thesis.

\section{Cartesian closed (bi)clones}

We shall follow the procedure of the previous two chapters, synthesising our 
type theory from the 
construction of a free biclone. The 1-categorical setting remains an 
enlightening starting point: in this setting, the type theory 
we synthesise ought to be the familiar simply-typed lambda calculus. To show 
this is indeed the case, we shall 
extend the diagram of adjunctions~(\ref{eq:cart-clone-cat-multigraph-diagram})
on page~\pageref{eq:cart-clone-cat-multigraph-diagram} to 
the cartesian closed setting. The ideas involved are not especially novel; 
however, to the best of my knowledge they have not been presented in 
this style elsewhere (although Jacobs'~\cite{Jacobs1992} shares many of the 
same basic insights).

\subsection{Cartesian closed clones}

Lambek~\cite{Lambek1989} defines a \Def{(right) internal hom} in a 
multicategory $\mCat$ to be a choice of object $\expobj{A}{B}$ for every $A, B 
\in \mCat$, together with a family of multimaps 
$\eval_{A,B} : (\expobj{A}{B}), A \to B$ inducing isomorphisms
\begin{align*}
\mCat(\Gamma; \expobj{A}{B}) &\xra{\iso} \mCat(\Gamma, A; B) \\
(h : \Gamma \to \expobj{A}{B}) &\mapsto 
	(\Gamma, A \xra{\ms{\eval_{A,B}}{h, \id_A}} B)
\end{align*}
for every $\Gamma$, $A$ and $B$. This suggests the following definition for 
clones~(\cf~Definition~\ref{def:representable-and-cartesian-clones}).

\begin{mydefn}
A  clone $(S, \clone)$ \Def{has a (right) internal hom} 
if the corresponding multicategory $\mCatOf\clone$ has a right internal hom. If 
$\clone$ is also cartesian, we say $\clone$ is \Def{cartesian closed}. 
\end{mydefn}

\begin{myexmp} \label{ex:cartesian-closed-clone-from-cc-cat}
The cartesian clone $\cloneOf\catC$ constructed from a cartesian closed category 
$\ccBicat{\catC}$ (recall Example~\ref{ex:fp-cat-to-cartesian-clone} on page~\pageref{ex:fp-cat-to-cartesian-clone}) is cartesian closed. The exponential of $A, B \in \catC$ is $\expobj{A}{B}$, the evaluation multimap is the evaluation map of $\catC$, and the currying of $f : \prodop_{n+1}(A_1, \dots, A_n, X) \to Y$ is the exponential transpose of
\[
\prodop_2{\left( \prodop_n(A_1, \dots, A_n), X \right)}
\xra{\iso} 
\prodop_{n+1}(A_1, \dots, A_n, X)
\xra{f} 
Y
\]
\end{myexmp}

Since every cartesian clone is representable, for any cartesian closed clone
$\ccClone{S}{\clone}$ one 
obtains the following chain of natural isomorphisms 
for every $A_1, \,\dots\, , A_n, B, C \in S \:\: (n \in \Nat)$:
\begin{equation} \label{eq:cartesian-closed-clone}
\begin{aligned}
\clone{\left(\prodop_{n+1}(A_1, \,\dots,\, A_n, B); C\right)} 
	&\iso \clone(A_1, \,\dots,\, A_n, B; C)
			& \text{by representability} \\ 
	&\iso \clone(A_1, \,\dots,\, A_n; \expobj{B}{C})
			& \text{by cartesian closure} \\
	&\iso \clone{\left(\prodop_n (A_1, \,\dots,\, A_n);  \expobj{B}{C}\right)}
			& \text{by representability} 
\end{aligned}
\end{equation}
Thus, for any multimap 
$t : A_1, \dots, A_n, B \to C$ in a 
cartesian closed clone $\ccClone{S}{\clone}$ 
there exists a multimap
$\lambda t : A_1, \,\dots\, , A_n \to (\exptype{B}{C})$ 
(called the \Def{currying} of $t$), which is the unique 
$g : A_1, \,\dots\, , A_n \to (\exptype{B}{C})$ satisfying
\[
t = \cslr
		{\eval_{A,B}}
		{\cs{g}{\p{1}{\ind{A};B}, \,\dots\, , \p{n}{\ind{A};B}}, 
			\p{n+1}{\ind{A};B}}
\]
Observe in particular how the requirement that the isomorphisms are defined on 
$\mCatOf{\clone}$---rather than on $\clone$---abstractly enforces the use of 
the \Def{weakening} operation taking $h : X_1, \,\dots\, , X_n \to Z$ to the 
multimap
$\cslr{h}{\p{1}{\ind{X}, Y}, \,\dots\, , \p{n}{\ind{X}, Y}} : 
	X_1, \,\dots\, , X_n, Y \to Z$.

\begin{myremark} \label{rem:exponential-structure-not-preserved}
For any cartesian closed clone $\ccClone{S}{\clone}$ the 
isomorphisms~(\ref{eq:cartesian-closed-clone}) entail that the nucleus 
$\overline{\clone}$ is also cartesian closed. 
Thus products are given as in $(S, \clone)$, and exponentials are given by the 
composite natural isomorphism
\begin{equation} \label{eq:exponentials-in-restriction}
\overline{\clone}(X \times A, B) =
\clone(X \times A, B)  \iso 
\clone(X, A; B) \iso
\clone(X, \expobj{A}{B}) =
\overline{\clone}(X, \expobj{A}{B})
\end{equation}
However, the evaluation map 
$\eval_{A,B} : (\expobj{A}{B}), A \to B$
witnessing exponentials in $\clone$ 
is not a morphism in $\overline\clone$. Chasing through the 
isomorphism~(\ref{eq:exponentials-in-restriction}), one sees that the 
evaluation map 
$(\exp{A}{B}) \times A \to B$ in 
$\overline{\clone}$ 
is 
$\cs{\eval_{A,B}}{\pi_1, \pi_2}$ and the currying of
$f : X \times A \to B$ 
is the 1-cell
$\lambda\big( 
	X, A \xra{\pair{\p{1}{X, A}, \p{2}{X,A}}} X \times A \xra{f} B
\big)$. 
%
%
%
%
To see this is the case, observe first that for any $u : X \to (\expobj{A}{B})$
one has: 
\begin{align*}
\csthree
	{\eval_{A,B}}
	{\cs{u}{\p{1}{X, A}}, \p{2}{X, A}}
	{\pi_1, \pi_2}
&= 
\cslr
	{\eval_{A,B}}
	{\csthreesmall
		{u}{\p{1}{X,A}}{\pi_1, \pi_2}, \cs{\p{2}{X, A}}{\pi_1, \pi_2}} \\
&=
\cs
	{\eval_{A,B}}
	{\cs{u}{\pi_1}, \pi_2} 
\end{align*}
Next recall that for any $u : X \to Y$ in 
$\overline\clone$ the corresponding morphism 
$u \times A : X \times A \to Y \times A$
is $\pair{\cs{u}{\pi_1}, \pi_2}$.
Putting these components together, one sees that for any 
$f : X \times A \to B$, 
\begin{align*}
{\cslr{\eval_{A,B}}{\pi_1, \pi_2}}
	&\uncs{\pairlr
		{\cslr
			{\lambda{\big( 
				\cs
					{f}
					{\pair{\p{1}{X,A}, \p{2}{X,A}}}		
				\big)}}
					{\pi_1}, 
			\pi_2}}  \\
&=
\cslr
	{\eval_{A,B}}
	{\cs
		{\lambda{{\big( 
				\cs
					{f}
					{\pair{\p{1}{X,A}, \p{2}{X,A}}}		
				\big)}}}{\pi_1},
	\pi_2} & \text{cartesian structure of $\clone$} \\
&=
\csthree
	{\eval_{A,B}}
	{
		\cs{\lambda{{\big( 
				\cs
					{f}
					{\pair{\p{1}{X,A}, \p{2}{X,A}}}		
				\big)}}}
			{\p{1}{X, A}},
		\p{2}{X, A}
	}
	{\pi_1, \pi_2} \\
&=
\cs
	{\cs{f}{\pair{\p{1}{X,A}, \p{2}{X,A}}}}
	{\pi_1, \pi_2} & \text{exponentials in $\clone$} \\
&= f
\end{align*}
The final line follows by 
Lemma~\ref{lem:product-structure-from-representable-arrow}. On the other hand, 
for any $u : X \to (\expobj{A}{B})$, 
\enlargethispage*{3\baselineskip}
\begin{align*}
\lambda{\big( 
		\cs
			{\eval_{A,B}}
			{\pi_1, \pi_2}
		\uncs
			{\pair{\cs{u}{\pi_1}, \pi_2}}
		\uncssmall
			{\pair{\p{1}{X, A}, \p{2}{X, A}}}
	\big)} 
&=	
\lambda{\left( 
		\csthree
			{\eval_{A,B}}
			{\cs{u}{\pi_1}, \pi_2}
			{\pair{\p{1}{X, A}, \p{2}{X, A}}}
	\right)} \\
&= 
\lambda{\left( 
		\cslr
			{\eval_{A,B}}
			{\cs{u}{\p{1}{X, A}}, \p{2}{X, A}}
	\right)} \\
&= 
	u
\end{align*}
where the final line follows again from the cartesian closed structure in 
$(S, \clone)$. 
It follows that $\cs{\eval_{A,B}}{\pi_1, \pi_2}$ is the universal arrow 
defining exponentials, as claimed. 

This structure is not surprising: it corresponds to the cartesian closed 
structure on the syntactic model of the simply-typed lambda calculus, 
restricted to unary contexts~(\eg~\cite[Theorem~4.8.4]{Crole1994}). 
\end{myremark}

The following two definitions follow the schema of 
Chapters~\ref{chap:biclone-lang} and~\ref{chap:fp-lang}.

\begin{mydefn} \label{def:stlc-sig}
A \Def{$\stlc$-signature} $\sig = (\baseTypes, \graph)$ consists of 
\begin{enumerate}
\item A set of base types $\baseTypes$, 
\item A multigraph $\graph$ with nodes generated by the grammar
\begin{equation} \label{eq:stlc-sig-grammar}
A_1, \,\dots\, , A_n, C, D ::= 
	B \st 
	\prodop_n(A_1, \,\dots\, , A_n)  \st
	\exptype{C}{D} \qquad (B \in \baseTypes, n \in \Nat)
\end{equation}
\end{enumerate}
If the multigraph $\graph$ is a graph we call the signature \Def{unary}. A 
\Def{homomorphism} of $\stlc$-signatures
$h : \sig \to \sig'$ 
is a morphism $h : \graph \to \graph'$ of the underlying multigraphs
such that, additionally,
\begingroup
\addtolength{\jot}{.5em}
\begin{align*}
h{\left( \prodop_n(A_1, \,\dots\, , A_n) \right)} 
&= \prodop_n{\left( hA_1, \,\dots\, , hA_n \right)} \\
h(\exptype{C}{D}) &= (\exptype{hC}{hD})
\end{align*}
\endgroup
We denote the category of $\stlc$-signatures and their homomorphisms by
$\stlcSigCat$, and the full subcategory of unary $\stlc$-signatures by
$\stlcUnSigCat$. 
\end{mydefn}

\begin{mynotation}[{\cf~Notation~\ref{not:stlc-times-alltypes-not}}] 
\label{not:stlc-alltypes-not}
For any $\stlc$-signature $\sig = (\baseTypes, \graph)$ we write 
$\allTypes\baseTypes$ for the set generated from $\baseTypes$ by the 
grammar~(\ref{eq:stlc-sig-grammar}).
In particular, when the signature is just a set (\ie~the graph $\graph$ has no edges) we denote the 
signature $\sig = (\baseTypes, \sig)$ simply by $\allTypes\baseTypes$. 
\end{mynotation}

\begin{mydefn}
A \Def{cartesian closed clone homomorphism} 
$h : \ccClone{S}{\clone} \to \ccClone{T}{\altClone}$ 
is a cartesian clone homomorphism 
$\cartClone{S}{\clone} \to \cartClone{T}{\altClone}$
such that the canonical map
$
\lambda{\left( h(\eval_{A,B}) \right)} :
h(\exptype{A}{B}) \to (\expobj{hA}{hB})$ 
is invertible. 
We call $h$ \Def{strict} if
\begin{align*}
h(\expobj{A}{B}) &= (\exptype{hA}{hB}) \\
h(\eval_{A,B}) &= \eval_{hA,hB}
\end{align*}
for every $A, B \in S$.
\end{mydefn}

In a similar fashion, we call a cartesian closed functor \Def{strict} if it 
strictly preserves exponentials and the evaluation map.

\newpage
We now construct the following diagram of 
adjunctions, in which $\CCCatCat$ denotes the category of cartesian closed 
categories 
and strict cartesian closed functors and $\CCClone$ denotes the category of 
cartesian closed clones and strict homomorphisms. As in the preceding chapter, 
we implicitly restrict to cartesian structure in which $\prodop_1(-)$ 
is the identity functor. 
\begin{equation} 
\label{eq:cc-clone-cat-multigraph-diagram}
\begin{tikzcd}[column sep = 5.5em, row sep = 3em]
\: &
\CCClone 
\arrow[bend left = 18]{dr}{\overline{(-)}}
\arrow[bend right = 18]{dl}[swap]{\text{forget}} &
\: \\
\stlcSigCat
\arrow[phantom]{ur}[description]{\adjDown}
\arrow[bend right = 18]{ur}[swap]{\freeCartClosedClone{-}}
\arrow[bend right = 18]{dr}[swap]{\widetilde{\lin}} &
\: &
\CCCatCat
\arrow[phantom]{ul}[description]{\adjDown}
\arrow[bend left = 18]{ul}[near end]{\prom}
\arrow[bend left = 18]{dl}{\text{forget}} \\
\: &
\stlcUnSigCat 
\arrow[phantom]{ur}[description]{\adjUp{}}
\arrow[bend left = 18]{ur}{\text{free}}
\arrow[phantom]{ul}[description]{\adjUp{}}
\arrow[bend right = 18, hookrightarrow]{ul} &
\:
\end{tikzcd} 
\end{equation}

The right adjoint to the inclusion 
$\inc : \stlcUnSigCat \hookrightarrow \stlcSigCat$
is defined by
$\widetilde{\lin}(\baseTypes, \graph) = (\baseTypes, \lin\graph)$ for 
$\lin : \MultiGraph \to \Graph$ the right adjoint to the inclusion
$\Graph \hookrightarrow \MultiGraph$ 
(\cf~Lemma~\ref{lem:stlcTimesSig-reflective-subcat}). The free-forgetful 
adjunction between 
cartesian closed categories and $\stlc$-signatures is the classical construction 
of the syntactic model of the simply-typed lambda calculus over a 
signature~\cite{Lambek1980}. There are two adjunctions left to construct.

\begin{mylemma} \quad \label{lem:free-cc-clone-on-sig}
The forgetful functor $\CCClone \to \stlcSigCat$ 
has a left adjoint. 
\begin{proof}
Define a clone 
$\freeCartClosedClone{\sig}$ over a signature $(\baseTypes, \graph)$ as follows. The 
sorts are
generated by the grammar
\[
A_1, \,\dots\, , A_n, C, D ::= 
	B \st 
	\prodop_n(A_1, \,\dots\, , A_n)  \st
	\exptype{C}{D} \qquad (B \in \baseTypes, n \in \Nat)
\] 
The operations are those of Construction~\ref{constr:free-cartesian-clone}
(page~\pageref{constr:free-cartesian-clone}) 
together with two additional rules:
\begin{center}
\unaryRule
	{\phantom{\freeCartClosedClone{\sig}}}
	{\evalterm_{B,C} \in \freeCartClosedClone{\sig}(\exptype{B}{C}, B; C)}
	{}
\quad
\unaryRule
	{t \in \freeCartClosedClone{\sig}(A_1, \,\dots\, , A_n, B; C)}
	{\lambda t \in \freeCartClosedClone{\sig}(A_1, \,\dots\, , A_n; \exptype{B}{C})}
	{$(n \in \Nat)$}
\vspace{-\treeskip}
\end{center}
Similarly, one extends the equational theory $\equiv$ by requiring that
\begin{itemize}
\item $\cslr
			{\evalterm_{B,C}}
			{\cs
				{(\lambda t)}
				{\p{1}{\ind{A}, B}, \,\dots\, , \p{n}{\ind{A}, B}},
			\p{n+1}{\ind{A}, B}
			} 
		\equiv  t$ 
for any $t : A_1, \,\dots\, , A_n, B \to C$,
\item $\lambda{\left(
		\cslr
			{\evalterm_{B,C}}
			{\cs
				{u}
				{\p{1}{\ind{A}, B}, \,\dots\, , \p{n}{\ind{A}, B}},
			\p{n+1}{\ind{A}, B}
			}\right)} 
		\equiv u$ 
for any $u : A_1, \,\dots\, , A_n \to (\exptype{B}{C})$. 
\end{itemize}
It is clear $\freeCartClosedClone{\sig}$ is cartesian closed. To see that it is also 
free, let 
$h : \sig \to \altClone$ be any $\stlc$-signature homomorphism from $\sig$ to 
the underlying $\stlc$-signature of a cartesian closed clone
$\ccClone{T}{\altClone}$. Define a cartesian closed clone homomorphism 
$\ext{h} : \freeCartClosedClone{\sig} \to \altClone$
by extending the definition of Lemma~\ref{lem:free-cartesian-clone}
(page~\pageref{lem:free-cartesian-clone}) as follows:
\begin{align*}
\ext{h}(\exptype{A}{B}) &:= (\exptype{\ext{h}A}{\ext{h}B}) \\
\ext{h}(\evalterm_{A,B}) &:= \eval_{(\ext{h}A, \ext{h}B)} \\
\ext{h}(\lambda t) &:= \lambda(\ext{h}t) 
\end{align*}
For uniqueness, we already know from Lemma~\ref{lem:free-cartesian-clone} and the definition of a cartesian closed clone homomorphism that any cartesian clone homomorphism strictly preserves all the structure, except for currying. So it suffices to show that any cartesian clone homomorphism preserves the $\lambda(-)$ mapping. 
%
Since $\lambda t $ is the unique multimap
$g: A_1, \,\dots\, , A_n \to (\exptype{B}{C})$ such that
$t = \cslr
			{\eval_{B,C}}
			{\cs
				{g}
				{\p{1}{\ind{A}, B}, \,\dots\, , \p{n}{\ind{A}, B}},
			\p{n+1}{\ind{A}, B}
			}$,
for any cartesian clone homomorphism 
$f : \freeCartClosedClone{\sig} \to \altClone$
one has
\begin{align*}
f(t) &= 
		f{\left( \cslr{\evalterm_{B,C}}
					{\cs
						{\big(\lambda t \big)}
						{\p{1}{\ind{A}, B}, \,\dots\, , \p{n}{\ind{A}, B}},
					\p{n+1}{\ind{A}, B}
					} \right)} \\
	&= \cslr{\eval_{fB,fC}}
			{\cslr
				{f(\lambda t)}
				{\p{1}{f\ind{A}, fB}, \,\dots\, , \p{n}{f\ind{A}, fB}},
			\p{n+1}{f\ind{A}, fB}		
			} 
\end{align*}
it follows that $f(\lambda t) = \lambda f(t)$ for every 
$t : A_1, \,\dots\, , A_n, B \to (\exptype{B}{C})$, as required.
\end{proof}
\end{mylemma}

It remains to construct the adjunction 
$\CCClone \leftrightarrows \CCCatCat$. 

\begin{mylemma} \label{lem:free-cc-clone-on-cc-cat}
The functor 
$\overline{(-)} : \CCClone \to \CCCatCat$ 
restricting a cartesian closed clone to its nucleus has a left adjoint.
\begin{proof}
Consider the functor 
$\prom : \CartCatCat \to \CartClone$ 
defined in
Lemma~\ref{lem:free-cart-clone-on-cart-cat}. This restricts to a functor
$\CCCatCat \to \CCClone$. Explicitly, the evaluation map in $\prom\catC$ is the 
evaluation map $\eval_{A,B}$ in 
$\catC$ and for any $f : X_1, \,\dots\, , X_n \to (\exp{A}{B})$ the composite
$\cslr{\eval_{A,B}}
	{
	\cs
		{f}
		{\p{1}{\ind{X}, A}, \,\dots\, , \p{n}{\ind{X}, A}}, 
	\p{n+1}{\ind{X}, A}
	}$
in $\prom\catC$ is the composite
$\eval_{A,B} \circ \seqlr{ f \circ \seq{\pi_1, \,\dots\, , \pi_n}, \pi_{n+1} }
	= \eval_{A,B} \circ (f \times A) \circ 
		\seqlr{\seq{\pi_1, \,\dots\, , \pi_n}, \pi_{n+1}}$
in $\catC$. The currying of $g : X_1, \,\dots\, , X_n, A \to B$ is the currying 
(in 
$\catC$) of the 
morphism
\[
\lambda{\big( 
	\prodop_{i=1}^n X_i \times A \xra{\iso}
	X_1 \times \cdots \times X_n \times A \xra{g}
	B
\big)}
\]
Now suppose that $F : \catC \to \overline{\altClone}$ is a strict cartesian 
closed functor. Define $\ext{F}$ as the free cartesian extension of $F$ from 
Lemma~\ref{lem:free-cart-clone-on-cart-cat}:
\[\ext{F}(X_1, \,\dots\, , X_n \xra{t} Y) :=
	\big( FX_1, \,\dots\, , FX_n 
			\xra{\psi_{F\ind{X}}(\p{1}{}, \,\dots\, , \p{n}{})}
		\prodop_{i=1}^n FX_i = F(\prodop_{i=1}^n X_i) 
			\xra{Ft}
		FY \big)
\]
To see that $\ext{F}$ preserves the evaluation map, note that---since $F$ is a 
strict cartesian closed functor---the equation
$F(\eval_{A,B}) = \cs{\eval_{FA,FB}}{\pi_1, \pi_2}$ must hold by 
Remark~\ref{rem:exponential-structure-not-preserved}. It follows that
\begin{align*}
\ext{F}(\eval_{A,B}) &= 
	\csthree
		{\eval_{FA,FB}}
		{\pi_1, \pi_2}
		{\psi_{F\ind{X}}(\p{1}{}, \,\dots\, , \p{n}{})} \\
	&= \cslr
		{\eval_{FA,FB}}
		{\p{1}{\scriptsizeexpobj{FA}{FB}, FA}, \p{2}{\scriptsizeexpobj{FA}{FB}, FA}}
		&\text{by equation~(\ref{eq:equations-from-universal-arrow}) 
				on page~\pageref{eq:equations-from-universal-arrow}}
		\\
	&= \eval_{FA,FB}
\end{align*}
as required. The proof of uniqueness is exactly as in the cartesian case.
\end{proof}
\end{mylemma}

This completes the construction of the diagram of 
adjunctions~(\ref{eq:cc-clone-cat-multigraph-diagram}). As for the diagram of 
adjunctions~(\ref{eq:cart-clone-cat-multigraph-diagram}) for cartesian strucure, 
it is easy to see that 
the outer edges of~(\ref{eq:cc-clone-cat-multigraph-diagram}) commute and
that $\nucleus{(-)} \circ \prom = \id_{\CCCatCat}$. One thereby obtains the 
following chain of natural isomorphisms~(\cf~equation~(\ref{eq:free-cartcat-suffices-to-restrict})),
in which we write $\mathbb{FC}\mathrm{at}^{{\times},{\to}}(\sig)$ for the free cartesian
closed category on a unary signature $\sig$:
\begin{equation} \label{eq:free-cartclosedcat-suffices-to-restrict}
\CCCatCat(\mathbb{FC}\mathrm{at}^{{\times}, {\to}}(\sig), \catC)
	= 
\CCCatCat{\left(\nucleus{\prom(\mathbb{FC}\mathrm{at}^{{\times},{\to}}(\sig))}, \catC\right)}
	\iso 
\CCCatCat{\left(\nucleus{\freeCartClosedClone{\inc\sig})}, \catC\right)}
\end{equation}
It follows 
that the free cartesian closed category on a $\stlc$-signature is described by 
restricting the deductive system of Lemma~\ref{lem:free-cc-clone-on-sig} to 
unary contexts. 

\begin{myremark} \label{rem:exponentials-preserved-proof-fails}
In the preceding lemma we rely on the equation
$ \cs{\eval_{FA,FB}}
	{\p{1}{(\scriptsizeexpobj{A}{B}, A)}, 
	\p{2}{(\scriptsizeexpobj{A}{B}, A)}} = 
	\eval_{FA,FB}$
to show that $\ext{F}$ is strictly cartesian closed. In the bicategorical 
setting, where this equality is generally only an isomorphism, the 
argument fails. 
As we shall see, the free cc-bicategory on a signature (in the strict sense of 
\emph{free} we have been using throughout) 
is not obtained by restricting the free cartesian biclone on the same signature.
\end{myremark}

\paragraph*{Cartesian closed clones and the simply-typed lambda calculus.} Let 
us examine how one extracts the simply-typed lambda calculus from the 
internal language of $\freeCartClosedClone{\sig}$ (defined in 
Lemma~\ref{lem:free-cc-clone-on-cc-cat}). The $\evalterm_{B,C}$ multimap 
becomes an application operation on variables:
\begin{center}
\unaryRule
	{\faketext}
	{f : \exptype{B}{C}, x : B \vdash \app{f}{x} : C}
	{}
\vspace{-\treeskip}
\end{center}
The weakening operation
$t \mapsto \cslr{t}{\p{1}{\ind{A}, B}, \,\dots\, , \p{n}{\ind{A}, B}}$
is the following form of the usual substitution lemma:
\begin{center}
\binaryRule
	{x_1 : A_1, \,\dots\, , x_n : A_n \vdash t : C}
	{x_1 : A_1, \,\dots\, , x_n : A_n, y : B \vdash t : C}
	{x_1 : A_1, \,\dots\, , x_n : A_n, y : B \vdash t[x_1/x_1, \,\dots\, , 
	x_n/x_n] : C}
	{}
\vspace{-\treeskip}
\end{center}
This mirrors the construction in $\langBiclone$ and its extensions, where 
weakening arises from explicit substitutions corresponding to inclusions of 
contexts. 

The $\lambda(-)$ mapping is the usual lambda abstraction operation, and the two 
equations become the following rules for every
$\context, x : A \vdash t :B$ and $\context \vdash u : \exptype{A}{B}$:
\begin{center}
$\app{(\lam{x}{t})[x_1/ x_1, \,\dots\, , x_n/x_n]}{x}$ \quad and \quad 
$\lam{x}{\app{u[x_1/ x_1, \,\dots\, , x_n/x_n]}{x}} = u$
\end{center}
As we saw in Section~\ref{sec:products-in-stlc}, these rules extend to 
rules on all terms in the presence of the meta-operation of capture 
avoiding substitution.
Thus, we recover the usual $\beta\eta$-laws of the simply-typed lambda calculus. The diagram of 
adjunctions~(\ref{eq:cc-clone-cat-multigraph-diagram}), together with 
the isomorphism~(\ref{eq:free-cartclosedcat-suffices-to-restrict}), 
then expresses the usual free property of the unary-context syntactic model~\cite[Chapter~4]{Crole1994}.

Our aim in what follows is to define cartesian closed biclones, 
construct the free instance to obtain a diagram 
matching~(\ref{eq:cc-clone-cat-multigraph-diagram}), and use this to extract a 
type theory in the same way as we have just sketched for the simply-typed 
lambda calculus. As for products, our insistence on strict universal properties makes 
the full diagram impossible to replicate (recall Example~\ref{exmp:nucleus-of-fp-not-free} on 
page~\pageref{exmp:nucleus-of-fp-not-free}). Nonetheless, we shall
see that a version of it exists up to biequivalence.

\subsection{Cartesian closed biclones}

The definitions of the previous section bicategorify in the way one would 
expect.

\begin{mydefn} \quad
\begin{enumerate}
\item A \Def{(right) closed} bi-multicategory is a bi-multicategory $\mbCat$ 
equipped with the following data for every $A, B \in \mbCat$:
\begin{enumerate}
\item A chosen object $\expobj{A}{B}$,
\item A chosen multimap $\eval_{A,B}: (\expobj{A}{B}), A \to B$, 
\item For every sequence of objects $\Gamma$ in $\mbCat$, an adjoint 
equivalence
\begin{equation*} 
\begin{tikzcd} 
\mbCat(\Gamma; \expobj{A}{B}) 
\arrow[bend left = 20]{r}
	{\ms{\eval_{A,B}}{(-), \Id_A}} 
\arrow[phantom]{r}[xshift=0em]{\adjUp{\simeq}} & 
\mbCat(\Gamma, A; B) 
\arrow[bend left = 20]{l}{\lambda} 
\end{tikzcd} 
\end{equation*}
specified by choosing a universal arrow with components
$\epsilon_t : \ms{\eval_{A,B}}{\lambda t, \Id_A} \To t$.
\end{enumerate} 
\item A \Def{(right) closed biclone} is a biclone $(S, \biclone)$ equipped with 
a choice of right-closed structure on the corresponding 
bi-multicategory $\mCatOf\biclone$.
\item A \Def{cartesian closed} biclone is a biclone equipped 
with a choice of both cartesian structure and right-closed structure. 
\qedhere  
\end{enumerate}
\end{mydefn}

Explicitly, a cartesian closed biclone is defined by the following 
universal property. For every sequence of objects 
$\Gamma := (A_1, \,\dots\, , A_n)$ and multimap
$t : \Gamma, A \to B$
there exists a multimap
$\lambda t : \Gamma \to (\exp{A}{B})$
and a 2-cell 
$\epsilon_t : \csbig
				{\eval_{A,B}}
				{\cs
					{(\lambda t)}
					{\p{1}{\ind{A}, B}, \,\dots\, , \p{n}{\ind{A}, B}},
				\p{n+1}{\ind{A}, B}
				} 
				\To 
				t$. 
This 2-cell is
universal in the sense that for every 
$u : \Gamma \to (\exp{A}{B})$
and 
$\alpha : \csbig
				{\eval_{A,B}}
				{\cs
					{u}
					{\p{1}{\ind{A}, B}, \,\dots\, , \p{n}{\ind{A}, B}},
				\p{n+1}{\ind{A}, B}
				}
		\To t$
there exists a 2-cell
$\transExp{\alpha} : u \To \lambda t$, 
unique such that
\small
\begin{equation} \label{eq:cc-biclone-ump}
\begin{tikzcd}[column sep = 1.5em]
%
\cslr
	{\eval_{A,B}}
	{\cs
		{u}
		{\p{1}{\ind{A}, B}, \,\dots\, , \p{n}{\ind{A}, B}},
		\p{n+1}{\ind{A}, B}
	} 
\arrow[swap]{dr}{\alpha}
\arrow{rr}[yshift=3mm]
	{\cslr{\eval_{A,B}}
	{\cs
		{\transExp{\alpha}}
		{\p{1}{\ind{A}, B}, \,\dots\, , \p{n}{\ind{A}, B}},
		\p{n+1}{\ind{A}, B}
	} } &
\: &
\cslr
	{\eval_{A,B}}
	{\cs
		{(\lambda t)}
		{\p{1}{\ind{A}, B}, \,\dots\, , \p{n}{\ind{A}, B}},
	\p{n+1}{\ind{A}, B}
	}
\arrow{dl}{\epsilon_t} \\
\: &
t &
\:
\end{tikzcd}
\end{equation}
\normalsize

Moreover, since every cartesian biclone is representable 
(Theorem~\ref{thm:representable-iff-cartesian-biclone}), one also obtains a 
sequence of 
pseudonatural adjoint equivalences lifting~(\ref{eq:cartesian-closed-clone}) to 
biclones:
\begin{equation} \label{eq:biclone-representable-iff-products} 
\begin{aligned}
\biclone\big(\prodop_{n+1}(A_1, \,\dots,\, A_n, B); C\big) 
	&\simeq \biclone(A_1, \,\dots,\, A_n, B; C) \\
	&\simeq \biclone(A_1, \,\dots,\, A_n; \expobj{B}{C}) \\
	&\simeq \biclone{\big(\prodop_n(A_1, \,\dots,\, A_n); \expobj{B}{C}\big)}
\end{aligned}
\end{equation}
It follows that, if
$(S, \biclone)$ is cartesian closed, then so is its nucleus
$\overline\biclone$. 

\begin{myremark} \label{rem:bicat-exp-structure-not-preserved}
We saw in Remark~\ref{rem:exponential-structure-not-preserved} that the
evaluation map witnessing cartesian closed structure in the  
nucleus $\overline\clone$ of a cartesian closed clone $\ccClone{S}{\clone}$ is not the evaluation multimap in 
$\clone$. Similarly, chasing through the 
equivalences~(\ref{eq:biclone-representable-iff-products}) one sees that the 
biuniversal arrow witnessing exponentials in the nucleus $\overline\biclone$
of a cartesian closed biclone $\ccClone{S}{\biclone}$ is
$\cs{\evalterm_{A,B}}{\pi_1, \pi_2} : 
	A \times (\expobj{A}{B}) \to B$ and the currying of
$f : X \times A \to B$ is 
$\lambda{\left( \cs{f}{\pair{\p{1}{X,A}, \p{2}{X,A}}} \right)}$. 
To see this defines 
an exponential, one can replace each of the equalities in the proof of 
Remark~\ref{rem:exponential-structure-not-preserved} to construct natural 
isomorphisms 
\begin{align*}
\csthree
	{\eval_{A,B}}
	{\cs{(-)}{\p{1}{X, A}}, \p{2}{X, A}}
	{\pi_1, \pi_2}
&\iso 
\id_{\biclone(X \times A, B)} \\
\lambda{\big( 
		\cs
			{\eval_{A,B}}
			{\pi_1, \pi_2}
		\uncs
			{\pair{\cs{(-)}{\pi_1}, \pi_2}}
		\uncssmall
			{\pair{\p{1}{X, A}, \p{2}{X, A}}}
	\big)} 
&\iso 
\id_{\biclone(X, \expobj{A}{B})}
\end{align*}
witnessing an equivalence, which may be promoted to the required adjoint 
equivalence without changing the functors (see~\eg~\cite[\S~IV.4]{cfwm}). 
\end{myremark}

\begin{myexmp}[{\cf~Example~\ref{ex:cartesian-closed-clone-from-cc-cat}}] 
\label{ex:cc-bicategory-to-cc-biclone}
The cartesian biclone $\bicloneFromProducts{\baseCat}$ constructed from a \mbox{cc-bicategory}
$\ccBicat{\baseCat}$~(recall Example~\ref{ex:fp-bicategory-defines-a-cartesian-biclone} on page~\pageref{ex:fp-bicategory-defines-a-cartesian-biclone}) is cartesian closed. The precise statement requires some juggling of products, for which we introduce the following notation. For any 
$A_1, \,\dots\, , A_n, B \in \baseCat \:\: (n \in \Nat)$ there 
exists a canonical equivalence
\begin{equation} \label{eq:hsem-canonical-equivalence}
e_{\ind{A}, B}
:
\prodop_{n+1}(A_1, \,\dots\, , A_n, B) 
	\leftrightarrows
\prodop_{2} \left( \prod_{n}(A_1, \,\dots\, , A_n), B\right) 
:
\psinv{e}_{\ind{A}, B} 
\end{equation}
where 
$
e_{\ind{A}, B} := 
\seq{\seq{\pi_1, \,\dots\, , \pi_n}, \pi_{n+1}}
$
and 
$
\psinv{e}_{\ind{A}, B} := 
\seq{\pi_1 \circ \pi_1, \,\dots\, , \pi_n \circ \pi_1, \pi_2}
$.
The witnessing 2-cells
\begin{equation} \label{eq:semantic-interpretation-witnessing-2-cells}
\begin{aligned}
\co_{\ind{A}, B} &: 
\psinv{e}_{\ind{A}, B} \circ e_{\ind{A}, B} 
\To 
\Id_{\prod_{n+1}(A_1, \,\dots\, , A_n, B)} \\
\un_{\ind{A}, B} &: \Id_{\prod_n(A_1, \,\dots\, , A_n) \times B} \To 
{e_{\ind{A}, B}  \circ \psinv{e}_{\ind{A}, B}} 
\end{aligned}
\end{equation}
are defined by the two diagrams below:
\begin{td}[column sep = 0em]
\seq{\pi_1 \circ \pi_1, \,\dots\, , \pi_n \circ \pi_1, \pi_2} \circ 
	\seqlr{\seq{\pi_1, \,\dots\, , \pi_n}, \pi_{n+1}} 
\arrow{r}{\co_{\ind{A}, B}}
\arrow[swap]{d}{\postName} &
\Id_{\prod_{n+1}(A_1, \,\dots\, , A_n, B)} \\
\seqlr{\left(\pi_1 \circ \pi_1\right) \circ e_{\ind{A}, B}, \,\dots\, , 
	\left(\pi_n \circ \pi_1\right) \circ e_{\ind{A}, B}, 
	\pi_2 \circ e_{\ind{A}, B}} 
\arrow[swap]{d}{\iso}  &
\seq{\pi_1, \,\dots\, , \pi_n, \pi_{n+1}} 
\arrow[swap]{u}{\widehat{\etaTimes{\Id}}^{-1}} \\
\seqlr{\pi_1 \circ \left(\pi_1 \circ e_{\ind{A}, B}\right), \,\dots\, , 
	\pi_n \circ \left(\pi_1 \circ e_{\ind{A}, B}\right), 
	\pi_2 \circ e_{\ind{A}, B}} 
\arrow[swap]{r}[yshift=-2mm]
	{\seq{\pi_1 \circ \epsilonTimesInd{1}{}, \,\dots\, 
			\pi_n \circ \epsilonTimesInd{1}{}, 
			\epsilonTimesInd{2}{}}}	&
\seqlr{\pi_1 \circ \seq{\ind{\pi}}, \,\dots\, , \pi_n \circ \seq{\ind{\pi}}, 
		\pi_{n+1}} 
\arrow[swap]{u}
	{\seq{\epsilonTimesInd{1}{}, \,\dots\, , \epsilonTimesInd{n}{}, \pi_{n+1}}} 
	\\
\: &
\: &
\: \\
\Id_{\prod_n(A_1, \,\dots\, , A_n) \times B} 
\arrow{r}{\un_{\ind{A},B}}
\arrow[swap]{d}{\widehat{\etaTimes{}}_{\Id}} &
\seqlr{\seq{\pi_1, \,\dots\, , \pi_n}, \pi_{n+1}}  \circ \psinv{e}_{\ind{A}, B} 
\\
\seq{\pi_1, \pi_2} 
\arrow[swap]{d}{\iso} &
\seqlr{\seq{\pi_1, \,\dots\, , \pi_n} \circ \psinv{e}_{\ind{A}, B}, 
		\pi_{n+1} \circ \psinv{e}_{\ind{A}, B}} 
\arrow[swap]{u}{\postName^{-1}} \\
\seq{\Id_{\prod_n(A_1, \,\dots\, , A_n)} \circ \pi_1, \pi_2}
\arrow[swap]{d}{\seq{\widehat{\etaTimes{}}_{\Id} \circ \pi_1, \pi_2}} &
\seq{\seq{\pi_1 \circ \psinv{e}_{\ind{A}, B}, \,\dots\, , 
			\pi_n \circ \psinv{e}_{\ind{A}, B}}, 
			\pi_{n+1} \circ \psinv{e}_{\ind{A}, B}} 
\arrow[swap]{u}{\seq{\postName^{-1}, \pi_{n+1} \circ \psinv{e}}} \\
\seqlr{\seq{\pi_1, \,\dots\, , \pi_n} \circ \pi_1, \pi_2} 
\arrow[swap]{r}{\seq{\postName, \pi_2}} &
\seqlr{\seq{\ind{\pi} \circ \pi_1}, \pi_2}
\arrow[swap]{u}
	{\seqlr{\seq{\epsilonTimesInd{-1}{}, \,\dots\, , \epsilonTimesInd{-n}{}},
			\epsilonTimesInd{-(n+1)}{}}}
\end{td}
Here $\widehat{\etaTimes{}}_{\Id_X}$ abbreviates the following composite:
\begin{equation} \label{eq:def-of-widehat-eta}
\widehat{\etaTimes{}}_{\Id_X} :=
\Id_X 
\XRA{\etaTimes{\Id_X}} 
\seq{\pi_1 \circ \Id_X, \,\dots\, , \pi_n \circ \Id_X}
\XRA{\iso}
\seq{\pi_1, \,\dots\, , \pi_n}
\end{equation}

The exponential 
of $A, B \in \baseCat$ is $\exptype{A}{B}$, the evaluation multimap is the evaluation map of $\baseCat$, and 
the currying of $f : \prodop_{n+1}(A_1, \dots, A_n, X) \to Y$ is  the exponential transpose of
\[
\prodop_2{\left( \prodop_n(A_1, \dots, A_n), X \right)} 
	\xrightarrow[\simeq]{\psinv{e}_{\ind{A}, X}} 
		\prodop_{n+1}(A_1, \dots, A_n, X) \xra{f} Y 
\]
The counit $\epsilonExp_f$ is the following composite:
\begin{equation*}
\makebox[\textwidth]{
\begin{tikzcd}[column sep = 5em, ampersand replacement = \&]
\eval_{X, Y} \circ 
			\seqlr{\lambda {(f \circ 
				\psinv{e}_{\ind{A},X})} 
			\circ \seq{\pi_1, \,\dots\, , \pi_n}, \pi_{n+1}} 
\arrow{r}{\epsilonExp_f}
\arrow[swap]{d}{\iso} \&
f  \\
\eval_{X, Y} \circ 
			\seqlr{\lambda (f \circ 
			\psinv{e}_{\ind{A},X})
			\circ \seq{\pi_1, \,\dots\, , \pi_n}, \Id_{X} \circ 
			\pi_{n+1}} 
\arrow[swap]{d}{\eval \circ \fuse^{-1}} \&
f \circ \Id_{\prod(\ind{A}) \times X} 
\arrow[swap]{u}{\iso} \\
\eval_{X, Y} \circ 
		\left(	\big(\lambda {(f \circ 
					\psinv{e}_{\ind{A},X}) 
				\times X}\big) 
			\circ e_{\ind{A},X}\right) 
\arrow[swap]{d}{\iso} \&
f 
	\circ (\psinv{e}_{\ind{A},X}
	\circ e_{\ind{A},X}) 
\arrow[swap]{u}{f \circ \co_{\ind{A},X}} \\
\left(\eval_{X, Y} \circ 
			\big(\lambda {(f \circ 
					\psinv{e}_{\ind{A},X}) 
				\times X}\big)\right) 
			\circ e_{\ind{A},X}
\arrow{r}[swap]
	{\epsilonExp_{(f \circ \scriptsizepsinv{e})} \circ e_{\ind{A}, X}} \&
(f 
	\circ \psinv{e}_{\ind{A},X}) 
	\circ e_{\ind{A},X}
\arrow[swap]{u}{\iso} 
\end{tikzcd}
}
\end{equation*}

For any 1-cell $g : \prodop_n(A_1, \dots, A_n) \to (\expobj{X}{Y})$ and 2-cell 
$\alpha : {\eval_{X,Y} \circ \seq{g \circ \seq{\pi_1, \dots, \pi_n}, \pi_{n+1}} \To f}$
the corresponding  mediating 2-cell $g \To \lambda(f \circ \psinv{e}_{\ind{A}, X})$ 
is $\transExp{\alpha^{{\circ}}}$, for $\alpha^{{\circ}}$ defined by the diagram below.
\begin{equation*}
\makebox[\textwidth]{
\begin{tikzcd}[column sep = 2em, ampersand replacement = \&]
\eval_{X, Y} \circ 
	(g \times X) 
\arrow{r}{{\alpha}^{\circ}}
\arrow[swap]{d}{\iso} \&
f \circ \psinv{e}_{\ind{A},X} \\
\left(\eval_{X, Y} 
	\circ (g \times X)\right) 
	\circ \Id_{\prod_2((\prod_n \ind{A}), B)}
\arrow[swap]{d}
	{\eval 
		\circ (g \times X) 
		\circ \un_{\prod_2((\prod_n \ind{A}), B)}} \&
\: \\
\left(\eval_{X, Y} 
	\circ (g \times X)\right)
	\circ \left(e_{\ind{A},X} 
	\circ \psinv{e}_{\ind{A},X}\right)
\arrow[swap]{d}{\iso} \&
\: \\
\left(\eval_{X, Y} 
	\circ \left((g \times X)\right)
	\circ e_{\ind{A},X} \right)
	\circ \psinv{e}_{\ind{A},X}
\arrow[swap]{d}
	{\eval \circ 
		\fuse \circ \psinv{e}} \&
\: \\
\left(\eval_{X, Y} 
	\circ \seqlr{g \circ \seq{\pi_1, \dots, \pi_n}, 
			\Id_{X} \circ \pi_{n+1}}\right) 
	\circ \psinv{e}_{\ind{A},X}
\arrow[swap]{r}{\iso} \&
\left(\eval_{X, Y} 
	\circ \seqlr{g \circ \seq{\ind{\pi}}, \pi_{n+1}}\right) 
	\circ \psinv{e}_{\ind{A},X} 
\arrow[swap]{uuuu}{\alpha \circ \psinv{e}} 
\end{tikzcd}
}
\end{equation*}
\end{myexmp}

\paragraph*{The free cartesian closed biclone.} In Chapters~\ref{chap:biclone-lang}~and~\ref{chap:fp-lang} we synthesised 
the required type theory from two 
principles: first, an appropriate notion of biclone, and 
second, the fact that the internal language of those biclones---when each rule is restricted to unary contexts---gives
rise to an internal language for the corresponding bicategories. For the cartesian closed case, 
we cannot restrict every rule of the internal language to unary contexts without also discarding all curried morphisms (lambda abstractions). Nonetheless we can show that the nucleus of the free cartesian closed biclone is
the free cartesian closed bicategory \emph{up to biequivalence}. Thus, one obtains the internal language 
of cartesian closed bicategories (in a bicategorical sense) by synthesising the internal language of cartesian closed biclones. 

We shall begin by defining an appropriate notion of signature and (strict) 
pseudofunctors of cartesian closed biclones. Then we 
shall construct the adjunctions of the following diagram, in which we write
$\CCBicloneCat$ for the category of cartesian closed biclones and strict pseudofunctors 
and $\CCBicatCat$ for the category of cc-bicategories and strict pseudofunctors.
\begin{equation} 
\label{eq:cc-biclone-cat-multigraph-diagram}
\begin{tikzcd}[column sep = 5em, row sep = 2.5em]
\: &
\CCBicloneCat 
\arrow[bend right = 18]{dl}[font=\scriptsize, swap]{\text{forget}} &
\: \\
\langSigCat
\arrow[phantom]{ur}[description]{\adjDown}
\arrow[bend right = 18]{ur}[swap, near end]{\freeCartClosedBiclone{-}}
\arrow[bend right = 18]{dr}[swap]{\widetilde{\lin}} &
\: &
\CCBicatCat
\arrow[bend left = 18]{dl}[font=\scriptsize]{\text{forget}} \\
\: &
\langUnSigCat 
\arrow[phantom]{ur}[description]{\adjUp}
\arrow[bend left = 18]{ur}{\freeCartClosedBicat{-}}
\arrow[phantom]{ul}[description]{\adjUp}
\arrow[bend right = 18, hookrightarrow]{ul} &
\:
\end{tikzcd} 
\end{equation}
Thereafter we shall extract our type 
theory $\langCartClosed$ from the free cartesian closed biclone over a 
signature, and use this to show that the nucleus of the free cartesian closed biclone 
is biequivalent to the free cc-bicategory 
over the same (unary) signature. 


\begin{mydefn} \label{def:lang-cart-closed-sig}
A \Def{$\langCartClosed$-signature} $\sig = (\baseTypes, \graph)$ consists of 
\begin{enumerate}
\item A set of base types $\baseTypes$, 
\item A 2-multigraph $\graph$, with nodes 
generated by the grammar
\begin{equation} \label{eq:langcartclosed-sig-grammar}
A_1, \,\dots\, , A_n, C, D ::= 
	B \st 
	\prodop_n(A_1, \,\dots\, , A_n)  \st
	\exptype{C}{D} \qquad 
	(B \in \baseTypes, n \in \Nat)
\end{equation}
\end{enumerate}
If $\graph$ is a 2-graph we call the signature \Def{unary}. A 
\Def{homomorphism} of $\langCartClosed$-signatures
$h : \sig \to \sig'$ 
is a morphism $h : \graph \to \graph'$ of the underlying multigraphs such that
\begin{center}
$h{\left( \prodop_n(A_1, \,\dots\, , A_n) \right)} 
= \prodop_n{\left( hA_1, \,\dots\, , hA_n \right)}$
\:\: and \:\:
$h(\exptype{C}{D}) = (\exptype{hC}{hD})$
\end{center}
for all $A_1, \dots, A_n, C, D \in \nodes\graph \:\: (n \in \Nat)$.
\hide{
\begingroup
\addtolength{\jot}{.3em}
\begin{align*}
h{\left( \prodop_n(A_1, \,\dots\, , A_n) \right)} 
&= \prodop_n{\left( hA_1, \,\dots\, , hA_n \right)} \\
h(\exptype{C}{D}) &= (\exptype{hC}{hD})
\end{align*}
\endgroup
}
We denote the category of $\langCartClosed$-signatures and their homomorphisms 
by
$\langSigCat$, and the full subcategory of unary $\langCartClosed$-signatures by
$\langUnSigCat$. 
\end{mydefn}

\begin{mynotation}[{\cf~Notation~\ref{not:stlc-alltypes-not}}] 
For a $\langCartClosed$-signature $\sig = (\baseTypes, \graph)$, we write 
$\allTypes\baseTypes$ for the set generated from $\baseTypes$ by the 
grammar~(\ref{eq:langcartclosed-sig-grammar}).
In particular, when the signature is just a set (\ie~the graph $\graph$ has no edges) we denote the 
signature $\sig = (\baseTypes, \graph)$ simply by $\allTypes\baseTypes$. 
\end{mynotation}

The embedding $\inc : \stlcTimesUnSigCat \hookrightarrow \stlcTimesSigCat$ 
has a right adjoint by an argument similar to that for 
Lemma~\ref{lem:stlcTimesSig-reflective-subcat} 
(\cf~also Lemma~\ref{lem:inc-adjoint-cartesian-structure-signatures}). 

The definition of cartesian closed pseudofunctor follows the template given by 
cartesian pseudofunctors of biclones, while the construction of the free 
cartesian closed biclone on a $\langCartClosed$-signature echoes that 
for the free cartesian closed clone on a $\stlc$-signature 
(Lemma~\ref{lem:free-cc-clone-on-sig}).

\begin{mydefn}
Let $\ccClone{S}{\biclone}$ and $\ccClone{T}{\altBiclone}$ be cartesian closed 
biclones. A \Def{cartesian closed pseudofunctor} 
$(F, \prodPres, \expPres) : \ccClone{S}{\biclone} \to \ccClone{T}{\altBiclone}$
is a cartesian pseudofunctor 
$(F, \prodPres) : \cartClone{S}{\biclone} \to \cartClone{T}{\biclone}$
equipped with a choice of equivalence
$\evBar_{A,B} : 
F(\expobj{A}{B}) \leftrightarrows \expobj{FA}{FB} : 
\expPres_{A,B}$
for every $A, B \in S$, where 
$\evBar_{A,B}
:= \lambda{\big(F\eval_{A,B} \big)}$.
We 
call $(F, \prodPres, \expPres)$ \Def{strict} if $(F, \prodPres)$ is a strict 
cartesian pseudofunctor such that
\begin{align*}
F(\expobj{A}{B}) &= (\expobj{FA}{FB}) \\
F(\eval_{A,B}) &=\eval_{FA,FB} \\
F(\lambda t) &=\lambda(Ft) \\
F(\epsilonExp_t) &=\epsilonExp_{Ft} \\
\expPres_{A,B} &=\Id_{\scriptsizeexpobj{FA}{FB}}
\end{align*}
and the isomorphisms witnessing the adjoint equivalences are the canonical 
$2$-cells
\[
\Id_{(\scriptsizeexpobj{FA}{FB})}
   \XRA{\etaExp{\Id}}
\lambda{\left( 
	\cslr
		{\eval_{FA, FB}}
		{\cs
			{\Id_{(\scriptsizeexpobj{FA}{FB})}}
			{\p{1}{(\scriptsizeexpobj{FA}{FB}), FA}}, 
		\p{2}{(\scriptsizeexpobj{FA}{FB}), FA}} 
\right)}
	\XRA\iso
   \lambda(\eval_{FA,FB})
\]
obtained from the unit and the canonical structural isomorphism. 
\end{mydefn}


For the construction of the free cc-biclone, it will be useful to introduce 
some notation. 
For 
$t : A \to B$ 
we define 
$t \times X := \pair{\cs{t}{\pi_1}, \cs{\Id_X}{\pi_2}} : \prodop_2(A, X) \to \prodop_2(B, X)$, and similarly on 2-cells.

\begin{myconstr} \label{constr:free-cc-biclone}
For any $\langCartClosed$-signature $\sig$, define a cartesian closed biclone 
$\freeCartClosedBiclone{\sig}$ with sorts generated by the grammar
\[
A_1, \,\dots\, , A_n, C, D ::= 
	B \st 
	\prodop_n(A_1, \,\dots\, , A_n)  \st
	\exptype{C}{D} \qquad (B \in \baseTypes, n \in \Nat)
\] 
by extending
Construction~\ref{constr:free-cart-biclone} (page~\pageref{constr:free-cart-biclone}) 
with the following rules:
\begin{center}
\unaryRule
	{\phantom{\freeCartClosedBiclone{\sig}}}
	{\evalterm_{B,C} \in \freeCartClosedBiclone{\sig}(\exptype{B}{C}, B; C)}
	{}
\qquad
\unaryRule
	{t \in \freeCartClosedBiclone{\sig}(A_1, \,\dots\, , A_n, B; C)}
	{\lambda t \in \freeCartClosedBiclone{\sig}(A_1, \,\dots\, , A_n; \exptype{B}{C})}
	{}

\unaryRule
	{t \in \freeCartClosedBiclone{\sig}(A_1, \,\dots\, , A_n, B; C)}
	{\epsilonExp_t \in 
		\freeCartClosedBiclone{\sig}(A_1, \,\dots\, , A_n, B; C)
			\left(\cslr
				{\evalterm_{B,C}}
				{\cs{(\lambda t)}{\p{1}{\ind{A}, B}, \,\dots\, , \p{n}{\ind{A}, 
				B}}
					, \p{n+1}{\ind{A}, B}} , 
				t \right)}
	{} \vspace{-0.5\treeskip}

\begin{prooftree}
\AxiomC{$u \in \freeCartClosedBiclone{\sig}(A_1, \,\dots\, , A_n; \exptype{B}{C})$}
\noLine
\UnaryInfC{$\alpha \in 
		\freeCartClosedBiclone{\sig}(A_1, \,\dots\, , A_n, B; C)
			\left( \cslr
				{\evalterm_{B,C}}
				{\cs{u}{\p{1}{\ind{A}, B}, \,\dots\, , \p{n}{\ind{A}, B}}
					, \p{n+1}{\ind{A}, B}} , 
			   t\right)$}
\UnaryInfC{$\transExp{\alpha} \in 	
		\freeCartClosedBiclone{\sig}(A_1, \,\dots\, , A_n; \exptype{A}{B})(u, \lambda t)$}
\end{prooftree}
\end{center}
The equational theory $\equiv$ is that of
Construction~\ref{constr:free-cart-biclone}, extended by requiring that
\begin{itemize}
\item 	
For every 
		$\alpha :  \cslr
				{\evalterm_{B,C}}
				{\cs{u}{\p{1}{\ind{A}, B}, \,\dots\, , \p{n}{\ind{A}, B}}
					, \p{n+1}{\ind{A}, B}} \To t
				: A_1, \,\dots\, , A_n, B \to C$, 
\[
\alpha \equiv 
	\epsilon_t \vert 
		\cslr
				{\evalterm_{B,C}}
				{\cs{\transExp{\alpha}}
					{\p{1}{\ind{A}, B}, \,\dots\, , \p{n}{\ind{A}, B}}
					, \p{n+1}{\ind{A}, B}}
\]

\item For every $\gamma : u \To \lambda t : A_1, \,\dots\, , A_n \to 
(\exptype{A}{B})$,
\[
\gamma \equiv 
	\transExplr{
			\epsilon_t \vert 
				\cslr
				{\evalterm_{B,C}}
				{\cs{\gamma}
					{\p{1}{\ind{A}, B}, \,\dots\, , \p{n}{\ind{A}, B}}
					, \p{n+1}{\ind{A}, B}}
			}
\]
\item 
If $\alpha \equiv \alpha' : 
				\cs
				{\evalterm_{B,C}}
				{u \times B} \To t
				: X_1, \,\dots\, , X_n, B \to C$ then
	$\transExp{\alpha} \equiv \transExp{\alpha'}$. 
\end{itemize}
Finally we require that every $\epsilonExp_t$ and
$\transExp{\id_{\cs
				{\evalterm}
				{\prod_2(u, B)}}}$ 
is invertible.
\end{myconstr}

It follows that for any 2-cell
$\alpha : \cslr
				{\evalterm_{B,C}}
				{\cs{u}
					{\p{1}{\ind{A}, B}, \,\dots\, , \p{n}{\ind{A}, B}}
					, \p{n+1}{\ind{A}, B}}
			\To {t : A_1, \,\dots\, , A_n, B \to C}$,
$\transExp\alpha$ is the unique 2-cell $\gamma$ of type
$u \To \lambda t$ 
such that 
$\alpha \equiv 
	\epsilon_t \vert 
		\cslr
				{\evalterm_{B,C}}
				{\cs{\gamma}
					{\p{1}{\ind{A}, B}, \,\dots\, , \p{n}{\ind{A}, B}}
					, \p{n+1}{\ind{A}, B}}$.
Existence is the first equation and uniqueness follows by the latter two
(\cf~Lemma~\ref{lem:free-cart-biclone-cartesians-structure}).

The required universal property extends that for cartesian biclones. 

\begin{mylemma} \label{lem:free-cc-biclone}
For any $\langCartClosed$-signature $\sig$, 
cartesian closed biclone $\ccClone{T}{\altBiclone}$
and $\langCartClosed$-signature homomorphism 
$h : \sig \to \altBiclone$
from $\sig$ to the $\langCartClosed$-signature underlying $\altBiclone$, 
there exists a unique strict cartesian closed pseudofunctor
$\ext{h} : \freeCartClosedBiclone{\sig} \to \altBiclone$
such that $\ext{h} \circ \inc = h$, for 
$\inc : \sig \hookrightarrow \freeCartClosedBiclone{\sig}$
the inclusion.
\begin{proof}
We extend the strict cartesian pseudofunctor $\ext{h}$ defined in 
Lemma~\ref{lem:free-cart-biclone-proved} (page~\pageref{lem:free-cart-biclone-proved}) 
with the following rules:
\begin{align*}
\ext{h}(\exptype{B}{C}) &:= (\exp{\ext{h}A}{\ext{h}B}) \\[5pt]
\ext{h}(\evalterm_{B,C}) &:= \eval_{\ext{h}B, \ext{h}C} \\
\ext{h}(\lambda t) &:= \lambda \big( \ext{h} t \big) \\[5pt]
\ext{h}(\epsilonExp_t) &:= \epsilonExp_{\ext{h}t} \\
\ext{h}(\transExp{\alpha}) &:= \transExp{\ext{h}\alpha}
\end{align*}
For uniqueness, it suffices to show that any strict cartesian closed 
pseudofunctor commutes with the $\transExp{-}$ operation. For this we use the 
universal property. Let $F : \freeCartClosedBiclone{\sig} \to \altBiclone$ be any cartesian closed 
pseudofunctor. Then, for any 
$\alpha : \cslr
				{\evalterm_{B,C}}
				{\cs{u}
					{\p{1}{\ind{A}, B}, \,\dots\, , \p{n}{\ind{A}, B}}
					, \p{n+1}{\ind{A}, B}}
			\To t : {A_1, \,\dots\, , A_n, B \to C}$
in $\freeCartClosedBiclone{\sig}$, 
\begin{align*}
\epsilon_{Ft} &\vert 
		\cslr
				{\eval_{FB,FC}}
				{\cslr{\left(F\transExp{\alpha}\right)}
					{\p{1}{F\ind{A}, FB}, \,\dots\, , \p{n}{F\ind{A}, FB}}
					, \p{n+1}{F\ind{A}, FB}} \\
	&= F(\epsilon_{t}) \vert 
		F\left(
			\cslr
				{\evalterm_{B,C}}
				{\cslr{\transExp{\alpha}}
					{\p{1}{\ind{A}, B}, \,\dots\, , \p{n}{\ind{A}, B}}
					, \p{n+1}{\ind{A}, B}} 
		\right) 
		&\text{by strict preservation} \\
	&= F\left(
			\epsilon_t \vert 
		\cslr
				{\evalterm_{B,C}}
				{\cs{\transExp{\alpha}}
					{\p{1}{\ind{A}, B}, \,\dots\, , \p{n}{\ind{A}, B}}
					, \p{n+1}{\ind{A}, B}}
		\right) \\
	&= F\alpha
\end{align*} 
Hence $\transExp{F\alpha}$ must equal $F{\left(\transExp\alpha\right)}$.
\end{proof}
\end{mylemma}

We saw in Example~\ref{exmp:nucleus-of-fp-not-free} (page~\pageref{exmp:nucleus-of-fp-not-free}) 
that the free fp-bicategory on a 
$\langCart$-signature cannot arise as the nucleus of the free cartesian biclone over the same signature. 
We can now see that the addition of exponentials introduces a further obstacle
(\cf~Remark~\ref{rem:exponentials-preserved-proof-fails}). Let 
$\sig$ be a unary $\langCartClosed$-signature and 
$\overline{\freeCartClosedBiclone{\sig}}$ 
be its nucleus. Just as in the categorical case, the 
maps $\pi_i$ in $\overline{\freeCartClosedBiclone{\sig}}$ are the biuniversal arrows 
defining 
products in $\freeCartClosedBiclone{\sig}$, but the evaluation map in 
$\overline{\freeCartClosedBiclone{\sig}}$ is
$\cs{\evalterm_{B,C}}{\pi_1,\pi_2}$~(recall~Remark~\ref{rem:bicat-exp-structure-not-preserved}).
It follows that for any
cc-bicategory $\ccBicat{\baseCat}$ and strict cc-pseudofunctor 
$F : \overline{\freeCartClosedBiclone{\sig}} \to \baseCat$ one must have
\begin{equation} \label{eq:strict-preservation-on-eval}
\begin{aligned}
\eval_{FB,FC} &= F{\left(\cs{\evalterm_{B,C}}{\pi_1, \pi_2} \right)} \\
	&= F{\left({\evalterm_{B,C}} \circ \seq{\pi_1, \pi_2} \right)} 
		\qquad &\text{\small by def. of products in 
		$\overline{\freeCartClosedBiclone{\sig}}$} \\
	&= F(\evalterm_{B,C}) \circ F\seq{\pi_1, \pi_2}  \\
	&=	F(\evalterm_{B,C}) \circ \seq{\pi_1, \pi_2} 
		&\text{\small by strict preservation}
\end{aligned}
\end{equation}
In particular, since
$\ext{h}(\evalterm_{B,C}) = \eval_{\ext{h}B, \ext{h}C}$, 
the restriction $\overline{\ext{h}}$ of $\ext{h}$ to unary multimaps cannot be 
strictly cartesian closed whenever 
$\eval_{\ext{h}B, \ext{h}C} \circ \seq{\pi_1, \pi_2} \neq
	\eval_{\ext{h}B, \ext{h}C}$
in the target cc-bicategory. This occurs, for instance, in the cc-bicategories 
of generalised species~\cite{FioreSpecies} and concurrent 
games~\cite{PaquetThesis}. 

One way to diagnose the problem is the chain of 
equivalences~(\ref{eq:biclone-representable-iff-products}). The product 
structure in a cartesian closed biclone arises via the $\prodop_n(-)$ 
operation, but exponentials are defined with respect to context extension. This 
mismatch makes it impossible for $\overline{\ext{h}}$ to strictly preserve both 
products and exponentials. To construct the free cc-bicategory over a unary 
signature, one must define 
exponentials directly with respect to products, resulting in a construction 
similar to that given in~\cite{Ouaknine1997}. 

\paragraph*{The free cc-bicategory.}
As for Construction~\ref{constr:free-cc-biclone}, we write $t \times B$ for the 
(derived) arrow
$\pair{\cs{t}{\pi_1}, \cs{\Id}{\pi_2}}$, and likewise on 2-cells.

\begin{myconstr} \label{constr:free-cc-bicat}
For any unary $\langCartClosed$-signature $\sig = (\baseTypes, \graph)$, define 
a cc-bicategory
$\freeCartClosedBicat{\sig}$ as follows. The objects are generated by the grammar
\[
A_1, \,\dots\, , A_n, C, D ::= 
	B \st 
	\prodop_n(A_1, \,\dots\, , A_n)  \st
	\exptype{C}{D} \qquad (B \in \baseTypes, n \in \Nat)
\] 
For 1-cells and 2-cells, one takes 
the deductive system defining the free fp-bicategory on $\sig$ 
(Lemma~\ref{lem:free-fp-bicat-from-free-biclone}, page~\pageref{lem:free-fp-bicat-from-free-biclone}), extended as follows.
For 1-cells:
\begin{center}
\unaryRule
	{\phantom{\freeCartClosedBicat{\sig}(X \times B; C)}}
	{\evalterm_{B,C} \in \freeCartClosedBicat{\sig}(\exptype{B}{C} \times B; C)}
	{} \quad
\unaryRule
	{t \in \freeCartClosedBicat{\sig}(X \times B; C)}
	{\lambda t \in \freeCartClosedBicat{\sig}(X, \exptype{B}{C})}
	{}
\vspace{-\treeskip}
\end{center}
For 2-cells:
\begin{figure}[!h]
\centerfloat
\begin{bprooftree}
\AxiomC{$\phantom{\big(\big)}$}
\noLine
\UnaryInfC{$t \in \freeCartClosedBicat{\sig}(X \times B, C)$}
\UnaryInfC{$\epsilonExp_t \in 
		\freeCartClosedBicat{\sig}(X \times B, C)
			\big(\cs
				{\evalterm_{B,C}}
				{\lambda t \times B} , 
				t \big)$}
\end{bprooftree}
\begin{bprooftree}
\AxiomC{$u \in \freeCartClosedBicat{\sig}(X, \exptype{B}{C})$}
\noLine
\UnaryInfC{$\alpha \in 
		\freeCartClosedBicat{\sig}(X \times B, C)
			\big( \cs
				{\evalterm_{B,C}}
				{u \times B} , 
				t \big)$}
\UnaryInfC{$\transExp{\alpha} \in 	
		\freeCartClosedBicat{\sig}(X, \exptype{A}{B})(u, \lambda t)$}
\end{bprooftree}
\end{figure}

Moreover, we extend the equational theory of 
Lemma~\ref{lem:free-fp-bicat-from-free-biclone} with the following three rules:
\begin{itemize}
\item 	
For every 
		$\alpha :  
				\cs
				{\evalterm_{B,C}}
				{u \times B} \To t
				: X \times B \to C$, 
\[
\alpha \equiv 
	\epsilon_t \vert 
		\cs
			{\evalterm_{B,C}}
			{\transExp{\alpha} \times B}
\]

\item For every $\gamma : u \To \lambda t : X \to (\exptype{A}{B})$,
\[
\gamma \equiv 
	\transExplr{
			\epsilon_t \vert 
		\cs
			{\evalterm_{B,C}}
			{\gamma \times B}}
\]
\item 
If $\alpha \equiv \alpha' : 
				\cs
				{\evalterm_{B,C}}
				{u \times B} \To t
				: X \times B \to C$ then
	$\transExp{\alpha} \equiv \transExp{\alpha'}$. 
\end{itemize} \enlargethispage*{2\baselineskip}
Finally we require that every $\epsilonExp_t$ and
$\transExp{\id_{	\cs
				{\evalterm}
				{u \times B}}}$ is invertible.
\end{myconstr}

The bicategory $\freeCartClosedBicat{\sig}$ is cartesian closed by exactly the same 
argument as for the 
biclone $\freeCartClosedBiclone{\sig}$. The associated free property is similarly 
straightforward.

\begin{mylemma} \label{lem:free-cc-bicat-on-sig}
For any unary $\langCartClosed$-signature $\sig$, 
cc-bicategory $\ccBicat{\altCat}$
and \mbox{$\langCartClosed$-signature} homomorphism 
$h : \sig \to \altCat$
from $\sig$ to the $\langCartClosed$-signature underlying $\altCat$, 
there exists a unique strict cartesian closed pseudofunctor
$\ext{h} : \freeCartClosedBicat{\sig} \to \altCat$
such that $\ext{h} \circ \inc = h$, for 
$\inc : \sig \hookrightarrow \freeCartClosedBicat{\sig}$
the inclusion.
\begin{proof}
We extend the strict cartesian pseudofunctor $\ext{h}$ defined in 
Lemma~\ref{lem:free-fp-bicat-from-free-biclone} (page~\pageref{lem:free-fp-bicat-from-free-biclone}) 
as follows:
\begin{align*}
\ext{h}(\exptype{B}{C}) &:= (\exp{\ext{h}A}{\ext{h}B}) \\[5pt]
\ext{h}(\evalterm_{B,C}) &:= \eval_{\ext{h}B, \ext{h}C} \\
\ext{h}(\lambda t) &:= \lambda \big( \ext{h} t \big) \\[5pt]
\ext{h}(\epsilonExp_t) &:= \epsilonExp_{\ext{h}t} \\
\ext{h}\left(\transExp{\alpha}\right) &:= \transExp{\ext{h}\alpha}
\end{align*}
For uniqueness, it suffices to show that any strict cartesian closed 
pseudofunctor commutes with the $\transExp{-}$ operation. The proof is as in 
Lemma~\ref{lem:free-cc-biclone} (or, more abstractly, follows from
Lemma~\ref{lem:strict-preservation-strict-pres-UMP}).
\end{proof}
\end{mylemma}

The preceding lemma entails that one may construct a type 
theory for cartesian 
closed bicategories by synthesising the internal language of 
$\freeCartClosedBicat{\sig}$. 
Within this `bicategorical' (rather than \emph{biclone-theoretic}) type theory
the variables play 
almost no role. 
For instance, the lambda abstraction rule takes on the 
following form:
\begin{center}
\binaryRule
	{p : A \times B \vdash t : C}
	{q \text{ fresh }}
	{q : A \vdash \lambda(q, p \bind t) : \exptype{B}{C}}
	{lam}\vspace{-\treeskip}
\end{center}
The variable $p$ is bound, but $q$ is free. 
It is possible to place such rules within the general framework of 
binding signatures,
and the syntactic model of the resulting type theory is 
biequivalent to the syntactic model of the type theory  
extracted from the construction of $\freeCartClosedBiclone{\sig}$, 
restricted to unary 
contexts. However, the result is rather
alien to the usual conception of a type theory. 
We therefore call the 
internal language of $\freeCartClosedBiclone{\sig}$ the `type theory for cartesian closed 
bicategories'. In Section~\ref{sec:syntactic-model-free-property} we shall show 
that this terminology is warranted. 

The freeness universal property of $\freeCartClosedBicat{\sig}$ also entails an up-to-equivalence uniqueness property we 
shall employ later. We begin by stating a result for the case where the 
signature is just a set; 
thereafter we employ slightly stronger hypotheses to handle constants.
We write $t : A_1, \,\dots\, , A_n \to B$ and 
$\tau : t \To t' : A_1, \,\dots\, , A_n \to B$
for 1-cells and 2-cells in $\freeCartClosedBicat{\sig}$. 
%
%
%
%

\begin{mylemma} \label{lem:syntactic-model-uniqueness-up-to-equivalence}
Let $\sig = (\baseTypes, \graph)$ be a unary $\langCartClosed$-signature 
for which $\graph$ is a set, 
$\ccBicat{\baseCat}$ be a 
cc-bicategory and $h : \sig \to \altCat$ be a 
$\langCartClosed$-signature homomorphism. Then, for any \mbox{cc-pseudofunctor} 
$(F, \prodPres, \expPres)$ such that the following diagram commutes,
\begin{equation} \label{eq:termcatCom-free-up-to-equivalence}
\begin{tikzcd}
\freeCartClosedBicat{\sig} 
\arrow{r}{F} &
\altCat \\
\sig
\arrow[hookrightarrow]{u}
\arrow[swap]{ur}{h} &
\:
\end{tikzcd}
\end{equation}
there exists an equivalence $F \simeq \ext{h}{}$ between $F$ and the 
canonical cc-pseudofunctor extending $h$.
\begin{proof}
We construct a pseudonatural transformation 
$(\natTrans, \natCell) : F \To \ext{h}{}$ 
whose components are all equivalences. We define the components $\natTrans_X$ 
and their pseudo-inverses $\psinv{\natTrans}_X$ by mutual induction as follows:
\begin{gather*}
\begin{aligned}
\natTrans_{B} &:= FB \xra{=} hB \xra{\Id_{hB}} hB \xra{=} \ext{h}{B} 
	\qquad \text{ for } B \in \baseTypes \\
\psinv{\natTrans}_B &:= \ext{h}{B} \xra{=} hB \xra{\Id_{hB}} hB \xra{=} FB
\end{aligned} 
\: \\[8pt]
\begin{aligned}
\natTrans_{(\prod_n \ind{A})} &:=
	F{\left(\prodop_n \ind{A}\right)} 
	\xra{\seq{F\pi_1, \dots, F\pi_n}}
	\prodop_{i=1}^n F(A_i) 
	\xra{\prod_{i=1}^n \natTrans_{A_i}}
	\prodop_{i=1}^n \ext{h}{A_i} \\
\psinv{\natTrans}_{(\prod_n \ind{A})} &:=
	\prodop_{i=1}^n \ext{h}{A_i}
	\xra{\prod_{i=1}^n \psinv{\natTrans}_{A_i}}
	\prodop_{i=1}^n F(A_i) 
	\xra{\prodPres_{\ind{A}}}
	F{\left(\prodop_n \ind{A}\right)} 
\end{aligned}
\: \\[8pt]
\begin{aligned}
\natTrans_{(\scriptsizeexpobj{X}{Y})} &:= 
	F(\exptype{X}{Y}) 
	\xra{\evBar_{X,Y}} 
	\left(\exptype{FX}{FY}\right) 
	\xra{\scriptsizeexpobj{\psinv{\natTrans}_{X}}{\natTrans_{Y}}} 
	\left(\exptype{\ext{h}{X}}{\ext{h}{Y}}\right) \\
\psinv{\natTrans}_{(\scriptsizeexpobj{X}{Y})} &:=
	\left(\exptype{\ext{h}{X}}{\ext{h}{Y}}\right)
	\xra{\scriptsizeexpobj{{\natTrans}_{X}}{\psinv{\natTrans}_{Y}}}
	\left(\exptype{FX}{FY}\right)
	\xra{\expPres_{X,Y}}
	F(\exptype{X}{Y}) 
\end{aligned}
\end{gather*}
We denote the 
unit and counit of the equivalence 
\[
\natTrans_X: FX \leftrightarrows \ext{h} X : \psinv{\natTrans}_X
\]
by
$\un_X : \Id_{FX} \To
	 \psinv{\natTrans}_X \circ \natTrans_X$ 
and 
$\co_X : \natTrans_X \circ \psinv{\natTrans}_X \To 
	\Id_{\ext{h}{X}}$, 
respectively, and assume without loss of generality that they satisfy the two 
triangle laws.

We now
construct the witnessing 2-cells 
$\natCell_t : \natTrans_B \circ Ft \To \ext{h}(t) \circ \natTrans_A$
by induction.

For identities, the definition is forced upon us by the unit law of a 
pseudonatural transformation. We define 
\[
\natCell_{\Id_A} := 
\natTrans_{A} \circ F(\Id_A) 
\XRA{\natTrans_{A} \circ (\psi^F_{A})^{-1}}
\natTrans_{A} \circ \Id_{F(A)} 
\XRA{\iso}
\Id_{\ext{h}{(A)}} \circ \natTrans_{A} 
\]
For the product structure, we define 
$\natCell_{\pi_k}$ and 
$\natCell_{\pair{t_1, \,\dots\, , t_n}}$
by the commutativity of the 
following diagrams:
\begin{td}[column sep = 1em]
\natTrans_{A_k} \circ 
	F\pi_k 
\arrow{rr}{\natCell_{\pi_k}}
\arrow[swap]{d}
	{\natTrans_{A_k} \circ \epsilonTimesInd{-k}{}} &
\: &
\ext{h}{(\pi_k)} \circ 
	\natTrans_{(\prod_n \ind{A})} \\
\natTrans_{A_k} \circ \left(\pi_k \circ \seq{F\ind{\pi}}\right)
\arrow[swap, bend right = 8]{dr}{\iso} &
\: &
\left(\pi_k \circ \prod_{i=1}^n \natTrans_{A_i}\right) 
	\circ \seq{F\ind{\pi}}
\arrow[swap]{u}{\iso} \\

\: &
\left(\natTrans_{A_k} \circ \pi_k\right) \circ \seq{F\ind{\pi}}
\arrow[swap, bend right = 8]{ur}	
	{\epsilonTimesInd{-k}{} \circ \seq{F\ind{\pi}}} &
\: \\

\left(\prod_{i=1}^m \natTrans_{A_i} \circ \seq{F\ind{\pi}}\right) \circ
	F{\left(\pair{t_1, \,\dots\, , t_m} \right)} 
\arrow{rr}[yshift=0mm]{\natCell_{\pair{t_1, \,\dots\, , t_m}}}
\arrow[swap]{d}{\iso} &
\: &
\ext{h}{\left(\pair{t_1, \,\dots\, , t_m}\right)} \circ \natTrans_{X} \\
\left(\prod_{i=1}^m \natTrans_{A_i}\right) \circ \left(\seq{F\ind{\pi}} \circ
	F{\left(\pair{t_1, \,\dots\, , t_m} \right)}\right) 
\arrow[swap]{d}{(\prod_i \natTrans_{A_i}) \circ \unpack} &
\: &
\seq{\ext{h}{(\ind{t})}} \circ \natTrans_{X} 
\arrow[equals]{u} \\
\left(\prod_{i=1}^m \natTrans_{A_i}\right) \circ  
	\seq{F(\ind{t})}
\arrow[swap]{r}{\fuse} &
\seq{ \natTrans_{\ind{A}} \circ 
	F(\ind{t})}
\arrow[swap]{r}{\seq{\natCell_{t_1}, \,\dots\, , \natCell_{t_m}}} &
\seq{\ext{h}{(\ind{t})} \circ \natTrans_{X}} 
\arrow[swap]{u}{\postName^{-1}} 
\end{td}
The \rulename{eval} and \rulename{lam} cases require more work, but are in a 
similar spirit. 

\vspace{1em}
\subproof{\rulename{eval} case.} 
We are required to give an invertible 2-cell filling the diagram
\begin{td}
F\big((\exptype{A}{B}) \times A\big) 
\arrow[
swap,
rounded corners,
to path=
{ -- ([xshift=1ex]\tikztostart.west)
-| ([xshift=-1cm]\tikztotarget.west)
-- ([xshift=-1cm]\tikztotarget.west)
-- (\tikztotarget.west)}, 
]{dd}
\arrow[phantom]{dd}[font=\scriptsize, xshift=-3.7cm, description]
	{(\natTrans_{(\scriptsizescriptsizeexpobj{A}{B})} \times \natTrans_{A}) 
		\circ \seq{F\pi_1, F\pi_2}}
\arrow[phantom]{ddrr}[description]{\twocell{\natCell_{\eval}}}
\arrow[swap]{d}{\seq{F\pi_1, F\pi_2}}
\arrow{rr}{F\evalterm_{A,B}} &
\: &
FB 
\arrow{dd}{\natTrans_{B}} \\

F(\exptype{A}{B}) \times F(A) 
\arrow[swap]{d}
	{\natTrans_{(\scriptsizescriptsizeexpobj{A}{B})} \times \natTrans_{A} } &
\: &
\: \\

\ext{h}{(\exptype{A}{B})} \times \ext{h}{A} 
\arrow[equals]{r}{} &
(\exp{\ext{h}{A}}{\ext{h}{B}}) \times \ext{h}{A} 
\arrow[swap]{r}{\eval} &
\ext{h}{B}
\end{td}
To this end, first define an invertible 2-cell $\delta_{A,B}$ 
applying the counit $\epsilonExp$ as far as possible:
\begin{td}
\eval_{\ext{h}{A}, \ext{h}{B}} \circ \left( 
	\natTrans_{(\scriptsizeexpobj{A}{B})} \times 
			\natTrans_{A}\right) 
\arrow[equals]{d} 
\arrow[bend left = 40]{ddddddr}{\delta_{A,B}} &
\: \\

\eval_{\ext{h}{A}, \ext{h}{B}} \circ \left( 
	(\exp{\psinv{\natTrans}_{A}}{\natTrans_{B}}) 
	\circ 			\evBar^F_{A, B} \times 
			\natTrans_{A}\right)
\arrow[swap]{d}{\iso} &
\: \\

\left(\eval_{\ext{h}{A}, \ext{h}{B}} 
	\circ \big( (\exp{\psinv{\natTrans}_{A}}{\natTrans_{B}}) \times 
			\ext{h}{A}\big)\right)
	\circ (\evBar^F_{A, B} \times \natTrans_{A})
\arrow[swap]{d}
	{\epsilonExp_{(\natTrans \circ \eval 
	\circ (\Id \times \psinv{\natTrans}))} 
		\circ (\evBar^F_{A, B} \times \natTrans_{A})} &
\: \\

\left( \left(\natTrans_{B} 
	\circ \eval_{FA, FB}\right) 
	\circ \big(\Id_{(\scriptsizeexpobj{FA}{FB})} \times 
					\psinv{\natTrans}_{A}\big)\right) 
	\circ (\evBar^F_{A, B} \times \natTrans_{A}) 
\arrow[swap]{d}{\iso} &
\: \\

\left(\natTrans_{B} 
	\circ \left(\eval_{FA, FB} 
	\circ \big(\evBar^F_{A, B} \times FA \big)\right)\right) 
	\circ \big( \Id_{(\scriptsizeexpobj{FA}{FB})} \times 
					\psinv{\natTrans}_{A}\natTrans_{A} 
				\big) 
\arrow[swap]{d}{\natTrans \circ 
	\epsilonExp_{(F(\eval) \circ \prodPres)} \circ 
	( \Id \times 
				\psinv{\natTrans}\natTrans)
			} &
\: \\

\left( \natTrans_B 
	\circ \left( F(\eval_{A,B})
	\circ \prodPres_{\scriptsizeexpobj{A}{B}, A}\right)\right)
	\circ \big( \Id_{(\scriptsizeexpobj{FA}{FB})} \times 
					\psinv{\natTrans}_{A}\natTrans_{A} 
				\big) 
\arrow[swap]{d}
	{\natTrans \circ F\eval \circ \prodPres \circ (\Id \times \un_A^{-1})} &
\: \\

\left( \natTrans_B 
	\circ \left( F(\eval_{A,B})
	\circ \prodPres_{\scriptsizeexpobj{A}{B}, A}\right)\right)
	\circ \big( \Id_{(\scriptsizeexpobj{FA}{FB})} \times 
					\Id_{FA} 
				\big) 
\arrow[swap]{r}{\iso} &
\left( \natTrans_B \circ F(\evalterm_{A,B}) \right) \circ 
	\prodPres_{\scriptsizeexpobj{A}{B}, A}
\end{td}
\vspace{2mm}
Then define $\natCell_{\evalterm}$ to be the composite
\begin{td}
\natTrans_{B} \circ F(\eval_{A,B}) 
\arrow[swap]{d}{\iso}
\arrow{r}{\natCell_{\evalterm}} &
\eval_{\ext{h}{A}, \ext{h}{B}}
	 \circ \left(\left( \natTrans_{(\scriptsizeexpobj{A}{B})} \times 
			\natTrans_{A}\right) 
	\circ \seq{F\pi_1, F\pi_2}\right) \\

\left(\natTrans_{B} \circ F(\eval_{A,B})\right) \circ 
	\Id_{F{\left((\scriptsizeexpobj{A}{B}) \times A\right)}}
\arrow[swap]{d}
	{\left(\natTrans_{B} 
		\circ F(\eval_{A,B})\right) 
		\circ (\coTimes_{\scriptsizescriptsizeexpobj{A}{B}, A})^{-1}} &
\: \\

\left(\natTrans_{B} 
		\circ F(\eval_{A,B})\right) 
		\circ \left( \prodPres_{\scriptsizeexpobj{A}{B}, A} 
		\circ \seq{F\pi_1, F\pi_2} 
	\right)
\arrow[swap]{d}{\iso} &
\: \\

\left(\natTrans_{B} 
		\circ \big(F(\eval_{A,B}) 
		\circ \prodPres_{\scriptsizeexpobj{A}{B}, A} \big)\right)
		\circ \seq{F\pi_1, F\pi_2} 
\arrow[swap, bend right = 8]{dr}{\delta_{A,B}^{-1} \circ \seq{F\pi_1, F\pi_2}} &
\: \\

\: &
\left(\eval_{\ext{h}{A}, \ext{h}{B}}
	 \circ \left( \natTrans_{(\scriptsizeexpobj{A}{B})} \times 
			\natTrans_{A}\right)\right) 
	\circ \seq{F\pi_1, F\pi_2}
\arrow[swap]{uuuu}{\iso}
\end{td}

\newpage
\subproof{\rulename{lam} case.} Suppose $ t: Z \times A \to B$. By induction we 
are given $\natCell_t$ 
filling
\begin{td}
F(Z \times A) 
\arrow[
swap,
rounded corners,
to path=
{ -- ([xshift=1ex]\tikztostart.west)
-| ([xshift=-1cm]\tikztotarget.west)
-- ([xshift=-1cm]\tikztotarget.west)
-- (\tikztotarget.west)}, 
]{dd}
\arrow[phantom]{dd}[font=\scriptsize, xshift=-3.1cm, description]
	{(\natTrans_Z \times \natTrans_A) \circ \seq{F\pi_1, F\pi_2}}
\arrow[swap]{d}{\seq{F\pi_1, F\pi_2}}
\arrow[phantom]{ddrr}[description]{\twocell{\natCell_t}}
\arrow{rr}{Ft} &
\: &
FB
\arrow{dd}{\natTrans_{B}} \\

FZ \times FA 
\arrow[swap]{d}{\natTrans_{Z} \times \natTrans_{A}} &
\: &
\: \\

\ext{h}{(Z)} \times \ext{h}{(A)} 
\arrow[equals]{r}{} &
\ext{h}{(Z \times A)} 
\arrow[swap]{r}[yshift=0mm]{\ext{h}{t}} &
\ext{h}{B}
\end{td}
and we are required to fill the diagram
\begin{td}[column sep = 8em]
FZ 
\arrow[swap]{dd}{\natTrans_{Z}}
\arrow[phantom]{ddr}[description, xshift=-3mm]{\twocell{\natCell_{\lambda t}}}
\arrow{r}{F(\lambda t)} &
F(\exptype{A}{B}) 
\arrow[
swap,
rounded corners,
to path=
{ -- ([xshift=1ex]\tikztostart.east)
-| ([xshift=1cm]\tikztotarget.east)
-- ([xshift=1cm]\tikztotarget.east)
-- (\tikztotarget.east)}, 
]{dd}
\arrow[phantom]{dd}[font=\scriptsize, xshift=2.8cm, description]
	{(\scriptsizeexpobj{\psinv{\natTrans}_{A}}{\natTrans_{B}}) 
		\circ \evBar^F_{A,B}}
\arrow{d}{\evBar^F_{A,B}} \\

\: &
\big(\exp{FA}{FB} \big)
\arrow{d}{ (\scriptsizeexpobj{\psinv{\natTrans}_{A}}{\natTrans_{B}}) } \\

\ext{h}{Z} 
\arrow[swap]{r}{\ext{h}{\left(\lambda t\right)}} &
(\exp{\ext{h}{A}}{\ext{h}{B}})
\end{td}
Our strategy is the following. Writing $\cl$ for the clockwise 
composite 
around the preceding diagram, we define a 2-cell 
\[
\zeta_{A,B} : \eval_{\ext{h}{A}, \ext{h}{B}} \circ 
	(\cl \times \ext{h}{A})
	\To 
	\ext{h}{(t)} \circ 
		(\natTrans_{Z} \times \ext{h}A)
\]
so that 
$\transExp{\zeta_{A,B}} : \cl \To 
	\lambda {\left(\ext{h}{(t)} \circ 
			(\natTrans_{Z} \times \ext{h}A)\right)}$. We then define
$\natCell_{\lambda t}$ as the composite
\[
\cl \XRA{\transExp{\zeta_{A,B}}} 
	\lambda {\left(\ext{h}{(t)} \circ 
			(\natTrans_{Z} \times \ext{h}A)\right)} 
	\XRA{\pushName^{-1}}
	\lambda{\left(\ext{h}{t}\right)} \circ \natTrans_{Z} 
	=
	\ext{h}{\left(\lambda t\right)} \circ \natTrans_Z
\]
The 2-cell $\zeta_{A,B}$ is defined in stages. 
First we set $\upsilon_{A,B}$ to be 
the following
composite, where we write $\iso$ for 
composites of $\phiTimes$ and structural isomorphisms:
\begin{td}
\eval_{\ext{h}{A}, \ext{h}{B}} \circ (\cl \times \ext{h}{A}) 
\arrow[swap]{d}{\iso} \\

\left(\eval_{\ext{h}{A}, \ext{h}{B}} 
	\circ 
	\big((\exp{\psinv{\natTrans}_{A}}{\natTrans_{B}}) 
		\times \ext{h}{A}\big)\right) 
	\circ 
	\left(\left(\evBar^F_{A, B} \circ F(\lambda t)\right) \times \ext{h}{A} 
	\right)
\arrow[swap]{d}
	{\epsilonExp_{\natTrans \circ \eval \circ (\Id \times \psinv{\natTrans})} 
		\circ 
	(\evBar^F_{A, B} F(\lambda t) \times \ext{h}{A})} \\
	
\left( \left(\natTrans_{B} 
	\circ \eval_{FA, FB}\right) 
	\circ (\Id_{(\scriptsizeexpobj{FA}{FB})} 
		\times \psinv{\natTrans}_{A})\right) 
	\circ \left(\left(\evBar^F_{A, B} \circ F(\lambda t)\right) 
		\times \ext{h}{A} 
	\right) 
\arrow[swap]{d}{\iso} \\

\left(\natTrans_{B} 
	\circ \left(\eval_{FA, FB} 
	\circ \big(\evBar^F_{A, B} \times F(A)\big)\right)\right)
	\circ \big( F(\lambda t) \times \psinv{\natTrans}_{A} \big) 
\arrow[swap]{d}
	{\natTrans_{B} \circ \epsilonExp_{(F(\eval) \circ \prodPres)} \circ 
		( F(\lambda t) \times \psinv{\natTrans} )} 
		\\
	
\left(\natTrans_{B} 
	\circ \left(F(\eval_{A,B}) 
	\circ \prodPres_{\scriptsizeexpobj{A}{B}, A}\right)\right) 
	\circ \left( F(\lambda t) \times \psinv{\natTrans}_{A} \right)
\end{td}
Next we define $\theta_{A,B}$ to be the composite
\begin{td}[column sep = 6em]
F(\eval_{A,B}) 
	\circ \left(\prodPres_{\scriptsizeexpobj{A}{B}, A} 
	\circ \left( F\lambda t \times FA \right)\right)
\arrow{r}{\theta_{A,B}}
\arrow[swap]{d}
	{F(\eval) \circ \prodPres \circ 
		( F(\lambda t) \times \psi^F_{A} )} &
Ft \circ \prodPres_{Z,A} \\

F(\eval_{A,B}) 
	\circ \left(\prodPres_{\scriptsizeexpobj{A}{B}, A} 
	\circ \big( \lambda t \times F\Id_{A} \big) \right)
\arrow[swap]{d}{F(\eval) \circ \nat} &
F{\left(\eval_{A,B} \circ (\lambda t \times A) \right)} \circ 
	\prodPres_{Z,A} 
\arrow[swap]{u}{F(\epsilonExp_{t}) \circ \prodPres} \\
	
F(\eval_{A,B}) 
	\circ \left(F{\left( \lambda t \times A \right)} 
	\circ \prodPres_{Z,A} \right)
\arrow[swap]{r}{\iso} &
\left(F(\eval_{A,B}) 
	\circ F{\left( \lambda t \times A \right)}\right)
	\circ \prodPres_{Z,A} 
\arrow[swap]{u}[yshift=-0mm]
	{\phi^F_{(\eval, \lambda t \times A)} 
		\circ 	\prodPres}
\end{td}
We can now define $\zeta_{A,B}$ as follows:
\begin{td}[column sep = 1em]
\eval_{\ext{h}{A}, \ext{h}{B}} \circ (\cl \times \ext{h}{A})  
\arrow[swap]{d}{\upsilon_{A,B}}
\arrow{r}{\zeta_{A,B}} &
\ext{h}{(t)} \circ 
	\big(\natTrans_{Z} \times A\big) \\

\left(\natTrans_{B} 
	\circ \left(F\eval_{A,B}
	\circ \prodPres_{\scriptsizeexpobj{A}{B}, A}\right)\right) 
	\circ \left( F(\lambda t) \times \psinv{\natTrans}_{A} \right)
\arrow[swap]{d}{\iso} &
\: \\

\left(\natTrans_{B} 
	\circ \left(F\eval_{A,B} 
	\circ \left(\prodPres_{\scriptsizeexpobj{A}{B}, A} 
	\circ \left( F(\lambda t) \times FA\right)\right)\right)\right) 
	\circ  \left(FZ \times \psinv{\natTrans}_{A} \right) 
\arrow[swap]{d}{\natTrans_{B} \circ \theta_{A,B} \circ 
	(FZ \times \psinv{\natTrans}_{A})} &
\: \\

\left(\natTrans_{B} 
	\circ \left(Ft 
	\circ \prodPres_{Z,A}\right)\right)
	\circ \left(FZ \times \psinv{\natTrans}_{A} \right) 
\arrow[swap]{d}{\iso} &
\: \\

\left(\natTrans_B
	\circ Ft\right)
	\circ \left( \prodPres_{Z,A}
	\circ \left(FZ \times \psinv{\natTrans}_{A} \right)\right)
\arrow[swap]{d}
	{\natCell_{t} \circ \prodPres \circ 
		(FZ \times \psinv{\natTrans}_{A} ) } &
\: \\

\left(\ext{h}{(t)} 
	\circ \left((\natTrans_{Z} \times \natTrans_{A}) 
	\circ \seq{F\pi_1, F\pi_2} \right)\right)
	\circ \left(\prodPres_{Z,A}  
	\circ \left(FZ \times \psinv{\natTrans}_{A} \right)\right) 
\arrow[swap]{d}{\iso} &
\: \\

\left(\left(\ext{h}{(t)} 
	\circ (\natTrans_{Z} \times \natTrans_{A})\right) 
	\circ \left(\seq{F\pi_1, F\pi_2} 
	\circ \prodPres_{Z,A}\right)\right)
	\circ \left(FZ \times \psinv{\natTrans}_{A} \right)
\arrow[swap]{d}
	{\ext{h}{(t)} \circ (\natTrans_{Z} \times \natTrans_{A}) 
	\circ (\unTimes_{Z,A})^{-1} \circ 
		(FZ \times \psinv{\natTrans}_{A})} &
\: \\

\ext{h}{(t)} \circ (\natTrans_{Z} \times \natTrans_{A}) \circ 
	\Id_{FZ \times FA} \circ 
	\big(FZ \times \psinv{\natTrans}_{A} \big) 
\arrow[swap]{dr}{\iso} &
\: \\

\: &
\ext{h}{(t)} \circ 
	(\natTrans_{Z} \times \natTrans_{A}\psinv{\natTrans}_{A})
\arrow[swap]{uuuuuuuu}
	{\ext{h}{(t)} \circ 
		(\natTrans_{Z} \times \co_{A})}
\end{td}
This completes the definition of $\natCell_{\lambda t}$. The only  
remaining case is horizontal composition. 

\vspace{1em}
\subproof{\rulename{hcomp} case.} As was the case for identities, 
the definition for multimaps of the form 
$t \circ u :  Z \to  B$ 
is forced by the axioms of a pseudonatural transformation. Using 
that $\ext{h}{}$ is a strict pseudofunctor, we define	
\begin{td}[column sep = 4em]
\natTrans_{B} \circ F(t \circ u)
\arrow[swap]{d}{\natTrans_{B} \circ (\phi^F_{t, u})^{-1}} 
\arrow{rr}{\natCell_{t \circ u}} &
\: &
\left(\ext{h}{(t)} 
	\circ \ext{h}{(u)}\right) 
	\circ \natTrans_{Z} \\

\natTrans_{B} \circ \left(F(t) \circ F(u)\right)
\arrow[swap]{d}{\iso} &
\: &
\ext{h}{(t)} \circ \left(\ext{h}(u) \circ \natTrans_{Z}\right)
\arrow[swap]{u}{\iso} \\

\left( \natTrans_B \circ Ft\right) \circ Fu
\arrow[swap]{r}{\natCell_t \circ F(u)} &
\left(\ext{h}{(t)} \circ \natTrans_{A}\right) \circ Fu
\arrow[swap]{r}{\iso} &
\ext{h}{(t)} \circ \left(\natTrans_{A} \circ Fu\right)
\arrow[swap]{u}{\ext{h}(t) \circ \natCell_u}
\end{td}

To show that $(\natTrans, \natCell)$ is indeed a pseudonatural transformation, 
we need to check the naturality condition and two axioms. Naturality is a 
straightforward check for each case outlined above. The two 
axioms---corresponding to the identity and \rulename{hcomp} cases---hold 
by construction.
\end{proof} 
\end{mylemma}

Examining the construction of the pseudonatural transformation just given, one 
extracts the following result.

\begin{mycor} \label{cor:termcatccc-unique-up-to-equivalence}
For any unary $\langCartClosed$-signature $\sig = (\baseTypes, \graph)$, 
cc-bicategory $\ccBicat{\baseCat}$, 
$\langCartClosed$-signature homomorphism
$h : \sig \to \altCat$,
and cc-pseudofunctor 
$(F, \prodPres, \expPres)$ such that 
\begin{enumerate}
\item Diagram~(\ref{eq:termcatCom-free-up-to-equivalence}) commutes,~\ie:
\begin{td}
\freeCartClosedBicat{\sig} 
\arrow{r}{F} &
\altCat \\
\sig
\arrow[hookrightarrow]{u}
\arrow[swap]{ur}{h} &
\:
\end{td}
\item For every $A_1, \,\dots\, , A_n, A, B \in \freeCartClosedBicat{\sig}$, the 1-cells
$\seq{F\pi_1, \,\dots\, , F\pi_n}$ and 
$\evBar_{A,B}$ are isomorphic to the 
identity,
\end{enumerate}
there exists an equivalence $F \simeq \ext{h}{}$ between $F$ and the 
canonical cc-pseudofunctor extending $h$.
\begin{proof}
One only needs to extend the pseudonatural equivalence $(\natTrans, \natCell)$ 
constructed in the proof 
of Lemma~\ref{lem:syntactic-model-uniqueness-up-to-equivalence} to cover 
constants. For these, one employs the second hypothesis. For any constant
$c \in \graph(A, B)$, condition~(1) requires that 
$F(c) = h(c) = \ext{h}(c)$. Condition~(2), on the other hand, entails that the 
components of
$(\natTrans, \natCell)$ are, inductively, each isomorphic to the identity. For 
the 2-cell filling
\begin{td}
FA 
\arrow[phantom]{dr}[description]{\twocell{\natCell_c}}
\arrow{r}{Fc} 
\arrow[swap]{d}{\natTrans_A} &
FB
\arrow{d}{\natTrans_B} \\

\ext{h}(A) 
\arrow[swap]{r}{\ext{h}(c)} &
\ext{h}(B)
\end{td}
one may therefore take the composite
$
\natTrans_B \circ Fc 
\XRA{\iso} 
Fc 
= \ext{h}(c) 
\XRA{\iso}
\ext{h}(c) \circ \natTrans_A
$
This definition is natural in $c$, and the two axioms of a 
pseudonatural transformation continue to hold. The claim follows. 
\end{proof}
\end{mycor}

\section{The type theory \texorpdfstring{$\langCartClosed$}{for cartesian 
closed bicategories}}

Fix a $\langCartClosed$-signature $\sig$. The type theory 
$\langCartClosed(\sig)$ is constructed as the internal language of 
$\freeCartClosedBiclone{\sig}$, with rules matching those of 
Construction~\ref{constr:free-cc-biclone}. These are collected together in 
Figures~\ref{r:ccc:terms}--\ref{r:ccc:invertibility}. Recall that for a context 
renaming $r$ we write 
$\hcomp{t}{r}$ to denote the term $\hcomp{t}{x_i \mapsto 
r(x_i)}$~(Figure~\ref{fig:context-renaming}), and that we write 
 $\mathrm{inc}_x$ for the inclusion of contexts
 $\Gamma \hookrightarrow \Gamma, x : A$ extending $\Gamma$ with a fresh 
 variable $x$.

The lambda abstraction operation extends to a (functorial) mapping on
rewrites, and the unit is derived as the mediating map corresponding to the 
identity (\cf~the discussion following~Definition~\ref{def:cc-bicat}). 

\begin{mydefn} \quad
\begin{enumerate}
\item For any derivable rewrite 
\mbox{$(\Gamma, x : A \vdash \tau : \rewrite{t}{t'} :B)$} we define 
$\lam{x}{\tau} : 
\rewrite{\lam{x}{t}}{\lam{x}{t'}}$ to be the rewrite $\transExp{ x \bind \tau 
\vert \epsilonExpRewr{t}}$ in context $\Gamma$.
\item For any derivable term $(\Gamma \vdash u : \exptype{A}{B})$ we define the 
unit $\etaExp{u} : \rewrite{u}{\lam{x}{\genevalterm{\wkn{u}{x}, x}}}$ to be the 
rewrite $\transExp{x \bind \id_{\hcomp{\evalterm}{\wkn{u}{x}, x}}}$ in context 
$\Gamma$. \qedhere
\end{enumerate}
\end{mydefn}

The usual application operation becomes a derived rule:
\begin{center}
\begin{prooftree}
\AxiomC{$\Gamma \vdash t : \exptype{A}{B}$}
\AxiomC{$\Gamma \vdash u : A$}
\BinaryInfC{$\Gamma \vdash \hcomp{\evalterm}{t,u} : B$}
\end{prooftree}
\end{center}	
The \mbox{$\epsilonExpRewr{}$-introduction} rule only relates  
lambda abstractions and variables, but the general form of (explicit) 
$\beta$-reduction is derivable. In the definition we use the following 
notation. For a context 
$\Gamma := (x_i : A_i)_{i=1, \dots, n}$ 
and terms 
$\Gamma, x : A \vdash t : B$ 
and
$\Gamma \vdash u : A$,
we write
$\hcomp{t}{\id_\Gamma, x \mapsto u}$
to denote the term
$\hcomp{t}{x_1 \mapsto x_1, \, \dots \,, x_n \mapsto x_n, x \mapsto u}$
in context $\Gamma$.

\begin{mydefn} \label{def:beta-rewrite}
For derivable terms \mbox{$\Gamma, x : A \vdash t : B$} and \mbox{$\Gamma 
\vdash u : A$} we define the \Def{$\beta$-reduction rewrite} 
$\genEpsilonExp{x \smallbind t, u} 
		:  
		\rewrite
			{\hcomp{\evalterm}{\lam{x}{t}, u}}
			{\hcomp{t}{\id_{\Gamma}, x 	\mapsto u}}$ to be $ 
\hcomp{\epsilonExpRewr{t}}{\id_{\Gamma}, x \mapsto u} \vert \tau$ in context 
$\Gamma$, where $\tau$ is the following composite of structural isomorphisms:
\begin{align*}
\hcomp{\evalterm}{\lam{x}{t}, u} &\iso 
\hcomp{\evalterm}{\hcomp{(\lam{x}{t})}{\mathrm{inc}_x}, u} \\
&\iso 
\hcompbig
		{\evalterm}
		{\hcomp
			{(\lam{x}{t})}
			{\hcomp{\mathrm{inc}_x}{\id_\Gamma, x \mapsto u}} , u
		}  \\
&\iso 
\hcompbig{\evalterm}
	{\hcomp
		{
			\hcomp{(\lam{x}{t})}{\mathrm{inc}_x}
		}
		{
			{\id_\Gamma, x \mapsto u}
		}, 
	 \hcomp{x}{{\id_\Gamma, x \mapsto u}}
	} \\
&\iso 
	\hcompthree
		{\evalterm}
		{\hcomp{(\lam{x}{t})}{\mathrm{inc}_x}, x}
		{\id_\Gamma, x \mapsto u} 
\qedhere
\end{align*}
\end{mydefn}

In a similar vein, one may wish to introduce the counit via the 
following more explicit rule: 
\begin{equation*} 
\unaryRule
		{\Gamma, x : A \vdash t : B}
		{\Gamma, y : A \vdash 
			\epsilonExpRewr{x \bind t} : 
				\rewrite
					{\hcomp{\evalterm}{\wkn{(\lam{x}{t})}{y},y}}
					{\hcomp{t}{\id_\Gamma, x \mapsto y}} : B}
		{}
\vspace{-\treeskip}
\end{equation*}
In the presence of the structural rewrites, this definition is equivalent to 
that given in Figure~\ref{r:ccc:rewrites}.

We continue to work up to $\alpha$-equivalence of terms and rewrites. Unlike 
the extension from $\langBiclone$ to $\langCart$, the type theory 
$\langCartClosed$ has new
\Def{binding operations}: alongside the usual binding rules for lambda 
abstraction, we require that the variable $x$ is bound in the rewrite 
$\transExp{x \bind \alpha}$. This is reflected in the definition of 
$\alpha$-equivalence.

\setlength{\floatsep}{5pt plus 1.0pt minus 2.0pt} 

\begin{figure*}[!h]
{\small
\begin{minipage}{\textwidth}
\begin{mdframed}
\centering
\input{rules/ccc/exponential-terms}
\vspace{-\treeskip}
\caption{Terms for cartesian closed structure \label{r:ccc:terms}}
\end{mdframed}
\end{minipage}

\begin{minipage}{\textwidth}
\begin{mdframed}
\centering
\input{rules/ccc/exponential-rewrites}
\caption{Rewrites for cartesian closed structure \label{r:ccc:rewrites}}
\end{mdframed}
\end{minipage}

\begin{minipage}{\textwidth}
\begin{mdframed}
\centering
\input{rules/ccc/trans-ump}
\input{rules/ccc/extra-congruence}
\caption{Universal property and congruence laws for  $\transExp{\alpha}$\label{r:ccc:umpTrans}}
\end{mdframed}
\end{minipage}

\begin{minipage}{\textwidth}
\begin{mdframed}
\centering
\input{rules/ccc/invertibility-intro}

\input{rules/ccc/invertibility-congruences}

\caption{Inverses for the unit and counit\label{r:ccc:invertibility}}
\end{mdframed}
\end{minipage}

}

\begin{manyfigcap}
Rules for $\langCartClosed(\sig)$.
\end{manyfigcap}

\end{figure*}

\clearpage

\paragraph*{$\alpha$-equivalence and free variables}

For 
$\lambda$-abstraction we follow the usual conventions of the simply-typed 
lambda calculus 
(\cf~\cite{Barendregt1985}).

\begin{mydefn}
For any $\langCartClosed$-signature $\sig$ define the 
\Def{$\alpha$-equivalence relation} 
$\aeq$ on terms by extending Definition~\ref{def:fp:alpha-equivalence} with the 
rules
\begin{center}
\binaryRule
	{t[y /x] \aeq t'[y /x']}
	{y \text{ fresh }}
	{\lam{x}{t} \aeq \lam{x'}{t'}}
	{} \vspace{-\treeskip} \quad
\unaryRule
	{t \aeq t'}
	{\epsilonExpRewr{t} \aeq \epsilonExpRewr{t'}}
	{} \vspace{-\treeskip} \quad
\unaryRule
	{\sigma[y / x] \aeq \sigma[y / x'] \qquad y \text{ fresh }}
	{\transExp{x \bind \sigma} \aeq \transExp{x' \bind \sigma}}
	{} 
	\vspace{-\treeskip}
\end{center}
Similarly, the meta-operation of capture-avoiding substitution is that 
of~Definition~\ref{def:fp:alpha-equivalence}, 
extended by the rules 
\begin{center}
$\evalterm(f,x)[t / f, u / x] := \heval{t, u}$ \quad and \quad 
$(\lam{x}{t})[u_i / x_i] := \lam{z}{(t[z/ x, u_i/x_i])}$ for $z$ fresh
\end{center}
and
\begin{center}
$\epsilonExpRewr{t}[u_i/x_i] := \epsilonExpRewr{t[u_i/x_i]}$ \:\:\:\: and \:\:\:\:
$\transExp{y \dot \alpha}[u_i / x_i] := \transExp{z \dot \alpha[z / y, u_i / x_i]}$ for $z$ fresh
\end{center}
These rules extend to the inverses of rewrites in the obvious fashion.
\end{mydefn}

\begin{prooflesslemma} 
Let $\sig$ be a $\langCartClosed$-signature. Then in $\langCartClosed(\sig)$:
\begin{enumerate} 
\item If $\Gamma \vdash t : B$ and $t =_\alpha t'$ then $\Gamma \vdash t' :B$,
\item If $\Gamma \vdash \tau : \rewrite{t}{t'} : B$ and $\tau =_\alpha \tau'$ then 
$\Gamma \vdash \tau' : \rewrite{t}{t'} : B$. \qedhere
\end{enumerate} 
\end{prooflesslemma} 

The $\aeq$ relation is a congruence on the derived structure. In particular, 
one obtains the expected equality for the 
induced lambda abstraction operation on rewrites. 

\begin{prooflesslemma}
Let $\sig$ be a $\langCartClosed$-signature. Then in $\langCartClosed(\sig)$:
\begin{enumerate}
\item If $\tau[y / x] \aeq \tau'[y/x']$ (for $y$ fresh) then 
$\lam{x}{\tau} 
\aeq 
\lam{x'}{\tau'}$, 
\item If $u \aeq u'$ then $\etaExp{u} \aeq \etaExp{u'}$,
\item If $t[y / x] \aeq t'[y/x']$ and $u \aeq u'$ then  
$\beta_{x \smallbind t,\, u} \aeq \beta_{x' \smallbind t', \, u'}$. \qedhere
\end{enumerate}
\end{prooflesslemma} 

As for $\langCart$, the type theory 
$\langCartClosed$ satisfies all the expected type-theoretic well-formedness 
properties.

\begin{mydefn} 
Fix a $\langCartClosed$-signature $\sig$. We define the \Def{free variables in 
a term} $t$ in 
$\langCartClosed(\sig)$ by extending Definition~\ref{def:fp:free-vars} by setting 
$\fv{\left(\lam{x}{t}\right)} := \fv(t) - \{ x \}$ and
$\fv{\left(\heval{p}\right)} := \{ p \}$.
Similarly, we define the \Def{free variables in a rewrite $\tau$} in 
$\langCartClosed(\sig)$ by extending Definition~\ref{def:fp:free-vars} as 
follows:
$\fv(\epsilonExpRewr{t}) = \fv(t)$ and 
$\fv{\left(\transExp{x \bind \alpha}\right)} = \fv(\alpha) - \{ x\}$.
We define the free variables of a specified inverse $\sigma^{-1}$ to be exactly 
the free variables of $\sigma$. An occurrence of a variable in a term or 
rewrite is \Def{bound} if it is not free. 
\end{mydefn} 

\begin{prooflesslemma}
Let $\sig$ be a $\langCartClosed$-signature. For any derivable judgements 
$\Gamma \vdash u : 
B$ and \mbox{$\Gamma \vdash \tau : \rewrite{t}{t'} : B$} in 
$\langCartClosed(\sig)$, 
\begin{enumerate} 
\item $\fv(u) \subseteq \dom(\Gamma)$, 
\item $\fv(\tau) \subseteq \dom(\Gamma)$, 
\item The judgements $\Gamma \vdash t : B$ and $\Gamma \vdash t' : B$ are both derivable. 
\end{enumerate} 
Moreover, whenever $(\Delta \vdash u_i : A_i)_{i= 1, \,\dots\, , n}$ and 
$\Gamma := (x_i  : A_i)_{i = 1, \,\dots\, , n}$, then
\begin{enumerate}
\item If $\Gamma \vdash t : B$, then $\Delta \vdash t[u_i / x_i] : B$, 
\item If $\Gamma \vdash \tau : \rewrite{t}{t'} : B$, then 
$\Delta \vdash \tau[u_i / x_i] : \rewrite{t[u_i / x_i]}{t'[u_i / x_i]} : B$. \qedhere
\end{enumerate}
\end{prooflesslemma}

\subsection{\texorpdfstring{The syntactic model of $\langCartClosed$}{The syntactic model}}

We now turn to constructing the syntactic model for $\langCartClosed(\sig)$ and 
proving it is the free cartesian closed biclone on $\sig$. The 
construction is a straightforward extension of 
Construction~\ref{constr:langcart-syntactic-model} (page~\pageref{constr:langcart-syntactic-model}).

\begin{myconstr} \label{constr:langcartclosed-syntactic-model}
For any $\langCartClosed$-signature $\sig = (\baseTypes, \graph)$, define the 
\Def{syntactic model} 
$\syncloneCartClosed{\sig}$ of $\langCartClosed(\sig)$ as follows. The 
sorts are nodes $A, B, \dots$ of $\graph$. The 1-cells are
$\alpha$-equivalence classes of terms 
$(x_1 : A_1, \,\dots\, , x_n : A_n \vdash t : B)$ 
derivable in 
$\langCartClosed(\sig)$. We 
assume a fixed enumeration $x_1, x_2, \dots$ of variables, and that the 
variable name in the $i$th position is determined by this enumeration. 
The 2-cells are 
$\alpha{\equiv}$-equivalence classes of rewrites 
\mbox{$(x_1 : A_1, \,\dots\, , x_n : A_n \vdash \tau : \rewrite{t}{t'} : B)$}. 
Composition is vertical composition 
and and the identity on $t$ is $\id_t$; the substitution operation is explicit 
substitution and the structural rewrites are $\assoc{}, \subid{}$ and 
$\indproj{i}{}$. 
\end{myconstr}

$\syncloneCartClosed{\sig}$ is a cartesian closed biclone. Products are as in 
$\syncloneCart{\sig}$ (Section~\ref{sec:fp:syntactic-model}) and for 
exponentials the biuniversal arrow is
$\evalterm(f,x) : (f : (\exptype{A}{B}), x : A) \to (y : B)$. 
Indeed, 
for any judgement 
$(\Gamma, x : A \vdash \alpha : 
		\rewrite{\hcomp{\evalterm}{\wkn{u}{x}, x}}{t} : B)$ 
in $\langCartClosed(\sig)$, the rewrite $\transExp{x \bind \alpha}$ is the 
unique $\gamma$ (modulo $\alpha{\equiv}$) such that
\begin{equation} \label{eq:exponentials-ump-with-gamma}
\Gamma, x : A \vdash \alpha \equiv \epsilonExpRewr{t} \vert 
\genevalterm{\wkn{\gamma}{x}, x} : 
	\rewrite{\hcomp{\evalterm}{\wkn{u}{x}, x}}{t} : B
\end{equation}
Existence is precisely rule~\rulename{U1}. For uniqueness, for 
any $\gamma$ 
satisfying~(\ref{eq:exponentials-ump-with-gamma}) one has 
\[
\gamma \overset{\rulename{U2}}{\equiv} \transExp{x \bind \epsilonExpRewr{t} 
\vert 
\genevalterm{\wkn{\gamma}{x}, x}}  \overset{\rulename{cong}}{\equiv} 
\transExp{x \bind \alpha}
\]
Moreover, $\syncloneCartClosed{\sig}$ is the free cartesian closed biclone on 
$\sig$, which validates our claim that $\langCartClosed(\sig)$ is the internal 
language of $\freeCartClosedBiclone{\sig}$.

\begin{mypropn} \label{propn:synclonecartclosed-free-property}
For any $\langCartClosed$-signature $\sig$, cartesian closed biclone 
$\ccClone{T}{\altBiclone}$, and $\langCartClosed$-signature homomorphism 
$h : \sig \to \altBiclone$, there 
exists a unique strict cartesian closed pseudofunctor 
$h\sem{-} : \syncloneCartClosed{\sig} \to \altBiclone$  such that 
$h\sem{-} \circ \iota = h$, for 
$\iota : \sig \hookrightarrow \syncloneCartClosed{\sig}$ 
the inclusion.
\begin{proof}
We extend the pseudofunctor $h\sem{-}$ of 
Proposition~\ref{propn:synclonecart-free-property} (page~\pageref{propn:synclonecart-free-property}) 
with the following rules.
\begin{align*} 
h\sem{\exptype{A}{B}} &:= \expobj{h\sem{A}}{h\sem{B}}  \\[5pt]
h\sem{f : \exptype{A}{B}, a : A \vdash \evalterm(f,a) : B} &:= \eval_{A,B} \\
h\sem{\Gamma \vdash \lam{x}{t} : \exptype{A}{B}} &:= 
	\lambda{\left(h\sem{\Gamma, x : A \vdash t : B}\right)} \\[5pt]
h\sem{\Gamma, x : A \vdash \epsilonExpRewr{t} : 
\rewrite{\genevalterm{\wkn{(\lam{x}{t})}{x}, x}}{t} : B} &:= 
\epsilon_{h\sem{\Gamma, x : A \vdash t : B}} \\
h\sem{\Gamma \vdash \transExp{x \bind \alpha} : \rewrite{u}{\lam{x}{t}} : 
\exptype{A}{B}} &:= 
	\transExplr{h\sem{\Gamma, x : A \vdash \alpha : 
		\rewrite{\genevalterm{\wkn{u}{x}, x}}{t} : B}} 
\end{align*} 
Uniqueness follows because any strict cc-pseudofunctor must strictly preserve 
the $\lambda(-)$ and $\transExp{-}$ operations 
(\cf~Lemma~\ref{lem:free-cc-biclone} and 
Lemma~\ref{lem:strict-preservation-strict-pres-UMP}). 
\end{proof}
\end{mypropn}

\begin{myremark}
As we saw for 
products~(Remark~\ref{rem:fp:biuniversal-universal-interleaving}), the 
universal property of the counit for exponentials 
gives rise to a nesting 
of (global) biuniversal arrows and (local) universal arrows. These are related 
by the following bijective correspondence, in which we write $(x : A)$ to 
indicate the variable $x$ of type $A$ is free in 
the context~(\cf~\cite{MartinLof}):
\begin{center}
\begin{bprooftree} 
\AxiomC{$(x : A)$}
\noLine
\UnaryInfC{$\hcomp{\evalterm}{\wkn{u}{x}, x} \To t : B$} 
\doubleLine 
\UnaryInfC{$u \To \lam{x}{t} : \exptype{A}{B}$} 
\end{bprooftree} 
\end{center}
We conjecture that 
a calculus for cartesian 
closed \emph{tri}categories (cartesian closed $\infty$-categories) would have 
three (a countably infinite tower) of such correspondences. 
\end{myremark}

For a unary $\langCartClosed$-signature $\sig$, the nucleus 
$\nucleus{\syncloneCartClosed\sig}$ of 
$\syncloneCartClosed\sig$ is cartesian closed with exponentials as described in
Remark~\ref{rem:bicat-exp-structure-not-preserved}. We make this explicit in 
the next construction, which mirrors the syntactic model of the simply-typed lambda 
calculus~(\eg~\cite[Chapter~4]{Crole1994}). 

\begin{myconstr}  \label{constr:restrict-of-langcartclosed-biclone}
For any $\langCartClosed$-signature $\sig$, define a bicategory 
$\increstrict{\syncloneCartClosed{\sig}}$ as follows. The objects are unary 
contexts 
with a single \emph{fixed} variable name. The 1-cells $(x : A) \to (x : B)$ 
are 
\mbox{$\alpha$-equivalence} classes of terms $(x : A \vdash t :B)$ derivable 
in 
$\langCartClosed(\sig)$. The \mbox{2-cells} are 
\mbox{$\alpha{\equiv}$-equivalence} classes of rewrites 
\mbox{$(x : A \vdash \tau : \rewrite{t}{t'} : B)$}. Vertical composition is 
given by the $\vert$ operation, and horizontal composition is given by explicit 
substitution. 
\end{myconstr} 

As we have seen, we 
cannot hope for 
$\increstrict{\syncloneCartClosed{\sig}}$ to satisfy a strict universal 
property~(recall the discussion following Lemma~\ref{lem:free-cc-biclone} 
on page~\pageref{lem:free-cc-biclone}, as well as Example~\ref{exmp:nucleus-of-fp-not-free} on 
page~\pageref{exmp:nucleus-of-fp-not-free}). 
Nonetheless, we shall see in Section~\ref{sec:syntactic-model-free-property} 
that it is \Def{weakly initial}: any morphism of 
$\langCartClosed$-signatures may 
be extended to a pseudofunctor out of 
$\increstrict{\syncloneCartClosed{\sig}}$, but this may not be unique. Hence, 
$\langCartClosed$ may be soundly 
interpreted in any cc-bicategory. We shall also see that 
$\increstrict{\syncloneCartClosed{\sig}}$ is biequivalent to the free
cc-bicategory $\freeCartClosedBicat{\sig}$ on $\sig$, yielding a bicategorical 
universal property. 
Before proceeding to these results, we 
first establish a series of lemmas that will simplify their proofs. 

\subsection{\texorpdfstring{Reasoning within $\langCartClosed$}{Reasoning with type-theoretic exponentials}}

We begin by recovering the unit-counit presentation of 
exponentials (\cf~\cite{Seely1987, Hilken1996}) as a series of admissible 
rules. These are 
collected together in  
Figure~\ref{fig:ccc:admissible-rules}, below. The proofs are similar to the 
case for 
products, so we omit them. 

\begin{prooflesslemma}
For any $\langCartClosed$-signature $\sig$, the rules of 
Figure~\ref{fig:ccc:admissible-rules} 
are admissible in $\langCartClosed(\sig)$. 
\end{prooflesslemma}

A direct corollary is that the $\beta$-reduction rewrite of 
Definition~\ref{def:beta-rewrite} is natural. 

\begin{prooflesscor}
For any $\langCartClosed$-signature $\sig$,
if the judgements
$(\Gamma, x : A \vdash \tau : \rewrite{t}{t'} : B)$ and 
$(\Gamma \vdash \sigma : \rewrite{u}{u'} : A)$ are derivable in
$\langCartClosed(\sig)$, 
then the following diagram of rewrites commutes:
\begin{td}[column sep = 5em]
\heval{\lam{x}{t}, u} 
\arrow[Rightarrow]{r}{\heval{\lam{x}{\tau}, \sigma}} 
\arrow[swap, Rightarrow]{d}{\genEpsilonExp{x \smallbind t, u}} &
\heval{\lam{x}{t'}, u'} 
\arrow[Rightarrow]{d}{\genEpsilonExp{x \smallbind t', u'}} \\
\hcomp{t}{\id_\Gamma, x \mapsto u}  
\arrow[swap, Rightarrow]{r}{\hcomp{\tau}{\id_\Gamma, x \mapsto \sigma}} &
\hcomp{t'}{\id_\Gamma, x \mapsto u'}
\end{td}
\end{prooflesscor}

\setlength{\floatsep}{5pt plus 1.0pt minus 2.0pt} 

\begin{figure*}[h]
\begin{minipage}{\textwidth}
{
\begin{mdframed}
\begin{minipage}{\textwidth}
\centering 


\unaryRule {\Gamma, x : A \vdash t : B} 
{\Gamma \vdash \lam{x}{\id_t} \equiv \id_{\lam{x}{t}} : \rewrite{\lam{x}{t}}{\lam{x}{t}} : \exptype{A}{B}} 
{} 
\binaryRule {\Gamma, x : A \vdash \tau' : \rewrite{t'}{t''} : B} 
{\Gamma, x : A \vdash \tau : \rewrite{t}{t'} : B} 
{\Gamma \vdash \lam{x}{(\tau' \vert \tau)} \equiv (\lam{x}{\tau'}) \vert (\lam{x}{\tau}) : \rewrite{\lam{x}{t}}{\lam{x}{t''}} : \exptype{A}{B}} 
{} 
\end{minipage}

%

\begin{minipage}{\textwidth}
\centering
\unaryRule {\Gamma \vdash \sigma : \rewrite{u}{u'} : \exptype{A}{B}}
{\Gamma \vdash \etaExp{u'} \vert \sigma \equiv \lam{x}{\genevalterm{\wkn{\sigma}{x}, x}} \vert \etaExp{u} : \rewrite{u}{\lam{x}{\hcomp{\evalterm}{\wkn{u'}{x}, x}}} : \exptype{A}{B}} 
{$\etaExp{}$-nat} 

\unaryRule {\Gamma, x : A \vdash \tau : \rewrite{t}{t'} : B} 
{\Gamma, x : A \vdash \tau \vert \epsilonExpRewr{t} \equiv \epsilonExpRewr{t'} \vert \genevalterm{\wkn{(\lam{x}{\tau})}{x}, x} : \rewrite{ \genevalterm{\wkn{(\lam{x}{t})}{x}, x}}{t'} : B} 
{$\epsilonExpRewr{}$-nat} 
\end{minipage}

\begin{minipage}{\textwidth}
\centering
\begin{mbox}
\centering 
\unaryRule {\Gamma, x : A \vdash t : B} 
{\Gamma \vdash (\lam{x}{\epsilonExpRewr{t}}) \vert \etaExp{t} \equiv \id_{\lam{x}{t}} : \rewrite{\lam{x}{t}}{\lam{x}{t}} : \exptype{A}{B}} 
{triangle-law-1} 

{\small
\begin{bprooftree} 
\AxiomC{$\Gamma \vdash u : \exptype{A}{B}$} 
\RightLabel{{\scriptsize triangle-law-2}}
\UnaryInfC{$\Gamma, x : A \vdash \epsilonExpRewr{\hcomp{\evalterm}{\wkn{u}{x}, 
x}} \vert \hcomp{\evalterm}{\wkn{\etaExp{u}}{x}, x} \equiv 
\id_{\hcomp{\evalterm}{\wkn{u}{x}, x}}$\hspace{22mm}}
\noLine
\UnaryInfC{\hspace{50mm}$: \rewrite{\hcomp{\evalterm}{\wkn{u}{x}, x}}{\hcomp{\evalterm}{\wkn{u}{x}, x}} : B$}
\end{bprooftree}\vspace{\treeskip} 
} 
\end{mbox}
\end{minipage}

{
\caption{Admissible rules for $\langCartClosed(\graph)$\label{fig:ccc:admissible-rules}}
}
\end{mdframed}
}
\end{minipage}
\end{figure*}

Recall that for products we constructed a rewrite $\postName$ of type 
\[
{\hcomp{\pair{t_1, \,\dots\, , t_m}}{u_1, \,\dots\, , u_n}} \To 
	\pair{\hcomp{t_1}{u_1, \,\dots\, , u_n}, \,\dots\, , \hcomp{t_m}{u_1, 
	\,\dots\, , u_n}}
\]
For exponentials we call the corresponding rewrite 
$\pushName$~(\cf~Construction~\ref{constr:push-semantically}).  
Just as $\postName$ witnesses that explicit substitutions and the tupling 
operation commute (up to isomorphism), so $\pushName$ witnesses 
that explicit substitutions and lambda abstractions can be permuted (up to 
isomorphism). 
Precisely, $\pushName$ 
relates the following two derivations 
(where $\Gamma := (x_i : A_i)_{i=1,\dots,n}$):
\begin{prooftree}
\AxiomC{$\Gamma, x : A \vdash t : B$}
\UnaryInfC{$\Gamma \vdash \lam{x}{t} : \exptype{A}{B}$}
\AxiomC{$(\Delta \vdash u_i : A_i)_{i=1,\dots,n}$}
\BinaryInfC{$\Delta \vdash \hcomp{(\lam{x}{t})}{x_i \mapsto u_i} : 
\exptype{A}{B}$}
\end{prooftree}
and
\begin{prooftree}
\AxiomC{$\Gamma, x : A \vdash t : B$}
\AxiomC{$(\Delta \vdash u_i : A_i)_{i=1,\dots,n}$}
\UnaryInfC{$(\Delta, x : A \vdash \wkn{u_i}{x} : A_i)_{i=1,\dots,n}$}
\AxiomC{$\Delta, x : A \vdash x : A$}
\TrinaryInfC{$\Delta, x : A \vdash 
	\hcomp{t}{x_i \mapsto \wkn{u_i}{x}, x \mapsto x} : B$}
\UnaryInfC{$\Delta \vdash 
	\lam{x}{\hcomp{t}{x_i \mapsto \wkn{u_i}{x}, x \mapsto x}} : \exptype{A}{B}$}
\end{prooftree}
From the perspective of the simply-typed lambda calculus, the rewrite 
\[
\pushName : 
	\hcomp{(\lam{x}{t})}{x_i \mapsto u_i} 
	\To 
	\lam{x}{\hcomp{t}{x_i \mapsto \wkn{u_i}{x}, x \mapsto x}}
\]	
is an explicit 	version of the usual rule 
$(\lam{x}{t})[u_i/x_i] = \lam{z}{t[u_i/x_i, z /x]}$ for the meta-operation of 
capture-avoiding substitution~(\cf~\cite[Definition~4]{Ritter1997}, where a similar 
operation is constructed for a version of the simply-typed lambda calculus 
with explicit substitution).

We construct $\pushName$ by emulating 
Construction~\ref{constr:push-semantically} within $\langCartClosed$.

\begin{myconstr}
For any $\langCartClosed$-signature $\sig$ we construct a rewrite 
$\push{t}{\ind{u}}$ in $\langCartClosed(\sig)$ making the following rule is 
admissible:
\begin{center}
\begin{prooftree}
\AxiomC{$\Gamma, x : A \vdash t : B$}
\AxiomC{$(\Delta \vdash u_i : A_i)_{i=1,\dots,n}$}
\BinaryInfC{$\Delta \vdash \push{t}{\ind{u}} 
	: 
	\rewrite
		{\hcomp{(\lam{x}{t})}{x_i \mapsto u_i}}
		{\lam{x}{
			\hcomp{t}
				{x_i \mapsto \wkn{u_i}{x}, x \mapsto x}}} 
	: 		
	\exptype{A}{B}$}
\end{prooftree}
\end{center}
Following Construction~\ref{constr:push-semantically}, we first need to 
construct the 2-cell $\phiCell{2}{}$ witnessing the pseudofunctorality of the 
product-former. From the judgements 
$\Gamma \vdash t : B$ and $(\Delta \vdash u_i : A_i)_{i=1, \,\dots\, , n}$ one 
obtains the terms
\begin{center}
$\hcompthree{t}{\mathrm{inc}_x}{x_i \mapsto \wkn{u_i}{x}, x \mapsto x}$ \quad 
and \quad $\wkn{\hcomp{t}{x_i \mapsto u_i}}{x}$
\end{center}
of type $B$ in context $\Delta, x : B$ by either performing explicit 
substitution or weakening first. These terms are related by the following 
composite, which we 
call $\phiCell{(2)}{t, \ind{u}}$:
\begin{align*}
\hcompthree{t}{\mathrm{inc}_x}{x_i \mapsto \wkn{u_i}{x}, x \mapsto x} 
&\overset{\assoc{}}{\iso} 
	\hcomp{t}{\hcomp{\mathrm{inc}_x}{x_i \mapsto \wkn{u_i}{x}, x \mapsto x}} \\
	&\overset{\hcomp{t}{\indproj{\bullet}{}}}{\iso}  \hcomp{t}{x_i \mapsto 
	\wkn{u_i}{x}}  \\
	&\overset{\assoc{}^{-1}}{\iso}  \wkn{\hcomp{t}{x_i \mapsto u_i}}{x}
\end{align*}

We therefore set $\push{t}{\ind{u}}$ to be $\transExp{x \bind \tau}$, for 
$\tau$ 
the composite
\begin{align*}
\evalterm \big\{ \wkn{\hcomp{(\lam{x}{t})}{x_i \mapsto u_i}}{x}, x \big\} 
\hspace{-35mm}& 
\\
	&\iso \evalterm 
	\big\{ 
		\hcompthree
			{(\lam{x}{t})}
			{\mathrm{inc}_x}
			{x_i \mapsto \wkn{u_i}{x}, x \mapsto x}, 
		\hcomp{x}{{x_i \mapsto \wkn{u_i}{x}, x \mapsto x}} 
	\big\} \\
	&\iso \evalterm \big\{  
	\hcomp{(\lam{x}{t})}{\mathrm{inc}_x}, x \big\} 
	\big\{ x_i \mapsto \wkn{u_i}{x}, 
			x \mapsto x  
	\big\}  \\ 
	&\iso \hcomp{t}{x_i \mapsto \wkn{u_i}{x}, x \mapsto x}
\end{align*}
where the first isomorphism is 
$\heval
	{(\phiCell{2}{\lam{x}{t}, \ind{x}})^{-1}, 
		\indproj{-(\len{\Delta}+1)}{\wkn{\ind{u}}{x}, x}}$, the second 
is $\assoc{}^{-1}$ and the third is
$\hcomp{\epsilonExpRewr{t}}{\wkn{u_i}{x}, x}$.
\end{myconstr}

Thinking of rewrites in $\langCartClosed$ as witnesses 
for equalities in the simply-typed lambda calculus, the following lemma is as 
expected 
(\cf~Lemma~\ref{lem:PseudoCCCCanonical2CellsLaws}).

\begin{mylemma} \label{lem:properties-of-push}
For any $\langCartClosed$-signature $\sig$, if 
$\Gamma := (x_i : A_i)_{i=1,\dots,n}$ and the judgement
$({\Delta \vdash \sigma_i : \rewrite{u_i}{u_i'} : A_i})$
is derivable in $\langCartClosed(\sig)$, then:
\begin{enumerate}
%
{\tikzcdset{arrow style=tikz, arrows={Rightarrow}}
\item (Naturality). If $\Gamma, x : A \vdash \tau : \rewrite{t}{t'} : B$, then
\begin{td}
\hcomp{(\lam{x}{t})}{\ind{u}} 
\arrow{r}{\pushName} 
\arrow[swap]{d}{\hcomp{(\lam{x}{\tau})}{\ind{\sigma}}} &
\lam{x}{\hcomp{t}{\wkn{\ind{u}}{x}, x}} 
\arrow{d}{\lam{x}{\hcomp{\tau}{\wkn{\ind{\sigma}}{x}, x}}} \\
\hcomp{(\lam{x}{t'})}{\ind{u}'} 
\arrow[swap]{r}{\pushName} &
\lam{x}{\hcomp{t'}{\wkn{\ind{u}'}{x}, x}} 
\end{td}

\item (Compatibility with $\subid{}$). \label{c:push-and-subid} If 
$\Gamma, x : A \vdash t : B$, then 
\begin{td}[column sep = 6em]
\lam{x}{t} \arrow{r}{\subid{}} 
\arrow[swap]{d}{\lam{x}{\subid{}}}  &
\hcomp{(\lam{x}{t})}{\ind{x}} 
\arrow{d}{\pushName} \\
\lam{x}{\hcomp{t}{\ind{x}}} &
\lam{x}{\hcomp{t}{\wkn{\ind{x}}{x}, x}} 
\arrow{l}{\lam{x}{\hcompsmall{t}{x, \indproj{\bullet}{}}}}
\end{td}

\newpage
\item (Compatibility with $\assoc{}$). \label{c:push-and-assoc} If $\Gamma, x : 
A \vdash t : C$, $\Delta := (y_j : B_j)_{j = 1,\dots, m}$ and \mbox{$(\Sigma 
\vdash v_j : B_j)_{j=1,\dots, m}$}, then
\small
\begin{td}[column sep = tiny]
\: &
\hcomp{\big(\lam{x}{\hcomp{t}{\wkn{\ind{u}}{x}, x}}\big)}{\ind{v}} 
\arrow{dr}{\pushName} &
\: \\

\hcompthree{(\lam{x}{t})}{\ind{u}}{\ind{v}} 
\arrow{ur}{\hcomp{\pushName}{\ind{v}}} \arrow[swap]{d}{\assoc{}}  &
\:  &
\lam{x}{\hcompthree{t}{\wkn{\ind{u}}{x}, x}{\wkn{\ind{v}}{x}, x}} 
\arrow{d}{\lam{x}{\assoc{}}} \\

\hcomp{(\lam{x}{t})}{\hcomp{\ind{u}}{\ind{v}}} \arrow[swap]{d}{\pushName} &
\: &
\lam{x}
	{\hcompbig{t}
	{\hcomp{\wkn{\ind{u}}{x}}{\wkn{\ind{v}}{x}, x}, 
	 \hcomp{x}{\wkn{\ind{v}}{x}, x}}
	} 
\arrow{d}{\lam{x}{\hcomp{t}{\assoc{}, \indproj{m+1}{}}}} \\

\lam{x}{\hcompbig{t}{\wkn{\hcomp{\ind{u}}{\ind{v}}}{x}, x}} 
\arrow[swap]{dr}[xshift=-1mm, yshift =-.5mm]{\assoc{}} &
\: &
\lam{x}{\hcompbig{t}{\hcomp{\ind{u}}{\hcomp{\ind{y}}{\wkn{\ind{v}}{x}, x}}, x}} 
\arrow{dl}
	{\lam{x}
		{
			\hcompsmall{t}{\hcompsmall{\ind{u}}{\indproj{\bullet}{}}, x}
		}
	} \\

\: &
\lam{x}{\hcompbig{t}{\hcomp{\ind{u}}{\wkn{\ind{v}}{x}}, x}}  &
\: 
\end{td}
\normalsize

\item (Compatibility with $\etaExp{}$). \label{c:push-with-etaexp} If $\Gamma, 
x : A \vdash t : B$ then 
\small
\begin{td}[column sep = 3em]
\hcomp{t}{\ind{u}} \arrow{r}{\hcomp{\etaExp{}}{\ind{u}}} 
\arrow[swap]{dd}{\etaExp{}} &
\hcomp{\big(\lam{x}{\heval{\wkn{t}{x}, x}}\big)}{\ind{u}} 
\arrow{d}{\pushName} 
\\

\: &
\lam{x}{\hcompthree{\evalterm}{\wkn{t}{x}, x}{\wkn{\ind{u}}{x}, x}} 
\arrow{d}{\lam{x}{\assoc{}}} \\

\lam{x}{\hcompbig{\evalterm}{\wkn{\hcomp{t}{\ind{u}}}{x}, x}} &
\lam{x}
	{\hcompbig
		{\evalterm}
		{\hcomp{\wkn{t}{x}}{\wkn{\ind{u}}{x}, x}, 
			 \hcomp{x}{\wkn{\ind{u}}{x}, x}}
	} 
\arrow{l}[yshift=-2mm]
	{\lam{x}{\hcompsmall{\evalterm}{\phiCell{2}{t; \ind{u}}, \indproj{m+1}{}}}}
\end{td}
\normalsize
}
\end{enumerate}
\begin{proof}
Long but direct calculations using the universal property of 
$\transExp{x \bind \alpha}$. 
\end{proof}
\end{mylemma}

The rewrite $\pushName$ is also compatible with the $\genEpsilonExp{}$-rewrite. 
In the simply-typed lambda 
calculus, for any terms $\Gamma, x : A \vdash t : B$ and $\Gamma \vdash u : A$
and any family
$(\Delta \vdash v_i : A_i)_{i=1, \,\dots\, , n}$, then 
\begin{equation} \label{eq:stlc-beta-eta-and-comp}
{\left(\app{\lam{x}{t}}{u}\right)}[v_i / x_i] 
	=_{\beta\eta} t[u/x][v_i/x_i] 
	= t\big[ u[v_i/x_i] / x, v_i/x_i \big]
\end{equation}
In $\langCartClosed$ this corresponds to the two derivations
\begin{prooftree}
\AxiomC{$\Gamma, x : A \vdash t : B$}
\UnaryInfC{$\Gamma \vdash \lam{x}{t} : \exptype{A}{B}$}
\AxiomC{$\Gamma \vdash u : A$}
\BinaryInfC{$\Gamma \vdash \heval{\lam{x}{t}, u} : B$}
\AxiomC{$(\Delta \vdash v_i : A_i)_{i=1,\dots,n}$}
\BinaryInfC{$\Delta \vdash \hcomp{\heval{\lam{x}{t}, u}}{x_i \mapsto v_i} : B$}
\end{prooftree}
and 
\begin{prooftree}
\AxiomC{$\Gamma, x : A \vdash t : B$}
\AxiomC{$(\Delta \vdash v_i : A_i)_{i=1,\dots,n}$}
\AxiomC{$\Gamma \vdash u : A$}
\BinaryInfC{$\Delta \vdash \hcomp{u}{x_i \mapsto v_i} : A$}
\AxiomC{$(\Delta \vdash v_i : A_i)_{i=1,\dots,n}$}
\TrinaryInfC{$\Delta \vdash 
	\hcomp{t}{x_i \mapsto v_i, x \mapsto \hcomp{u}{x_i \mapsto v_i}} : B$}
\end{prooftree}

Continuing the equalities-as-rewrites perspective---which we make precise 
in Proposition~\ref{prop:stlc-up-to-iso}---the 
equation~(\ref{eq:stlc-beta-eta-and-comp}) 
becomes the following lemma.

\begin{mylemma}
Let $\sig$ be any $\langCartClosed$-signature and 
$\Gamma := (x_i : A_i)_{i=1,\dots,n}$ and
$\Delta := {(y_j : B_j})_{j=1, \,\dots\, , m}$ 
be contexts.
If the judgements
$(\Gamma, x : A \vdash t : B)$ and 
$(\Gamma \vdash u : A)$ and 
$(\Delta \vdash v_i : A_i)_{i=1,\dots,n}$
are derivable in $\langCartClosed(\sig)$, 
then 
%
{\tikzcdset{arrow style=tikz, arrows={Rightarrow}}
\begin{td}
\hcompthree{\evalterm}{\lam{x}{t}, u}{\ind{v}} \arrow{r}{\assoc{}} 
\arrow[swap]{d}{\hcomp{\genEpsilonExp{x \smallbind t, u}}{\ind{v}}} &
\hcompbig{\evalterm}{\hcomp{(\lam{x}{t})}{\ind{v}}, \hcomp{u}{\ind{v}}} 
\arrow{d}{\heval{\pushName, \hcomp{u}{\ind{v}}}} \\

\hcompthree{t}{\id_\Gamma, x \mapsto u}{\ind{v}} \arrow[swap]{d}{\iso} &
\heval{\lam{x}{\hcomp{t}{\wkn{\ind{v}}{x}, x}}, \hcomp{u}{\ind{v}}} 
\arrow{d}{\genEpsilonExp
			{x \smallbind \hcomp{t}{\wkn{\ind{v}}{x}, x}, \hcomp{u}{\ind{v}}}} 
			\\

\hcomp{t}{\wkn{\ind{v}}{x}, \hcomp{u}{\ind{v}}} &
\hcomp{\hcomp{t}{\wkn{\ind{v}}{x}, x}}{\id_\Delta, x \mapsto 
\hcomp{u}{\ind{v}}} \arrow{l}{\iso} 
\end{td}
where the unlabelled isomorphisms are defined by commutativity of the following 
two diagrams:
\begin{td}[column sep = 4em]
\hcompthree{t}{\id_\Gamma, u}{\ind{v}} 
\arrow{r}
\arrow[swap]{d}{\assoc{}} &
\hcomp{t}{\wkn{\ind{v}}{x}, \hcomp{u}{\ind{v}}} \\

\hcomp{t}{\hcomp{\id_\Gamma}{\ind{v}}, \hcomp{u}{\ind{v}}} 
\arrow[swap]{r}{\hcompsmall{t}{{\indproj{\bullet}{}}, \hcomp{u}{\ind{v}}}} &
\hcomp{t}{\ind{v}, \hcomp{u}{\ind{v}}}
\arrow[swap]{u}{\hcomp{t}{\subid{}, \hcomp{u}{\ind{v}}}}
\end{td}
\begin{td}[column sep = 6em] 
\hcomp{\hcomp{t}{\wkn{\ind{v}}{x}, x}}{\id_\Delta, \hcomp{u}{\ind{v}}} 
\arrow[swap, Rightarrow]{d}{\assoc{}} 
\arrow[Rightarrow]{r} 
& 
\hcomp{t}{\wkn{\ind{v}}{x}, \hcomp{u}{\ind{v}}} \\

\hcompbig
	{t}
	{\hcomp
		{\wkn{\ind{v}}{x}}
		{\id_\Delta, \hcomp{u}{\ind{v}}}, 
	 \hcomp{x}{\id_\Delta, \hcomp{u}{\ind{v}}}
	} 
\arrow[swap, Rightarrow]{r}{\hcomp{t}{\assoc{}, \indproj{1}{}}} &
\hcompbig
	{t}
	{\hcomp
		{\ind{v}}
		{\hcomp{\ind{y}}{\id_\Delta, \hcomp{u}{\ind{v}}}}, 
     \hcomp{u}{\ind{v}}
	} 
\arrow[swap, Rightarrow]{u}
	{\hcompsmall{t}{\hcompsmall{\ind{v}}{\indproj{\bullet}{}}, 
	\hcomp{u}{\ind{v}}}}
\end{td} }
\begin{proof}
Unfold the definitions and apply coherence.
\end{proof}
\end{mylemma}

\subsection{\texorpdfstring{The free property of 
$\increstrict{\syncloneCartClosed{\sig}}$}{The free property of the 
syntactic model}}
\label{sec:syntactic-model-free-property}

In this section we shall make precise the relationship between  
$\increstrict{\syncloneCartClosed{\sig}}$ and the free cc-bicategory $\freeCartClosedBicat\sig$ on $\sig$ (Construction~\ref{constr:free-cc-bicat}). 
We establish two related results. 
First, we shall show that for any cc-bicategory $\ccBicat{\baseCat}$ and
$\langCartClosed$-homomorphism $h : \sig \to \baseCat$, there exists a
\Def{semantic intepretation} cc-pseudofunctor
$h\sem{-} : \increstrict{\syncloneCartClosed{\sig}} \to \baseCat$. Along the 
way, 
we shall observe that such an interpretation extends to the cc-bicategory 
defined by extending $\termCatContextExt(\sig)$ 
(Construction~\ref{constr:context-cart-termcat}) with exponentials. This 
cc-bicategory, in which every context appears as an object, will 
play an important role in the normalisation-by-evaluation proof of 
Chapter~\ref{chap:nbe}. 
Second, we shall show that 
$\increstrict{\syncloneCartClosed{\sig}}$ is biequivalent $\freeCartClosedBicat\sig$. Thus, one does not obtain a 
strict universal property in the style of 
Theorem~\ref{thm:unary-syntactic-model-of-lang-biclone-free-bicat} 
(for $\langBicat$) or
Theorem~\ref{thm:unary-contexts-fp-bicat} (for $\langCart$), but one does 
obtain such a universal property up to biequivalence.

\paragraph*{Semantic interpretation.}

The semantic
interpretation of $\langCartClosed$ follows the tradition of 
semantic interpretation of the simply-typed lambda 
calculus~\cite{Lambek1980, Lambek1985}. For a fixed cartesian closed category 
$\ccBicat{\catC}$ and $\stlc$-signature homomorphism 
$h : \sig \to \catC$, the interpretation of a judgement 
$(\Gamma \vdash t : B)$ in the simply-typed lambda calculus over $\sig$
is $h\sem{\Gamma \vdash t : B}$, where $h\sem{-}$ is the unique cartesian closed clone homomorphism
extending $h$ (so $h\sem{-}$ has domain the free cartesian closed clone on $\sig$---namely, the syntactic model of the simply-typed lambda calculus---and codomain the cartesian closed clone 
$\cloneFromProducts\catC$ constructed in Example~\ref{ex:cartesian-closed-clone-from-cc-cat} 
(page~\pageref{ex:cartesian-closed-clone-from-cc-cat})). 

\hide{
For STLC, it is common to interpret 
contexts by induction on their length, using only the terminal object and 
binary products (\eg~\cite{Pitts2000}).  Because we take $n$-ary products 
as our primitive, however, we interpret a context of length $n$ directly as an 
$n$-ary 
product. This enables us to interpret the \rulename{var} rule
$x_1 : A_1,\, \dots,\, x_n : A_n \vdash x_i : A_i$ as a straightforward 
projection
$\sem{A_1} \times \cdots \times \sem{A_n} \xra{\pi_i} \sem{A_i}$, but entails 
some work translating between the $(n+1)$-ary product
$\prodop_{n+1}(A_1, \dots, A_n, B)$ and the binary product
$\prodop_2{\left( \prodop_n(A_1, \dots, A_n), B \right)}$.
}

\begin{mypropn} \label{prop:semantic-interpretation}
For any unary $\langCartClosed$-signature $\sig$, cartesian closed bicategory 
$\ccBicat\baseCat$, and unary $\langCartClosed$-signature homomorphism 
$h : \sig \to \baseCat$, there exists a \Def{semantic interpretation}
$h\sem{-}$ assigning to every term $(\Gamma \vdash t : B)$ a 1-cell in 
$\baseCat$
and to every rewrite $(\Gamma \vdash \tau : \rewrite{t}{t'} : B)$ a 2-cell
in $\baseCat$. Moreover, this interpretation is sound in the sense that if
$(\Gamma \vdash \tau \equiv \tau' : \rewrite{t}{t'} : B)$ then
$h\sem{\Gamma \vdash \tau : \rewrite{t}{t'} : B} 
	= h\sem{\Gamma \vdash \tau' : \rewrite{t}{t'} : B}$. 
\begin{proof}
The $\langCartClosed$-signature homomorphism $h$ also defines a $\langCartClosed$-signature
homomorphism $\sig \to \bicloneFromProducts\baseCat$ from $\sig$ to the cartesian closed biclone
arising from the cartesian closed structure of $\baseCat$ (recall Example~\ref{ex:cc-bicategory-to-cc-biclone} on
page~\pageref{ex:cc-bicategory-to-cc-biclone}). It follows from the universal 
property of 
$\syncloneCartClosed{\sig}$ (Proposition~\ref{propn:synclonecartclosed-free-property}) that there exists a strict cartesian closed pseudofunctor of biclones 
$h\sem{-} : \syncloneCartClosed{\sig} \to \bicloneFromProducts\baseCat$. We take this to be the semantic interpretation. Soundness is then automatic. 
\end{proof}
\end{mypropn}

To avoid obstructing the flow of our discussion we leave the full description of the semantic interpretation to an appendix (Section~\ref{sec:semantic-interpretation}).

The following observation entails a weak universal property for $\nucleus{\syncloneCartClosed\sig}$.

\begin{mylemma} 
\label{lem:cc-pseudofunctor-from-cartesian-closed-biclone-pseudofunctor}
Let $\ccBicat\baseCat$ be a cc-bicategory and 
$\ccClone{ob(\baseCat)}{\bicloneFromProducts{\baseCat}}$ the associated 
cartesian closed biclone. Then, for any cartesian closed biclone 
$\cartClone{S}{\biclone}$ and cartesian closed pseudofunctor of biclones
$(F, \prodPres, \expPres) : \biclone \to 
	\bicloneFromProducts{\baseCat}$
such that $\prodPres_{\ind{X}} \iso \Id_{\prod_{i=1}^n FX_i}$ for all 
$X_1, \dots, X_n \in S \:\: (n \in \Nat)$,
the restriction to unary multimaps
$(\restrictedPseudofunctorToProduct{F}, \prodPres, \expPres) : 
	\nucleus\biclone \to \baseCat$ 
is a cc-pseudofunctor of bicategories. 
\begin{proof}
Define $\restrictedPseudofunctor{F}(X) := FX$ and
$\restrictedPseudofunctor{F}_{X, Y} :=
	F_{X;Y} : \nucleus{\biclone}(X,Y) = \biclone(X;Y)
			\to \baseCat(X, Y)$. 
The 2-cells $\phi^{\restrictedPseudofunctor{F}}$ and
$\psi^{\restrictedPseudofunctor{F}}$ are defined by restricting the 
2-cells 
$\phi$ and $\psi^{(i)}$ of $F$ to linear multimaps. The three axioms to check 
then follow from the three laws of a biclone pseudofunctor, restricted to 
linear multimaps.

For 
preservation of products, we are already given an equivalence
\[\seq{F\pi_1, \dots, F\pi_n} 
	: F{\big( \prodop_n(X_1, \dots, X_n) \big)}
		\leftrightarrows
	  \prodop_n(FX_1, \dots, FX_n) 
	: \prodPres_{\ind{X}}
\]
for every $X_1, \dots, X_n \in S \:\: (n \in \Nat)$
because tupling in $\bicloneFromProducts{\baseCat}$ is tupling in 
$\baseCat$. It follows that $(\restrictedPseudofunctorToProduct{F}, \prodPres)$ 
is an 
fp-pseudofunctor.

For preservation of exponentials, the cartesian closure of $F$ provides an 
equivalence
\[
\lambda{\big( F(\evalterm_{A,B}) \circ \seq{\pi_1, \pi_2} \big)} 
	: F(\expobj{A}{B}) 
		\leftrightarrows
	  (\expobj{FA}{FB}) 
	: \expPres_{A,B}
\]
for every $A, B \in S$ 
(recall from 
Example~\ref{ex:cc-bicategory-to-cc-biclone} the definition of currying in $\bicloneFromProducts{\baseCat}$). On the other hand, 
\begin{align*}
\evBar^{\restrictedPseudofunctorToProduct{F}}_{A,B} &:=
	\lambda{\big( 
		\restrictedPseudofunctorToProduct{F}(\evalterm_{A,B}) \circ 
		\prodPres_{A,B} 
	\big)} \\
	&\iso
\lambda{\big( 
		\restrictedPseudofunctorToProduct{F}(\evalterm_{A,B}) \circ \Id_{FA 
		\times FB} 
	\big)} 
		&\text{by assumption on $\prodPres$} \\
	&\iso 
\lambda{\big( 
		\restrictedPseudofunctorToProduct{F}(\evalterm_{A,B}) \circ \seq{\pi_1, 
		\pi_2} 
	\big)}
\end{align*}
Since $(f, \psinv{g})$ is an equivalence whenever $(g, \psinv{g})$ is an 
equivalence and $f \iso g$, it follows that 
$(\evBar^{\restrictedPseudofunctorToProduct{F}}_{A,B}, \expPres_{A,B})$ is an 
equivalence for every $A, B \in S$. Hence, $(F, \prodPres, \expPres)$ is a 
cc-pseudofunctor.
\end{proof}
\end{mylemma}

Applying this lemma to the semantic interpretation $h\sem{-}$ of Proposition~\ref{prop:semantic-interpretation}
immediately yields the following weak universal property of 
$\increstrict{\syncloneCartClosed{\sig}}$.

\begin{prooflesscor} \label{cor:weak-ump-for-restricted-contexts}
For any unary $\langCartClosed$-signature $\sig$, 
cc-bicategory $\ccBicat\baseCat$,
and \mbox{$\langCartClosed$-signature} homomorphism
$h : \sig \to \baseCat$,
there exists a cc-pseudofunctor 
$h\sem{-} : \increstrict{\syncloneCartClosed{\sig}} \to \baseCat$ 
such that 
$h\sem{-} \circ \inc = h$, for 
$\inc : \sig \hookrightarrow \increstrict{\syncloneCartClosed{\sig}}$ the 
inclusion. 
\end{prooflesscor}

For the normalisation-by-evaluation argument in Chapter~\ref{chap:nbe}
we shall work with sets of terms indexed by types and contexts. We shall therefore require 
a syntactic model in which all contexts appear. For this purpose we extend $\termCatContextExt(\sig)$
(Construction~\ref{constr:context-cart-termcat} on page~\pageref{constr:context-cart-termcat}) with 
exponentials. Recall from Section~\ref{sec:products-from-context-extension} 
that the resulting bicategory has two product structures: one from 
context extension, and the other from the type theory. We emphasise this 
fact in our notation.

\begin{myconstr} \label{constr:ccc-termcat}  
For any $\langCartClosed$-signature $\sig$, define a bicategory 
$\syncloneAtClosed{\sig}$ as 
follows. The objects are contexts $\Gamma, \Delta, \dots$. The 1-cells 
$\Gamma \to (y_j : B_j)_{j = 1, \,\dots\, , m}$ are $m$-tuples of 
\mbox{$\alpha$-equivalence}
classes of terms $(\Gamma \vdash t_j :B_j)_{j= 1, \,\dots\, , m}$ derivable in 
$\langCartClosed(\sig)$, and the \mbox{2-cells} 
$(\Gamma \vdash t_j :B_j)_{j= 1, \,\dots\, , m} 
	\To 
(\Gamma \vdash t_j' :B_j)_{j= 1, \,\dots\, , m}$ 
are 
$m$-tuples of $\alpha{\equiv}$-equivalence classes of rewrites \mbox{$(\Gamma 
\vdash \tau : \rewrite{t_j}{t'_j} : B_j)_{j= 1, \,\dots\, , m}$}. Vertical 
composition is given pointwise by the $\vert$ operation, and horizontal 
composition  
\begin{align*} 
(t_1, \,\dots\, , t_l), (u_1, \,\dots\, , u_m) &\mapsto (\hcomp{t_1}{x_i 
\mapsto 
u_i}, 
\dots, \hcomp{t_m}{x_i \mapsto u_i}) \\ 
(\tau_1, \,\dots\, , \tau_l), (\sigma_1, \,\dots\, , \sigma_m) &\mapsto 
(\hcomp{\tau_1}{x_i \mapsto \sigma_i}, \,\dots\, , \hcomp{\tau_m}{x_i \mapsto 
\sigma_i})
\end{align*} 
by explicit substitution. The identity on $\Delta = (y_j : B_j)_{j = 1, 
\,\dots\, , 
m}$ is $(\Delta \vdash y_j : B_j)_{j = 1, \,\dots\, , m}$. The structural 
isomorphisms $\l, \r$ and $\a$ are given pointwise by $\proj{}$, 
$\subid{}^{-1}$ and $\assoc{}$, 
respectively. 
\end{myconstr} 

We define exponentials in a similar way to the type-theoretic product structure 
on $\termCatContextExt(\sig)$ (Lemma~\ref{lem:context-product-structure}): 
following 
Remark~\ref{rem:exponentials-up-to-equiv}, the exponential 
$\exp{\Gamma}{\Delta}$ is defined to be
\[
\exp{\left( p : \prodop_n(A_1, \,\dots\, , A_n) \right)}
	{\left(q : \prodop_m(B_1, \,\dots\, , B_m)\right)}
\]
for 
$\Gamma := (x_i : A_i)_{i=1, \,\dots\, , n}$ 
and
$\Delta := (y_j : B_j)_{j=1,\dots,m}$.

\begin{myremark} \label{rem:cc:syntactic-cc-bicategories-biequivalent}

Since Lemma~\ref{lem:fp-bicat:equivalences} extends verbatim to 
$\syncloneAtClosed{\sig}$, one sees that 
$\syncloneAtClosed{\sig} \simeq \increstrict{\syncloneCartClosed{\sig}}$ for 
every unary
$\langCartClosed$-signature $\sig$
(\cf~Remark~\ref{rem:syntactic-cc-bicategories-biequivalent}). 
Indeed, it is plain from the two definitions that the full sub-bicategory
of $\syncloneAtClosed{\sig}$ consisting of just the unary contexts is exactly 
$\increstrict{\syncloneCartClosed{\sig}}$.
\end{myremark}

$\termCatContextExt(\sig)$ satisfies a weak universal property akin to Corollary~\ref{cor:weak-ump-for-restricted-contexts}. However, since this bicategory does not arise from $\syncloneCartClosed\sig$ we must define the interpretation pseudofunctor by hand.

\begin{mypropn} \label{propn:weak-universal-property}
For any unary \mbox{$\langCartClosed$-signature} $\sig$, 
cc-bicategory $\ccBicat\baseCat$,
and $\langCartClosed$-signature homomorphism
$h : \sig \to \baseCat$,
there exists a 
cc-pseudofunctor $h\sem{-} : \syncloneAtClosed{\sig} \to \baseCat$ (for the 
type-theoretic product structure of  
Lemma~\ref{lem:context-extension-bicat-type-product-structure}), such 
that 
$h\sem{-} \circ \inc = h$, for 
$\inc : \sig \hookrightarrow \syncloneAtClosed{\sig}$ the inclusion. 
\begin{proof}
As the notation suggests, we extend the interpretation
$h\sem{-}$ of Proposition~\ref{prop:semantic-interpretation} to $\syncloneAtClosed{\sig}$ by setting
\begingroup
\addtolength{\jot}{.5em}
\begin{align*}
h\semlr{(\Gamma \vdash t_j : B_j)_{j=1, \,\dots\, , m}} &:=
	\seqlr{h\sem{\Gamma \vdash t_1 : B_1}, 
			\dots, 
			h\sem{\Gamma \vdash t_m : B_m}} \\
h\semlr{(\Gamma \vdash \tau_j : \rewrite{t_j}{t_j'} : B_j)_{j=1, \,\dots\, ,m}} 
&:=
	\seqlr{h\sem{\Gamma \vdash \tau_1 : \rewrite{t_1}{t_1'} : B_1}, 
			\dots, 
			h\sem{\Gamma \vdash \tau_m : \rewrite{t_m}{t_m'} : B_m}}
\end{align*}
\endgroup
This is 
well-defined on $\alpha{\equiv}$-equivalence classes of rewrites by the soundness of the semantic interpretation. For  
preservation of composition, we define $\phi^{h\sem{-}}$ as follows 
(where $\Gamma := {(x_i : A_i)}_{i=1, \,\dots\, , n}$):
\begin{td}
h\semlr{(\Gamma \vdash t_j : B_j)_{j=1, \,\dots\, , m}} \circ
	h\semlr{(\Delta \vdash u_i : A_i)_{i=1, \dot, n}} 
\arrow{r}{\phi^{h\sem{-}}}
\arrow[equals]{d} &
h\semlr{(\Delta \vdash \hcomp{t_j}{x_i \mapsto u_i} : B_j)_{j=1, \,\dots\, , 
m}} \\
\seqlr{h\sem{t_j}^\Gamma}_j \circ \seqlr{h\sem{u_i}^\Delta}_i 
\arrow[swap]{r}{\postName} &
\seqlr{ h\sem{t_j}^\Gamma \circ \seq{h\sem{u_i}^\Delta}_i }_j
\arrow[equals]{u}
\end{td}
For preservation of identities, we take
\[
\psi^{h\sem{\Gamma}} := \Id_{h\sem{\Gamma}} 
	\XRA{\widehat{\etaTimes{}}_{\Id_{h\sem{\Gamma}}}}
	\seq{\pi_1, \,\dots\, , \pi_n} = 
	h\sem{(\Gamma \vdash x_i : A_i)_{i=1,\dots, n}}
\]
where $\widehat{\etaTimes{}}$ is defined in~(\ref{eq:def-of-widehat-eta}) on 
page~\pageref{eq:def-of-widehat-eta}. 
We check the three axioms of a pseudofunctor. For the left unit law, one 
derives the commutative diagram below, then applies the triangle law relating 
the unit $\etaTimes{}$ and counit $\epsilonTimes{}{}$ for products:
\begin{td}[row sep = 3em]
\Id_{h\sem{\Gamma}} \circ \seq{h\sem{u_i}^\Gamma}_i 
\arrow[phantom]{ddrr}[description, xshift=32mm]{\equals{nat.}}
\arrow{rr}{\iso}
\arrow[swap, bend right = 80]{dd}{\widehat{\etaTimes{}}_{\Id_{h\sem{\Gamma}}}} 
\arrow{d}[description]
	{\etaTimes{\Id} \circ \seq{h\sem{u_i}^\Gamma}_i} 
\arrow[bend left = 15]{dr}[description, yshift=-2mm]
	{\etaTimes{(\Id \circ \seq{h\sem{u_i}}_i)}} &
\: &
\seq{h\sem{u_i}^\Gamma}_i 
\arrow{dd}{\etaTimes{\seq{h\sem{u_i}}_i }} \\

\seqlr{\ind{\pi} \circ \Id_{h\sem{\Gamma}}} \circ \seq{h\sem{u_i}^\Gamma}_i 
\arrow{d}[description]{\iso} 
\arrow{dr}[description]{\postName} &
\seqlr{\ind{\pi} 
		\circ \left(\Id_{h\sem{\Gamma}} 
		\circ \seq{h\sem{u_i}^\Gamma}_i\right)}
\arrow{d}[description]{\iso} &
\: \\

\seq{\pi_1, \,\dots\, , \pi_n} \circ \seqlr{h\sem{u_i}^\Gamma}_i 
\arrow[phantom]{r}[description]{\equals{nat.}}
\arrow[swap]{d}{\postName} &
\seqlr{\left(\ind{\pi} \circ \Id_{h\sem{\Gamma}}\right) 
	\circ \seq{h\sem{u_i}^\Gamma}_i} 
\arrow{dl}[description]{\iso} 
\arrow{r}{\iso} &
\seqlr{\ind{\pi} \circ \seq{h\sem{u_i}^\Gamma}_i} 
\arrow[equals]{dll} \\

\seq{\ind{\pi} \circ \seq{h\sem{u_i}^\Gamma}_i} 
\arrow[swap]{d}{\seq{\epsilonTimesInd{\bullet}{} }} &
\: \\

\seq{h\sem{u_i}^\Gamma}_i
\end{td}
The 
unlabelled triangular shape is an easily-verified property of $\postName$ 
\big({\cf~Lemma~\ref{lem:PseudoproductCanonical2CellsLaws}, diagram (\ref{c:PostAndEta})}\big).
The right unit law is similar, and the associativity law follows directly from 
the naturality of $\postName$ and the observation 
that the following 
commutes~\big({\cf~Lemma~\ref{lem:PseudoproductCanonical2CellsLaws}(\ref{c:DoublePostLaw})}\big):
\begin{td}[column sep = 4em]
\left(\seq{\ind{f}} 
	\circ g\right) 
	\circ h 
\arrow{rr}[yshift=0mm]{\postName \circ h} 
\arrow[swap]{d}{\iso} &
\: &
\seq{\ind{f} \circ g} \circ h 
\arrow{d}{\postName} \\
\seq{\ind{f}} \circ (g \circ h)
\arrow[swap]{r}{\postName} &
\seq{\ind{f} \circ (g \circ h)}
\arrow[swap]{r}{\seq{\iso, \,\dots\,, \iso}} &
\seq{(\ind{f} \circ g) \circ h}
\end{td}

Now we want to show that $h\sem{-}$ is a cc-pseudofunctor. We start with 
products. It is immediate from the definition that, for any family of 
unary contexts
$(x_1 : A_1), \,\dots\, , (x_n : A_n) \:\: (n \in \Nat)$,
the pseudofunctor $h\sem{-}$ strictly preserves the data making
${\left(p : \prodop_n(A_1, \,\dots\, , A_n)\right)} = \prod_{i=1}^n (x_i : A_i)$
an $n$-ary product.
More generally, for contexts 
$\Gamma^{(i)} := (x_j^{(i)} : A_j^{(i)})_{j=1,\dots, \len{\Gamma^{(i)}}}
	\:\: (i=1,\dots, n)$, 
the $n$-ary product 
$\Gamma^{(1)} \times \dots \times \Gamma^{(n)}$
is interpreted as 
\[
h\semlr{p : \prodop_n \big( \prod_{\len{\Gamma^{(1)}}} \ind{A}^{(1)}, \,\dots\, 
, 
		\prodop_{\len{\Gamma^{(n)}}} \ind{A}^{(n)}\big)} =
		\prodop_{i=1}^n \prod_{j=1}^{\len{\Gamma^{(j)}}} h\sem{A_j^{(i)}} =
		\prodop_{i=1}^n h\sem{\Gamma^{(i)}}
\]
and the $i$th projection 
\[
\left(p : \prodop_n \big( \prodop_{\len{\Gamma^{(1)}}} \ind{A}^{(1)}, 
\,\dots\, , 
		\prodop_{\len{\Gamma^{(n)}}} \ind{A}^{(n)} \big) \vdash 
		\hpi{j}{\pi_i(p)} : A^{(i)}_j\right)_{j=1, \dots, \len{\Gamma^{(i)}}}
\]
is interpreted as 
$\prodop_{i=1}^n h\sem{\Gamma^{(i)}} 
	\xra{
		\seqlr
		{\pi_1 \circ \pi_i, \dots, \pi_{\len{\Gamma^{(i)}}} \circ \pi_i}
	} 
	\prod_{j=1}^{\len{\Gamma^{(i)}}} h\sem{A_j^{(i)}}
	= h\sem{\Gamma^{(i)}}
$. 
To witness that $h\sem{-}$ preserves products, then, one can take
$\prodPres_{\Gamma^{(\bullet)}}$ to be the identity, with witnessing 2-cell
\newpage
\begin{align*}
\seqlr{\seqlr{\ind{\pi} \circ \pi_1}, \dots, \seqlr{\ind{\pi} \circ \pi_n}} 
	&\XRA{\seqlr{\postName^{-1}, \dots, \postName^{-1}}} 
	\seqlr{ \seq{\pi_1, \dots, \pi_{\len{\Gamma^{(1)}}}} \circ \pi_1, \dots, 
	\seq{\pi_1, \dots, \pi_{\len{\Gamma^{(n)}}}} \circ \pi_n } \\
	&\XRA{\seq{\widehat{\etaTimes{}}^{-1}, \dots, \widehat{\etaTimes{}}^{-1}}} 
	\seq{ \Id_{h\sem{\Gamma^{(1)}}} \circ \pi_1, \dots, 	
	\Id_{h\sem{\Gamma^{(n)}}} \circ \pi_n } \\
	&\iso \seq{\pi_1, \dots, \pi_n} \\
	&\XRA{\widehat{\etaTimes{}}^{-1}} \Id_{h\sem{\prod_i \Gamma^{(i)}}}
\end{align*}
Note we once again use the 2-cell $\widehat{\etaTimes{}}$ defined 
in~(\ref{eq:def-of-widehat-eta}) on 
page~\pageref{eq:def-of-widehat-eta}.

For exponentials, one sees that (where $\Delta := (y_j : B_j)_{j=1, \dots, m}$):
\begingroup
\addtolength{\jot}{.3em}
\begin{align*}
h\sem{\exp{\Gamma}{\Delta}} &= 
	h\semlr{
		\exp{\big(p : \prodop_n(A_1, \,\dots\, , A_n)\big)}
			{\big(q : \prodop_m(B_1, \,\dots\, , B_m)\big)}} \\
	&= h\sem{f : \exptype{\prodop_n(A_1, \,\dots\, , A_n)}
							{\prodop_m(B_1, \,\dots\, , B_m)}} \\
	&= \exp{\big(\prodop_{i=1}^n h\sem{A_i}\big)}
			{\big(\prodop_{j=1}^n h\sem{B_j}\big)}
\end{align*}
\endgroup
and
\begingroup
\addtolength{\jot}{.3em}
\begin{align*}
h\sem{(\exp{\Gamma}{\Delta}) \times \Gamma} 
	&= h\semlr{p :\prodop_2\big(\exptype{\prodop_n \ind{A}}
								{\prodop_m \ind{B}}, 
								\prodop_n \ind{A} \big)} \\
	&= \big(\exp{\prodop_{i=1}^n h\sem{A_i}}{\prodop_{j=1}^n h\sem{B_j}}\big) 
	\times \prodop_{i=1}^n h\sem{A_i}
\end{align*}
\endgroup
It follows that
$\evBar^{h\sem{-}}_{\Gamma, \Delta}$
is the currying of
\begin{equation*} 
\begin{aligned}
h\big\llbracket p :\prodop_2\big(\exptype{\prodop_n \ind{A}}
								{\prodop_m \ind{B}}, 
								\prodop_n \ind{A} \big) &\vdash 
					\heval{\pi_1(p), \pi_2(p)}  : \prodop_m \ind{B} 
					\big\rrbracket 
					\circ \Id_{(h\sem{\scriptsizeexpobj{\Gamma}{\Delta}} \times 
					h\sem{\Gamma})} 
					\\
		&= \left(\eval_{h\sem{\Gamma}, h\sem{\Delta}} 
			\circ \seq{\pi_1, \pi_2}\right) 
			\circ \Id_{(h\sem{\scriptsizeexpobj{\Gamma}{\Delta}} 
				\times h\sem{\Gamma})}
\end{aligned}
\end{equation*}
Hence, $\evBar^{h\sem{-}}_{\Gamma, \Delta}$ is naturally isomorphic to the 
identity via the composite
\begin{align*}
\lambda \big(  &\left(\eval_{(h\sem{\Gamma}, h\sem{\Delta})} 
			\circ \seq{\pi_1, \pi_2}\right) 
			\circ \Id_{(h\sem{\scriptsizeexpobj{\Gamma}{\Delta}} 
					\times h\sem{\Gamma})} 
			\big)  \\
		&\iso 
		\lambda \big(  \eval_{(h\sem{\Gamma}, h\sem{\Delta})} 
		\circ \seqlr
			{\pi_1 \circ 
				\Id_{(h\sem{\scriptsizeexpobj{\Gamma}{\Delta}} 
					\times h\sem{\Gamma})},  
			\Id_{\pi_2\circ 
				(h\sem{\scriptsizeexpobj{\Gamma}{\Delta}} 
					\times h\sem{\Gamma})}}\big) 
			\\
		&\iso
		\lambda \big(  \eval_{(h\sem{\Gamma}, h\sem{\Delta})} 
		\circ (\Id_{h\sem{\scriptsizeexpobj{\Gamma}{\Delta}}} \times 
		\prodop_m 
		h\sem{\ind{B}})\big) \\
		&\overset{\etaExp{}}{\iso} 
		\Id_{h\sem{\scriptsizeexpobj{\Gamma}{\Delta}}}
\end{align*}
and $h\sem{-}$ is a cc-pseudofunctor.
\end{proof}
\end{mypropn}

Our aim now is to prove that 
$\increstrict{\syncloneCartClosed{\sig}}$ is biequivalent to the free 
cc-bicategory on the unary $\langCartClosed$-signature $\sig$ (defined in 
Construction~\ref{constr:free-cc-bicat}), and hence that 
$\langCartClosed$ is the internal language for cc-bicategories up to biequivalence.

\paragraph*{$\increstrict{\syncloneCartClosed{\sig}}$ is biequivalent to 
$\freeCartClosedBicat{\sig}$.}
Fix a unary $\langCartClosed$-signature $\sig$. We shall show that the 
canonical cc-pseudofunctors
$\ext{\iota} : \freeCartClosedBicat{\sig} \to \syncloneAtClosed{\sig}$ 
and
$\iota\sem{-} : \syncloneAtClosed{\sig} \to \freeCartClosedBicat{\sig}$
extending the respective inclusions
$\sig \hookrightarrow \freeCartClosedBicat{\sig}$
and
$\sig \hookrightarrow \syncloneAtClosed{\sig}$
induce a biequivalence 
$\syncloneAtClosed{\sig} \simeq \freeCartClosedBicat{\sig}$. 
(These cc-pseudofunctors are defined in Lemma~\ref{lem:free-cc-bicat-on-sig} and
Proposition~\ref{propn:weak-universal-property}, respectively.)  
One then obtains the required 
biequivalence by restricting 
$\syncloneAtClosed{\sig}$ to 
unary contexts (recall 
Remark~\ref{rem:cc:syntactic-cc-bicategories-biequivalent}). 

\begin{myremark}
Because the pseudofunctor $\ext{\iota}$ is defined inductively using the cartesian 
closed structure of $\syncloneAtClosed{\sig}$, we must be explicit about which cartesian 
closed structure we choose. We take the \emph{type-theoretic} product 
structure, so that the composite $\ext{\iota} \circ \iota\sem{-}$ takes an arbitrary 
context $\Gamma$ to an (equivalent) unary context. Because the restriction of 
$\syncloneAtClosed{\sig}$ to unary contexts is exactly 
$\nucleus{\syncloneCartClosed{\sig}}$, this ensures that the biequivalence 
we construct will restrict to $\nucleus{\syncloneCartClosed{\sig}}$ with its 
canonical cartesian closed structure (namely, that of Remark~\ref{rem:bicat-exp-structure-not-preserved}). 
Of course, up to biequivalence of the underlying bicategories, 
the uniqueness of products and exponentials ensures 
that the choice of cc-bicategory is immaterial
(recall Remark~\ref{rem:cc-pseudofunctor-for-different-struct} and 
Lemma~\ref{lem:bicat-to-cartclosedbicat-equivalence-lift}).   
\end{myremark}

Our two-step approach reflects two intended applications. In this 
chapter we wish to 
prove a free property, so restrict to unary contexts, but in 
Chapter~\ref{chap:nbe} 
we wish to interpret the syntax of $\langCartClosed$ varying over a 
(2-)category of contexts, and so require all contexts. 

\begin{myremark}
Although we present the argument indirectly here, 
it is also possible to prove directly that
the canonical cc-pseudofunctors induce a biequivalence 
$\nucleus{\syncloneCartClosed{\sig}} \simeq \freeCartClosedBicat{\sig}$. 
The calculations involved are similar to those we shall see below.
\end{myremark}

We begin by showing that 
$\inc\sem{-} \circ \ext{\inc} \simeq \id_{\freeCartClosedBicat{\sig}}$. Recall from 
Proposition~\ref{propn:weak-universal-property} 
that $\inc\sem{-}$ preserves products and exponentials up to equivalence in a 
particularly strong way, in 
the sense that 
$\seq{\inc\sem{\pi_1}, \dots, \inc\sem{\pi_n}} \iso \id$ and 
$\evBar^{\inc\sem{-}} \iso \id$. One 
may therefore apply Corollary~\ref{cor:termcatccc-unique-up-to-equivalence}.

\begin{mypropn} \label{propn:first-biequivalence}
For any unary $\langCartClosed$-signature $\sig$, the composite
$\inc\sem{-} \circ \ext{\inc} : \freeCartClosedBicat{\sig} \to \freeCartClosedBicat{\sig}$ 
induced by the following diagram is equivalent to $\id_{\freeCartClosedBicat{\sig}}$:
\begin{td}[column sep = 3em]
\freeCartClosedBicat{\sig} 
\arrow{r}{\ext{\inc}} &
\syncloneAtClosed{\sig} 
\arrow{r}{\inc\sem{-}} &
\freeCartClosedBicat{\sig} \\

\sig 
\arrow[equals]{r}
\arrow[hookrightarrow]{u}{\inc} &
\sig 
\arrow[equals]{r}
\arrow[hookrightarrow]{u}{\inc} &
\sig 
\arrow[hookrightarrow, swap]{u}{\inc}
\end{td}
\begin{proof}
The diagram commutes, and the composite $\inc\sem{-} \circ \ext{\inc}$ is 
certainly a cc-pseudofunctor. Since $\ext{\inc}$ is strict and 
$\inc\sem{-}$ has $\prodPres$ and $\expPres$ both given by the identity, 
Corollary~\ref{cor:termcatccc-unique-up-to-equivalence} applies. Hence 
$\inc\sem{-} \circ \ext{\inc}$ is equivalent to the 
unique strict cc-pseudofunctor $\freeCartClosedBicat{\sig} \to \freeCartClosedBicat{\sig}$ extending
the inclusion $\sig \hookrightarrow \freeCartClosedBicat{\sig}$. Since the identity is 
such 
a strict cc-pseudofunctor, it follows that 
$\inc\sem{-} \circ \ext{\inc} \simeq \id_{\freeCartClosedBicat{\sig}}$, as required.
\end{proof}
\end{mypropn}

We shall see in Chapter~\ref{chap:nbe} that this result is 
crucial to the normalisation-by-evaluation proof. Roughly speaking, it plays 
the same role as the 1-categorical observation that the canonical map from the 
free cartesian closed category to itself is the identity.

We now turn to showing that $\usem{}{\inc} \circ \inc\sem{-}$ is equivalent to 
the identity. To this end, observe that for any context
 $\Gamma := (x_i : A_i)_{i=1, \,\dots\, , n}$, 
\[
\usem{\left(\inc\sem{\Gamma}\right)}{\inc} = 
	\ext{\inc}{\left(\prodop_n(A_1, \,\dots\, , A_n)\right)} =
	\left( p : \prodop_n(A_1, \,\dots\, , A_n) \right)
\]
We define a pseudonatural transformation 
$(\altNat, \altCell) : \usem{}{\inc} \circ \inc\sem{-} \To 
							\id_{\syncloneAtClosed{\sig}}$
with components
$\altNat_\Gamma : \usem{{\left(\inc\sem{\Gamma}\right)}}{\inc} \to \Gamma$ 
given by the equivalence 
\begin{equation*} \label{eq:context-product-equivalence}
\begin{tikzcd}[column sep = 14em]
\Gamma 
\arrow[yshift=1mm]{r}
	{(\Gamma \vdash \pair{x_1, \,\dots\, , x_n} : \prod_n \ind{A})} &
\big( p : \prod_n(A_1, \,\dots\, , A_n) \big)
\arrow[yshift=-1mm]{l}
	{(p : \prod_n(A_1, \,\dots\, , A_n) \vdash 
			\pi_i(p) : A_i)_{i=1, \,\dots\, , n}}
\end{tikzcd}
\end{equation*}
constructed in Lemma~\ref{lem:fp-bicat:equivalences} (page~\pageref{lem:fp-bicat:equivalences}).
We are therefore required to provide an invertible 2-cell filling the diagram 
below for every judgement $(\Gamma \vdash t : B)$:
\begin{equation} \label{eq:usem-inc-to-id}
\begin{tikzcd}[column sep = 6em] 
\usem{\left(\inc\sem{\Gamma}\right)}{\inc}
\arrow[phantom]{dr}[description]{\twocell{\altCell_t}}
\arrow[swap]{d}{\altNat_\Gamma}
\arrow{r}{\usem{{\left(\inc\sem{\Gamma \vdash t : B}\right)}}{\inc}} &
\usem{{\left(\inc\sem{y : B}\right)}}{\inc}
\arrow{d}{\altNat_B} \\
\Gamma 
\arrow[swap]{r}{(\Gamma \vdash t : B)} &
(y : B)
\end{tikzcd}
\end{equation}

\begin{myconstr} \label{constr:altcell-defined}
For any $\langCartClosed$-signature $\sig$, we define a family of 2-cells
$\altCell_t$ filling~(\ref{eq:usem-inc-to-id}) in $\syncloneAtClosed{\sig}$.
Unfolding the anticlockwise composite, one sees that
\begin{align*}
(\Gamma \vdash t : B) \circ \altNat_\Gamma &= 
		(\Gamma \vdash t : B) \circ 
		\big( p : \prodop_n \ind{A} \vdash 
			\pi_i(p) : A_i\big)_{i=1,\dots,n} \\
		&= \big( p : \prodop_n (A_1, \,\dots\, , A_n) \vdash 
			\hcomp{t}{x_i \mapsto \pi_i(p)} : B \big)
\end{align*}
Thus, it suffices to define 
2-cells $\natCell_t$ of type 
$(p : \prod_n \ind{A} \vdash \rewrite{\overline{t}}
	{\hcomp{t}{x_i \mapsto \pi_i(p)}} : B)$, where $\overline{t}$ is 
the term in the judgement 
$\ext{\inc}{\left(\inc\sem{\Gamma \vdash t : B}\right)}$. Since 
$\altNat_B$ is simply 
$(y : B \vdash y :B)$, one
may then define 
the required 2-cell 
$\altCell_t$ 
to be
\[
\altCell_t := \hcomp{y}{\overline{t}} \XRA{\indproj{1}{\overline{t}}} 
		\overline{t} \XRA{\natCell_t} \hcomp{t}{x_i \mapsto \pi_i(p)}
\]
We define $\natCell_t$ by induction on the derivation of $t$.

\vspace{1em}
\subproof{\rulename{var} case.} For
$(\Gamma \vdash x_k : A_k)$ the corresponding term $\overline{x_k}$ is 
$\big( p : \prodop_n \ind{A} \vdash \pi_k(p) : A_k \big)$, so 
we define
\[
\natCell_{x_i} := \big( p : \prodop_n \ind{A} \vdash 
		\indproj{-k}{\ind{\pi}(p)}
		: \rewrite{\pi_k(p)}{\hcomp{x_k}{x_i \mapsto \pi_i(p)}} : A_k \big)
\]
\subproof{\rulename{const} case.} For any constant 
$c \in \graph(A,B)$, the judgement 
$\ext{\inc}\inc\sem{x : A \vdash c (x) : B}$
is simply $( x : A \vdash c(x) : B)$. Since the context is unary, 
$\altNat_\Gamma$ is the identity and we may take 
$\natCell_{c(x)}$ to be canonical structural isomorphism. 

\vspace{1em}
\subproof{\rulename{proj} case.}  Observing that
$\usem{}{\inc} \circ \inc\sem{-}$ is the identity on
$\left( {p : \prodop_n (A_1, \,\dots\, , A_n)} \vdash \pi_i(p) : A_i \right)$, 
we 
take the canonical isomorphism
\begin{td}[column sep = 8em, row sep = 5em ] 
\big( p : \prodop_n(A_1, \,\dots\, , A_n) \big) 
\arrow{dr}[description]{(p : \prod_n \ind{A} \vdash \pi_i(p) : A_i)}
\arrow[swap]{d}{(p : \prod_n \ind{A} \vdash p : \prod_n \ind{A})}
\arrow{r}{(p : \prod_n \ind{A} \vdash \pi_i(p) : A_i)} &
(x_i : A_i) 
\arrow[phantom]{dl}[description, near start, xshift=6mm, yshift=1mm]{\iso} 
\arrow{d}{(x_i : A_i \vdash x_i : A_i)} \\

\big( p : \prodop_n(A_1, \,\dots\, , A_n) \big) 
\arrow[phantom]{ur}[description, near start, xshift=-6mm, yshift=-1mm]{\iso}
\arrow[swap]{r}{(p : \prod_n \ind{A} \vdash \pi_i(p) : A_i)} &
(x_i : A_i) 
\end{td}

\vspace{1em}
\subproof{\rulename{tup} case.} 
From the induction hypothesis one obtains 
$
\left(p : \prodop_n \ind{A} \vdash \natCell_{t_i} : 
	\rewrite{\overline{t_j}}{\hcomp{t_j}{x_i \mapsto \pi_i(p)}} : B_j\right) 
$ 
for 
$j=1,\dots,m$. So for $\natCell_{\pair{t_1, \,\dots\, , t_m}}$ we take the 
composite rewrite
\[
\pair{\overline{t_1}, \,\dots\, , \overline{t_m}} 
\XRA{\pair{\natCell_{t_1}, \,\dots\, , \natCell_{t_m}}}
\pair{\hcomp{t_1}{\ind{\pi}(p)}, \,\dots\, , \hcomp{t_m}{\ind{\pi}(p)}} 
\XRA{\postName^{-1}}
\hcomp{\pair{t_1, \,\dots\, , t_m}}{\ind{\pi}(p)} 
\]
of type $\prodop_m(B_1, \,\dots\, , B_m)$ in context 
$\big( p : \prodop_n (A_1, \,\dots\, , A_n) \big)$.

\vspace{1em}
\enlargethispage{3\parskip}
\subproof{\rulename{eval} case.} 
The evaluation 1-cell 
$(f : \exptype{A}{B}) \times (x : A) \to (y : B)$ 
in $\syncloneAtClosed{\sig}$ with the type-theoretic product structure is
$\big( p : (\exptype{A}{B}) \times A \vdash 
		\heval{\pi_1(p), \pi_2(p)} : B \big) $, so one obtains
\begin{align*}
\usem{{\left(\inc
		{
			\sem{f : \exptype{A}{B}, x : A \vdash \evalterm(f,x) : B}
		 }
	\right)}}{\inc} 
	&= \ext{\inc}(\eval_{\inc\sem{A}, \inc\sem{B}}) \\
	&= \big( p : (\exptype{A}{B}) \times A \vdash 
			\heval{\pi_1(p), \pi_2(p)} : B \big)
\end{align*}
\vspace{-3\parskip}
We therefore define $\natCell_{\evalterm(f,x)}$ to be the identity. 

\vspace{1em}
\subproof{\rulename{lam} case.}
The exponential 
transpose of a term
$({p : Z \times B \vdash t : C})$ in 
$\syncloneAtClosed{\sig}$
is
\[
(z : Z \vdash \lam{x}{(\hcomp{t}{p \mapsto \pair{z, x}})} : \exptype{B}{C}) 
\] 
It follows that
\begin{align*}
\usem{\left(\inc\sem{\Gamma \vdash \lam{x}{t} : \exptype{B}{C}}\right)}{\inc} 
		&= \lambda {\left( q  : \prodop_2 (\prodop_n \ind{A}, B) 
					\vdash \hcompbig{\overline{t}}
						{\pair{  \hcomp{\ind{\pi}}{\pi_1(q)}, \pi_2(q) }} : 
						C\right)} \\
		&= \big( p : \prodop_n \ind{A} \vdash 
				\lam{x}{\hcompthreebigfst{\overline{t}}
							{\pair{\hcomp{\ind{\pi}}{\pi_1(q)}, 				
													\pi_2(q) }}
							{\pair{p, x}}} : 
							\exptype{B}{C} \big)
\end{align*}
Now, the induction hypothesis provides the 2-cell
$\big( s : \prodop_n(A_1, \,\dots\, , A_n, B) \vdash 
		\natCell_t : {\rewrite{\overline{t}}{\hcomp{t}{x_i \mapsto \pi_i(s)}} : 
		C} \big)$
so for $\natCell_{\lam{x}{t}}$ we begin by defining a composite 
$\vartheta_t$ by
{\tikzcdset{arrow style=tikz, arrows={Rightarrow}}
\small
\begin{td}
\hcompthreebigfst{\overline{t}}
	{\pair{\hpi{1}{\pi_1(q)}, \,\dots\, , \hpi{n}{\pi_1(q)}, \pi_2(q)}}
	{\pair{p,x}}
\arrow[swap]{d}{\assoc{}}
\arrow[bend left=23]{dddr}{\vartheta_t} &
\: \\

\hcompbig{\overline{t}}
	{\hcomp{\pairbig{\hpi{1}{\pi_1(q)}, \,\dots\, , \hpi{n}{\pi_1(q)}, 
	\pi_2(q)}}
	{\pair{p,x}}} 
\arrow[swap]{d}{\hcomp{\overline{t}}{\postName}} &
\: \\

\hcompbig{\overline{t}}
	{\pairbig{ \hcompthree{\pi_1}{\pi_1(q)}{\pair{p,x}}, \,\dots\, , 
		\hcompthree{\pi_n}{\pi_1(q)}{\pair{p,x}}, 
		\hpi{2}{\pair{p,x}} }} 
\arrow[swap, bend right = 9]{dr}{
		\hcomp
			{\overline{t}}
			{\pair{\gamma_1, \,\dots\, , \gamma_n, \epsilonTimesInd{2}{p,x}}}} &
\: \\

\: &
\hcompbig{\overline{t}}{\pair{\hcomp{\pi_1}{p} , \,\dots\, , \hcomp{\pi_n}{p} , 
x}}
\end{td}
\normalsize
in context $\big( p : \prodop_n (A_1, \,\dots\, , A_n), x : B \big)$, where 
$\gamma_k$ is defined, in the same context, to be
\[
\gamma_k :=
\hcompthree{\pi_k}{\pi_1(q)}{\pair{p,x}}
\XRA{\assoc{}}
\hcompbig{\pi_k}{\hpi{1}{\pair{p,x}}}
\XRA{\hcompsmall{\pi_k}{\epsilonTimesInd{1}{p,x}}}
\hcomp{\pi_k}{p} 
\]
%
for $k=1, \,\dots\, , n$. We then define $\natCell_{\lam{x}{t}}$ to be the 
composite
\begin{td}
\lam{x}{\hcompthreebigfst{\overline{t}}
						{\pair{\hcomp{\ind{\pi}}{\pi_1(q)}, \pi_2(q) }}
						{\pair{p, x}}} 
\arrow{r}{\natCell_{\lam{x}{t}}}
\arrow[swap]{d}{\lam{x}{\vartheta_t}} &
\hcomp{(\lam{x}{t})}{\pi_1(p), \,\dots\, , \pi_n(p)} \\

\lam{x}
	{\hcompbig
		{\overline{t}}
		{\pair{\hcomp{\pi_1}{p} , \,\dots\, , \hcomp{\pi_n}{p} , x}}
	} 
\arrow[swap]{d}
	{\lam{x}
		{\hcomp
			{\natCell_t}
			{\pair{\hcomp{\pi_1}{p} , \,\dots\, , \hcomp{\pi_n}{p} , x}}
		}
	} &
\: \\

\lam{x}{\hcompthree
		{t}
		{\pi_1(s), \,\dots\, , \pi_n(s), \pi_{n+1}(s)}
		{\pair{\hcomp{\pi_1}{p} , \,\dots\, , \hcomp{\pi_n}{p} , x}}
		} 
\arrow[swap]{d}{\lam{x}{\assoc{}}} &
\: \\

\lam{x}{\hcompbig
		{t}
		{\hcomp
			{\ind{\pi}}
			{\pair{\hcomp{\pi_1}{p} , \,\dots\, , \hcomp{\pi_n}{p} , x}}
		}} 
\arrow[swap]{r}
		{
		\lam{x}{\hcompsmall
				{t}
				{\epsilonTimesInd{\bullet}{}}
				}
		} &
\lam{x}{\hcomp
	{t}
	{\hcomp{\pi_1}{p} , \,\dots\, , \hcomp{\pi_n}{p} , x}
	} 
\arrow[swap]{uuu}{\pushName^{-1}}
\end{td}
}
It remains to consider the cases of explicit substitutions and $n$-tuples of 
terms. We take the latter first and then put it to work for explicit 
substitutions.

\vspace{1em}
\subproof{$n$-tuples case.}
For contexts 
$\Gamma := (x_i : A_i)_{i=1, \,\dots\, ,n}$ and
$\Delta := (z_j : Z_j)_{j=1, \,\dots\, , m}$
and an $n$-tuple 
$(\Delta \vdash t_i : A_i)_{i=1, \,\dots\, ,n} : \Delta \to \Gamma$, we 
directly 
define the rewrite
$\altCell_{(t_j)_{j=1, \,\dots\, ,m}}$ filling 
\begin{td}[column sep = 8em]
\big( q : \prodop_m (Z_1, \,\dots\, , Z_m) \big)
\arrow{r}[yshift=2mm]{\left(q : \prod_m \ind{Z} \vdash 
		\pair{\overline{t_1}, \,\dots\, , \overline{t_n}} : \prod_n \ind{A} 
		\right)}
\arrow[phantom]{dr}[description]{\twocell{\altCell_{(t_i)_{i=1, \,\dots\, ,n}}}}
\arrow[swap]{d}{\simeq} &
\big( p : \prodop_n(A_1, \,\dots\, , A_n) \big) 
\arrow{d}{\simeq} \\

\Delta 
\arrow[swap]{r}{(\Delta \vdash t_i : A_i)_{i=1, \,\dots\, ,n}} &
\Gamma
\end{td}
to be the $n$-tuple with components
\[
\altCell_{(t_i)_{i=1, \,\dots\, ,n}} :=
\hpi{k}{\pair{\overline{t_1}, \,\dots\, , \overline{t_n}}} 
\XRA{\epsilonTimesInd{k}{}}
\overline{t_k} 
\XRA{\natCell_{t_k}}
\hcomp{t_k}{\pi_1(q), \,\dots\, , \pi_m(q)}
\]
for $k=1, \,\dots\, , n$. 

\vspace{1em}
\subproof{\rulename{hcomp} case.}
For explicit substitutions 
$(\Delta \vdash \hcomp{t}{x_i \mapsto u_i} : B) = 
	(\Gamma \vdash t : B) \circ (\Delta \vdash u_i : A_i)_{i=1,\dots,n}$
we take the definition from the associativity law of a pseudonatural 
transformation. Thus, we define $\altCell_{\hcomp{t}{x_i \mapsto u_i}}$ to be 
the pasting diagram
\begin{td}[column sep = 5em, row sep = 3em]
\big( q : \prodop_m (B_1, \,\dots\, , B_m) \big) 
\arrow{rr}{\left(q : \prod_m \ind{B} \vdash 
	\hcomp{\overline{t}}{\pair{\overline{u_1}, \,\dots\, , \overline{u_n}}} : 
	C\right)}
\arrow[swap]{ddd}{\simeq}
\arrow{dr}[description]{(q : \prod_m \ind{B} \vdash 
		\pair{\overline{u_1}, \,\dots\, , \overline{u_n}} : \prod_n \ind{A})} &
\: &
(z : C) 
\arrow{ddd}{(z : C \vdash z : C)} \\

\: &
\big( p : \prodop_n(A_1, \,\dots\, , A_n) \big) 
\arrow[phantom]{dl}[description]{\twocell{\altCell_{(u_i)_{i=1, \,\dots\, , 
n}}}}
\arrow[phantom]{dr}[description]{\twocell{\altCell_t}}
\arrow{d}[description]{\simeq}
\arrow{ur}[description]{\left(p : \prod_n(A_1, \,\dots\, , A_n) \vdash 
\overline{t} 
	: C\right)} 
&
\: \\

\: &
\Gamma 
\arrow{dr}[description]{(\Gamma \vdash t : C)} &
\: \\

\Delta 
\arrow[swap]{rr}{(\Delta \vdash \hcomp{t}{x_i \mapsto u_i} : C)}
\arrow{ur}[description]{(\Delta \vdash u_i : A_i)_{i=1, \,\dots\, , n}} &
\: &
(z : C)
\end{td}
\end{myconstr}

The preceding construction does indeed define a pseudonatural transformation. 
It is clear that each $\altCell_t$ is natural, so it remains to check the 
unit and associativity laws. For the unit law, we are required to show the 
following equality of pasting diagrams 
for every context $\Gamma := (x_i : A_i)_{i=1,\dots, n}$:
\begin{equation*}
\begin{tikzcd}[column sep = 12em, row sep = 5em]
( p : \prodop_n \ind{A} )
\arrow[phantom]{dr}[description, yshift=-2mm]
	{\twocell{\altCell_{(x_i)_{i=1, \,\dots\, ,n}}}}
\arrow[phantom]{r}[description]
	{\twocellIso{\psi^{\usem{}{\inc} \circ \inc\sem{-}}}}
\arrow[swap]{d}{\simeq}
\arrow[bend right=12, swap]{r}
	{\left( p : \prod_n \ind{A} \vdash 
		\pair{\hpi{1}{p} ,\dots, \hpi{n}{p})} : \prod_n \ind{A} \right)}
\arrow[bend left=12]{r}
	{\left( p : \prod_n \ind{A} \vdash p : \prod_n \ind{A}\right)} 
&
( p : \prodop_n \ind{A} ) 
\arrow{d}{\simeq} \\
\Gamma 
\arrow[swap]{r}{(\Gamma \vdash x_i : A_i)_{i=1, \,\dots\, , n}} &
\Gamma
\end{tikzcd}
\: = \: 
\begin{tikzcd}[column sep = 4.5em, row sep = 5em]
( p : \prodop_n \ind{A} )
\arrow{dr}[description]{\simeq}
\arrow[swap]{d}{\simeq}
\arrow[bend left=10]{r}
	{\left( p : \prod_n \ind{A} \vdash p : \prod_n \ind{A}\right)} &
( p : \prodop_n \ind{A} ) 
\arrow{d}{\simeq}
\arrow[phantom]{dl}[near start, description]{\iso} \\
\Gamma 
\arrow[phantom]{ur}[near start, description]{\iso}
\arrow[swap]{r}{(\Gamma \vdash x_i : A_i)_{i=1, \,\dots\, , n}} &
\Gamma
\end{tikzcd}
\end{equation*}
Applying the definition of $\psi^{\inc\sem{-}}$ given in 
Proposition~\ref{propn:weak-universal-property}, this entails checking the 
outer 
edges of the following diagram commute for $k=1, \,\dots\, , n$:
\begin{td}[row sep = 4em, column sep = 4em]
\hpi{k}{p} 
\arrow{rr}{\subid{\pi_k(p)}^{-1}}
\arrow[swap]{d}{\hpi{k}{\etaTimes{p}}} 
\arrow[equals, bend left = 20]{dr}
\arrow[phantom]{dr}[description, xshift=-2mm, yshift=-2mm]
	{\equals{triang. law}} &
\: &
\pi_k(p) 
\arrow[phantom]{ddl}[description]{=}
\arrow{dd}{\indproj{-k}{\pi_k(p)}} \\

\hpi{k}{\pair{\hpi{1}{p} ,\dots, \hpi{n}{p}}} 
\arrow[swap]{d}
	{\hpi{k}{\pair{\subid{\pi_1(p)}^{-1}, \,\dots\, , \subid{\pi_n(p)}^{-1}}}} 
\arrow[swap]{r}{\epsilonTimesInd{k}{\hcomp{\ind{\pi}}{p}}}
\arrow[phantom]{dr}[description]{\equals{nat.}} &
\hpi{k}{p} 
\arrow{d}[description]{\subid{\pi_k(p)}^{-1}} &
\: \\

\hpi{k}{\pair{\pi_1(p), \,\dots\, , \pi_n(p)}} 
\arrow[swap]{r}{\epsilonTimesInd{k}{\ind{\pi}(p)}} &
\pi_k(p) 
\arrow[swap]{r}{\indproj{-k}{\pi_k(p)}} &
\hcomp{x_k}{x_i \mapsto \pi_i(p)} 
\end{td}
Hence, the unit law does indeed hold. The associativity law 
holds by construction 
for composites of terms in unary contexts. For the general case, one 
instantiates the definition of $\phi^{\inc\sem{-}}$ from 
Proposition~\ref{propn:weak-universal-property} and applies the definition of 
$\postName$ to get exactly the required composite. This completes the proof of 
the next lemma.

\begin{prooflesslemma} 
For any unary $\langCartClosed$-signature $\sig$, the composite
$\ext{\inc} \circ \inc\sem{-} : \syncloneAtClosed{\sig} \to 
\syncloneAtClosed{\sig}$ 
induced by the following diagram is equivalent to 
$\id_{\syncloneAtClosed{\sig}}$:
\begin{equation} \label{eq:usemh-inc-sem-typing}
\begin{tikzcd}
\syncloneAtClosed{\sig}
\arrow{r}{\inc\sem{-}} &
\freeCartClosedBicat{\sig} 
\arrow{r}{\usem{}{\inc}} &
\syncloneAtClosed{\sig} \\
\sig
\arrow[equals]{r}
\arrow[hookrightarrow]{u} &
\sig 
\arrow[equals]{r}
\arrow[hookrightarrow]{u} &
\sig
\arrow[hookrightarrow]{u}
\end{tikzcd}
\end{equation}
\end{prooflesslemma}

Putting this lemma together with
Proposition~\ref{propn:first-biequivalence}, 
one obtains the biequivalence between 
$\syncloneAtClosed{\sig}$ and $\freeCartClosedBicat{\sig}$:

\begin{prooflesspropn} \label{prop:ccc-syntactic-models-biequivalent}
For any unary $\langCartClosed$-signature $\sig$, the 
cc-pseudofunctors $\inc\sem{-}$ and $\usem{}{\inc}$ extending the inclusion as 
in the diagram
\begin{td}
\freeCartClosedBicat{\sig} 
\arrow{r}{\ext{\inc}} &
\syncloneAtClosed{\sig}
\arrow{r}{\inc\sem{-}} &
\freeCartClosedBicat{\sig} \\

\sig 
\arrow[equals]{r}
\arrow[hookrightarrow]{u}{\inc} &
\sig 
\arrow[equals]{r}
\arrow[hookrightarrow]{u}{\inc} &
\sig 
\arrow[hookrightarrow]{u}{\inc}
\end{td}
form a biequivalence 
$\freeCartClosedBicat{\sig} \simeq \syncloneAtClosed{\sig}$.
\end{prooflesspropn}

It is not hard to see that the pseudonatural transformation 
$(\altNat, \altCell)$ defined in Construction~\ref{constr:altcell-defined} 
restricts to a pseudonatural transformation 
$\inc\sem{-} \circ \usem{}{\inc} \simeq 
\id_{\increstrict{\syncloneCartClosed{\sig}}}$ 
for $\inc\sem{-}$ the restriction of the interpretation pseudofunctor
of Proposition~\ref{propn:weak-universal-property} to 
$\increstrict{\syncloneCartClosed{\sig}}$. Since the proof of 
Proposition~\ref{propn:first-biequivalence} also restricts to the unary case, 
one 
obtains the following. 

\begin{prooflesscor}  \label{cor:uniqueness-up-to-bieq-of-restricted-contexts}
For any unary $\langCartClosed$-signature $\sig$, the 
cc-pseudofunctors $\inc\sem{-}$ and $\usem{}{\inc}$ extending the inclusion as 
in the diagram
\begin{td}
\freeCartClosedBicat{\sig} 
\arrow{r}{\ext{\inc}} &
\increstrict{\syncloneCartClosed{\sig}}
\arrow{r}{\inc\sem{-}} &
\freeCartClosedBicat{\sig} \\

\sig 
\arrow[equals]{r}
\arrow[hookrightarrow]{u}{\inc} &
\sig 
\arrow[equals]{r}
\arrow[hookrightarrow]{u}{\inc} &
\sig 
\arrow[hookrightarrow]{u}{\inc}
\end{td}
form a biequivalence 
$\freeCartClosedBicat{\sig} \simeq \increstrict{\syncloneCartClosed{\sig}}$. 
\end{prooflesscor}

Hence, up to canonical biequivalence, the syntactic model of $\langCartClosed(\sig)$ is 
the free cc-bicategory on the $\langCartClosed$-signature $\sig$. 
We are 
therefore justified in calling $\langCartClosed$ the internal language of 
cartesian closed bicategories. 

It further follows 
that the canonical pseudofunctor is unique up to equivalence.

\begin{mycor} \label{cor:uniqueness-up-to-equiv-of-ext-pseudofun}
For any cc-bicategory $\ccBicat{\baseCat}$, unary $\langCartClosed$-signature 
$\sig$ and \mbox{$\langCartClosed$-signature} homomorphism
$h : \sig \to \baseCat$, there exists a strict cc-pseudofunctor
$h\sem{-} : \increstrict{\syncloneCartClosed{\sig}} \to \baseCat$. Up to 
equivalence, this is the unique strict cc-pseudofunctor 
$F : \increstrict{\syncloneCartClosed{\sig}} \to \baseCat$ such that
$F \circ \inc = h$, for $\inc$ the inclusion.
\begin{proof}
Existence is Corollary~\ref{cor:weak-ump-for-restricted-contexts} so it 
suffices to show uniqueness. To this 
end, consider the diagram
\begin{td}
\freeCartClosedBicat{\sig} 
\arrow{r}{\ext{\inc}} &
\increstrict{\syncloneCartClosed{\sig}}
\arrow{r}{F} &
\baseCat \\
\: &
\sig
\arrow[hookrightarrow]{u}{\inc}
\arrow[hookrightarrow]{ul}{\inc}
\arrow{ur}[swap]{h} &
\:
\end{td}
where $F$ is any strict cc-pseudofunctor. By the free property of 
$\freeCartClosedBicat{\sig}$ (Lemma~\ref{lem:free-cc-bicat-on-sig}), 
$\ext{h} = F \circ \ext{\inc}$. Then, applying 
Corollary~\ref{cor:uniqueness-up-to-bieq-of-restricted-contexts}, one sees that
\[
F 	\simeq F \circ (\ext{\inc} \circ \inc\sem{-})
	\simeq (F \circ \ext{\inc}) \circ \inc\sem{-}
	= \ext{h} \circ \inc\sem{-}
\]
It follows that any strict cc-pseudofunctor extending $h$ is equivalent to 
$\ext{h} \circ \inc\sem{-}$. Hence, $h\sem{-}$ is unique up to 
equivalence.
\end{proof}
\end{mycor} 


We finish this section with a corollary relating the semantic interpretation of 
Proposition~\ref{prop:semantic-interpretation} to the free property 
of the free cc-bicategory (Lemma~\ref{lem:free-cc-bicat-on-sig}). 

\begin{mycor} \label{cor:interpretations-related}
For any cc-bicategory 
$\ccBicat\bicatX$,
set of base types $\baseTypes$, 
and \mbox{$\langCartClosed$-signature} homomorphism 
$h : \sig \to \bicatX$,
there exists an equivalence 
$\ext{h} \circ \inc\sem{-} \simeq h\sem{-} 
	: \syncloneAtClosed{\allTypes\baseTypes} \to \bicatX$. 
\begin{proof}
Observe that the composite
$
\allTypes\baseTypes
	\hookrightarrow
\freeCartClosedBicat{\allTypes\baseTypes} 
	\xra{\ext{\inc}}
\syncloneAtClosed{\allTypes\baseTypes}
	\xra{h\sem{-}}
\bicatX
$
is equal to simply $h$. Thus, applying
Lemma~\ref{lem:syntactic-model-uniqueness-up-to-equivalence}, there exists an 
equivalence $\ext{h} \simeq h\sem{-} \circ \ext{\inc}$. But by 
Proposition~\ref{prop:ccc-syntactic-models-biequivalent} there also exists an 
equivalence 
$\ext{\inc} \circ \inc\sem{-} 
\simeq 
\id_{{\freeCartClosedBicat{\allTypes\baseTypes}}}$. Hence,
\[
\ext{h} \circ \inc\sem{-} 
	\simeq
(h\sem{-} \circ \ext{\inc}) \circ \inc\sem{-}
	\simeq
h\sem{-}
\]
as claimed.
\end{proof}
\end{mycor}

\section{\texorpdfstring{Normal forms in $\langCartClosed$}{Normal forms}} 
\label{sec:STLC-vs-pseudoSTLC}

In this final section we shall make precise the sense in which 
$\langCartClosed$ is 
the simply-typed lambda calculus `up to isomorphism', which will enable us to 
port the notion of (long-$\beta\eta$) normal form from the simply-typed lambda 
calculus into 
$\langCartClosed$. Our approach is to extend 
the mappings defined in Section~\ref{sec:coherence-for-biclones} for 
$\langBiclone$ to include cartesian closed structure. One could go further, and 
prove that the syntactic model of $\langCartClosed$ is biequivalent to the 
syntactic model of the strict language $\hir$ extended with pseudo cartesian 
closed structure. Such a result provides a constructive proof that the free 
cartesian closed bicategory on a $\langCartClosed$-signature $\sig$ is 
biequivalent to the free 2-category with bicategorical products and 
exponentials on $\sig$. Since this follows from the Mac~Lane-Par{\'e} coherence 
theorem~\cite{MacLane1985}, together with fact that biequivalences 
preserve bilimits and 
biadjunctions, we restrict ourselves to mappings on terms. However, we shall 
present certain results one requires in order to construct this biequivalence, 
as they turn out 
to be of importance in the proof of our main theorem in Chapter~\ref{chap:nbe}. 

To fix notation, let $\stlc(\sig)$ denote the simply-typed lambda calculus 
with constants and base types specified by a $\stlc$-signature 
$\sig = (\baseTypes, \graph)$. This is 
defined in Figure~\ref{r:stlc-rules} below. As for $\langCartClosed$, we 
present products in an $n$-ary style which is equivalent to the usual 
presentation in terms of binary 
products and a terminal object. The equational theory is the usual 
$\alpha\beta\eta$-equality 
for the simply-typed lambda calculus (\eg~\cite{Barendregt1985, 
Crole1994}).

\setlength{\floatsep}{5pt plus 1.0pt minus 2.0pt} 

\begin{figure*}[!h]
{\small
\begin{mdframed}
\centering

\unaryRule	{\faketext}
			{x_1 : A_1, \dots, x_n : A_n \vdash x_k : A_k}
			{var $(1 \leq k \leq n)$}
			
\binaryRule	{c \in \graph(A_1, \dots, A_n; B)}
			{(\Delta \vdash	u_i : A_i)_{i=1,\dots,n}}
			{\Delta \vdash c(u_1, \dots, u_n) : B}
			{const}

\unaryRule	{\Gamma \vdash t_1 : A_1 \quad \dots \quad \Gamma \vdash t_n : A_n}
			{\Gamma \vdash \seq{t_1, \dots, t_n} : \prodop_n (A_1, \dots, A_n)}
			{$n$-tuple}
\quad
\unaryRule	{\Gamma \vdash t : \prodop_n(A_1, \dots, A_n)}
			{\Gamma \vdash \pi_k(t) : A_k}
			{$k$-proj ($1 \leq k \leq n$)}
			
\unaryRule 	{\Gamma, x : A \vdash t : B} 
			{\Gamma \vdash \lam{x}{t} : \exptype{A}{B}} 
			{lam} 
\qquad
\binaryRule	{\Gamma \vdash t : \exptype{A}{B}} 
			{\Gamma \vdash u : A}
			{\Gamma \vdash \app{t}{u} : B} 
			{app} \vspace{-\treeskip}

\caption{Rules for $\stlc(\sig)$. \label{r:stlc-rules}}

\end{mdframed}
}

\end{figure*}

%

We shall not distinguish notationally between the 
type theory $\stlc$ (resp. $\langCartClosed$) and its set of terms (or set of 
terms and rewrites) up to $\alpha$-equivalence. 
We employ the following 
notation: 
\begin{align*}
\stlc(\sig)(\Gamma; B) &:= 
	\slice{\{ t \st \Gamma \vdash_{\text{STLC}} t : B  \}\:}{\:\aeq} 
\\
\langCartClosed(\sig)(\Gamma; B) 
	&:= \slice{\{ t \st \Gamma \vdash_{\langCartClosed} t : B  \}\:}{\:\aeq} 
\end{align*}
Similarly, we write $\stlc(\sig)$ to denote the set 
of all $\stlc$-terms modulo $\alpha$-equivalence, and 
$\langCartClosed(\sig)$ to denote the set of all 
$\langCartClosed$-terms modulo $\alpha$-equivalence. (Precisely, these are 
sets indexed by (context, type) pairs.) 
We drop the 
decorations on the turnstile symbol unless the type theory in question is 
ambiguous.


\paragraph*{Relating $\langCartClosed$ and $\stlc$.} We define a pair of maps 
${\into{-} : 
\stlc(\sig) \leftrightarrows \langCartClosed(\sig) : \out{({-})}}$ for a fixed
$\stlc$-signature 
$\sig$. 
These maps extend those constructed 
in Section~\ref{sec:coherence-for-biclones} for biclones; indeed, 
the terms of $\hir(\sig)$ are exactly the 
variables and constants in $\stlc(\sig)$.

\begin{myconstr} \label{constr:cart-closed-out}
For any $\stlc$-signature $\sig$, define a mapping  
$\overline{(-)}: \langCartClosed(\sig) \to \stlc(\sig)$ as follows: 
\begin{align*}
\begin{aligned}[c]
\out{x_i} &:= x_i \\
\out{\pi_k(p)} &:= \pi_k(p) \\
\out{\evalterm(f,a)} &:= \app{f}{a}
\end{aligned}
\qquad\qquad\qquad
\begin{aligned}[c]
\out{c(x_1, \,\dots\, , x_n)} &:= c(x_1, \,\dots\, , x_n) \\
\out{\pair{t_1, \,\dots\, , t_n}} &:= \seq{\out{t_1}, \,\dots\, , \out{t_n}} \\
\out{\lam{x}{t}} &:= \lam{x}{\out{t}}
\end{aligned}
\end{align*} 
\end{myconstr}

It is elementary to check this definition respects $\alpha$-equivalence and the 
equational theory~$\equiv$.

\begin{prooflesslemma} \label{lem:out-respects-typing}
For any $\stlc$-signature $\sig$,
\begin{enumerate}
\item For all derivable terms $t, t'$ in $\langCartClosed(\sig)$, if $t \aeq 
t'$ then $\out{t} \aeq \out{t'}$, 
\item If $\Gamma \vdash t : B$ in $\langCartClosed(\sig)$ then $\Gamma \vdash 
\out{t} : B$ in $\stlc(\sig)$,~\ie~one obtains maps of indexed sets. \qedhere
\end{enumerate}
\end{prooflesslemma}

As we did for biclones, we think of $\out{t}$ as the \emph{strictification} 
of a term in $\langCartClosed$. The map $\into{-}$ 
interprets $\stlc$-terms in $\langCartClosed$.

\begin{myconstr} \label{constr:ccc:into-terms}
For any $\stlc$-signature $\sig$, define a mapping 
$\into{-} : \stlc(\sig) \to \langCartClosed(\sig)$ as follows:
\begin{align*}
\begin{aligned}[c]
\into{x_k} &:= x_k \\
\into{\pi_k(t)} &:= \hcomp{\pi_k}{\into{t}} \\
\into{\seq{t_1, \,\dots\, , t_n}} &:= \pair{\into{t_1}, \,\dots\, , \into{t_n}} 
\end{aligned}
\qquad\qquad\qquad
\begin{aligned}[c]
\into{c(u_1, \,\dots\, , u_n)} &:= \hcomp{c}{\into{u_1}, \,\dots\, , 
\into{u_n}} \\
\into{\app{t}{u}} &:= \hcomp{\evalterm}{\into{t}, \into{u}} \\
\into{\lam{x}{t}} &:= \lam{x}{\into{t}}
\end{aligned}
\end{align*} 
\end{myconstr}

This mapping also respects typing and $\alpha$-equivalence. 

\newpage
\begin{prooflesslemma}  \label{lem:into-respects-typing} 
For any $\stlc$-signature $\sig$,
\begin{enumerate}
\item For all derivable terms $t, t'$ in $\stlc(\sig)$, if $t \aeq t'$ then 
$\into{t} \aeq \into{t'}$, 
\item If $\Gamma \vdash t : B$ in $\stlc(\sig)$ then $\Gamma \vdash 
\into{t} : B$ in $\langCartClosed(\sig)$,~\ie~one obtains maps of indexed 
sets. \qedhere
\end{enumerate}
\end{prooflesslemma}

As in Section~\ref{sec:coherence-for-biclones}, strictifying a $\stlc$-term 
does nothing. 

\begin{mylemma} \label{lem:out-into-is-id}
The composite mapping $\out{({-})} \circ \into{-}$ is exactly the 
identity on $\stlc(\sig)$. 
\begin{proof}
The claim holds by induction, using the usual laws of capture-avoiding 
substitution for the simply-typed lambda calculus:
\begin{gather*}
x_k \mapsto x_k \mapsto x_k \\
c(u_1, \,\dots\, , u_n) \mapsto \hcomp{c}{\into{u_1}, \,\dots\, , \into{u_n}} 
\mapsto c(x_1, \,\dots\, , x_n)[\out{\into{u_i}}/x_i] \\[8pt]
\pi_k(t) \mapsto \hpi{k}{\into{t}} \mapsto \pi_k(p)[\out{\into{t}} / p] \\
\seq{t_1, \,\dots\, , t_n} \mapsto \pair{\into{t_1}, \,\dots\, , \into{t_n}} 
\mapsto \seq{\out{\into{t_1}}, \,\dots\, , \out{\into{t_n}}} \\[8pt]
\app{t}{u} \mapsto \heval{\into{t}, \into{u}} \mapsto (\app{f}{a})[\out{\into{t}} / f, \out{\into{u}} / a] \\
\lam{x}{t} \mapsto \lam{x}{\into{t}} \mapsto \lam{x}{\out{\into{t}}}
\end{gather*}
\end{proof}
\end{mylemma}

We shall require a rewrite reducing explicit 
substitutions to the meta-operation of capture-avoiding substitution.
As in 
the biclone case, this is the extra data required to 
make $\into{-}$ into a pseudofunctor.  Unlike the biclone case, however, we 
must now deal with variable binding. 
This entails an extra step in our 
construction.
To inductively prove a lemma about substitution in the simply-typed lambda 
calculus, it is 
common to first prove a lemma about weakening. This auxiliary result allows one 
to deal with the fresh variable appearing in the lambda abstraction step. We 
shall do something similar. First, we shall define a rewrite reducing context 
renamings (in particular, weakenings) to actual syntactic substitutions. Then, 
we shall use this to construct our rewrite handling arbitrary substitutions.

We call the auxiliary rewrite $\contName$ for 
\emph{context renaming}.

\begin{myconstr} \label{constr:def-of-cont}
For any $\stlc$-signature $\sig$ and context renaming $r$, we construct a 
rewrite 
$\cont{t}{r}$ 
making the following rule admissible:
\begin{center}
\binaryRule{\Gamma \vdash \into{t} : B}{r : \Gamma \to \Delta}{\Delta \vdash \cont{t}{r} : \rewrite{\hcomp{\into{t}}{x_i \mapsto r(x_i)}}{\into{t[r(x_i)/x_i]}} : B}{} \vspace{-\treeskip}
\end{center}
The definition is by induction on the derivation of $t$:
\small
\begin{gather*}
\cont{x_k}{r} := \hcomp{x_k}{x_i \mapsto r(x_i)} 
	\XRA{\indproj{r(x_i)}{}} \into{r(x_i)} \\
\cont{c(\ind{u})}{r} 
	:= \hcompthree{c}{\into{u_1}, \,\dots\, , \into{u_n}}{r} 
	\XRA{\assoc{}} 
	\hcomp{c}{\hcomp{\into{\ind{u}}}{r}} 
	\XRA{\hcomp{c}{\contName, \,\dots\, , \contName}} 
	\hcomp{c}{\into{\ind{u}[r(x_i) / x_i]}} \\
\cont{\pi_k(t)}{r} := 
	\hcomp{\hpi{k}{\into{t}}}{r} 
	\XRA{\assoc{}} 
	\hpi{k}{\hcomp{\into{t}}{r}} 
	\XRA{\hpi{k}{\contName}} 
	\hpi{k}{\into{t[r(x_i)/x_i]}} \\
\cont{\seq{t_1, \,\dots\, , t_n}}{\ind{u}} := 
	\hcomp{\pair{\into{t_1}, \,\dots\, , \into{t_n}}}{\into{\ind{u}}} 
	\XRA{\postName } \pair{\hcomp{\into{\ind{t}}}{\into{\ind{u}}}} 
	\XRA{\pair{\contName, \,\dots\, , \contName}} 
	\pair{\into{\ind{t}[u_i/x_i]}} 
\end{gather*}
\begin{align*}
\cont{\app{t}{u}}{r} 
	:= 
\hcompbig{\heval{\into{t}, \into{u}}}{r} 
	&\XRA{\assoc{}} 
\hcompbig
	{\evalterm}
	{\hcomp
		{\into{t}}
		{r}, 
	 \hcomp
	 	{\into{u}}
	 	{r}
	} \\ 
	&\XRA{\heval{\contName, \contName}} 
	\hcompbig{\evalterm}{\into{t[r(x_i)/x_i]}, {\into{u[r(x_i)/x_i]}}}  
\end{align*}
\begin{align*}
\cont{\lam{x}{t}}{r} := \hcomp{\left(\lam{x}{\into{t}}\right)}{r} 
		&\XRA{\pushName} 
		\lam{x}{\hcomp{\into{t}}{x \mapsto x, x_i \mapsto \wkn{r(x_i)}{x}}} \\ 
		&\XRA{\lam{x}{\hcomp{\into{t}}{x, \cont{r(x_i)}{\mathrm{inc}_x}}}} 
		\lam{x}{\hcomp{\into{t}}{x \mapsto x, x_i \mapsto r(x_i)}} 
		\\ 					
		&\XRA{\lam{x}{\contName}} \lam{x}{\into{t[x / x, r(x_i)/x_i]}}
\hfill\qedhere
\end{align*}
\normalsize
\end{myconstr}

We can now define $\subName$. The construction extends its biclone counterpart,
Construction~\ref{constr:def-of-sub}.

\begin{myconstr}
For any $\stlc$-signature $\sig$, we construct a rewrite $\sub{t}{\ind{u}}$ so 
that the following rule is admissible:
\begin{center}
\binaryRule{x_1 : A_1, \,\dots\, , x_n : A_n \vdash \into{t} : B}{(\Delta 
\vdash \into{u_i} : A_i)_{i=1,\dots, n}}{\Delta \vdash \sub{t}{\ind{u}} : 
\rewrite{\hcomp{\into{t}}{x_i \mapsto \into{u_i}}}{\into{t[u_i/x_i]}} : B}{} \vspace{-\treeskip}
\end{center}
The definition is by induction on the derivation of $t$:
\small
\begin{gather*}
\sub{x_k}{\ind{u}} := 
	\hcomp{x_k}{x_i \mapsto \into{u_i}} \XRA{\indproj{k}{}} \into{u_k} \\
\sub{c(\ind{u})}{\ind{v}} := 
	\hcompthree{c}{\into{u_1}, \,\dots\, , \into{u_n}}{\into{\ind{v}}} 
	\XRA{\assoc{}} 
	\hcomp{c}{\hcomp{\into{\ind{u}}}{\into{\ind{v}}}} 
	\XRA{\hcomp{c}{\subName, \,\dots\, , \subName}} 
	\hcomp{c}{\into{\ind{u}[v_j/y_j]}} \\
\sub{\pi_k(t)}{\ind{u}} := 
	\hcomp{\hpi{k}{\into{t}}}{\into{\ind{u}}} 
	\XRA{\assoc{}} 
	\hpi{k}{\hcomp{\into{t}}{\into{\ind{u}}}} 
	\XRA{\hpi{k}{\subName}} 
	\hpi{k}{\into{t[u_i/x_i]}} \\
\sub{\seq{t_1, \,\dots\, , t_n}}{\ind{u}} := 
	\hcomp{\pair{\into{t_1}, \,\dots\, , \into{t_n}}}{\into{\ind{u}}} 
	\XRA{\postName } 
	\pairlr{\hcomp{\into{\ind{t}}}{\into{\ind{u}}}} 
	\XRA{\pair{\subName, \,\dots\, , \subName}} 
	\pairlr{\into{\ind{t}[u_i/x_i]}}
\end{gather*}
\begin{align*}
\sub{\app{t}{u}}{\ind{v}} := 
\hcomp{\heval{\into{t}, \into{u}}}{\into{\ind{v}}} 
	&\XRA{\assoc{}} 
\hcompbig
	{\evalterm}
	{\hcomp
		{\into{t}}
		{\into{\ind{v}}}, 
	 \hcomp{\into{u}}{\into{\ind{v}}}
	 } \\ 
	&\XRA{\heval{\subName, \subName}} 
\hcompbig
	{\evalterm}
	{\into{t[v_j/y_j]} {\into{u[v_j/y_j]}}}  \\[10pt]
\sub{\lam{x}{t}}{\ind{u}} := 
	\hcomp{\left(\lam{x}{\into{t}}\right)}{\into{\ind{v}}} 	
	&\XRA{\pushName} 
\lam{x}{\hcompbig{\into{t}}{x, \wkn{\into{u}}{x}}} \\ 
	&\XRA{\lam{x}{\hcomp{\into{t}}{x, \cont{u}{\mathrm{inc}_x}}}}
\lam{x}{\hcomp{\into{t}}{x, \into{u}}} \\ 					
	&\XRA{\lam{x}{\subName}} 
\lam{x}{\into{t[x / x, u_i/x_i]}}
\end{align*}
\normalsize
\end{myconstr}

Note the use of $\contName$ in the lambda abstraction step. As one would 
expect, $\subName$ and $\contName$ coincide where the terms being 
substituted are all variables. 

\begin{mylemma} \label{lem:out-cont-coincide}
For any $\stlc$-signature $\sig$, 
judgement
$(\Gamma \vdash \into{t} : B)$ in $\langCartClosed(\sig)$,
and context renaming
$r : \Gamma \to \Delta$, then 
\[
\Delta \vdash \sub{t}{r(\ind{x})} \equiv \cont{t}{r} : 
	\hcomp{\into{t}}{x_i \mapsto r(x_i)} \To \into{t} : B
\]
\begin{proof}
By induction on the derivation of $t$: comparing the cases one-by-one, the equality is immediate.
\end{proof}
\end{mylemma}

Let us note some of other the ways in which $\contName$ and $\subName$ behave 
as expected~(\cf~Lemma~\ref{lem:pseudofun-laws-for-sub}). We shall not need 
these results immediately, but they will play an 
important role in the normalisation-by-evaluation proof of 
Chapter~\ref{chap:nbe}. 

\begin{mylemma} \label{lem:properties-of-cont-and-sub} \quad
For any $\stlc$-signature $\sig$  and any contexts 
$\Gamma := (x_i : A_i)_{i=1, \,\dots\, ,n}$ and
$\Delta := (y_j : B_j)_{j=1,\dots,m}$,
{\tikzcdset{arrow style=tikz, arrows={Rightarrow}}
\begin{enumerate}
\item \label{c:unit-law} If $\Gamma \vdash \into{t} : B$ then
\begin{equation} \label{eq:sub-unit-law}
\begin{tikzcd}[column sep = 4em]
\into{t}
\arrow[swap]{d}
	{\subid{\into{t}}}
\arrow[bend left=10, equals]{dr} &
\: \\
\hcomp{\into{t}}{x_i \mapsto x_i}
\arrow[swap]{r}[yshift=-2mm]
	{\cont{t}{\id_\Gamma}}  &
\into{t[x_i /x_i]}
\end{tikzcd}
\end{equation}

\item \label{c:pre-assoc-law} If $\Gamma \vdash \into{t} : B$ and 
$(\Delta \vdash \into{u_i} : A_i)_{i=1, \,\dots\, , n}$ then
\begin{equation} \label{eq:assoc-law-left-unit}
\begin{tikzcd}
\hcompthree{\into{t}}{x_i \mapsto \into{u_i}}{\id_\Delta} 
\arrow[swap]{d}{\hcomp{\sub{t}{\ind{u}}}{\id_\Delta}}
\arrow{r}{\assoc{}} &
\hcompbig{\into{t}}{x_i \mapsto \hcomp{\into{u_i}}{\id_\Delta}} 
\arrow{r}[yshift=2mm]{\hcomp{\into{t}}{\sub{u_i}{\id_\Delta}}} &
\hcomp{\into{t}}{x_i \mapsto \into{u_i}} 
\arrow{d}{\sub{t}{\ind{u}}} \\
\hcomp{\into{t[u_i / x_i]}}{\id_\Delta} 
\arrow[swap]{rr}{\sub{t[u_i / x_i]}{\id_\Delta}} &
\: &
\into{t[u_i / x_i]}
\end{tikzcd}
\end{equation}

\item \label{c:assoc-law} If $(\Gamma \vdash \into{t} : B)$, 
$(\Delta \vdash \into{u_i} : A_i)_{i=1, \,\dots\, , n}$ and 
$(\Sigma \vdash \into{v_j} : B_j)_{j=1, \,\dots\, , m}$, then
\begin{equation} \label{eq:sub-assoc-law}
\begin{tikzcd}
\hcompthree{\into{t}}{\into{\ind{u}}}{\into{\ind{v}}} 
\arrow[swap]{d}{\hcomp{\sub{t}{\ind{u}}}{\ind{v}}}
\arrow{r}{\assoc{}} &
\hcompbig{\into{t}}{\hcomp{\into{\ind{u}}}{\into{\ind{v}}}} 
\arrow{r}[yshift=2mm]{\hcomp{\into{t}}{\sub{u_i}{\ind{v}}}} &
\hcomp{\into{t}}{\into{\ind{u}[v_j/y_j]}} 
\arrow{d}{\sub{t}{\ind{u}}} \\
\hcomp{\into{t[u_i / x_i]}}{\into{\ind{v}}} 
\arrow[swap]{rr}{\sub{t[u_i / x_i]}{\ind{v}}} &
\: &
\into{t\left[u_i[v_j/y_j] / x_i\right]}
\end{tikzcd}
\end{equation}
\end{enumerate}
}
\vspace{\parskip}
\begin{proof}
{\tikzcdset{arrow style=tikz, arrows={Rightarrow}}
Each of the claims is proven by induction. Most of 
the cases for~(\ref{c:unit-law}) are almost immediate, except for lambda 
abstraction. There one 
uses Lemma~\ref{lem:properties-of-push}(\ref{c:push-and-subid}).

For~(\ref{c:pre-assoc-law}) and~(\ref{c:assoc-law}), all the cases except for 
lambda abstraction are relatively simple. One can 
prove~(\ref{c:assoc-law}) and derive~(\ref{c:pre-assoc-law}) as a special 
case. For lambda 
abstraction,~\ie~for judgements of the form
$(\Gamma \vdash t : \exptype{A}{B})$, one 
must deal with fresh variables. For this we take the claims in order. 

To prove 
the 
\rulename{lam} case of~(\ref{c:pre-assoc-law}) one first proves 
three further lemmas building towards the target result. The first is that 
whenever
$(\Delta \vdash \into{u_i} : A_i)$, then
\begin{equation} \label{eq:assoc-law-variables}
\begin{tikzcd}[column sep = 5em]
\hcompthree{\into{u_i}}{\id_\Delta}{\id_\Delta} 
\arrow[swap]{d}{\hcomp{\sub{t}{\id_\Delta}}{\id_\Delta}}
\arrow{r}{\assoc{}} &
\hcomp{\into{u_i}}{\hcomp{y_j}{\id_\Delta}} 
\arrow{r}[yshift=0mm]{\hcompsmall{\into{u_i}}{\indproj{\bullet}{\ind{y}}}} &
\hcomp{\into{u_i}}{\id_\Delta} 
\arrow{d}{\sub{u_i}{\ind{y}}} \\
\hcomp{\into{u_i}}{\id_\Delta} 
\arrow[swap]{rr}{\sub{t}{\id_\Delta}} &
\: &
\into{u_i}
\end{tikzcd}
\end{equation}
To show this diagram commutes, one inducts on the derivation of $\into{t}$; all 
the cases but \rulename{lam} follow as for~(\ref{c:assoc-law}). For the 
\rulename{lam} case one uses the inductive hypothesis, the coherence of 
$\langBiclone$, and 
Lemma~\ref{lem:properties-of-push}(\ref{c:push-and-assoc}). 

Next we show that, whenever 
$(\Gamma \vdash \into{t} : B)$ and
$(\Delta \vdash \into{u_i} : A_i)_{i=1, \,\dots\, , n}$, then
\begin{equation} \label{eq:assoc-law-terms-variables}
\begin{tikzcd}[column sep = 5em]
\hcompthree{\into{t}}{\into{\ind{u}}}{\id_\Delta} 
\arrow[swap]{d}{\hcomp{\sub{t}{\ind{u}}}{\id_\Delta}}
\arrow{r}{\assoc{}} &
\hcompbig{\into{t}}{x_i \mapsto \hcomp{\into{u_i}}{\id_\Delta}} 
\arrow{r}[yshift=2mm]{\hcomp{\into{t}}{\sub{u_i}{\id_\Delta}}} &
\hcomp{\into{t}}{\into{\ind{u}}} 
\arrow{d}{\sub{t}{\ind{u}}} \\
\hcomp{\into{t[u_i/x_i]}}{\id_\Delta} 
\arrow[swap]{rr}{\sub{t[u_i/x_i]}{\id_\Delta}} &
\: &
\into{t[u_i/x_i]}
\end{tikzcd}
\end{equation}
Once again all the cases but \rulename{lam} follow from the generality 
of~(\ref{c:assoc-law}). For the lambda abstraction case the proof is similar to 
that for~(\ref{eq:assoc-law-variables}): one applies the 
inductive hypothesis, 
Lemma~\ref{lem:properties-of-push}(\ref{c:push-and-assoc}) 
and~(\ref{eq:assoc-law-variables}).  

The final lemma required is the following. For any judgements
$(\Gamma \vdash \into{t} : B)$, 
$({\Delta \vdash \into{u_i} : A_i})_{i=1, \,\dots\, , n}$ and 
$(\Sigma, x : A \vdash \into{v_j} : B_j)_{j=1,\dots, m}$, 
one shows that
\begin{equation} \label{eq:assoc-law-free-variable}
\begin{tikzcd}[column sep = 5em]
\hcompthree{\into{t}}{\into{u_i}}{\into{\id_\Delta}}
\arrow[swap]{d}{\hcomp{\sub{t}{\ind{u}}}{\id_\Delta}}
\arrow{r}{\assoc{}} &
\hcomp{\into{t}}{x_i \mapsto \hcomp{\into{u_i}}{\id_\Delta}} 
\arrow{r}[yshift=0mm]{\hcomp{\into{t}}{\sub{u_i}{\id_\Delta}}} &
\hcompbig{\into{t}}{\into{\ind{u}}} 
\arrow{d}{\sub{t}{\ind{u}}} \\
\hcomp{\into{t[u_i/x_i]}}{\id_\Delta} 
\arrow[swap]{rr}{\sub{t[u_i/x_i]}{\id_\Delta}} &
\: &
\into{t[u_i/x_i]}
\end{tikzcd}
\end{equation}

We are finally in a position to prove the \rulename{lam} case 
of~(\ref{c:assoc-law}). Unwinding the clockwise route around the claim, one 
obtains the left-hand edge of Figure~\ref{fig:big-diagram} below~(page~\pageref{fig:big-diagram}), in which we 
abbreviate the term
\[
 \lam{x}{\hcomp{\into{t}^{\Gamma, x : A}}{
			\hcompbig
				{\wkn{\into{\ind{u}}}{x}^{\Delta, x : A}}
				{\wkn{\into{\ind{v}}}{x}^{\Sigma, x : A}, 
								x^{\Sigma, x : A}},
			\hcompbig
				{x^{\Delta, x : A}}
				{\wkn{\into{\ind{v}}}{x}^{\Sigma, x : A}, 
								x^{\Sigma, x : A}}
		} }
\]
by $\lam{x}{\hcomp{\into{t}}{(\ast)}}$ and write $\indproj{x}{\ind{u}, x}$ 
for the rewrite 
$\indproj{x}{\ind{u}, x} : \rewrite{\hcomp{x}{x_i \mapsto u_i, x \mapsto v}}
								{v}$ 
taking the projection at the variable $x$. One then unfolds the anticlockwise 
route and applies the inductive hypothesis to obtain the outer edge of 
Figure~\ref{fig:big-diagram}, completing the proof.
}
\end{proof}
\end{mylemma}

\paragraph*{STLC up to isomorphism.}

One approach in the field of game semantics is to quotient a (putative) 
cc-bicategory 
to obtain a cartesian closed category~(see~\eg~\cite[Chapter~2]{PaquetThesis}). 
Doing so 
loses intensional information, but makes calculations simpler. This suggests 
that one ought to be able to quotient $\langCartClosed$ (up to the existence of 
an invertible rewrite) to obtain $\stlc$ (up to $\beta\eta$-equality). 

We begin by making precise the sense in which the $\into{-}$ mapping respects 
$\beta\eta$-equality up to 
isomorphism. 

\begin{mylemma} \label{lem:existence-of-BE} \quad
Let $\sig$ be a $\stlc$-signature.
\begin{enumerate}
\item \label{c:out-well-defined} If $\Gamma \vdash \tau : \rewrite{t}{t'} : A$ 
in $\langCartClosed(\sig)$, then $\out{t} =_{\beta\eta} \out{t'}$. 
\item \label{c:BE-constructed} If $t =_{\beta\eta} t'$ for 
$t, t' \in \stlc(\sig)(\Gamma; A)$, then 
there exists a rewrite 
$\Gamma \vdash \BE{t,t'} : \rewrite{\into{t}}{\into{t'}} : A$ in 
$\langCartClosed(\sig)$. 
\end{enumerate}
\begin{proof}
For~(\ref{c:out-well-defined}) we induct on the derivation of $\tau$. For the 
structural rewrites and the identity the result is trivial, while for $\tau' 
\vert \tau$ it follows immediately from the inductive hypothesis. For 
$\epsilonTimesInd{k}{}$ one obtains 
$\out{\hpi{k}{\pair{t_1, \,\dots\, , t_n}}} 
	= \pi_k(\seq{\out{t_1}, \,\dots\, , \out{t_n}}) 
	=_{\beta\eta} \out{t_k}$, 
while for 
$\transTimes{\alpha_1, \,\dots\, , \alpha_n}$ one has 
$\out{u} 
	=_{\beta\eta} \seq{\pi_1(\out{u}), \,\dots\, , \pi_n(\out{u})} 
	\overset{\text{IH}}{=}_{\beta\eta} \seq{\out{t_1}, \,\dots\, , \out{t_n}}$. 
The cases for exponential structure are 
similar: for $\epsilonExpRewr{t}$ one sees that 
$\out{\heval{\wkn{(\lam{x}{t})}{x}, x}} 
	= \app{\lam{x}{\out{t}}}{x} 
	=_{\beta\eta} \out{t}$, while for 
$\transExp{x \bind \tau}$ one finds that
$\out{u} =_{\beta\eta} \lam{x}{\app{u}{x}} 
	\overset{\text{IH}}{=}_{\beta\eta} \lam{x}{\out{t}}$.
	
For~(\ref{c:BE-constructed}) we induct on the definition of $\beta\eta$-equality (\eg~\cite[Figure 4.2]{Crole1994}). 

\subproof{$\beta$-rules}
For the $\pi_k(\seq{t_1, \,\dots\, , t_n}) =_{\beta\eta} t_k$ rule one takes 
$\hpi{k}{\pair{\into{t_1}, \,\dots\, , \into{t_n}}}
\XRA{\epsilonTimesInd{k}{}} 
\into{t_k}$. 
For $\app{\lam{x}{t}}{u} =_{\beta\eta} t[u/x]$ one takes 
$\heval{\lam{x}{\into{t}}, \into{u}} 
\XRA{\genEpsilonExp{}} 
\hcomp{\into{t}}{\id_\Gamma, x \mapsto \into{u}} 
\XRA{\subName} 
\into{t[u/x]}$. 	

\subproof{$\eta$-rules}
In a similar fashion, for $t =_{\beta\eta} \seq{\pi_1(t), \,\dots\, , 
\pi_n(t)}$ one 
takes 
$\into{t} 
\XRA{\etaTimes{}} 
\pairlr{\hpi{1}{\into{t}}, \,\dots\, , \hpi{n}{\into{t}}}$ while for 
$t =_{\beta\eta} \lam{x}{\app{t}{x}}$ one takes 
$\into{t} 
\XRA{\etaExp{}} 
\lam{x}{\heval{\wkn{\into{t}}{x}, x}} 
\XRA{\lam{x}{\heval{\subName, x}}} 
\lam{x}{\heval{\into{t}, x}}$.

The rules for an equivalence relation hold by the categorical rules on vertical 
composition. The congruence rules hold by the functoriality of explicit 
substitution and the functoriality of the $\pair{-, \,\dots\, , =}$ and 
$\lam{x}{(-)}$ operations.
\end{proof}
\end{mylemma}

The preceding lemma motivates the following definition.

\begin{mydefn}
Fix a $\stlc$-signature $\sig$. 
For every context $\Gamma$ and type $A$, define an equivalence relation 
$\iso^\Gamma_A$ on $\langCartClosed(\sig)(\Gamma; A)$ by setting 
$t \iso^\Gamma_A t'$ if and only if there exists a (necessarily invertible) 
rewrite $\tau$ such that \mbox{$\Gamma \vdash \tau : \rewrite{t}{t'} : A$}.
\end{mydefn}

We can therefore rephrase Lemma~\ref{lem:existence-of-BE} as follows. For any 
pair of terms 
$t, t' \in \stlc(\Gamma;A)$ 
such that 
$t =_{\beta\eta} t'$, then $\into{t} \iso_A^\Gamma \into{t'}$; moreover, if 
$t \iso_A^\Gamma t'$ then $\out{t} =_{\beta\eta} \out{t'}$. To show that 
$\stlc(\sig)(\Gamma; A)$-terms modulo-$\beta\eta$ are in 
bijection with 
$\langCartClosed(\sig)(\Gamma; A)$-terms modulo-$\iso_A^\Gamma$, 
it remains to show how to \Def{reduce} a term of the form 
$\into{\out{t}}$ to the original term $t$.

%


\begin{myconstr}
Define an invertible rewrite $\reduceName$ with typing 
\begin{center}
\unaryRule	{\Gamma \vdash t : A}
{\Gamma \vdash \reduce{t} : \rewrite{t}{\into{\out{t}}} :  A}
{} \vspace{-\treeskip}
\end{center}
by extending Construction~\ref{constr:biclone-reduce} with the following rules:
\enlargethispage{2\baselineskip}
\begin{gather*}
\reduce{\pi_k(p)} := \pi_k(p) \XRA{\subid{}} \hpi{k}{p} \\
\reduce{\pair{t_1, \,\dots\, , t_n}} := \pair{t_1, \,\dots\, , t_n} 
\XRA{\pair{\reduceName, \,\dots\, , \reduceName}} \pair{\into{\out{t_1}}, 
\,\dots\, , \into{\out{t_n}}} \\
\reduce{\evalterm(f,x)} := \evalterm(f,x) \XRA{\subid{}} \heval{f,x} \\
\reduce{\lam{x}{t}} := \lam{x}{t} \XRA{\lam{x}{\reduce{t}}} \lam{x}{\into{\out{t}}}
\end{gather*}
\end{myconstr}

Thought of as syntax trees, the term $\into{\out{t}}$ is constructed
by evaluating explicit substitutions as far as possible and 
pushing them as far as possible to the left. 
The $\reduceName$ rewrites reach a fixpoint on terms of form $\into{\out{t}}$, 
thereby providing a notion of normalisation in the sense of 
abstract rewriting systems (\eg~\cite{Baader1998}).

\begin{mylemma}
For any $\stlc$-signature $\sig$ and any term $(\Gamma \vdash t : A)$ derivable 
in $\stlc(\sig)$, the judgement
$
\left(\Gamma \vdash \reduce{\into{t}} \equiv \id_{\into{t}} : 
\rewrite{\into{t}}{\into{t}} : A\right)
$
is derivable in $\langCartClosed(\sig)$. 
\begin{proof}
Induction on the structure of $t$.
\end{proof}
\end{mylemma}

We are now in a position to make precise the sense in which $\langCartClosed$ 
is $\stlc$ up to isomorphism.

\begin{mypropn} \label{prop:stlc-up-to-iso}
For any $\stlc$-signature $\sig$, the maps 
$\into{-} : \stlc(\sig) \leftrightarrows 
			\langCartClosed(\sig) : 
\out{({-})}$ 
descend to a bijection 
\[
\slice{\stlc(\sig)(\Gamma; A)}{\beta\eta} \iso 
\slice{\langCartClosed(\sig)(\Gamma; 
A)}{\iso^\Gamma_A }
\]
\normalsize
between $\alpha\beta\eta$-equivalence classes of $\stlc(\sig)$-terms and 
$\alpha{\iso^\Gamma_A}$-equivalence classes of 
$\langCartClosed(\sig)$-terms.
\begin{proof}
The maps are well-defined on equivalence classes by 
Lemma~\ref{lem:existence-of-BE} and respect typing by 
Lemmas~\ref{lem:out-respects-typing} and~\ref{lem:into-respects-typing}, so it 
suffices to check the isomorphism. By Lemma~\ref{lem:out-into-is-id}, the 
composite $\out{({-})} \circ \into{-}$ is the identity. For the other 
composite, one needs to construct an invertible rewrite $\into{\out{t}} \iso t$ 
for every derivable term $t$: we take $\reduceName$. 
\end{proof}
\end{mypropn}

In particular, every typeable term $(\Gamma \vdash t : A)$ in 
$\langCartClosed(\sig)$ 
has 
a natural choice of \emph{normal form}, namely the long-$\beta\eta$ normal form 
(\eg~\cite{Huet1976}) of $\out{t}$ as an $\stlc$-term. 

\newpage
\begin{mycor} 
Let $\sig$ be a $\stlc$-signature. For any derivable term $\Gamma \vdash 
t : B$ in $\langCartClosed(\sig)$, there exists a unique long-$\beta\eta$ 
normal form term $N$ in $\stlc(\sig)$ such that $t \iso^\Gamma_B 
\into{N}$ and $\reduce{\into{N}} \equiv \id_{\into{N}}$. 
\begin{proof}
We take $N$ to be the long-$\beta\eta$ normal form of $\out{t}$. Then 
$N =_{\beta\eta} \out{t}$ so, by Proposition~\ref{prop:stlc-up-to-iso}, 
\[
\into{N} \iso_B^\Gamma 
\into{\out{t}} \iso_B^\Gamma 
t
\] 
For uniqueness, suppose that $N$ 
and $N'$ are long-$\beta\eta$ normal terms such that $\into{N} \iso_B^\Gamma t 
\iso_B^\Gamma \into{N'}$. Then $\out{\into{N}} =_{\beta\eta} \out{\into{N'}}$, 
so that $N =_{\beta\eta} N'$, and hence $N = N'$ by the uniqueness of long 
$\beta\eta$-normal forms.
\end{proof}
\end{mycor}

We end this chapter by recording the bicategorical statement of the work in 
this section. 

\begin{prooflessthm}
Fix a unary $\langCartClosed$-signature $\sig$.
The mappings $\into{-}$ and $\out{(-)}$ extend to pseudofunctors between the 
free cartesian closed bicategory on $\sig$ and the free 2-category with 
bicategorical cartesian closed structure on $\sig$. Together with the 
pseudonatural transformation $(\Id, \reduceName)$, they form a biequivalence. 
\end{prooflessthm}

\newgeometry{margin=.2cm}
\begin{landscape}
\thispagestyle{empty}
\begin{figure}
\small
\centering
\begin{td}[column sep = 1.2em, row sep = 2.2em]
\hcompthree
	{\left(\lam{x}{\into{t}}\right)^{\Gamma}}
	{\into{\ind{u}}^{\Delta}}
	{\into{\ind{v}}^\Sigma} 
\arrow[swap]{dd}{\assoc{}}
\arrow{r}{\hcomp{\pushName}{\ind{v}}}
\arrow[phantom]{dddrr}[description]
	{\equals{Lemma~\ref{lem:properties-of-push}(\ref{c:push-and-assoc})}}
	&
\hcomp{\lam{x}{\hcomp{\into{t}^{\Gamma, x : A}}{
	\wkn{\into{\ind{u}}}{x}^{\Delta, x : A}, 
		x^{\Delta, x : A}}}}
	{\into{\ind{v}}^\Sigma} 
\arrow{r}[yshift=0mm]{\pushName} &
\lam{x}{\hcompthree{\into{t}^{\Gamma, x : A}}
		{\wkn{\into{\ind{u}}}{x}^{\Delta, x : A}, 
				x^{\Delta, x : A}}
		{\wkn{\into{\ind{v}}}{x}^{\Sigma, x : A}, 
				x^{\Sigma, x : A}}} 
\arrow{d}[description]{\lam{x}{\assoc{}}} \\

\: &
\: &
\lam{x}{\hcomp{\into{t}^{\Gamma, x : A}}{(\ast)}} 
\arrow[bend left = 90, phantom]{ddddd}
	[description, near end, ""{coordinate, name=Z}, xshift=14mm]{}
\arrow[
		swap,
		rounded corners,
		to path=
			{ -- ([xshift=2ex]\tikztostart.east)
			-| (Z) [near end]\tikztonodes
			|- ([xshift=2ex]\tikztotarget.east)
			-- (\tikztotarget)}
	]{ddddddl}{}
\arrow{d}[description]
	{\lam{x}{\hcomp{\into{t}}{
				\hcompsmall{\into{\ind{u}}}{\indproj{\bullet}{}} \vert 
				\assoc{}, 
				\indproj{x}{}}
			} }  \\

\hcomp
	{\left(\lam{x}{\into{t}}\right)^{\Gamma}}
	{\hcomp{\into{\ind{u}}}{\into{\ind{v}}}^{\Sigma}} 
\arrow[swap]{d}{\hcomp{(\lam{x}{\into{t}})}{\sub{u_i}{\ind{v}}}} 
\arrow{dr}[description]{\pushName} &
\: &
\lam{x}{\hcomp
			{\into{t}^{\Gamma, x : A}}
			{
				\hcomp{\into{\ind{u}}}{\wkn{\into{\ind{v}}}{x}}^{\Sigma, 
				x 			
				: A}, 
				 x^{\Sigma, x : A}}
		} 
\arrow[bend left = 10]{dddl}[description]{\lam{x}{
				\hcomp{\into{t}}
						{\hcomp{
							\into{\ind{u}}}
							{\sub{\into{\ind{v}}}{\mathrm{inc}_x}}, 
						\id_x
						}}} \\

\hcomp{(\lam{x}{\into{t}^{\Gamma, x : A}})}{\into{u_i[v_j/y_j]}^\Sigma}
\arrow[swap]{d}{\pushName} 
\arrow[phantom]{r}[description]{\equals{nat.}} &
\lam{x}{\hcomp{\into{t}^{\Gamma, x : A}}{
	\wkn{\hcomp{\into{\ind{u}}}{\into{\ind{v}}}}{x}^{\Sigma, x : A}, 
	 x^{\Sigma, x : A}}} 
\arrow{dl}[description]
	{\lam{x}{\hcomp{\into{t}}
		{\wkn{\sub{u_i}{\ind{v}}}{x}, \id_x}}}
\arrow{ur}[description]{\lam{x}{\hcomp{\into{t}}{\assoc{}, \id_x}}}
\arrow[phantom]{dd}[description]{\equals{(\ref{eq:assoc-law-terms-variables})}} 
&
\:
\\

\lam{x}
	{\hcomp{\into{t}^{\Gamma, x : A}}
		{\wkn{\into{u_i[v_j/y_j]}}{x}^{\Sigma, x : A}, 
			 x^{\Sigma, x : A}}} 
\arrow[swap]{dd}
	{\lam{x}{\hcomp{\into{t}}{\sub{u_i[v_j/y_j]}{\mathrm{inc}_x}, \id_x}}} &
\: &
\: \\

\: &
\lam{x}{\hcomp{\into{t}^{\Gamma, x : A}}
			{\hcomp{\into{\ind{u}}}{\into{\ind{v}}}^{\Sigma, x : A},  
			x^{\Sigma, x : A}}}
\arrow{dl}[description]{\lam{x}{\hcomp{\into{t}}
			{\sub{\ind{u}}{\ind{v}}, 
			\id_x}}} &
\: \\

\lam{x}{\hcomp{\into{t}^{\Gamma, x :A}}
	{\into{u_i[v_j/y_j]}^{\Sigma, x : A}, 
		 x^{\Sigma, x : A}}}
\arrow[swap]{ddd}{\lam{x}{\sub{t}{\ind{u}[v_j/y_j], x}}} &
\lam{x}{\hcomp{\into{t}}
		{\hcomp
			{\into{\ind{u}}}
			{\hcomp{\ind{y}}{\into{\ind{v}}, x}}^{\Sigma, x : A} , 
		x^{\Sigma, x : A}}}
\arrow{u}[swap]{\lam{x}{\hcomp{\into{t}}
					{\hcompsmall{\into{\ind{u}}}{\indproj{\bullet}{}} , 
							\id_x}}}
\arrow[phantom]{ur}[description, xshift=10mm]{\equals{nat.}} &
\: \\

\: &
\lam{x}{\hcomp{\into{t}^{\Gamma, x : A}}
			{\hcomp{\wkn{\into{\ind{u}}}{x}}{\into{\ind{v}}, x}^{\Sigma, x:A}, 
			x^{\Sigma, x :A}}}
\arrow[swap]{u}{\lam{x}{\hcomp{\into{t}}
			{\assoc{}, \id_x}}} 
\arrow{d}{\lam{x}{\hcompbig{\into{t}^{\Gamma, x : A}}
			{\hcomp{\sub{\into{\ind{u}}}{\mathrm{inc}_x}}
				{\into{\ind{v}}, x}, \id_x}}}
\arrow[phantom]{r}[yshift=-1em, font=\scriptsize, at end]
		{\lam{x}{
			\hcomp{\into{t}}{\hcomp
				{\wkn{\into{\ind{u}}}{x}}
				{\sub{\into{\ind{v}}}{\mathrm{inc}_x}, \id_x}, 
			\indproj{x}{}} }} &
\: \\
							
\: &
\lam{x}{\hcomp{\into{t}^{\Gamma, x : A}}
			{\hcomp{\into{\ind{u}}}{\into{\ind{v}}, x}^{\Sigma, x:A}, 
			x^{\Sigma, x: A} }}
\arrow[phantom]{uuul}
	[description, xshift=4mm]{\equals{(\ref{eq:assoc-law-variables})}}
\arrow[bend left = 15]{uul}[description, xshift=-4mm]{\lam{x}
	{\hcomp{\into{t}}{\sub{\ind{u}}{\into{\ind{v}}, x}}, \id_x}} \\

\lam{x}{\into{t\left[ u_i[v_j/y_j] / x_i \right]}^{\Sigma, x : A}} &
\: 
\end{td}
\normalsize
\caption{Diagram for the proof of 
Lemma~\ref{lem:properties-of-cont-and-sub}(\ref{c:assoc-law}) 
\label{fig:big-diagram}}
\end{figure}
\end{landscape}
\restoregeometry

\part{Glueing and normalisation-by-evaluation}
\label{part:nbe}

%

\chapter{Indexed categories as bicategorical presheaves} 
\label{chap:calculations}

Categories of (pre)sheaves are often useful as a kind of `completion', allowing 
one to employ extra structure 
that may not exist in the original category. The aim of 
this chapter is to show that 
bicategorical versions of some of these properties extend to the 
bicategory $\Hom(\baseCat, \Cat)$ of pseudofunctors from a bicategory 
$\baseCat$ to the 2-category $\Cat$. (Pseudofunctors $\op\baseCat \to \Cat$ are 
also called \Def{indexed categories}~\cite{MacLane1985}.) Recall 
that, since $\Cat$ is a 
2-category, so is $\Hom(\baseCat, \Cat)$, and that we write $\Cat$ for the 
2-category of small categories 
(Notation~\ref{not:size-issues}). 

Specifically, we shall prove three results 
which will be used in later chapters:
\begin{enumerate}
\item $\Hom(\baseCat, \Cat)$ has all small bilimits, which are given pointwise,
\item $\Hom(\baseCat, \Cat)$ is cartesian closed, and the value of the 
exponential 
$\altexp{P}{Q}$ 
at $X \in \baseCat$ can be taken to be 
$\Hom(\baseCat, \Cat)(\Yon X 
\times P, Q) : \baseCat \to \Cat$, for $\Yon X := \baseCat(X, -)$ the covariant 
Yoneda embedding, 
\item For any $X \in \baseCat$ the exponential $\altexp{\Yon X}{P}$ in 
$\Hom(\baseCat, \Cat)$ may be given by $P(- \times X)$. 
\end{enumerate}
The proofs are rather technical. The reader willing to take these three 
statements on trust---for example, by analogy with the case of presheaves---may 
safely skip this chapter. For reference, the cartesian closed structures we construct here are 
summarised in an appendix  
(Tables~\ref{table:ccc-structure-on-presheaves} 
and~\ref{table:exponentiating-by-representable}).

Our first result is that $\Hom(\baseCat, \Cat)$ is bicomplete. For brevity, we 
provide an abstract argument which 
relies on the notions of \Def{pseudolimit}~\cite{Street1980} and 
\Def{flexible limit}~\cite{Blackwell1981}. We 
will not use these concepts anywhere else, so do not delve into the details 
here: an excellent overview of the various forms of limit and their 
relationship is 
available in~\cite{Lack2010}. 

\newpage
\begin{mypropn} \label{prop:hom-basecat-cat-bicomplete}
For any bicategory $\baseCat$, the 2-category $\Hom(\baseCat, \Cat)$ is 
bicomplete, with bilimits given pointwise.
\begin{proof}
We may assume without loss of generality that $\baseCat$ is a 2-category.  
To see this is the case, observe that if 
$\catV \simeq \catV'$ are biequivalent bicategories then
$\Hom(\catV, \Cat) \simeq \Hom(\catV', \Cat)$ (see Lemma~\ref{lem:ccc-and-biequivalent-implies-ccc}), 
and hence 
$\Hom(\catV, \Cat)$ has all small bilimits if and only if
$\Hom(\catV', \Cat)$ does. By 
the coherence theorem for bicategories~\cite{MacLane1985} every bicategory is 
biequivalent to a 2-category, so the claim follows.

Now, by~\cite[Proposition~3.6]{Power1989bilimit} for any 2-category $\altCat$ 
the 2-category $\Hom(\altCat, \Cat)$ admits all flexible limits, calculated 
pointwise. The so-called 
`PIE limits' are flexible (\cite[Proposition~4.7]{Bird1989})
and suffice to construct all pseudolimits (\cite[Proposition~5.2]{Kelly1989}), 
so $\Hom(\baseCat, \Cat)$ has all pseudolimits. But, as explained 
in~\cite[\S 6.12]{Lack2010}, a 2-category with all pseudolimits has all 
bilimits, completing the proof. 
\end{proof}
\end{mypropn}

This result may also be obtained directly, in a manner similar to the 
categorical argument, as a corollary of the following proposition. 
We do not pursue the point any further here for reasons of space.

\begin{prooflesspropn} \label{prop:hom-bicategory-biuniversal-arrows-pointwise}
Let $F : \baseCat \to \catW$ and $D : \catV \to \catW$ 
($D$ for `diagram') be 
pseudofunctors equipped with a chosen biuniversal 
arrow $(LB, u_B : D(LB) \to FB)$ from $D$ to $FB$ for every $B \in \baseCat$. 
Then 
\begin{enumerate}
\item \label{c:pseudofunctor-pointwise} The mapping $L : ob(\baseCat) 
\to ob(\catV)$ extends canonically to a pseudofunctor 
$\baseCat \to \catV$, and
\item \label{c:biuniversal-arrow-pointwise} The biuniversal arrows $u_B$ are 
the components of a biuniversal arrow $DL \To F$
from $D \circ (-) : \Hom(\baseCat, \catV) \to \Hom(\baseCat, \catW)$ to 
$F$. \qedhere
\end{enumerate}
\end{prooflesspropn}

\section{\texorpdfstring{$\Hom(\baseCat, \Cat)$ is cartesian closed}{Cartesian 
closure}} 
\label{sec:presheaves-cartesian-closed}

It follows immediately from Proposition~\ref{prop:hom-basecat-cat-bicomplete} 
that,
for any bicategory $\baseCat$, the 2-category $\Hom(\baseCat, \Cat)$ has all 
finite products. In this section we confront the construction 
of exponentials. The usual Yoneda argument 
(see~\eg~\cite[\S 8.7]{awodey}), expressed bicategorically, gives us a 
canonical 
choice of exponential to check. For any pseudofunctors $P, Q : \baseCat \to 
\Cat$, putative exponential $\altexp{P}{Q}$ and object $X \in \baseCat$ one 
must have
\begin{align*}
\altexp{P}{Q}(X) &\simeq \Hom(\baseCat, \Cat)(\Yon X, \altexp{P}{Q}) \qquad 
\text{by the Yoneda lemma} \\
		&\simeq \Hom(\baseCat, \Cat)(\Yon X \times P, Q) \qquad \text{by 
		definition 
		of an exponential}
\end{align*}
So it remains to show that the pseudofunctor 
$\Hom(\baseCat, \Cat)\big(\Yon(-) \times P, Q \big) : \baseCat \to \Cat$ is 
indeed 
the exponential 
$\altexp{P}{Q}$ in $\Hom(\baseCat, \Cat)$, where $\Yon X := \baseCat(X, -)$ 
denotes the 
covariant Yoneda embedding.

To simplify the presentation we assume throughout this section that $\baseCat$ 
is a 2-category. The following lemma guarantees that this entails no loss of 
generality.

\begin{mylemma} \label{lem:ccc-and-biequivalent-implies-ccc}
Suppose that $\baseCat \simeq \baseCat'$ are biequivalent 
bicategories and $\catV$ is any bicategory. Then:
\begin{enumerate}
\item \label{c:biequivalence-to-hom-category} The hom-bicategories 
$\Hom(\baseCat, \catV)$ and 
$\Hom(\baseCat', \catV)$ are 
biequivalent, and
\item \label{c:biequivalence-to-ccc} If $\baseCat$ is cartesian closed, so is 
$\baseCat'$. 
\end{enumerate}
\begin{proof}
For~(\ref{c:biequivalence-to-hom-category}), suppose the biequivalence is given 
by pseudofunctors $P : \baseCat \leftrightarrows \baseCat' : Q$. Define 
pseudofunctors 
$Q_\ast : \Hom(\baseCat, \catV) \leftrightarrows \Hom(\baseCat' , \catV) : 
P_\ast$ by setting $Q_\ast(H) := H \circ Q$ and $P_\ast(F) := F \circ P$. From 
the biequivalence $\baseCat \simeq \baseCat'$ one obtains equivalences $PQ 
\simeq \id_{\baseCat'}$ and $QP \simeq \id_{\baseCat}$ and hence equivalences 
$P_\ast Q_\ast \simeq \id_{\Hom(\baseCat, \catV)}$ and $Q_\ast P_\ast \simeq 
\id_{\Hom(\baseCat', \catV)}$, as required.

For~(\ref{c:biequivalence-to-ccc}), one applies 
Lemma~\ref{lem:biequivalences-preserve-biuniversal-arrows} to carry the 
required biuniversal arrows from $\baseCat$ to $\baseCat'$ 
(\cf~also~Corollary~\ref{cor:biequivalences-preserve-bilimits}). 
\end{proof}
\end{mylemma}

We now turn to the construction of exponentials in $\Hom(\baseCat, \Cat)$. This 
entails constructing an adjoint equivalence 
$\Hom(\baseCat, \Cat)(R, \altexp{P}{Q}) 
	\simeq 
\Hom(\baseCat, \Cat)(R \times P, Q)$ 
for every triple of pseudofunctors 
$P, Q, R : \baseCat \to \Cat$. 
Since the 
definition of $\altexp{P}{Q}$ is also in terms of hom-categories, working with 
the 1- and 2-cells in $\Hom(\baseCat, 
\Cat)(R, \altexp{P}{Q})$ and $\Hom(\baseCat, \Cat)(R \times P, Q)$ quickly 
becomes complex, with several layers of data to consider. We therefore take 
the time to 
unwind some of the definitions we shall be using; as well as serving as a 
quick-reference on the details of the various definitions, this will fix 
notation for what follows. 

\subsection{A quick-reference summary}

\paragraph*{The pseudofunctor 
$\Hom(\baseCat, \Cat)\big( \Yon(-) \times P, Q\big)$.}

Suppose $f : X \to X'$ in $\baseCat$. The functor 
$	
\Hom(\baseCat, \Cat)(\Yon f \times P, Q) : 
	\Hom(\baseCat, \Cat)(\Yon X \times P , Q) 
		\to 
	\Hom(\baseCat, \Cat)(\Yon X' \times P, Q)$ 
takes a pseudonatural 
transformation $(\natTrans, \natCell) : \Yon X \times P \to Q$ to the 
pseudonatural transformation with components $\natTrans(- \circ f, =)$ and 
witnessing 2-cell given by the following composite for every $g : B \to B'$:
\begin{td}[column sep = 7em]
\baseCat(X', B) \times PB 
\arrow[swap]{d}{\baseCat(f, B) \times PB}
\arrow{r}{\baseCat(X', g) \times Pg}
\arrow[phantom]{dr}[description]{=} &
\baseCat(X', B') \times PB' 
\arrow{d}{\baseCat(f, B') \times PB'} \\

\baseCat(X, B) \times PB \arrow{r}{\baseCat(X, g) \times Pg}
\arrow[swap]{d}{\natTrans_B}
\arrow[phantom]{dr}[description, xshift=-4mm]{\twocell{\natCell_g}} &
\baseCat(X, B') \times PB'
\arrow{d}{\natTrans_{B'}} \\

QB \arrow[swap]{r}{Qg} &
QB'
\end{td}
The top square commutes because products in $\Cat$ are strict and we 
have assumed that $\baseCat$ is a 2-category. 

\begin{myremark}
We shall write both 
$\natTrans_B$ and 
$\natTrans(B, -, {=})$ to denote the component of a pseudonatural 
transformation $(\natTrans, \natCell)$ at an object $B$. These are just 
two notations 
for the same concept: the choice in any particular context is only dependent on 
which is clearest for exposition. Similar remarks apply to the 2-cells 
$\natCell$ and to 
modifications. 
\end{myremark}

\paragraph*{Pseudonatural transformations $R \To \altexp{P}{Q}$.}  To give a 
pseudonatural transformation $(\natTrans, \natCell) : R \To \Hom(\baseCat, 
\Cat)\big( \Yon(-) \times P , Q \big)$ is to give 
\begin{itemize}
\item For every $X \in \baseCat$ a functor 
$\natTrans_X : RX \to \Hom(\baseCat, \Cat)(\Yon X \times P , Q)$, 

\item For every $f : X \to X'$ in $\baseCat$ an invertible 2-cell (that is, a 
natural isomorphism) 
$\natCell_f$ as in the following diagram:
\begin{td}
RX \arrow{r}{Rf} \arrow[swap]{d}{\natTrans_X} 
\arrow[phantom]{dr}[description, xshift=4mm]{\twocell{\natCell_f}} &
RX' \arrow{d}{\natTrans_{X'}} \\

\Hom(\baseCat, \Cat)(\Yon X \times P , Q)
\arrow[swap]{r}[yshift=-2mm]{\Hom(\baseCat, \Cat)(\Yon f \times P , Q)} &
\Hom(\baseCat, \Cat)(\Yon X' \times P , Q)
\end{td} 
\end{itemize}
Thus, for every $r \in RX$ one obtains a pseudonatural 
transformation $\natTrans(X, r,-) : \Yon X \times P \To Q$ and an invertible 
2-cell (modification) 
$\natCell(f, r) : \natTrans(X', (Rf)(r), -) \to
	{\Hom(\baseCat, \Cat)(\Yon f \times P , Q)\big(\natTrans(X, r, -)\big)}$. 
The components of this modification are natural isomorphisms 
$\natCell(f, r, B)$, 
with components
\begin{equation}\label{eq:natTransCell-f-r-B} 
\lambda(h,x)^{\baseCat(X',B) \times PB} \bind 
\natTrans(X', (Rf)(r), B)(h,x) \xra{\natCell(f, r, B)(h,x)} 
\natTrans(X, r, B)(h \circ f, x)
\end{equation}
indexed by $B \in \baseCat$. (Note that we use the $\lambda$-notation 
$\lambda(h,x)^{\baseCat(X',B) \times PB} \bind \natTrans(X,r,B)(h,x)$ to 
anonymously refer to the action on objects 
$(h, x) \in \baseCat(X', B) \times PB$.) The modification axiom on $\natCell(f, r)$ requires 
\enlargethispage{\baselineskip}
that the diagram below commutes for every 
$(h, p) \in \baseCat(X, B) \times PB$, 
$g : B \to B'$ and 
$f : X \to X'$ in $\baseCat$:
\begin{equation} \label{eq:natTrans-f-r-modification-axiom}
\begin{tikzcd}[column sep = 6em, row sep = 2.5em]
\natTrans(X', (Rf)(r), B')\big(gh, (Pf)(p)\big) 
\arrow{r}[yshift=2mm]{\natCell(X', (Rf)(r), g)(h, (Pf)(p))}
\arrow[swap]{d}{\natCell(f,r)(gh, (Pf)(p))} &
(Qg)\big( \natTrans(X', (Rf)(r), B)(h,p) \big) 
\arrow{d}{(Qg)(\natCell(f,r)(h,p))} \\
\natTrans(X, r, B')\big(ghf, (Pf)(p)\big) 
\arrow[swap]{r}[yshift=-0mm]{\natCell(X,r,g)(hf, (Pf)(p))} &
(Qg)\big( \natTrans(X,r,B)(hf, p) \big)
\end{tikzcd}
\end{equation}

We can unfold the pseudonatural transformation $\natTrans(X, r,-)$ further. It 
has components given by functors
$\natTrans(X,r,B) : \baseCat(X, B) \times PB \to QB$  
(for $B \in \baseCat$), and 
for every $g : B \to B'$ one obtains an invertible 2-cell (that is, a natural 
isomorphism) 
$\cellOf{\natTrans}(X, r, g)$ as in
\begin{equation} \label{eq:cellOf-k-X-r-g}
\begin{tikzcd}[column sep = 5em, row sep = 2.5em]
\baseCat(X, B) \times PB \arrow{r}{\baseCat(X, g) \times Pg} 
\arrow[swap]{d}{\natTrans(X, r, B)}
\arrow[phantom]{dr}[description,xshift=-3mm]
	{\twocell{\cellOf{\natTrans}(X,r,g)}} &
\baseCat(X, B') \times PB' \arrow{d}{\natTrans(X,r,B')} \\
QB \arrow[swap]{r}{Qg} &
QB'
\end{tikzcd}
\end{equation}
Examining the components of this 2-cell, one sees that for each 
$(h, p) \in \baseCat(X, B) \times PB$ one obtains an invertible 
1-cell 
$\cellOf{\natTrans}(X,r,g)(h, p) : 
	\natTrans(X,r,B')\big( g \circ h, (Pg)(p) \big) 
		\to 
	(Qg)\big( \natTrans(X,r,B)(h,p) \big)$. 

There are then two levels of naturality at play, related 
via~(\ref{eq:natTrans-f-r-modification-axiom}). The naturality condition  
making 
$\cellOf{k}(X,r,-)$ a pseudonatural transformation requires that for every 
2-cell $\tau : g \To g' : {B \to B'}$ the following commutes:
\begin{td}[column sep = 10em, row sep = 2.5em]
\natTrans(X,r,B')\big( g \circ h, (Pg)(p) \big)  
\arrow[swap]{d}{\cellOf{\natTrans}(X,r,g)(h,p)}
\arrow{r}{\natTrans(X,r,B')( \tau \circ h, (P\tau )(p) )} &
\natTrans(X,r,B')\big( g' \circ h, (Pg)(p) \big) 
\arrow{d}{\cellOf{\natTrans}(X,r,g')(h,p)} \\

(Qg)\big( \natTrans(X,r,B)(h,p) \big)
\arrow[swap]{r}{(Q\tau)( \natTrans(X,r,B)(h,p) )} &
(Qg')\big( \natTrans(X,r,B)(h,p) \big)
\end{td}
On the other hand, the naturality condition making $\cellOf{\natTrans}(X,r,g)$ 
a natural 
transformation requires that for every $\rho : h \To h'$ in $\baseCat(X,B)$ and 
$t : p \to p'$ in $PB$, the following commutes:
\begin{td}[column sep = 10em, row sep = 2.5em]
\natTrans(X,r,B')\big( g \circ h, (Pg)(p) \big)  
\arrow[swap]{d}{\cellOf{\natTrans}(X,r,g)(h,p)}
\arrow{r}{\natTrans(X,r,B')( g \circ \rho, (Pg)(t) )} &
\natTrans(X,r,B')\big( g \circ h', (Pg)(p') \big)  
\arrow{d}{\cellOf{\natTrans}(X,r,g)(h',p')} \\

(Qg)\big( \natTrans(X,r,B)(h,p) \big)
\arrow[swap]{r}{(Qg)( \natTrans(X,r,B)(\rho,t) )} &
(Qg)\big( \natTrans(X,r,B)(h',p') \big)
\end{td}

\paragraph*{Modifications 
$(\altNat, \altCell) \to 
(\altaltNat, \altaltCell) 
: R \To \altexp{P}{Q}$.} 

To give a modification $\altModif : (\altNat, \altCell) \to 
(\altaltNat, \altaltCell) $ between pseudonatural transformations $R \To 
\altexp{P}{Q}$ is to give a natural transformation 
$\altModif_X : \altNat_X \To \altaltNat_X$ between functors of type 
$RX \to \Hom(\baseCat, \Cat)(\Yon X \times P, Q)$ for every $X \in \baseCat$, 
such that the whole 
$X$-indexed family of natural transformations satisfies the modification axiom. 

Unwinding the definition of natural transformation, $\altModif_X$ is a family 
of 2-cells (that is, modifications) 
$\altModif(X,r,-) : \altNat(X,r,-) \To \altaltNat(X,r,-)$, natural in $r \in 
\baseCat$ and such that every $\altModif(X,r,-)$ satisfies the modification 
axiom.
In 
particular, since every $\altModif(X,r,-)$ is a modification between 
pseudonatural 
transformations $\Yon X \times P \To Q$, for every 
$B \in \baseCat$ we have a natural transformation 
$\altModif(X,r,B) : 
\altNat(X,r,B) \To \altaltNat(X,r,B) : \baseCat(X, B) \times PB \to QB$.

\subsection{The cartesian closed structure of $\Hom(\baseCat, \Cat)$} 

To construct exponentials in $\Hom(\baseCat, \Cat)$ we are required to give:
\begin{itemize}
\item A biuniversal arrow $\eval_{P,Q}: \altexp{P}{Q} \times P \to Q$ for each 
$P, Q : \baseCat \to \Cat$, 
\item A mapping $\Lambda : 
ob\big(\Hom(\baseCat, \Cat)(R \times P, Q)\big) \to 
ob\big( \Hom(\baseCat, \Cat)(R, \altexp{P}{Q}) \big)$,
\item An invertible universal 2-cell 
$\eval_{P,Q} \circ \Lambda(\altNat, 
\altCell) \To 
(\altNat, \altCell)$ defining the counit, such that the unit is also invertible.
\end{itemize}
We take these components in turn. The 
main difficulty of the proof is maintaining a clear view of what one is 
required to 
construct, and ensuring that all the relevant axioms have been checked.

\paragraph*{The biuniversal arrow.} Our first step is the construction of the 
biuniversal arrow $\eval_{P,Q} : \altexp{P}{Q} \times P \to Q$. To be a 1-cell 
in $\Hom(\baseCat, \Cat)$, this needs to be a pseudonatural transformation for 
which each component is a functor
$e_X : \Hom(\baseCat, \Cat)(\Yon X \times P, Q) \times PX \to QX$.

Let $X \in \baseCat$ be fixed; we define $e_X$. Consider a pair 
$\big((\natTrans, \natCell), p\big) 
	\in \Hom(\baseCat, \Cat)(\Yon X \times P, Q)$ 
consisting of a pseudonatural transformation 
$(\natTrans, \natCell) : \Yon X \times P \To Q$ and an element $p \in PX$. 
Noting that, in particular, the 
component of $(\natTrans, \natCell)$ at $X \in \baseCat$ has type
$\baseCat(X, X) \times PX \to QX$, one obtains a functor 
$\natTrans(X,\Id_X, -) : PX \to QX$. We therefore define $e_X\big( (\natTrans, 
\natCell), p \big) := \natTrans(X,\Id_X, p)$. 

To extend this to 
morphisms, we need to define a morphism $\natTrans(X, \Id_X, p) \to 
\natTrans'(X, \Id_X, p')$ for every pair $(\modif, f)$ consisting of a
modification 
$\modif : (\natTrans, \natCell) \to (\natTrans', \natCell')$
 and morphism $f : p \to p'$.  The modification $\modif$  
is a family of natural transformations 
$\modif_X: \natTrans(X, -, =) \To \natTrans'(X, -, =)$ 
for $X \in \baseCat$, where naturality amounts to the following commutative 
diagram for every $\tau : h \To h' : X \to B$ and $f : p \to p'$ in $PB$:
\begin{td}[column sep = 4em]
\natTrans(X,h, p) \arrow{r}{\natTrans(X, \tau, f)} 
\arrow[swap]{d}{\modif_X(h,p)}  &
\natTrans(X, h', p') \arrow{d}{\modif_X(h',p')} \\
\natTrans'(X,h, p) \arrow[swap]{r}{\natTrans(X, \tau, f)}  &
\natTrans'(X, h', p') 
\end{td}
We 
define $e_X(\modif, f)$ to be the composite
\[
e_X(\modif, f) := \natTrans(X, \Id_X, p) \XRA{\modif_X(\Id_X, p)} \natTrans'(X, 
\Id_X, p) \XRA{\natTrans'(X, \Id_X, f)} \natTrans'(X, \Id_X, p')
\]
This definition is functorial.

Next we need to provide invertible 2-cells witnessing that the mappings $e_X$ 
are pseudonatural. That is, for every
$f : X \to X'$ in $\baseCat$ we need 
to provide a natural isomorphism as in the 
following diagram:
\begin{td}[column sep = 7em]
\Hom(\baseCat, \Cat)(\Yon X \times P, Q) \times PX
\arrow[phantom]{dr}[description, xshift=-8mm]{\twocell{\cellOf{e}_f}}
\arrow{r}[yshift=2mm]{\Hom(\baseCat, \Cat)(\Yon f \times P, Q) \times Pf} 
\arrow[swap]{d}{e_X} &
\Hom(\baseCat, \Cat)(\Yon X' \times P, Q) \times PX'
\arrow{d}{e_{X'}} \\
QX
\arrow[swap]{r}{Qf} &
QX'
\end{td}
Chasing an arbitrary element 
$\big( (\natTrans, \natCell), p\big) \in 
\Hom(\baseCat, \Cat)(\Yon X \times P, Q) \times PX$
through this diagram,
one sees that we need to provide an isomorphism 
$\natTrans\big(X',f, (Pf)(p)\big) \iso (Qf)(\natTrans(X, \Id_X, p))$
in $QX'$. We take 
\[
\cellOf{e}_f\big( (\natTrans, \natCell), p\big) 
:=  \natTrans(X', f, (Pf)(p)) = 
\natTrans(X', f \circ \Id_X, (Pf)(p)) 
\XRA{\cellOf{\natTrans}(X,r,f)(\Id_X, p)} 
(Qf)\big(\natTrans(X,r,B)(\Id_X, p) \big)
\]
using the natural isomorphism provided by 
diagram~(\ref{eq:cellOf-k-X-r-g}). 

\begin{mylemma} \label{lem:eval-defn}
The pair $(e, \cellOf{e})$ defined above is a pseudonatural 
transformation 
$\altexp{P}{Q} \times P \To Q$.
\begin{proof}
The naturality condition follows directly from that for $\natCell$. 
Similarly, the unit and associativity and unit laws hold immediately because 
they hold for $(\natTrans, \natCell)$. 
\end{proof}
\end{mylemma} 

We now have a candidate for the biuniversal arrow 
$\eval_{P,Q}$ defining 
exponentials. The next step is to define a mapping 
$\Lambda : 
ob\big(\Hom(\baseCat, \Cat)(R \times P, Q)\big) \to 
ob\big( \Hom(\baseCat, \Cat)(R, \altexp{P}{Q}) \big)$.

\paragraph*{The mapping $\Lambda$.} Let $(\altNat, \altCell)$ be a 
pseudonatural transformation $R \times P \To Q$. 
We define 
$\Lambda(\altNat, \altCell) : R \To \altexp{P}{Q}$ in stages. For the 1-cell
components 
we need to define a functor 
$RX \to \Hom(\baseCat, \Cat)(\Yon X \times P, Q)$ for every $X \in \baseCat$. 
We do this first. 

Fix some $X \in 
\baseCat$ and $r \in RX$. We define a pseudonatural transformation 
$(\Lambda\altNat)(X, r, -) : \Yon X \times P \To Q$. For every 
$B \in \baseCat$ we take the functor 
\begin{align*}
\baseCat(X, B) \times PB &\to QB \\
(h, p) &\mapsto \altNat\big(X, (Rh)(r), p\big)
\end{align*}
This is well-defined because $\altNat_X : RX \times PX \to QX$, so 
$(Rh)(r) \in RB$. We take the evident functorial action on 2-cells: 
$(\Lambda \altNat)(X, r, B)(\tau, f) 
:= \altNat\big(X, (R\tau)(r), f\big)$. 

To extend these 1-cells to a pseudonatural transformation we need to provide a 
natural isomorphism 
$\cellOf{(\Lambda \altNat)}(X, r, g)$ as in  
\begin{td}[column sep = 7em, row sep = 2.5em]
\baseCat(X, B) \times PB \arrow{r}{\baseCat(X, g) \times Pg} 
\arrow[swap]{d}{(\Lambda\altNat)(X,r,B)}  
\arrow[phantom]{dr}[description, xshift=-4mm]
	{\twocell{\cellOf{(\Lambda \altNat)}(X, r, g)}}  &
\baseCat(X, B') \times PB' \arrow{d}{(\Lambda\altNat)(X,r)_{B'}} \\
QB \arrow[swap]{r}{Qg} &
QB'
\end{td}
for every $g : B \to B'$ in $\baseCat$. 
So for every $(h, p) \in \baseCat(X, B) \times PB$ 
we need to give an isomorphism 
$\altNat\big( X, (Rgh)(r), (Pg)(p) \big) \iso 
(Qg)\big(\altNat\big( X, (Rh)(r), p \big)\big)$, for which we take the 
composite defined by commutativity of
\begin{td}
\altNat\big( X, (Rgh)(r), (Pg)(p) \big)  
\arrow{rr}{\cellOf{(\Lambda \altNat)}(X, r, g)}
\arrow[swap]{dr}{\altNat\left( X, (\phi_{g,h}^R)^{-1}(r), (Pg)(p) \right)} &
\: &
(Qg)\big(\altNat\big( X, (Rh)(r), p \big)\big) \\

\: &
\altNat\big( X, (Rg)(Rh)(r), (Pg)(p) \big)
\arrow[swap]{ur}{\altCell(g, (Rh)(r), p)} &
\: 
\end{td}
This definition is natural 
in $g$ because $\phi^R_{g,h}$ and $\altCell_g$ both are. The unit and 
associativity laws follow easily from those of $(\altNat, \altCell)$, 
yielding the following.

\begin{prooflesslemma} \label{lem:def-of-Lambda-part-one}
For every $X \in \baseCat$, $r \in RX$ and pseudonatural transformation 
$(\altNat, \altCell) : R \times P \To Q$, the pair 
$\big( (\Lambda\altNat)(X,r,-), \cellOf{(\Lambda \altNat)}(X,r,-)\big)$ is a 
pseudonatural transformation $\Yon X \times P \To Q$.
\end{prooflesslemma}

The preceding lemma defines a mapping 
$ob(RX) \to ob\big(\Hom(\baseCat, \Cat)(\Yon X \times P, Q)\big)$. Our next 
task is to extend this to a functor. So suppose $f : r \to r'$ in $RX$. To give 
a modification $(\Lambda \altNat)(X, f, -) : (\Lambda \altNat)(X, r, -) \to 
(\Lambda \altNat)(X, r', -)$, one must provide a family of natural 
transformations $(\Lambda \altNat)(X, r, B) \To (\Lambda \altNat)(X, r', B)$ 
indexed by $B \in \baseCat$. For a fixed choice of $B$ and 
$(h, p) \in 
\baseCat(X,B) \times PB$, we take the 1-cell
\[
(\Lambda \altNat)(X, f, B)(h,p) := 
\altNat( X, (Rh)(r), p) 
\XRA{\altNat( X, (Rh)(f), p)} 
\altNat( X, (Rh)(r'), p) 
\]
This is natural in $h$ and $p$ by functoriality. The modification law for
$(\Lambda \altNat)(X, f, -)$ is a consequence of the naturality 
properties. For $(h,p)$ as above and $f : r \to r'$, one has
\begin{td}[column sep = 7.5em, row sep = 2.5em]
\altNat\big(X', (Rgh)(r), (Pg)(p)\big) 
\arrow[swap]{d}{\altNat\left(X', (\phi_{g,h}^R)^{-1}(r), (Pg)(p)\right)}
\arrow{r}{\altNat(X', (Rgh)(f), (Pg)(p))} &
\altNat\big(X', (Rgh)(r'), (Pg)(p)\big) 
\arrow{d}{\altNat\left(X', (\phi_{g,h}^R)^{-1}(r'), (Pg)(p)\right)} \\

\altNat\big( X', (Rg)(Rh)(r), (Pg)(p) \big)
\arrow{r}[yshift=2mm]{\altNat( X', (Rg)(Rh)(f), (Pg)(p) )}
\arrow[swap]{d}{\altCell(g, (Rh)(r), p)} &
\altNat\big( X', (Rg)(Rh)(r'), (Pg)(p) \big) 
\arrow[]{d}{\altCell(g, (Rh)(r'), p)} \\

(Qg)\big( \altNat(X, (Rh)(r), p) \big) 
\arrow[swap]{r}{(Qg)(\altNat(X, (Rh)(f), p))} &
(Qg)\big( \altNat(X, (Rh)(r'), p) \big)
\end{td}
in which the top square commutes by naturality of $\phi^R$ and the bottom 
square by the fact that $\altCell_g$ is a natural transformation.

We have now defined a functor $(\Lambda \altNat)(X, -, =) : RX \to 
\Hom(\baseCat, \Cat)\big(  \Yon X \times P, Q \big)$ for each $X \in \baseCat$. 
It remains to show these functors are the components of a pseudonatural 
transformation. Thus, for every $f : X \to X'$ we need to provide invertible 
2-cells $\cellOf{(\Lambda\altNat)}(f, -, =)$ as in
\begin{td}[column sep=3em, row sep = 2.4em]
RX \arrow{r}{Rf}  
\arrow[swap]{d}{(\Lambda \altNat)(X, -, =)}
\arrow[phantom]{dr}[description, xshift=4mm]
	{\twocell{\cellOf{(\Lambda\altNat)}(f, -, =)}} &
RX' 
\arrow[]{d}{(\Lambda \altNat)(X', -, =)} \\

\Hom(\baseCat, \Cat)(\Yon X \times P, Q) 
\arrow[swap]{r}[yshift=-2mm]{\Hom(\baseCat, \Cat)(\Yon f \times P, Q)} &
\Hom(\baseCat, \Cat)(\Yon X' \times P, Q)
\end{td}
This diagram requires an isomorphism 
\begin{equation} \label{eq:lambda-j-nat-witness}
\lambda B^\baseCat \bind
\lambda(h,p)^{\baseCat(X', B) \times PB} \bind 
\altNat(X, (Rh)(Rf)(r), p) \iso \altNat(X, (Rhf)(r), p) 
\end{equation}
for each $r \in RX$, for which we 
take simply 
$\lambda B^\baseCat \bind 
\lambda(h,p)^{\baseCat(X' B) \times PB} \bind 			
\altNat(X, \phi^R_{h,f}(r), p)$. 
The unit and 
associativity laws then follow from the unit and associativity laws of the 
pseudofunctor $R$. 

We record our progress in the following lemma.

\begin{prooflesslemma} \label{lem:lem:def-of-Lambda-part-two}
The pair $\big( (\Lambda \altNat)(X, -, =), \cellOf{(\Lambda \altNat)}(f, -, =) 
\big)$ is a pseudonatural transformation $R \To \Hom(\baseCat, \Cat)(\Yon X 
\times P, Q)$. 
\end{prooflesslemma}

We therefore define the required mapping as follows:
\begin{align*}
\Lambda : 
ob\big(\Hom(\baseCat, \Cat)(R \times P, Q)\big) &\to 
ob\big( \Hom(\baseCat, \Cat)(R, \altexp{P}{Q}) \big) \\
(\altNat, \altCell) &\mapsto 
\big( (\Lambda \altNat)(X, -, =), 
\cellOf{(\Lambda \altNat)}(f, -, =) 
\big)
\end{align*}
Our next task is to define the universal arrow, which will act as the 
counit. 

\paragraph*{The counit $\evalMod$.} We begin by calculating $\eval_{P,Q} \circ 
\big((\natTrans, 
\natCell) \times P \big) : R \times P \To Q$ for any 
$(\natTrans, \natCell) : 
R \To \altexp{P}{Q}$. 
The component at $X \in \baseCat$ is the 
functor 
acting on $(r, p) \in RX \times PX$ by
\begin{align*}
\big(e_X \circ (\natTrans_X \times PX)\big)(X, r, p) &= e_X\big( 
\natTrans(X,r,-), p \big) \\
	&= e_X\big( \lambda B^\baseCat \bind 
				\lambda(h,x)^{\baseCat(X, B) \times PB} \bind 
				\natTrans(X,r,B)(h,x), p  
	\big)  \\
	&= \natTrans(X,r,X)(\Id_X, p) 
\end{align*}
For any $f : X \to X'$ and $(r, p) \in RX \times PX$, the witnessing 2-cell is 
defined by the following commutative diagram: 
\begin{equation} \label{eq:eval-circ-natTrans-times-P-twocell}
\begin{tikzcd}[column sep = 9em, row sep = 2.5em]
\natTrans(X',(Rf)(r),X')\big(\Id_{X'}, (Pf)(p)\big)
\arrow{r}
	{\overline{(\eval_{P,Q} \circ ((\natTrans, \natCell) \times P))}_f(r,p)}
\arrow[swap]{d}{\natCell(f,r)(\Id_{X'}, (Pf)(p))} &
(Qf)\big(  \natTrans(X,r,X)(\Id_X, p) \big) \\
\natTrans(X, r, X')\big(\Id_{X'} \circ f, (Pf)(p)\big) 
\arrow[equals]{r} &
\natTrans(X, r, X')\big(f \circ \Id_X, (Pf)(p) \big) 
\arrow[swap]{u}{ \cellOf{\natTrans}(X,r,f)(\Id_X, p) }
\end{tikzcd}
\end{equation}
Note that both levels 
of naturality appear in this definition: the first arrow arises
from the components of the modification 
$\natCell(f,r)$ given in~(\ref{eq:natTransCell-f-r-B}), while 
the second arises from the 2-cell witnessing
the naturality of $\natTrans_X$ in diagram~(\ref{eq:cellOf-k-X-r-g}).

Now suppose that $(\altNat, \altCell) : R \times P \To Q$ and consider 
$\eval_{P,Q} \circ 
	\big(\Lambda(\altNat, \altCell) \times P \big): R \times P \To Q$. 
The 1-cell components of this pseudonatural transformation act by
\begin{equation} \label{eq:eval-circ-lambda-altNat-times-P}
\begin{aligned}
RX \times PX &\to QX \\
(r, p) &\mapsto \altNat\big(X, (R\Id_X)(r), p\big)
\end{aligned}
\end{equation}
and for $f : X \to X'$ and $(r, p) \in RX \times PX$ the witnessing 2-cell is 
the composite
\begin{td}[column sep = 7em]
\altNat\big(X', (R\Id_{X'})(Rf)(r), (Pf)(p)\big)
\arrow{r}
	{\overline{(\eval_{P,Q} \circ (\Lambda(\altNat, \altCell) \times P))}_f}
\arrow[swap]{d}{\altNat\left(X', \phi_{\Id, f}^R(r), (Pf)(p)\right)} &
(Qf)\big( \altNat(X, R(\Id_X)(r), p) \big) \\

\altNat\big(X', R(\Id_{X'} \circ f)(r), (Pf)(p)\big)
\arrow[equals]{d} &
\: \\

\altNat\big(X', R(f\circ\Id_{X})(r), (Pf)(p)\big)
\arrow[swap]{r}[yshift = -2mm]
	{\altNat\left(X', (\phi_{f, \Id}^R)^{-1}(r), (Pf)(p)\right)} &
\altNat\big(X', R(f)R(\Id_{X})(r), (Pf)(p)\big)
\arrow[swap]{uu}{\altCell(f, (R\Id_X)(r), p)}
\end{td}

By the identification~(\ref{eq:eval-circ-lambda-altNat-times-P}), to define the 
counit modification 
$\evalMod : \eval_{P,Q} \circ 
\big(\Lambda(\altNat, \altCell) \times P \big) \to 
(\altNat, \altCell)$ we need to provide a 
natural transformation 
$\evalMod_X : \altNat\big(X, (R\Id_X)(-), {=}\big) 
\To 
\altNat(X,-,{=}) :
RX \times PX \to QX$
for every $X \in \baseCat$. We take the obvious choice, namely 
$\lambda(r, p)^{RX \times PX} \bind \altNat\big(X, (\psi_X^R)^{-1}(r), p\big)$. 
Since $\psi^R_X : \Id_{RX} \To R\Id_X$ is a 2-cell in $\Cat$,~\ie~a natural 
transformation, it only remains to check the 
modification axiom. 

\begin{mylemma} \label{lem:counit-defn}
The family of 2-cells 
$\evalMod_X := \altNat\big(X, (\psi_X^R)^{-1}(-), {=}\big)$ 
(for $X \in \baseCat$) 
form a modification 
$\eval_{P,Q} \circ \Lambda(\altNat, \altCell) 
\to (\altNat, \altCell)$.
\begin{proof}
We need to verify that the following diagram commutes for every $f : X \to X'$ 
in $\baseCat$:
\begin{equation} \label{eq:evalMod-a-modification}
\begin{tikzcd}[column sep = 1.5em, row sep = 1.9em]
\altNat\big(X', (Rf)(r), (Pf)(p) \big) 
\arrow[swap]{d}[]{\evalMod_{X'}((Rf)(r),(Pf)(p))}
\arrow{r}{\altCell(f, r,p)} &
(Qf)\big(\altNat(X, r,p)\big)  
\arrow{ddddd}[]{(Qf)(\evalMod_X(r,p))) } \\
\altNat\big(X', (R\Id_{X'})(Rf)(r), (Pf)(p)\big) 
\arrow[swap]{d}{\altNat(X', \phi_{\Id, f}^R(r), (Pf)(p)) }
\arrow[swap, bend right = 90]{dddd}
	[description, yshift=8.5mm, xshift=-1mm]
	{\overline{(\eval_{P,Q} \circ (\Lambda(\altNat, \altCell) \times P))}
		(f, r,p)} &
\: \\
\altNat\big(X', R(\Id_{X'} \circ f)(r), (Pf)(p)\big) 
\arrow[equals]{d}{} &
\: \\
\altNat\big(X', R(f \circ \Id_X)(r), (Pf)(p)\big)  
\arrow[swap]{d}{\altNat\left(X', (\phi_{f, \Id}^R)^{-1}(r), (Pf)(p)\right) } &
\: \\
\altNat\big(X', R(f)R(\Id_X)(r), (Pf)(p)\big) 
\arrow[swap]{d}{\altCell(f, R(\Id_X)(r), p))} &
\: \\
(Qf)\big( \altNat(X, R(\Id_X)(r), p) \big) 
\arrow[equals]{r}{} &
(Qf)\big( \altNat(X, R(\Id_X)(r), p) \big) &
\end{tikzcd}
\end{equation}
To this end, one uses the two unit laws of a 
pseudofunctor to see that the following commutes:
\begin{td}[column sep = 2em, row sep = 2em]
\: &
\altNat_{X'} \circ (Rf \times Pf)
\arrow[swap]{dl}[]{\altNat_{X'} \circ (\psi^R_{X'} \times Pf)}
\arrow[]{ddddr}
	{\altNat_{X'} \circ ((Rf \circ \psi^R_X) \times Pf)}
\arrow[equals]{ddd} &
\: \\
\altNat_{X'} \circ \big( (R\Id_{X'} \circ Rf) \times Pf \big)
\arrow[swap]{d}{\altNat_{X'} \circ (\phi^R_{\Id, f} \times Pf) } &
\: &
\: \\
\altNat_{X'} \circ \big( R(\Id_{X'} \circ f) \times Pf \big)
\arrow[equals]{dr}{} 
\arrow[equals, bend right]{ddr} &
\: &
\: \\
\: &
\altNat_{X'} \circ \big( Rf \times Pf \big)
\arrow[equals]{d}{} &
\: \\
\: &
\altNat_{X'} \circ \big( R(f \circ \Id_X) \times Pf \big)
\arrow[swap]{r}[yshift=-2mm]
	{\altNat\left(X', (\phi_{f, \Id}^R)^{-1}, Pf\right) } &
\altNat_{X'} \circ \big( (Rf \circ R\Id_X) \times Pf \big) 
\end{td}
Diagram~(\ref{eq:evalMod-a-modification}) therefore reduces to 
\begin{td}[row sep = 2.5em]
\altNat\big(X', (Rf)(r), (Pf)(p) \big) 
\arrow[swap]{d}[]{\altNat(X', (Rf)(\psi^R_X)(r), (Pf)(p) ) }
\arrow{r}{\altCell(f,r,p)} &
(Qf)\big(\altNat(X, r,p)\big)  
\arrow{dd}[]{(Qf)(\altNat(X, \psi_X^R(r),p)) } \\

\altNat\big(X', R(f)R(\Id_X)(r), (Pf)(p)\big) 
\arrow[swap]{d}{\altCell(f, R(\Id_X)(r), p))} &
\: \\

(Qf)\big( \altNat(X, R(\Id_X)(r), p) \big) 
\arrow[equals]{r}{} &
(Qf)\big( \altNat(X, R(\Id_X)(r), p) \big) &
\end{td}
which commutes by the naturality of $\altCell(f, -, {=})$ in $r$. 
\end{proof}
\end{mylemma}

We have constructed our candidate counit $\evalMod$; now we
need to show it is universal. 
For the existence part of this claim, we need to construct a modification 
$\trans{\modif} : (\natTrans, \natCell)  \to 
\Lambda(\altNat, \altCell)$
for every pair of pseudonatural transformations
$(\altNat, \altCell) : R \times P \To Q$ and 
$(\natTrans, \natCell) : R \To \altexp{P}{Q}$
and every modification 
$\modif : 
	\eval_{P,Q} \circ \big((\natTrans, \natCell) \times P\big) 
		\to 
	(\altNat, \altCell)$.

\paragraph*{The modification $\trans{\modif}$.}
We begin by unwinding the definition of a modification 
${\eval_{P,Q} \circ 
\big((\natTrans, \natCell) \times P\big)} 
\to (\altNat, \altCell)$. 
For every $X \in \baseCat$ and $(r, p) \in RX \times PX$, we are given a 1-cell 
$\modif(X,r,p) : 
	\natTrans(X,r,X)(\Id_X, p) 
		\to 
	\altNat(X,r,p)$ 
in $QX$. These are 
natural in the sense that, for any 
$g : r \to r'$ and $h : p \to p'$ in $RX \times PX$, the following commutes:
\begin{equation*} \label{eq:eq:transModif-def-naturality}
\begin{tikzcd}[column sep = 6em]
\natTrans(X,r,X)(\Id_X, p) 
\arrow{r}{\natTrans(X,g,X)(\Id_X, h)} 
\arrow[swap]{d}{\modif(X,r,p)} & 
\natTrans(X,r',X)(\Id_X, p') 
\arrow{d}{\modif(X,r',p')} \\
\altNat(X,r,p) \arrow[swap]{r}{\altNat(X,g,h)} &
\altNat(X,r',p')
\end{tikzcd}
\end{equation*}
The $X$-indexed family of natural transformations
$\modif(X, -, {=})$
is subject to the modification axiom, which requires that the following 
commutes for every $f : X \to X'$ in $\baseCat$ (recall the definition of 
$\overline{(\eval_{P,Q} \circ ( (\natTrans, \natCell) \times P )}_f$  
from~(\ref{eq:eval-circ-natTrans-times-P-twocell})):
\begin{equation} \label{eq:transModif-def-modification-axiom}
\begin{tikzcd}
\natTrans\big( X', (Rf)(r), X'\big)\big(\Id_{X'}, (Pf)(p)\big) 
\arrow[swap]{d}{\natCell(f,r, B)(\Id_{X'}, (Pf)(p))}
\arrow{r}[yshift=2mm]{\modif(X', (Rf)(r), (Pf)(r))} &
\altNat\big(X', (Rf)(r), (Pf)(p) \big) 
\arrow{ddd}{\altCell(f, r, p)} \\
\natTrans(X, r, X')\big( \Id_{X'} \circ f, (Pf)(p) \big)
\arrow[equals]{d} &
\: \\
\natTrans(X, r, X')\big( f \circ \Id_X, (Pf)(p) \big)
\arrow[swap]{d}{\natCell(X,r,f)(\Id_X, p)} &
\: \\
(Qf)\big( \natTrans(X,r,X)(\Id_X,p) \big) 
\arrow[swap]{r}{(Qf)(\modif(X,r,p))} &
(Qf)\big( \altNat(X,r,p) \big)
\end{tikzcd}
\end{equation}

Now, to define $\trans\modif$ we are required to provide a 2-cell
$\trans{\modif}_X : \natTrans_X \to (\Lambda \altNat)_X$ for every
$X \in \baseCat$, subject to the modification axiom. Since
$\natTrans_X$ and $(\Lambda \altNat)_X$ are functors 
$RX \to \altexp{P}{Q}X$, such a natural transformation consists of a family of 
1-cells (modifications) 
$\trans{\modif}(X, r, -) : 
	\natTrans(X, r, -) \to (\Lambda \altNat)(X, r, -)$
that is natural in $r$. We build this data in stages.

Fix $X \in \baseCat$ and $r \in RX$. We begin  
by defining the modifications $\trans{\modif}(X, r, -)$. For the components, we 
define a natural transformation 
$\trans{\modif}(X, r, B) : \natTrans(X,r,B) 
\To (\Lambda \altNat)(X,r,B)$ for each $B \in \baseCat$ as follows. For 
$(h, p) \in \baseCat(X, B) \times PB$, we take the 1-cell defined by
commutativity of the diagram below, 
where the bottom arrow arises from the fact that each $\natCell_f$ is a 
modification with type given 
in~(\ref{eq:natTransCell-f-r-B}):
\begin{equation} \label{eq:def-of-transModif}
\begin{tikzcd}[column sep = 7em, row sep=2.5em]
\natTrans(X,r,B)(h, p) 
\arrow{r}{\trans{\modif}(X,r,B)(h,p)}
\arrow[equals]{d} &
\altNat\big(B, (Rh)(r), p \big) \\
\natTrans(X,r, B)(\Id_B \circ h, p) 
\arrow[swap]{r}{\natCell(h, r, B)(\Id_B, p)^{-1}} &
\natTrans\big(B, (Rh)(r), B\big)(\Id_B, p) 
\arrow[swap]{u}{\modif(B, (Rh)(r), p)} &
\end{tikzcd}
\end{equation}
The family of 1-cells thus defined is 
natural in $(h,p)$ because each component is. We claim that the family of 
natural transformations 
$\trans{\modif}(X, r, -)$ is a modification. This 
entails checking that the 
following commutes for every $f : B \to B'$ in $\baseCat$: 
\begin{td}[column sep = 4em, row sep=2.5em]
\natTrans(X,r,B) \circ \big(\baseCat(X, f) \times Pf\big)
\arrow{r}[yshift=2mm]{\trans{\modif}(X,r,B) \circ  (\baseCat(X, f) \times Pf)}
\arrow[swap]{d}{\natCell(X,r,f)} &
(\Lambda \altNat)(X,r,B) \circ \big(\baseCat(X, f) \times Pf\big) 
\arrow{d}{\cellOf{(\Lambda \altNat)}(X,r,f)} \\

(Qf)\big( \natTrans(X,r,B) \big) 
\arrow[swap]{r}{(Qf)(\trans{\modif}(X,r,B))} &
(\Lambda \altNat)(X,r,B)
\end{td}
To prove this, fix some 
$(h, p) \in \baseCat(X, B) \times PB$. 
Applying the naturality of $\modif$ with 
respect to the map 
$\phi^R_{f,h}(r) : (Rf)(Rh)(r) \to R(f \circ h)(r)$, and the modification 
axiom~(\ref{eq:transModif-def-modification-axiom}), one reduces the claim to 
showing that
\begin{equation*}
\makebox[\textwidth]{
\begin{tikzcd}[column sep = -1.7em, row sep=2em, ampersand replacement = \&]
\: \&
\natTrans(X, r, B')(\Id_{B'} \circ f \circ h, (Pf)(p))
\arrow[equals]{dr} \&
\: \\
\natTrans(B', R(fh)(r), B')(\Id_{B'}, (Pf)(p)) 
\arrow{ur}{\natCell(f \circ h, r)(\Id_{B'}, (Pf)(p))} \&
\: \&
\natTrans(X,r,B')(f \circ h, (Pf)(p)) 
\arrow{dd}{\natCell(X,r,f)(h,p)} \\
\natTrans(B', (Rf)(Rh)(r), B')(\Id_{B'}, (Pf)(p))
\arrow{u}{\natTrans\left(B', \phi_{f,h}^R(r), B')(\Id_{B'}, (Pf)(p)\right)}
\arrow[swap]{d}{\natCell(B, R(h)(r), f)(\Id_{B'}, (Pf)(p))} \&
\: \&
\: \\
\natTrans(B, (Rh)(r), B')(\Id_{B'} \circ f, (Pf)(p))
\arrow[equals]{d}{} \&
\: \&
(Qf)\big( \natTrans(X,r,B)(h,p) \big) 
\arrow[equals]{d} \\
\natTrans(B, (Rh)(r), B')(f \circ \Id_B, (Pf)(p)) 
\arrow[swap]{dr}{\natCell(B, (Rh)(r), f)(\Id_B, p)} \&
\: \&
(Qf)\big( \natTrans(X,r,B)(\Id_B \circ h,p) \big) \\
\: \&
(Qf)\big( \natTrans(B', (Rh)(r), B')(\Id_B, p) \big) 
\arrow[swap]{ur}{(Qf)(\natCell(h,r)(\Id_B,p))} \&
\: 
\end{tikzcd}
}
\end{equation*}
This commutes by an application of the associativity law for $R$ and the 
modification axiom~(\ref{eq:natTrans-f-r-modification-axiom}) for 
$\natCell(f,r)$. 

Thus, $\trans{\modif}(X,r)$ is a modification 
$\big( \natTrans(X,r, -), \natCell(X,r,-) \big) \to 
\big( (\Lambda\altNat)(X,r,-), \cellOf{(\Lambda\altNat)}(X,r,-) \big)$ 
for every 
$X \in \baseCat$ and $r \in RX$. Moreover, since each of the components in the 
definition of 
$\trans{\modif}(X,r)$ 
is natural in $r$, this $r$-indexed family of 1-cells forms a natural 
transformation 
$\trans{\modif}_X : \natTrans_X \To (\Lambda \altCell)_X$. 

To show that 
$\trans\modif$ is a modification 
$( \natTrans, \natCell) \to 
( \Lambda\altNat, \cellOf{\Lambda\altNat} )$, it 
remains to check the following modification law for every $f : X \to X'$ and 
$(h, p) \in \baseCat(X', B) \times PB$:
\begin{equation} \label{eq:transModif-modification-axiom}
\begin{tikzcd}
\natTrans\big(X', (Rf)(r), B\big)(h,p) 
\arrow{r}{\natCell(f,r)} 
\arrow[swap]{d}{\trans{\modif}(X, (Rf)(r), B)(h,p)} &
\natTrans\big(X, r, B\big)(h \circ f,p)
\arrow{d}{\trans{\modif}(X, r, B)(hf,p)} \\
(\Lambda \altNat)\big(X', (Rf)(r), B\big)(h,p) 
\arrow[swap]{r}{(\Lambda \altCell)(f)} &
(\Lambda \altNat)\big(X, r, B\big)(h \circ f,p)
\end{tikzcd}
\end{equation}
This follows from the  associativity law for 
${\eval_{P,Q} \circ 
\big((\natTrans, \natCell) \times P\big)}$, 
namely
\begin{td}
\natTrans\big(B, (Rh)(Rf)(r), B\big)(\Id_B, p) 
\arrow[swap]{d}{\natCell(h, (Rf)(r))(\Id_B, p)}
\arrow{rr}[yshift=0mm]{\natTrans\left(B, \phi_{h,f}^R(r), B\right)(\Id_B, p)} &
\: &
\natTrans\big(B, R(hf)(r), B\big)(\Id_B, p) 
\arrow{dd}{\natCell(h \circ f, r)(\Id_B, p)}  \\

\natTrans\big( X', (Rf)(r), B \big)(\Id_B \circ h, p) 
\arrow[equals]{d} &
\: &
\: \\

\natTrans\big( X', (Rf)(r), B \big)(h, p)
\arrow[swap]{r}{\natCell(f, r)(h, p)} &
\natTrans\big(X,r, B\big)(h \circ f, p)
\arrow[equals]{r}  &
\natTrans\big(X,r, B\big)(\Id_B \circ h \circ f, p) 
\end{td}
together with the naturality of $\modif_X$ 
with respect to the morphism 
$\phi^R_{h,f}(r) : (Rh)(Rf)(r) \to R(hf)(r)$. We summarise the result:

\begin{prooflesslemma} \label{lem:transModif-a-modification}
The family of natural transformations $\trans{\modif}(X,-,{=})$ defined 
in~(\ref{eq:def-of-transModif}) forms a modification $( \natTrans, \natCell) 
\to ( \Lambda\altNat, \cellOf{\Lambda\altNat} )$.
\end{prooflesslemma}
The final part of the proof is showing that $\trans{\modif}$ is the unique 
modification $\altModif$ such that
\begin{equation} \label{eq:transModif-ump}
\begin{tikzcd}
\eval_{P,Q} \circ \big( (\natTrans, \natCell) \times P \big) 
\arrow[swap]{dr}{\modif}
\arrow{rr}{\eval_{P,Q} \circ (\altModif \times P)} &
\: &
\eval_{P,Q} \circ \big( \Lambda(\altNat, \altCell) \times P\big)
\arrow{dl}{\evalMod} \\
\: &
(\altNat, \altCell) &
\:
\end{tikzcd}
\end{equation}
We turn to this next.

\paragraph*{The universal property of $\evalMod$.}

The existence part of the claim follows from the unit law of a 
pseudonatural transformation and the fact that $\modif(X,r,p)$ is a natural 
transformation: 
\begin{td}
\natTrans(X,r,X)(\Id_X, p) 
\arrow[phantom]{ddr}[yshift=-5mm]{\overset{\text{def}}{=}}
\arrow[swap]{dd}{\trans{\modif}(X,r,X)(\Id_X, p)}
\arrow[equals]{dr} &
\: &
\: \arrow[bend left, dd, phantom, ""{coordinate, name=Z}] \\

\: &
\natTrans(X, r, X)(\Id_X \circ \Id_X, p) 
\arrow{d}{\natCell(\Id_X, r)(\Id_X, p)^{-1}}
\arrow[
equals,
swap,
rounded corners,
to path=
{ -- ([xshift=2ex]\tikztostart.east)
-| (Z) [near end]\tikztonodes
|- ([xshift=2ex]\tikztotarget.east)
-- (\tikztotarget)}
]
{dd}[yshift=0em]{} &
\: \\

\altNat\big(X, R(\Id_X)(r), p\big)
\arrow[swap]{dd}{\altNat\left(X, (\psi_X^R)^{-1}(r), p\right)}
\arrow[phantom]{dr}[yshift=-5mm]{\overset{\text{nat}}{=}} &
\natTrans\big(X, R(\Id_X)(r), X\big)(\Id_X, p)
\arrow{l}[yshift=-2mm]{\modif(X, R(\Id_X)(r), p)}
\arrow{d}{\natTrans\left(X, (\psi_X^R)^{-1}(r), X\right)(\Id_X, p)}
\arrow[phantom]{r}
	[description, font=\scriptsize, xshift=1mm]{\overset{\text{unit law}}{=}} 
&
\: \\

\: &
\natTrans(X, r, X)(\Id_X, p) 
\arrow{dl}{\modif(X,r,p)} &
\: \\

\altNat(X,r,p) &
\: &
\:
\end{td}
For uniqueness, suppose that $\altModif$ is a modification 
filling~(\ref{eq:transModif-ump}). Then, applying the definition 
of $(\Lambda \altNat)(f, -, {=})$ from~(\ref{eq:lambda-j-nat-witness}), one 
obtains the diagram below, in which one uses the modification 
axiom \big(\cf~(\ref{eq:transModif-modification-axiom})\big), the assumption on 
$\altModif$ and the unit law of a pseudofunctor:
\begin{td}[row sep = 4em]
\natTrans(X,r,B)(h,p) 
\arrow[equals]{d}{} 
\arrow[phantom]{dr}[description]{=}
\arrow[bend left = 20]{dr}{\altModif(X,r,B)(h, p)} &
\: &
\: \arrow[bend left = 38, dd, phantom, ""{coordinate, name=Z}]  \\

\natTrans(X,r,B)(\Id_B \circ h, p)
\arrow[phantom]{dr}[description, yshift = 1mm]
	{\overset{\text{modif. law}}{=}}
\arrow{r}[yshift=2mm]{\altModif(X,r,B)(\Id_B \circ h, p)} 
\arrow[swap]{d}{\natCell(h,r)(\Id_B, p)^{-1}}  &
\altNat\big(B, R(\Id_B \circ h)(r), p\big) 
\arrow{d}{\altNat\left(B, (\phi_{\Id, h}^R)^{-1}(r), p\right) } 
\arrow[
equals,
swap,
rounded corners,
to path=
{ -- ([xshift=2ex]\tikztostart.east)
-| (Z) [near end]\tikztonodes
|- ([xshift=2ex]\tikztotarget.east)
-- (\tikztotarget)}
]
{dd}[yshift=0em]{} &
\: \\

\natTrans\big( B, (Rh)(r), B \big)(\Id_B, p) 
\arrow[phantom]{dr}[description]
{\overset{\text{(\ref{eq:transModif-ump})}}{=}}
\arrow{r}[yshift=2mm]{\altModif(B, (Rh)(r), B)(\Id_B, p)}
\arrow[bend right = 30, swap]{dr}{\modif(B, (Rh)(r), B)(\Id_B, p)} &
\altNat\big( B, (R\Id_B)(Rh)(r), p \big) 
\arrow{d}{\altNat\left(B, (\psi_B^R)^{-1}(Rh)(r), p\right) } 
\arrow[phantom]{r}[description, xshift=1mm, yshift=3mm]
	{\overset{\text{unit law}}{=}} &
\: \\

\: &
\altNat\big( B, (Rh)(r), p \big) &
\:
\end{td}
Since the left-hand leg of this diagram is the definition of 
$\trans\modif$~(\ref{eq:def-of-transModif}), one  
obtains the required universal property:

\begin{prooflesslemma}
For any modification 
$\modif : \eval_{P,Q} \circ \big( (\natTrans, \natCell) \times P \big) 
\to (\altNat, \altCell)$ 
the modification 
$\trans\modif$ of Lemma~\ref{lem:transModif-a-modification} is the unique such 
filling~(\ref{eq:transModif-ump}). 
\end{prooflesslemma}

Putting together everything we have seen in this section, for every 
$P, Q : \baseCat \to \Cat$ 
the 
pseudofunctor 
$\altexp{P}{Q} := \Hom(\baseCat, \Cat){\big(\Yon(-) \times P, Q\big)}$ 
satisfies an adjoint equivalence 
\[
\Lambda : 
\big(\Hom(\baseCat, \Cat)(R \times P, Q)\big) \leftrightarrows 
\big( \Hom(\baseCat, \Cat)(R, \altexp{P}{Q}) \big) : \eval_{P,Q} \circ (- 
\times P)
\] with evaluation map defined as in Lemma~\ref{lem:eval-defn} and 
counit $\evalMod$ defined as in Lemma~\ref{lem:counit-defn}. The 
universality of the 
counit is witnessed 
by the mapping $\trans{(-)}$ of Lemma~\ref{lem:transModif-a-modification}. 
Moreover, it is clear that $\trans\modif$ is invertible if $\modif$ is, so in 
particular the unit is invertible.
Thus, $\altexp{P}{Q}$ is an 
exponential in 
$\Hom(\baseCat, \Cat)$. 

\begin{prooflesspropn} \label{prop:cc-structure-of-Hom(baseCat,Cat)}
For any 2-category $\baseCat$ and pseudofunctors $P, Q : \baseCat \to \Cat$, 
the 
exponential $\altexp{P}{Q}$ exists and may be given by 
$\Hom(\baseCat, \Cat){\big(\Yon(-) \times P, Q\big)}$. 
\end{prooflesspropn}

Hence, $\Hom(\baseCat, \Cat)$ is cartesian closed for any 2-category 
$\baseCat$. Applying Lemma~\ref{lem:ccc-and-biequivalent-implies-ccc} 
yields our final result.

\begin{prooflessthm} \label{thm:2-presheaves-cartesian-closed}
For any bicategory $\baseCat$, the 2-category $\Hom(\baseCat, \Cat)$ is 
cartesian closed. 
\end{prooflessthm}

\section{Exponentiating by a representable} \label{sec:exponentiating-by-Yon}

For any 2-category $\baseCat$ with pseudo-products, object $X \in \baseCat$ and 
pseudofunctor
$P : \op\baseCat \to \Cat$, the exponential $\altexp{\Yon X}{P}$ may be given 
as $P(- \times X)$. This follows immediately from the the uniqueness 
of exponentials up to equivalence (Remark~\ref{rem:exponentials-up-to-equiv}), 
together with the following chain of equivalences:
\begin{equation} \label{eq:exp-by-rep}
\begin{aligned}
\altexp{\Yon X}{P} 
	&\simeq \Hom(\baseCat, \Cat){\big(\Yon(-) \times \Yon X, P\big)} 
		&\text{by Proposition~\ref{prop:cc-structure-of-Hom(baseCat,Cat)}}  \\ 
	&\simeq \Hom(\baseCat, \Cat){\big(\Yon(- \times X), P\big)} \\
	&\simeq P(- \times X) 
		&\text{by the Yoneda Lemma}
\end{aligned}
\end{equation}
For the second line we use the fact that birepresentables preserve bilimits 
(Lemma~\ref{lem:representables-and-adjoints-preserve-bilimits}). 

In the 
normalisation-by-evaluation argument (Chapter~\ref{chap:nbe}) we shall require 
an explicit description 
of the evaluation map 
witnessing $P(- \times X)$ as the exponential $\altexp{\Yon X}{P}$. In this 
section, therefore, we outline the exponential structure of $P(- \times X)$ and 
briefly show that it satisfies the required universal property. Since this 
structure may be extracted from the work of the preceding section by chasing 
through the equivalences~(\ref{eq:exp-by-rep}), our presentation will be less 
detailed than before.

Note that, for the rest of this chapter, we work contravariantly. Since we are 
assuming $\baseCat$ is a 2-category, the Yoneda pseudofunctor is now both 
strict (in fact, a 2-functor) and contravariant:  
$\Yon X = \op{\baseCat}(X, -) = \baseCat(-, X)$.

\paragraph{The evaluation map.}
We begin with the pseudonatural transformation 
${P(- \times X)} \times 
\Yon X \To P$ that will act as the evaluation map. For the component at 
$B \in \baseCat$ we take the functor
\begin{align*}
e_B : P(B \times X) \times \baseCat(B, X) &\to PB \\
(p, h) &\mapsto P(\seq{\Id_B, h})(p)
\end{align*}
with the evident action on 2-cells. To turn this into a pseudonatural 
transformation we need to provide an invertible 2-cell 
$\overline{e}_f$ as in the diagram below for every $f : B' \to B$ in $\baseCat$:
\begin{td}[column sep = 7em]
P(B \times X) \times \baseCat(B, X) 
\arrow[phantom]{dr}[description, xshift=-6mm]{\twocell{\overline{e}_f}}
\arrow[swap]{d}{e_B}
\arrow{r}{P(f \times X) \times \baseCat(f, X)} &
P(B' \times X) \times \baseCat(B', X)
\arrow{d}{e_{B'}} \\

PB \arrow[swap]{r}{Pf} &
PB'
\end{td}
At $h : B \to X$ we define 
$\overline{e}_f(h, -)$ to be the 
composite
\begin{td}[column sep = 4em]
P(\seq{\Id_B, h \circ f}) \circ P(f \times X)
\arrow[swap]{d}{\phi^P_{\seq{\Id, hf}, f \times X}}
\arrow{r}{\overline{e}_f(h, -)} &
P(f) \circ P\seq{\Id_B, h} \\

P\big( (f \times X)\seq{\Id_B, hf}\big)
\arrow{r}[swap]{P\swap_{h,f}} &
P(\seq{\Id_{B'},h} \circ f)
\arrow[swap]{u}{(\phi^P_{\seq{\Id, h}, f})^{-1}}
\end{td}
where the isomorphism $\swap_{h,f}$ is 
$(f \times X) \circ \seq{\Id_B, hf} 
\XRA{\fuse} 
\seq{f, hf} 
\XRA{\postName^{-1}} 
\seq{\Id_{B'}, h} \circ f$.
\nom{\swap_{h,f}}{The 2-cell of type
$(f \times X) \circ \seq{\Id_B, hf} \To \seq{\Id_{B'}, h} \circ f$, defined as 
the composite
$(f \times X) \circ \seq{\Id_B, hf} 
\XRA{\fuse} 
\seq{f, hf} 
\XRA{\postName^{-1}} 
\seq{\Id_{B'}, h} \circ f$
} The 
whole composite is a natural isomorphism because each component is, so 
it remains to check the two axioms of a pseudonatural transformation. 
The unit law is a short diagram chase using the unit law for $P$
and the fact that
\[
\Id_{B\times X} \circ \seq{\Id_B, h}
	\XRA{\etaTimes{\Id} \circ \seq{\Id, h}}
\seq{\Id_B, h} \circ \Id_{B}
	\XRA{\swap}
\Id_{B\times X} \circ \seq{\Id_B, h}
\]
is the identity.
To prove the associativity law, on the other hand, one uses the naturality of 
the $\phi^P$ 2-cells and the 
associativity law of a pseudofunctor to reduce the problem to a diagram 
in the image of $P$, whereupon one can apply standard properties of the product 
structure
(recall Lemma~\ref{lem:PseudoproductCanonical2CellsLaws}).

\begin{prooflesslemma} \label{lem:YX-exponential-eval-map}
For any $X \in \baseCat$ and pseudofunctor $P : \op\baseCat \to \Cat$, the  
pair $(e, \overline{e})$ defined above forms a pseudonatural transformation 
$P(- \times X) \times \Yon X \To P$. 
\end{prooflesslemma}

\paragraph{The mapping $\Lambda$.} 
Next we define the mapping
\[
\Lambda : ob\big(\Hom(\op\baseCat, \Cat)( R \times \Yon X, P)\big) \to 
ob\big(\Hom(\op\baseCat, \Cat)( R , {P(- \times X)})\big)
\] 
Let $(\natTrans, \natCell) : R \times \Yon X \To P$ be a pseudonatural 
transformation. We define 
$\Lambda(\natTrans, \natCell) := 
	(\Lambda \natTrans, \cellOf{\Lambda \natTrans}) : 
	R \To P(- \times X)$ as 
follows. For 
$B \in \baseCat$ we take the functor
\begin{align*}
(\Lambda \natTrans)_B : RB &\to P(B \times X) \\
r &\mapsto \natTrans_{B \times X}\big(R(\pi_1)(r), \pi_2 \big)
\end{align*}
Thus, $(\Lambda \natTrans)_B$ is the 
composite 
$RB 	\xra{R\pi_1} 
	R(B \times X) 
		\xra{\natTrans_{B \times X}(-, \pi_2)} 
	P(B \times X)$. 
To define 
$\cellOf{(\Lambda \natTrans)}_f$, where $f : B' \to B$, we need to give an 
invertible 2-cell as in 
\begin{td}[column sep = 3.5em]
RB \arrow{r}{Rf}
\arrow[phantom]{dr}[description, xshift=2mm]
	{\twocell{\cellOf{(\Lambda \natTrans)}_f}}
\arrow[swap]{d}{(\Lambda \natTrans)_{B}} &
RB'
\arrow{d}{(\Lambda \natTrans)_{B'}} \\
P(B \times X) 
\arrow[swap]{r}{P(f \times X)} &
P(B' \times X) 
\end{td}
This must be a natural isomorphism
$\natTrans_{B' \times X}\big(R(\pi_1)R(f)(-), \pi_2 \big) 
\XRA{\iso} 
P(f \times X)\big( \natTrans_{B \times X}(R(\pi_1)(-), \pi_2 ) \big)$, for 
which we take 
the following composite:
\begin{td}[column sep = 3em]
\natTrans_{B' \times X}(R(\pi_1) \circ R(f), \pi_2 ) 
\arrow{r}{(\cellOf{\Lambda \natTrans})_f} 
\arrow[swap]{d}{\natTrans_{B' \times X}(\phi^R_{f, \pi_1}, \pi_2)} &
P(f \times X)\big( \natTrans_{B \times X}(R\pi_1, \pi_2 ) \big) \\

\natTrans_{B' \times X}\big( R(f \circ \pi_1), \pi_2\big)
\arrow[swap]{d}
	{\natTrans_{B' \times X}
		(R\epsilonTimesInd{-1}{}, 
			\epsilonTimesInd{-2}{})} &
\: \\

\natTrans_{B' \times X}\big(R\left(\pi_1(f \times X)\right), \pi_2(f \times 
X)\big) 
\arrow[swap]{r}[yshift=-2mm]
	{\natTrans_{B' \times X}
		((\phi_{\pi_1, f \times X}^R)^{-1}, \pi_2(f \times 	X))} &
\natTrans_{B' \times X}\big(R(f \times X) \circ R(\pi_1), 
							\pi_2(f \times X)\big) 
\arrow[swap]{uu}{\natCell_{f \times X}(R\pi_1, \pi_2)}
\end{td}
To see that this is a pseudonatural transformation, observe that we have 
actually defined 
$\Lambda (\natTrans, \natCell)$ as a composite
\begin{equation} \label{eq:lamda-nat-trans-as-composite}
\begin{tikzcd}[row sep = 3em]
RB 
\arrow[swap]{d}{\natTransAux_B}
\arrow{r}{Rf}
\arrow[phantom]{dr}[description, xshift=2.5mm]{\twocell{\natCellAux_f}} &
RB'
\arrow{d}{\natTransAux_{B'}} \\
R(B \times X) \times \baseCat(B \times X, X)
\arrow[swap]{r}[yshift=-2mm]{R(f \times X) \times \baseCat(f \times X, X)} 
\arrow[swap]{d}{\natTrans_{B \times X}}
\arrow[phantom]{dr}[description, yshift =-2mm]{\twocell{\natCell_{f \times X}}} 
&
R(B' \times X) \times \baseCat(B' \times X, X) 
\arrow{d}{\natTrans_{B' \times X}} \\
P(B \times X)  
\arrow[swap]{r}{P(f \times X)} &
P(B' \times X)
\end{tikzcd}
\end{equation}
where $\natTransAux_B(r) := \big( R(\pi_1)(r), \pi_2 \big)$ and $\natCellAux_f$ 
has first component
\begin{equation} \label{eq:r-two-cell-defined} 
R\pi_1 \circ Rf  
\XRA{\phi^R_{f, \pi_1}}
R(f \circ \pi_1) 
\XRA{R\epsilonTimesInd{-1}{}} 
R\big( \pi_1 \circ (f \times X) \big) 
\XRA{(\phi^R_{\pi_1, f \times X})^{-1}} 
R(f \times X) \circ R\pi_1
\end{equation}
and second component 
$\pi_2 \XRA{\epsilonTimesInd{-2}{}} \pi_2 \circ (f \times X)$. So it suffices 
to show that $(\natTransAux, \natCellAux)$ defines a pseudonatural 
transformation 
$R \To R(- \times X) \times \baseCat(- \times X, X)$. Naturality follows 
immediately from the fact each component in the definition is natural. For the 
unit law, the first component is the triangle law for products, and the second 
component is a short diagram chase. 
%
%

For the associativity law, 
it is once again the second component that is more difficult. 
%
%
%
%
As for $(e, \overline{e})$~(Lemma~\ref{lem:YX-exponential-eval-map}), the proof 
consists of using the associativity axiom 
of a pseudofunctor and the naturality of $\phi^R$. Once the calculation has 
been pushed `inside' $R$, what remains is a relatively easy diagram chase. This 
completes 
the proof that $(\natTransAux, \natCellAux)$ is a pseudonatural transformation, 
and 
hence the definition of the mapping $\Lambda$. 

\begin{prooflesslemma} 
The pair $(\natTransAux, \natCellAux)$ defined in~(\ref{eq:r-two-cell-defined}) 
forms a pseudonatural transformation 
$R \To {R(- \times X) \times {\baseCat(- \times X, X)}}$.
\end{prooflesslemma}

\begin{prooflesscor} \label{cor:YX-exponentiating-lambda-mapping}
The pair $(\Lambda \natTrans, \cellOf{\Lambda \natTrans})$ defined 
in~(\ref{eq:lamda-nat-trans-as-composite}) forms a 
pseudonatural transformation $R \To P(- \times X)$ for every 
$(\natTrans, \natCell) : R \times \Yon X \To P$. 
\end{prooflesscor}

\paragraph{The counit $\evalMod$.}
For every 
$(\natTrans, \natCell) : R \times \Yon X \To P$ we need to provide an 
invertible modification 
$\evalMod^{(\natTrans, \natCell)} : 
(e, \overline{e}) \circ \big( \Lambda(\natTrans, \natCell) \times \Yon X \big) 
\to (\natTrans, \natCell)$. 

Unwrapping the definition of 
$(e, \overline{e}) \circ \big( \Lambda(\natTrans, \natCell) \times \Yon X 
\big)$ 
at $B \in \baseCat$ and $(r, h) \in RB \times \baseCat(B, X)$, one sees that
\begin{align*}
\Big(e_B \circ \big( (\Lambda\natTrans)_B \times \Yon X \big)\Big)(r,h) &= 
		e_B\big( \natTrans_{B \times X}( R(\pi_1)(r), \pi_2 ) , h \big) \\
		&= P\left(\seq{\Id_B, h}\right)
			\big(\natTrans_{B \times X}( R(\pi_1)(r), \pi_2 )\big)
\end{align*}
Furthermore, for $f : B' \to B$ the corresponding 2-cell 
$\overline{\big(e_B \circ \big( (\Lambda\natTrans)_B 
\times \Yon X \big)\big)}_f$ is defined by
\begin{equation*}
\makebox[\textwidth]{
\begin{tikzcd}[column sep = 1.8em, row sep = 2.5em, ampersand replacement = \&]
P(\seq{\Id_B, hf})
\big(\natTrans_{B' \times X}( R(\pi_1) R(f)(r), \pi_2 )\big) 
\arrow{r}[yshift=2mm]
	{\overline{\left(e_B \circ 
		\left( (\Lambda\natTrans)_B \times \Yon X \right)\right)}_f(r,h)}
\arrow[swap]{d}{P(\seq{\Id_B, hf})(\altCell_f(r))} \&
P(f)P(\seq{\Id_B, h})\big(\natTrans_{B \times X}(R(\pi_1)(r), \pi_2)\big) \\
P(\seq{\Id_B, hf})
	\big(\natTrans_{B \times X}(R(f \times X)R(\pi_1)(r),\pi_2(f\times X)\big) 
\arrow[swap]{r}[yshift=-2mm]
	{P(\seq{\Id_B, hf})(\natCell_{f \times X}(R(\pi_1)(r), \pi_2))} \&
P(\seq{\Id_B, hf})P(f \times X)
	\big(\natTrans_{B \times X}(R(\pi_1)(r), \pi_2)\big) 
\arrow[swap]{u}{\overline{e}_f(h, \natTrans_{B \times X}(R(\pi_1)(r), \pi_2))}
\end{tikzcd}
}
\end{equation*}
We therefore take the component at $B \in \baseCat$ of 
$\evalMod_B^{(\natTrans, \natCell)}$ to be 
the natural isomorphism defined by
\begin{equation} \label{eq:Yon-X-exponential-counit}
\begin{tikzcd}[column sep = 1.5em, row sep = 2.5em]
P(\seq{\Id_B, h})
\big(\natTrans_{B' \times X}( R(\pi_1)(r), \pi_2 )\big) 
\arrow[swap]{d}{\natCell_{\seq{\Id, h}}^{-1}(R(\pi_1)(r), \pi_2) } 
\arrow{rr}{\evalMod_B^{(\natTrans, \natCell)}(r,h)} &
\:  &
\natTrans_B(r, h) \\
\natTrans_B\big(  R(\seq{\Id_B, h}) R(\pi_1)(r) , \pi_2\seq{\Id_B, h} \big)
\arrow[swap]{r}[yshift=-2mm]
	{\natTrans_B(\phi^R_{\pi_1, \seq{\Id, h}}(r), \epsilonTimesInd{2}{})} &
\natTrans_B\big( R(\pi_1 \seq{\Id_B, h})(r) , h \big) 
\arrow[swap]{r}[yshift=-2mm]{\natTrans_B(R\epsilonTimesInd{1}{}, h)} &
\natTrans_B\big(R(\Id_B)(r), h\big)
\arrow[swap]{u}{\natTrans_B((\psi^R_B)^{-1}, h)}
\end{tikzcd}
\vspace{2mm}
\end{equation}
We need to check the $B$-indexed family of 2-cells 
$\evalMod^{(\natTrans, \natCell)}$ satisfies the
modification axiom, namely that
\begin{td}[column sep = 9em, row sep = 2.5em]
P(\seq{\Id_B, hf}) \big( \natTrans_{B \times X}(R(\pi_1)R(f)(r), \pi_2) \big) 
\arrow{r}[yshift=0mm]{\evalMod^{(\natTrans, \natCell)}_{B}(R(\pi_1)R(f)(r), 
\pi_2)} 
\arrow[swap]{d}{\overline{\big(e_B \circ \big( (\Lambda\natTrans)_B 
\times \Yon X \big)\big)}_f(r,h)} & 
\natTrans_B \big( R(f) (r), hf \big) 
\arrow{d}{\natCell_f(r,h)} \\

P(f)P(\seq{\Id_B, h})\big(\natTrans_{B \times X}(R(\pi_1)(r), \pi_2)\big) 
\arrow[swap]{r}{P(f)(\evalMod^{(\natTrans, \natCell)}_{B}(r,h))} &
P(f)\big( \natTrans_B(r,h) \big)
\end{td}
Unfolding all the data results in a long exercise in diagram chasing. The 
second component is relatively straightforward. For the first component, one 
applies the naturality properties and associativity law of a pseudofunctor to 
reduce the claim to the following:
\begin{td}[column sep = 1.3em]
R(\seq{\Id_{B'}, hf} ) \circ R(\pi_1) \circ R(f) 
\arrow{r}[yshift=2mm]{\phi^R_{\pi_1, \seq{\Id, hf}} \circ R(f)}
\arrow[swap]{d}{R(\seq{\Id_{B'}, hf} ) \circ  \phi^R_{f, \pi_1}} &
R(\pi_1 \seq{\Id_{B'}, hf}) \circ R(f) 
\arrow{r}[yshift=2mm]{R(\epsilonTimesInd{1}{}) \circ R(f)} &
R(\Id_{B'}) \circ R(f) & \\

R(\seq{\Id_{B'}, hf} ) \circ R(f \circ \pi_1) 
\arrow[swap]{d}{R(\seq{\Id_{B'}, hf} ) \circ  R(\epsilonTimesInd{-1}{})} &
\: &
R(f) 
\arrow[swap]{u}{\psi^R_{B'} \circ R(f)} \\

R(\seq{\Id_{B'}, hf} ) \circ R\big( \pi_1 \circ (f \times X) \big) 
\arrow[swap]{d}{R(\seq{\Id_{B'}, hf} ) \circ  (\phi^R_{\pi_1, f \times 
X})^{-1}} &
\: &
R(\Id_B \circ f) 
\arrow[equals]{u} \\

R(\seq{\Id_{B'}, hf} ) \circ R(f \times X) \circ R(\pi_1) 
\arrow[swap]{d}{\phi^R_{f \times X, \seq{\Id, hf}} \circ R(\pi_1)} &
\: &
R(\pi_1 \circ \seq{\Id_B, h} \circ f) 
\arrow[swap]{u}{R(\epsilonTimesInd{1}{} \circ f)} \\

R\big((f \times X) \circ \seq{\Id_{B'}, hf} \big) \circ R(\pi_1)
\arrow[swap]{r}[yshift=-2mm]{R(\fuse) \circ R(\pi_1)} &
R(\seq{f, hf}) \circ R(\pi_1) 
\arrow[swap]{r}[yshift=-2mm]{R(\postName^{-1}) \circ R(\pi_1)} &
R(\seq{\Id_B, h} \circ f) \circ R(\pi_1) 
\arrow[swap]{u}{\phi^R_{\pi_1, \seq{\Id_B, h} \circ f}}
\end{td}
The strategy is now familiar: one applies naturality and the associativity law 
to bring together all the morphisms in the image of $R$, and then unwraps the 
definition of $\postName$ and $\fuse$ to reduce the long anticlockwise claim to 
the top row.

We have therefore constructed a modification to act as the counit.

\begin{prooflesslemma} \label{lem:YX-exponentiating-counit-2-cell}
The 2-cells $\evalMod_B^{(\natTrans, \natCell)} \:\: (B \in \baseCat)$ defined 
in~(\ref{eq:Yon-X-exponential-counit}) 
form an invertible modification 
$(e, \overline{e}) \circ \big( \Lambda(\natTrans, \natCell) \times \Yon X \big) 
\to (\natTrans, \natCell)$.
\end{prooflesslemma}

All that remains is to show the modification $\evalMod^{(\natTrans, \natCell)}$ 
is a universal arrow. 

\paragraph{The modification $\trans\modif$.} 
We aim to construct a modification $\trans\modif$ for every 
pseudonatural transformation $(\altNat, \altCell) : R \To P(- \times X)$ and 
modification 
$\modif : (e, \overline{e}) \circ \big( (\altNat, \altCell) \times \Yon X \big) 
\to (\natTrans, \natCell)$, such that $\trans\modif$ is the unique modification 
filing
\begin{equation} \label{eq:trans-modif-ump-YX-exponential}
\begin{tikzcd}
(e, \overline{e}) \circ \big( (\altNat, \altCell) \times \Yon X \big)
\arrow[swap]{dr}{\modif}
\arrow{rr}{(e, \overline{e}) \circ ( \trans\modif \times \Yon X )} &
\: &
(e, \overline{e}) \circ \big( \Lambda(\natTrans, \natCell) \times \Yon X \big) 
\arrow{dl}{\evalMod^{(\natTrans, \natCell)}} \\
\: &
(\natTrans, \natCell) &
\:
\end{tikzcd}
\end{equation}
Because the definitions of $(e, \overline{e})$, $\Lambda(\natTrans, \natCell)$ 
and $\evalMod^{(\natTrans, \natCell)}$ are all composites, the proof requires 
working with a large 
accumulation of data. Nonetheless the diagram chases---although long---are not 
especially difficult.  

Suppose that 
$\modif : 
(e, \overline{e}) \circ \big( (\altNat, \altCell) \times \Yon X \big) \to 
(\natTrans, \natCell)$. Since 
\[
\big(e_B \circ ( \altNat_B \times \Yon X)\big)(r, h) =  
e_B(\altNat_B(r), h) 
= P(\seq{\Id_B, h})(\altNat_B(r))
\]
for every $B \in \baseCat$ we are provided 
with a 
natural transformation with components 
$\modif_B(r,h) : (P\seq{\Id_B, h})(\altNat_B(r)) \to \natTrans_B(r,h)$ for 
$(r, h) \in RB \times \baseCat(B, X)$.  
We define $\trans{\modif}_B$ to be the 
composite
\begin{equation}  \label{eq:Yon-X-exponential-trans-modif}
\begin{tikzcd}[column sep = 1.9em]
\altNat_B 
\arrow[swap]{d}{\psi^P_{B \times X} \circ \altNat_B}
\arrow{rr}{\trans{\modif}_B} &
\: &
\natTrans_{B\times X}(R\pi_1, \pi_2) \\
P(\Id_{B\times X}) \circ \altNat_B 
\arrow[swap]{d}{P(\etaTimes{\Id}) \circ \altNat_B} &
\: &
P(\seq{\Id_{B\times X}, \pi_2}) \circ \altNat_{B \times X} \circ R\pi_1 
\arrow[swap]{u}{\modif_B(R\pi_1, \pi_2)} \\
P(\seq{\pi_1, \pi_2}) \circ \altNat_B 
\arrow[swap]{r}[yshift=-2mm]{P(\fuse^{-1}) \circ \altNat_B} &
P\big( (\pi_1 \times X) \circ \seq{\Id_{B\times X}, \pi_2}) \circ \altNat_B
\arrow[swap]{r}[yshift=-2mm]
	{(\phi^P_{\pi_1 \times X, \seq{\Id, \pi_2}})^{-1} \circ \altNat_B} &
P(\seq{\Id_{B\times X}, \pi_2}) \circ P(\pi_1 \times X) \circ \altNat_{B}  
\arrow[swap]{u}{P(\seq{\Id_{B\times X}, \pi_2}) \circ \altCell^{-1}_{\pi_1}}
\end{tikzcd}
\vspace{2mm}
\end{equation}
and claim this does indeed define a modification. We therefore need to 
verify the following diagram of functors commutes for every $f : B' \to B$ in 
$\baseCat$:
\begin{td}[column sep = 6em]
\altNat_{B'}\big(R(f)\big) 
\arrow{r}{\trans{\modif}_{B'}(R(f))}
\arrow[swap]{d}{\altCell_f} &
\natTrans_{B' \times X}\big( R(\pi_1) R(f) , \pi_2 \big) 
\arrow{d}{(\cellOf{\Lambda\natTrans})_f} \\

P(f \times X)\big(\altNat_B\big) 
\arrow[swap]{r}{P(f \times X)(\trans{\modif}_B)} &
P(f \times X)\big( \natTrans_{B \times X} (R(\pi_1), \pi_2) \big) 
\end{td}
Unfolding all the various composites results in a very large diagram. We give 
the strategy for proving it commutes. One begins by using naturality until one 
can apply the modification axiom 
for $\modif$ to relate the 
final term in the composite defining 
$\cellOf{(\Lambda \natTrans)}_f$ with 
$P(f \times X)\big(\modif_{B \times X}(R(\pi_1)(r), \pi_2)\big)$. Next one 
applies the associativity law for $(\altNat, \altCell)$ in order to push the 
2-cells $\phi^P$ as early as possible. One then observes that the 
following diagram commutes, and hence that its image under 
$P$ commutes:
\begin{td}[column sep = 1em]
f \times X 
\arrow[swap]{d}{(f \times X) \circ \etaTimes{\Id}}
\arrow{r}{\etaTimes{\Id} \circ (f \times X)} &
\seq{\pi_1, \pi_2} \circ (f \times X) \\

(f \times X) \circ \seq{\pi_1, \pi_2} 
\arrow[swap]{d}{(f \times X) \circ \seq{\pi_1, \epsilonTimesInd{-2}{}}} &
\: \\

(f \times X) \circ \seq{\pi_1, \pi_2(f \times X)}
\arrow[swap]{d}{(f \times X) \circ \fuse^{-1}}  &
(\pi_1 \times X) \circ \seq{\Id_{B \times X}, \pi_2} \circ (f \times X)
\arrow[swap]{uu}{\fuse \circ (f \times X)} \\

(f \times X) \circ (\pi_1 \times X) \circ 
\seq{\Id_{B' \times X}, \pi_2(f \times X)} 
\arrow[swap]{d}{\phiTimes_{f, \pi_1; \Id_X} \circ \seq{\Id, \pi_2(f \times X)}} 
&
(\pi_1 \times X) \circ \big((f \times X) \times X \big) \circ 
\seq{\Id_{B' \times X}, \pi_2(f \times X)} 
\arrow[swap]{u}{(\pi_1 \times X) \circ \swap}  \\

\big((f \circ \pi_1) \times X \big) \circ 
\seq{\Id_{B' \times X}, \pi_2(f \times X)} 
\arrow[swap]{r}[yshift=-2mm]{(\epsilonTimesInd{-1}{} \times X) \circ \seq{\Id, 
\pi_2(f \times X)}} &
\big((\pi_1(f \times X)) \times X \big) \circ 
\seq{\Id_{B' \times X}, \pi_2(f \times X)} 
\arrow[swap]{u}
	{\phiTimes_{\pi_1, f \times X; \Id_X}^{-1} \circ 
		\seq{\Id, \pi_2(f \times X)}}
\end{td}
From this point the rest of the proof is a manageable diagram chase. Hence, 
$\trans\modif$ is a modification. 

\begin{prooflesslemma}
For every modification $\modif : 
(e, \overline{e}) \circ \big( (\altNat, \altCell) \times \Yon X \big) \to 
(\natTrans, \natCell)$ between pseudonatural transformations $R \times \Yon X 
\To P$, the 2-cells $\trans{\modif}_B$ form a modification $(\altNat, \altCell) 
\to \Lambda(\natTrans, \natCell)$. 
\end{prooflesslemma}

The last part of the proof is checking that $\trans\modif$ is the unique 
modification filling the diagram~(\ref{eq:trans-modif-ump-YX-exponential}). 

\paragraph{The universal property of $\evalMod$.}

The existence and uniqueness parts of~(\ref{eq:trans-modif-ump-YX-exponential}) 
also entail long but not especially difficult diagram chases. 
In each case one unfolds the various composites and applies the modification 
axiom for $\modif$. The rest of the proof is an 
exercise in applying the various naturality properties and the two laws of a 
pseudofunctor. 

Putting together all the work of this section, one obtains the following. 

\newpage
\begin{prooflesspropn}
For any 2-category $\baseCat$ with pseudo-products, pseudofunctor $P : 
\op\baseCat \to \Cat$ and 
object $X \in \baseCat$, the modification $\evalMod$ of 
Lemma~\ref{lem:YX-exponentiating-counit-2-cell} is the counit of an adjoint 
equivalence
\[\Lambda : \Hom(\op\baseCat, \Cat)(R \times \Yon X, P) 
\leftrightarrows \Hom(\op\baseCat, \Cat)\big(R, P(- \times X) \big) : 
(e, \overline{e}) \circ (- \times \Yon X)
\]
in which the pseudonatural transformation $(e, \overline{e})$ and mapping 
$\Lambda$ are as in Lemma~\ref{lem:YX-exponential-eval-map} 
and Corollary~\ref{cor:YX-exponentiating-lambda-mapping}, 
respectively. 
\end{prooflesspropn}

\begin{prooflessthm} \label{thm:exponentiating-by-Yoneda-identified}
For any 2-category $\baseCat$ with pseudo-products, pseudofunctor $P : 
\op\baseCat \to \Cat$ and 
object $X \in \baseCat$, the pseudofunctor $P(- \times X)$ is 
(up to equivalence) the exponential 
$\altexp{\Yon X}{P}$ in $\Hom(\op\baseCat, \Cat)$.
\end{prooflessthm}

Setting $\altCat := \op\baseCat$ recovers the covariant statement.

%

\chapter{Bicategorical glueing} 
\label{chap:glueing}

Glueing is a powerful technique which may be used to leverage semantic 
arguments in order to prove syntactic results.  
Intuitively, one `glues together' syntactic and semantic 
information, allowing one to extract proofs of syntactic properties from 
semantic arguments. The breadth and 
utility of this approach has led to its being discovered in various forms, 
with correspondingly various names: the notions of logical 
relation~\cite{Plotkin1973, Statman1985}, 
sconing~\cite{Freyd1990}, Freyd covers and glueing (\eg~\cite{LambekAndScott}) 
are all closely 
related (see~\eg~\cite{Mitchell1993} for an overview of the connections). 
Taylor  
identifies the basic apparatus as going back to 
Groethendieck~\cite[Section~7.7]{Taylor1999}, while versions of logical 
relations appear as early as Gandy's thesis (who, in turn, attributes some of 
the theory to Turing)~\cite{Gandy1953}. Originally presented 
in the 
set-theoretic setting, the technique was quickly given categorical 
expression~\cite{Ma1992, Mitchell1993}, for which Hermida provided an account 
in terms of fibrations in his thesis~\cite{Hermida1993}. Such techniques are 
now a standard component of the armoury for studying type theories.

In this chapter we define a notion of glueing for bicategories and 
prove a bicategorical version of the fundamental result establishing mild 
conditions for the glueing category to be cartesian closed. (For reference, the 
construction 
is summarised in the appendix on 
page~\pageref{table:glueing-ccc-structure}.) 
This will form the core of our normalisation-by-evaluation proof in the next 
chapter. 

We begin by recalling the categorical glueing construction and giving a precise 
statement of the cartesian closure result we wish to prove. These will provide 
a template for our bicategorical work.

\section{Categorical glueing}

The most succinct description of 
categorical glueing is as a special kind of comma 
category.

\begin{mydefn} \label{def:1-categorical-glueing}
\quad
\begin{enumerate}
\item
Let $F : \catA \to \catC$ and $G : \catB \to \catC$ be functors. The 
\Def{comma category} 
$(F \downarrow G)$ has objects triples $(A, f, B)$, where  
$A \in \catA$ and $B \in \catB$ are objects and $f : FA \to GB$ is a morphism 
in $\catC$. Morphisms 
$(A, f, B) \to (A', f', B')$
are pairs of morphisms $(p,q)$ such that the following square commutes:
\begin{equation} \label{eq:glueing-1-categorical-morphism-square}
\begin{tikzcd}
FA
\arrow{r}{Fp}
\arrow[swap]{d}{f} &
FA' 
\arrow{d}{f'} \\
GB
\arrow[swap]{r}{Gq} &
GB' 
\end{tikzcd}   
\end{equation}

\item
The \Def{glueing} $\gl{\glueFun}$ of $\catB$ to $\catC$ along a functor 
$\glueFun : \catB \to \catC$ is the comma 
category $(\id_\catC \downarrow \glueFun)$. We denote the objects and 
morphisms following the 
vertical order 
of their appearance in 
diagram~(\ref{eq:glueing-1-categorical-morphism-square}), as 
$(C \in \catC, c : C \to \glueFun B, B \in \catB)$ and 
$(q : C \to C', p : B \to B')$. \qedhere
\end{enumerate}
\end{mydefn}

There are evident \Def{projection functors} 
$\catB \xla{\pi_{\mathrm{dom}}} \gl{\glueFun} \xra{\pi_{\mathrm{cod}}} \catC$.  
We wish to bicategorify the following folklore result~(\cf~\cite[Proposition 
2]{Ma1992}):

\begin{mypropn} \label{prop:glueing-cat-ccc}
Let $\glueFun : \catB \to \catC$ be a functor between cartesian closed 
categories, such that $\glueFun$ preserves products and $\catC$ has all 
pullbacks. Then the glueing category $\gl{\glueFun}$ is cartesian closed, and 
the 
projection $\pi_{\mathrm{dom}}$ strictly preserves the cartesian closed 
structure.
\begin{proof}
For $n \in \Nat$ the $n$-ary product of objects 
$(C_i, c_i, B_i) \:\: (i=1, \,\dots\, , n)$ is 
the composite 
\[
\prodop_{i=1}^n C_i \xra{\prod_i c_i} \prodop_{i=1}^n( \glueFun B_i) \xra{\iso} 
\glueFun \big(\prodop_{i=1}^n B_i\big)
\]
Projections are given pointwise, as $(\pi_i^\catC, \pi_i^\catB)$, and  
the $n$-ary tupling of a family of 1-cells
$(f_i, g_i) : (X, x, Y) \to (C_i, c_i, B_i) \:\: (i=1, \,\dots\, , n)$  
is the pair $(\seq{f_1, \,\dots\, , f_n}, \seq{g_1, \,\dots\, , g_n})$. Hence 
both 
$\pi_{\mathrm{dom}}$ and $\pi_{\mathrm{cod}}$ strictly preserve products. 

The exponential 
$\exp{(C, c, B)}{(C', c', B')}$ is defined to be the left-hand vertical map in 
the pullback diagram
\begin{equation} \label{eq:exponential-in-glued-cat}
\begin{tikzcd}[column sep = 3.5em]
\glexp{C}{C'} 
\arrow[phantom]{dr}[at start, description, yshift=2mm, xshift=-1mm]{\sidecorner}
\arrow[swap]{d}{p_{c,c'}}
\arrow{rr}{q_{c,c'}} &
\: &
(\exp{C}{C'}) 
\arrow{d}{\scriptsizeexpobj{C}{c'}} \\
\glueFun(\exp{B}{B'})
\arrow[swap]{r}{\evBar_{B,B'}} &
(\exp{\glueFun B}{\glueFun B'}) 
\arrow[swap]{r}{(\scriptsizeexpobj{c}{\glueFun B'})} &
(\exp{C}{\glueFun B'})
\end{tikzcd}
\vspace{1mm}
\end{equation}
where 
$\evBar_{B,B'}$ is the exponential transpose of
$\big(\glueFun (\exp{B}{B'}) \times \glueFun B 
	\xra\iso 
	\glueFun ((\exp{B}{B'}) \times B) 
	\xra{\glueFun \eval_{B,B'}} 
	\glueFun B' \big)$.
The evaluation map has first component
$(\glexp{C}{C'}) \times C \xra{q_{c,c'} \times C} (\exp{C}{C'}) \times C 
\xra{\eval_{C,C'}} C'$
and second component simply $\eval_{B,B'}$. The currying operation is given by 
the universal property of pullbacks. 
\end{proof}
\end{mypropn}

The rest of the chapter is dedicated to proving a bicategorical version of this 
proposition.

%
%

\section{Bicategorical glueing}

We bicategorify Definition~\ref{def:1-categorical-glueing} in the usual way: by 
replacing commuting squares with invertible 2-cells, subject to coherence 
conditions.

\begin{mydefn} \label{def:comma-bicategory}
Let $F  : \bicatA \to \altCat$ and $G : \baseCat \to \altCat$ be pseudofunctors 
of 
bicategories. The 
\Def{comma bicategory} $(F \downarrow G)$ has objects triples 
$(A \in \bicatA, f : FA \to GB, B \in \baseCat)$. The 1-cells 
$(A, f, B) \to (A', f', B')$ are triples 
$(p, \alpha, q)$, where $p : A \to A'$ and $q : B \to B'$ are 1-cells and 
$\alpha$ is an invertible 
2-cell $\alpha : f' \circ Fp \To Gq \circ f$ witnessing the commutativity 
of~(\ref{eq:glueing-1-categorical-morphism-square}):
\begin{equation} \label{eq:glueing-bicategorical-morphism-square}
\begin{tikzcd}
FA
\arrow[phantom]{dr}[description]{\twocell{\alpha}}
\arrow{r}{Fp}
\arrow[swap]{d}{f} &
FA' 
\arrow{d}{f'} \\
GB
\arrow[swap]{r}{Gq} &
GB' 
\end{tikzcd}   
\end{equation}
The 2-cells $(p,\alpha,q) \To (p', \alpha', q')$ are pairs of 2-cells 
$(\sigma : p \To p', \tau : q \To q')$ such that the following diagram 
commutes:
\begin{equation} \label{eq:glueing-2-cell-condition}
\begin{tikzcd}[column sep = 3.5em]
f' \circ F(p) 
\arrow[swap]{d}{\alpha}
\arrow{r}{f' \circ F(\sigma)} &
f' \circ F(p') \arrow{d}{\alpha'}  \\
G(q) \circ f 
\arrow[swap]{r}{G(\tau) \circ f} &
G(q') \circ f
\end{tikzcd}   
\end{equation}
The horizontal composite of 
$(A, f, B) \xra{(p, \alpha, q)} 
(A', f', B') \xra{(r, \beta, s)}
(A'', f'', B'')$ 
is
$(r \circ p, \iso, s \circ q)$, where the isomorphism is the composite on the 
left below:
\begin{center}
\begin{tikzcd}
f'' \circ F(r \circ p)
\arrow{r} 
\arrow[swap]{d}{f'' \circ (\phi^F_{r,p})^{-1}} &
G(s \circ q) \circ f \\
f'' \circ (Fr \circ Fp) 
\arrow[swap]{d}{\iso} &
(Gs \circ Gq) \circ f 
\arrow[swap]{u}{\phi^G_{s, q} \circ f} \\
(f'' \circ Fr) \circ Fp 
\arrow[swap]{d}{\beta \circ Fp} &
Gs \circ (Gq \circ f) 
\arrow[swap]{u}{\iso} \\
(Gs \circ f') \circ Fp 
\arrow[swap]{r}{\iso} &
Gs \circ (f' \circ Fp) 
\arrow[swap]{u}{Gs \circ \alpha} 
\end{tikzcd}
\quad\qquad
\begin{tikzcd}
f \circ F\Id_A
\arrow{r}
\arrow[swap]{d}{f \circ (\psi^F_A)^{-1}} &
G\Id_B \circ f \\
f \circ \Id_{FA} 
\arrow[swap]{r}{\iso} &
\Id_{GB} \circ f 
\arrow[swap]{u}{\psi^G_B \circ f}
\end{tikzcd}
\end{center}
In a similar fashion, the identity 1-cell on $(A, f, B)$ is
$(\Id_A, \iso, \Id_B)$ with isomorphism $\iso$ as on the right above.

Vertical composition and 
the identity 2-cell are given component-wise, as are the structural 
isomorphisms $\a, \l$ and $\r$. 
\end{mydefn} 
The identities and composition may be expressed as the following pasting 
diagrams:
\begin{equation*} 
\begin{tikzcd}
FA 
\arrow{dr}[description]{f}
\arrow{r}{F\Id_A} 
\arrow[swap]{d}{f} &
FA
\arrow[phantom]{dl}[description, font=\scriptsize, very near start]{\iso}
\arrow{d}{f} \\
GB 
\arrow[phantom]{ur}[description, font=\scriptsize, very near start]{\iso}
\arrow[swap]{r}{G\Id_B} &
GB
\end{tikzcd}  
\qquad\qquad\qquad
\begin{tikzcd}
\: &
\: &
\: \\
FA 
\arrow[bend left = 50]{rr}{F(r \circ p)}
\arrow[phantom]{dr}[description]{\twocell{\alpha}}
\arrow{r}{Fp}
\arrow[swap]{d}{f}  &
FA'
\arrow[phantom]{u}[description]{\twocellIso{\phi^F}}
\arrow[phantom]{dr}[description]{\twocell{\beta}}
\arrow{r}{Fr}
\arrow[swap]{d}{f'} &
FA'' 
\arrow{d}{f''} \\
GB 
\arrow[bend right = 50, swap]{rr}{G(s \circ q)}
\arrow[swap]{r}{Gq} &
GB'
\arrow[phantom]{d}[description]{\twocellIso{\phi^G}}
\arrow[swap]{r}{Gs} &
GB'' \\
\: &
\: &
\: 
\end{tikzcd}
\end{equation*}
We call axiom~(\ref{eq:glueing-2-cell-condition}) the 
\Def{cylinder condition} due to its shape when viewed as a (3-dimensional) 
pasting diagram (\cf~the \Def{cylinders} of ~\cite[\S~8]{Benabou1967}). 
From this perspective, the axiom requires that if one passes across the top of 
the cylinder and then down the front, the result is the same as passing first 
down the back of the cylinder and then the bottom~(\cf~the definition of 
\Def{transformation} between $T$-algebra morphisms in 2-dimensional universal 
algebra~\cite[\S~4.1]{Lack2010}):
\begin{equation*}
\begin{tikzcd}[row sep = 2.5em]
FA 
\arrow[phantom]{r}[description, font=\small]{\Downarrow}
\arrow[phantom]{dr}[yshift=-1.5mm]{\twocell{\alpha}}
\arrow[bend left = 20]{r}
\arrow[bend right = 20]{r} 
\arrow{d} &
FA'
\arrow{d} \\
GA 
\arrow[bend left = 20, phantom]{r}
\arrow[bend right = 20]{r} &
GB'
\end{tikzcd}
\:\: = \:\:
\begin{tikzcd}[row sep = 2.5em]
FA 
\arrow[phantom]{dr}[yshift=1.5mm]{\twocell{\alpha'}}
\arrow[bend left = 20]{r}
\arrow[bend right = 20, phantom]{r} 
\arrow{d} &
FA'
\arrow{d} \\
GA 
\arrow[phantom]{r}[description, font=\small]{\Downarrow}
\arrow[bend left = 20, dashed]{r}
\arrow[bend right = 20]{r} &
GB'
\end{tikzcd}
\end{equation*}

The following lemma, which mirrors the categorical statement, helps assure 
us the preceding definition is correct. For the proof one simply unwinds the 
two universal properties. 

\begin{prooflesslemma}
For any pseudofunctor $F : \baseCat \to \altCat$ and $C \in \altCat$, the 
following are equivalent:
\begin{enumerate}
\item $(R, u)$ is a biuniversal arrow from $F$ to $C$, 
\item $(FR \xra{u} C)$ is the terminal object in $(F \downarrow \constPseudofunctor_C)$, 
where $\constPseudofunctor_C$ denotes the constant pseudofunctor at $C$. \qedhere
\end{enumerate}
\end{prooflesslemma}
The glueing construction is an instance of the comma construction.

\begin{mydefn} \label{def:glueing-bicat}
The \Def{glueing bicategory} $\gl{\glueFun }$ of bicategories $\baseCat$ and 
$\altCat$ 
along a pseudofunctor $\glueFun  : \baseCat \to \altCat$ is the comma 
bicategory 
$(\id_{\altCat} \downarrow \glueFun )$.
\end{mydefn}

\newpage
As in 
Definition~\ref{def:1-categorical-glueing}, we order the tuples in a comma 
bicategory as they are read 
down the page. In the particular case of a glueing bicategory, therefore, the 
objects, 1-cells and 2-cells have the following form:
\begin{align*}
\text{objects} &: (C \in \altCat, c : C \to \glueFun B, B \in \baseCat) \\
\text{1-cells} &: (q : C \to C', \alpha : c' \circ q \To \glueFun (p) \circ c, 
p : B \to B') \\
\text{2-cells} &: (\tau : q \To q', \sigma : p \To p')
\end{align*} 
One now obtains projection 
\emph{pseudo}functors 
$\baseCat 
	\xla{\pi_{\mathrm{dom}}} 
\gl{\glueFun } 
	\xra{\pi_{\mathrm{cod}}} 
\altCat$.
Note also that there is a `weakest link' property at play: the bicategory 
$\gl{\glueFun}$ 
is a 2-category only if $\baseCat$, $\altCat$ and $\glueFun$ are all strict.

\begin{myremark}
The preceding definitions are \emph{pseudo}. One obtains a \Def{lax} comma 
bicategory (and hence lax glueing bicategory) by dropping the requirement that 
the 2-cells 
filling~(\ref{eq:glueing-bicategorical-morphism-square}) are invertible. 
\end{myremark}

\section{Cartesian closed structure on 
\texorpdfstring{$\gl{\glueFun}$}{the glueing bicategory}}

We now turn to a bicategorical version of 
Proposition~\ref{prop:glueing-cat-ccc}. The construction for products is 
relatively easy. 

\subsection{Finite products in \texorpdfstring{$\gl{\glueFun}$}{the glueing 
bicategory}}
\label{sec:glued-products-constructed}

Recall from Definition~\ref{def:fp-bicat} that a bicategory with finite 
products---an \Def{fp-bicategory}---is a bicategory $\baseCat$ equipped with a 
chosen object $\prodop_n(A_1, \dots, A_n)$ and a  
biuniversal arrow
$(\pi_1, \,\dots\, , \pi_n) : 
\Delta\big(\prodop_n(A_1, \,\dots\, , A_n)\big) \to (A_1, \,\dots\, , A_n)$ for 
every 
$A_1, \,\dots\, , A_n \in \baseCat \:\: (n \in \Nat)$. An 
\Def{fp-pseudofunctor} is then a pseudofunctor of the underlying bicategories 
that preserves these biuniversal arrows (Definition~\ref{def:fp-pseudofunctor}).

We claim the following:

\begin{prooflesspropn} \label{prop:glueing-bicat-has-products}
Let $\fpBicat{\baseCat}$ and $\fpBicat{\altCat}$ be fp-bicategories and 
$(\glueFun , \prodPres) : \baseCat \to \altCat$ an fp-pseudofunctor. Then 
$\gl{\glueFun }$ is an fp-bicategory with both projection 
pseudofunctors $\pi_{\mathrm{dom}}$ and $\pi_{\mathrm{cod}}$
strictly preserving products.
\end{prooflesspropn}

\newpage
We construct the data in stages and then verify the required 
equivalence on hom-categories. Recall that we denote 
the 2-cells witnessing the fact that $\glueFun$ preserves products by
\begin{align*}
\unTimes_{\ind{B}} : \Id_{(\prod_i \glueFun B_i)} \To \seq{\glueFun \pi_1, 
\dots, 
\glueFun \pi_n} \circ 
\prodPres_{\ind{B}} \\
\coTimes_{\ind{B}} : \prodPres_{\ind{B}} \circ \seq{\glueFun \pi_1, \,\dots\, , 
\glueFun \pi_n} \To 
\Id_{\glueFun (\prod_i B_i)}
\end{align*}
We begin with the product mapping. For a family of objects 
$(C_i, c_i, B_i)_{i=1, \,\dots\, , n}$ we define the $n$-ary product 
$\prod_{i=1}^n (C_i, c_i, B_i)$ to be the tuple 
$\big( \prod_{i=1}^n C_i, 
	\prodPres_{\ind{B}} \circ \prod_{i=1}^n c_i, 
\prod_{i=1}^n B_i \big)$.
We set the $k$-th projection $\glued{\pi}_k$ 
to be
$(\pi_k, \mu_k, \pi_k)$, where 
$\mu_k$ is defined
%
%
by commutativity of the following diagram:
\begin{equation} \label{eq:def-of-mu-k}
\begin{tikzcd}[column sep = 4em]
c_k \circ \pi_k
\arrow[swap]{d}{\epsilonTimesInd{-k}{}} 
\arrow{r}{\mu_k} &
\glueFun (\pi_k) \circ \left(\prodPres_{\ind{B}} \circ \prod_i c_i\right) 
\\
\pi_k \circ \prod_i c_i 
\arrow[swap]{d}{\iso} &
(\glueFun\pi_k \circ \prodPres_{\ind{B}}) \circ \prod_i c_i
\arrow[swap]{u}{\iso} \\
(\pi_k \circ \Id_{(\prod_i \glueFun B_i)}) \circ \prod_i c_i
\arrow[swap, bend right = 12]{dr}[yshift=0mm]
	{\pi_k \circ \unTimes_{\ind{B}} \circ \prod_i c_i}  
&
\left( (\pi_k \circ \seq{\glueFun \pi_1, \,\dots\, , \glueFun \pi_n}) \circ 
	\prodPres_{\ind{B}} \right) \circ \prod_i c_i
\arrow[swap]{u}
	{\epsilonTimesInd{k}{} \circ \prodPres_{\ind{B}} \circ \prod_i c_i} \\
\: &
\left(\pi_k \circ (\seq{\glueFun \pi_1, \,\dots\, , \glueFun \pi_n} \circ 
\prodPres_{\ind{B}})\right) \circ 
\prod_i c_i 
\arrow[swap]{u}{\iso}
\end{tikzcd}
\end{equation}

Next we define the $n$-ary tupling map. For an $n$-ary family of 1-cells 
$(g_i, \alpha_i, f_i) : (Y, y, X) \to (C_i, c_i, B_i) \:\: (i=1,\dots, n)$, 
we set the $n$-ary tupling to be 
$(\seq{g_1, \,\dots\, , g_n}, 
\glTup{\alpha_1, \,\dots\, , \alpha_n},
\seq{f_1, \,\dots\, , f_n})$, 
where $\glTup{\alpha_1, \,\dots\, , \alpha_n}$ is the composite
\begin{equation} \label{eq:def-of-glTup}
\begin{tikzcd}[column sep = 3.5em]
\left(\prodPres_{\ind{B}} \circ \prod_i c_i\right) \circ \seq{g_1, \,\dots\, , 
g_n} 
\arrow[swap]{d}{\iso}
\arrow{r}{\glTup{\alpha_1, \,\dots\, , \alpha_n}} &
\glueFun\seq{f_1, \,\dots\, , f_n} \circ y \\
\prodPres_{\ind{B}} \circ \left( \prod_i c_i \circ \seq{g_1, \,\dots\, , 
g_n}\right)
\arrow[swap]{d}{\prodPres_{\ind{B}} \circ \fuse}  &
\Id_{\glueFun (\prod B_i)} 
	\circ \left(\glueFun\seq{f_1, \,\dots\, , f_n} \circ y\right)  
\arrow[swap]{u}{\iso} \\
\prodPres_{\ind{B}} \circ \seqlr{c_1 \circ g_1, \,\dots\, , c_n \circ g_n} 
\arrow[swap]{d}{\prodPres_{\ind{B}} \circ \seq{\alpha_1, \,\dots\, , \alpha_n}} 
&
\left( \prodPres_{\ind{B}} 
	\circ \seq{\glueFun \pi_1, \,\dots\, , \glueFun \pi_n}\right) 
	\circ \left(\glueFun\seq{f_1, \,\dots\, , f_n} \circ y\right)
\arrow[swap]{u}
	{\coTimes_{\ind{B}} \circ \glueFun\seq{f_1, \,\dots\, , f_n} ) \circ y} \\
\prodPres_{\ind{B}} \circ 
	\seqlr{\glueFun f_1 \circ y, \,\dots\, , \glueFun f_n \circ y} 
\arrow[swap, bend right = 12]{dr}{\prodPres_{\ind{B}} \circ \postName^{-1}} &
\prodPres_{\ind{B}} 
	\circ \left( ( \seq{\glueFun \pi_1, \,\dots\, , \glueFun \pi_n} 
				\circ \glueFun\seq{f_1, \,\dots\, , f_n} ) \circ y \right)
\arrow[swap]{u}{\iso}\\
\: &
\prodPres_{\ind{B}} \circ 
	\left(\seqlr{\glueFun f_1, \,\dots\, , \glueFun f_n} \circ y\right)
\arrow[swap]{u}{\prodPres_{\ind{B}} \circ \unpack_{\ind{f}}^{-1} \circ y} 
\end{tikzcd}
\vspace{2mm}
\end{equation}

Finally, we are required to provide a universal arrow to act as the 
counit. For every family of 1-cells 
$(g_i, \alpha_i, f_i) : (Y, y, X) \to (C_i, c_i, B_i) \:\: (i=1, \,\dots\, , 
n)$ we 
require a glued 2-cell 
\[
\glued{\pi}_k \circ (\seq{g_1, \,\dots\, , g_n}, 
\glTup{\alpha_1, \,\dots\, , \alpha_n}, 
\seq{f_1, \,\dots\, , f_n}) \To (g_k, \alpha_k, f_k)
\] 
for which we take simply 
$(\epsilonTimesInd{k}{\ind{g}}, \epsilonTimesInd{k}{\ind{f}})$. The next lemma 
establishes that this is a 2-cell in $\gl{\glueFun }$.

\begin{mylemma}
For every family of 1-cells 
$(g_i, \alpha_i, f_i) : (Y, y, X) \to (C_i, c_i, B_i) \:\: 
(i=1, \,\dots\, , n)$, the cylinder condition holds for 
$(\epsilonTimesInd{k}{\ind{g}}, \epsilonTimesInd{k}{\ind{f}})$. That is, the 
following diagram commutes:
\begin{td}
c_k \circ (\pi_k \circ \seq{g_1, \,\dots\, , g_n})
\arrow[swap]{d}{\iso}
\arrow{r}{c_k \circ \epsilonTimesInd{k}{}} &
c_k \circ g_k
\arrow{r}{\alpha_k} &
\glueFun(f_k) \circ y \\

(c_k \circ \pi_k) \circ \seq{g_1, \,\dots\, , g_n} 
\arrow[swap]{d}{\mu_k \circ \seq{g_1, \,\dots\, , g_n}} &
\: &
\glueFun\left( \pi_k \circ \seq{f_1, \,\dots\, , f_n} \right) \circ y
\arrow[swap]{u}{\glueFun(\epsilonTimesInd{k}{}) \circ y} \\

\left( \glueFun(\pi_k) 
	\circ \left( \prodPres_{\ind{B}} \circ \prod_i c_i \right) \right) 
	\circ \seq{g_1, \,\dots\, , g_n} 
\arrow[swap]{d}{\iso} &
\: &
\left( \glueFun\pi_k \circ \glueFun\seq{f_1, \,\dots\, , f_n} \right) \circ y
\arrow{u}[swap]{\phi^{\glueFun}_{\pi_k; \seq{\ind{f}}} \circ y} \\

\glueFun\pi_k 
	\circ \left( \left( \prodPres_{\ind{B}} \circ \prod_i c_i \right) 
	\circ \seq{g_1, \,\dots\, , g_n} \right)
\arrow[swap]{rr}{\glueFun(\pi_k) \circ \glTup{\alpha_1, \,\dots\, , \alpha_n}} &
\: &
\glueFun\pi_k 
	\circ \left(\glueFun\seq{f_1, \,\dots\, , f_n} \circ y \right)
\arrow[swap]{u}{\iso}
\end{td}
\begin{proof}
Unfolding the definition of $\fuse$ and applying the functoriality of 
composition as far as possible, the claim reduces to commutative 
diagram below, in which the unlabelled cells are all instances of functoriality 
of composition or naturality. To improve readability we neglect the bracketing 
and 
corresponding 
associativity constraints; the 
coherence theorem for bicategories guarantees that one can translate to the 
`fully bicategorical' version as required.
\begin{td}[column sep = 2.5em, row sep = 2em]
\pi_k \circ \seq{\ind{c} \circ \ind{g}}
\arrow[swap]{d}{\iso}
\arrow{rr}{\epsilonTimesInd{k}{}}
\arrow{dr}[description]{\pi_k \circ \seq{\alpha_1, \,\dots\, , \alpha_n}} &
\: &
c_k \circ g_k
\arrow{dddddddd}{\alpha_k}  \\

\pi_k \circ \Id_{(\prod_i \glueFun B_i)} \circ 
\seq{\ind{c} \circ \ind{g}}
\arrow[swap]{d}{\pi_k \circ \unTimes_{\ind{B}} \circ \seq{\ind{c} \circ 
\ind{g}}}
\arrow[bend left = 20]{ddr}[description, xshift=1mm]{\pi_k \circ \Id \circ 
\seq{\alpha_1, \,\dots\, , \alpha_n}}
\arrow[phantom]{ddr}[description, yshift=-4mm]{\equals{triang. law}} &
\pi_k \circ \seq{\glueFun (\ind{f}) \circ y}
\arrow{dd}[description]{\iso}
\arrow[bend left = 60]{ddddd}[description]{\pi_k \circ \postName^{-1}} 
\arrow[phantom]{dddddddr}
	[description, yshift=-22mm, xshift=6mm]{\equals{$\postName$ def.}} &
\: \\

\pi_k \circ \seq{\glueFun \ind{\pi}} \circ \prodPres_{\ind{B}} \circ  
\seq{\ind{c} 
\circ \ind{g}}
\arrow[swap]{d}{\pi_k \circ \seq{\glueFun \ind{\pi}} \circ \prodPres_{\ind{B}} 
\circ 
\seq{\alpha_1, \,\dots\, , \alpha_n}} &
\: &
\: \\

\pi_k \circ \seq{\glueFun \ind{\pi}} \circ \prodPres_{\ind{B}} \circ 
\seq{\glueFun (\ind{f}) \circ y}
\arrow[swap]{d}{\pi_k \circ \postName^{-1}}
\arrow{r}[yshift=2mm]{\pi_k \circ \coTimes_{\ind{B}} \circ \seq{\glueFun 
(\ind{f}) 
\circ y}} 
&
\pi_k \circ \Id_{(\prod_i \glueFun B_i)} \circ \seq{\glueFun (\ind{f}) \circ 
y}
\arrow{d}[description]{\pi_k \circ \Id \circ \postName^{-1}}  &
\: \\

\pi_k \circ \seq{\glueFun \ind{\pi}} \circ \prodPres_{\ind{B}} \circ 
\seq{\glueFun \ind{f}} \circ y
\arrow[swap]{d}{\pi_k \circ \seq{\glueFun \ind{\pi}} \circ \prodPres_{\ind{B}} 
\circ 
\unpack_{\ind{f}}^{-1} \circ y }
\arrow{r}[yshift=2mm]{\pi_k \circ \coTimes_{\ind{B}} \circ \seq{\glueFun 
\ind{f}} 
\circ y} &
\pi_k \circ \Id_{(\prod_i \glueFun B_i)} \circ \seq{\glueFun \ind{f}} 
\arrow[bend left = 20]{ddl}[description]{\pi_k \circ \seq{\glueFun \ind{\pi}} 
\circ 
\unpack_{\ind{f}}^{-1} \circ y}
\arrow{dd}[description]{\iso{}}
\circ y  &
\: \\

\pi_k \circ \seq{\glueFun \ind{\pi}} \circ \prodPres_{\ind{B}} \circ 
\seq{\glueFun \ind{\pi}} \circ \glueFun \seq{\ind{f}} \circ y
\arrow[swap]{d}{\pi_k \circ \seq{\glueFun \ind{\pi}} \circ \coTimes_{\ind{B}} 
\circ 
\glueFun \seq{\ind{f}} \circ y } &
\:  &
\: \\

\pi_k \circ \seq{\glueFun \ind{\pi}} \circ \Id_{(\glueFun \prod_i B_i)} \circ 
\glueFun \seq{\ind{f}} 
\circ y 
\arrow[swap]{d}{\iso} &
\pi_k \circ \seq{\glueFun \ind{f}} \circ y
\arrow{ddr}[description]{\epsilonTimesInd{k}{} \circ y}
\arrow[bend left = 10]{dl}[description]
	{\pi_k \circ \unpack_{\ind{f}}^{-1} \circ y} 
\arrow[phantom]{dd}[description, yshift=-3mm]{\equals{$\unpack$ def.}}
 &
\: \\

\pi_k \circ \seq{\glueFun \ind{\pi}} \circ \glueFun \seq{\ind{f}} \circ y 
\arrow[swap]{d}{\epsilonTimesInd{k}{} \circ \glueFun \seq{\ind{f}} \circ y} &
\: &
\: \\

\glueFun (\pi_k) \circ \glueFun \seq{\ind{f}} \circ y
\arrow[swap]{r}{\phi^\glueFun _{(\pi_k, \seq{\ind{f}})} \circ y} &
\glueFun (\pi_k \circ \seq{\ind{f}}) \circ y
\arrow[swap]{r}{\glueFun (\epsilonTimesInd{k}{}) \circ y} &
\glueFun (f_k) \circ y
\end{td}
\end{proof}
\end{mylemma}
It remains to check the universal property. Taking arbitrary 1-cells 
\begin{align*}
(v, \gamma, u) : (Y, y, X) &\to \prodop_{i=1}^n (C_i, c_i, B_i) \\
(t_i, \tau_i, s_i) : (Y, y, X) &\to (C_i, c_i, B_i) \qquad\qquad\quad 
	(i=1, \,\dots\, , n)
\end{align*}
related by 2-cells 
\[(\beta_i, \alpha_i) : 
	\glued{\pi_i} \circ (v, \gamma, u) 
		\To 
	(t_i, \tau_i, s_i) \qquad\qquad\quad (i=1,\dots,n)
\]
we observe that
$\beta_i : \pi_i \circ v \To t_i$ and
$\alpha_i : \pi_i \circ u \To s_i$ 
for each $i$. We therefore claim that 
$\left(\transTimes{\beta_1, \,\dots\, , \beta_n}, 
	\transTimes{\alpha_1, \,\dots\, , \alpha_n}\right)$ 
is the unique 2-cell in 
$\gl{\glueFun}$ such that the following commutes 
for $i=1, \,\dots\, , n$:
\begin{td}
\glued{\pi_i} \circ (v, \gamma, u)
\arrow{rr}
	{\glued{\pi_i} \circ (\transTimes{\ind{\beta}},\transTimes{\ind{\alpha}})}
\arrow[swap]{dr}{(\beta_i, \alpha_i)} &
\: &
\glued{\pi_i} \circ (\seq{\ind{t}}, \glTup{\ind{\tau}}, \seq{\ind{s}}) 
\arrow{dl}{(\epsilonTimesInd{i}{\ind{t}}, \epsilonTimesInd{i}{\ind{s}})} \\
\: &
(t_i, \tau_i, s_i) &
\:
\end{td}
Of course, it suffices to show that
$\left(\transTimes{\ind{\beta}}, \transTimes{\ind{\alpha}}\right)$ is a 2-cell 
in 
$\gl{\glueFun}$: the rest of the claim 
follows from the (bi)universality of products in $\baseCat$ and 
$\altCat$. 

\begin{mylemma} \label{lem:products-cylinder-condition}
For any 1-cells
$(v, \gamma, u)$ and
$(t_i, \tau_i, s_i)$ and any 2-cells 
$(\beta_i, \alpha_i) : 
	\glued{\pi_i} \circ (v, \gamma, u) 
		\To 
	(t_i,\tau_i, s_i) 
	\:\: (i=1,\dots,n)$
as above, 
the pair 
$\left(\transTimes{\beta_1, \,\dots\, , \beta_n}, 
	\transTimes{\alpha_1, \,\dots\, , \alpha_n}\right)$ is a 2-cell in 
	$\gl{\glueFun}$. 
\begin{proof}
We need to check the cylinder condition, which in this case is the following:
\begin{td}[column sep = 11em]
\left(\prodPres_{\ind{B}} \circ  \prod_i c_i\right) \circ v
\arrow[swap]{d}{\gamma}
\arrow{r}
	{\prodPres_{\ind{B}} \circ \left(\prod_i c_i\right) \circ 	
	\transTimes{\beta_1, \dots, \beta_n}} &
\left(\prodPres_{\ind{B}} \circ \prod_i c_i\right) \circ \seq{t_1, \,\dots\, , 
t_n}
\arrow{d}{\glTup{\tau_1, \,\dots\, , \tau_n}} \\
\glueFun (u) \circ y 
\arrow[swap]{r}{\glueFun (\transTimes{\alpha_1, \,\dots\, , \alpha_n}) \circ y} 
&
\glueFun (\seq{s_1, \,\dots\, , s_n}) \circ y
\end{td}
For this, one begins by observing that the following commutes for every
$k=1, \,\dots\, , n$:
\begin{td}[row sep = 3em]
\pi_k \circ \left( \prod_i c_i \circ v \right) 
\arrow{r}{\iso}
\arrow[swap]{d}
	{\pi_k \circ \left(\prod_i c_i\right) 
		\circ \transTimes{\beta_1, \,\dots\, , \beta_n}} &
\left( \pi_k \circ \prod_i c_i \right) \circ v
\arrow{r}{\epsilonTimesInd{k}{} \circ v} 
\arrow{d}[description]
	{\pi_k \circ \prod_i c_i \circ \transTimes{\beta_1, \,\dots\, , \beta_n}} &
(c_k \circ \pi_k) \circ v 
\arrow[]{d}{\iso} \\

\pi_k \circ \left(\prod_i c_i \circ \seq{\ind{t}}\right)
\arrow[swap]{dd}{\pi_k \circ \fuse}
\arrow[swap]{r}{\iso} &
\left(\pi_k \circ \prod_i c_i\right) \circ \seq{\ind{t}}
\arrow{d}[description]{\epsilonTimesInd{k}{} \circ \seq{\ind{t}}} &
c_k \circ (\pi_k \circ v)
\arrow{d}[description]
	{c_k \circ \pi_k \circ \transTimes{\beta_1, \,\dots\, , \beta_n}}
\arrow[bend left = 90]{dd}{c_k \circ \beta_k} \\

\: &
(c_k \circ \pi_k) \circ \seq{\ind{t}}
\arrow[swap]{r}{\iso} 
\arrow[phantom]{dl}[description]{\equals{def. of $\fuse$}} &
\c_k \circ \left(\pi_k \circ \seq{\ind{t}}\right) 
\arrow{d}[description]{c_k \circ \epsilonTimesInd{k}{}} \\

\pi_k \circ \seq{\ind{c} \circ \ind{t}}
\arrow[swap]{rr}{\epsilonTimesInd{k}{}}
\arrow[swap]{d}{\pi_k \circ \seq{\ind{\tau}}} &
\: &
c_k \circ t_k 
\arrow{d}{\tau_k} \\

\pi_k \circ \seqlr{\glueFun(\ind{s}) \circ y}
\arrow[phantom]{drr}[description]{\equals{def. of $\postName$}}
\arrow[swap]{d}{\pi_k \circ \postName^{-1}}
\arrow[swap]{rr}{\epsilonTimesInd{k}{}} &
\: &
\glueFun(s_k) \circ y
\arrow[equals]{d} \\

\pi_k \circ \left(\seqlr{\glueFun\ind{s}} \circ y\right)
\arrow[swap]{r}{\iso} &
(\pi_k \circ \seqlr{\glueFun\ind{s}}) \circ y
\arrow[swap]{r}{\epsilonTimesInd{k}{} \circ y} &
\glueFun(s_k) \circ y 
\end{td}
and that the following commutes:
\begin{td}
c_k \circ (\pi_k \circ v )
\arrow[swap]{dd}{c_k \circ \beta_k}
\arrow{rr}{\iso} &
\: &
(c_k \circ \pi_k) \circ v
\arrow{d}{\mu_k \circ v} \\

\: &
\: &
\left( \glueFun(\pi_k) \circ (\prodPres_{\ind{B}} \circ \prod_i c_i ) \right)
	\circ v
\arrow{d}{\iso} \\

\pi_k \circ t_k 
\arrow[swap]{dd}{\tau_k} &
\:  &
\glueFun(\pi_k) 
	\circ \left( \left( \prodPres_{\ind{B}} \circ \prod_i c_i \right)
	\circ v \right)
\arrow{d}{\glueFun(\pi_k) \circ \gamma} \\

\: &
\glueFun(\pi_k \circ u) \circ y 
\arrow[phantom]{uuu}[description]{\equals{cylinder condition}}
\arrow[phantom]{ddl}[description]{\equals{def.}}
\arrow[phantom]{ddr}[xshift=-7mm, description]{\equals{nat.}}
\arrow[bend left]{ddl}{\glueFun (\pi_k \circ \transTimes{\ind\alpha}) \circ y}
\arrow{dl}[description]{\glueFun (\alpha_k) \circ y} &
\glueFun (\pi_k) \circ (\glueFun(u) \circ y)
\arrow{d}{\iso} \\

\glueFun(s_k) \circ y
\arrow[swap]{d}{\glueFun(\epsilonTimesInd{-k}{}) \circ y} &
\: &
\left( \glueFun\pi_k \circ \glueFun u \right) \circ y
\arrow{ul}[description]{\phi_{\pi_k, u}^\glueFun  \circ y}
\arrow{d}{\glueFun(\pi_k) \circ \glueFun (\transTimes{\ind{\alpha}}) \circ y} \\

\glueFun (\pi_k \circ \seq{\ind{s}}) \circ y
\arrow[swap]{rr}{(\phi^\glueFun _{\pi_k, \seq{\ind{s}}})^{-1} \circ y} &
\: &
\left(\glueFun(\pi_k) \circ \glueFun\seq{\ind{s}}\right) \circ y
\end{td}
Putting these two together and applying the definition of $\unpack$, one 
obtains the following commuting diagram:
\begin{td}[column sep = 1.5em]
\pi_k \circ \left(\prod_i c_i \circ v\right)
\arrow[swap]{d}{\pi_k \circ \prod_i c_i \circ \transTimes{\ind{\beta}}}
\arrow{rr}{\iso} &
\: &
(\pi_k \circ \Id_{(\prod_i \glueFun B_i)}) 
	\circ \left(\prod_i c_i \circ v\right)
\arrow{d}{\pi_k \circ \unTimes_{\ind{B}} \circ \prod_i c_i \circ v} \\

\pi_k \circ \left(\prod_i c_i \circ \seq{\ind{t}}\right)
\arrow[swap]{d}{\pi_k \circ \fuse} &
\: &
\left(\pi_k 
	\circ \left(\seq{\glueFun \ind{\pi}} \circ \prodPres_{\ind{B}}\right)\right)
	\circ \left(\prod_i c_i \circ v\right)
\arrow{d}{\iso}  \\

\pi_k \circ \seq{\ind{c} \circ \ind{t}}
\arrow[swap]{d}{\pi_k \circ \seq{\ind{\tau}}} &
\: &
\left(\pi_k \circ \seq{\glueFun \ind{\pi}}\right) 
	\circ \left( \left(\prodPres_{\ind{B}} \circ \prod_i c_i\right) \circ v	
	\right)
\arrow{d}
	{\epsilonTimesInd{k}{} \circ \prodPres_{\ind{B}} \circ \prod_i c_i \circ v} 
	\\

\pi_k \circ \seq{\glueFun(\ind{s}) \circ y} 
\arrow[swap]{d}{\pi_k \circ \postName^{-1}} &
\: &
\glueFun(\pi_k) 
	\circ \left( \left(\prodPres_{\ind{B}} \circ \prod_i c_i\right) 
	\circ v	\right) 
\arrow{d}{\glueFun(\pi_k) \circ \gamma} \\

\pi_k \circ \left(\seq{\glueFun \ind{s}} \circ y\right) 
\arrow[swap]{d}{\pi_k \circ \unpack_{\ind{s}}^{-1} \circ y} &
\: &
\glueFun(\pi_k) \circ \left(\glueFun(u) \circ y\right)
\arrow{d}
	{\glueFun(\pi_k) \circ \glueFun(\transTimes{\ind{\alpha}}) \circ y} \\

\pi_k \circ 
	\left( \left(\seq{\glueFun \ind{\pi}} \circ \glueFun\seq{\ind{s}} \right) 
	\circ y\right)
\arrow[swap]{r}{\iso} &
\left(\pi_k \circ \seq{\glueFun\ind{\pi}}\right) \circ 
	\left(\glueFun\seq{\ind{s}} \circ y\right)
\arrow[swap]{r}[yshift=0mm]{\epsilonTimesInd{k}{} \circ \glueFun 
\seq{\ind{s}} \circ y} &
\glueFun(\pi_k) \circ \left(\glueFun\seq{\ind{s}} \circ y\right)
\end{td}
With this lemma in hand, the rest of the proof is a diagram 
chase applying naturality and the definition of $\postName$. 
\end{proof}
\end{mylemma}

Lemma~\ref{lem:products-cylinder-condition} completes the proof that 
$\gl{\glueFun}$ does indeed have finite 
products, and hence the proof of 
Proposition~\ref{prop:glueing-bicat-has-products}. For the construction of 
exponentials we will require morphisms of the form $f \times A$. We briefly 
check that such morphisms appear in $\gl{\glueFun }$ in the way one would 
expect, namely as pasting diagrams of the form
\begin{td}[column sep = 6em, row sep = 5em]
C \times Y
\arrow[phantom]{drr}[description, yshift=-1.5mm, near end]
{\twocellIso{\alpha \times y}}
\arrow[bend left = 10]{drr}[description]{c'g \times y}
\arrow[bend right = 20]{drr}[description]{\glueFun (f)c \times y}
\arrow{rr}{g \times Y} 
\arrow[swap]{d}{c \times y}
\arrow[bend right = 60, swap]{dd}{\prodPres_{B, X} \circ (c \times y)} &
\: &
C' \times Y
\arrow[phantom]{dl}[description, near start]{\twocellIso{\phiTimes}}
\arrow{d}{c' \times y}
\arrow[bend left = 60]{dd}{\prodPres_{B', X} \circ (c' \times y)} \\

\glueFun B \times \glueFun X
\arrow[phantom]{rr}
	[description, xshift=-4mm, yshift=-2mm]
	{\twocellIso{\phiTimes}}
\arrow[phantom]{drr}[description, yshift=-5mm]{\twocellIso{\nat}}
\arrow[bend right = 20, swap]{rr}{\glueFun f \times \glueFun X}
\arrow[swap]{d}{\prodPres_{B,X}} &
\: &
\glueFun B' \times \glueFun X 
\arrow{d}{\prodPres_{B',X}} \\

\glueFun (B \times X) 
\arrow[swap]{rr}{\glueFun (f \times X)} &
\: &
\glueFun (B' \times X)
\end{td}
In particular, when the bicategories
$\baseCat$ and $\altCat$ are 2-categories with strict products and
$\glueFun  : \baseCat \to \altCat$ is a strict fp-pseudofunctor, this 2-cell is 
simply $\alpha \times y$. 

\begin{mylemma} \label{lem:f-times-Y-in-glued-bicat}
For every 1-cell $\glued{g} := (g, \alpha, f) : (C, c, B) \to (C', c', B')$ and 
object 
$\glued{Y} := (Y, y, X)$ in $\gl{\glueFun }$, the 1-cell 
$\glued{g} \times \glued{Y} : 
	(C, c, B) \times (Y, y, X) 
		\to 
	(C', c', B') \times (Y, y, X)$ 
is equal to 
$(g \times Y, \glued{\alpha}_{\glued{Y}}, f \times Y)$, 
where $\glued{\alpha}_{\glued{Y}}$ is the composite
\begin{equation} \label{eq:glued-alpha-times-Y}
\begin{tikzcd}[column sep = 7em]
\left(\prodPres_{B', X} \circ (c' \times y)\right) \circ (g \times Y)
\arrow{r}{\glued{\alpha}_{\glued{Y}}}
\arrow[swap]{d}{\iso} &
\glueFun(f \times X) 
	\circ \left(\prodPres_{B, X} \circ (c \times y)\right) \\
\prodPres_{B', X} \circ \left( (c' \times y) \circ (g \times Y) \right)
\arrow[swap]{d}{\prodPres_{B', X} \circ \phiTimes_{c', g; y, \Id}} &
\left(\glueFun(f \times X) \circ \prodPres_{B, X}\right) \circ (c \times y)
\arrow[swap]{u}{\iso} \\
\prodPres_{B', X} \circ \left( (c' \circ g) \times (y \circ \Id_Y) \right) 
\arrow[swap]{d}{\iso} &
\left( \prodPres_{B', X} \circ (\glueFun f \times \glueFun \Id_X) \right) 
	\circ (c \times y)
\arrow[swap]{u}{\nat_{f, \Id_X} \circ (c \times y)} \\
\prodPres_{B', X} \circ 
	\left( (c' \circ g) \times (\Id_{\glueFun X} \circ y) \right)
\arrow[swap]{d}{\prodPres_{B', X} \circ 
	(\alpha \times \left(\Id_{\glueFun X} \circ y)\right)} &
\prodPres_{B', X} \circ 
\left((\glueFun f \times \glueFun \Id_X) \circ (c \times y)\right)
\arrow[swap]{u}{\iso} \\
\prodPres_{B', X} \circ 
	\left((\glueFun f \circ c) \times (\Id_{\glueFun X} \circ y)\right)
\arrow[swap]{r}[yshift=-2mm]
	{\prodPres_{B', X} \circ ((\glueFun f \circ c) \times 
		(\psi^\glueFun _X \circ y))} &
\prodPres_{B', X} \circ 
	\left((\glueFun f \circ c) \times (\glueFun\Id_X \circ y)\right) 
\arrow[swap]{u}
	{\prodPres_{B', X} \circ 
	\phiTimes^{-1}_{\glueFun f, c; \glueFun\Id, y}} 
\end{tikzcd}
\end{equation}
\begin{proof}
The proof amounts to unfolding the definition and checking that 
it does indeed equal the composite given in the claim. Let 
$\tau_1$ and $\tau_2$ respectively denote the 2-cells defined by the 
pasting diagrams on the left and right below:
\begin{center}
\makebox[\textwidth]{\parbox{1.1\textwidth}{
\centering
\begin{tikzcd}[ampersand replacement = \&]
C \times Y 
\arrow[swap]{d}{c \times y}
\arrow[phantom]{ddr}[description, xshift=3mm]{\twocell{\mu_1}}
\arrow{r}{\pi_1}
\arrow[bend left]{rr}{g \circ \pi_1}
\arrow[bend right = 70, swap]{dd}{\prodPres_{B,X} \circ (c \times y)} \&
C 
\arrow[phantom]{ddr}[description]{\twocell{\alpha}}
\arrow{r}{g} 
\arrow{dd}[description]{c} \&
C' 
\arrow{dd}{c'} \\
\glueFun B \times \glueFun X 
\arrow[swap]{d}{\prodPres_{B,X}} \&
\: \&
\: \\
\glueFun (B \times X)
\arrow[bend right = 60, swap]{rr}{\glueFun (f \circ \pi_1)}
\arrow[swap]{r}{\glueFun \pi_1} \&
\glueFun B
\arrow[phantom]{d}[description]{\twocellIso{\phi^\glueFun _{f, \pi_1}}} 
\arrow[swap]{r}{\glueFun f} \&
\glueFun B' \\
\: \&
\: \&
\: 
\end{tikzcd}
\qquad\quad
\begin{tikzcd}[ampersand replacement = \&]
C \times Y 
\arrow[swap]{d}{c \times y}
\arrow[phantom]{ddr}[description, xshift=3mm]{\twocell{\mu_2}}
\arrow{r}{\pi_2}
\arrow[bend left]{rr}{\Id_Y\circ \pi_2}
\arrow[bend right = 70, swap]{dd}{\prodPres_{B,X} \circ (c \times y)} \&
Y 
\arrow{ddr}[description]{y}
\arrow{r}{\Id_Y} 
\arrow{dd}[description]{y} \&
Y
\arrow[phantom]{dl}[description, near start, font=\small, xshift=1mm]{\iso}
\arrow{dd}{y} \\
\glueFun B \times \glueFun X 
\arrow[swap]{d}{\prodPres_{B,X}} \&
\: \&
\: \\
\glueFun (B \times X)
\arrow[bend right = 75, swap]{rr}{\glueFun (\Id_X \circ \pi_2)}
\arrow[swap]{r}{\glueFun \pi_2} \&
\glueFun X
\arrow[phantom]{d}[description,  xshift=-2mm, 
yshift=-2mm]{\twocellIso{\phi^\glueFun _{\Id, \pi_2}}} 
\arrow[phantom]{uu}[near start, description, font=\small, xshift=5mm, 
yshift=2mm]{\iso}
\arrow[phantom]{r}[description]{\twocellIso{\psi^\glueFun _X}}
\arrow[bend left]{r}{\Id_{\glueFun X}}
\arrow[bend right = 50, swap]{r}{\glueFun \Id_X} \&
\glueFun X \\
\: \&
\: \&
\: 
\end{tikzcd}
}}
\end{center}
By definition, the 1-cell $\glued{g} \times \glued{Y}$ has a 
witnessing 2-cell given by the following composite, in which we write
$(\ast)$ for
$\prodPres_{B', X} \circ 
	\seqlr{\left(\glueFun (f \circ \pi_1) 
		\circ \prodPres_{B', X}\right) \circ (c 	\times 	y), 
		\left(\glueFun(\Id_X \circ \pi_2) 
		\circ \prodPres_{B', X}\right) \circ (c \times y)}$:
\begin{center}
\makebox[\textwidth]{\parbox{1.1\textwidth}{
\centering
\begin{tikzcd}[column sep = -4em, ampersand replacement = \&]
\left(\prodPres_{B', X} \circ (c' \times y)\right) \circ 
	\seqlr{g \circ \pi_1, \Id_Y \circ \pi_2}
\arrow[swap]{d}{\iso}
\arrow{r}{\glTup{\tau_1, \tau_2}} \&
\glueFun(f \times B) \circ \left(\prodPres_{B, X} \circ (c \times y)\right) \\
\prodPres_{B',X} 
	\circ \left( (c' \times y) \circ \seqlr{g \circ \pi_1, \Id_Y \circ \pi_2} 
	\right) 
\arrow[swap]{d}{\prodPres_{B', X} \circ \fuse} \&
\: \\
\prodPres_{B', X} \circ 
	\seqlr{c' \circ (g \circ \pi_1), y \circ (\Id_Y \circ \pi_2)} 
\arrow[swap]{d}{\prodPres_{B', X} \circ \seq{\tau_1, \tau_2}} \&
\: \\
(\ast)
\arrow[swap]{d}{\prodPres_{B', X} \circ \postName^{-1}} \&
\Id_{\glueFun(B' \times X)} \circ \glueFun (f \times X) \circ (c \times y) 
\arrow[swap]{uuu}{\iso}	 \\
\prodPres_{B', X} \circ 
	\left(\seqlr{\glueFun (f \circ \pi_1), 
						\glueFun(\Id_X \circ \pi_2)} 
	\circ \left(\prodPres_{B',X} \circ (c \times y) \right) \right)
\arrow[swap]{d}{\prodPres_{B', X} \circ 
		\unpack^{-1}_{f \circ \pi_1, \Id \circ 	\pi_2} 
		\circ (c \times y)}	 \&
\:	\\
\prodPres_{B', X} \circ \left(\left(\seq{\glueFun \pi_1, \glueFun \pi_2} 
	\circ \glueFun(f \times X)\right) 
	\circ \left( \prodPres_{B' \times X} \circ (c \times y)\right) \right)
\arrow[swap, bend right = 8]{dr}{\iso} \&
\: \\
\: \&
\left(\prodPres_{B', X} \circ \seq{\glueFun \pi_1, \glueFun \pi_2}\right)
	\circ \left(\glueFun(f \times X) 
	\circ \left( \prodPres_{B' \times X} \circ (c \times y)\right) \right)
\arrow[swap]{uuu}
	{\coTimes_{B', X} \circ \glueFun(f \times X) 
		\circ (c \times y)}
\end{tikzcd}
}}
\end{center}
Applying naturality and the lemma relating $\unpack$ with $\unTimes$ 
(Lemma~\ref{lem:unpack-prod-pres}), a long diagram 
chase transforms this to the composite in the claim.
\end{proof}
\end{mylemma}

\subsection{Exponentials in \texorpdfstring{$\gl{\glueFun}$}{the glueing 
bicategory}}
\label{sec:exponentials-in-glueing-constructed}

As in the 1-categorical case, the 
definition of currying in $\gl{\glueFun}$ employs pullbacks. We therefore 
take a brief diversion to spell out their universal property.

\paragraph{Pullbacks in a bicategory.}
The notion of pullback we employ is sometimes referred to as a 
\Def{bipullback}~(\eg~\cite{Lack2010}) to 
distinguish it from pullbacks defined as a pseudolimit. Since the only limits 
we work with in this thesis are bilimits, we omit the prefix. 

\begin{mydefn} \label{def:pullback}
Let $\cospanCat$ (for `cospan') denote the category $\cospan$ and $\baseCat$ be 
any 
bicategory. A \Def{pullback} of the cospan
$(X_1 \xra{f_1} X_0 \xla{f_2} X_2)$ in $\baseCat$ is a bilimit for the strict
pseudofunctor $\cospanCat \to \baseCat$ determined by this cospan.
\end{mydefn}

This characterisation of pullbacks, while precise, 
must be unfolded to obtain a universal property one can use for calculations. 
The next lemma establishes such a property. The proof is not especially hard, 
and the result appears to be known---although not explicitly proven---in the 
literature, so we leave it for an 
appendix
(Appendix~\ref{chap:extra-proofs}).

\begin{prooflesslemma} \label{lem:pullback-ump}
For any bicategory $\baseCat$ and cospan $(X_1 \xra{f_1} X_0 \xla{f_2} X_2)$ 
in $\baseCat$, the pullback of  $(X_1 \xra{f_1} X_0 \xla{f_2} X_2)$ is 
determined, up to 
equivalence, by the following universal property: there exists a chosen
object $P \in \baseCat$, span 
$(X_1 \xla{\gamma_1} P \xra{\gamma_2} X_2)$ and 
invertible 2-cell $\overline{\gamma}$ filling the diagram on the left below
\begin{equation}  \label{eq:pullback-biuniversal-arrow}
\begin{tikzcd}[column sep = small, row sep = small]
\: &
P 
\arrow[phantom]{dd}[description]{\twocellIso{\overline{\gamma}}}
\arrow[swap]{dl}{\gamma_1}
\arrow{dr}{\gamma_2} &
\: \\
X_1 
\arrow[swap]{dr}{f_1} &
\: &
X_2 
\arrow{dl}{f_2} \\
\: &
X_0 &
\:
\end{tikzcd}
\qquad\qquad 
\begin{tikzcd}[column sep = small, row sep = small]
\: &
Q 
\arrow[phantom]{dd}[description]{\twocellIso{\overline{\mu}}}
\arrow[swap]{dl}{\mu_1}
\arrow{dr}{\mu_2} &
\: \\
X_1 
\arrow[swap]{dr}{f_1} &
\: &
X_2 
\arrow{dl}{f_2} \\
\: &
X_0 &
\:
\end{tikzcd}
\end{equation}
such that for any other such square as on the right above
there exists an invertible \Def{fill-in} 
$(u, \modif_1, \modif_2)$~(\cf~\cite{Vitale2010}), namely a 1-cell 
$u : Q \to P$ and invertible 2-cells 
$\modif_i : \gamma_i \circ u \To \mu_i \:\: (i=1,2)$ such that
\begin{equation} \label{eq:fill-in-actual}
\begin{tikzcd}[column sep = 3.5em]
\left(f_2 \circ \gamma_2\right) \circ u
\arrow{r}{\iso} 
\arrow[swap]{d}{\overline{\gamma} \circ u}&
f_2 \circ (\gamma_2 \circ u) 
\arrow{r}{f_2 \circ \modif_2} &
f_2 \circ \mu_2 
\arrow{d}{\overline{\mu}} \\
\left(f_1 \circ \gamma_1\right) \circ u
\arrow[swap]{r}{\iso} &
f_1 \circ \left(\gamma_1 \circ u\right)
\arrow[swap]{r}{f_1 \circ \modif_1} &
f_1 \circ \mu_1
\end{tikzcd}
\end{equation}
This fill-in is universal in the following sense. For any other fill-in 
\[(v : Q \to P, 
\altModif_1 : {\gamma_1 \circ v \To \mu_1}, 
\altModif_2 : \gamma_2 \circ v \To \mu_2)\] 
there exists a 2-cell 
$\trans{\altModif} : v \To u$, unique such that 
\begin{equation} \label{eq:fill-in-mediating-property}
\begin{tikzcd}
\gamma_i \circ v 
\arrow[swap]{dr}{\altModif_i}
\arrow{rr}{\gamma_i \circ \trans\altModif} &
\: &
\gamma_i \circ u \arrow{dl}{\modif_i} \\
\: &
\mu_i &
\:
\end{tikzcd}
\end{equation}
for $i=1,2$. Finally, it is required that for any $w : Q \to P$ the 
2-cell $\trans{\id}$ obtained by applying the universal property to 
$(w, \id_{\gamma_1 \circ w}, \id_{\gamma_2 \circ w})$ is invertible.
\end{prooflesslemma}

\begin{myremark} \quad \label{rem:pullback-remarks}
The universal property of pullbacks can be 
stated in a slightly different way, which is more useful for some calculations. 
The pullback of a cospan $(X_1 \xra{f_1} X_0 \xla{f_2} X_2)$ is determined by a 
biuniversal arrow $(\gamma, \overline{\gamma}) : \Delta P \To F$, for $F$ the 
pseudofunctor determined 
by the cospan, $P$ the pullback, and $(\gamma, \overline{\gamma})$ an 
iso-commuting 
square as in~(\ref{eq:pullback-biuniversal-arrow}). It follows that the functor
$(\gamma, \overline{\gamma}) \circ \Delta(-) : 
	\baseCat(Z, P) 
	\to 
	\Hom(\cospanCat, \baseCat)(\Delta Z, F)$
is fully-faithful and essentially surjective for every $Z \in \baseCat$. Being 
essentially surjective is 
exactly the existence of a fill-in for every iso-commuting square, as in the 
preceding lemma. Being full and faithful entails that, for every pair of
1-cells $t, u : Z \to P$ equipped with 2-cells
$\Gamma_i : \gamma_i \circ t \To \gamma_i \circ u \:\: (i=1,2)$ 
satisfying the fill-in law~(\ref{eq:fill-in-actual}), there exists a
unique 2-cell $\trans\Gamma : t \To u$ such that
$\gamma_i \circ \trans\Gamma = \Gamma_i$ for $i=1,2$. \qedhere
\end{myremark}

The following is an example of where it is convenient to use the 
universal property of 
Remark~\ref{rem:pullback-remarks}. The lemma 
guarantees that 
one may define objects in a glueing bicategory (up to equivalence) by pullback.

\begin{mylemma} \label{lem:pullback-equivalent-glued-objects}
For any pseudofunctor $\glueFun : \baseCat \to \altCat$ and any pullbacks 
\begin{center}
\begin{tikzcd}
P
\arrow[phantom]{dr}[very near start, yshift=4mm, xshift=-1mm]{\sidecorner}
\arrow[phantom]{dr}[description]{\twocellIso{\pi}}
\arrow[swap]{d}{p}
\arrow{r}{q} &
B 
\arrow{d}{b} \\
\glueFun A  
\arrow[swap]{r}{a} &
C
\end{tikzcd}
\qquad\qquad
\begin{tikzcd}
X
\arrow[phantom]{dr}[very near start, yshift=4mm, xshift=-1mm]{\sidecorner}
\arrow[phantom]{dr}[description]{\twocellIso{\chi}}
\arrow[swap]{d}{x}
\arrow{r}{y} &
B 
\arrow{d}{b} \\
\glueFun A  
\arrow[swap]{r}{a} &
C
\end{tikzcd}
\end{center}
in $\altCat$, the objects $(P \xra{p} \glueFun A)$ and $(X \xra{x} \glueFun A)$ 
are equivalent in $\gl{\glueFun}$.
\begin{proof}
It is immediate from the uniqueness of bilimits that there exists a canonical 
equivalence
$P \simeq X$. The 
only question is whether this equivalence lifts to a 1-cell
in $\gl{\glueFun}$. If one constructs the equivalence using the 
universal property of 
Remark~\ref{rem:pullback-remarks}, this follows 
immediately.
\end{proof}
\end{mylemma}

Preliminaries complete, we can now give the data for defining exponentials 
in the glueing bicategory. Precisely, we extend 
Proposition~\ref{prop:glueing-bicat-has-products} to the following. Recall that 
a cartesian closed bicategory---a \Def{cc-bicategory}---is an fp-bicategory 
equipped with a right biadjoint to $(-) \times A$ for every object $A$ 
(Definition~\ref{def:cc-bicat}).

\begin{prooflessthm} \label{thm:glueing-bicat-cartesian-closed}
Let $\ccBicat{\baseCat}$ and $\ccBicat{\altCat}$ be cc-bicategories and suppose 
that $\altCat$ has all pullbacks. Then for any fp-pseudofunctor
$(\glueFun , \prodPres) : \fpBicat{\baseCat} \to \fpBicat{\altCat}$
the glueing bicategory 
$\gl{\glueFun}$ is cartesian closed with forgetful 
pseudofunctor 
$\pi_{\mathrm{dom}} : \gl{\glueFun} \to \baseCat$ strictly preserving products 
and exponentials.
\end{prooflessthm}

Much of the complication in the 
definitions that follow arises from 
the 
invertible 2-cells moving 1-cells in and out of products; where the 
product structure is strict, the exponentials in $\gl{\glueFun}$ are given 
similarly to the 1-categorical case. The reader 
happy to employ Power's coherence result for fp-bicategories
(Proposition~\ref{prop:power-coherence}) may therefore greatly simplify the 
definitions just 
given and the calculations to come. Because we wish to prove an independent 
coherence result, we do not take this approach. 

We begin by defining the mapping $\expobj{(-)}{(=)}$ and the evaluation 1-cell 
$\glued{\eval}$. 

\paragraph*{Defining $\expobj{(-)}{(=)}$ and $\glued{\eval}$.}

For 
$\glued{C} := (C, c, B)$ and $\glued{C'} := (C', c', B')$ in $\gl{\glueFun }$ 
we set 
the exponential $\expobj{\glued{C}}{\glued{C'}}$ to be the left-hand vertical 
leg of the 
following pullback diagram, in which 
$\evBar_{B,B'}$ is the exponential transpose of 
$\glueFun(\eval_{B,B'}) \circ \prodPres$
\big(\cf~the definition in the 
1-categorical 
case~(\ref{eq:exponential-in-glued-cat})\big):
\begin{equation}  \label{eq:exponential-in-glued-bicat}
\begin{tikzcd}[column sep = 8em]
\glexp{C}{C'} 
\arrow[phantom]{drr}[description]{\twocell{\omega_{c,c'}}}
\arrow[phantom]{dr}[at start, description, yshift=2mm, xshift=-1mm]{\sidecorner}
\arrow[swap]{d}{p_{c,c'}}
\arrow{rr}{q_{c,c'}} &
\: &
(\exp{C}{C'}) 
\arrow{d}{\lambda(c' \circ \eval_{C,C'})} \\
\glueFun (\exp{B}{B'})
\arrow[
swap,
rounded corners,
to path=
{ -- ([yshift=0ex]\tikztostart.south)
|- ([yshift=-4.5ex]\tikztostart.south)
-| ([yshift=-4.5ex]\tikztotarget.south)
-- (\tikztotarget.south)}, 
]{rr}
\arrow[phantom, swap]{rr}[font=\scriptsize, yshift=-8.1ex]
	{\lambda(\eval_{\glueFun B,\glueFun B'} 
		\circ ((\scriptsizeexpobj{\glueFun B}{\glueFun B'}) \times c)) 
		\circ \evBar_{B,B'}}
\arrow[swap]{r}{\evBar_{B,B'}} &
(\exp{\glueFun B}{\glueFun B'}) 
\arrow[swap]{r}[yshift=-2mm]
	{\lambda(\eval_{\glueFun B,\glueFun B'} 
		\circ 
	((\scriptsizescriptsizeexpobj{\glueFun B}{\glueFun B'}) \times c))} &
(\exp{C}{\glueFun B'})
\end{tikzcd}
\end{equation}
We use $\lambda(c' \circ \eval_{C,C'})$ and 
$\lambda(\eval_{\glueFun B, \glueFun B'} 
	\circ ((\expobj{\glueFun B}{\glueFun B'}) \times c))$ 
instead of 
$(\expobj{\glueFun B}{c})$ and $(\expobj{C}{c'})$ as a simplifying measure: 
doing so 
avoids the need to apply the isomorphisms 
$(\expobj{\glueFun B}{c}) \iso \lambda(c' \circ \eval_{C,C'})$ and
$(\expobj{C}{c'}) 
\iso \lambda(\eval_{\glueFun B,\glueFun B'} 
	\circ ((\expobj{\glueFun B}{\glueFun B'}) \times c))$ removing the 
redundant identities in the left-hand side (recall the comment after 
Notation~\ref{not:exponentials}). 

\begin{mynotation}
For reasons of space---particularly for fitting pasting 
diagrams onto a single page---we 
will sometimes write
$\compact{c} := 
	\eval_{\glueFun B, \glueFun B'} \circ ((\expobj{\glueFun B}{\glueFun B'}) 
	\times c)$ 
where $c : C \to \glueFun B$ in $\altCat$ (see, for 
example,~(\ref{eq:def-of-eval-map-witnessing-2-cell})). 
\end{mynotation}
For the evaluation 1-cell $\glued{\eval}$ we take the 1-cell with components
\begin{gather*}
(\glexp{C}{C'}) \times C \xra{q_{c,c'} \times C} 
(\expobj{C}{C'}) \times C \xra{\eval_{C,C'}} C' \\
(\expobj{B}{B'}) \times B \xra{\eval_{B,B'}} B'
\end{gather*}
The witnessing 2-cell $\evalMod_{\glued{C}, \glued{C'}}$ is given by the 
following 
pasting diagram.
\begin{equation} \label{eq:def-of-eval-map-witnessing-2-cell}
\makebox[\textwidth]{\parbox{1.1\textwidth}{
\centering
\begin{tikzcd}[column sep = 1em, row sep = 2.5em, ampersand replacement = \&]
(\glexp{C}{C'}) \times C \arrow{rrr}{q_{c,c'} \times C} 
\arrow{dr}[description]{p_{c,c'} \times C} 
\arrow[swap]{dd}[description]{p_{c,c'} \times c} 
\arrow[bend left = 15]{rrrr}{\eval_{C,C'} \circ (q_{c,c'} \times C)}
\arrow[bend right = 70, swap]{ddd}[description]
	{\prodPres_{(\scriptsizescriptsizeexpobj{B}{B'}, B)} 
		\circ (p_{c,c'} \times c)}
\arrow[phantom]{drrr}[description, xshift=13ex]
	{\twocellIso{\omega_{c,c'} \times C}} \&
\: \&
\: \&
(\exp{C}{C'}) \times C 
\arrow{r}{} 
\arrow{d}[description]{\lambda(c' \circ \eval_{C,C'}) \times C} \&
C' \arrow{ddd}{c'} \\
\arrow[phantom]{r}[description]{\twocellIso{\phiTimes}} \&
\glueFun (\exp{B}{B'}) \times C 
\arrow{r}[yshift=0.5ex]{\evBar_{B,B'} \times C} 
\arrow{dl}[description]{\glueFun(\scriptsizeexpobj{B}{B'}) \times c} 
\arrow[phantom]{dr}[description, xshift=-2ex]{\twocellIso{}} \&
(\exp{\glueFun B}{\glueFun B'}) \times C 
\arrow{r}[yshift=0.5ex]{\lambda \compact{c} \times C} 
\arrow{d}[description]{(\scriptsizeexpobj{\glueFun B}{\glueFun B'}) \times c} 
\arrow[phantom]{ddrr}[yshift=2ex]{\twocellIso{\epsilonExp}} \&
(\exp{C}{\glueFun B'}) \times C 
\arrow{ddr}[description]{\eval_{C, \glueFun B'}} 
\arrow[phantom]{r}[description]{\twocellIso{\epsilonExp}} \&
\: \\
\glueFun (\exp{B}{B'}) \times \glueFun B 
\arrow[swap]{rr}{\evBar_{B,B'} \times \glueFun B} 
\arrow[swap]{d}[description]{\prodPres_{(\scriptsizeexpobj{B}{B'}, B)}} 
\arrow[phantom]{drrr}[description]{\twocellIso{\epsilonExp}} \&
\: \&
(\exp{\glueFun B}{\glueFun B'}) \times \glueFun B 
\arrow{drr}[description]{\eval_{\glueFun B, \glueFun B'}} \&
\:  \&
\: \\
\glueFun \left((\exp{B}{B'}) \times B\right) 
\arrow[swap]{rrrr}{\glueFun \eval_{B,B'}} \&
\: \&
\: \&
\: \&
\glueFun B'
\end{tikzcd}
}}
\vspace{2mm}
\end{equation}
Here we omit the canonical 2-cells for the product 
structure: thus, the shape labelled $\omega_{c,c'} \times C$ is actually the 
composite 
\begin{td}
\big(\lambda(c' \circ \eval_{C,C'}) \times C\big) \circ 
\big( q_{c,c'} \times C \big) 
\arrow{r}
\arrow[swap]{d}{\phiTimes_{\lambda(c' \circ \eval), q; \Id}} &
(\lambda \compact{c} \times C) 
	\circ \left((\evBar_{B,B'} \times C) \circ (p_{c,c'} \times C)\right) \\

\big(\lambda(c' \circ \eval_{C,C'}) \circ q_{c,c'}\big) \times 
	(\Id_C \circ \Id_C)
\arrow[swap]{d}{\iso} &
\: \\

\big(\lambda(c' \circ \eval_{C,C'}) \circ q_{c,c'}\big) \times C 
\arrow[swap]{r}{\omega_{c,c'} \times C} &
(\lambda\compact{c} \circ \evBar_{B,B'} \circ p_{c,c'}) \times C
\arrow[swap]{uu}{\iso}
\end{td}
in which the unlabelled isomorphism employs two applications of 
$\phiTimes^{-1}$, together with the evident structural isomorphisms. 

\begin{mynotation}
For the rest of this chapter we will 
adopt the convention just employed, and write simply $\iso$ for
instances of either $\phiTimes$ or 
its inverse, composed with structural isomorphisms. Power's coherence result 
guarantees that this is valid as an explanatory shorthand: of course, the 
masochistic reader could work explicitly with all the instances of $\phiTimes$ 
and prove exactly the same set of diagrams commute. Thus, while Power's result 
is useful 
\emph{for reasons of exposition and presentation}, the proofs we present do 
not rely on it.  
\end{mynotation}

With this convention, $\evalMod_{\glued{C}, \glued{C'}}$ is the following 
composite:
\begin{equation} \label{eq:def-of-glued-E}
\makebox[\textwidth]{
\begin{tikzcd}[ampersand replacement = \&]
c' \circ \left(\eval_{C,C'} \circ (q_{c,c'} \times C)\right)
\arrow[swap]{d}{\iso}
\arrow{r}{\evalMod_{\glued{C}, \glued{C'}}} \& 
\glueFun (\eval_{B,B'}) \circ 
	\big(\prodPres_{(\scriptsizeexpobj{B}{B'}, B)} 
		\circ (p_{c,c'} \times c)\big)
\\
\left( c' \circ \eval_{C,C'} \right) \circ (q_{c,c'} \times C)
\arrow[swap]{d}{\epsilonExp_{(c' \circ \eval)}^{-1} \circ (q_{c,c'} \times C)} 
\&
\:
	\\
\left(\eval_{C, C'} \circ \left( \lambda(c' \circ \eval_{C,C'}) \times C 
		\right)\right) \circ (q_{c,c'} \times C)
\arrow[swap]{d}{\iso} \&
\big(\glueFun(\eval_{B,B'}) \circ 
	\prodPres_{(\scriptsizeexpobj{B}{B'}, B)}\big) \circ (p_{c,c'} \times c)
\arrow[swap]{uu}{\iso}  \\
\eval_{C, C'} \circ \big( \lambda(c' \circ \eval_{C,C'}) \circ q_{c,c'}\big) 
\times C 
\arrow[swap]{d}{\eval \circ (\omega_{c,c'} \times C)} \&
\left(\eval_{\glueFun B,\glueFun B'} 
	\circ (\evBar_{B,B'} \times \glueFun B)\right) 
	\circ (p_{c,c'} \times c) 
\arrow[swap]{u}{\epsilonExp_{(\glueFun \eval \circ \prodPres)} 
	\circ (p_{c,c'} \times c)} \\ 
\eval_{C, C'} \circ 
\left( \left(\lambda\compact{c} \circ \evBar_{B,B'}\right) 
	\circ p_{c,c'}\right) \times C 
\arrow[swap]{d}{\iso} \&
\:  \\
\left(\eval_{C,C'} \circ \left( \lambda\compact{c} \times C\right) \right)
	\circ \left(\evBar_{B,B'}p_{c,c'} \times C \right) 
\arrow[swap]{r}[yshift=0mm]
	{\epsilonExp_{\compact{c}} \circ (\evBar_{B,B'}p_{c,c'} \times C)} \&
\compact{c}
	\circ \left(\evBar_{B,B'}p_{c,c'} \times C\right)  
\arrow[swap]{uu}{\iso}
\end{tikzcd}
}
\end{equation}

\paragraph*{The mapping $\glued{\lambda}$.} 
Next we need to provide a mapping $\glued{\lambda}$ assigning a 
1-cell 
of type $\glued{R} \to (\expobj{\glued{C}}{\glued{C'}})$ to every 
1-cell $\glued{R} \times \glued{C} \to \glued{C'}$.
Let
$\glued{R} := (R, r, Q)$, $\glued{C} := (C, c, B)$ and 
$\glued{C'} := (C', c', B')$.
As our starting point, suppose given a 1-cell 
$(t, \alpha, s) : \glued{R} \times \glued{C} \to \glued{C'}$, as on the left 
below: 
\begin{center}
\begin{tikzcd}
R \times C 
\arrow[phantom]{ddr}[description, xshift=3mm]{\twocell{\alpha}}
\arrow[bend right = 70]{dd}[description, xshift=-4mm]
	{\prodPres_{Q, B} \circ (r \times c)}
\arrow[swap]{d}{r \times c}
\arrow{r}{t} &
C'
\arrow{dd}{c'} \\
\glueFun Q \times \glueFun B
\arrow[swap]{d}{\prodPres_{Q,B}} &
\: \\
\glueFun (Q \times B) 
\arrow[swap]{r}{\glueFun s} &
\glueFun B'
\end{tikzcd}
\quad\:
\begin{tikzcd}
R 
\arrow[bend right = 65, swap]{dd}[description]{\glueFun(\lambda s) \circ r}
\arrow[phantom]{ddrr}[description, xshift=-1mm]{\twocell{\L_\alpha}}
\arrow[swap]{d}{r}
\arrow{rr}{\lambda t} &
\: &
\expobj{C}{C'}
\arrow{dd}{\lambda(c' \circ \eval_{C,C'})} \\
\glueFun Q 
\arrow[swap]{d}{\glueFun \lambda s} & 
\: &
\: \\
\glueFun (\expobj{B}{B'})
\arrow[bend right = 20, swap]{rr}{\lambda{\compact{c}} \circ \evBar_{B,B'}}
\arrow[swap]{r}{\evBar_{B,B'}} &
\expobj{\glueFun B}{\glueFun B'} 
\arrow[swap]{r}[yshift=0mm]{\lambda\compact{c}} &
\expobj{C}{\glueFun B'}
\end{tikzcd}
\end{center}
We construct a 2-cell $\L_\alpha$ as on 
the right above and apply the universal property of the 
pullback~(\ref{eq:exponential-in-glued-bicat}). To this end, let us define two 
invertible composites, which we denote by $\T_\alpha$ and $\U_\alpha$. For 
$\T_\alpha$ we take
\begin{td}[column sep=-1em]
\eval_{C, \glueFun B'} \circ 
\big( \lambda(c' \circ \eval_{C,C'}) \circ \lambda t \big) \times C 
\arrow{rr}{\T_\alpha}
\arrow[swap]{d}{\iso} &
\: & 
c' \circ t \\

\left(\eval_{C, \glueFun B'} \circ 
\left(\lambda(c' \circ \eval_{C,C'}) \times C \right)\right) \circ 
(\lambda t \times C) 
\arrow[swap, bend right = 8]{dr}[yshift=0mm]
{\epsilonExp_{(c' \circ \eval)} \circ (\lambda t \times C)} &
\: &
c' \circ \left(\eval_{C,C'} \circ (\lambda t \times C)\right)
\arrow[swap]{u}{c' \circ \epsilonExp_t} \\

\: &
\left(c' \circ \eval_{C,C'}\right) \circ (\lambda t \times C) 
\arrow[swap, bend right = 10]{ur}{\iso} &
\:
\end{td}
and for $\U_\alpha$ we take 
\begin{center}
\makebox[\textwidth]{
\begin{small}
\begin{tikzcd}[column sep = -0.5em, ampersand replacement = \&]
\eval_{C, \glueFun B} \circ 
\left(\left(\lambda\compact{c} \circ \evBar_{B,B'}\right) 
	\circ \left(\glueFun(\lambda s) \circ r \right) \right) 
	\times C
\arrow[swap]{d}{\iso}
\arrow{r}{\U_\alpha} \&
\glueFun s \circ  \left(\prodPres_{Q,B} \circ (r \times c)\right) \\
\eval_{C, \glueFun B} 
	\circ \left(\lambda\compact{c} \times C\right) 
	\circ \left(
		\evBar_{B,B'} \circ \left(\glueFun(\lambda s) \circ r\right) 
		\right) 
	\times C 
\arrow[swap]{d}
{\epsilonExp_{\compact{c}} \circ 
	(\evBar_{B,B'} \circ \glueFun (\lambda s) \circ r ) \times C} 
\&
\: \\
\compact{c} 
	\circ \left(
		\evBar_{B,B'} \circ \left(\glueFun(\lambda s) \circ r\right) 
	\right) \times C 
\arrow[swap]{d}{\iso} \&
\glueFun{\left(\eval_{B,B'} \circ (\lambda s \times B)\right)} 
	\circ \left(\prodPres_{Q,B} \circ (r \times c)\right)
\arrow[swap]{uu}{\glueFun \epsilonExp_s \circ \prodPres \circ (r \times c)} \\
\left(\eval_{\glueFun B, \glueFun B'} 
	\circ (\evBar_{B,B'} \times \glueFun B) \right)
	\circ \left((\glueFun (\lambda s) \times \glueFun B) 
	\circ (r \times c) \right)
\arrow[swap]{d}
	{\epsilonExp_{(\glueFun \eval \circ \prodPres)} 
		\circ (\glueFun(\lambda s) \times \glueFun B) 
		\circ (r \times c)} \&
\left(\glueFun(\eval_{B,B'}) \circ \glueFun(\lambda s \times B)\right) 
	\circ \left(\prodPres_{Q,B} \circ (r \times c)\right)
\arrow[swap]{u}
	{\phi^\glueFun_{\eval, \lambda s \times B} 
		\circ \prodPres 
		\circ (r \times c)} \\
\left(\glueFun(\eval_{B,B'}) 
	\circ \prodPres_{(\scriptsizescriptsizeexpobj{B}{B'}, B)}\right) 
	\circ \left((\glueFun (\lambda s) \times \glueFun \Id_B) 
			\circ (r \times c)\right)
\arrow[swap, bend right = 8]{dr}{\iso} \&
\glueFun(\eval_{B,B'}) 
	\circ \left(\left(\glueFun(\lambda s \times B) 
	\circ \prodPres_{Q,B}\right) 
	\circ (r \times c)\right)
\arrow[swap]{u}{\iso} \\
\: \&
\glueFun(\eval_{B,B'}) 
	\circ \big(\big(\prodPres_{(\scriptsizescriptsizeexpobj{B}{B'}, B)}
	\circ (\glueFun (\lambda s) \times \glueFun \Id_B)\big) 
			\circ (r \times c)\big)
\arrow[swap]{u}
	{\glueFun (\eval_{B,B'}) \circ \nat \circ (r \times c)} 
\end{tikzcd}
\end{small}
}
\end{center}
We may therefore define a 2-cell $\K_\alpha$ as the composite
\begin{td}[column sep = 5em]
\eval_{C, \glueFun B'} 
	\circ \left( \lambda(c' \circ \eval_{C, C'}) \circ \lambda t \right) 
		\times C
\arrow[swap]{d}{\T_\alpha}
\arrow{r}{\K_\alpha} &
\eval_{C, \glueFun B} \circ 
\left(\left(\lambda\compact{c} \circ \evBar_{B,B'}\right) 
	\circ \left(\glueFun(\lambda s) \circ r \right) \right) 
	\times C \\
c' \circ t
\arrow[swap]{r}{\alpha} &
\glueFun s \circ \left(\prodPres_{Q,B} \circ (r \times c) \right)
\arrow[swap]{u}{\U_\alpha^{-1}}
\end{td}
and, finally, $\L_\alpha$ as
\begin{equation*} 
\begin{tikzcd}[column sep = 0.8em]
\lambda(c' \circ \eval_{C,C'}) \circ \lambda t 
\arrow{rr}{\L_\alpha} 
\arrow[swap, bend right = 8]{dr}[yshift=1mm, xshift=-1mm]{\transExp{\K_\alpha}} 
&
\: &
\left(\lambda\compact{c} \circ \evBar_{B,B'}\right) 
	\circ \left(\glueFun(\lambda s) \circ r\right)
\\
\: &
\lambda\big( 
	\eval_{C, \glueFun B} 
		\circ \left(\left(\lambda\compact{c} \circ \evBar_{B,B'}\right) 
		\circ \left(\glueFun(\lambda s) \circ r \right) \right) \times C \big)
\arrow[swap, bend right = 8]{ur}[yshift=1mm, xshift=1mm]{\etaExp{}^{-1}} &
\:
\end{tikzcd}
\end{equation*}
Since we work in the pseudo setting, $\U_\alpha$,  
$\T_\alpha$, $\K_\alpha$---and hence $\L_\alpha$---are all invertible. 

Now, $\L_\alpha$ fills the following diagram:
\begin{equation} \label{eq:L-alpha-diagram}
\begin{tikzcd}[column sep = 7em]
R 
\arrow[swap]{d}{\glueFun(\lambda s) \circ r}
\arrow[phantom]{dr}[description, yshift=-1mm]{\twocellIso{\L_\alpha}}
\arrow{r}{\lambda t} &
(\expobj{C}{C'})
\arrow{d}{\lambda(c' \circ \eval_{C,C'})} \\
\glueFun(\expobj{B}{B'})
\arrow[swap]{r}{\lambda\compact{c} \circ \evBar_{B,B'}} &
(\expobj{C}{\glueFun B'})
\end{tikzcd}
\end{equation}
Hence, by the universal property of the 
pullback~(\ref{eq:exponential-in-glued-bicat}), one obtains a 1-cell 
$\glLam(t)$ and a pair of 
invertible 2-cells $\Gamma_{c,c'}$ and $\Delta_{c,c'}$ 
filling the diagram
\begin{equation} \label{eq:L-alpha-ump-diagram}
\begin{tikzcd}[column sep = 4em, row sep = 2.3em]
R
\arrow[bend right, swap]{ddr}{\glueFun(\lambda s) \circ r}
\arrow[bend left = 10]{drrr}{\lambda t}
\arrow[phantom]{drrr}[description, yshift=-1mm]{\twocellRight{\Delta_{c,c'}}}
\arrow[phantom]{ddr}[description, yshift=-1mm]{\twocell{\Gamma_{c,c'}}}
\arrow[dashed]{dr}[description]{\glLam(t)} &
\: &
\: &
\: \\
\: &
\glexp{C}{C'} 
\arrow[phantom]{drr}[description]{\twocell{\omega_{c,c'}}}
\arrow[phantom]{dr}[at start, description, yshift=2mm, xshift=1mm]{\sidecorner}
\arrow[swap]{d}[description]{p_{c,c'}}
\arrow{rr}[description]{q_{c,c'}} &
\: &
(\exp{C}{C'}) 
\arrow{d}{\lambda(c' \circ \eval_{C,C'})} \\
\: &
\glueFun (\exp{B}{B'})
\arrow[swap]{rr}{\lambda\compact{c} \circ \evBar_{B,B'}} &
\: &
(\exp{C}{\glueFun B'})
\end{tikzcd}
\end{equation}
such that the pasting diagrams~(\ref{eq:L-alpha-diagram}) 
and~(\ref{eq:L-alpha-ump-diagram}) are 
equal,~\ie~the following commutes:
\begin{equation} \label{eq:lam-coherence-condition}
\begin{tikzcd}[column sep = -1em, row sep = 2.5em]
\: &
\lambda(c' \circ \eval_{C, C'}) 
	\circ \left(q_{c,c'} \circ \glLam(t)\right)
\arrow{dr}{\lambda(c' \circ \eval_{C,C'}) \circ \Delta_{c,c'}} &
\: \\
\left(\lambda(c' \circ \eval_{C,C'}) \circ q_{c,c'}\right) \circ \glLam(t)
\arrow{ur}{\iso}
\arrow[swap]{d}{\omega_{c,c'} \circ \glLam(t)} &
\: &
\lambda(c' \circ \eval_{C,C'}) \circ \lambda t
\arrow{d}{\L_\alpha} \\
\left(\left(\lambda\compact{c} \circ \evBar_{B,B'}\right) 
	\circ p_{c,c'}\right) \circ \glLam(t)
\arrow[swap]{dr}{\iso} &
\: &
\left(\lambda\compact{c} \circ \evBar_{B,B'}\right) 
	\circ \left(\glueFun(\lambda s) \circ r \right) \\
\: &
\left(\lambda\compact{c} \circ \evBar_{B,B'}\right)
	\circ \left( p_{c,c'} \circ \glLam(t) \right)
\arrow[swap]{ur}{\lambda \compact{c} \circ \evBar_{B,B'} \circ \Gamma_{c,c'}} &
\:
\end{tikzcd}
\end{equation}
Moreover, $\Gamma_{c,c'}$ and $\Delta_{c,c'}$ are universal in the sense of 
Lemma~\ref{lem:pullback-ump}. We define 
$\glued{\lambda}(t, \alpha, s) := 
	\big( \glLam(t) , \Gamma_{c,c'}, \lambda s \big)$.

\newpage
\paragraph*{The counit $\glued{\epsilonExp}$.}
Finally we come to the counit. 
Let us first calculate 
$\glued{\eval} \circ \big( \glued{\lambda}(t, \alpha, s) \times (C,c,B)\big)$ 
for a 1-cell 
$\glued{t} := (t, \alpha, s) : (R, r, Q) \times (C, c, B) \to (C', c', B')$. 
Using Lemma~\ref{lem:f-times-Y-in-glued-bicat}, one unwinds this 1-cell to the 
following 
pasting diagram, in which we omit the canonical 
isomorphisms for the product structure as well as the structural isomorphisms:
\begin{td}[column sep = 6em, row sep = 4em]
R \times C 
\arrow[
swap,
rounded corners,
to path=
{ -- ([yshift=0ex]\tikztostart.north)
|- ([yshift=7ex]\tikztostart.east)
-| ([yshift=7ex]\tikztotarget.south)
-- (\tikztotarget.north)}
]{rr}
\arrow[phantom]{rr}[yshift=9ex, font = \scriptsize]
	{\left( \eval_{C,C'} \circ (q_{c,c'} \times C) \right) \circ 
		(\glLam(t) \times C)}
\arrow[bend right = 50]{dd}[swap]{\prodPres_{Q,B} \circ (r \times c)}
\arrow[phantom]{dr}[description, yshift=5mm]{\twocell{\Gamma_{c,c'} \times c}}
\arrow[swap]{d}[description]{r \times c}
\arrow{r}{\glLam(t) \times C} &
(\glexp{C}{C'}) \times C 
\arrow[phantom]{ddr}[description, xshift=2mm]
	{\twocellIso{\evalMod_{\glued{C}, \glued{C'}}}}
\arrow{d}[description]{q_{c,c'} \times c}
\arrow{r}{\eval_{C,C'} \circ (q_{c,c'} \times C)} &
C' 
\arrow{dd}{c'} \\

\glueFun Q \times \glueFun B 
\arrow[phantom]{dr}[description, yshift=-5mm]{\twocellIso{\nat}}
\arrow[phantom]{r}[description]{\twocellIso{\glueFun (\lambda s) \times 
\psi^\glueFun _B}}
\arrow[bend right, swap]{r}{\glueFun (\lambda s) \times \glueFun \Id_B}
\arrow[bend left]{r}{\glueFun (\lambda s) \times \glueFun B}
\arrow[swap]{d}[description]{\prodPres_{Q,B}} &
\glueFun (\expobj{B}{B'}) \times \glueFun B 
\arrow{d}[description]{\prodPres_{(\expobj{B}{B}',B)}} &
\: \\

\glueFun (Q \times B)
\arrow[phantom]{rr}[yshift=-9ex, font=\scriptsize]
	{\glueFun(\eval_{B,B'} \circ (\lambda s \times B))}
\arrow[
swap,
rounded corners,
to path=
{ -- ([yshift=0ex]\tikztostart.south)
|- ([yshift=-5ex]\tikztostart.south)
-| ([yshift=-5ex]\tikztotarget.north)
-- (\tikztotarget.south)}
]{rr}

\arrow[phantom]{rr}[description, yshift=-6mm]{\twocellIso{\phi^\glueFun }}
\arrow[swap]{r}{\glueFun (\lambda s \times B)} &
\glueFun \big( (\expobj{B}{B'}) \times B\big) 
\arrow[swap]{r}{\glueFun \eval_{B,B'}} &
\glueFun B' 
\end{td}
For the counit $\glued{\epsilonExp}_{\glued{t}}$ we therefore take the 2-cell 
with first 
component $\gluedCo{t}$ defined by
\begin{equation} \label{eq:def-of-counit}
\begin{tikzcd}
\left(\eval_{C,C'} \circ (q_{c,c'} \times C)\right) \circ (\glLam(t) \times C) 
\arrow{r}{\gluedCo{t}}
\arrow[swap]{d}{\iso}
&
t \\
\eval_{C,C'} \circ (q_{c,c'} \circ \glLam(t)) \times C  
\arrow[swap]{r}[yshift=-2mm]{\eval_{C,C'} \circ (\Delta_{c,c'} \times C)} &
\eval_{C,C'} \circ (\lambda t \times C) 
\arrow[swap]{u}{\epsilonExp_t}
\end{tikzcd}
\end{equation}
and second component simply
\[
\eval_{B,B'} \circ (\lambda s \times B) 
\XRA{\epsilonExp_s} 
s
\]
We need to check that this to be a legitimate 2-cell in 
$\gl{\glueFun}$,~\ie~that the cylinder condition holds.

\newpage
\begin{mylemma}
For any objects
$\glued{R} := (R, r, Q)$,
$\glued{C} := (C, c, B)$ and
$\glued{C'} := (C', c', B')$ and
1-cell 
$\glued{t} := (t, \alpha, s) : \glued{R} \times \glued{C} \to \glued{C'}$ 
in $\gl{\glueFun }$, the pasting diagram
\begin{td}[column sep = 7em, row sep = 4.3em]
R \times C 
\arrow[
swap,
rounded corners,
to path=
{ -- ([yshift=0ex]\tikztostart.north)
|- ([yshift=7ex]\tikztostart.east)
-| ([yshift=7ex]\tikztotarget.south)
-- (\tikztotarget.north)}
]{rr}
\arrow[phantom]{rr}[yshift=9ex, font = \scriptsize]
	{\left( \eval_{C,C'} \circ (q_{c,c'} \times C) \right) \circ 
		(\glLam(t) \times C)}
\arrow[bend right = 50]{dd}[swap]{\prodPres_{Q,B} \circ (r \times c)}
\arrow[phantom]{dr}[description, yshift=5mm]{\twocell{\Gamma_{c,c'} \times c}}
\arrow[swap]{d}[description]{r \times c}
\arrow{r}{\glLam(t) \times C} &
(\glexp{C}{C'}) \times C 
\arrow[phantom]{ddr}[description]{\twocellIso{\evalMod_{\glued{C}, \glued{C'}}}}
\arrow{d}[description]{q_{c,c'} \times c}
\arrow{r}{\eval_{C,C'} \circ (q_{c,c'} \times C)} &
C' 
\arrow{dd}{c'} \\

\glueFun Q \times \glueFun B 
\arrow[phantom]{dr}[description, yshift=-6mm]{\twocellIso{\nat}}
\arrow[phantom]{r}[description]{\twocellIso{\glueFun (\lambda s) \times 
\psi^\glueFun _B}}
\arrow[bend right, swap]{r}{\glueFun (\lambda s) \times \glueFun \Id_B}
\arrow[bend left]{r}{\glueFun (\lambda s) \times \glueFun B}
\arrow[swap]{d}[description]{\prodPres_{Q,B}} &
\glueFun (\expobj{B}{B'}) \times \glueFun B 
\arrow{d}[description]{\prodPres_{(\expobj{B}{B}',B)}} &
\: \\

\glueFun (Q \times B)
\arrow[phantom]{r}[yshift=-9ex, font=\scriptsize]
	{\glueFun(\eval_{B,B'} \circ (\lambda s \times B))}
\arrow[phantom]{r}[yshift=-15ex, font=\scriptsize]
	{\glueFun s}
\arrow[
swap,
rounded corners,
to path=
{ -- ([yshift=0ex]\tikztostart.south)
|- ([yshift=-5ex]\tikztostart.south)
-| ([yshift=-5ex]\tikztotarget.north)
-- (\tikztotarget.south)}
]{rr}
\arrow[
swap,
rounded corners,
to path=
{ -- ([yshift=0ex]\tikztostart.south)
|- ([yshift=-12ex]\tikztostart.south)
-| ([yshift=-12ex]\tikztotarget.north)
-- (\tikztotarget.south)}, 
]{rr}{\glueFun s}

\arrow[phantom]{rr}[description, yshift=-6mm]{\twocellIso{\phi^\glueFun }}
\arrow[phantom]{rr}[description, yshift=-17mm]
	{\twocellIso{\glueFun \epsilonExp_s}}
\arrow[phantom]{rr}[description, yshift=-6mm]{\twocellIso{\phi^\glueFun }}
\arrow[swap]{r}{\glueFun (\lambda s \times B)} &
\glueFun \big( (\expobj{B}{B'}) \times B\big) 
\arrow[swap]{r}{\glueFun \eval_{B,B'}} &
\glueFun B' 
\end{td}
is equal to 
\begin{td}[column sep = 3em, row sep = 2em]
\:
\arrow[phantom]{rr}[yshift=7ex, font=\scriptsize]
	{(q_{c,c'} \times C) \circ (\glLam(t) \times C)} &
(\glexp{C}{C'}) \times C
\arrow[phantom]{l}[description, xshift=0mm, yshift=0mm]{\iso}
\arrow[phantom]{d}
	[description, near start, yshift=-1mm, xshift=1mm, font=\scriptsize]
	{\twocellDown{\Delta_{c,c'} \times C}}
\arrow{r}{q_{c,c'} \times C} &
(\exp{C}{C'}) \times C 
\arrow[phantom]{r}[yshift=7mm, xshift=-3mm, font=\small]{\iso}
\arrow[phantom]{d}[description]{\twocellIso{\epsilonExp_t}}
\arrow{dr}[description]{\eval_{C,C'}} &
\: \\

R \times C
\arrow[
swap,
rounded corners,
to path=
{ -- ([yshift=0ex]\tikztostart.north)
|- ([yshift=14ex]\tikztostart.east)
-| ([yshift=4ex]\tikztotarget.south)
-- (\tikztotarget.north)}, 
]{urr}
\arrow[
swap,
rounded corners,
to path=
{ -- ([yshift=0ex]\tikztostart.north)
|- ([yshift=20ex]\tikztostart.east)
-| ([yshift=20ex]\tikztotarget.south)
-- (\tikztotarget.north)}, 
]{rrr}
\arrow[phantom]{rrr}[yshift=22ex, font = \scriptsize]
	{\left( \eval_{C,C'} \circ (q_{c,c'} \times C)\right)
		\circ (\glLam(t) \times C)}
\arrow[bend right = 70, swap]{dd}{\prodPres_{Q,B} \circ (r \times c)}
\arrow[phantom]{ddrrr}[description]{\twocell{\alpha}}
\arrow[swap]{d}{r \times c}
\arrow[bend right = 10, swap]{urr}[description]{\lambda t \times C}
\arrow{rrr}[description]{t}
\arrow{ur}[description]{\glLam(t) \times C} &
\: &
\: &
C' 
\arrow{dd}{c'} \\

\glueFun Q \times \glueFun B 
\arrow[swap]{d}{\prodPres_{Q,B}} &
\: &
\: &
\: \\

\glueFun (Q \times B) 
\arrow[swap]{rrr}{\glueFun s} &
\: &
\: &
\glueFun B'
\end{td}
Hence $\glued{\epsilonExp}_{\glued{t}} := (\gluedCo{t}, \epsilon_s)$ is a 
2-cell in $\gl{\glueFun }$.
\begin{proof}
Unfolding the first diagram, one sees that it is equal to the composite 
\begin{td}[column sep = 1em]
c' 
	\circ \left( \left(\eval_{C, C'} \circ (q_{c,c'} \times C)\right) 
	\circ (\glLam(t) \times C)\right)
\arrow[swap]{d}{(\ast)}
\arrow{r} &
\glueFun s \circ \left(\prodPres_{Q,B} \circ (r \times c)\right) \\

\eval_{\glueFun B, \glueFun B'} 
	\circ \left(\evBar_{B,B'} 
	\circ \left(p_{c,c'} 
	\circ \glLam(t)\right)\right) 
	\times c 
\arrow[swap]{d}
	{\eval_{\glueFun B, \glueFun B'} 
		\circ (\evBar_{B,B'} \circ \Gamma_{c,c'}) \times C} &
\eval_{C, \glueFun B} 
	\circ \left( 
		\left(\lambda\compact{c} \circ	\evBar_{B,B'}
		\right) 
	\circ \left(\glueFun (\lambda s) \circ r\right) 
	\right) 
	\times 	C 
\arrow[swap]{u}{\U_\alpha} \\

\eval_{\glueFun B, \glueFun B'}
	 \circ \left( \evBar_{B,B'} 
				\circ \left(\glueFun(\lambda s) 
				\circ r \right)\right) 
	\times c 
\arrow[swap, bend right = 8]{dr}{\iso} &
\left(\eval_{C, \glueFun B} 
	\circ \left(\lambda\compact{c} \times C\right)\right) 
	\circ \left(\evBar_{B,B'} \circ \left(\glueFun(\lambda s) \circ r 
			\right)\right) 
	\times 	C 
\arrow[swap]{u}{\iso}
 \\

\: &
\left(\eval_{B,B'} \circ \compact{c} \right)
	\circ \left(\evBar_{B,B'} 
		\circ \left(\glueFun(\lambda s) \circ r \right)\right) \times C  
\arrow[swap]{u}
	{\epsilonExp_{\compact{c}}^{-1} 
		\circ (\evBar_{B,B'} \circ \glueFun(\lambda s) \circ r) \times C}
\end{td}
where the arrow labelled $(\ast)$ arises by composing the following with 
structural isomorphisms and $\phiTimes$:
\begin{td}[column sep = 3em]
c' 
	\circ \left(\eval_{C, C'} \circ (q_{c,c'} \times C)\right) 
\arrow[swap]{d}{\iso} 
\arrow{r} &
\eval_{\glueFun B, \glueFun B'} 
	\circ \left((\evBar_{B,B'} \circ p_{c,c'}) \times c\right) \\

c' 
	\circ \big( \eval_{C,C'} \circ (q_{c,c'} \times C)\big)
\arrow[swap]{d}{\evalMod_{C,C'}} &
\left( 
	\eval_{\glueFun B, \glueFun B'} 
	\circ (\evBar_{B,B'} \times \glueFun B)\right) 
	\circ (p_{c,c'} \times c)
\arrow[swap]{u}{\iso} \\

\glueFun(\eval_{B,B'}) 
	\circ \big( \prodPres_{(\scriptsizeexpobj{B}{B'}, B)} 
	\circ (p_{c, c'} \times c) \big)
\arrow[swap]{r}{\iso} &
\big( \glueFun(\eval_{B,B'}) \circ 
	\prodPres_{(\scriptsizeexpobj{B}{B'}, B)}\big) 
	\circ (p_{c, c'} \times c)
\arrow[swap]{u}
	{\epsilonExp_{\eval \circ (\evBar \times \glueFun B)}^{-1} 
	\circ (p_{c,c'} \times c)}
\end{td}
Applying the coherence 
condition~(\ref{eq:lam-coherence-condition}), the first diagram in the claim 
reduces further to
\begin{equation} \label{eq:checking-counit-cylinder-diagram}
\begin{small}
\begin{tikzcd}[column sep = 1em]
c' 
	\circ \left( \left(\eval_{C, C'} \circ (q_{c,c'} \times C)\right) 
	\circ (\glLam(t) \times C)\right)
\arrow[swap]{d}{\iso} 
\arrow{r} &
\glueFun s \circ \left(\prodPres_{Q,B} \circ (r \times c)\right) \\
c' 
	\circ 
		\left(\eval_{C, C'} 
			\circ \big((q_{c,c'} \circ \glLam(t)\big) \times C\right) 
\arrow[swap]{d}{c' \circ \eval_{C,C'} \circ (\Delta_{C,C'} \times C)} &
\: \\
c' \circ \left(\eval_{C,C'} \circ (\lambda t \times C)\right)
\arrow[swap]{d}{\iso}
& 
\eval_{\glueFun B,C'} 
	\circ \left( \left(\lambda \compact{c} \circ \evBar_{B,B'}\right)
	\circ \left(\glueFun(\lambda s) \circ r\right) \right) 
\times C 
\arrow[swap]{uu}{\U_\alpha} \\
\left( c' \circ \eval_{C,C'} \right) \circ 
	\left( \lambda t \times C \right)
\arrow[swap]{d}
	{\epsilonExp_{(c' \circ \eval)}^{-1} \circ (\lambda t \times C)}  &
\: \\
\left(\eval_{\glueFun B, C'} 
	\circ (\lambda(c' \circ \eval_{C, C'}) \times C)\right) 
	\circ (\lambda t \times C) 
\arrow[swap]{r}{\iso} &
\eval_{\glueFun B,C'} \circ \big( \lambda(c' \circ \eval_{C, C'}) \circ 
\lambda t \big) 
\times C 
\arrow[swap]{uu}{\eval_{\glueFun B,C'} \circ (\L_\alpha \times C)}
\end{tikzcd}
\end{small}
\end{equation}
Next, by the 
definition of $\L_\alpha$ and the triangle law relating $\etaExp{}$ and 
$\epsilonExp$, one sees that
\begin{td}
\: &
\eval_{\glueFun B,C} \circ (\lambda h \times C) 
\arrow{dr}[description]{\eval_{\glueFun B,C} \circ \lambda (\etaExp{h}^{-1} 
\times C)}  
\arrow{r}{\epsilonExp_h} &
h 
\arrow[equals, bend left = 40]{d} \\

\eval_{\glueFun B,C} \circ \big( \lambda(c' \circ \eval_{C, C'}) \circ \lambda 
t \big) 
\times C  
\arrow[bend left = 40]{urr}[near end]{\K_\alpha}
\arrow{ur}[description]{\eval_{\glueFun B,C} \circ (\transExp{\K_\alpha} \times 
C)}
\arrow[swap]{rr}{\eval_{\glueFun B,C'} \circ (\L_\alpha \times C)} &
\: &
h
\end{td}
for $h :=  \eval_{C, \glueFun B} \circ 
\left( \left(\lambda\compact{c} \circ \evBar_{B,B'}\right) 
	\circ \left(\glueFun(\lambda s) \circ r \right)\right) \times C $. 
Hence, the composite~(\ref{eq:checking-counit-cylinder-diagram}) is equal to 
the 
anti-clockwise route around the diagram below, in which $(\dagger)$ 
abbreviates
\[
\left(c' \circ \eval_{C,C'}\right) \circ (\lambda t \times C) 
\XRA\iso
c' \circ \left( \eval_{C,C'} \circ (\lambda t \times C) \right)
\XRA{c' \circ \epsilonExp_t}
c' \circ t
\]
and the bottom two shapes commute by definition:
\begin{td}[column sep = 1em]
c' 
	\circ \left( \left(\eval_{C, C'} \circ (q_{c,c'} \times C)\right) 
	\circ (\glLam(t) \times C)\right)
\arrow[swap]{d}{\iso} &
\: \\

\left(c' \circ \eval_{C, C'}\right) 
	\circ \big(q_{c,c'} \circ \glLam(t)\big) \times C 
\arrow[swap]{d}{c' \circ \eval_{C,C'} \circ (\Delta_{c,c'} \times C)} &
\: \\

\left( c' \circ \eval_{C,C'} \right) \circ (\lambda t \times C) 
\arrow[swap]{d}{\epsilonExp_{(c' \circ \eval)}^{-1} \circ (\lambda t \times C)}
\arrow[equals, bend left = 17]{dr} &
\: \\

\left(
	\eval_{\glueFun B, C'} 
	\circ (\lambda(c' \circ \eval_{C, C'}) \times C)\right) 
	\circ (\lambda t \times C) 
\arrow{r}[swap, yshift=-1.5mm]{\epsilonExp_{(c' \circ \eval)} \circ (\lambda t 
\times C)}
\arrow[swap]{d}{\iso} &
\left(c' \circ \eval_{C,C'}\right) \circ (\lambda t \times C) 
\arrow{r}{(\dagger)} &
c' \circ t  
\arrow{dd}{\alpha}\\

\eval_{\glueFun B,C'} \circ 
	\big( \lambda(c' \circ \eval_{C, C'}) \circ \lambda t \big) \times C 
\arrow[bend right = 5]{urr}[swap, description]{\T_\alpha}
\arrow[swap]{d}{\K_\alpha} &
\:  &
\: \\

\eval_{\glueFun B,C'} 
	\circ \left( \left( \lambda \compact{c} \circ \evBar_{B,B'}\right) 
		\circ \left(\glueFun(\lambda s) \circ r \right) \right) \times C 
\arrow[swap]{rr}{\U_\alpha} &
\: &
\glueFun s \circ \left(\prodPres_{Q,B} \circ (r \times c)\right)
\end{td}
The clockwise route around this diagram is equal to the 2-cell given by the 
second diagram in the claim, so the proof is complete.
\end{proof}
\end{mylemma} 

We have now constructed all the data we shall require. It 
remains to show 
that, together, it defines an adjoint equivalence
\[
\glued{\lambda} : 
\gl{\glueFun }{\big(\glued{R} \times \glued{C}, \glued{C'}\big)} 
\leftrightarrows 
\gl{\glueFun }\big(\glued{R}, \expobj{\glued{C}}{\glued{C'}}\big) 
: \glued{\eval}_{\glued{C}, \glued{C'}} \circ (- \times \glued{C})
\]
Thus, we need to check that for every pair of 1-cells 
$\glued{g} : \glued{R} \to (\exp{\glued{C}}{\glued{C'}})$ and 
$\glued{t} : \glued{R} \times \glued{C} \to \glued{C'}$ 
related by a 2-cell 
$\glued{\tau} := (\tau, \sigma) : 
\glued{\eval}_{\glued{C}, \glued{C'}} \circ (\glued{g} \times \glued{C}) \To 
\glued{t}$, there exists a 2-cell 
$\transExp{\glued{\tau}} : \glued{g} \To \glued{\lambda}\glued{t}$, unique such 
that
\begin{equation} \label{eq:ump-for-exponentials}
\begin{tikzcd}[column sep = 4em]
\glued{\eval}_{\glued{C}, \glued{C'}} \circ (\glued{g} \times \glued{C})
\arrow[swap]{dr}{\glued{\tau}}
\arrow{rr}
{\glued{\eval}_{\glued{C}, \glued{C'}} \circ (\transExp{\glued{\tau}} \times 
\glued{C})} &
\: &
\glued{\eval}_{\glued{C}, \glued{C'}} \circ (\glued{\lambda}\glued{t} \times 
\glued{C}) 
\arrow{dl}{\glued{\epsilonExp}_t} \\
\: &
\glued{t} &
\:
\end{tikzcd}
\end{equation}
We turn to this next.

\paragraph*{Universality of $\glued{\epsilon} = (\gluedCo{}, \epsilonExp)$.}
We begin with the existence part of the claim. Let 
$\glued{g} := (g, \gamma, f) : (R, r, Q) \to (\glexp{C}{C'}, p_{c,c'}, 
\expobj{B}{B'})$ 
and 
$\glued{t} := (t, \alpha, s) : 
(R \times C, \prodPres_{Q,B} \circ (r \times c), Q \times B) \to
(C', c', B')$ be 1-cells and suppose that 
$\glued{\tau} := (\tau, \sigma) : 
\glued{\eval}_{\glued{C}, \glued{C'}} \circ (\glued{g} \times \glued{C}) \To 
\glued{t}$. Thus, $\tau$ and $\sigma$ have type
\begin{gather*}
\tau : \left(\eval_{C, C'} \circ (q_{c, c'} \times C)\right) \circ (g \times C) 
	\To t 
\\
\sigma : \eval_{B,B'} \circ (f \times B) \To s
\end{gather*}
and we are required to provide 2-cells $\altTrans{\tau}$ and 
$\altTrans{\sigma}$ of type
\begin{align*}
\altTrans\tau : g &\To \glLam(t) \\
\altTrans\sigma : f &\To \lambda s
\end{align*}
satisfying the cylinder condition.
For the second component we can simply take 
$\transExp{\sigma}$. For the first component we use the universal property of 
pullbacks. We aim to construct a pair of 2-cells 
\begin{align*}
p_{c, c'} \circ g &\To \glueFun (\lambda s) \circ r \\
q_{c, c'} \circ g &\To \lambda t
\end{align*}
such that the coherence condition~(\ref{eq:lam-coherence-condition}) holds. We 
claim that the following 2-cells suffice
\begin{equation} \label{eq:def-of-Sigma-1-and-Sigma-2}
\begin{aligned} 
\Sigma_1 &:= p_{c, c'} \circ g 
\XRA{\gamma} \glueFun (f) \circ r 
\XRA{\glueFun (\transExp{\sigma}) \circ r} \glueFun (\lambda s) \circ r \\
\Sigma_2 &:= 
	q_{c,c'} \circ g 
		\XRA{\transExplr{\chi}} 
	\lambda t
\end{aligned}
\end{equation}
where
$\chi :=
	\eval_{C,C'} \circ \left((q_{c,c'} \circ g) \times C \right)
	\XRA\iso 
	\left(\eval_{C,C'} \circ (q_{c,c'} \times C)\right) \circ (g \times c) 
	\XRA{\tau} 
	\lambda t
$. 
The required coherence condition is the subject of the following lemma.

\begin{mylemma} \label{lem:fill-in-condition-for-sigmas}
Consider a pair of 1-cells
\begin{align*}
\glued{g} &:= 
	(g, \gamma, f) : 
		(R, r, Q) 
			\to 
		(\glexp{C}{C'}, p_{c,c'}, \expobj{B}{B'})  \\
\glued{t} &:= 
	(t, \alpha, s) : 
		(R \times C, \prodPres_{Q,B} \circ (r \times c), Q \times B) 
			\to
	(C', c', B')
\end{align*} 
in $\gl{\glueFun}$ related by a 2-cell 
$\glued{\tau} 
	:= (\tau, \sigma) : 
		\glued{\eval}_{\glued{C}, \glued{C'}} \circ 
			(\glued{g} \times \glued{C}) 
		\To 
		\glued{t}$.
Then, where  $\Sigma_1$ and $\Sigma_2$ are defined 
in~(\ref{eq:def-of-Sigma-1-and-Sigma-2}), the following diagram commutes:
\begin{equation*}
\makebox[\textwidth]{
\begin{tikzcd}[ampersand replacement = \&]
\left(\lambda(c' \circ \eval_{C,C'}) \circ q_{c,c'}\right) \circ g
\arrow{r}{\iso} 
\arrow[swap]{d}{\omega_{c,c'} \circ g} \&
\lambda(c' \circ \eval_{C,C'}) \circ \left( q_{c,c'} \circ g\right)
\arrow{r}[yshift=2mm]{\lambda(c' \circ \eval_{C,C'}) \circ \Sigma_2} \&
\lambda(c' \circ \eval_{C,C'}) \circ \lambda t
\arrow{d}{\L_\alpha} \\
\left( \left(\lambda\compact{c} \circ \evBar_{B,B'}\right) 
	\circ p_{c,c'}\right) 
	\circ g
\arrow[swap]{r}{\iso} \&
 \left(\lambda\compact{c} \circ \evBar_{B,B'}\right)
	\circ \left( p_{c,c'} \circ g \right)
\arrow[swap]{r}[yshift=-2mm]{\lambda \compact{c} \circ \evBar_{B,B'} \circ 
\Sigma_1} \&
\left(\lambda\compact{c} 
	\circ \evBar_{B,B'}\right)
	\circ \left(\glueFun(\lambda s) \circ r \right)
\end{tikzcd}
}
\end{equation*}
\begin{proof}
Straightforward manipulations and an application of the cylinder condition on 
$\glued{\tau}$ 
unfolds the clockwise route to the 
following composite:
\begin{equation} \label{eq:mediating-2-cell-calculation}
\begin{small}
\begin{tikzcd}[column sep = 0.5em]
\left(\lambda(c' \circ \eval_{C, C'}) \circ q_{c,c'}\right) \circ g 
\arrow[swap]{dd}{\etaExp{}}
\arrow{r} &
\left(\lambda\compact{c} \circ \evBar_{B,B'}\right) 
	\circ \left(\glueFun(\lambda s) \circ r \right)
\\
\: &
\lambda {\left( 
	\eval_{C, \glueFun B} 
	\circ \left( \left(\lambda\compact{c} \circ \evBar_{B,B'}\right) 
	\circ  \left(\glueFun(\lambda s) \circ r \right)\right) 
	\times C 
	\right)} 
\arrow[swap]{u}{\etaExp{}^{-1}} \\
\lambda\left( 
	\eval_{C, \glueFun B'} 
	\circ \left( \left(\lambda(c' \circ \eval_{C, C'})
 		\circ q_{c,c'}\right) \circ g\right) \times C\right) 
\arrow[swap]{r}{\lambda \zeta} &
\lambda {\left(\glueFun s 
			\circ \left(\prodPres_{Q,B} \circ (r \times c)\right)\right)}
\arrow[swap]{u}{\lambda\U_\alpha^{-1}}
\end{tikzcd}
\end{small}
\vspace{3mm}
\end{equation}
Here 
$\zeta : 
	\eval_{C, \glueFun B'} 
	\circ \left( \left(\lambda(c' \circ \eval_{C, C'})
 		\circ q_{c,c'}\right) \circ g\right) \times C 
	\to 
 	\glueFun s \circ \left(\prodPres_{Q,B} \circ (r \times c)\right)$ 
is the composite defined by  commutativity of the following diagram:
\begin{equation*}
\makebox[\textwidth]{
\begin{tikzcd}[column sep = 0em, ampersand replacement = \&]
\eval_{C, \glueFun B'} 
	\circ \left( \left(\lambda(c' \circ \eval_{C, C'})
 		\circ q_{c,c'}\right) \circ g\right) \times C
\arrow{r}{\zeta}
\arrow[swap]{d}{\iso} \&
\glueFun s \circ \left(\prodPres_{Q,B} \circ (r \times c)\right) \\
\left(\eval_{C, \glueFun B'} 
	\circ \left(\lambda(c' \circ \eval_{C, C'}) \times C\right)\right)
	\circ \left( (q_{c,c'} \circ g) \times C \right)
\arrow[swap]{d}
	{\epsilonExp_{(c' \circ \eval)} \circ (qg \times C)} \&
\glueFun\left(\eval_{B,B'} \circ (f \times B)\right) 
	\circ \left(\prodPres_{Q,B} 
	\circ (r \times c) \right)
\arrow[swap]{u}{\glueFun (\sigma) \circ \prodPres \circ (r \times c)} \\
(c' \circ \eval_{C,C'}) \circ \left( (q_{c,c'} \circ g) \times C \right) 
\arrow[swap]{d}{\iso} \&
\left(\glueFun(\eval_{B,B'}) 
	\circ \glueFun(f \times B)\right) 
	\circ \left(\prodPres_{Q,B} 
	\circ (r \times c)\right)
\arrow[swap]{u}
	{\phi^\glueFun _{\eval, f \times B} \circ \prodPres \circ (r \times c)}\\
\left(c' 
	\circ \left(\eval_{C,C'} \circ (q_{c,c'} \times C)\right)\right)
	\circ (g \times C)
\arrow[swap]{d}{\evalMod_{\glued{C}, \glued{C'}} \circ (g \times C) } \&
\left(\glueFun(\eval_{B,B'}) 
	\circ \left(\glueFun(f \times B) 
	\circ \prodPres_{Q,B}\right)\right) 
	\circ (r \times c)
\arrow[swap]{u}{\iso} \\
\big(\glueFun(\eval_{B,B'}) 
	\circ \big(\prodPres_{(\scriptsizescriptsizeexpobj{B}{B'}, B)} 
	\circ (p_{c,c'} \times c) \big)\big)
	\circ (g \times C)
\arrow[swap]{d}{\iso} \&
\:
\\
\big(\glueFun(\eval_{B,B'}) 
	\circ \prodPres_{(\scriptsizescriptsizeexpobj{B}{B'}, B)}\big) 
	\circ \big((p_{c,c'} \circ g) \times c\big) 
\arrow[swap]{d}
	{\glueFun(\eval) \circ \prodPres \circ (\gamma \times c)} \&
\big(\glueFun(\eval_{B,B'}) 
	\circ \big(\prodPres_{(\scriptsizescriptsizeexpobj{B}{B'}, B)} 
	\circ (\glueFun f \times \glueFun\Id_B)\big)\big)
	\circ (r \times c) 
\arrow[swap]{uu}{\glueFun(\eval) \circ \nat_{f, \Id} \circ (r \times c)}  \\
\big(\glueFun(\eval_{B,B'}) 
	\circ \prodPres_{(\scriptsizescriptsizeexpobj{B}{B'}, B)}\big)
	\circ \big((\glueFun f \circ r) \times c\big) 
\arrow[swap]{r}{\iso} \&
\big(\glueFun(\eval_{B,B'}) 
	\circ \big(\prodPres_{(\scriptsizescriptsizeexpobj{B}{B'}, B)} 
	\circ (\glueFun f \times \glueFun B)\big)\big)
	\circ (r \times c) 
\arrow[swap]{u}
	{\glueFun(\eval) 
		\circ \prodPres 
		\circ (\glueFun f \times \psi^\glueFun _B) 
		\circ (r \times c) }
\end{tikzcd}
}
\end{equation*}
A short calculation shows that the following also commutes:
\begin{td}
\eval_{C, \glueFun B'} 
	\circ \left( \left(\lambda(c' \circ \eval_{C, C'})
 		\circ q_{c,c'}\right) \circ g\right) \times C
\arrow[swap]{d}{\eval \circ (\omega_{c,c'} \circ g) \times C}
\arrow{r}{\zeta} &
\glueFun s \circ \left(\prodPres_{Q,B} \circ (r \times c)\right) \\
\eval_{C, \glueFun B'} 
	\circ \left(\left( \left(\lambda\compact{c} \circ \evBar_{B,B'}\right) 
		\circ p_{c,c'}\right) \circ g\right)
	\times C 
\arrow[swap]{d}{\iso} &
\: \\
\eval_{C, \glueFun B'} 
	\circ \left( \left(\lambda\compact{c} \circ \evBar_{B,B'}\right) 
		\circ \left(p_{c,c'}\circ g\right)\right)
	\times C 
\arrow[swap]{r}[yshift=-2mm]
	{\eval \circ 
			(\lambda\compact{c} \circ \evBar \circ \Sigma_1) 
			\times C } &
\eval_{C, \glueFun B'} 
	\circ \left(\left(\lambda\compact{c} \circ \evBar_{B,B'}\right) 
		\circ \left(\glueFun(\lambda s) \circ r\right) \right) \times C  
\arrow[swap]{uu}{\U_\alpha}
\end{td}
Substituting this back into~(\ref{eq:mediating-2-cell-calculation}) and 
applying the naturality of $\etaExp{}$, one obtains the anticlockwise route 
around 
the claim, as required. 
\end{proof}
\end{mylemma}

It follows that $(g, \Sigma_1, \Sigma_2)$ is a fill-in.
By the universality of the fill-in $(\glLam(t), \Gamma, \Delta)$, therefore, 
one obtains a 2-cell 
$\trans\Sigma : g \To \glLam(t)$, unique such that the following two diagrams 
commute \big(\cf~(\ref{eq:fill-in-mediating-property})\big):
\begin{equation} \label{eq:ump-of-trans-Sigma}
\begin{tikzcd}
p_{c, c'} \circ g 
\arrow{r}{\gamma}
\arrow[swap]{d}{p_{c,c'} \circ \trans\Sigma} &
\glueFun (f) \circ r 
\arrow{d}{\glueFun (\transExp{\sigma}) \circ r} \\
p_{c,c'} \circ \glLam(t)
\arrow[swap]{r}{\Gamma_{c,c'}} &
\glueFun (\lambda s) \circ r
\end{tikzcd}
\qquad\qquad
\begin{tikzcd}
q_{c,c'} \circ g 
\arrow[swap]{d}{q_{c,c'} \circ \trans\Sigma} 
\arrow{dr}{\transExp{\chi}} &
\: \\
q_{c,c'} \circ \glLam(t) 
\arrow[swap]{r}{\Delta_{c,c'}} &
\lambda t
\end{tikzcd}
\end{equation}
We therefore define the components of $\transExp{\glued{\tau}}$ as follows:
\begin{equation} \label{eq:def-of-trans-tau-trans-sigma}
\begin{aligned}
\altTrans\tau &:= \trans\Sigma :  g \To \glLam(t) \\
\altTrans\sigma &:= \transExp{\sigma} : f \To \lambda s
\end{aligned}
\end{equation}
Note that the left-hand diagram of~(\ref{eq:ump-of-trans-Sigma}) establishes 
this pair is a 2-cell in $\gl{\glueFun }$. We need to show that this 2-cell 
makes~(\ref{eq:ump-for-exponentials}) commute. For 
the second component, this holds by assumption. For the first component, we 
observe that $\gluedCo{t}$ is the right-hand leg of the following diagram:
\begin{equation*}
\makebox[\textwidth]{
\begin{tikzcd}[column sep = 0em, ampersand replacement = \&]
\left(\eval_{C,C'} \circ (q_{c,c'} \times C)\right) \circ (g \times C)
\arrow[phantom]{drr}[description]{\equals{nat.}}
\arrow[bend right]{dddrr}[swap]{\tau}
\arrow{dr}[description]{\iso}
\arrow{rr}[yshift=0mm]{\eval_{C,C'} \circ (q_{c,c'} \times C) \circ 
(\trans\Sigma \times C)} \&
\: \&
\left(\eval_{C,C'} 
	\circ (q_{c,c'} \times C) \right)
	\circ (\glLam(t) \times C)
\arrow{d}{\iso}  \\
\: \&
\eval_{C,C'} \circ \big((q_{c,c'} \circ g\big) \times C) 
\arrow[bend right, swap]{ddr}{\chi}
\arrow[phantom]{dd}[description]{\equals{def.}}
\arrow[bend right = 12, swap]{dr}
\arrow[phantom]{ddr}[description, yshift=-4mm, xshift=2mm]{\equals{UMP}}
\arrow{r}[yshift=2mm]{\eval_{C,C'} \circ (\trans\Sigma \times C)} \&
\eval_{C,C'} \circ \big((q_{c,c'} \circ \glLam(t)\big) \times C) 
\arrow{d}{\eval_{C,C'} \circ (\Delta_{c,c'} \times C)} \\
\: \&
\: \&
\eval_{C, C'} \circ (\lambda t \times C)
\arrow{d}{\epsilonExp_t} \\
\: \&
\: \&
t
\end{tikzcd}
}
\end{equation*}
The unlabelled inner arrow is 
$\eval_{C,C'} \circ (\transExp{\chi} \times C)$ (where 
$\chi$ is defined just after~(\ref{eq:def-of-Sigma-1-and-Sigma-2})), 
so the triangular shape commutes by~(\ref{eq:ump-of-trans-Sigma}). This 
completes the existence part of the universality claim; we record 
our progress so far in the following lemma.

\begin{prooflesslemma}
For any triple of 1- and 2-cells as in 
Lemma~\ref{lem:fill-in-condition-for-sigmas}, the pair
$\transExp{\glued{\tau}} := (\trans\Sigma, \transExp{\sigma})$ defined 
in~(\ref{eq:def-of-trans-tau-trans-sigma}) is a 2-cell 
$\glued{g} \To \glued{\lambda} \, \glued{t}$ 
in $\gl{\glueFun }$ satisfying (\ref{eq:ump-for-exponentials}).
\end{prooflesslemma}

It remains to show uniqueness. Suppose given a 2-cell 
$\glued{\theta} : \glued{g} \To \glued{\lambda} \, \glued{t}$ in 
$\gl{\glueFun}$ with 
components
\begin{align*}
\theta : g &\To \glLam(t) \\
\vartheta : f &\To \lambda s
\end{align*}
such that $\glued{\theta}$ fills~(\ref{eq:ump-for-exponentials}). Examining the 
second component, it is immediate from the universal property of 
$\transExp{\sigma}$ that $\transExp{\sigma} = \vartheta$. For the first 
component, 
we show that $\theta = \trans\Sigma$ by showing that $\theta$ satisfies the two 
diagrams of~(\ref{eq:ump-of-trans-Sigma}). For the left-hand diagram, the 
cylinder condition on $\glued{\theta}$ requires that
\begin{td}[column sep = 3em]
p_{c, c'} \circ g 
\arrow{r}{\gamma}
\arrow[swap]{d}{p_{c,c'} \circ \theta} &
\glueFun (f) \circ r 
\arrow{d}{\glueFun (\vartheta) \circ r} \\
p_{c,c'} \circ \glLam(t) 
\arrow[swap]{r}{\Gamma_{c,c'}} &
\glueFun (\lambda s) \circ r
\end{td}
But we already know that $\vartheta = \transExp{\sigma}$, so the required 
diagram commutes. For the right-hand diagram, it follows 
from~(\ref{eq:ump-for-exponentials}) and the definition of $\gluedCo{t}$ that 
the following commutes:
\begin{td}[column sep = 1.5em]
\eval_{C,C'} \circ \left((q_{c,c'} \circ g) \times C\right)
\arrow[swap]{d}{\eval_{C,C'} \circ (q_{c,c'} \circ \theta) \times C} 
\arrow{rr}{\iso} &
\: &
\left(\eval_{C,C'} \circ (q_{c,c'} \times C)\right) \circ (g \times C) 
\arrow{d}{\tau}  \\
\eval_{C,C'} \circ \left(\left(q_{c,c'} \circ \glLam(t)\right) \times C\right) 
\arrow[swap]{r}[yshift=-2mm]{\eval_{C,C'} \circ (\Delta_{c,c'} \times C)} &
\eval_{C,C'} \circ (\lambda t \times C)
\arrow[swap]{r}{\epsilonExp_t} &
t
\end{td}
The claim then holds by the universal property of 
$\transExp{\vartheta}$. Thus:

\begin{prooflesslemma}
For any triple of 1- and 2-cells as in 
Lemma~\ref{lem:fill-in-condition-for-sigmas}, the pair
$\transExp{\glued{\tau}} := (\trans\Sigma, \transExp{\sigma})$ defined 
in~(\ref{eq:def-of-trans-tau-trans-sigma}) is the unique 2-cell 
$\glued{g} \To \glued{\lambda} \, \glued{t}$ 
in $\gl{\glueFun }$ satisfying (\ref{eq:ump-for-exponentials}).
\end{prooflesslemma}
This completes the proof that for any
$\glued{R}, \glued{C}$ and $\glued{C'}$ in $\gl{\glueFun}$ the diagram
\[
\glued{\lambda} : 
\gl{\glueFun }\big(\glued{R} \times \glued{C}, \glued{C'}\big) 
\leftrightarrows 
\gl{\glueFun }\big(\glued{R}, \expobj{\glued{C}}{\glued{C'}}\big) 
: \glued{\eval}_{\glued{C}, \glued{C'}} \circ (- \times \glued{C})
\]
is an adjoint equivalence, and hence the proof of 
Theorem~\ref{thm:glueing-bicat-cartesian-closed}.

%

\chapter{\texorpdfstring{Normalisation-by-evaluation for 
$\langCartClosed$}{Coherence via normalisation-by-evaluation}} 
\label{chap:nbe}

We now turn to the main result of this thesis, namely the coherence result 
for cartesian closed bicategories. Our strategy is to employ a bicategorical 
treatment of the
\Def{normalisation-by-evaluation} proof technique. It is well-known that the 
na\"ive strategy for proving strong normalisation of the simply-typed lambda 
calculus---by a straightforward structural induction on terms---fails because 
an application $\app{t}{u}$ may contain redexes that do not occur in 
either $t$ or $u$. One classical solution, originally due to Tait~\cite{Tait}, 
is to strengthen the inductive hypothesis using \Def{reducibility predicates}. 
This approach was refined by Girard~\cite{Girard1972}, who introduced the 
notion of \Def{neutral terms}. These can be viewed as the obstructions to the 
normalisation proof: they are the 
terms whose introduction rules may introduce new $\beta$-redexes.

Normalisation-by-evaluation provides an alternative strategy: as a slogan, one 
`inverts the evaluation functional' to construct a mapping from 
neutral to normal terms. 
Loosely speaking, one constructs a 
model with enough intensional information to pass back and forth 
between semantics and syntax. One \emph{quotes} a morphism $f$ to a (normal) 
term in the 
syntax, and \emph{unquotes} a term $t$ to a morphism in 
the semantics 
(these operations are also known as \emph{reify} and \emph{reflect}). 

The intuition is---very roughly---as follows. Consider a 
semantics $\sem{-}$ for the simply-typed lambda calculus, determined by a 
choice of cartesian closed category and an interpretation of the base types, 
and suppose that one has constructed mappings $\quote$ and $\unquote$ between 
the syntax and semantics, as indicated above. 
For a term $(x : A \vdash t : B)$ one has an 
interpretation $\sem{t} : \sem{A} \to \sem{B}$. Now, where $x$ is a generic 
fresh variable, $\unquote(x) : \sem{A}$. So one may evaluate $\sem{t}$ at 
$\unquote(x)$
to obtain a normal term 
$\quote\left(\sem{t}\left(\unquote(x)\right)\right)$ of 
type $B$. The 
normal form of $\lam{x}{t}$ is then 
$\lam{x}{\,\quote\big( {\sem{t}\left(\unquote(x)\right)} \big)}$. 

First introduced by  
Berger \& Schwichtenberg~\cite{Berger1991} for the simply-typed lambda 
calculus, normalisation-by-evaluation has 
become a standard tool for tackling normalisation problems. It has been 
extended to a number of richer calculi, including the
simply-typed lambda calculus with sum 
types~\cite{Altenkirch2001}, versions of Martin-L{\"o}f type 
theory (\eg~\cite{Abel2007, Altenkirch2016dep, Altenkirch2017}),  and even to type theories with 
algebraic effects~\cite{Staton2013}. Moreover, the normalisation algorithm one 
extracts from normalisation-by-evaluation is generally highly efficient, which 
has led to significant study for applications in interactive proof systems 
(see~\eg~\cite{Berger1998}).

Here we follow in the vein of \emph{categorical} reconstructions of the 
normalisation-by-evaluation argument 
(\eg~\cite{Altenkirch1995, Coquand1997, Cubric1998, Fiore2002}). 
In particular, the argument we present closely follows~\cite{Fiore2002}; 
the reliance on categorical properties there lends itself 
especially to bicategorical translation.

The chapter is arranged as follows. We begin in Section~\ref{sec:fiore-proof} 
by briefly recapitulating the argument of~\cite{Fiore2002}. 
In Sections~\ref{sec:syntax-as-pseudofunctors}--\ref{sec:nbe-glueing} we show 
how the crucial elements of this argument can be lifted to the bicategorical 
setting. Section~\ref{sec:main-result} presents the main result of this 
thesis: $\langCartClosed$ is locally coherent.

\section{Fiore's categorical normalisation-by-evaluation proof}
\label{sec:fiore-proof}

We extract the bare bones of Fiore's argument~\cite{Fiore2002}. The intention 
is not to 
provide the reader with the full proof, but to waypoint the key steps in the 
bicategorical argument we present thereafter.

\paragraph*{Syntax as presheaves.}
For any set of base types $\baseTypes$, let 
$\Con_{\allTypes\baseTypes}$ 
denote the free strict
cocartesian category on the set $\allTypes\baseTypes$
generated by the grammar
\[
X_1, \,\dots\, , X_n, Y, Z ::= 
	B \st
	\prodop_n(X_1, \,\dots\, , X_n) \st
	\exptype{Y}{Z}
	\qquad
	(B \in \baseTypes, n \in \Nat)
\]
Explicitly, this is the 
comma category 
$(\mathbb{F} \downarrow \allTypes{\baseTypes})$, where $\mathbb{F}$ is a 
skeleton of the category of 
finite sets and all set-theoretic functions. For our purposes, however, we 
identify it with the \Def{category of contexts}, in which the objects are 
contexts (defined by Figure~\ref{fig:nbe:contexts}, below) and the 
morphisms are context renamings. Note that we index from 0 to avoid awkward 
off-by-one manipulations.
\begin{rules}
\centering
\unaryRule	{\faketext}
			{\diamond \mathrm{\: ctx}}
			{\qquad\qquad} 
\binaryRule	{\Gamma \mathrm{\: ctx}}
			{\len{\Gamma} = n}
			{\Gamma, x_{n} : A \mathrm{\: ctx}} 
			{$\big(A \in \allTypes{\baseTypes}\big)$} \vspace{-\treeskip}
\caption{Rules for contexts\label{fig:nbe:contexts}}
\end{rules}

To ensure that that $\Con_{\allTypes{\baseTypes}}$ is strict cocartesian, we 
stipulate that variables are named in order according to a fixed 
enumeration. 
However, following our standing 
abuse~(Notation~\ref{not:multiple-variable-names}), we 
shall freely employ 
more indicative variable names, such as using $f$ to denote a variable of 
exponential type. 

An object
$\gamma : [n] \to \allTypes\baseTypes$ 
(for $[n] = \{0, \,\dots\, , n-1\} \in \mathbb{F}$) in 
$(\mathbb{F} \downarrow \allTypes{\baseTypes})$ 
corresponds to the context
$\left(x_i : \gamma(i)\right)_{i=1, \,\dots\, , n}$. A morphism 
$h : \gamma \to \delta$, namely a set map 
$[n] \to [m]$ such that the diagram below commutes,
corresponds 
to the context 
renaming $x_i \mapsto x_{hi}$.
\begin{equation*}
\begin{tikzcd} 
{[n]}
\arrow[swap]{dr}{\gamma}
\arrow{rr}{h} &
\: &
{[m]}
\arrow{dl}{\delta} \\
\: &
\allTypes\baseTypes &
\: 
\end{tikzcd}
\end{equation*}
The coproduct $\Gamma + \Delta$ is the concatenated context 
$\Gamma \concat \Delta$. 

We denote the universal embedding of $\allTypes\baseTypes$ into 
$\Con_{\allTypes\baseTypes}$ by $\coerce{-}$; thus, 
$\coerce{A}$ \Def{coerces} the type $A$ into the unary context $(x_1 : A)$, and 
the coproduct $\Gamma + \coerce{A}$ is the weakening of $\Gamma$ by a 
variable of type $A$. The notation is chosen to suggest a list of length one.

In the tradition of \Def{algebraic type theory} (\eg~\cite{Fiore1999, Fiore2012types}), the 
category 
$\Psh{\op{\Con_{\allTypes\baseTypes}}}$ 
of covariant presheaves $\Con_{\allTypes\baseTypes} \to \Set$
provides a semantic universe for 
the study of abstract syntax. For example, for the simply-typed lambda calculus 
$\stlc(\baseTypes)$ over 
$\baseTypes$, the set of terms-in-context of a 
given type $B$ (modulo $\alpha$-equivalence) define a presheaf 
$\langTerms(-; B)$ by 
$\langTerms(\Gamma;B) := \slice{\{ {t \st \Gamma \vdash t : B} \}\:}{\aeq}$. 
The functorial action is given by context renamings: for contexts 
$\Gamma := ({x_i : A_i})_{i=1, \,\dots\, , n}$ and
$\Delta := (y_j : B_j)_{j=1, \,\dots\, , m}$ and a 
context renaming $r : \Gamma \to \Delta$, one obtains a mapping
\begin{align*}
\langTerms(\Gamma; B) &\to \langTerms(\Delta; B) \\
t &\mapsto t[r(x_i) / x_i]
\end{align*}
by the admissibility of the rule
\begin{center}
\binaryRule{\Gamma \vdash t : B}
			{r : \Delta \to \Gamma}
			{\Delta \vdash t[r(x_i) / x_i] : B}
			{} \vspace{-\treeskip}
\end{center}
The Yoneda embedding $\yon$ yields a \Def{presheaf of variables}:
for any type $A \in \allTypes\baseTypes$ and context $\Gamma$,
$\yon(\coerce{A})(\Gamma) = \yon(x : A)(\Gamma) = \Con_{\allTypes\baseTypes}((x: A), \Gamma)$ 
corresponds to the set of inclusions of contexts 
$(x : A) \hookrightarrow \Gamma$.
This determines a presheaf 
$\varTerms(-; A)$ defined by
$\varTerms(\Gamma; A) = \{ x \st \Gamma \vdash x : A \}$. The well-known fact 
that $\altexp{\yon X}{P} \iso P(- \times X)$ in any presheaf 
category over a cartesian category corresponds to the observation that the exponential presheaf
$\altexp{\yon A}{\langTerms(-; B)}$ consists of terms of type $B$ in the 
extended context
$\Gamma + \coerce{A}$
(note that, since $\Con_{\allTypes\baseTypes}$ is strict cocartesian, its 
opposite category is strict cartesian).

\paragraph*{Intensional Kripke relations}
We extend the \Def{Kripke logical relations of varying arity} 
of~\cite{Jung1993, Alimohamed1995} to a category of 
\Def{intensional Kripke relations}. Encoding this extra intensional information 
allows one to extract a normalisation algorithm from the proof. Abstractly, the 
key to this construction is 
the \Def{relative hom-functor} (also known as the \Def{nerve} functor). 
For any functor $ \glueFun  : \catB \to \catX$ the left Kan extension 
$\lanext{\glueFun } := \lan{\glueFun }{(\yon)}$ exists as in the following 
diagram, in which 
$\Psh{\catB}$ denotes the presheaf category:
\begin{equation} \label{eq:glueing-left-kan-diagram}
\begin{tikzcd}
\catB 
\arrow[swap]{dr}{\glueFun }
\arrow[hookrightarrow]{rr}{\yon} &
\: 
\arrow[phantom]{d}[description, font=\scriptsize, 
yshift=1mm]{\twocellDown{\mathrm{lan}}} &
\Psh{\catB} \\
\: &
\catX 
\arrow{ur}[swap]{\lanext{\glueFun }} &
\:
\end{tikzcd}
\end{equation}
Explicitly, 
$\lanext{\glueFun }(X) := \catX\big( \glueFun (-), X\big) : \op\catB \to 
\Set$ and 
$\mathrm{lan}_B : \catB(-, B) \To \catX\big(\glueFun (-), \glueFun B\big)$ is 
just the functorial action of $\glueFun $.  This construction is particularly 
well-known in the context of profunctors 
(distributors), since 
$\catB\big(\glueFun (-), X\big)$ and $\catB\big(X, \glueFun (-)\big)$ provide 
canonical (indeed, adjoint) profunctors 
$\catX \nrightarrow \catB$ for every functor 
$\glueFun  : \catB \to \catX$ (\eg~\cite[Example 7.8.3]{Borceux1994}).  

\begin{mydefn} \quad \label{def:relative-hom-functor}
\begin{enumerate}
\item For $\glueFun  : \catB\to \catX$ a functor, the 
\Def{relative hom-functor} 
is the functor $\lanext{\glueFun } : \catX \to \Psh\catB$ defined above.
\item For a category $\catB$ and a functor $\glueFun  : \catB\to \catX$, the 
category of
\Def{$\catB$-intensional Kripke relations of arity $\glueFun $} is the glueing 
category 
$\gl{\lanext{\glueFun }}$ associated to the relative hom-functor. \qedhere
\end{enumerate}
\end{mydefn}

The relative hom-functor preserves limits so, when $\catX$ is cartesian closed, 
the glueing category $\gl{\lanext{\glueFun}}$ is cartesian closed and the 
forgetful functor to $\catX$ strictly preserves products and exponentials. 
Moreover, the 
Yoneda embedding extends to an embedding 
$\glued{\yon} : \catB \to \gl{\lanext{\glueFun }}$
by 
$\glued{\yon}(B) := 
	\left(\yon(B), \yon(B) \XRA{\mathrm{lan}_B} \lanext{\glueFun }(\glueFun 
	B), \glueFun B  \right)$.

Consider now the following situation. Fix a set of base types $\baseTypes$ and 
an interpretation $h : \baseTypes \to \catX$ in a cartesian closed category 
$\catX$. By the cartesian closed structure, this extends to a map
$\allTypes\baseTypes \to \catX$ we also denote by $h$. Applying the
universal property, $h$ extends in turn  
to a cartesian functor 
$\prodext{h}{} : \op{\Con_{\allTypes\baseTypes}} \to \catX$ 
interpreting all 
contexts within $\catX$. Moreover, writing $\mathcal{F}(\allTypes\baseTypes)$
for the free cartesian closed category on $\allTypes\baseTypes$, namely the 
syntactic model of the simply-typed lambda calculus $\stlc(\baseTypes)$, the 
coercion 
$\coerce{-} : 
	\allTypes{\baseTypes} \hookrightarrow \Con_{\allTypes{\baseTypes}}$
extends to a cartesian functor 
$\Con_{\allTypes{\baseTypes}} \to \mathcal{F}(\allTypes\baseTypes)$. By the 
various uniqueness properties, this factors the semantic interpretation 
$h\sem{-} : \mathcal{F}(\allTypes{\baseTypes}) \to \catX$ 
extending $h$. 
The situation is summarised in the following diagram.
\begin{equation}
\begin{tikzcd} \label{eq:h-to-interpret}
\: &
\mathcal{F}(\allTypes{\baseTypes})
\arrow[dashed]{dr}{h\sem{-}} &
\: 
\\
\op{\Con_{\allTypes\baseTypes}}
\arrow[dashed]{ur}
\arrow[swap, dashed]{rr}{\prodext{h}{}} &
\: &
\catX &
\: \\
\allTypes{\baseTypes}
\arrow[hookrightarrow]{u}{\coerce{-}}
\arrow[swap]{urr}{} &
\: \\
\baseTypes
\arrow[hookrightarrow]{u}
\arrow[swap]{uurr}{h}
\end{tikzcd}
\end{equation}
Note in particular that $\prodext{h}{\Gamma} = h\sem{\Gamma}$ for every context 
$\Gamma \in \Con_{\allTypes\baseTypes}$, and that for any type 
$A \in \allTypes{\baseTypes}$ the interpretation $h\sem{A}$ is equal to
$\prodext{h}{\coerce{A}}$. (Here we use the assumption that
$\prodop_1(X) = X$ to identify $h\sem{x: A}$ with $h\sem{A}$.)

An object in the category $\gl{\lanext{\prodext{h}{}}}$ of 
$\Con_{\allTypes\baseTypes}$-intensional Kripke relations of arity 
$\prodext{h}{}$ 
then consists of a presheaf $P :\Con_{\allTypes\baseTypes} \to \Set$ (which one 
might think of as \emph{syntactic} intensional information),
an object $X \in \catX$, and a natural transformation
$\pi : P \To \catX{\left(\prodext{h}{(-)}, X\right)}$
(which one might think of as \emph{semantic} information). One may think of 
this category as internalising the relationship between syntax and semantics 
required for the normalisation-by-evaluation argument.

\paragraph*{Neutral and normal terms as glued objects.}
The definitions of neutral and (long-$\beta\eta$) normal terms for the 
simply-typed lambda 
calculus, given in Figure~\ref{fig:neutral-and-normal-terms} below, are 
standard 
(\eg~\cite[Chapter 4]{Girard1989}). 
We define a family of judgements 
$\Gamma \vdash_M t : B$ and
$\Gamma \vdash_N t : B$  
characterising neutral and normal terms, respectively, 
by mutual induction. 

\begin{rules}
\unaryRule	{\faketext}
			{x_1 : A_1, \dots, x_n : A_n \vdash_M x_i : A_i}
			{var} 			

\unaryRule 	{\Gamma \vdash_M t : \prodop_n(A_1, \dots, A_n)}
			{\Gamma \vdash_M \pi_k (t) : A_k}
			{proj $(k=1,\dots, n)$}
			\:\:
\binaryRule	{\Gamma \vdash_M t : \exptype{A}{B}}
			{\Gamma \vdash_N u : A}
			{\Gamma \vdash_M \app{t}{u} : B}
			{app}
			
\unaryRule	{\Gamma \vdash_N t_i : A_i \quad (i = 1,\dots, n)}
			{\Gamma \vdash_N \seq{t_1, \dots, t_n} : 
				\prodop_n(A_1, \dots, A_n)}
			{tuple}
			\quad
\unaryRule	{\Gamma, x : A \vdash_N t : B}
			{\Gamma \vdash_N \lam{x}{t} : \exptype{A}{B}}
			{lam}

\unaryRule	{\Gamma \vdash_M t : B}
			{\Gamma \vdash_N t : B}
			{inc ($B$ a base type)}

\caption{Neutral terms and normal 
terms in the simply-typed lambda calculus\label{fig:neutral-and-normal-terms}}
\end{rules}

Crucially, the sets of neutral and normal terms are 
invariant under renaming, so for every type 
$A \in \allTypes\baseTypes$ one now obtains four presheaves
$\Con_{\allTypes\baseTypes} \to \Set$, defined at
$\Gamma \in \Con_{\allTypes\baseTypes}$ as follows:
\begin{equation} \label{eq:syntax-presheaves}
\begin{aligned}
\varTerms(\Gamma; A) &:= \yon \coerce{A} = 
	\slice{\{ x \st \Gamma \vdash x : A \}\:}{\: \aeq}  \\
\neutTerms(\Gamma; A) &:= \slice{\{ t \st \Gamma \vdash_M t : A \}\:}{\: \aeq} 
\\
\normTerms(\Gamma; A) &:= \slice{\{ t \st \Gamma \vdash_N t : A \}\:}{\: \aeq} 
\\
\langTerms(\Gamma;A) &:= \slice{\{ t \st \Gamma \vdash t : A \}\:}{\: \aeq}
\end{aligned}
\end{equation}

Each rule of Figure~\ref{fig:neutral-and-normal-terms} defines a 
morphism on these indexed families of presheaves. For the lambda abstraction 
case we employ the coproduct structure on $\Con_{\allTypes{\baseTypes}}$.

\begin{mylemma} \label{lem:typing-rules-are-nat-trans}
The rules of  
Figure~\ref{fig:neutral-and-normal-terms} 
give rise to natural transformations, as follows:
\begin{align*}
\mathsf{var}(-; A_i) : \varTerms(-; A_i) &\To 
	\neutTerms(-; A_i) \\
\mathsf{inc}\big(-; B \big) : 
	\neutTerms(-; B) 
	&\To \normTerms(-; B) &\text{($B$ a base type)} \\
\mathsf{proj}_k(-; \ind{A}) : 
	\neutTerms{\left(-; \prodop_n (A_1, \,\dots\, , A_n)\right)} &\To 
		\neutTerms(-; A_k) 
		&(k = 1, \,\dots\, , n) \\
\mathsf{app}(-; A, B) : 
	\neutTerms(-; \exptype{A}{B}) \times \normTerms(-; A) 
	&\To \neutTerms(-; B) \\
\mathsf{tuple}(-; \ind{A} ) : 
	\prodop_{i=1}^n \normTerms(-; A_i) 
	&\To \normTerms{\left(-; \prodop_n(A_1, \,\dots\, , A_n)\right)} \\
\mathsf{lam}(-; \exptype{A}{B} ) : 
	 \normTerms{\big(- + \coerce{A}; B\big)} &\To 
	 	\normTerms(-; \exptype{A}{B}) 
\end{align*}
\begin{proof}
The mappings are just the operations on terms.
In each case naturality follows from the definition of the 
meta-operation of 
capture-avoiding substitution, in particular the fact that substitution passes 
through the various constructors, and that it respects $\alpha$-equivalence. 
\end{proof}
\end{mylemma}

Returning to the development described by the 
diagram~(\ref{eq:h-to-interpret}), and noting 
that 
$\lanext{\prodext{h}{}}{\left(\prodext{h}{\coerce{A}}\right)} = 
	\catX{\left(\prodext{h}{(-)}, \prodext{h}{\coerce{A}} \right)} = 
	\catX{\left(h\sem{-}, h\sem{A} \right)}$
for every type $A$, one obtains the following glued 
objects in $\gl{\lanext{\prodext{h}{}}}$ for every 
$A \in \allTypes\baseTypes$:
\begin{equation} \label{eq:cat-glued-data}
\begin{aligned}
\underline{\varTerms}_A &:=
	\left(\varTerms(-; A), 
		\varTerms(-; A) \To
		\lanext{\prodext{h}{}}{\left(h\sem{A}\right)},  
		h\sem{A} 
		\right) = 
		\glued{\yon}(\coerce{A})
	 \\
\underline{\neutTerms}_A &:=
	\left(\neutTerms(-; A), 
		\neutTerms(-; A) \To
			\lanext{\prodext{h}{}}{\left(h\sem{A}\right)},  
		h\sem{A} \right) \\
\underline{\normTerms}_A &:=
	\left(\normTerms(-; A), 
		\normTerms(-; A) \To
	\lanext{\prodext{h}{}}{\left(h\sem{A}\right)},  
	h\sem{A} \right) \\
\underline{\langTerms}_A &:=
	\left(\langTerms(-; A), 
		\langTerms(-; A) \To
		\lanext{\prodext{h}{}}{\left(h\sem{A}\right)},  
		h\sem{A} 
		\right)				
\end{aligned} 
\end{equation}
In each case, the natural transformation is the canonical interpretation of 
$\stlc(\baseTypes)$-terms in $\catX$. 
Moreover, extending the natural 
transformations induced from the rules of 
Figure~\ref{fig:neutral-and-normal-terms} in a similar fashion, one obtains a 
morphism in $\gl{\lanext{\prodext{h}{}}}$ for each rule.

\paragraph*{Normalisation-by-evaluation.}
We paste together the various elements seen thus far. Since 
$\gl{\lanext{\prodext{h}{}}}$ is cartesian closed, one may consider the 
interpretation $B \mapsto \underline{\neutTerms}_B$ of base types in 
$\gl{\lanext{\prodext{h}{}}}$. This extends to an interpretation
$\overline{h}\sem{-} : \mathcal{F}(\allTypes\baseTypes) \to  
\gl{\lanext{\prodext{h}{}}}$. 
Write
$\overline{h}\sem{A} 
	:=
	\left( G_A, 
		\gamma_A, 
		h\sem{A}\right)$
and
$\overline{h}\sem{\Gamma \vdash t : A}
	:= \left( h'\sem{\Gamma \vdash t  : A},
				h\sem{\Gamma \vdash t : A} \right)$. 
Since the forgetful functor 
$\pi_{\mathrm{dom}} : \gl{\lanext{\prodext{h}{}}} \to \catX$ 
is strictly 
cartesian closed, the final component in each case is exactly the 
interpretation in 
$\catX$ extending $h$.

One then 
employs the cartesian closed structure of $\gl{\lanext{\prodext{h}{}}}$, and 
the 1-cells in $\gl{\lanext{\prodext{h}{}}}$ induced from the rules of
Figure~\ref{fig:neutral-and-normal-terms}, to 
inductively define \emph{quote} and \emph{unquote} as
$\allTypes{\baseTypes}$-indexed maps of the 
following type: 
\begin{align*}
\unquote_A : \underline{\neutTerms}_A \to \overline{h}\sem{A} \\
\quote_A : \overline{h}\sem{A} \to \underline{\normTerms}_A 
\end{align*}
For every 
$\stlc(\baseTypes)$-term $\Gamma \vdash t : A$ 
(where $\Gamma := (x_i : A_i)_{i=1, \,\dots\, , n}$), one thereby obtains the 
following commutative diagram in 
$\Psh{\op{\Con_{\allTypes\baseTypes}}}$, in which the unlabelled arrows 
are the canonical interpretations of terms inside $\catX$:
\begin{equation} \label{eq:1-cat-glueing-diagram}
\begin{tikzcd}[column sep = 2.6em, row sep = 2.1em]
\prod_{i=1}^n \neutTerms(-; A_i) 
\arrow{r}{\prod_{i=1}^n \unquote_{A_i}}
\arrow[bend right = 15]{dr} &
\prod_{i=1} G_{A_i}
\arrow[swap]{d}{\prod_{i=1}^n \gamma_{A_i}}
\arrow{r}{h'\sem{\Gamma \vdash t : A}}
\arrow{d} &
G_A
\arrow{r}{\quote_A} 
\arrow{dd}{\gamma_A}
\arrow{dd} &
\normTerms(-; A) 
\arrow[bend left = 8]{ddl} \\
\: &
\prod_{i=1}^n \catX\left(h\sem{-}, h\sem{{A_i}}\right)
\arrow[swap]{d}{\iso} &
\: &
\: \\
\: &
\catX\left(h\sem{-}, h\sem{\Gamma}\right)
\arrow[swap]{r}{h\sem{\Gamma \vdash t : A} \circ (-)} &
\catX\left(h\sem{-}, h\sem{{A}}\right) &
\:
\end{tikzcd}
\end{equation}
Chasing the $n$-ary variable-projection tuple 
$(\Gamma \vdash x_i : A_i)_{i=1, \,\dots\, , n}$ 
through this diagram, one obtains a normal term $\mathsf{nf}(t)$ for which the 
semantic interpretation $h\sem{\mathsf{nf}(t)}$ is equal to $h\sem{t}$. 
Moreover, for every type $A$ the projections
$\pi_{\mathrm{dom}}(\quote_A)$ and $\pi_{\mathrm{dom}}(\unquote_A)$ are both 
the identity. 
It follows that, for $\catX = \mathcal{F}(\allTypes\baseTypes)$ the 
syntactic model of 
$\stlc(\baseTypes)$, one obtains a normal form $\mathsf{nf}(t)$ for $t$
such that $t =_{\beta\eta} \mathsf{nf}(t)$. Hence, every 
$\stlc(\baseTypes)$-term has a long-$\beta\eta$ normal form, which can be 
explicitly calculated. This yields a normalisation algorithm. 

Our aim in what follows is to leverage as much of this proof as possible as we 
lift it to the bicategorical setting. We follow the strategy just outlined 
stage-by-stage, with the aim of building up a version 
of~(\ref{eq:1-cat-glueing-diagram}) in which each of the commuting shapes is 
filled by a witnessing 2-cell. Throughout we shall assume that 
$\baseTypes$ is a fixed set of base types. 

\section{Syntax as pseudofunctors}
\label{sec:syntax-as-pseudofunctors}

\paragraph*{The locally discrete 2-category of contexts.}

The notion of context in $\langCartClosed$ is the same as that in the 
simply-typed lambda calculus. We therefore require the same categorical 
structure on the category of contexts $\Con_{\allTypes{\baseTypes}}$, which we 
now wish to treat as a degenerate 2-category. Keeping track of such 
degeneracies will help identify instances where we can apply the 1-categorical 
theory.

\begin{mynotation} \quad
\begin{enumerate}
\item For $S$ a set, write $\oneDisc S$ for the 
discrete category with objects 
the elements of $S$. Similarly, write 
$\oneDisc f$ for the discrete functor $\oneDisc S \to \oneDisc S'$
induced by the set map $f : S \to S'$. 
\item 
\begin{enumerate}
\item For $\catC$ a category, write $\twoDisc{\catC}$ for the 
\Def{locally discrete} 2-category with objects those of $\catC$
and hom-categories
$(\twoDisc{\catC})(X, Y) := \oneDisc{\left(\catC(X, Y) \right)}$.

\item Write $\twoDisc{F}$ for the \Def{locally discrete} 2-functor
$\twoDisc{\catC} \to \twoDisc{\catD}$ 
induced from the functor 
$F : \catC \to \catD$
by setting 
$(\twoDisc{F})X := FX$ 
and 
$(\twoDisc{F})_{X,Y} := \oneDisc{\left(F_{X,Y}\right)}$. 

\item Write $\twoDisc{\mu}$ for the \Def{locally discrete} 2-natural 
transformation 
$\twoDisc{F} \To \twoDisc{G}$ induced from the natural transformation
$\mu : F \To G : \catC \to \catD$
by setting
$(\twoDisc{\mu})_C := \mu_C$ for every $C \in \catC$. \qedhere
\end{enumerate}
\end{enumerate}
\end{mynotation}

The $\oneDisc{(-)}$ and $\twoDisc{(-)}$ constructions will be our 
main technical tool for constructing 
(degenerate) 
bicategorical structure from 1-categorical data. The next lemma collects 
together some of their important properties. The proofs are not especially 
difficult, but stating all the details precisely requires some care. Since we employ the 
notation $\exp{-}{=}$ for exponentials in $\Hom(\baseCat, \Cat)$ we  
denote the usual categorical functor category by $\funCat{\catC}{\catD}$. 

\begin{mylemma} \label{lem:cat-to-2-cat}
Let $\catC$ and $\catD$ be 1-categories. Then:
\begin{enumerate}
\item \label{c:op-and-discrete} 
	$\op{(\twoDisc\catC)} = \twoDisc{(\op\catC)}$.
\item \label{c:Cat-to-Hom} There exists an isomorphism of 2-categories 
$\twoDisc{\big(\funCat{\catC}{\catD}\big)} 
		\iso \Hom(\twoDisc\catC, \twoDisc\catD)$.
\item \label{c:lifting-Yoneda} 
There exists an injective-on-objects, locally isomorphic 2-functor 
$\iota : \twoDisc{\funCat{\catC}{\Set}} \hookrightarrow 
	\Hom(\twoDisc\catC, \Cat)$, 
	which induces a commutative diagram
\begin{equation} \label{eq:Yon-as-iota-factored}
\begin{tikzcd}
\twoDisc{\big(\funCat{\catC}{\Set} \big)}
\arrow[hookrightarrow]{r}{\iota} &
\Hom(\twoDisc{\catC}, \Cat) \\
\twoDisc{\catC} 
\arrow{u}{\twoDisc{\yon}}
\arrow[swap]{ur}{\Yon} &
\:
\end{tikzcd}
\end{equation}
In particular, $\Yon(C) = (\twoDisc{\yon})C$ for all $C \in \catC$.
\item \label{c:1-cat-to-2-cat-products-and-exponentials} If $\catC$ is 
cartesian (resp. cartesian closed) as a 1-category, then 
$\twoDisc\catC$ has finite products (resp. is cartesian closed) as a 
2-category.
\item \label{c:presheaves-to-pseudofunctors} Let $P, Q : \catC \to \Set$. The 
exponential 
$\altexp{\iota{P}}{\iota{Q}}$ in $\Hom(\twoDisc\catC, \Cat)$ is given up 
to equivalence by
$\iota{\big(\funCat{\catC}{\Set}\big( \yon (-) \times P, Q\big)\big)}$, 
for 
$\yon : \catC \to \funCat{\catC}{\Set}$ the 1-categorical Yoneda embedding.  
\end{enumerate}
\begin{proof}
(\ref{c:op-and-discrete}) is immediate from the definitions.

For~(\ref{c:Cat-to-Hom}), consider the mapping
$\twoDisc{(-)} : \twoDisc{\big( \funCat{\catC}{\catD} \big)}
	\to \Hom(\twoDisc\catC, \twoDisc\catD)$
taking 
$F : \catC\to\catD$
to the locally discrete 2-functor $\twoDisc{F}$
and $\mu : F \to G$ 
to the locally discrete pseudonatural transformation
$\twoDisc\mu$. Since $\twoDisc{\big( \funCat{\catC}{\catD} \big)}$ is locally 
discrete, 
this extends canonically to a 2-functor. 

Now suppose that 
$F : \twoDisc\catC \to \twoDisc\catD$ is a pseudofunctor. By definition, this 
is a set map 
$F : ob(\twoDisc\catC) \to ob(\twoDisc\catD)$ 
with functors
$F_{X,Y} : (\twoDisc\catC)(X, Y) \to (\twoDisc\catD)(FX, FY)$ 
for every $X, Y \in \twoDisc\catC$.
Since every $(\twoDisc\catC)(X, Y)$ is a discrete category, 
every
$F_{X,Y}$ is discrete, and so $F = \twoDisc H$ for a unique
functor $H : \catC\to\catD$. So $\twoDisc{(-)}$ is bijective on objects.

Next fix functors 
$F, G : \catC \to \catD$ and consider the hom-category
$\Hom(\twoDisc\catC, \twoDisc\catD)(\twoDisc{F}, \twoDisc{G})$. A pseudonatural 
transformation
$(\natTrans, \natCell) : \twoDisc{F} \To \twoDisc{G}$ consists of  a family of 
1-cells
$\natTrans_X : FX \to GX \:\: (X \in \twoDisc\catC)$, together with a 2-cell
$\natCell_f : \natTrans_Y \circ Ff \To Gf \circ \natTrans_X$ in 
$\twoDisc\catD$ 
for every $f : X \to Y$ in $\twoDisc\catC$. Since $\twoDisc\catD$ is locally 
discrete, 
the only choice of such a 2-cell is the identity. So $(\natTrans, \natCell)$ 
is a 2-natural transformation, and is of the form 
$\twoDisc \mu$ for a unique natural transformation $\mu : F \To G$. Similarly, 
every modification 
$\modif : (\natTrans, \natCell) \to (\altNat, \altCell) : \twoDisc{F} \To 
\twoDisc{G}$
consists of a family of 2-cells, and must therefore be the identity. 
It follows that 
$\twoDisc{(-)}_{F,G} : \twoDisc{\big( \funCat{\catC}{\catD} \big)}(F,G)
	\to \Hom(\twoDisc\catC, \twoDisc\catD)(\twoDisc{F}, \twoDisc{G})$
is an isomorphism for every $F$ and $G$, as required.

For~(\ref{c:lifting-Yoneda}), we define $\iota$ 
by setting $\iota P$ to be the composite
$\catC \xra{P} \Set \xra{\oneDisc(-)} \Cat$, so that
$\iota P := \lambda C^{\catC} \bind \oneDisc (PC)$ 
and $(\iota \mu)_C := \oneDisc{\left(\mu_C\right)}$ for
every $\mu : P \To Q$ and $C \in \catC$. 
It is clear that $\iota$ is injective 
on objects. To see that
$\iota_{P,Q} :  \twoDisc{\big(\funCat{\catC}{\Set}\big)}(P, Q)
	\to \Hom(\twoDisc\catC, \Cat)(\iota{P}, \iota{Q})$
is an isomorphism for every $P$ and $Q$, one reasons as above: 
since $(\iota P)C$ is a discrete category for every $C \in \catC$, every 
pseudonatural transformation $\iota P \To \iota Q$ must be of the form 
$\iota(\mu)$ for a unique natural transformation $\mu : P \To Q$, 
and there can be no non-identity modifications between such transformations.

To relate the 1-categorical and bicategorical Yoneda embeddings, one notes that
\begin{align*}
(\iota \circ \twoDisc{\yon})(C) 
	&= \iota \big( \catC(C, -) \big) \\
	&= \lambda X^{\catC} \bind \oneDisc{\left(  \catC(C, X) \right)} \\
	&= \lambda X^{\catC} \bind (\twoDisc{\catC})(C, X) \\
	&= \Yon C
\end{align*}
as claimed.

For~(\ref{c:1-cat-to-2-cat-products-and-exponentials}), one simply observes 
that the natural isomorphisms 
$\catC{\left(X, \prod_{i=1}^n A_i\right)} \iso \prod_{i=1}^n \catC(X, A_i)$ 
immediately provide 2-natural isomorphisms of hom-categories
$
(\twoDisc\catC){\left( X, \prod_{i=1}^n A_i \right)} 
\iso 
\prod_{i=1}^n (\twoDisc\catC)(X, A_i)
$, and similarly for exponentials. 

For~(\ref{c:presheaves-to-pseudofunctors}), recall from 
Theorem~\ref{thm:2-presheaves-cartesian-closed} that for pseudofunctors 
$G, H : \twoDisc{\catC} \to \Cat$, the exponential
$\altexp{G}{H}$ may be given by the pseudofunctor
$\Hom(\twoDisc{\catC}, \Cat){\left( \Yon(-) \times G, H\right)} 
: \twoDisc{\catC} \to \Cat$. Next observe that the embedding $\iota$ 
of~(\ref{c:lifting-Yoneda}) preserves products:
\begin{align*}
\big(\iota (P \times Q) \big)C 
	&= \oneDisc{\big((P \times Q)(C)\big)}  \\
	&= \oneDisc{(PC \times QC)} \\ 
	&= \oneDisc{(PC)} \times \oneDisc{(QC)} \\
	&= (\oneDisc{P} \times \oneDisc{Q})C \\
	&= \big(\iota(P) \times \iota(Q)\big)C 
\end{align*}
Hence:
\begin{align*}
\Hom(\twoDisc{\catC}, \Cat)
	{\left( \Yon X \times \iota{P}, \iota{Q}\right)} 
	&= \Hom(\twoDisc{\catC}, \Cat){\left( 
			(\iota \circ \twoDisc{\yon})X \times \twoDisc{P}, \twoDisc{Q}
		\right)} \qquad &\text{ by diagram}~(\ref{eq:Yon-as-iota-factored}) \\
	&=  \Hom(\twoDisc{\catC}, \Cat){\left( 
				\iota (\yon X) \times \iota(P), \iota(Q)
			\right)} \\
	&= \Hom(\twoDisc{\catC}, \Cat){\left( 
					\iota (\yon X \times P), \iota(Q)
				\right)} \\
	&\iso \left(\twoDisc{\funCat{\catC}{\Set}}\right)(\yon X \times P, Q)
		&\text{ by }~(\ref{c:lifting-Yoneda}) \\
	&= \oneDisc{{\left( \funCat{\catC}{\Set}(\yon X \times P, Q) \right)}} 
		&\text{by definition of } \twoDisc{(-)}
\end{align*}
completing the proof.
\end{proof}
\end{mylemma}

The preceding lemma provides a framework for treating the category of contexts 
$\Con_{\allTypes{\baseTypes}}$ as a 2-category. Next we show how to 
extend from an interpretation of (base) types to an interpretation of all 
contexts, that is, to an fp-pseudofunctor out of 
$\twoDisc{\op{\Con_{\allTypes\baseTypes}}}$. In the categorical setting, one 
merely uses the fact that $\op{\Con_{\allTypes\baseTypes}}$ is the free strict 
cartesian category on $\allTypes\baseTypes$. The pseudo nature of bicategorical 
products and exponentials entails a little more work, but the construction is 
essentially the same.

Note that any interpretation 
$s : \baseTypes \to \bicatX$
of base types in a cc-bicategory $\ccBicat\bicatX$ extends canonically
to an interpretation
$\allTypes\baseTypes \to \bicatX$
by the cartesian closed structure, which we also denote by $s$.
%

\begin{mylemma} \label{lem:contexts-pseudo-free-prop}
For any set of base types $\baseTypes$, 
cc-bicategory $\ccBicat{\bicatX}$, and
set map
$s : \baseTypes \to \bicatX$, 
there exists an fp-pseudofunctor 
$\prodext{s}{} : \twoDisc{\op{\Con_{\allTypes\baseTypes}}} \to \bicatX$ 
making the following diagram commute:
\begin{td}
\: &
\twoDisc{\op{\Con_{\allTypes\baseTypes}}}
\arrow{dr}{\prodext{s}{}} &
\: \\

\allTypes\baseTypes
\arrow{ur}{\coerce{-}}
\arrow{rr}{} &
\: &
\bicatX \\

\baseTypes
\arrow[hookrightarrow]{u}
\arrow[swap]{urr}{s} &
\: &
\:
\end{td}
\begin{proof}
We define $\prodext{s}{}$ on types by 
$\prodext{s}{A} := sA$ 
and extend to contexts in the usual manner: 
$\prodext{s}{\big((x_i : A_i)_{i=1, \,\dots\, , n}\big)} 
		:= \prod_{i=1}^n \prodext{s}{A_i}$ 
and
$\prodext{s}{(\diamond)} := \prodop_0 ()$. In particular, for a unary context 
$(x : A)$ we define $\prodext{s}{(x:A)} = sA$, so that
$\prodext{s}{\coerce{A}} = sA$. 

The action on 1-cells is the following. For contexts
$\Gamma := (x_i : A_i)_{i=1, \,\dots\, , n}$ and
$\Delta := ({y_j : B_j})_{j=1, \,\dots\, , m}$ 
and a context renaming $r : \Gamma \to \Delta$, we define 
$\prodext{s}{r} : \prodop_{j=1}^m \prodext{s}{B_j} \to \prodop_{i=1}^n 
\prodext{s}{A_i}$ to be
$\seqlr{\pi_{r(1)}, \,\dots\, , \pi_{r(n)}}$, where we write
$r(i)$ to indicate the index of $r(x_i)$ within $(y_1, \,\dots\, , y_m)$. The 
action 
on 2-cells is trivial since 
$\twoDisc{\op{\Con_{\allTypes\baseTypes}}}$ 
is locally discrete.

For the 2-cell 
$\psi^{\prodext{s}{}}_\Gamma : 
	\Id_{\prodext{s}{\Gamma}} 
		\To 
	\prodext{s}{(\Id_\Gamma)}$ 
we take
\[
\widehat{\etaTimes{}}_{\Id_{\prodext{s}{\Gamma}}} :=
	\Id_{\prodext{s}{\Gamma}}
		\XRA{\etaTimes{\Id_{\prodext{s}{\Gamma}}}}
	\seqlr{\pi_1 \circ \Id_{\prodext{s}{\Gamma}}, \,\dots\, , 
		\pi_n \circ \Id_{\prodext{s}{\Gamma}}}
		\XRA{\iso}	
	\seqlr{\pi_1, \,\dots\, , \pi_n}
\]
For a composable pair of context renamings 
$\Sigma \xra{r} \Gamma \xra{r'} \Delta$,
we define $\phi^{\prodext{s}{}}_{r', r}$ to be the composite
\begin{td}
\seqlr{\pi_{r(1)}, \,\dots\, , \pi_{r(n)}} \circ 
	\seqlr{\pi_{r'(1)}, \,\dots\, , \pi_{r'(m)}} 
\arrow[swap]{d}{\postName} 
\arrow[bend left = 10]{dr}{\phi^{\prodext{s}{}}_{r', r}} &
\: \\

\seqlr{\pi_{r(1)} \circ \seq{\pi_{r'(\bullet)}}, \,\dots\, , 
	\pi_{r(n)} \circ \seq{\pi_{r'(\bullet)}}} 
\arrow[swap]{r}[yshift=-2mm]
	{\seq{\epsilonTimesInd{r(1)}{}, \,\dots\, , \epsilonTimesInd{r(n)}{}}}  &
\seqlr{\pi_{r'r(1)}, \,\dots\, , \pi_{r'r(n)}}
\end{td}
The three axioms to check are diagram chases using the product 
structure, along with the properties of
Lemma~\ref{lem:PseudoproductCanonical2CellsLaws}. For the associativity 
law one uses naturality and the commutativity 
of the following diagram, in which we abbreviate 
$\seq{\pi_{r(1)}, \,\dots\, , \pi_{r(n)}}$ by $\seq{\pi_r}$:
\begin{td}
\seq{\pi_{r}} \circ \seq{\pi_{r'}} \circ \seq{\pi_{r''}}
\arrow[bend left = 10]{dr}{\postName}
\arrow[swap]{d}{\postName \circ \seq{\pi_{r''}}} &
\: \\

\seqlr{\pi_r \circ \seq{\pi_{r'}}} \circ \seq{\pi_{r''}} 
\arrow[swap]{r}{\postName} &
\seqlr{\pi_r \circ \seq{\pi_{r'}} \circ \seq{\pi_{r''}}} 
\end{td}
For the left and right unit laws, one respectively uses the diagrams on the 
left and right below:
\begin{center}
\begin{tikzcd}
\Id_{\prodext{s}{\Sigma}} \circ \seq{\pi_r} 
\arrow[swap]{d}{\etaTimes{\Id} \circ \seq{\pi_r}}
\arrow[bend left = 10]{dr}{\etaTimes{\Id \circ \seq{\pi_r}}} &
\: \\
\seqlr{\ind{\pi} \circ \Id_{\prodext{s}{\Sigma}}}
	\circ \seq{\pi_r}
\arrow[swap]{r}{\postName} &
\seqlr{\ind{\pi} \circ \Id_{\prodext{s}{\Sigma}} \circ \seq{\pi_r}}
\end{tikzcd}
\:\:
\begin{tikzcd}
\seq{\pi_r} \circ \Id_{\prodext{s}{\Gamma}} 
\arrow[swap]{d}{\postName}
\arrow[bend left = 10]{dr}{\iso} &
\: \\
\seqlr{\pi_{r(1)} \circ \Id_{\prodext{s}{\Gamma}}, \,\dots\, , \pi_{r(n)}
	\circ \Id_{\prodext{s}{\Gamma}}} 
\arrow[swap]{r}{\iso} &
\seq{\pi_r}
\end{tikzcd}
\end{center}
It remains to show that $\prodext{s}{}$ preserves products. For $n$ contexts
$\Gamma_1, \,\dots\, , \Gamma_n \:\: (n \in \Nat)$ of the form
$\Gamma_i := {(x_j^{(i)} : A_j^{(i)})}_{j=1, \,\dots\, , \len{\Gamma_i}}$, note 
that
\begin{align*}
\prodext{s}{\left(\prodop_{i=1}^n \Gamma_i\right)} 
 &= \prodext{s}{(\Gamma_1 \concat \cdots \concat \Gamma_n)} 
 =  \prodop_{\substack{j_i=1,\dots, \len{\Gamma_i} \\ i=1, \,\dots\, ,n}} 
 s(A_i) \\
\prodop_{i=1}^n \prodext{s}(\Gamma_i) &= 
	\prodop_{i=1}^n \prodop_{j=1}^{\len{\Gamma_i}} s(A_j^{(i)})
\end{align*}
and that 
$\prodext{s}{(\pi_k)} = 
	\prodext{s}
		{(\Gamma_k \hookrightarrow \Gamma_1 \concat \cdots \concat \Gamma_n)}$ 
is the 1-cell 
$\seq{\pi_{1+ \sum_{i=1}^{k-1} \len{\Gamma_i}}, \,\dots\, , 
	\pi_{\sum_{i=1}^k \len{\Gamma_i}} }$. 
One therefore obtains the required equivalence
$\prod_{i=1}^n \prod_{j=1}^{\len{\Gamma_i}} s(A_i^{(j)}) \simeq
	 \prod_{\substack{j=1,\dots, \len{\Gamma_i} \\ i=1, \,\dots\, ,n}} 
	 s(A_i^{(j)})$ 
by taking $\prodPres_{\ind{\Gamma}}$ to be the 1-cell
$\prod_{i=1}^n \prod_{j=1}^{\len{\Gamma_i}} s(A_j^{(i)}) \to
\prod_{\substack{j=1,\dots, \len{\Gamma_i} \\ i=1, \,\dots\, ,n}} s(A_j^{(i)})$
given by
\begin{equation} \label{eq:prodext-prod-pres}
\seqlr{\pi_1 \circ \pi_1, \,\dots\, , \pi_{\len{\Gamma_1}} \circ \pi_1, 
\,\dots\, , 
		\pi_1 \circ \pi_k, \,\dots\, , \pi_{\len{\Gamma_k}} \circ \pi_k, 
		\dots, 
		\pi_1 \circ \pi_n, \,\dots\, , \pi_{\len{\Gamma_n}} \circ \pi_n}
\end{equation}
This defines an equivalence with witnessing 2-cells defined by the 
commutativity of the following two diagrams:
\begin{td}[column sep = 5em]
\seq{\pi_1 \circ \pi_1, \,\dots\, , 
	\pi_{\len{\Gamma_n}} \circ \pi_n} \circ 
\seq{\prodext{s}{\ind{\pi}}}
\arrow{r}
\arrow[swap]{d}{\postName} &
\Id_{\prodext{s}{(\prod_i \Gamma_i)}}  \\

\seqlr{\dots, \pi_1 \circ \pi_k \circ 
		\seq{\prodext{s}{\ind{\pi}}}, \,\dots\, , 
	\pi_{\len{\Gamma_k}} \circ \pi_k \circ 
		\seq{\prodext{s}{\ind{\pi}}}, \dots} 
\arrow[swap]{d}
	{\seq{\dots, \pi_1 \circ \epsilonTimesInd{k}{}, \,\dots\, , 
		\pi_{\len{\Gamma_k}} \circ \epsilonTimesInd{k}{}, \dots}} &
\seq{\pi_1, \,\dots\, , \pi_{\sum_{i=1}^n\sum_{j=1}^{\len{\Gamma_i}} j}}
\arrow[swap]{u}
	{\widehat{\etaTimes{}}_{\Id_{\prodext{s}{(\prod_i \Gamma_i)}}}^{-1}} \\

\seqlr{\dots, 
	\pi_j \circ 
		\seq{\pi_{1+ \sum_{i=1}^{k-1} \len{\Gamma_i}}, \,\dots\, , 
			\pi_{\sum_{i=1}^k \len{\Gamma_i}} }, \dots} 
\arrow[swap]{r}{
	\seq{\dots, \epsilonTimesInd{j}{}, \dots}} &
\seqlr{\dots, \pi_{j+ \sum_{i=1}^{k-1} \len{\Gamma_i}}, \dots}
\arrow[swap, equals]{u} \\

\seqlr{\prodext{s}{(\ind{\pi})}} \circ \seq{\pi_1 \circ \pi_1, \,\dots\, , 
	\pi_{\len{\Gamma_n}} \circ \pi_n} 
\arrow{r}
\arrow[swap]{d}{\postName} &
\Id_{(\prod_i \prodext{s}{\Gamma_i})} \\

\seqlr{\dots, \prodext{s}{(\pi_k)} \circ \seq{\pi_1 \circ \pi_1, \,\dots\, , 
	\pi_{\len{\Gamma_n}} \circ \pi_n}, \dots} 
\arrow[swap]{d}{\iso} &
\seq{\pi_1, \,\dots\, , \pi_n} 
\arrow[swap]{u}
	{\widehat{\etaTimes{}}_{\Id_{(\prod_i \prodext{s}{\Gamma_i})}}^{-1}} \\

\seqlr{\dots, \seq{\pi_1, \,\dots\, , \pi_{\len{\Gamma_k}}} \circ \pi_k, \dots}
\arrow[swap]{r}[yshift=-2mm]
	{\seq{\widehat{\etaTimes{}}_{\Id_{(\prodext{s}{\Gamma_n})}}^{-1} 
		\circ \ind{\pi}}} &
\seqlr{\Id_{(\prodext{s}{\Gamma_1})} \circ  \pi_1, \,\dots\, , 
	\Id_{(\prodext{s}{\Gamma_n})} \circ \pi_n}
\arrow[swap]{u}{\iso}
\end{td}
The downwards arrow labelled $\iso$ is the $n$-ary tupling of 
\begin{td}
\seq{\pi_{1+ \sum_{i=1}^{k-1} \len{\Gamma_i}}, \,\dots\, , 
	\pi_{\sum_{i=1}^k \len{\Gamma_i}} } \circ 
\seq{\pi_1 \circ \pi_1, \,\dots\, , 
	\pi_{\len{\Gamma_n}} \circ \pi_n}
\arrow{r}
\arrow[swap]{d}{\postName} &
\seq{\pi_1, \,\dots\, , \pi_{\len{\Gamma_k}}} \circ \pi_k \\

\seq{\dots, \pi_{j+ \sum_{i=1}^{k-1} \len{\Gamma_i}} \circ 
		\seq{\pi_1 \circ \pi_1, \,\dots\, , 
			\pi_{\len{\Gamma_n}} \circ \pi_n}, \dots}_{j=1, 
			\dots,\len{\Gamma_k}} 
\arrow[swap]{r}[yshift=-4mm]
	{\seq{\dots, \epsilonTimesInd{j+ \sum_{i=1}^{k-1}}{}, \dots}} &
\seqlr{\dots, \pi_j \circ \pi_k, \dots}_{j=1, \,\dots\, ,\len{\Gamma_k}}
\arrow[swap]{u}{\postName^{-1}}
\end{td}
for $k=1, \,\dots\, , n$. Hence $\prodext{s}{}$ is an fp-pseudofunctor, as 
claimed.
\end{proof}
\end{mylemma}

\begin{myremark} \label{rem:prod-pres-special-case}
We shall need the following special case of the fact that the pseudofunctor 
$\prodext{s}{}$ preserves products. For a 
context $\Gamma = (x_i : A_i)_{i=1, \,\dots\, ,n}$ and type $A$, the 
1-cell~(\ref{eq:prodext-prod-pres}) becomes simply
$
\seq{\pi_1 \circ \pi_1, \,\dots\, , \pi_n \circ \pi_1, \pi_2} : 
\prodext{s}{\Gamma} \times \prodext{s}{\coerce{A}} \to
\prodext{s}{\left(\Gamma \concat \coerce{A}\right)} 
$.
\end{myremark}


One also obtains the following version of 
Proposition~\ref{propn:weak-universal-property} by taking the 
\emph{context extension} product structure of the syntactic model 
instead of 
the type-theoretic product structure (recall 
Section~\ref{sec:products-from-context-extension}). 

\begin{mypropn}
For any $\langCartClosed$-signature $\sig$, 
cc-bicategory $\ccBicat{\bicatX}$,
and $\langCartClosed$-signature homomorphism
$s : \sig \to \bicatX$, 
there exists a cc-pseudofunctor 
$s\sem{-} : \syncloneAtClosed{\sig} \to \bicatX$ 
with respect to the context extension product structure,
such that $s\sem{-} \circ \inc = s$, for 
$\inc : \sig \hookrightarrow \syncloneAtClosed{\sig}$ the inclusion.
\begin{proof}
Define $s\sem{-}$ as in 
Proposition~\ref{propn:weak-universal-property}, except that for preservation 
of 
products one takes $\prodPres$ as in the 
preceding lemma. Preservation of exponentials then takes the following form. 
For $\Gamma := (x_i : A_i)_{i=1, \,\dots\, , n}$ and
$\Delta := (y_j : B_j)_{j=1, \,\dots\, , m}$, the evaluation map
is the $m$-tuple with components 
\[
f : \exptype{\prodop_n \ind{A}}{\prodop_m \ind{B}}, 
	x_1 : A_1, \,\dots\, , x_n : A_n \vdash 
\hpi{j}{\heval{f, \pair{x_1, \,\dots\, , x_n}}} : B_j
\]
for $j=1, \,\dots\, , m$. One then obtains the following chain of natural 
isomorphisms: 
\begin{align*}
s\sem{\eval_{\Gamma, \Delta}} 
	&\circ \prodPres_{\scriptsizeexpobj{\Gamma}{\Delta}, \Gamma} \\
	&= \seqlr{\ind{\pi} 
		\circ \eval_{s\sem{\prod_n \ind{A}}, s\sem{\prod_m \ind{B}} } 
		\circ \seqlr{\pi_1, \seq{\pi_2, \,\dots\, , \pi_{n+1}}}} \circ
		\seq{\pi_1, \pi_1 \circ \pi_2, \,\dots\, , \pi_n \circ \pi_2} \\
	&\iso \seqlr{
		\ind{\pi} 
		\circ  \eval_{s\sem{\prod_n \ind{A}}, s\sem{\prod_m \ind{B}}}
		\circ \seqlr{\pi_1, \seq{\pi_1 \circ \pi_2, \,\dots\, , \pi_n \circ 
		\pi_2}}} 
		\\
	&\iso  \seqlr{
			\ind{\pi} 
			\circ  \eval_{s\sem{\prod_n \ind{A}}, s\sem{\prod_m \ind{B}}}
			\circ \seqlr{\pi_1, \seq{\pi_1, \,\dots\, , \pi_n} \circ \pi_2}} \\
	&\iso \seqlr{
				\ind{\pi} 
				\circ  \eval_{s\sem{\prod_n \ind{A}}, s\sem{\prod_m \ind{B}}}
				\circ \seqlr{\pi_1, \pi_2}} \\
	&\iso \seqlr{\ind{\pi} 
					\circ  \eval_{s\sem{\prod_n \ind{A}}, 
						s\sem{\prod_m \ind{B}}}} \\
	&\iso \seq{\pi_1, \,\dots\, , \pi_m} \circ \eval_{s\sem{\prod_n \ind{A}}, 
							s\sem{\prod_m \ind{B}}} \\
	&\iso \eval_{s\sem{\prod_n \ind{A}}, 
							s\sem{\prod_m \ind{B}}}
\end{align*}
It follows that
$\evBar_{\Gamma, \Delta} 
	= \lambda(s\sem{\eval_{\Gamma, \Delta}} ) 
	\iso \lambda {\left(\eval_{s\sem{\prod_n \ind{A}}, 
								s\sem{\prod_m \ind{B}}}\right)} 
	\iso \id_{s\sem{\scriptsizeexpobj{\Gamma}{\Delta}}}$, so
$s\sem{-}$ preserves exponentials. 
\end{proof}
\end{mypropn}

While the interpretation of Proposition~\ref{propn:weak-universal-property} is 
useful for proving uniqueness properties, the interpretation of the preceding 
proposition is the natural choice when working with the {(2-)category} of 
contexts. Of course, the two pseudofunctors are canonically equivalent. 
Throughout this chapter, we shall work with the version just defined.

For any interpretation of base types 
$s : \baseTypes \to \bicatX$
in a cc-bicategory $\ccBicat\bicatX$, one therefore obtains the following 
diagram lifting~(\ref{eq:h-to-interpret}) to the bicategorical setting:
\begin{equation*}
\begin{tikzcd} 
\: &
\syncloneAtClosed{\allTypes\baseTypes}
\arrow[]{dr}{s\sem{-}} &
\: 
\\
\twoDisc{\op{\Con_{\allTypes\baseTypes}}}
\arrow[]{ur}{\prodext{\iota}{}}
\arrow[swap]{rr}{\prodext{s}{}} &
\: &
\bicatX &
\: \\
\allTypes{\baseTypes}
\arrow[hookrightarrow]{u}{\coerce{-}}
\arrow[swap]{urr}{} &
\: \\
\baseTypes
\arrow[hookrightarrow]{u}
\arrow[swap]{uurr}{s}
\end{tikzcd}
\end{equation*}
Note in particular that, just 
as in the 1-categorical case, the equality 
$s\sem{\Gamma} = \prodext{s}{\Gamma}$ holds for every context $\Gamma$. 

\paragraph*{Syntactic presheaves for $\langCartClosed$.}

Lemma~\ref{lem:contexts-pseudo-free-prop} provides a way to interpret contexts 
whenever one has an interpretation of base types, 
while Lemma~\ref{lem:cat-to-2-cat} guarantees that, in order to interpret the 
syntax of $\langCartClosed$ as a pseudofunctor 
$\twoDisc{\Con_{\allTypes{\baseTypes}}} \to \Cat$, 
it suffices to a define a presheaf
${\Con_{\allTypes{\baseTypes}}} \to \Set$
on the underlying category. There 
remains the question of what it means to 
be a neutral or normal term in $\langCartClosed$. The answer is provided by the 
embedding of $\stlc$ into $\langCartClosed$ constructed in
Section~\ref{sec:STLC-vs-pseudoSTLC}. Thus, for every 
$A \in 
{\allTypes{\baseTypes}}$ 
we define four presheaves 
$\dvarTerms(-; A), 
\dneutTerms(-; A), 
\dnormTerms(-; A), 
\dlangTerms(-; A) :
 \Con_{\allTypes\baseTypes} \to \Set$ by setting
\begin{equation} \label{eq:langCartClosed-syntactic-presheaves}
\begin{aligned}
\dvarTerms(\Gamma; A) := \{ \into{t} \st t \in \varTerms(\Gamma; A) \} \\
\dneutTerms(\Gamma; A) := \{ \into{t} \st t \in \neutTerms(\Gamma; A) \} \\
\dnormTerms(\Gamma; A) := \{ \into{t} \st t \in \normTerms(\Gamma; A) \} \\
\dlangTerms(\Gamma; A) := \{ \into{t} \st t \in \langTerms(\Gamma; A) \} 
\end{aligned}
\end{equation}
where $\into{-}$ is defined in Construction~\ref{constr:ccc:into-terms}
on page~\pageref{constr:ccc:into-terms}
and the presheaves $\varTerms(; A), \neutTerms(-; A), \normTerms(-;A)$ and 
$\langTerms(-;A)$ are defined 
in~(\ref{eq:syntax-presheaves}) on~page~\pageref{eq:syntax-presheaves}. 
Since $\into{-}$ respects $\alpha$-equivalence 
(Lemma~\ref{lem:into-respects-typing}), these definitions are well-defined on 
$\alpha$-equivalence classes. To see that these definitions are invariant under 
variable renamings, recall from Construction~\ref{constr:def-of-cont} that the 
following rule is admissible in $\langCartClosed$:
\begin{center}
\binaryRule{\Gamma \vdash \into{t} : B}{r : \Gamma \to \Delta}{\Delta \vdash 
\cont{t}{r} : \rewrite{\hcomp{\into{t}}{x_i \mapsto 
r(x_i)}}{\into{t[r(x_i)/x_i]}} : B}{} \vspace{-\treeskip}
\end{center}
\newpage
Since a rewrite $\tau : t \To t'$ is typeable in context $\Gamma$ only if 
both $t$ and $t'$ are also typeable in $\Gamma$, it follows that the following 
rule is admissible:
\begin{center}
\binaryRule
	{\Gamma \vdash \into{t} : B}
	{r : \Gamma \to \Delta}
	{\Delta \vdash \into{t[r(x_i)/x_i]} 
		: B}
	{} \vspace{-\treeskip}
\end{center}
Since the presheaves~(\ref{eq:syntax-presheaves}) are invariant under 
renamings, it follows that those 
of~(\ref{eq:langCartClosed-syntactic-presheaves}) are too, as required. 
 
The functorial action is the unique choice
such that the following diagram
commutes, where 
$\terms(-; A) \in \{ \varTerms(-; A), \neutTerms(-; A), \normTerms(-; A) \}$ 
and 
$\dterms(-; A)$ denotes the 
image of $\terms(-; A)$ under $\into{-}$:
\begin{equation} \label{eq:extended-presheaf-nat-square}
\begin{tikzcd}[column sep = 3em]
\terms(\Gamma; A) 
\arrow[swap]{d}{\into{-}^\Gamma_A}
\arrow{r}{\terms(r; A)} &
\terms(\Delta; A) 
\arrow{d}{\into{-}^\Delta_A} \\
\dterms(\Gamma; A)
\arrow[swap]{r}{\dterms(r;A)} &
\dterms(\Delta; A)
\end{tikzcd}
\end{equation}
Explicitly, for a context renaming 
$r : \Gamma \to \Delta$ we define 
$\dterms(-; A)(r)(\into{t}^\Gamma_A) := \into{t[r(x_i)/ x_i]}^\Delta_A$.

This formulation is particularly 
convenient as it allows one to make use of 
standard facts about the simply-typed lambda calculus. Moreover, we can employ 
many of the details of Fiore's proof via the following observation.

\begin{mylemma} \label{lem:into-a-natural-isomorphism}
For any 
type 
$A \in \allTypes{\baseTypes}$, let
$\terms(-; A) \in \{ \varTerms(-; A), \neutTerms(-; A), 
	\normTerms(-; A), \langTerms(-; A) \}$ 
and let 
$\dterms(-; A) \in 
	\{ \dvarTerms(-; A),  \dneutTerms(-; A), 
		\dnormTerms(-; A), \dlangTerms(-; A)\}$ denote the 
image of $\terms_A$ under $\into{-}$.  Then the mappings 
$\into{-}_A^{(=)} : \terms_A \To \dterms_A$ form a natural 
isomorphism. 
\begin{proof}
Since $\into{-}_A^{(=)}$ respects the typings, it is clear from the definition 
that it is an injection, hence  
a bijection onto its image. Naturality is 
exactly~(\ref{eq:extended-presheaf-nat-square}). 
\end{proof}
\end{mylemma}

For example, one may immediately extend the natural 
transformations of
Lemma~\ref{lem:typing-rules-are-nat-trans} 
to $\langCartClosed$. One therefore obtains the following natural 
transformations:
\begin{equation} \label{eq:typing-rules-as-pseudonat-trans}
\begin{aligned}
\mathsf{var}(-; A_i) : \dvarTerms(-; A_i) &\To 
	\dneutTerms(-; A_i) \\
\mathsf{inc}(-; B) : \dneutTerms(-; B) &\To 
	\dnormTerms(-; B) &\text{($B$ a base type)}  \\
\mathsf{proj}_k(-; \ind{A}) : 
	\dneutTerms{\left(-; \prodop_n (A_1, \,\dots\, , A_n)\right)} &\To 
		\dneutTerms(-; A_k) \qquad &(k = 1, \,\dots\, , n) \\
\mathsf{app}(-; A, B) : 
	\dneutTerms(-; \exptype{A}{B}) \times \dnormTerms(-; A) 
	&\To \dneutTerms(-; B) \\
\mathsf{tuple}(-; \ind{A} ) : 
	\prodop_{i=1}^n \dnormTerms(-; A_i) 
	&\To \dnormTerms{\left(-; \prodop_n(A_1, \,\dots\, , A_n)\right)} \\
\mathsf{lam}(-; A, B ) : 
	 \dnormTerms{\big(- + \coerce{A}; B\big)} &\To 
	 	\dnormTerms(-; \exptype{A}{B}) 
\end{aligned}
\end{equation}
Explicitly, the action on terms is the following:
\begin{align*}
x_k &\mapsto x_k \\
\into{t} &\mapsto \into{t} \\
\into{t} &\mapsto \into{\pi_k(t)} = \hpi{k}{\into{t}}
		&(k = 1, \,\dots\, , n) \\
\left(\into{t}, \into{u}\right) &\mapsto 
	\into{\app{t}{u}} = \heval{\into{t}, \into{u}} \\
\left(\into{t_1}, \,\dots\, , \into{t_n}\right) &\mapsto 
\into{\seq{t_1, \,\dots\, , t_n}} = 
		\pair{\into{t_1}, \,\dots\, , \into{t_n}} \\
\into{t} &\mapsto \into{\lam{x}{t}} = \lam{x}{\into{t}} 
\end{align*}

The presheaves~(\ref{eq:langCartClosed-syntactic-presheaves}) and natural 
transformations~(\ref{eq:typing-rules-as-pseudonat-trans})---viewed as 
locally discrete pseudofunctors and locally discrete pseudonatural 
transformations---describe the syntax of $\langCartClosed$ within
$\Hom(\twoDisc{\op{\Con_{\allTypes\baseTypes}}}, \Cat)$. As we saw in 
Chapter~\ref{chap:calculations}, this bicategory shares many of the important 
features of the presheaf category 
$\Psh{\op{\Con_{\allTypes\baseTypes}}}$. Our next task, therefore, is to 
construct the bicategorical correlate to the category of intensional Kripke 
relations.

\subsection{Bicategorical intensional Kripke relations}

\paragraph*{The relative hom-pseudofunctor.}
We start by constructing the pseudo correlate of the relative hom-functor and 
establishing its key properties. Precisely, we show that 
diagram~(\ref{eq:glueing-left-kan-diagram}) on 
page~\pageref{eq:glueing-left-kan-diagram} lifts to the bicategorical setting, 
and that the relative hom-pseudofunctor preserves bilimits. 

The construction is the natural bicategorification of 
Definition~\ref{def:relative-hom-functor}.

\begin{myconstr}
For any pseudofunctor $\glueFun  : \baseCat \to \bicatX$ one obtains a 
\Def{relative hom-pseudofunctor}
$\lanext{\glueFun } : \bicatX \to \Hom(\op\baseCat, \Cat)$ as follows. 

On objects, we set $\lanext{\glueFun }X := 
\bicatX{\left(\glueFun (-), X\right)}$. 
On morphisms, we define a pseudonatural transformation
$\lanext{\glueFun }f : \lanext{\glueFun }X \To \lanext{\glueFun }X'$
for every $f : X \to X'$ in $\bicatX$. The 1-cell components are 
\[
(\lanext{\glueFun }f)_B := 
	\bicatX(\glueFun B, X) \xra{f \circ (-)} \bicatX(\glueFun B, X')
\]
and for $g : B' \to B$ in $\baseCat$ the witnessing 2-cell
$\cellOf{(\lanext{\glueFun }f)}_g$
filling
\begin{td}[column sep = 4em]
\bicatX(\glueFun B, X) 
\arrow[phantom]{dr}[description]{\twocell{\cellOf{\left(\lanext{\glueFun 
}f\right)}_g}}
\arrow{r}{(\lanext{\glueFun }X)(g)}
\arrow[swap]{d}{f \circ (-)} &
\bicatX(\glueFun B', X) 
\arrow{d}{f \circ (-)} \\

\bicatX(\glueFun B, X') 
\arrow[swap]{r}{(\lanext{\glueFun }X')(g)} &
\bicatX(\glueFun B', X')
\end{td}
is the structural isomorphism
$\lambda h^{\bicatX(\glueFun B,X)} \bind \a^{-1}_{f, h, \glueFun g}$. Finally, 
for a 2-cell
$\tau : f \To f'$ in $\bicatX$, 
we define a modification
$\lanext{\glueFun }f \to \lanext{\glueFun }f'$
by setting 
$\lanext{\glueFun }\tau  := \tau \circ (-)$. The modification axiom holds by 
the naturality of the associator $\a$. 

It remains to give the extra data witnessing preservation of units and 
composition. For 
$\psi^{\lanext{\glueFun }}_X : \Id_{\lanext{\glueFun }X} \To \lanext{\glueFun 
}(\Id_X)$
we take the modification with components given by the structural isomorphisms
$\id_{\bicatX(\glueFun B, X)} \XRA{\iso} \Id_X \circ (-)$. Similarly, for a 
composable 
pair $X \xra{g} X' \xra{f} X''$ in $\catX$, the modification
$\phi^{\lanext{\glueFun }}_{f,g} : 
	\lanext{\glueFun }(f) \circ \lanext{\glueFun }(g) 
		\To
	\lanext{\glueFun }(f \circ g)$
has components 
$f \circ \left( g \circ (-)\right) \XRA\iso (f \circ g) \circ (-)$. 
\end{myconstr}

The preceding construction leads us to the following definition
(\cf~Definition~\ref{def:relative-hom-functor}).

\begin{mydefn}
For a category $\baseCat$ and pseudofunctor $\glueFun  : \baseCat \to 
\bicatX$, 
the bicategory of \Def{$\baseCat$-intensional Kripke relations of arity 
$\glueFun $}
is the glueing bicategory 
$\gl{\lanext{\glueFun }}$ associated
to the relative hom-pseudofunctor.
\end{mydefn}

To  
bicategorify~(\ref{eq:glueing-left-kan-diagram})
we employ the canonical 
equivalences
$
\Hom(\altCat \times \baseCat, \catV) \simeq
\Hom(\baseCat \times \altCat, \catV) \simeq 
\Hom\big(\baseCat, \Hom(\altCat, \catV)\big)$ of~\cite[\S1.34]{Street1980}.

\begin{mylemma} \label{lem:composition-with-yoneda-pseudonat-trans}
For any pseudofunctor $\glueFun  : \baseCat \to \bicatX$ there exists a 
pseudonatural 
transformation $(\lanNat, \lanCell)$ as in the diagram
\begin{equation} \label{eq:two-argument-lanext-diagram}
\begin{tikzcd}
\op\baseCat \times \baseCat
\arrow[swap]{dr}{\op{\glueFun } \times \glueFun }
\arrow{rr}{\Hom(-, =)} &
\arrow[phantom]{d}[description, 
yshift=1mm,font=\scriptsize]{\twocellDown{(\lanNat, \lanCell)}} &
\Cat \\
\: &
\op\bicatX \times \bicatX
\arrow[swap]{ur}{\Hom(-, =)} &
\:
\end{tikzcd}
\end{equation}
where 
\begin{gather*}
\op{\glueFun } := \glueFun  : ob(\op{\baseCat}) \to ob(\op{\bicatX}) \\
\op{(\glueFun_{B,C})} :=  
\op{\baseCat}(B, C) = \baseCat(C,B) \xra{\glueFun _{C,B}} \bicatX(C, B) = 
						\op{\bicatX}(C,B)
\end{gather*}
\begin{proof}
For the functors $\lanNat_{(B,C)} : \baseCat(B, C) \to \bicatX(\glueFun B, 
\glueFun C)$ we take 
$\glueFun _{B,C}$. For $f : B' \to B$ and $g : C \to C'$, the witnessing 
isomorphism 
$\lanCell_{(f,g)}$ in the diagram below
\begin{td}[column sep = 5em]
\baseCat(B,C) 
\arrow{r}{\baseCat(f,g)} 
\arrow[swap]{d}{\glueFun _{B,C}}
\arrow[phantom]{dr}[description]{\twocell{\lanCell_{(f,g)}}} &
\baseCat(B', C')
\arrow{d}{\glueFun _{B',C'}} \\

\bicatX(\glueFun B, \glueFun C) 
\arrow[swap]{r}{\bicatX(\op{\glueFun }f, \glueFun g)}&
\bicatX(\glueFun B', \glueFun C')
\end{td}
is defined to be the composite natural isomorphism
\begin{equation} \label{eq:lanext-2-cell-definition}
\glueFun \big(g \circ (h \circ f)\big) 
\XRA{(\phi^\glueFun _{g, h \circ f})^{-1}}
\glueFun (g) \circ \glueFun (h \circ f) 
\XRA{\glueFun (g) \circ (\phi^\glueFun _{h,f})^{-1}}
\glueFun (g) \circ \big( \glueFun h \circ \glueFun f\big)  
\end{equation}
This composite is natural in $g$ and $f$; the unit and associativity laws 
follow from the corresponding laws of a pseudofunctor.
\end{proof}
\end{mylemma}

\begin{mycor} \label{cor:lanext-transformation-construction}
For any pseudofunctor $\glueFun  : \baseCat \to \bicatX$ there exists a 
pseudonatural 
transformation 
$(\lanNat, \lanCell) : \Yon \To \lanext{\glueFun } \circ \glueFun  : 
	\baseCat \to \Hom(\op\baseCat, \Cat)$, which is 
given by the functorial action of $\glueFun $ on hom-categories.
\begin{proof}
Passing~(\ref{eq:two-argument-lanext-diagram}) through 
the equivalences 
$\Hom(\op\baseCat \times \baseCat, \Cat) \simeq
\Hom(\baseCat \times \op\baseCat, \Cat) \simeq
\Hom\big(\baseCat, \Hom(\op\baseCat, \Cat) \big)$ 
at an arbitrary  $P : \op\baseCat \times \baseCat \to \Cat$
yields the following:
\begin{align*}
\lambda (B, C)^{\op\baseCat \times \baseCat} \bind P(B, C) \mapsto 
\lambda (C, B)^{\baseCat \times \op\baseCat} \bind P(B, C) \mapsto
\lambda C^{\baseCat} \bind \lambda B^{\op\baseCat} \bind P(B, C)
\end{align*}
so that $\Hom(-, =) \mapsto \lambda C^\baseCat \bind \Yon C$ and 
$\Hom(\glueFun (-), \glueFun (=)) \mapsto \lambda C^\baseCat \bind 
\lanext{\glueFun }(C)$. By the 
preceding lemma, these are related by the pseudonatural transformation with 
components
$\lanNat_C := \glueFun _{(-), C} : \baseCat(-, C) \to \bicatX\big(\glueFun (-), 
\glueFun C\big)$ and 
witnessing 2-cells given as in~(\ref{eq:lanext-2-cell-definition}).
\end{proof}
\end{mycor}

We may now extend the Yoneda pseudofunctor $\Yon$ to its glued counterpart
$\glued{\Yon}$. 

\begin{myconstr}
For any pseudofunctor 
$\glueFun  : \baseCat \to \bicatX$,
define the \Def{extended Yoneda pseudofunctor}
$\glued{\Yon} : \baseCat \to \gl{\lanext{\glueFun }}$
as follows.

On objects, we set
\begin{equation} \label{eq:def-of-glued-Yon}
\glued{\Yon}B := \big( \Yon B, (\lanNat, \lanCell)_{(-, B)}, \glueFun B \big)
\end{equation}
where $(\lanNat, \lanCell)_{(-, B)}$ is pseudonatural since 
$(\lanNat, \lanCell)$ 
is pseudonatural in both arguments.

For a 1-cell $f : B \to B'$ in $\baseCat$, we define
$\glued{\Yon}f$ to be the 1-cell 
$(\Yon f, (\phi^\glueFun _{-, f})^{-1}, \glueFun f)$ 
as in the diagram
\begin{td}[column sep = 4em]
\baseCat(-, B) 
\arrow[swap]{d}{\glueFun _{-, B}} 
\arrow{r}{f \circ (-)} 
\arrow[phantom]{dr}[description]{\twocell{(\phi^\glueFun _{-, f})^{-1}}} &
\baseCat(-, B') 
\arrow{d}{\glueFun _{-, B'}} \\

\bicatX\big( \glueFun (-), \glueFun B\big) 
\arrow[swap]{r}{\glueFun (f) \circ (-)} &
\bicatX\big( \glueFun (-), \glueFun B'\big) 
\end{td}
On 2-cells, we set 
$\glued{\Yon}(\tau : f \To f' : B \to B')$ to 
be the pair $(\Yon \tau, \glueFun \tau)$, which satisfies the cylinder 
condition by the 
naturality of $\phi^\glueFun $. 

Finally we need to define $\psi^{\glued{\Yon}}$ and $\phi^{\glued{\Yon}}$. 
Since $\Yon \Id_X = (\Yon \Id_X, \glueFun \Id_X)$, we may 
take simply 
$\psi^{\glued{\Yon}} := (\psi^{\Yon}, \psi^\glueFun )$. This forms a 2-cell in 
$\gl{\lanext{\glueFun }}$ by the unit law on $(\lanNat, \lanCell)$. Similarly, 
for 
$\phi^{\glued{\Yon}}$ we take $(\phi^{\Yon}, \phi^\glueFun )$, which satisfies 
the 
cylinder condition by the associativity law on $(\lanNat, \lanCell)$. The three 
axioms to check then hold pointwise.
\end{myconstr}

In the next section we shall provide an explicit presentation of 
exponentials
$\expobj{\glued{\Yon}B}{\glued{X}}$ in the glueing bicategory, which will 
provide a bicategorical, glued correlate of the 
identification $\altexp{\yon B}{P} \iso P(- \times X)$ for presheaves. First, 
however, we finish 
our examination of the relative hom-pseudofunctor by showing that it 
preserves bilimits.

\begin{mylemma} \label{lem:lanext-preserves-bilimits}
For any pseudofunctor $\glueFun  : \baseCat \to \bicatX$ the relative 
hom-pseudofunctor
$\lanext{\glueFun } : \bicatX\to \Hom(\op\baseCat, \Cat)$ preserves all 
bilimits that 
exist in $\bicatX$. 
\begin{proof}
Let $ H : \catJ \to \bicatX$ be a pseudofunctor and suppose the bilimit 
$(\bilim_{j \in \catJ} Hj, \lambda_j)$ exists in $\bicatX$. By
Proposition~\ref{prop:hom-basecat-cat-bicomplete}, the bilimit 
$\bilim (\lanext{\glueFun } \circ H)$ exists in
$\Hom(\op\baseCat, \Cat)$
and is given pointwise. 

Now, since representable pseudofunctors preserve 
bilimits~(Lemma~\ref{lem:representables-and-adjoints-preserve-bilimits}), the 
canonical map 
$e_B : \bilim_{j \in \catJ} \bicatX( \glueFun B, Hj) \to 
	\bicatX( \glueFun B, \bilim_{j \in \catJ} Hj)$
is an equivalence for 
every $B \in \baseCat$. These extend canonically to a pseudonatural 
transformation, yielding 
the required equivalence
$\bilim (\lanext{\glueFun } \circ H) \XRA\simeq \lanext{\glueFun }\left( \bilim 
H \right)$.
\end{proof}
\end{mylemma}

It will be useful to have an explicit description of how $\lanext{\glueFun}$ 
preserves products. For this we rely on the $\postName$ 2-cells. 

\begin{mylemma} \label{lem:seq-post-pseudonatural-trans}
For any fp-bicategory $\fpBicat{\baseCat}$, the $n$-ary tupling operation and 
2-cells $\postName$ together form a pseudonatural transformation
$\prod_{i=1}^n \baseCat(-, B_i) \To 
	\baseCat{\left(-, \prodop_{i=1}^n B_i\right)}$, and 
hence 
an equivalence of pseudofunctors
$\prod_{i=1}^n \baseCat(-, B_i) \simeq 
	\baseCat{\left(-, \prodop_{i=1}^n B_i\right)}$ 
in 
$\Hom(\op\baseCat, \Cat)$. 
\begin{proof}
For every $X \in \baseCat$ the $n$-ary tupling operation defines a functor 
$\seq{-, \,\dots\, , =} : 
\prod_{i=1}^n \baseCat(X, B_i) \to \baseCat\left(X, \prod_{i=1}^n B_i\right)$ 
which, by the definition of an fp-bicategory~(Definition~\ref{def:fp-bicat}), 
is an equivalence in $\Cat$. For these 
functors to be the components of a pseudonatural transformation, we need to 
provide an invertible 2-cell filling the diagram below 
for every $f : Y \to X$:
\begin{td}[column sep = 6em]
\prod_{i=1}^n \baseCat(X, B_i)
\arrow[phantom]{dr}[description]{\twocell{}}
\arrow[swap]{d}{\seq{-, \,\dots\, , =}} 
\arrow{r}{\prod_{i=1}^n \baseCat(f, B_i)} &
\arrow{d}{\seq{-, \,\dots\, , =}}
\prod_{i=1}^n \baseCat(Y, B_i) \\
\baseCat(X, \prod_{i=1}^n B_i) 
\arrow[swap]{r}{\baseCat(f, \prod_{i=1}^n B_i)} &
\baseCat(Y, \prod_{i=1}^n B_i)
\end{td}
Thus, we require a natural isomorphism 
$\seq{h_1 \circ f, \,\dots\, , h_n \circ f} 
	\To \seq{h_1, \,\dots\, , h_n} \circ f$, for 
which we take $\post{\ind{h}, f}^{-1}$. The two axioms are exercises in using
Lemma~\ref{lem:PseudoproductCanonical2CellsLaws}.
\end{proof}
\end{mylemma}

\begin{prooflesscor} \label{cor:lanext-prod-pres}
For any pseudofunctor $\glueFun : \baseCat \to \bicatX$, the relative 
hom-pseudofunctor
$\lanext{\glueFun}$ extends to an fp-pseudofunctor
$(\lanext{\glueFun}, \prodPres)$
with
$\prodPres_{\ind{X}}$ given by the pseudonatural transformation 
$(\seq{-, \,\dots\, , =}, \postName)$ 
defined in the preceding lemma.
\end{prooflesscor}

\begin{myremark}
From the perspective of 
biuniversal arrows, Lemma~\ref{lem:seq-post-pseudonatural-trans} is an 
instance of 
Lemma~\ref{lem:preservation-of-biadjoints-to-equivalence}.
\end{myremark}

\subsection{Exponentiating by glued representables} 
In order to emulate Fiore's construction of the 1-cells $\quote$ and $\unquote$ 
in the glueing bicategory, we require a correlate of the following categorical 
fact:

\begin{prooflesslemma}[{\cite{Fiore2002}}]
For any cartesian category $\catB$, cartesian closed category
$\catX$ and cartesian functor 
$\glueFun  : \catB \to \catX$, 
the exponential 
$\altexp{\glued{\yon}B}{(P,p,X)}$
in $\gl{\lanext{\glueFun}}$
may be described explicitly as
\[
\altexp{\yon B}{P}
	\xra{\altexp{\yon B}{p}} 
\altexp{\yon B}{\lanext{\glueFun }(X)}
	\xra{\iso}
\lanext{\glueFun }\left(\expobj{\glueFun B}{X}\right)
\]
Here the unlabelled isomorphism is the composite
\[
\altexp{\yon B}{\lanext{\glueFun }(X)}	
	\xra{\iso}
\catX\left( \glueFun(- \times B), X \right)
	\xra\iso
\catX\left(\glueFun(-) \times \glueFun B, X\right) 
	\xra{\iso}
\catX\left(\glueFun(-), \expobj{\glueFun B}{X}\right)
\]
arising from the canonical isomorphism 
$[\yon B, P] \iso P(- \times X)$, the product-preservation of 
$\glueFun$,
and the cartesian closed structure on $\catX$.
\end{prooflesslemma} 

For the bicategorical version of this lemma 
we note that,
since products in $\Cat$ are strict, one 
obtains
$\id_P \times \id_Q = \id_{P \times Q}$ for every $P, Q : \op\baseCat \to 
\Cat$, so that 
$\altexp{\id_P}{(\natTrans, \natCell)} : 
	\altexp{P}{Q} \To \altexp{P}{Q'}$
is equal to $\Lambda{\big( (\natTrans, \natCell) \circ (e, \cellOf{e}) \big)}$
(recall from Section~\ref{sec:presheaves-cartesian-closed} that 
$(e, \cellOf{e})$ denotes the evaluation 1-cell in 
$\Hom(\op\baseCat, \Cat)$). With our (locally discrete) use-case in mind, we 
shall simplify what follows by assuming the bicategory $\baseCat$ to be a 
2-category.

\begin{prooflesspropn} \label{prop:explicit-glued-yon-exp}
For any 2-category $\baseCat$ with pseudo-products, 
cc-bicategory   
$\ccBicat\bicatX$
and fp-pseudofunctor 
$(\glueFun, \prodPres) : 
	\fpBicat\baseCat \to \fpBicat{\bicatX}$,
the exponential
$\expobj{\glued{\Yon}B}{\left(K, (\natTrans, \natCell), X\right)}$ 
in $\gl{\lanext{\glueFun}}$
may be given explicitly
by the following composite in 
$\Hom(\op\baseCat, \Cat)$:
\begin{equation} \label{eq:glueing-exponential-identified}
\altexp{\Yon B}{K} 
	\xra{{[\Yon B, (\natTrans, \natCell)]}}
\altexp{\Yon B}{\lanext{\glueFun}X} 
	\xra{u_{B,X}}
\lanext{\glueFun }(\expobj{\glueFun B}{X})
\end{equation} 
where
$u_{B,X}$ 
is the composite of equivalences
\begin{equation} \label{eq:def-of-u-x}
\altexp{\Yon B}{\lanext{\glueFun }X} \xra{\overset{(1)}{\simeq}} 
\bicatX\big(\glueFun (- \times B), X\big) \xra{\overset{(2)}{\simeq}}
\bicatX\big(\glueFun (-) \times \glueFun B, X\big) \xra{\overset{(3)}{\simeq}}
\bicatX\big(\glueFun (-), \expobj{\glueFun B}{X}\big)
\end{equation}
arising from the following, respectively:
\begin{enumerate}
\item The canonical 
equivalence arising 
from the identification of  
$(\lanext{\glueFun }X)(- \times B)$ as 
$\altexp{\Yon B}{\lanext{\glueFun}X}$ 
(Theorem~\ref{thm:exponentiating-by-Yoneda-identified}),
\item The fact that $\glueFun$ preserves products,
\item The definition of exponentials in $\bicatX$. \qedhere
\end{enumerate}
\end{prooflesspropn}

Our strategy is to show that the 
composite~(\ref{eq:glueing-exponential-identified}) is the left-hand leg of a 
pullback diagram in $\Hom(\op\baseCat, \Cat)$; by 
Lemma~\ref{lem:pullback-equivalent-glued-objects}, this is sufficient to prove 
an equivalence in the glueing bicategory.
We prove this using the following 
fact, which generalises the 1-categorical 
situation. 

\begin{mylemma} \label{lem:pullback-from-equivalence}
Let $\baseCat$ be a bicategory and 
$e : B \leftrightarrows C : f$ be any adjoint equivalence
in $\baseCat$, 
with witnessing invertible 2-cells 
$\un : \Id_C \XRA\iso e \circ f$ and 
$\co : f \circ e \XRA\iso \Id_B$. Then for any $r : A \to C$ 
the pullback of the cospan 
$(B \xra{e} C \xla{r} A)$ exists and is given by
\begin{equation} \label{eq:pullback-given-explicitly}
\begin{tikzcd}[row sep = 2.3em]
A 
\arrow[bend right = 50]{dd}[swap]{f \circ r}
\arrow[phantom]{dr}[description, yshift=0mm]{\iso}
\arrow[swap]{d}{r}
\arrow{r}{\Id_A} &
A 
\arrow{dd}{r} \\
C 
\arrow[phantom]{dr}[description, yshift=-1mm]{\twocellIso{\un}}
\arrow[bend left]{dr}[near start, yshift=-1mm]{\Id_C}
\arrow[swap]{d}{f} &
\: \\
B 
\arrow[swap]{r}{e} &
C
\end{tikzcd}
\end{equation}
where the top isomorphism is a composite of structural isomorphisms.
\begin{proof}
Suppose given any other iso-commuting square
\begin{td}
X
\arrow[phantom]{dr}[description]{\twocellIso{\rho}} 
\arrow{r}{p} 
\arrow[swap]{d}{q} &
A 
\arrow{d}{r} \\

B 
\arrow[swap]{r}{e} &
C
\end{td}
We take the mediating map $X \to A$ to be $p$. For the 2-cells we take 
$\Gamma := \Id_A \circ p \XRA\iso p$ and $\Delta$ to be defined by the 
following diagram:
\begin{td}
(f \circ r) \circ p
\arrow[swap]{d}{\iso}
\arrow{r}{\Delta} &
q \\
f \circ (r \circ p)
\arrow[swap]{d}{f \circ \rho} &
\Id_B \circ q
\arrow[swap]{u}{\iso} \\
f \circ (e \circ q)
\arrow[swap]{r}{\iso} &
(f \circ e) \circ q
\arrow[swap]{u}{\co \circ q}
\end{td}
A short diagram chase using the triangle law relating $\un$ and $\co$ shows 
this is a fill-in.

Next we claim that $(p, \Gamma, \Delta)$ is universal. To this end, let 
$(v, \Sigma_1, \Sigma_2)$ be any other fill-in, so that the following diagram 
commutes:
\begin{equation} \label{eq:fill-in-property}
\begin{tikzcd}
(r \circ \Id_A) \circ v
\arrow{d}
\arrow{r}{\iso} &
r \circ (\Id_A \circ v)
\arrow{r}{r \circ \Sigma_1} &
r \circ p
\arrow{d}{\rho} \\
\left( e \circ (f \circ r) \right) \circ v
\arrow[swap]{r}{\iso} &
e \circ \left( (f \circ r) \circ v\right)
\arrow[swap]{r}{e \circ \Sigma_2} &
e \circ q
\end{tikzcd}
\end{equation}
The unlabelled arrow is the composite~(\ref{eq:pullback-given-explicitly}) 
given in the claim.

We define 
$\trans\Sigma
:= v \XRA{\iso} \Id_A \circ v \XRA{\Sigma_1} p$, and claim that 
both the following equations hold:
\begin{equation} \label{eq:fill-in-universality-conditions}
\begin{tikzcd}
\Id_A \circ v 
\arrow{rr}{\Id_A \circ \trans\Sigma}
\arrow[swap]{dr}{\Sigma_1} &
\: &
\Id_A \circ p 
\arrow{dl}{\Gamma} \\
\: &
p  &
\:
\end{tikzcd}
\qquad\qquad
\begin{tikzcd}
(f \circ r) \circ v 
\arrow{rr}{(f \circ r) \circ \trans\Sigma}
\arrow[swap]{dr}{\Sigma_2} &
\: &
(f \circ r) \circ p
\arrow{dl}{\Delta} \\
\: &
q  &
\:
\end{tikzcd}
\end{equation}
The right-hand diagram is an relatively easy check. The left-hand diagram 
follows by naturality, the triangle law relating $\un$ and $\co$, and the 
assumption~(\ref{eq:fill-in-property}).

It remains to check the uniqueness condition for $\trans\Sigma$. For 
any other 
$\Theta : v \To p$ satisfying the two diagrams 
of~(\ref{eq:fill-in-universality-conditions}), one sees that
\begin{td}[row sep = 2.5em]
v
\arrow[phantom]{dr}[description]{\equals{nat.}}
\arrow{r}{\Theta}
\arrow[swap]{d}{\iso} &
p 
\arrow[bend left = 60, equals]{dd}
\arrow{d}[description]{\iso} \\

\Id_A \circ v
\arrow[swap]{dr}{\Sigma_1}
\arrow{r}{\Id_A \circ \Theta} &
\Id_A \circ p 
\arrow{d}[description]{\iso} \\

\: &
p
\end{td}
where the bottom triangle commutes by the right-hand diagram 
of~(\ref{eq:fill-in-universality-conditions}), and the left-hand leg is 
exactly the definition of $\trans\Sigma$. Hence $\Theta = \trans\Sigma$ as 
required. Finally we observe that $\trans\id$ is certainly invertible.
\end{proof}
\end{mylemma}

The requirement for an adjoint equivalence in the preceding lemma is, by the 
usual argument, no stronger than requiring just an equivalence 
(\eg~\cite[Proposition~1.5.7]{Leinster2004}). Importantly, the adjoint 
equivalence one constructs from an equivalence has the same 1-cells. 

In the light of the lemma, if we can show that the equivalence 
$u_{B,X}$ defined 
in~(\ref{eq:def-of-u-x}) has a pseudo-inverse given by the composite
$\altexp{(\lanNat, \lanCell)_{(-, B)}}{\lanext{\glueFun}X}
 	\circ \evBar_{\glueFun B, X}$,
then the following is a pullback diagram:
\begin{td}[column sep = 5em]
\altexp{\Yon B}{K}
\arrow[swap]{d}{[\Yon B, (\natTrans, \natCell)]}
\arrow[phantom]{drr}[description]{\iso}
\arrow{rr}{\Id_{\altexp{\Yon B}{K}}} &
\: &
\altexp{\Yon B}{K} 
\arrow{dd}{\Lambda( (\natTrans, \natCell) \circ (e, \overline{e}))} \\

\altexp{\Yon B}{\lanext{\glueFun }X} 
\arrow[phantom]{dr}[description]{\iso}
\arrow[bend left=10]{drr}[description]
	{\Id_{\scriptsizescriptsizeexpobj{\Yon B}{\lanext{\glueFun}X}}}
\arrow[swap]{d}{u_{B,X}} &
\: &
\: \\

\lanext{\glueFun }(\expobj{\glueFun B}{X}) 
\arrow[swap]{r}{\evBar_{\glueFun B,X}} &
\altexp{\lanext{\glueFun }(\glueFun B)}{\lanext{\glueFun }X} 
\arrow[swap]{r}[yshift=-2mm]
{\Lambda ( (e, \overline{e}) \circ 
(\altexp{\lanext{\glueFun }(\glueFun B)}{\lanext{\glueFun }X} \times (\lanNat, 
\lanCell)) )} 
&
\altexp{\Yon B}{\lanext{\glueFun}X}
\end{td}
It will then follow that for any
$\glued{K} := \left(K, (\natTrans, \natCell), X\right)$
the composite~(\ref{eq:glueing-exponential-identified})---the 
left-hand leg of the above diagram---is an explicit description of the 
exponential 
$(\expobj{\glued{\Yon}X}{\glued{K}})$. The difficulty, therefore, is not in 
showing that $u_{B,X}$ is an equivalence, but in checking whether it has a 
pseudo-inverse of the form we require. We turn to this next. (The cartesian 
closed structures we employ are summarised in 
Appendix~\ref{chap:cartesian-closed-structures}).

\subsubsection{The equivalence 
$\altexp{\Yon B}{\lanext{\glueFun }X} \simeq 
	\lanext{\glueFun}(\expobj{\glueFun B}{X})$: calculating the 1-cells}
\label{sec:1-cells-calculated}

In this section we shall calculate the action of the maps
$u_{B,X}$ and 
$\altexp{(\lanNat, \lanCell)_{(-, B)}}{\lanext{\glueFun}X}
	\circ \evBar_{\glueFun B, X}$;
in the next section we shall show these form an equivalence. To shorten 
notation, let us introduce the following abbreviation:
\begin{align*}
\wholecell_{B,X} &:= 	
	\altexp{(\lanNat, \lanCell)_{(-, B)}}{\lanext{\glueFun}X} \circ 
	\evBar_{\glueFun B,X} 
\end{align*}
Our first task is to unfold each of the equivalences in the definition of 
$u_{B,X}$ to determine the action of the whole composite.

\paragraph*{Calculating the composite $u_{B,X}$.} 
If $\altexp{X}{Y}$ and $\exp{X}{Y}$ are both the exponential of $X$ and $Y$ in 
a bicategory $\baseCat$,  with associated currying operation and evaluation 
maps $\lambda, \eval_{X,Y}$ and $\widehat{\lambda}, \widehat{\eval}_{X,Y}$, 
respectively, 
then 
$\widehat{\lambda}{\left((\altexp{X}{Y}) \times X \xra{\eval_{X,Y}} Y \right)} 
: 
\altexp{X}{Y} \to (\exp{X}{Y})$
is canonically an equivalence.

Now let $\fpBicat\baseCat$ be a 2-category with pseudo-products, 
$B \in \baseCat$, and 
$P : \op\baseCat \to \Cat$ be any pseudofunctor. We calculate the 
equivalence
$\altexp{\Yon B}{P} 
	= \Hom(\op\baseCat, \Cat)\left( \Yon(-) \times \Yon B, P\right)
	\xra{\simeq} {P(- \times B)}$
arising from 
Theorem~\ref{thm:exponentiating-by-Yoneda-identified}.
The evaluation 1-cell
$\eval_{\Yon B, P} : \altexp{\Yon B}{P} \times \Yon B \to P$ is the 
pseudonatural transformation $(e, \overline{e})$ with components
\begin{align*}
\Hom(\op\baseCat, \Cat)(\Yon C \times \Yon B, P) \times \baseCat(C,B) 
	&\xra{e_C} 
	PC \\
\big( (\natTrans, \natCell), h \big) &\mapsto \natTrans_C (\Id_C, h)
\end{align*}
On the other hand, the currying operation 
\[
\widehat{\Lambda} : \Hom(\op\baseCat, \Cat)(R \times \Yon B, P) \to 
\Hom( \op\baseCat, \Cat)\big(R, P(- \times B)\big)
\]
witnessing $P(- \times  X)$ as an exponential  
takes a pseudonatural transformation 
$(\altNat, \altCell)$ 
to the pseudonatural transformation with components
$RC 
	\xra{R\pi_1} 
R(C \times B) 
	\xra{\altNat_{C \times B}(-, \pi_2)} 
{P(C \times B)}$. 
Using the assumption that $\baseCat$ is a 2-category, the 
component of the canonical equivalence
$\altexp{\Yon B}{P} \xra{\simeq} P(- \times B)$
at $C \in \baseCat$ is therefore
\begin{equation} \label{eq:canonical-equivalence-of-exponentials}
\begin{aligned}
\Hom(\op\baseCat, \Cat)\big(\Yon C \times \Yon B, P\big)  &\to P(C \times B) \\
(\natTrans, \natCell) &\mapsto  
\natTrans_{C \times B}(\pi_1, \pi_2)
\end{aligned}
\end{equation}
It follows that $u_{B,X}(C)$ is the following composite:
\begin{equation} \label{eq:u-identified}
\begin{aligned}
\altexp{\Yon B}{\lanext{\glueFun }X}(C) \xra{\simeq} 
\bicatX\big(\glueFun (C \times B), X\big) \xra{\simeq}
\bicatX\big(\glueFun C \times \glueFun B, X\big) \xra{\simeq}
\bicatX\big(\glueFun C, \expobj{\glueFun B}{X}\big) \\
(\natTrans,\natCell) \mapsto 
\natTrans_{C \times B}(\pi_1, \pi_2) \mapsto
\natTrans_{C \times B}(\pi_1, \pi_2) \circ \prodPres_{C,B} \mapsto
\lambda {\big( \natTrans_{C \times B}(\pi_1, \pi_2) \circ 
\prodPres_{C,B} \big)}
\end{aligned}
\end{equation}
Next we turn to calculating
$\wholecell_{B,X} := 	
	\altexp{(\lanNat, \lanCell)_{(-, B)}}{\lanext{\glueFun}X} \circ 
		\evBar_{\glueFun B, X}$. 

\paragraph*{Calculating 
	$\altexp{(\lanNat, \lanCell)_{(-, B)}}{\lanext{\glueFun}X}$.}
We begin by calculating the composite
\begin{equation} \label{eq:e-circ-id-times-l}
\begin{tikzcd}[column sep = 4em]
\altexp{\lanext{\glueFun }(\glueFun B)}{\lanext{\glueFun }(X)} \times \Yon B
\arrow{r}[yshift=2mm]{\altexp{\lanext{\glueFun }(\glueFun B)}{\lanext{\glueFun 
}(X)} \times 
(\lanNat, \lanCell)_{(-, B)}} &
\altexp{\lanext{\glueFun }(\glueFun B)}{\lanext{\glueFun }(X)} \times 
\lanext{\glueFun }{\glueFun B} 
\arrow{r}{(e, \overline{e})} &
\lanext{\glueFun }(X)
\end{tikzcd}
\end{equation}
Applying the definition of $(e, \widebar{e})$ again, the component of the 
composite~(\ref{eq:e-circ-id-times-l}) at 
$C \in \baseCat$ is
\begin{align*}
\Hom(\op\baseCat, \Cat)
\big( \baseCat(-, C) \times \bicatX(\glueFun (-), \glueFun B), \bicatX(\glueFun 
(-), X) \big) \times 
\baseCat(C, B) &\to
\bicatX(\glueFun C, X) \\
\big( (\natTrans, \natCell), h \big) &\mapsto 
\natTrans(C, \Id_C, \glueFun h)
\end{align*}
Naturality in $C$ is witnessed by the following 2-cell, where $r : C' \to C$
is any 1-cell in $\baseCat$:
\begin{td}
\natTrans\big(C', \Id_{C'} \circ r, \glueFun (h \circ r)\big)
\arrow{r}
\arrow[swap]{d}
	{\natTrans\left(C', \Id_{C'} \circ r, (\phi^\glueFun _{h,r})^{-1} \right)} &
\natTrans(C, \Id_C, \glueFun h) \circ \glueFun r \\

\natTrans(C', \Id_{C'} \circ r, \glueFun h \circ \glueFun r ) 
\arrow[equals]{r} &
\natTrans(C', r \circ \Id_{C}, \glueFun h \circ \glueFun r ) 
\arrow[swap]{u}{\natCell(r, \Id_C, \glueFun h)}
\end{td}
Instantiating this with the cartesian closed structure constructed in 
Section~\ref{sec:presheaves-cartesian-closed}, 
one may identify
$\altexp{(\lanNat, \lanCell)_{(-, B)}}{\lanext{\glueFun}X} : 
	\altexp{\lanext{\glueFun }(\glueFun B)}{\lanext{\glueFun}(X)} 
		\to 
	\altexp{\Yon B}{\lanext{\glueFun }(X)}$
as in the following lemma. 

\begin{prooflesslemma} \label{lem:action-of-second-composite}
For any 2-category with pseudo-products 
$\fpBicat\baseCat$, 
cc-bicategory  
$\ccBicat\bicatX$,
and fp-pseudofunctor 
$(\glueFun, \prodPres) : 
	\fpBicat\baseCat \to \fpBicat{\bicatX}$,
the pseudonatural transformation
$\altexp{(\lanNat, \lanCell)_{(-, B)}}{\lanext{\glueFun}X} 
: 
	\altexp{\lanext{\glueFun }(\glueFun B)}{\lanext{\glueFun}(X)} 
		\To 
	\altexp{\Yon B}{\lanext{\glueFun }(X)}$ 
(where ${B \in \baseCat}$ and ${X \in \bicatX}$)
has functorial components
\begin{align*}
\altexp{\lanext{\glueFun }(\glueFun B)}{\lanext{\glueFun}(X)}(C) 
	&\xra{\altexpsmall{(\lanNat, \lanCell)_{(-, B)}}{\lanext{\glueFun}X}(C)}
\altexp{\Yon B}{\lanext{\glueFun }(X)}(C) \\
(\natTrans, \natCell) &\mapsto 
\lambda A^\baseCat \bind
\lambda h^{A \to C} \bind
\lambda p^{A \to B} \bind
\natTrans(A, h, \glueFun p)
\end{align*}
For $s : A' \to A$, the witnessing 2-cell of 
${\altexp{(\lanNat, \lanCell)_{(-, B)}}{\lanext{\glueFun}X}(C)}
	{\big( (\natTrans, \natCell) \big)}$ 
as in the diagram
\begin{td}[column sep = 6em]
\baseCat(A, C) \times \baseCat(A, B) 
\arrow[phantom]{dr}[description]{\twocell{\iso}}
\arrow{r}{\baseCat(s, C) \times \baseCat(s, B)}
\arrow[swap]{d}{\natTrans(A, -, \glueFun (=))} &
\baseCat(A', C)  \times \baseCat(A', B)
\arrow{d}{\natTrans(A', -, \glueFun (=))} \\
\bicatX(\glueFun A, X)
\arrow[swap]{r}{\bicatX(\glueFun s, X)} &
\bicatX(\glueFun A', X)
\end{td}
is given by
\begin{equation*} 
\natTrans\big(A', ({-}) \circ s, \glueFun ({=} \circ s)\big) 
\XRA{\natTrans(A', ({-}) \circ s, (\phi^\glueFun_{({=}),s})^{{-}1} ) }
\natTrans\big(A', ({-}) \circ s, \glueFun ({=}) \circ \glueFun s\big) 
\XRA{\natCell(s, {-}, \glueFun ({=}))}
\natTrans(A, {-}, \glueFun ({=})) \circ \glueFun s
\end{equation*}
\end{prooflesslemma}

\paragraph*{Calculating $\evBar_{\glueFun B, X}$.}
By Lemma~\ref{lem:seq-post-pseudonatural-trans}, the 
pseudonatural transformation 
$\lanext{\glueFun}(\eval_{\glueFun B, X}) \circ
\prodPres_{\glueFun B, X}$
has components
defined by 
$\lambda C^\baseCat \bind 
\lambda h^{\glueFun C \to (\scriptsizeexpobj{\glueFun B}{X})} \bind
\lambda g^{\glueFun C \to \glueFun B} \bind
\eval_{\glueFun B, X} \circ \seq{h, g}$ 
and witnessing \mbox{2-cells} of the form
\begin{td}
\bicatX\big(\glueFun C, \expobj{\glueFun B}{X}\big) \times 
	\bicatX\big(\glueFun C, \glueFun B\big)
\arrow{r}[yshift=2mm]{\bicatX(\glueFun f, \expobj{\glueFun B}{X}) \times 
\bicatX(\glueFun f, \glueFun B)} 
\arrow[swap]{d}{\eval_{\glueFun B,X} \circ \seq{-, =}}
\arrow[phantom]{dr}[description]{\twocell{\iso}}
&
\bicatX\big(\glueFun C', \expobj{\glueFun B}{X}\big) \times 
\bicatX\big(\glueFun C', \glueFun B\big) 
\arrow{d}{\eval_{\glueFun B,X} \circ \seq{-, =}} \\
\bicatX\big(\glueFun C, X\big) 
\arrow[swap]{r}{\bicatX(\glueFun f, X)} &
\bicatX\big(\glueFun C', X\big)
\end{td}
given by
\[
\eval_{\glueFun B, X} \circ \seq{h \circ \glueFun f, g \circ \glueFun f} 
\XRA{\eval_{\glueFun B,X} \circ \postName^{-1}}
\eval_{\glueFun B,X} \circ (\seq{h,g} \circ \glueFun f)
\XRA\iso
(\eval_{\glueFun B,X} \circ \seq{h,g}) \circ \glueFun f
\]
for every $f : C' \to C$ in $\baseCat$.
Applying the currying operation defined in  
Section~\ref{sec:presheaves-cartesian-closed},
one obtains the following characterisation of $\evBar_{\glueFun B, X}$.

\newpage
\begin{prooflesslemma} \label{lem:def-of-first-composite}
For any 2-category with pseudo-products 
$\fpBicat\baseCat$, 
cc-bicategory   
$\ccBicat\bicatX$,
and fp-pseudofunctor 
$(\glueFun, \prodPres) : 
	\fpBicat\baseCat \to \fpBicat{\bicatX}$,
the pseudonatural 
transformation 
$\evBar_{\glueFun B,X}$ 
has components
$\evBar_{\glueFun B,X}(C)$ given by the functors
\begin{align*}
\bicatX(\glueFun C, \expobj{\glueFun B}{X}) &\to 
\Hom(\op\baseCat, \Cat)\big( 
\Yon C \times \lanext{\glueFun}(\glueFun B),  
\lanext{\glueFun}X
\big) \\
f &\mapsto 
\lambda A^\baseCat \bind
\lambda (h^{A \to C}, g^{\glueFun A \to \glueFun B}) \bind
\big(\glueFun A \xra{\seq{f \circ \glueFun h, g}} (\exp{\glueFun B}{X}) \times 
\glueFun B \xra{\eval_{\glueFun B,X}} 
X\big)
\end{align*}
Moreover, for every $r : A' \to A$ the pseudonatural transformation 
$\evBar_{\glueFun B, X}(C)(f)$ has 
witnessing 2-cell
\begin{td}[column sep = 5em]
\baseCat(A, C) \times \bicatX(\glueFun A, \glueFun B) 
\arrow[phantom]{dr}[description]{\twocell{\overline{\evBar_{\glueFun B, 
X}(C)(f)}_r}}
\arrow{r}[yshift=2mm]{\baseCat(r, C) \times \bicatX(\glueFun r, \glueFun B)}
\arrow[swap]{d}{\eval_{\glueFun B,X} \circ \seq{f \circ \glueFun (-), =}} &
\baseCat(A', C) \times \bicatX(\glueFun A', \glueFun B) 
\arrow{d}{\eval_{\glueFun B,X} \circ \seq{f \circ \glueFun (-), =}} \\
\bicatX(\glueFun A, X) 
\arrow[swap]{r}{\baseCat(\glueFun r, X)} &
\bicatX(\glueFun A', X)
\end{td}
defined by
\begin{equation*}
\makebox[\textwidth]{
\begin{tikzcd}[column sep = 0em, ampersand replacement = \&]
\eval_{\glueFun B,X} \circ 
	\seqlr{f \circ \glueFun (h \circ r), g \circ \glueFun r} 
\arrow{rr}{\overline{\evBar_{\glueFun B, X}(C)(f)}_r}
\arrow[swap]{d}
	{\eval_{\glueFun B,X} \circ 
		\seqlr{f \circ (\phi^\glueFun _{h,r})^{-1}, g \circ \glueFun r }} \&
\: \&
\left(\eval_{\glueFun B,X} \circ \seqlr{f \circ \glueFun h, g}\right)  
	\circ \glueFun r \\
\eval_{\glueFun B,X} \circ 
	\seqlr{f \circ (\glueFun h \circ \glueFun r), g \circ \glueFun r}
\arrow[swap]{dr}{\iso} \&
\: \&
\eval_{\glueFun B,X} \circ \left(\seqlr{f \circ \glueFun h, g}  
	\circ \glueFun r\right)
\arrow[swap]{u}{\iso} \\
\: \&
\eval_{\glueFun B,X} \circ 
	\seqlr{(f \circ \glueFun h) \circ \glueFun r, g \circ \glueFun r} 
\arrow[swap]{ur}{\eval_{\glueFun B,X} \circ \postName^{-1}}  \&
\:
\end{tikzcd}
}
\end{equation*}
\end{prooflesslemma}

\paragraph*{Calculating $\wholecell_{B,X}$.}

Combining~Lemma~\ref{lem:action-of-second-composite} with 
Lemma~\ref{lem:def-of-first-composite}, one obtains the following 
identification of $\wholecell_{B,X}$.

\begin{prooflesslemma}
For any 2-category with pseudo-products 
$\fpBicat\baseCat$, 
cc-bicategory   
$\ccBicat\bicatX$,
and fp-pseudofunctor 
$(\glueFun, \prodPres) : 
	\fpBicat\baseCat \to \fpBicat{\bicatX}$,
the composite pseudonatural transformation
$\wholecell_{B,X} :  \lanext{\glueFun }(\expobj{\glueFun B}{X}) \to 
\altexp{\Yon 
B}{\lanext{\glueFun }X}$
has components
\begin{align*}
\bicatX(\glueFun C, \exp{\glueFun B}{X}) &\xra{\wholecell_{B,X}(C)} 
\Hom(\op\baseCat, \Cat)\big( \Yon C \times \Yon B, \bicatX(\glueFun (-), X) 
\big) \\
f &\mapsto 
\lambda A^\baseCat \bind
\lambda h^{A \to C} \bind
\lambda p^{A \to B} \bind
\big( \glueFun A \xra{\seq{f \circ \glueFun h, \glueFun p}} (\exp{\glueFun 
B}{X}) \times \glueFun B \xra{\eval_{\glueFun B, X}} 
X \big)
\end{align*}
The witnessing 2-cells for the pseudonatural transformation 
$\wholecell_{B,X}(C)(f)$
are defined by the following commutative diagram, where
$r : A' \to A$ 
is any 1-cell:
\begin{equation} \label{eq:b-circ-u-x-2-cell}
\makebox[\textwidth]{
\begin{tikzcd}[column sep = -1em, ampersand replacement = \&]
\eval_{\glueFun B,X} \circ 
	\seqlr{f \circ \glueFun (h \circ r), \glueFun(p \circ r)} 
\arrow{rr}{\overline{\wholecell_{B,X}(C)(f)}_r}
\arrow[swap]{d}
{\eval_{\glueFun B,X} \circ 
	\seq{f \circ (\phi^\glueFun _{h,r})^{-1}, (\phi^\glueFun _{p,r})^{-1}} } \&
\: \&
\eval_{\glueFun B,X} \circ 
	\seqlr{f \circ \glueFun h, \glueFun p}  \circ \glueFun r \\
\eval_{\glueFun B,X} \circ 
	\seqlr{f \circ \left(\glueFun h \circ \glueFun r\right), 
	\glueFun p \circ \glueFun r} 
\arrow[swap]{dr}
	{\iso} \&
\: \&
\eval_{\glueFun B,X} \circ 
	\left(\seqlr{f \circ \glueFun h, 
	\glueFun p} \circ \glueFun r\right)
\arrow[swap]{u}{\iso} \\
\: \&
\eval_{\glueFun B,X} \circ 
	\seqlr{\left(f \circ \glueFun h\right) \circ \glueFun r, 
	\glueFun p \circ \glueFun r}
\arrow[swap]{ur}[yshift=1mm, xshift=1mm]
	{\eval_{\glueFun B,X} \circ \postName^{-1}} \&
\:
\end{tikzcd}
}
\vspace{2mm}
\end{equation}
\end{prooflesslemma}

%

\subsubsection{The equivalence 
$\altexp{\Yon B}{\lanext{\glueFun }X} \simeq 
	\lanext{\glueFun}(\expobj{\glueFun B}{X})$}

We are finally in a position to prove that
$u_ X : \altexp{\Yon B}{\lanext{\glueFun }X} \leftrightarrows \lanext{\glueFun 
}(\exp{\glueFun B}{X}) :
\wholecell_{B,X} $
defines an equivalence of pseudofunctors
in $\Hom(\op\baseCat, \Cat)$. By 
Lemma~\ref{lem:equivalence-of-pseudofunctors-pointwise} it suffices to 
construct an equivalence of categories 
$
u_{B,X}(C) : \altexp{\Yon B}{\lanext{\glueFun }X}(C) \leftrightarrows 
\lanext{\glueFun }(\exp{\glueFun B}{X})(C) :
\wholecell_{B,X}(C)
$
for each $C \in \baseCat$. We deal with this in the following lemma.

\begin{mylemma}
For any 2-category with pseudo-products 
$\fpBicat\baseCat$, 
cc-bicategory  
$\ccBicat\bicatX$,
and fp-pseudofunctor 
$(\glueFun, \prodPres) : 
	\fpBicat\baseCat \to \fpBicat{\bicatX}$, 
the following composites are naturally isomorphic to the identity functor for 
every $B, C \in \baseCat$ and $X \in \bicatX$:
\begin{enumerate}
\item \label{c:isomorphic-to-id-functor}
\[
\bicatX(\glueFun C, \exp{\glueFun B}{X}) 
\xra{\wholecell_{B,X}(C)} 
\Hom(\op\baseCat, \Cat)(\Yon C \times \Yon B, \lanext{\glueFun }X) 
\xra{u_{B,X}(C)}
\bicatX(C, \exp{\glueFun B}{X}) 
\]

\item \label{c:isomorphic-to-id-pseudonat-trans} 
\begin{equation*}
\begin{tikzcd}[column sep = -1em]
\Hom(\op\baseCat, \Cat)(\Yon C \times \Yon B, \lanext{\glueFun }X) 
\arrow[swap]{dr}{u_{B,X}(C)}
\arrow{rr} &
\: &
\Hom(\op\baseCat, \Cat)(\Yon C \times \Yon B, \lanext{\glueFun }X) \\
\: &
\bicatX(\glueFun C, \exp{\glueFun B}{X}) 
\arrow[swap]{ur}{\wholecell_{B,X}(C)} &
\:
\end{tikzcd}
\end{equation*}
\end{enumerate}
Hence, $\wholecell_{B,X}$ is pseudo-inverse to 
$u_{B,X} : \altexp{\Yon B}{\lanext{\glueFun }X} \to \lanext{\glueFun 
}(\exp{\glueFun B}{X})$ in
$\Hom(\op\baseCat, \Cat)$. 
\begin{proof}
For~(\ref{c:isomorphic-to-id-functor}), we begin by calculating
\begin{align*}
\big(u_{B,X}(C) \circ \wholecell_{B,X}(C)\big)(f) &= 
u_{B,X}(C)\big( 
\lambda A^\baseCat \bind
\lambda h^{A \to C} \bind
\lambda p^{A \to B} \bind
\eval_{\glueFun B,X} \circ \seq{f \circ \glueFun h, \glueFun p}
\big) \\
&=  \lambda \big( 
\glueFun C \times \glueFun B 
\xra{\prodPres_{C,B}} 
\glueFun (C \times B) 
\xra{\eval_{\glueFun B,X} \circ \seq{f \circ \glueFun \pi_1, \glueFun \pi_2}} 
X \big) 
\end{align*}
for $f : \glueFun C \to (\exp{\glueFun B}{X})$.
For each such $f$, one 
obtains an invertible 2-cell 
${\big( u_{B,X} \circ \wholecell_{B,X}(C) \big)(f)} \XRA{\iso}f$ 
as the composite
\begin{td}[column sep = 5em, row sep = 2.5em]
\lambda 
\big( {\left(\eval_{\glueFun B, X} 
		\circ \seq{f \circ \glueFun \pi_1, \glueFun \pi_2}\right) 
		\circ \prodPres_{C,B}}\big) 
\arrow[swap]{d}
{\lambda ( \eval_{\glueFun B, X} \circ \fuse^{-1} \circ \prodPres_{C,B}) }
\arrow{r} &
f \\

\lambda 
\big( {\left(\eval_{\glueFun B, X} 
	\circ \left((f \times \glueFun B) 
	\circ \seq{\glueFun \pi_1, \glueFun \pi_2}\right)\right) 
	\circ \prodPres_{C,B}}\big) 
\arrow[swap]{d}{\iso} &
\lambda \big( \eval_{\glueFun B, X} 
	\circ (f \times \glueFun B)\big) 
\arrow[swap]{u}{\etaExp{f}^{-1}} \\

\lambda 
\big( {\left(\eval_{\glueFun B, X} 
	\circ (f \times \glueFun B)\right) 
	\circ \left(\seq{\glueFun \pi_1, \glueFun \pi_2}\right) 
	\circ \prodPres_{C,B}}\big) 
\arrow[swap]{r}[yshift=-2mm]
	{\lambda 
		( \eval_{\glueFun B, X} \circ (f \times \glueFun B) 
			\circ (\unTimes_{C,B})^{-1})} & 
\lambda \big( {\left(\eval_{\glueFun B, X} 
	\circ (f \times \glueFun B)\right) 
	\circ \Id_{\glueFun B \times \glueFun C}}\big) 
\arrow[swap]{u}{\iso} 
\end{td}
where the bottom isomorphism arises from the equivalence 
$\seq{\glueFun \pi_1, \glueFun \pi_2} : \glueFun (B \times C) \leftrightarrows 
{\glueFun B \times \glueFun C : \prodPres_{C,B}}$
witnessing $(\glueFun , \prodPres)$ as an fp-pseudofunctor. This composite is 
clearly 
natural in $f$, so one obtains the required natural isomorphism.

For~(\ref{c:isomorphic-to-id-pseudonat-trans}) one must work a little 
harder. We 
are required to
construct an invertible modification 
$\modif^{(\natTrans, \natCell)} : \big( \wholecell_{B,X}(C) 
\circ u_{B,X}(C)\big)
\big( (\natTrans, \natCell) \big)  \xra\iso (\natTrans, \natCell)$
for every pseudonatural transformation
$(\natTrans, \natCell) : 
	\Yon C \times \Yon B \To \bicatX\big( \glueFun (-), X \big)$, 
and this family
which must be
natural in the sense that, for any modification 
$\altModif : (\natTrans, \natCell) \to (\altNat, \altCell)$, 
the following diagram commutes:
\begin{equation} \label{eq:modification-natural-in-2-cell}
\begin{tikzcd}[column sep = 9em]
\big( \wholecell_{B,X}(C) \circ u_{B,X}(C)\big)
\big( (\natTrans, \natCell) \big) 
\arrow[swap]{d}{\iso}
\arrow{r}{( \wholecell_{B,X}(C) \circ u_{B,X}(C))(\altModif)} &
\big( \wholecell_{B,X}(C) \circ u_{B,X}(C)\big)
\big( (\altNat, \altCell) \big)
\arrow{d}{\iso} \\
(\natTrans, \natCell) 
\arrow[swap]{r}{\altModif} &
(\altNat, \altCell) 
\end{tikzcd}
\end{equation}
To this end, let us first unwind the data we are given. Applying the 
work of the preceding section, one sees that for 
$(\natTrans, \natCell) : 
	\Yon C \times \Yon B 
	\to 
	\bicatX\big( \glueFun (-), X \big)$
one has
\begin{align*}
\big( \wholecell_{B,X}(C) &\circ u_{B,X}(C)\big)\big( (\natTrans, \natCell) 
\big) \\
&=
\wholecell_{B,X}(C) 
\Big( \lambda\big(\natTrans_{C \times B}(\pi_1, \pi_2) \circ 
\prodPres_{C,B} 
\big)\Big) 
\\ 
&=  
\lambda A^\baseCat \bind
\lambda h^{A \to C} \bind
\lambda p^{A \to B} \bind
\eval_{\glueFun B, X} \circ 
\seqlr{\lambda\big(\natTrans_{C \times B}(\pi_1, \pi_2) \circ 
\prodPres_{C,B} 
\big) \circ \glueFun h, \glueFun p} 
\end{align*}
Moreover, writing
$L := \natTrans_{C \times B}(\pi_1, \pi_2) \circ \prodPres_{C,B}$, 
the 2-cell
required for the diagram below (in which $r : A' \to A$) is the composite 
defined in~(\ref{eq:b-circ-u-x-2-cell}) with 
$f := \lambda L$:
\begin{td}[column sep = 8em, row sep = 3em]
\baseCat(A, C) \times \baseCat(A,B) 
\arrow[swap]{d}
	{\eval_{\glueFun B,X} \circ \seq{\lambda L \circ \glueFun (-), 
		\glueFun ({=})}}
\arrow[phantom]{dr}[description]
	{\twocell
		{\overline{( \wholecell_{B,X}(C) \circ u_{B,X}(C))
			( (\natTrans, \natCell))}_r}}
\arrow{r}{\baseCat(r, C) \times \baseCat(r, B)} &
\baseCat(A', C) \times \baseCat(A', B) 
\arrow{d}
	{\eval_{\glueFun B,X} \circ \seq{\lambda L \circ \glueFun (-), 
		\glueFun ({=})}} \\

\bicatX(\glueFun A, X) 
\arrow[swap]{r}{\bicatX(\glueFun r, X)} &
\bicatX(\glueFun A', X)
\end{td}

We now turn to defining the modification 
$\modif^{(\natTrans, \natCell)}$. 
For $A \in \baseCat$ 
and $(h, p) \in \baseCat(A, C) \times \baseCat(A,B)$ there 
exists an evident choice of isomorphism 
$\modif^{(\natTrans, \natCell)}(A,h,p) : 
	\big( \wholecell_{B,X}(C) \circ u_{B,X}(C)\big)
		\big( (\natTrans, \natCell) \big)(A,h,p)
		\To 
	\natTrans(A,h,p)$, namely
\begin{equation*} 
\makebox[\textwidth]{
\begin{tikzcd}[column sep = 4em, ampersand replacement = \&]
\eval_{\glueFun B, X} \circ 
\seq{\lambda L \circ \glueFun h, \glueFun p}
\arrow{r}{\modif^{(\natTrans, \natCell)}_A(h,p)} 
\arrow[swap]{d}{\iso} \&
\natTrans_A(h,p) \\
\eval_{\glueFun B,X} \circ \seq{\lambda L \circ \glueFun h, \Id_{\glueFun B} 
\circ  \glueFun p} 
\arrow[swap]{d}
{\eval_{\glueFun B, X} \circ \fuse^{-1}}  \&
\natTrans_A(\pi_1\seq{p,q}, \pi_2\seq{p,q}) 
\arrow[swap]{u}{\natTrans_A(\epsilonTimesInd{1}{p,q}, 
\epsilonTimesInd{2}{p,q})}  \\
\eval_{\glueFun B,X} 
	\circ \left((\lambda L \times \glueFun \Id_B ) 
	\circ \seq{\glueFun h, \glueFun p}\right)
\arrow[swap]{d}{\eval_{\glueFun B,X} \circ ( \lambda L \times (\psi^\glueFun 
_B)^{-1} ) \circ 
\seq{\glueFun h, \glueFun p}}  \&
\natTrans_{C \times B}(\pi_1, \pi_2) \circ \glueFun \seq{h, p} 
\arrow[swap]{u}{\natCell_{\seq{p,h}}^{-1}(\pi_1, \pi_2)} \\
\eval_{\glueFun B,X} 
	\circ \left((\lambda L \times \glueFun B) 
	\circ \seq{\glueFun h, \glueFun p}\right)
\arrow[swap]{d}
	{\iso} \&
\left(\natTrans_{C \times B}(\pi_1, \pi_2) 
	\circ \Id_{\glueFun C \times \glueFun B}\right) 
	\circ \glueFun \seq{h, p}  
\arrow[swap]{u}{\iso} \\
\left(\eval_{\glueFun B,X} 
	\circ (\lambda L \times \glueFun B)\right) 
	\circ \seq{\glueFun h, \glueFun p}
\arrow[swap]{d}
	{\epsilonExp_L \circ \seq{\glueFun h, \glueFun p}} \&
\left(\natTrans_{C \times B}(\pi_1, \pi_2) 
	\circ \left(\prodPres_{C,B} 
	\circ \seq{\glueFun \pi_1, \glueFun \pi_2}\right)\right) 
	\circ \glueFun \seq{h, p}
\arrow[swap]{u}
{\natTrans_{C \times B}(\pi_1, \pi_2) \circ \coTimes_{C,B} \circ \glueFun 
\seq{h, p}} \\
\left(\natTrans_{C \times B}(\pi_1, \pi_2) \circ \prodPres_{C,B}\right) 
	\circ \seq{\glueFun h, \glueFun p}
\arrow[swap]{r}[yshift=-2mm]
{\natTrans_{C \times B}(\pi_1, \pi_2) \circ \prodPres_{C,B} \circ 
\unpack^{-1}}  \&
\left(\natTrans_{C \times B}(\pi_1, \pi_2) 
	\circ \prodPres_{C,B}\right) 
	\circ \left(\seq{\glueFun \pi_1, \glueFun \pi_2} 
	\circ \glueFun \seq{h, p}\right)
\arrow[swap]{u}{\iso}
\end{tikzcd}
}
\end{equation*}
It is clear from the definition that 
$\modif^{(\natTrans, \natCell)}_A := \modif^{(\natTrans, \natCell)}(A, -, {=})$ 
is natural in its two arguments and so a 2-cell
${{\big( \wholecell_{B,X}(C) \circ u_{B,X}(C)\big)}
\big( (\natTrans, \natCell)\big)}(A, -, {=})
\To \natTrans(A, -, {=})$ 
in $\Cat$. Moreover, the naturality 
condition~(\ref{eq:modification-natural-in-2-cell}) holds by naturality of each 
of the components defining $\modif^{(\natTrans, \natCell)}$ and the 
modification axiom on 
$\altModif : (\natTrans, \natCell) \to (\altNat, \altCell)$, which requires 
that the following diagram commutes for every $r : A' \to A$ in $\baseCat$ and
$(p, h) \in \baseCat(A, C) \times \baseCat(A, B)$:
\begin{td}
\natTrans_{A'}(pr, hr) 
\arrow[swap]{d}{\altModif_A'(pr,hr)}
\arrow{r}{\natCell_r(p,h)} &
\natTrans_A(p,h) \circ \glueFun r 
\arrow{d}{\altModif_A(p,h) \circ \glueFun r} \\

\altNat_{A'}(pr, hr) 
\arrow[swap]{r}{\altCell_r(p,h)} &
\altNat_A(p,h) \circ \glueFun r
\end{td}
It therefore remains to show that 
the family of 2-cells
$(\modif^{(\natTrans, \natCell)}_A)_{A \in \baseCat}$ satisfies the following 
instance of the
modification axiom 
for every $r : A' \to A$ in $\baseCat$:
\begin{td}[column sep = 5em, row sep = 2.5em]
\big( \wholecell_{B,X}(C) \circ u_{B,X}(C)\big)
\big( (\natTrans, \natCell) \big)(A, pr, hr) 
\arrow{r}{\modif^{(\natTrans, \natCell)}(A, pr, hr)} 
\arrow[swap]{d}
{\overline{( \wholecell_{B,X}(C) \circ u_{B,X}(C))
( (\natTrans, \natCell) ) }_r} &
\natTrans(A, pr, hr) 
\arrow{d}{\natCell_r(p, h)} \\

\big( \wholecell_{B,X}(C) \circ u_{B,X}(C)\big)
\big( (\natTrans, \natCell) \big)(A, p, h) \circ \glueFun r
\arrow[swap]{r}{\modif^{(\natTrans, \natCell)}(A, p, h) \circ \glueFun r} &
\natTrans(A, p, h) \circ \glueFun r
\end{td}
Unfolding the definitions around the anticlockwise composite and applying the 
lemma relating $\fuse$ and $\postName$ 
(Lemma~\ref{lem:PseudoproductCanonical2CellsLaws}),
the problem reduces to the following two lemmas:
\begin{equation} \label{eq:modif-law-lemma-1}
\begin{tikzcd}[column sep = -2em, row sep = 2em]
\: &
\natTrans_A\left(\pi_1\seq{p,h}, \pi_2\seq{p,h}\right) \circ \glueFun r 
\arrow{dr}[yshift=-1mm, xshift=1mm]{\natTrans_A(\epsilonTimesInd{1}{p,h}, 
\epsilonTimesInd{1}{p,h}) \circ 
\glueFun r} &
\: \\
\left(\natTrans_{B\times C}(\pi_1, \pi_2) 
	\circ \glueFun \seq{p, h}\right) 
	\circ \glueFun r 
\arrow{ur}[yshift=-1mm, xshift=-1mm]{\natCell_{\seq{p,h}}^{-1}(\pi_1, \pi_2) 
\circ \glueFun r} 
\arrow[swap]{d}{\iso} &
\: &
\natTrans_A(p,h) \circ \glueFun r  \\
\natTrans_{B\times C}(\pi_1, \pi_2) 
	\circ \left(\glueFun \seq{p, h} 
	\circ \glueFun r\right)
\arrow[swap]{d}
{\natTrans_{B\times C}(\pi_1, \pi_2) \circ \phi^\glueFun _{\seq{p, h}, r}}  &
\: &
\: \\
\natTrans_{B\times C}(\pi_1, \pi_2) \circ \glueFun (\seq{p, h} \circ r) 
\arrow[swap]{d}
{\natTrans_{B\times C}(\pi_1, \pi_2) \circ \glueFun \postName} &
\: &
\natTrans_{A'}(pr, hr) 
\arrow[swap]{uu}{\natCell_r(p,h)} \\
\natTrans_{B\times C}(\pi_1, \pi_2) \circ \glueFun \seq{pr, hr}
\arrow[swap]{rr}{\natCell_{\seq{pr, hr}}^{-1}(\pi_1, \pi_2)} &
\: &
\natTrans_{A'}\left(\pi_1\seq{pr, hr}, \pi_2\seq{pr, hr}\right) 
\arrow[swap]{u}
{\natTrans_{A'}(\epsilonTimesInd{1}{pr,hr}, \epsilonTimesInd{2}{pr,hr})}
\end{tikzcd}
\end{equation}
and
\begin{equation} \label{eq:modif-law-lemma-2}
\makebox[\textwidth]{
\begin{tikzcd}[column sep = 0.5em, row sep = 2em, ampersand replacement = \&]
\: \&
\prodPres_{C,B} 
	\circ \left(\left(\seqlr{\glueFun \pi_1, \glueFun \pi_2} 
	\circ \glueFun\seq{p, h}\right)
	\circ \glueFun r\right) 
\arrow[swap]{dl}{\prodPres_{C,B} \circ \unpack \circ \glueFun r}
\arrow{dr}{\iso} \&
\: \\
\prodPres_{C,B} \circ \seqlr{\glueFun p, \glueFun h} 
	\circ \glueFun r
\arrow[swap]{d}{\prodPres_{C,B} \circ \postName} \&
\: \&
\glueFun \seq{p, h} \circ \glueFun r 
\arrow{d}{\phi^\glueFun _{\seq{p,h}, r}}  \\
\prodPres_{C,B} \circ 
	\seqlr{\glueFun p \circ \glueFun r, \glueFun h \circ \glueFun r}
\arrow[swap]{d}{\prodPres_{C,B} \circ 
	\seqlr{\phi^\glueFun _{p,r}, \phi^\glueFun_{h,r}}} \&
\: \&
\glueFun (\seq{p, h} \circ r)
\arrow{d}{\glueFun \postName}  \\
\prodPres_{C,B} \circ \seqlr{\glueFun (pr), \glueFun (hr)}
\arrow[swap]{d}{\prodPres_{C,B} \circ \unpack^{-1}} \&
\: \&
\glueFun \seq{pr, hr} \\
\prodPres_{C,B} 
	\circ \left(\seqlr{\glueFun \pi_1, \glueFun \pi_2} 
	\circ \glueFun\seq{pr, hr}\right) 
\arrow[swap]{rr}{\iso} \&
\: \&
\left(\prodPres_{C,B} 
	\circ \seqlr{\glueFun \pi_1, \glueFun \pi_2}\right) 
	\circ \glueFun \seq{pr, hr}
\arrow[swap]{u}{\coTimes_{C,B} \circ h}
\end{tikzcd}
\vspace{6mm}
}
\end{equation}
Here the top unlabelled isomorphism is the composite
\begin{td}[column sep = 7em]
\prodPres_{C,B} 
	\circ \left(\left(\seqlr{\glueFun \pi_1, \glueFun \pi_2} 
	\circ \glueFun\seq{p, h}\right)
	\circ \glueFun r\right) 
\arrow[swap]{d}{\iso}
\arrow{r} &
\glueFun\seq{p, h}
	\circ \glueFun r \\

\left(\prodPres_{C,B} 
	\circ \seqlr{\glueFun \pi_1, \glueFun \pi_2}\right) 
	\circ \left(\glueFun\seq{p, h}
	\circ \glueFun r\right)
\arrow{r}[swap, yshift=-0mm]
	{\coTimes_{C,B} 
	\circ \glueFun\seq{p, h}
	\circ \glueFun r} &
\Id_{\glueFun(B \times C)}
	\circ \left(\glueFun\seq{p, h}
	\circ \glueFun r\right)
\arrow[swap]{u}{\iso}
\end{td}
applying the isomorphism $\coTimes_{C,B}$ witnessing that 
$\prodPres_{C,B} : \glueFun C \times \glueFun B 
	\leftrightarrows \glueFun(C \times B) : 
	\seqlr{\glueFun \pi_1, \glueFun \pi_2}$ 
forms an equivalence.
 
For~(\ref{eq:modif-law-lemma-1}), one applies the 
associativity law for $(\natTrans, \natCell)$ along with the 
definition of $\postName$ as part of a short diagram chase. 
For~(\ref{eq:modif-law-lemma-2}), one unwinds the 
definition of $\unpack$ in each of the two given composites and repeatedly 
applies naturality.
\end{proof}
\end{mylemma}

This lemma, together with Lemma~\ref{lem:pullback-from-equivalence}, completes 
the proof of Proposition~\ref{prop:explicit-glued-yon-exp}. 

\section{Glueing syntax and semantics}
\label{sec:nbe-glueing}

Our aim now is to show how the structure of $\langCartClosed$, together with 
the identification of neutral and normal terms in 
Section~\ref{sec:syntax-as-pseudofunctors}, determines data 
in the bicategory of intensional Kripke 
relations~(\cf~(\ref{eq:cat-glued-data}) on page~\pageref{eq:cat-glued-data}). 
Fix a cc-bicategory 
$\ccBicat\bicatX$ and consider an interpretation $\baseTypes \to \bicatX$ 
of base types in $\bicatX$ with canonical extension 
$s : \allTypes\baseTypes \to \bicatX$. We show that the terms of 
$\langCartClosed$ determine 
objects in the 
glueing bicategory, and that the typing rules determine 1-cells.

\paragraph*{From terms to glued objects.} On neutral and normal terms, the key 
observation is that the 
interpretation of $\langCartClosed$-terms in $\bicatX$ is pseudonatural.

\begin{myconstr}
Let $\baseTypes$ be a set of base types, 
$\ccBicat{\bicatX}$ be a cc-bicategory, and
$s : \allTypes{\baseTypes} \to \bicatX$ the canonical extension of a set map
$\baseTypes \to \bicatX$. By Proposition~\ref{propn:weak-universal-property} 
there exists a \mbox{cc-pseudofunctor} 
$s\sem{-} : \syncloneAtClosed{\allTypes\baseTypes} \to \bicatX$ 
interpreting $\langCartClosed(\allTypes\baseTypes)$ in $\bicatX$ (see 
Construction~\ref{constr:interpretation-of-langCartClosed} for the full definition). 
We define a pseudonatural transformation
$(s\sem{-}, \cellOf{s\sem{-}}) : \twoDisc{\dlangTerms(-;A)} \To 
			\bicatX\big( s\sem{-},  s\sem{A}\big)
			: 
			\twoDisc{{\Con_{\allTypes{\baseTypes}}}} \to \Cat$ 
for every $A \in \allTypes\baseTypes$. 

For the component at $\Gamma \in \Con_{\allTypes\baseTypes}$ we take 
the functor
\begin{align*}
\twoDisc{\dlangTerms(\Gamma; A)} 
	&\xra{s\sem{-}_{\Gamma, A}} 
	\bicatX(s\sem{\Gamma},  s\sem{A}) \\
\into{t} &\mapsto s\sem{\Gamma \vdash \into{t} : A}
\end{align*}
Next, for every context renaming 
$r : \Gamma \to \Delta$ 
we need to provide a 2-cell---\ie~natural isomorphism---as in
\begin{td}[column sep = 6em]
\twoDisc{\dlangTerms(\Gamma; A)}
\arrow{r}{\twoDisc{\dlangTerms(r; A)}}
\arrow[phantom]{dr}[description]{\twocell{(\cellOf{s\sem{-}})_r}}
\arrow[swap]{d}{s\sem{-}} &
\twoDisc{\dlangTerms(\Delta; A)}
\arrow{d}{s\sem{-}} \\

\bicatX(s\sem{\Gamma},  s\sem{A}) 
\arrow[swap]{r}
	{\bicatX(s\sem{r},  s\sem{A})} &
\bicatX(s\sem{\Delta},  s\sem{A})
\end{td}
Thus, for every $\into{t} \in \langTerms(\Gamma; A)$ we need to provide an 
isomorphism in $\bicatX$ of type
${s\sem{\Delta \vdash \into{t[r(x_i) / x_i]} : A}} \to 
	{s\sem{\Gamma \vdash \into{t} : A}} \circ 
		s\sem{r}$. 
Calculating, one sees that
\newpage
\begin{align*}
s\semlr{\into{\Gamma \vdash \into{t} : A}} \circ s\sem{r} 
	&= s\semlr{\into{\Gamma \vdash \into{t} : A}} 
		\circ \seqlr{\pi_{r(1)}, \,\dots\, , \pi_{r(n)}} \\
	&= s\semlr{\into{\Gamma \vdash \into{t} : A}} 
		\circ \seqlr{s\sem{(\Delta \vdash x_{r(i)} : A_{r(i)}})}_i 
	\\
	&= s\semlr{\into{\Gamma \vdash \into{t} : A}} \circ 
		s\semlr{(\Delta \vdash x_{r(i)} : A_{r(i)})_{i=1, \,\dots\, , n}} \\
	&= s\semlr{\Delta \vdash \hcomp{\into{t}}{r} : A }
\end{align*}
Now recall from Construction~\ref{constr:def-of-cont} that we have already 
constructed a rewrite typed by the rule
\begin{center}
\binaryRule{\Gamma \vdash \into{t} : A}
		{r : \Gamma \to \Delta}
		{\Delta \vdash \cont{t}{r} : 
			\rewrite{\hcomp{\into{t}}{x_i \mapsto 
				r(x_i)}}
			{\into{t[r(x_i)/x_i]}} : A}{} \vspace{-\treeskip}
\end{center}
We therefore define $(\cellOf{s\sem{-}})_r$ to be the interpretation of 
$\contName$:
\[
(\cellOf{s\sem{-}})_r(t) := 
s\sem{ \Delta \vdash \cont{t}{r} : 
			\rewrite{\hcomp{\into{t}}{x_i \mapsto 
				r(x_i)}}
			{\into{t[r(x_i)/x_i]}} : A }
\]
To see that this is a pseudonatural transformation, observe first that it is 
certainly natural: there are no non-trivial 2-cells in 
$\twoDisc{\langTerms(\Gamma;A)}$. 
For the unit law, we need to show that
\begin{equation} \label{eq:cont-name-pseudonat-unit-law}
\begin{tikzcd}[column sep = 7em]
s\semlr{\Gamma \vdash \into{t} : A} \circ \Id_{s\semlr{\Gamma}}
\arrow[swap]{d}
	{s\semlr{\Gamma \vdash \into{t} : A} \circ 
	\widehat{\etaTimes{}}_{\Id_{s\semlr{\Gamma}}}}
\arrow{r}{\iso} &
s\semlr{\Delta \vdash \into{t[x_i/x_i]} : A}
\arrow[equals]{d} \\
s\semlr{\Gamma \vdash \into{t} : A} \circ \seq{\pi_1, \,\dots\, , \pi_n}
\arrow[swap]{r}[yshift=-2mm]{s\semlr{\Gamma \vdash \cont{t}{\id_\Gamma} : 
	\rewrite{\hcomp{t}{x_i \mapsto x_i}}{\into{t[x_i/x_i]}} : A}} &
s\semlr{\Delta \vdash \into{t} : A}
\end{tikzcd}
\end{equation}
where 
$ \widehat{\etaTimes{}}_{\Id_{s\sem{\Gamma}}}
	:= \Id_{s\sem{\Gamma}}  
		\XRA{\etaTimes{\Id_{s\sem{\Gamma}}}}
		\seqlr{\pi_1 \circ \Id_{s\sem{\Gamma}}, \,\dots\, , \pi_n \circ 
			\Id_{s\sem{\Gamma}}}
		\XRA\iso
		\seq{\pi_1, \,\dots\, , \pi_n}$.
To see this commutes, note that
$s\sem{\Gamma \vdash \subid{\into{t}} : 
	\rewrite{\into{t}}{\hcomp{\into{t}}{x_i \mapsto x_i}} : A}$ is, by 
definition, the composite
\[
s\semlr{\Gamma \vdash \into{t} : A} 
\XRA{\iso}
s\semlr{\Gamma \vdash \into{t} : A} \circ \Id_{s\sem{\Gamma}} 
\XRA{s\sem{\Gamma \vdash \into{t} : A} \circ 
	\widehat{\etaTimes{}}_{\Id_{s\sem{\Gamma}}}}
s\sem{\Gamma \vdash \into{t} : A} \circ \seq{\pi_1, \,\dots\, , \pi_n}
\]
Hence~(\ref{eq:cont-name-pseudonat-unit-law}) commutes by 
Lemma~\ref{lem:out-cont-coincide} and
Lemma~\ref{lem:properties-of-cont-and-sub}(\ref{c:unit-law}). 

For the associativity law we need to show that, for any contexts
$\Gamma := (x_i : A_i)_{i=1, \,\dots\, , n}$
and 
$\Delta := (y_j : A_j)_{j=1, \,\dots\, , m}$,
and any context renamings 
$\Gamma \xra{r} \Delta \xra{r'} \Sigma$, the following diagram 
commutes:
\begin{equation*} \label{eq:cont-name-pseudonat-assoc-law}
\begin{tikzcd}[column sep = -1.5em, row sep = 3em]
s\semlr{\Gamma \vdash \into{t} : A} 
	\circ \left(\seq{\pi_{r}} 
	\circ \seq{\pi_{r'}}\right) 
\arrow{rr}{s\semlr{\Gamma \vdash \into{t} : A} \circ \postName} &
\: &
s\semlr{\Gamma \vdash \into{t} : A} \circ \seqlr{\pi_{r} \circ \seq{\pi_{r'}}} 
\arrow{d}{s\semlr{\Gamma \vdash \into{t} : A} \circ 
	\seq{\epsilonTimesInd{r}{}}} \\
\left(s\semlr{\Gamma \vdash \into{t} : A} \circ \seq{\pi_{r}}\right) 
	\circ \seq{\pi_{r'}} 
\arrow[]{u}{\iso}
\arrow{d}[swap]{s\semlr{\cont{\into{t}}{r}} \circ \seq{\pi_{r'}}} &
\: &
s\semlr{\Gamma \vdash \into{t} : A} \circ \seq{\pi_{r'r}} 
\arrow{d}{s\semlr{\cont{\into{t}}{r'r}}} \\
s\semlr{\Delta \vdash \into{t[r(x_i)/x_i]} : A} \circ \seq{\pi_{r'}} 
\arrow{dr}[swap]{s\semlr{\cont{\into{t[r(x_i) /x_i]}}{r'}}} &
\: &
s\semlr{\Sigma \vdash \into{t[r'r(x_i)/x_i]} : A} \\
\: &
s\semlr{\Sigma \vdash \into{t[r(x_i)/x_i][r'(y_j)/y_j]} : A} 
\arrow[equals]{ur} &
\:
\end{tikzcd}
\end{equation*}
We suppress the full typing judgement in the vertical arrows for 
reasons of space. By 
Lemma~\ref{lem:out-cont-coincide}, this diagram is 
exactly the image of 
Lemma~\ref{lem:properties-of-cont-and-sub}(\ref{c:assoc-law}) 
under $s\sem{-}$, and so it commutes.
\end{myconstr}

The preceding construction restricts to 
neutral and normal terms, giving pseudonatural transformations
\begin{align*}
\twoDisc{\dneutTerms(-; A)} \XRA{\restr{(s\sem{-}, \cellOf{s\sem{-}})}{M}} 
	\bicatX\big( s\sem{-},  s\sem{A} \big) \\
\twoDisc{\dnormTerms(-; A)} \XRA{\restr{(s\sem{-}, \cellOf{s\sem{-}})}{N}} 
	\bicatX\big( s\sem{-},  s\sem{A} \big)
\end{align*}
One thereby obtains the following glued objects for every type 
$A \in \allTypes\baseTypes$:
\begin{equation} \label{eq:dneut-and-dnorm-defined}
\begin{aligned}
\mu_A &:= \big(\twoDisc{\dneutTerms(-; A)} , 
				\restr{(s\sem{-}, \cellOf{s\sem{-}})}{M}, 
				s\sem{A}\big) \\
\eta_A &:= \big(\twoDisc{\dnormTerms(-; A)}, 
				\restr{(s\sem{-}, \cellOf{s\sem{-}})}{N}, 
				s\sem{A}\big)
\end{aligned}
\end{equation}
Finally, for variables, we take
\[
\nu_A := \glued{\Yon}(\coerce{A}) = 
\big(\twoDisc{\Con_{\allTypes\baseTypes}(-;A)} , 
			(\lanNat, \lanCell)_{(-, A)} , 
			s\sem{A} \big) 
\]
where $(\lanNat, \lanCell)_{(-, A)}$ is the pseudonatural transformation of 
Corollary~\ref{cor:lanext-transformation-construction}.

\paragraph*{From typing rules to glued 1-cells.}
We also lift the natural transformations 
of~(\ref{eq:typing-rules-as-pseudonat-trans})---viewed as locally discrete 
pseudonatural transformations---to morphisms in 
$\gl{\lanext{\prodext{s}}}$. 

For the lambda abstraction case we will use the following observation. 
For 
types $A, B \in \allTypes\baseTypes$ the exponential
$\altexp
	{\twoDisc{\dvarTerms(-;A)}}
	{\twoDisc{\dnormTerms(-;B)}} = 
\altexp
		{\twoDisc{(\yon {\coerce{A}})}}
		{\twoDisc{\dnormTerms(-;B)}} =
\altexp
		{\Yon \coerce{A}}
		{\twoDisc{\dnormTerms(-;B)}}$
in $\Hom(\twoDisc{\Con_{\allTypes\baseTypes}}, \Cat)$ is, by 
Theorem~\ref{thm:exponentiating-by-Yoneda-identified}, equivalent to 
$\twoDisc{\dnormTerms(- \concat \coerce{A};B)}$. One thereby obtains a composite
\begin{equation} \label{eq:exponential-context-ext}
\altexp
	{\twoDisc{\dvarTerms(-;A)}}
	{\twoDisc{\dnormTerms(-;B)}} \xra{\simeq}
\twoDisc{\dnormTerms(- \concat \coerce{A};B)} 
\xra{\twoDisc{\mathsf{lam}(-;A,B)}}
\twoDisc{\dnormTerms(- ;\exptype{A}{B})}
\end{equation}
We put this to work in the next result, which is the bicategorical 
version of Fiore's~\cite[Proposition 7 and Proposition 8]{Fiore2002}.

\begin{myremark}
Examining the equivalence 
$\altexp
	{\twoDisc{\dvarTerms(-;A)}}
	{\twoDisc{\dnormTerms(-;B)}} \simeq
\twoDisc{\dnormTerms(- \concat \coerce{A};B)}$, one sees 
that it is in fact an isomorphism. Since 
$\dnormTerms(\Gamma \concat \coerce{A};B)$ is a set for every context 
$\Gamma$, 
the composite 
$\dnormTerms(\Gamma \concat \coerce{A};B) 
	\to 
\altexp
	{\twoDisc{\dvarTerms(-;A)}}
	{\twoDisc{\dnormTerms(-;B)}}(\Gamma)
	\to 
\dnormTerms(\Gamma \concat \coerce{A};B)$ 
must be equal to the identity. On the other 
hand, by Lemma~\ref{lem:cat-to-2-cat}(\ref{c:presheaves-to-pseudofunctors}), 
the exponential
$\altexp
	{\twoDisc{\dvarTerms(-;A)}}
	{\twoDisc{\dnormTerms(-;B)}}$
may be given by 
$\twoDisc{\left(\funCat{\catC}{\Set}
	\left( \yon (-) \times \dvarTerms(=;A), 
			\dnormTerms(=;B)\right)\right)}$.
But   
$\Cat(\twoDisc\catC, \Set)
	\left( \yon \Gamma \times \dvarTerms(=;A), 
			\dnormTerms(=;B)\right)$ 
is also a set for every context $\Gamma$. Hence, the composite 
$\altexp
	{\twoDisc{\dvarTerms(-;A)}}
	{\twoDisc{\dnormTerms(-;B)}} \to 
\altexp
	{\twoDisc{\dvarTerms(-;A)}}
	{\twoDisc{\dnormTerms(-;B)}}
$ 
must also be the identity.
\end{myremark}

\newpage
\begin{mypropn} \label{prop:1-cells-in-glueing-bicat}
For every set of base types $\baseTypes$, 
cc-bicategory $\ccBicat{\bicatX}$, 
and set map
$s : \allTypes\baseTypes \to \bicatX$
canonically induced from an interpretation of base types
$\baseTypes \to \bicatX$, 
\begin{enumerate}
\item \label{c:map-in-gl-bicat-var} For every type 
$A_i \in\allTypes\baseTypes$, 
the triple 
$\glued{\mathsf{var}} := \big(\twoDisc{\mathsf{var}(-; A_i)}, \iso, 
\Id_{s\sem{A_i}}\big)$ 
is a 1-cell 
$\nu_{A_i} \to \mu_{A_i}$ in 
$\gl{\lanext{\prodext{s}}}$, where the 2-cell $\iso$ filling
\begin{td}[column sep = 4em]
\twoDisc{\dvarTerms(-;A_i)}
\arrow{r}{\twoDisc{\mathsf{var}(-; A_i)}}
\arrow[swap]{d}{s\sem{-}}
\arrow[phantom]{dr}[description]{\twocell{\iso}} &
\twoDisc{\dneutTerms(-; A_i)} 
\arrow{d}{s\sem{-}} \\

\bicatX\big( s\sem{-},  s\sem{A_i} \big) 
\arrow[swap]{r}[yshift=-2mm]
	{\bicatX\big( s\sem{-},  \Id_{s\sem{A_i}} \big)} &
\bicatX\big( s\sem{-},  s\sem{A_i} \big)
\end{td}
is the structural isomorphism 
$
s\sem{\Gamma \vdash x_i : A_i} 
\XRA{\iso}
\Id_{s\sem{A_i}} \circ s\sem{\Gamma \vdash x_i : A_i} 
$.
\item \label{c:map-in-gl-bicat-base-types} For any base type 
$B \in \baseTypes$, the 
triple
$\glued{\mathsf{inc}} :=
\left(\mathsf{inc}(-; B), \iso, \Id_{s\sem{B}} \right)$, in 
which $\iso$ is a 
structural isomorphism, is an isomorphism
$\mu_B \xra{\iso} \eta_B$ in $\gl{\lanext{\prodext{s}}}$.
\item \label{c:map-in-gl-bicat-proj} For every sequence of types 
$A_1, \,\dots\, , A_n \in\allTypes\baseTypes~~(n \in \Nat)$, the triple
$\glued{\mathsf{proj}}_k := 
\big(\twoDisc{\mathsf{proj}_k(-;\ind{A})}, \id, \pi_k \big)$ is 
a 1-cell $\mu_{\prod_n(A_1, \,\dots\, , A_n)} \to \mu_{A_k}$ 
in $\gl{\lanext{\prodext{s}}}$
for $k=1, \,\dots\, , n$.
\item \label{c:map-in-gl-bicat-app} For every pair of types 
$A, B \in\allTypes\baseTypes$, the triple
$\glued{\mathsf{app}} := 
	\big( 
		\twoDisc{\mathsf{app}}(-; A, B), \id, \eval_{s\sem{A}, s\sem{B}} 
	\big)$
is a 1-cell 
${\mu_{\scriptsizeexpobj{A}{B}} \times \eta_{A}} \to \mu_B$ in 
$\gl{\lanext{\prodext{s}}}$.
\item \label{c:map-in-gl-bicat-tupe} For every sequence of types 
$A_1, \,\dots\, , A_n \in\allTypes\baseTypes ~(n \in \Nat)$, the triple \newline
$\glued{\mathsf{tuple}} :=
{\big(\twoDisc{\mathsf{tuple}(-;\ind{A})}, \iso, 
\Id_{s\sem{\prod_n \ind{A}}} \big)}$ 
is 
a 1-cell 
$\prod_{i=1}^n \eta_{A_i} \to \eta_{\prod_n(A_1, \,\dots\, , A_n)}$
in $\gl{\lanext{\prodext{s}}}$, where the 
isomorphism filling
\begin{td}[column sep = 6em]
\prod_{i=1}^n \twoDisc{\dnormTerms(-; A_i)}
\arrow{r}{\twoDisc{\mathsf{tuple}(-;\ind{A})}}
\arrow[phantom]{ddr}[description]{\twocell{\iso}}
\arrow[swap]{d}{\prod_{i=1}^n s\sem{-}} &
\twoDisc{\dnormTerms\big(-; \prodop_n(A_1, \,\dots\, , A_n)\big)}
\arrow{dd}{s\sem{-}} \\

\prod_{i=1}^n \bicatX\big( s\sem{-}, s\sem{A_i}\big) 
\arrow[swap]{d}{\seq{-, \,\dots\, , =}} &
\: \\

\bicatX{\left(s\sem{-}, \prodop_{i=1}^n s\sem{A_i}\right)} 
\arrow[swap]{r}[yshift=-2mm]
	{\bicatX{\left( s\sem{-}, \Id_{s\sem{\prod_n \ind{A}}} \right)}} &
\bicatX{\left( s\sem{-}, 
	s\sem{{\prodop_n(A_1, \,\dots\, , A_n)}} \right)} 
\end{td}
is the structural isomorphism
\[
s\sem{\Gamma \vdash \pair{\into{t_1}, \,\dots\, , \into{t_n}} : \prodop_n 
\ind{A}} 
= \seq{s\sem{\Gamma \vdash \into{\ind{t}} : \ind{A} }}
\XRA{\iso}
\Id_{(\prod_i sA_i)} \circ \seq{s\sem{\Gamma \vdash \into{\ind{t}} : \ind{A} }} 
\] 
\enlargethispage*{3\baselineskip}
\item \label{c:map-in-gl-bicat-lam} For any pair of types 
$A, B \in \allTypes{\baseTypes}$, write
$\mathsf{L}_{A,B}$ for the 
composite 
\[
\altexp
	{\twoDisc{\dvarTerms(-;A)}}
	{\twoDisc{\dnormTerms(-;B)}} 
\xra{\simeq}
\twoDisc{\dnormTerms\big(- + \coerce{A}, B \big)}
\xra{\twoDisc{\mathsf{lam}(-; A, B)}}
\twoDisc{\dnormTerms\big(- , \exptype{A}{B} \big)}
\]
of~(\ref{eq:exponential-context-ext}). 
Then, where $\iso$ denotes a structural isomorphism, 
$\glued{\mathsf{lam}} := 
	(\mathsf{L}_{A,B}, \iso,\Id_{\scriptsizeexpobj{s\sem{A}}{s\sem{B}}} )$
is a 1-cell
${(\exp{\nu_A}{\eta_B})} \xra{\iso} \eta_{\scriptsizeexpobj{A}{B}}$
in $\gl{\lanext{\prodext{s}}}$.
\end{enumerate} 
\begin{proof}
(\ref{c:map-in-gl-bicat-var}) is immediate. 
For~(\ref{c:map-in-gl-bicat-base-types}), observe first that the only 
way to 
construct normal terms of base type is via the \rulename{inc} rule. Hence the 
natural transformation $\mathsf{inc}$ is a natural isomorphism. Next consider 
the diagram
\begin{td}[column sep = 6em]
\twoDisc{\dneutTerms(-; B)} 
\arrow[phantom]{dr}[description]{\twocell{\iso}}
\arrow[swap]{d}{s\sem{-}}
\arrow{r}{\mathsf{inc}{(-;B)}} &
\twoDisc{\dnormTerms(-;B)} 
\arrow{d}{s\sem{-}} \\

\bicatX\big( s\sem{-}, s\sem{B} \big) 
\arrow[swap]{r}{\bicatX( s\sem{-}, \Id_{s\sem{B}} )} &
\bicatX\big( s\sem{-}, s\sem{B} \big)
\end{td}
For a context $\Gamma$ and term $t \in \dneutTerms(\Gamma; B)$, the 
clockwise 
route returns $s\sem{\Gamma \vdash t : B}$ while the anticlockwise route 
returns $\Id_{s\sem{B}} \circ s\sem{\Gamma \vdash t : B}$. Hence 
the diagram is 
filled by a structural isomorphism, and 
$\big( \mathsf{inc}(-;B), \iso, \Id_{s\sem{B}}\big)$ is a 1-cell 
in
$\gl{\lanext{\prodext{s}{}}}$. To see that it is an isomorphism in 
$\gl{\lanext{\prodext{s}{}}}$, observe 
that the diagram 
\begin{td}[column sep = 6em]
\twoDisc{\dnormTerms(-; B)} 
\arrow[phantom]{dr}[description]{\twocell{\iso}}
\arrow[swap]{d}{s\sem{-}}
\arrow{r}{\mathsf{inc}{(-;B)}^{-1}} &
\twoDisc{\dneutTerms(-;B)} 
\arrow{d}{s\sem{-}} \\

\bicatX\big( s\sem{-}, s\sem{B} \big) 
\arrow[swap]{r}{\bicatX( s\sem{-}, \Id_{s\sem{B}} )} &
\bicatX\big( s\sem{-}, s\sem{B} \big)
\end{td}
is also filled by a structural isomorphism, giving a 1-cell
$\big( \mathsf{inc}(-;B)^{-1}, \iso, \Id_{s\sem{B}}\big)$. Then, by 
the coherence 
theorem for bicategories, the composite   
\begin{td}[column sep = 5em, row sep = 3em]
\twoDisc{\dneutTerms(-; B)} 
\arrow[bend left = 20]{rr}{\Id_{\twoDisc{\dneutTerms(-; B)} }}
\arrow[bend left = 10, phantom]{rr}[description]{=}
\arrow[phantom]{dr}[description]{\twocell{\iso}}
\arrow[swap]{d}{s\sem{-}}
\arrow{r}{\mathsf{inc}{(-;B)}} &
\twoDisc{\dnormTerms(-;B)} 
\arrow{d}[description]{s\sem{-}}
\arrow[phantom]{dr}[description]{\twocell{\iso}} 
\arrow{r}{\mathsf{inc}{(-;B)}^{-1}} &
\twoDisc{\dneutTerms(-; B)} 
\arrow{d}{s\sem{-}}  \\

\bicatX\big( s\sem{-}, s\sem{B} \big) 
\arrow[bend right = 12, phantom]{rr}[description]{\twocellDown{\iso}}
\arrow[bend right = 20, swap]{rr}{\bicatX( s\sem{-}, 
\Id_{s\sem{B}} )}
\arrow[swap]{r}{\bicatX( s\sem{-}, \Id_{s\sem{B}} )} &
\bicatX\big( s\sem{-}, s\sem{B} \big) 
\arrow[swap]{r}{\bicatX( s\sem{-}, \Id_{s\sem{B}} )} &
\bicatX\big( s\sem{-}, s\sem{B} \big) 
\end{td}
is equal to the identity 1-cell 
$\Id_{\mu_B}$ in $\gl{\lanext{\prodext{s}{}}}$, and similarly for the other 
composite. 

For~(\ref{c:map-in-gl-bicat-proj}) one needs to check that the following 
diagram commutes on the nose:
\begin{td}[column sep = 4em]
\twoDisc{\dneutTerms\left(-; \prodop_n(A_1,\dots, A_n) \right)}
\arrow[swap]{d}{s\sem{-}}
\arrow{r}{\twoDisc{\mathsf{proj}_k(-;\ind{A})}} &
\twoDisc{\dneutTerms(-; A_k)} 
\arrow{d}{s\sem{-}} \\

\bicatX\big( s\sem{-},  
	s\sem{\prodop_n(A_1, \,\dots\, , A_n)} \big) 
\arrow[swap]{r}
	{\bicatX( s\sem{-}, \pi_k)} &
\bicatX\big( s\sem{-},  
	s\sem{A_k} \big) 
\end{td}
For a fixed context $\Gamma$ and term 
$\into{t} \in \dneutTerms(\Gamma;B)$, 
\[
s\sem{\mathsf{proj}_k(\Gamma;\ind{A})(t)} = 
	s\sem{ \into{\pi_k(t)} } =
	s\sem{\hpi{k}{\into{t}}} =
	\pi_k \circ s\sem{\Gamma \vdash \into{t} : \prodop_n(A_1, \,\dots\, ,A_n)}
\]	
as required.

For~(\ref{c:map-in-gl-bicat-app}) one observes that
the product
$\mu_{\scriptsizeexpobj{A}{B}} \times \eta_A$ 
in $\gl{\lanext{\prodext{s}{}}}$
is the pseudonatural transformation $\kappa_{A,B}$ defined by the diagram below.
\begin{td}[column sep = -1em]
\: &
\bicatX\big( s\sem{-}, 
s\sem{\expobj{A}{B}} \big) \times 
\bicatX\big( s\sem{-},  s\sem{A} \big) 
\arrow{dr}{\seq{-, =}}  &
\: \\

\twoDisc{\dneutTerms(-; \exptype{A}{B})} \times
	\twoDisc{\dnormTerms(-; A)} 
\arrow{ur}{s\sem{-} \times s\sem{-}}
\arrow[swap]{rr}{\kappa_{A,B}} &
\: &
\bicatX\big( s\sem{-},  
	s\sem{\expobj{A}{B}} 
		\times s\sem{A} \big) 
\end{td}
Hence, the composite
$\bicatX\big(s\sem{-}, \eval_{sA, sB}\big) \circ  \kappa_{A,B}$
instantiated at a context $\Gamma$ and a pair of terms
$(\into{t}, \into{u})$ returns
\begin{align*}
	\eval_{sA, sB} \circ 
		\seqlr{s\semlr{\Gamma \vdash \into{t} : \exptype{A}{B}}, 
			s\semlr{\Gamma \vdash \into{u} : A}}
 &= s\semlr{\heval{\into{t}, \into{u}}} \\
 &= s\semlr{\twoDisc{\mathsf{app}}(\Gamma; A, B)(\into{t}, \into{u})}
\end{align*} as required. The calculation for~(\ref{c:map-in-gl-bicat-tupe}) is 
similar.

For~(\ref{c:map-in-gl-bicat-lam}) some calculations are required.
Since
$\nu_A = \Yon\coerce{A}$, the exponential $\exp{\nu_A}{\eta_B}$ may, by 
Proposition~\ref{prop:explicit-glued-yon-exp}, be given by the composite
\[
\altexp{\Yon\coerce{A}}{\twoDisc{\dnormTerms(-;B)}}
	\xra{[\Yon \coerce{A}, (s\sem{-}, \cellOf{s\sem{-}})]} 
\altexp{\Yon \coerce{A}}{\bicatX(s\sem{-}, s\sem{B})} 
	\xra{u_{\coerce{A}, s\sem{B}}}
\bicatX\big( s\sem{-}, \exp{s\sem{A}}{s\sem{B}}
\big)
\]
We therefore calculate the two routes around the diagram 
\begin{td}[column sep = 4em]
\altexp{\Yon\coerce{A}}{\twoDisc{\dnormTerms(-;B)}}
\arrow{r}{\simeq}
\arrow[swap]{d}
	{[\Yon \coerce{A}, (s\sem{-}, \cellOf{s\sem{-}})]} &
\twoDisc{\dnormTerms(- + \coerce{A};B)} 
\arrow{r}{\twoDisc{\mathsf{lam}(-;A,B)}} &
\twoDisc{\dnormTerms(- ;\exptype{A}{B})} 
\arrow{dd}{s\sem{-}} \\

\altexp{\Yon\coerce{A}}
	{\bicatX\big( s\sem{-}, s\sem{B} \big)} 
\arrow[swap]{d}{u_{\coerce{A}, s\sem{B}}} &
\: &
\: \\

\bicatX\big( s\sem{-}, 
\exp{s\sem{A}}{s\sem{B}} \big)
\arrow[swap]{rr}{\bicatX( s\sem{-}, 
\Id_{\scriptsizescriptsizeexpobj{s\sem{A}}{s\sem{B}}} )} &
\: &
\bicatX\big( s\sem{-}, 
\exp{s\sem{A}}{s\sem{B}} \big) &
\end{td}
We begin with the anticlockwise route, instantiated at a context $\Gamma$. 
For 
$(\altNat, \altCell) : 
	\Yon \Gamma \times \Yon\coerce{A} \To \twoDisc{\dnormTerms(-;B)}$
the pseudonatural transformation 
$[\Yon \coerce{A}, (s\sem{-}, \cellOf{s\sem{-}})]
	(\altNat, \altCell)$
is simply the composite 
\begin{equation} \label{eq:action-as-composition}
\Yon\Gamma \times \Yon\coerce{A} 
	\XRA{(\altNat, \altCell)}
\twoDisc{\dnormTerms(-;B)} 
	\XRA{(s\sem{-}, \cellOf{s\sem{-}})}
\bicatX\big( s\sem{-}, s\sem{B} \big)
\end{equation}
Moreover, from~(\ref{eq:u-identified}) on page~\pageref{eq:u-identified} we 
know that, at $\Gamma$, the equivalence 
$u_{s\sem{A}, s\sem{B}}$ 
takes a 
pseudonatural transformation 
$(\natTrans, \natCell) : \Yon\Gamma \times \Yon \coerce{A} \To 
		\bicatX(s\sem{-}, s\sem{B})$
to the 1-cell 
\[
\lambda \big( s\sem{\Gamma} \times s\sem{A}
				\xra{\prodPres_{\Gamma, \coerce{A}}} 
				s\sem{\Gamma \concat \coerce{A}}
				\xra{\natTrans_{\Gamma \concat \coerce{A}}(\iota_1, \iota_2) }
				 s\sem{B}\big)
\]
in $\bicatX$, where 
$\iota_1$ and $\iota_2$ denote the two inclusions
$\Gamma \hookrightarrow \Gamma + \coerce{A}$ and
$\coerce{A} \hookrightarrow \Gamma + \coerce{A}$. Instantiating in the case 
where 
$(\natTrans, \natCell)$ is 
given 
by~(\ref{eq:action-as-composition}), one obtains
\[
\big( 
u_{\coerce{A}, s\sem{B}} \circ \altexp{\Yon\coerce{A}}{s\sem{-}} 
\big)(\altNat, \altCell) =
	\lambda {\left( s\sem{\altNat_{\Gamma \concat \coerce{A}}(\iota_1, 
	\iota_2)} 
	\right)} \circ \prodPres_{\Gamma, \coerce{A}}
\]
It follows that the 
value of the whole anticlockwise route is 
$\Id_{\scriptsizeexpobj{sA}{sB}} \circ 
	\lambda(s\sem{\altNat_{\Gamma + \coerce{A}}(\iota_1, \iota_2)} 
	\circ \prodPres_{\Gamma, \coerce{A}})$.
	
Next we calculate the clockwise route. For a context $\Gamma$ and 
pseudonatural transformation 
$(\altNat, \altCell)$ as above, the unlabelled equivalence returns the 1-cell
$\altNat_{\Gamma \concat \coerce{A}}(\iota_1, \iota_2)$ 
(recall~(\ref{eq:canonical-equivalence-of-exponentials}) on 
page~\pageref{eq:canonical-equivalence-of-exponentials}). This is a 
normal term of type $B$ in context 
$\Gamma \concat \coerce{A} = (\Gamma, x_{\len{\Gamma} + 1} : A)$; let 
us write
$\mathsf{j}$ for this term. The clockwise 
composite therefore returns 
\begin{align*}
s\sem{\Gamma \vdash \lam{x}{\mathsf{j}} : \exptype{A}{B}}
	&= \lambda{\left(s\sem{\Gamma, x_{\len{\Gamma} + 1} : A \vdash \mathsf{j}: 
	B} \circ 
		\seq{\pi_1 \circ \pi_1, \,\dots\, , \pi_n \circ \pi_1, \pi_2}\right)} \\
	&= \lambda{\left(s\sem{\altNat_{\Gamma + \coerce{A}}(\iota_1, \iota_2)} 
		\circ 
		\seq{\pi_1 \circ \pi_1, \,\dots\, , \pi_n \circ \pi_1, \pi_2}\right)}
\end{align*}
Since the tupling of projections on the right is exactly 
$\prodPres_{\Gamma, \coerce{A}}$ (Remark~\ref{rem:prod-pres-special-case}), the 
required 2-cell is a structural isomorphism:
\begin{align*}
\Id_{\scriptsizeexpobj{sA}{sB}} \circ 
	\lambda(s\sem{\altNat_{\Gamma \concat \coerce{A}}(\iota_1, \iota_2)} 
	\circ \prodPres_{\Gamma, \coerce{A}})
&\iso 
\lambda(s\sem{\altNat_{\Gamma \concat \coerce{A}}(\iota_1, \iota_2)} 
	\circ \prodPres_{\Gamma, \coerce{A}}) \\
&= 
\lambda\left(s\sem{\altNat_{\Gamma \concat \coerce{A}}(\iota_1, \iota_2)} 
		\circ 
		\seq{\pi_1 \circ \pi_1, \,\dots\, , \pi_n \circ \pi_1, \pi_2}\right)
\end{align*} 
\end{proof}
\end{mypropn}

\section{\texorpdfstring{$\langCartClosed$ is locally coherent}{Cartesian 
closed bicategories are locally coherent}}
\label{sec:main-result}

We are finally in a position to prove the main result. To this end, 
let $\baseTypes$ be a set of base types, $\ccBicat{\bicatX}$ be a 
cc-bicategory, 
and
$s : \allTypes\baseTypes \to \bicatX$ be the canonical extension of a set map
$\baseTypes \to \bicatX$. This extends in turn to an interpretation
$s\sem{-} : \syncloneAtClosed{\allTypes\baseTypes} \to \bicatX$. 
From this interpretation one obtains the glued objects 
of~(\ref{eq:dneut-and-dnorm-defined})~(page~\pageref{eq:dneut-and-dnorm-defined})
 and hence a set map 
$\baseTypes \to \gl{\lanext{\prodext{s}{}}}$ sending
$B \mapsto \mu_B$. This extends via the cartesian closed structure of 
$\gl{\lanext{\prodext{s}{}}}$ to an interpretation
$\overline{s}\sem{-} : \syncloneAtClosed{\allTypes\baseTypes} \to 
\gl{\lanext{\prodext{s}{}}}$.
Since the forgetful functor
$\gl{\lanext{\prodext{s}{}}} \to \bicatX$
strictly preserves the cc-bicategorical structure, we may write
$\overline{s}\sem{A} := (G_A, \gamma_B, s\sem{A})$ for every type
 $A \in \allTypes{\baseTypes}$. Moreover, for every context 
$\Gamma := (x_i : A_i)_{i=1, \,\dots\, , n}$ and term
$\Gamma \vdash t : B$ in $\langCartClosed(\allTypes\baseTypes)$,
one obtains a 1-cell 
$\overline{s}\sem{\Gamma \vdash t : B} = \prod_{i=1}^n \overline{s}\sem{A_i} 
	\to \overline{s}\sem{B}$. 
Write 
$(\overline{s}'\sem{\Gamma \vdash t : B}, 
	\overline{\sigma}\sem{\Gamma \vdash t : B}, 
	s\sem{\Gamma \vdash t : B})$ 
for this 1-cell, which is described pictorially by the following 
pseudo-commutative diagram
in $\Hom( \twoDisc{{\Con_{\allTypes{\baseTypes}}}} , \Cat)$
(note that, since $\glued{s}$ is contravariant on 
$\Con_{\allTypes{\baseTypes}}$, the composite 
$\bicatX(\glued{s}(-), X) = \bicatX(s\sem{-}, X)$ is covariant):
\begin{equation} \label{eq:interpretation-in-glued-bicat}
\begin{tikzcd}[column sep = 6em]
\prod_{i=1}^n G_{A_i} 
\arrow[phantom]{ddr}[description, xshift=2mm]
	{\twocell{\overline{\sigma}\sem{\Gamma \vdash t : B}}}
\arrow[phantom]{ddr}[swap, description, xshift=2mm, yshift=-4mm]
	{\iso}
\arrow[swap]{d}{\prod_{i=1}^n \gamma_{A_i}}
\arrow{r}{\overline{s}'\sem{\Gamma \vdash t : B}} &
G_B 
\arrow{dd}{\gamma_B} \\
\prod_{i=1}^n \bicatX\big( s\sem{-}, s\sem{A_i}\big) 
\arrow[swap]{d}{\seq{-, \,\dots\, , =}} &
\: \\
\bicatX\big( s\sem{-}, \prod_{i=1}^n s\sem{A_i})\big) 
\arrow{r}[swap]
	{s\sem{\Gamma \vdash t : B} \circ (-)} &
\bicatX\big( s\sem{-}, s\sem{B}\big)
\end{tikzcd}
\end{equation}
Finally, for every rewrite 
$\Gamma \vdash \tau : \rewrite{t}{t'} : B$ one obtains a pair of 2-cells
\begin{align*}
\overline{s}'\sem{\Gamma \vdash \tau : \rewrite{t}{t'} : B} : 
	\overline{s}'\sem{\Gamma \vdash t : B} &\To 
	\overline{s}'\sem{\Gamma \vdash t' : B} \\
s\sem{\Gamma \vdash \tau : \rewrite{t}{t'} : B} : 
	s\sem{\Gamma \vdash t : B} &\To 
	s\sem{\Gamma \vdash t' : B}
\end{align*}
which, by the cylinder condition, satisfy the diagram below. 
Since 
$\Hom(\twoDisc{{\Con_{\allTypes{\baseTypes}}}} , \Cat)$ is a 
2-category, there is no need to distinguish between bracketings. 
\begin{equation} \label{eq:glueing-bicat-cylinder-condition}
\begin{tikzcd}[column sep = 6em, row sep = 2.5em]
\gamma_B \circ \overline{s}'\sem{\Gamma \vdash t : B} 
\arrow[swap]{d}{\overline{\sigma}\sem{\Gamma \vdash t : B}}
\arrow{r}
	{
	\gamma_B \circ \overline{s}'\sem{\Gamma \vdash \tau : \rewrite{t}{t'} : B} 
	} &
\gamma_B \circ \overline{s}'\sem{\Gamma \vdash t' : B} 
\arrow{d}{\overline{\sigma}\sem{\Gamma \vdash t' : B}} \\
s\sem{\Gamma \vdash t : B} \circ \seq{-, \,\dots\, , =} 
	\circ \prod_{i=1}^n 
\gamma_{A_i} 
\arrow{r}[swap, yshift=-2mm]
	{s\sem{\Gamma \vdash \tau : \rewrite{t}{t'} : B} 
		\circ \seq{-, \,\dots\, , =} \circ \prod_{i=1}^n \gamma_{A_i}}
&
s\sem{\Gamma \vdash t' : B} \circ \seq{-, \,\dots\, , =} 
	\circ \prod_{i=1}^n \gamma_{A_i}
\end{tikzcd}
\vspace{2mm}
\end{equation}
We now use Proposition~\ref{prop:1-cells-in-glueing-bicat} to define 1-cells 
$\unquote_A : \mu_A \to \overline{s}\sem{A}$ and
$\quote_A : \overline{s}\sem{A} \to \eta_A$ by induction on types. 
%
On base types $B$, we take
\begin{align*}
\unquote_B &:= \Id_{\mu_B} : \mu_B \to \mu_B = \overline{s}\sem{B} \\
\quote_B &:= (\mathsf{inc}(-;B)^{-1}, \iso, \Id_{sB}) : \overline{s}\sem{B} \to 
\eta_B
\end{align*}
where $(\twoDisc{\mathsf{inc}}(-;B)^{-1}, \iso, \Id_{sB})$ is defined in 
Proposition~\ref{prop:1-cells-in-glueing-bicat}(\ref{c:map-in-gl-bicat-base-types}).

\newpage
On product types $\prod_n(A_1, \,\dots\, , A_n)$, the 1-cell
$\unquote_{(\prod_n \ind{A})} : 
	\mu_{(\prod_n \ind{A})} \to \prod_{i=1}^n \overline{s}\sem{A_i}$ 
is the $n$-ary tupling of the composite
\[
\mu_{(\prod_n \ind{A})} \xra{(\twoDisc{\mathsf{proj}_k}, \id, \pi_k)}
	\mu_{A_k} \xra{\unquote_{A_k}}
	\overline{s}\sem{A_k}
\]
for $k=1, \,\dots\, , n$, where the first 1-cell is defined in
Proposition~\ref{prop:1-cells-in-glueing-bicat}(\ref{c:map-in-gl-bicat-proj}). 
For $\quote_{(\prod_n \ind{A})}$, we define 
\[
\quote_{(\prod_n \ind{A})} :=
\prodop_{i=1}^n \overline{s}\sem{A_i} \xra{\prod_{i=1}^n \quote_{A_i}} 
	\prodop_{i=1}^n \eta_{A_i} 
	\xra{(\twoDisc{\mathsf{tuple}},\iso, \Id_{s\sem{\prod_n \ind{A}}} )}
	\eta_{(\prod_n \ind{A})}
\]
where the second 1-cell is defined in
Proposition~\ref{prop:1-cells-in-glueing-bicat}(\ref{c:map-in-gl-bicat-tupe}).

Finally, for exponential types we define $\unquote_{\scriptsizeexpobj{A}{B}}$ 
to be the currying of 
$\left(  \unquote_B \circ \glued{\mathsf{app}} \right) 
	\circ (\mu_{\scriptsizeexpobj{A}{B}} \times \quote_A)$, thus:
\[
\lambda{\left(  \mu_{\scriptsizeexpobj{A}{B}} \times \overline{s}\sem{A} 
		\xra{\mu_{\scriptsizescriptsizeexpobj{A}{B}} \times \quote_A} 
		(\mu_{\scriptsizeexpobj{A}{B}}) \times \eta_A 
		\xra{(\twoDisc{\mathsf{app}}(-; A, B), \id, 
				\eval_{s\sem{A}, s\sem{B}})}
		\mu_B 
		\xra{\unquote_B}
		\overline{s}\sem{B} \right)}
\]
where we use 
Proposition~\ref{prop:1-cells-in-glueing-bicat}(\ref{c:map-in-gl-bicat-app}) 
for the second arrow. For $\quote_{\scriptsizeexpobj{A}{B}}$ we define
\[
\quote_{\scriptsizeexpobj{A}{B}} := 
	(\exp{\overline{s}\sem{A}}{\overline{s}\sem{B}}) \to 
	(\exp{\nu_A}{\eta_B}) 
	\xra{(\mathsf{L}_{A,B}, 
	\iso,\Id_{\scriptsizescriptsizeexpobj{s\sem{A}}{s\sem{B}}})}
	\eta_{\scriptsizeexpobj{A}{B}}
\]
where the second arrow is defined in 
Proposition~\ref{prop:1-cells-in-glueing-bicat}(\ref{c:map-in-gl-bicat-lam}) 
and the first arrow is the currying of
$
\left( \quote_B \circ \eval_{\overline{s}\sem{A}, {\overline{s}\sem{B}}}\right)
	\circ 
\big( 
	{\left(
	(\expobj{\overline{s}\sem{A}}{\overline{s}\sem{B}}) \times 	\unquote_A
	\right)
	\circ 
	\left(
	(\expobj{\overline{s}\sem{A}}{\overline{s}\sem{B}}) 
		\times 
	\glued{\mathsf{var}}
	\right)}
\big)
$; that is, the currying of the following composite:
\begin{td}[column sep = 5em]
(\exp{\overline{s}\sem{A}}{\overline{s}\sem{B}}) \times \nu_A 
\arrow[bend left = 25]{ddrr}
\arrow[swap]{d}
	{(\scriptsizeexpobj{\overline{s}\sem{A}}{\overline{s}\sem{B}}) \times 
		\glued{\mathsf{var}}} 
\arrow[
swap,
rounded corners,
to path=
{ -- ([xshift=-1ex]\tikztostart.west)
-| ([xshift=-2cm]\tikztotarget.west)
-- (\tikztotarget.west)}, 
]{dd} &
\: &
\: \\

(\exp{\overline{s}\sem{A}}{\overline{s}\sem{B}}) \times \mu_A 
\arrow[swap]{d}[yshift=-0mm]
	{(\scriptsizeexpobj{\overline{s}\sem{A}}{\overline{s}\sem{B}}) \times 
	\unquote_A} &
\: &
\: \\

(\exp{\overline{s}\sem{A}}{\overline{s}\sem{B}}) \times 
	\overline{s}\sem{A}
\arrow[swap]{r}[yshift=0mm]
	{\glued{\eval}_{\overline{s}\sem{A}, {\overline{s}\sem{B}}}} 
\arrow[phantom]{rr}[font=\scriptsize, yshift=-1.05cm]
	{\quote_B \circ \glued{\eval}_{\overline{s}\sem{A}, {\overline{s}\sem{B}}}}
\arrow[
swap,
rounded corners,
to path=
{ -- ([yshift=-1ex]\tikztostart.south)
|- ([yshift=-.5cm]\tikztotarget.south)
-- (\tikztotarget.south)}, 
]{rr} &
\overline{s}\sem{B}
\arrow[swap]{r}{\quote_B} &
\eta_B 
\end{td}
The morphism 
$\glued{\mathsf{var}} 
	:= 
	\big(\twoDisc{\mathsf{var}(-; A_i)}, \iso, \Id_{s\sem{A_i}}\big)$ 
is defined in 
Proposition~\ref{prop:1-cells-in-glueing-bicat}(\ref{c:map-in-gl-bicat-var}). 
Let us denote
$\unquote_B := (\widehat{u}_B, \overline{u}_B, u_B)$ and
$\quote_B := (\widehat{q}_B, \overline{q}_B, q_B)$, so that
$\pi_{\mathrm{dom}}(\unquote_B) = u_B$ and $\pi_{\mathrm{dom}}(\quote_B) = q_B$.

\newpage
\begin{mylemma}
For every type $B \in \allTypes\baseTypes$, there exist natural 
isomorphisms
$\pi_{\mathrm{dom}}(\unquote_B) \iso \Id_{s\sem{B}}$ and
$\pi_{\mathrm{dom}}(\quote_B) \iso \Id_{s\sem{B}}$.
\begin{proof}
We proceed inductively. On base types the claim holds trivially. For 
product 
types, we observe that, where 
$A_1, \,\dots\, , A_n \in \allTypes\baseTypes \:\: (n \in \Nat)$:
\begin{align*}
\pi_{\mathrm{dom}}(\unquote_{(\prod_n \ind{A})}) &= 
	\seq{ u_{A_1} \circ \pi_1, \,\dots\, , u_{A_n} \circ \pi_n } \\
	&\iso \left(\prodop_{i=1}^n u_{A_i}\right) \circ \seq{\pi_1, \,\dots\, , 
	\pi_n} 
	\\
	&\overset{\text{IH}}{\iso} \left(\prodop_{i=1}^n \Id_{A_i} \right) \circ  
		\seq{\pi_1, \,\dots\, , \pi_n} \\
	&\iso \Id_{s\sem{\prod_n \ind{A}}} \\[6pt]
\pi_{\mathrm{dom}}(\quote_{(\prod_n \ind{A})}) &= 
	\Id_{s\sem{\prod_n \ind{A}}} \circ \prodop_{i=1}^n q_{A_i} \\
	&\iso \prodop_{i=1}^n q_{A_i} \\
	&\overset{\text{IH}}{\iso} 
	\prodop_{i=1}^n \Id_{s\sem{A_i}} \\
	&\iso \Id_{s\sem{\prod_n \ind{A}}}
\end{align*}
Finally, for exponentials, one sees that
\begin{align*}
\pi_{\mathrm{dom}}(\unquote_{\scriptsizeexpobj{A}{B}}) &= 
	\lambda{\left( \left(u_B \circ \eval_{s\sem{A}, s\sem{B}}\right) \circ 
		(\Id_{s\sem{\scriptsizeexpobj{A}{B}}} \times q_A)  \right)} \\
	&\overset{\text{IH}}{\iso} 
		\lambda{\left( 
			\left(\Id_{s\sem{B}} \circ \eval_{s\sem{A}, s\sem{B}}\right) 
		\circ 
				(\Id_{s\sem{\scriptsizeexpobj{A}{B}}} \times \Id_{s\sem{A}})  
				\right)} \\
	&\iso 	\lambda{\left(\eval_{s\sem{A}, s\sem{B}} \circ 
					(\Id_{s\sem{\scriptsizeexpobj{A}{B}}} \times 
					\Id_{s\sem{A}})  
					\right)} \\
	&\overset{\etaExp{}}{\iso} \Id_{s\sem{\scriptsizeexpobj{A}{B}}} \\[6pt]
\pi_{\mathrm{dom}}(\quote_{\scriptsizeexpobj{A}{B}})  
	&\iso \lambda{\left(
				(q_B 
				\circ \eval_{{s}\sem{A}, {\sem{B}}})
				\circ \big(
					(\Id_{s\sem{\scriptsizeexpobj{A}{B}}} \times u_A) 
				\circ (\Id_{s\sem{\scriptsizeexpobj{A}{B}}} \times 	
						\Id_{s\sem{A}})\big)\right) } \\
	&\overset{\text{IH}}{\iso} 
		\lambda{\left(
					(\Id_{s\sem{B}} 
					\circ \eval_{{s}\sem{A}, {\sem{B}}})
					\circ \big(
					(\Id_{s\sem{\scriptsizeexpobj{A}{B}}} \times u_A) 
				\circ (\Id_{s\sem{\scriptsizeexpobj{A}{B}}} \times 	
						\Id_{s\sem{A}})\big)\right) } 		\\
	&\iso 
		\lambda{\left(
					\left(\Id_{s\sem{B}} 
					\circ \eval_{{s}\sem{A}, {\sem{B}}}\right)
					\circ \left(
						(\Id_{s\sem{\scriptsizeexpobj{A}{B}}} \times 
											\Id_{s\sem{A}})\right)\right) }
					\\
	&\iso \lambda{\left(
		\eval_{{s}\sem{A}, {\sem{B}}} \circ 
		\left(\Id_{s\sem{\scriptsizeexpobj{A}{B}}} \times 
						\Id_{s\sem{A}}\right)\right)} \\
	&\overset{\etaExp{}}{\iso} \Id_{s\sem{\scriptsizeexpobj{A}{B}}}
\end{align*}
In each case the isomorphisms are composites of structural isomorphisms or 
canonical isomorphisms for the cartesian closed structure, hence natural.
\end{proof}
\end{mylemma}

The definitions of $\unquote$ and $\quote$, together with the preceding lemma 
and the 2-cells 
$\psi_X^{s\sem{-}}$, 
give rise to diagrams of the 
following form for every type $B \in \allTypes\baseTypes$:
\begin{center}
\begin{tikzcd}[column sep = 4em, row sep = 3em]
\twoDisc{\dneutTerms(-;B)}
\arrow[phantom]{dr}[description]{\twocell{\overline{u}_B}}
\arrow[swap]{d}{s\sem{-}} 
\arrow{r}{\widehat{u}_B} &
G_B 
\arrow{d}{\gamma_B} \\
\bicatX\big( s\sem{-}, s\sem{B}\big)
\arrow[bend right = 11, phantom]{r}[description]{\iso} 
\arrow[bend right]{r}[swap]{\Id_{\bicatX( s\sem{-}, 
s\sem{B}) }}
\arrow{r}
	{\bicatX( s\sem{-}, u_B)} &
\bicatX\big( s\sem{-}, s\sem{B}\big)
\end{tikzcd}
\qquad
\begin{tikzcd}[column sep = 4em, row sep = 3em]
G_B
\arrow[swap]{d}{\gamma_B}
\arrow[phantom]{dr}[description]{\twocell{\overline{q}_B}}
\arrow{r}{\widehat{q}_B} &
\twoDisc{\dnormTerms(-;B)} 
\arrow{d}{s\sem{-}}  \\
\bicatX\big( s\sem{-}, s\sem{B}) 
\arrow{r}
	{\bicatX( s\sem{-}, q_B)} 
\arrow[bend right = 12, phantom]{r}[description]{\iso} 
\arrow[bend right]{r}[swap]{\Id_{\bicatX( s\sem{-}, 
s\sem{B}) }} &
\bicatX\big( s\sem{-}, s\sem{B})
\end{tikzcd}
\end{center}
Thus, for any sequence of types 
$A_1, \,\dots\, , A_n \in \allTypes\baseTypes \:\: (n \in \Nat)$, one obtains 
a diagram of shape
\begin{td}[column sep = 8em, row sep = 3em]
\prod_{i=1}^n \twoDisc{\dneutTerms(-;A_i)}
\arrow[phantom]{dr}[description, yshift=2mm]{\twocell{\iso}}
\arrow[swap]{d}{\prod_{i=1}^n s\sem{-}} 
\arrow{r}{\prod_{i=1}^n \widehat{u}_{A_i}} &
\prod_{i=1}^n G_{A_i} 
\arrow{d}{\prod_{i=1}^n \gamma_{A_i} } \\
\prod_{i=1}^n \bicatX\big( s\sem{-}, 
s\sem{A_i}\big) 
\arrow{r}
	{\prod_{i=1}^n \bicatX( s\sem{-}, u_{A_i})}
\arrow[bend right = 10, phantom]{r}[description]{\iso} 
\arrow[bend right=20]{r}[swap]
	{\Id_{\prod_i \bicatX( s\sem{-}, s\sem{A_i})}} &
\prod_{i=1}^n \bicatX\big( s\sem{-}, s\sem{A_i} \big)
\end{td}
by composing with the $\fuse$ 2-cells. Pasting these diagrams together 
with~(\ref{eq:interpretation-in-glued-bicat}), one obtains the following
diagram in $\Hom( \twoDisc{{\Con_{\allTypes{\baseTypes}}}} , \Cat)$ for 
every rewrite 
$(\Gamma \vdash \tau: \rewrite{t}{t'} : B)$ 
in $\langCartClosed(\allTypes\baseTypes)$.  We write
$\overline{s}'\sem{\tau}$ for 
$\overline{s}'\sem{\Gamma \vdash \tau : \rewrite{t}{t'} : B}$ and
$s\sem{\tau}$ for
$s\sem{\Gamma \vdash \tau : \rewrite{t}{t'} : B}$. Since there are no 
constants in $\langCartClosed(\allTypes\baseTypes)$, these rewrites are 
necessarily invertible.
\begin{equation} \label{eq:main-diagram}
\makebox[\textwidth]{
\begin{tikzcd}[column sep = 4.5em, row sep = 5em, ampersand replacement=\&]
\prod_{i=1}^n \twoDisc{\dneutTerms(-;A_i)}
\arrow[phantom]{dr}[description, yshift=2mm, yshift=1mm]{\twocell{\iso}}
\arrow[swap]{d}{\prod_{i=1}^n s\sem{-}} 
\arrow{r}{\prod_{i=1}^n \widehat{u}_{A_i}} \&
\prod_{i=1}^n G_{A_i} 
\arrow[phantom]{ddr}[description, xshift=3mm]
	{\twocell{\overline{\sigma}\sem{\Gamma \vdash t : B}}}
\arrow[phantom]{ddr}[description, xshift=3mm, yshift=-4mm]
	{\iso}
\arrow[swap]{d}{\prod_{i=1}^n \gamma_{A_i}}
\arrow[bend left=15]{r}{\overline{s}'\sem{\Gamma \vdash t' : B}}
\arrow[bend right=15, swap]{r}{\overline{s}'\sem{\Gamma \vdash t : B}}
\arrow[phantom]{r}[description, font=\scriptsize]
	{\overset{\overline{s}'\sem{\tau}}{\Uparrow\iso}} \&
G_B 
\arrow{dd}{\gamma_B}
\arrow{r}{\widehat{q}_B}
\arrow[phantom]{ddr}[description]{\twocell{\overline{q}_B}}
\arrow[phantom]{ddr}[description, yshift=-4mm]{\iso} \&
\twoDisc{\dnormTerms(-;B)} 
\arrow{dd}{s\sem{-}} \\
\prod_{i=1}^n \bicatX( s\sem{-}, s\sem{A_i}) 
\arrow{r}[yshift=1.5mm]
	{\prod_{i=1}^n \bicatX( s\sem{-}, u_{A_i})}
\arrow[bend right = 10, phantom]{r}[description]{\iso} 
\arrow[bend right=20]{r}[swap]
	{\Id_{\prod_i \bicatX( s\sem{-}, s\sem{A_i})}} \&
\prod_{i=1}^n \bicatX\big( s\sem{-}, s\sem{A_i}\big) 
\arrow[swap]{d}[description]{\seq{-, \,\dots\, , =}} \&
\: \\
\: \&
\bicatX{\left( s\sem{-}, \prodop_{i=1}^n s\sem{A_i}\right)} 
\arrow[bend left=15]{r}
	{s\sem{\Gamma \vdash t : B} \circ (-)}
\arrow[bend right=15]{r}[swap]
	{s\sem{\Gamma \vdash t' : B} \circ (-)}
\arrow[phantom]{r}[description, font=\scriptsize]
	{\overset{s\sem{\tau} \circ (-)}{\Downarrow\iso}} \&
\bicatX( s\sem{-}, s\sem{B})
\arrow{r}
	{\bicatX( s\sem{-}, q_B)} 
\arrow[bend right = 12, phantom]{r}[description]{\iso} 
\arrow[bend right]{r}[swap]{\Id_{\bicatX( s\sem{-}, s\sem{B}) }} \&
\bicatX( s\sem{-}, s\sem{B})
\end{tikzcd}
}
\vspace{2mm}
\end{equation}

The proof now hinges on two facts. Firstly, since 
$\dnormTerms(-;B)$ is a set, the composite 2-cell obtained by whiskering 
across the top row of the diagram above must be the identity. 

Secondly, the middle part of the diagram satisfies the cylinder condition. 
Precisely, writing $\mathrm{tup}$ for $\seq{-, \,\dots\, , =}$, let
$\kappa_t$ be the invertible 2-cell obtained from the front face:
\begin{equation} \label{eq:def-of-kappa-t}
\makebox[\textwidth]{
\begin{tikzcd}[ampersand replacement=\&]
s\sem{-} \circ \widehat{q}_B \circ  
	\overline{s}'\sem{\Gamma \vdash t : B} \circ 
	\prodop_{i=1}^n \widehat{u}_{A_i} 
\arrow{r}{\kappa_t}
\arrow[swap]{d}
	{\overline{q}_B \circ \overline{s}'\sem{\Gamma \vdash t : B} \circ 
		\prod_{i=1}^n \widehat{u}_{A_i}}
\arrow{d}{\iso} \&
s\sem{\Gamma \vdash t : B} 
		\circ \mathrm{tup} 
		\circ \prod_{i=1}^n s\sem{-} \\
\bicatX(s\sem{-}, q_B) \circ \gamma_B \circ 
	\overline{s}'\sem{\Gamma \vdash t : B} \circ 
	\prodop_{i=1}^n \widehat{u}_{A_i} 
\arrow[swap]{d}{\iso} \&
s\sem{\Gamma \vdash t : B} 
		\circ \mathrm{tup} 
		\circ \Id_{\bicatX(s\sem{-}, u_{A_i})}
		\circ \prod_{i=1}^n s\sem{-}
\arrow[swap]{u}{\iso} \\
\Id_{\bicatX( s\sem{-}, s\sem{B}) }  \circ \gamma_B \circ 
	\overline{s}'\sem{\Gamma \vdash t : B} \circ 
	\prodop_{i=1}^n \widehat{u}_{A_i} 
\arrow[swap]{d}{\iso} \&
s\sem{\Gamma \vdash t : B} 
		\circ \mathrm{tup} 
		\circ \prod_{i=1}^n \bicatX(s\sem{-}, u_{A_i}) 
		\circ \prod_{i=1}^n s\sem{-} 
\arrow[swap]{u}{\iso}  \\
\gamma_B \circ 
	\overline{s}'\sem{\Gamma \vdash t : B} \circ 
	\prodop_{i=1}^n \widehat{u}_{A_i} 
\arrow[swap]{d}
	{\overline{\sigma}\sem{\Gamma \vdash t : B} \circ 
		\prod_{i=1}^n \widehat{u}_{A_i}}
\arrow{d}{\iso} \&
s\sem{\Gamma \vdash t : B} 
		\circ \mathrm{tup} 
		\circ \prod_{i=1}^n \left(\bicatX(s\sem{-}, u_{A_i}) 
			\circ s\sem{-} \right) 
\arrow[swap]{u}
	{s\sem{\Gamma \vdash t : B} 
			\circ \mathrm{tup} 
			\circ \fuse^{-1}}
\arrow{u}{\iso} \\
s\sem{\Gamma \vdash t : B} 
		\circ \mathrm{tup} \circ \prod_{i=1}^n \gamma_{A_i} 
		\circ \prodop_{i=1}^n \widehat{u}_{A_i} 
\arrow[swap]{r}[yshift=-2mm]
	{s\sem{\Gamma \vdash t : B}  \circ \mathrm{tup} \circ 
		\fuse} \&
s\sem{\Gamma \vdash t : B} 
		\circ \mathrm{tup} 
		\circ \prod_{i=1}^n (\gamma_{A_i} \circ \widehat{u}_{A_i}) 
\arrow[swap]{u}
	{s\sem{\Gamma \vdash t : B}  \circ \mathrm{tup} \circ
		\prod_{i=1}^n \overline{u}_{A_i}}
\arrow{u}{\iso}
\end{tikzcd}
}
\vspace{2mm}
\end{equation}
The cylinder 
condition~(\ref{eq:glueing-bicat-cylinder-condition}) and the functorality of 
horizontal composition imply that $\kappa_t$ satisfies the following 
property in $\Hom( \twoDisc{{\Con_{\allTypes{\baseTypes}}}} , \Cat)$: 
\begin{td}[column sep = 5em]
s\sem{-} \circ \widehat{q}_B \circ  
	\overline{s}'\sem{\Gamma \vdash t : B} \circ 
	\prodop_{i=1}^n \widehat{u}_{A_i} 
\arrow{r}[yshift=2mm]
	{s\sem{-} \circ \widehat{q}_B \circ  
		\overline{s}'\sem{\Gamma \vdash \tau : \rewrite{t}{t'} : B} \circ 
		\prodop_{i=1}^n \widehat{u}_{A_i}}
\arrow[swap]{d}{\kappa_{t}}
\arrow{d}{\iso} &
s\sem{-} \circ \widehat{q}_B \circ  
	\overline{s}'\sem{\Gamma \vdash t' : B} \circ 
	\prodop_{i=1}^n \widehat{u}_{A_i} 
\arrow{d}{\kappa_{t'}} \\

s\sem{\Gamma \vdash t : B} 
		\circ \mathrm{tup} 
		\circ \prod_{i=1}^n s\sem{-}
\arrow[swap]{r}[yshift=-2mm]{s\sem{\Gamma \vdash \tau : \rewrite{t}{t'} : B} 
		\circ \mathrm{tup} 
		\circ \prod_{i=1}^n s\sem{-}}  &
s\sem{\Gamma \vdash t' : B} 
		\circ \mathrm{tup} 
		\circ \prod_{i=1}^n s\sem{-}
\arrow{u}{\iso}
\end{td}
Applying the first fact, this diagram degenerates to the following:
\begin{equation} \label{eq:degenerate-diagram}
\begin{tikzcd}
s\sem{-} \circ \widehat{q}_B \circ  
	\overline{s}'\sem{\Gamma \vdash t : B} \circ 
	\prodop_{i=1}^n \widehat{u}_{A_i} 
\arrow[equals]{r}
\arrow{d}{\iso}
\arrow[swap]{d}{\kappa_{t}} &
s\sem{-} \circ \widehat{q}_B \circ  
	\overline{s}'\sem{\Gamma \vdash t' : B} \circ 
	\prodop_{i=1}^n \widehat{u}_{A_i} 
\arrow{d}{\kappa_{t'}} \\
s\sem{\Gamma \vdash t : B} 
		\circ \mathrm{tup} 
		\circ \prod_{i=1}^n s\sem{-}
\arrow[swap]{r}[yshift=-2mm]{s\sem{\Gamma \vdash \tau : \rewrite{t}{t'} : B} 
		\circ \mathrm{tup} 
		\circ \prod_{i=1}^n s\sem{-}}  &
s\sem{\Gamma \vdash t' : B} 
		\circ \mathrm{tup} 
		\circ \prod_{i=1}^n s\sem{-}
\arrow{u}{\iso}
\end{tikzcd}
\end{equation}
Instantiating the bottom row of this diagram at the context 
$\Gamma := (x_i : A_i)_{i=1, \,\dots\, , n}$
and the $n$-tuple of terms
$(\Gamma \vdash x_i : A_i)_{i=1, \,\dots\, , n}$, one sees that
\begin{align*}
\left(s\sem{\Gamma \vdash t : B} 
		\circ \mathrm{tup} 
		\circ \prodop_{i=1}^n s\sem{-}\right)
		(\Gamma \vdash x_i : A_i)_{i=1, \,\dots\, ,n} 
&= s\sem{\Gamma \vdash t : B} \circ 
	\seq{s\sem{\Gamma \vdash {x}_i : {A}_i}}_i \\
&= s\sem{\Gamma \vdash t : B} \circ 
	\seq{\pi_1, \,\dots\, , \pi_n}
\end{align*}
We may now extend~(\ref{eq:degenerate-diagram}) downwards. Writing 
$T_t := 
	s\sem{-} \circ \widehat{q}_B \circ  
		\overline{s}'\sem{\Gamma \vdash t : B} \circ 
		\prodop_{i=1}^n \widehat{u}_{A_i}$ 
and instantiating at
$(\Gamma \vdash x_i : A_i)_{i=1, \,\dots\, , n}$, one obtains the following 
diagram. 
\begin{equation} \label{eq:reduction-to-interpretation}
\begin{tikzcd}[column sep = 4em, row sep = 3em]
T_t(\Gamma \vdash x_i : A_i)_{i=1, \,\dots\, , n} 
\arrow[equals]{r}
\arrow[swap]{d}{\kappa_{t}}
\arrow{d}{\iso} &
T_{t'}(\Gamma \vdash x_i : A_i)_{i=1, \,\dots\, , n}  
\arrow{d}{\kappa_{t'}}
\arrow[swap]{d}{\iso} \\
s\sem{\Gamma \vdash t : B} \circ 
	\seq{\pi_1, \,\dots\, , \pi_n}
\arrow[swap]{r}[yshift=-2mm]{
	s\sem{\Gamma \vdash \tau : \rewrite{t}{t'} : B} \circ 
			\seq{\pi_1, \,\dots\, , \pi_n}} 
\arrow[swap]{d}
	{\widehat{\etaTimes{}}_{\Id_{s\sem{\Gamma}}}^{-1}}
\arrow{d}{\iso}  &
s\sem{\Gamma \vdash t' : B} \circ 
	\seq{\pi_1, \,\dots\, , \pi_n} 
\arrow{d}
	{\widehat{\etaTimes{}}_{\Id_{s\sem{\Gamma}}}^{-1}}
\arrow[swap]{d}{\iso} \\
s\sem{\Gamma \vdash t : B} \circ 
	\Id_{s\sem{\Gamma}}
\arrow[swap]{r}[yshift=-2mm]{s\sem{\Gamma \vdash \tau : \rewrite{t}{t'} : B} 
		\circ \Id_{s\sem{\Gamma}}} 
\arrow[swap]{d}{\iso} &
s\sem{\Gamma \vdash t' : B} \circ 
	\Id_{s\sem{\Gamma}} 
\arrow{d}{\iso} \\
s\sem{\Gamma \vdash t : B}
\arrow[swap]{r}{s\sem{\Gamma \vdash \tau : \rewrite{t}{t'} : B}} &
s\sem{\Gamma \vdash t' : B}
\end{tikzcd}
\end{equation}
The bottom two squares commute by naturality. Hence, since each component is
invertible, it must be the case that 
$s\sem{\Gamma \vdash \tau : \rewrite{t}{t'} : B}$ 
is equal to the clockwise composite around this diagram. 
We record this result as the 
following proposition.

\begin{prooflesspropn} \label{prop:semantics-must-be-equal}
For any set of base types $\baseTypes$,
cc-bicategory $\ccBicat\bicatX$ and interpretation
$s : \baseTypes \to \bicatX$, the induced interpretation 
$s\sem{\Gamma \vdash \tau : \rewrite{t}{t'} : B}$
of any rewrite
$(\Gamma \vdash \tau : \rewrite{t}{t'} : B)$ 
in $\bicatX$
is equal to the 2-cell obtained by 
composing clockwise around~(\ref{eq:reduction-to-interpretation}). 
Moreover, 
this 2-cell depends only on the context $\Gamma$, the type $B$, and the 
terms 
$t$ and
$t'$. 
\end{prooflesspropn}

Hence, any
pair of parallel rewrites
$(\Gamma \vdash \tau : \rewrite{t}{t'} : B)$ and
$(\Gamma \vdash \tau' : \rewrite{t}{t'} : B)$
must be interpreted by the same 2-cell, namely the 2-cell obtained by 
composing clockwise around~(\ref{eq:reduction-to-interpretation}).

\begin{prooflessthm}
For any parallel pair of rewrites 
$\Gamma \vdash \tau : \rewrite{t}{t'} : B$ and
$\Gamma \vdash \tau' : \rewrite{t}{t'} : B$ in 
$\langCartClosed(\allTypes\baseTypes)$, the interpretations
$s\sem{\Gamma \vdash \tau : \rewrite{t}{t'} : B}$ and
$s\sem{\Gamma \vdash \tau' : \rewrite{t}{t'} : B}$ are equal.
\end{prooflessthm}

We wish to instantiate this theorem in the syntactic bicategory to see that
any parallel pair of rewrites must be equal in the equational theory of 
$\langCartClosed$. However, the cc-pseudofunctor 
$\inc\sem{-} : \syncloneAtClosed{\allTypes\baseTypes} \to \syncloneAtClosed{\allTypes\baseTypes}$ 
extending the inclusion
$\inc : \baseTypes \hookrightarrow \syncloneAtClosed{\allTypes\baseTypes}$ 
is not the identity: the definition for lambda abstractions requires an extra 
equivalence. Nonetheless, one can leverage the universal property to show that
$\inc\sem{-}$ is equivalent to the 
identity~(\cf~Corollary~\ref{cor:uniqueness-up-to-equiv-of-ext-pseudofun}). 

\newpage
\begin{mylemma} \label{lem:inc-sem-equivalent-to-id}
For any set of base types $\baseTypes$, the cc-pseudofunctor 
$\inc\sem{-} : \syncloneAtClosed{\allTypes\baseTypes} \to 
\syncloneAtClosed{\allTypes\baseTypes}$ extending the inclusion 
$\inc : \allTypes\baseTypes \hookrightarrow \syncloneAtClosed{\allTypes\baseTypes}$
is equivalent to the identity. Hence,
$\inc\sem{-}$ is a biequivalence.
\begin{proof}
By Proposition~\ref{prop:ccc-syntactic-models-biequivalent}, the canonical 
cc-pseudofunctor
$\usem{(-)}{\inc} : \freeCartClosedBicat{\allTypes\baseTypes} \to \syncloneAtClosed{\allTypes\baseTypes}$ 
(defined in Lemma~\ref{lem:free-cc-bicat-on-sig}) is part of a 
biequivalence; write 
$V_\inc$ for its 
pseudo-inverse. Moreover, considering the diagram
\begin{td}
\syncloneAtClosed{\allTypes\baseTypes}
\arrow{r}{\inc\sem{-}} &
\syncloneAtClosed{\allTypes\baseTypes} \\

\freeCartClosedBicat{\allTypes\baseTypes} 
\arrow[swap]{ur}{\usem{(-)}{\inc}}
\arrow{u}{\usem{(-)}{\inc}} &
\:
\end{td}
and applying 
Lemma~\ref{lem:syntactic-model-uniqueness-up-to-equivalence}, one sees that 
there exists an 
equivalence 
$\inc\sem{-} \circ \usem{(-)}{\inc} \simeq \usem{(-)}{\inc}$. One therefore 
obtains 
a chain of equivalences
\begin{align*}
\id_{\syncloneAtClosed{\allTypes\baseTypes}} &\simeq \usem{(-)}{\inc} \circ V_\inc \\ 
		&\simeq (\inc\sem{-} \circ \usem{(-)}{\inc}) \circ V_\inc \\
		&\simeq \inc\sem{-} \circ \id_{\syncloneAtClosed{\allTypes\baseTypes}} \\
		&\simeq \inc\sem{-}
\end{align*}
as required.
\end{proof}
\end{mylemma}

We can finally prove our theorem.

\begin{mythm} \label{thm:main-theorem}
For any set of base types $\baseTypes$ and
any rewrites
$(\Gamma \vdash \tau : \rewrite{t}{t'} : B)$ and
$(\Gamma \vdash \tau' : \rewrite{t}{t'} : B)$ 
in $\langCartClosed(\allTypes\baseTypes)$, the judgement
$(\Gamma \vdash \tau \equiv \tau' : \rewrite{t}{t'} : B)$
is derivable in $\langCartClosed(\allTypes\baseTypes)$.
Hence, $\langCartClosed(\allTypes\baseTypes)$ is locally coherent.
\begin{proof}
Consider the interpretation in the syntactic model
$\inc\sem{-} : \syncloneAtClosed{\allTypes\baseTypes} \to 
	\syncloneAtClosed{\allTypes\baseTypes}$
extending the inclusion of base types. Instantiating 
Proposition~\ref{prop:semantics-must-be-equal}, one sees that
$\inc\sem{\Gamma \vdash \tau : \rewrite{t}{t'} : B}
=
\inc\sem{\Gamma \vdash \tau' : \rewrite{t}{t'} : B}$
for every parallel pair of rewrites $\tau$ and $\tau'$. But 
biequivalences are locally fully faithful, so by the preceding lemma 
$\inc\sem{\Gamma \vdash \tau : \rewrite{t}{t'} : B}
	=
\inc\sem{\Gamma \vdash \tau' : \rewrite{t}{t'} : B}$ holds
if and only if 
$\tau$ and $\tau'$ are equal 2-cells in $\syncloneAtClosed{\allTypes\baseTypes}$; that 
is, 
$(\Gamma \vdash \tau \equiv \tau' : \rewrite{t}{t'} : B)$.
\end{proof}
\end{mythm}

\newpage
\begin{mythm} \label{thm:main-result-on-free-bicat}
Let $\baseTypes$ be any set and 
$\tau, \sigma : t \To t'$ be a parallel pair of 2-cells in the free 
cc-bicategory on $\baseTypes$. Then $\tau \equiv \sigma$. 
\begin{proof}
By Proposition~\ref{propn:first-biequivalence}, the syntactic bicategory 
$\syncloneAtClosed{\allTypes\baseTypes}$ is biequivalent to $\freeCartClosedBicat{\allTypes\baseTypes}$, 
the free 
cc-bicategory on 
$\baseTypes$. By the preceding theorem, the images of the 2-cells $\tau$ and 
$\sigma$
in $\syncloneAtClosed{\allTypes\baseTypes}$ must be equal. 
Since biequivalences are locally fully faithful, it follows that 
$\tau \equiv \sigma$.
\end{proof}
\end{mythm}

We can express this informally as follows. For any 
cc-bicategory 
$\ccBicat{\baseCat}$ and pair of parallel 2-cells 
$\sigma, \tau : f \To g$ in $\baseCat$, if $\sigma$ and $\tau$ are constructed 
from the cartesian closed structure using solely structural isomorphisms and 
the operations of vertical composition and horizontal composition, then 
$\sigma = \tau$. 
As a slogan: \emph{all pasting diagrams in the free cc-bicategory commute}.

\subsection{Evaluating the proof}
\label{sec:proof-failure}

It is worth examining where the proof of 
Theorem~\ref{thm:main-theorem} would fail if $\langCartClosed$ were not locally 
coherent. Our reasoning here is only 
informal, but it should provide a measure of confidence 
that the many pages of 
proof do not contain a fatal error, as well as throwing light on what makes the 
argument work. 

%

The normalisation-by-evaluation proof hinges 
crucially on two 
facts: (1) that any interpretation of $\langCartClosed$ induces an 
interpretation in the glueing bicategory, and (2) that the canonical 
interpretation of $\langCartClosed$ in
the syntactic model is biequivalent to the identity. The first fact 
entails that, whenever $\tau$ and $\sigma$ are parallel rewrites of type 
$\rewrite{t}{t'}$, 
their interpretations $s\sem{\tau}$ and $s\sem{\sigma}$ must coincide in 
every model. Then, writing $\mathsf{J}$ for the inverse to
$((\inc\sem{-})_{\Gamma, A})_{t,t'} : 
	\syncloneAtClosed{\allTypes\baseTypes}(\Gamma; A)(t, t')
		\to
	\syncloneAtClosed{\allTypes\baseTypes}(\Gamma; A)(\inc\sem{t}, \inc\sem{t'})$,
the second fact 
allows one to construct the chain of equalities
\[
\sigma \equiv 
\mathsf{J}(\inc\sem{\sigma}) \equiv 
\mathsf{J}(\inc\sem{\tau}) \equiv 
\tau
\]
witnessing local coherence. We give a small example showing how (1) fails if 
one adds extra structure that is not locally coherent.

Consider the $\langCartClosed$-signature $\sig$ consisting of a 
set of base types and a 
single constant rewrite $x : B \vdash \constrewr : \rewrite{x}{x} : B$ at 
a base type $B$. Since we add no extra equations, $\langCartClosed(\sig)$ is 
clearly not locally coherent. Now let $\ccBicat{\bicatX}$ be any 
cc-bicategory and $s : \baseTypes \to \bicatX$ an interpretation of base 
types.
Since variables are normal terms, the interpretation of our additional 
rewrite in the glueing bicategory as 
in~(\ref{eq:interpretation-in-glued-bicat}) on 
page~\pageref{eq:interpretation-in-glued-bicat} yields the diagram below,
for which we use the fact that the interpretation of the judgement
$(x : B \vdash x : B)$ is the identity: 
\begin{equation*} 
\begin{tikzcd}[column sep = 6em, row sep = 6em]
\twoDisc{\dneutTerms(-; B)}
\arrow[phantom]{dr}[description]
	{\twocell{\iso}}
\arrow[swap]{d}{s\sem{-}}
\arrow[bend left = 11]{r}{\id_{\twoDisc{\dneutTerms(-; B)}}}
\arrow[phantom]{r}[description, font=\scriptsize]
	{\twocellDown{\overline{s}'\sem{\constrewr}}}
\arrow[bend right = 11, swap]{r}{\id_{\twoDisc{\dneutTerms(-; B)}}}  &
\twoDisc{\dneutTerms(-; B)}
\arrow{d}{s\sem{-}} \\
\bicatX\big( s\sem{-}, s\sem{B})\big) 
\arrow[phantom]{r}[description, font = \scriptsize]
	{\twocellDown
		{s\sem{\constrewr} \circ (-)}}
\arrow[bend left = 11]{r}
	{s\sem{x : B \vdash x : B} \circ (-)}
\arrow[bend right = 11]{r}[swap]
	{s\sem{x : B \vdash x : B} \circ (-)} &
\bicatX\big( s\sem{-}, s\sem{B}\big)
\end{tikzcd}
\end{equation*}
Since $\twoDisc{\dneutTerms(-; B)}$ is locally discrete, the 
2-cell $\overline{s}'\sem{x: B \vdash \constrewr : \rewrite{x}{x} : B}$
can only be the identity.  Now consider a context $\Gamma$ and evaluate at a
neutral term $\into{t} \in \neutTerms(\Gamma; B)$. The isomorphism 
filling the central shape is the structural isomorphism
$s\sem{\Gamma \vdash t : B} \overset{\l_{s\sem{t}}}{\iso} 
	\Id_{s\sem{B}} \circ s\sem{\Gamma \vdash t : B}$,
so the cylinder condition requires that
\[
s\sem{x : B \vdash \constrewr : \rewrite{x}{x} : B} 
	= \l_{s\sem{t}} \vert \id_{\id_{\twoDisc{\dneutTerms(-; B)}}} \vert 
	\l_{s\sem{t}}^{-1}
	= \id_{s\sem{x}}
	= s\sem{x : B \vdash \id_x : \rewrite{x}{x} : B}
\]
Now, following the argument employed to prove Theorem~\ref{thm:main-theorem}, 
one sees that 
this equation is satisfied for the interpretation extending
$\inc : \baseTypes \hookrightarrow \syncloneAtClosed{\allTypes\baseTypes}$
if and only if the judgement
$(x : B \vdash \constrewr \equiv \id_x : \rewrite{x}{x} : B)$
is derivable. Since we assumed this not to be the case, 
\emph{the cylinder condition cannot hold}. Thus, 
the constant rewrite $\constrewr$ may not be soundly interpreted in every 
glueing bicategory $\gl{\lanext{s}}$, so one cannot rerun the 
normalisation-by-evaluation proof. 

\section{Another Yoneda-style proof of coherence}
\label{sec:yoneda-coherence}

Proposition~\ref{prop:yoneda-coherence-for-cc-bicats} proved a form of 
coherence for cc-bicategories. It turns out that this can be extended to an 
alternative proof of the main result just presented. The strategy is similar to 
that presented in Section~\ref{sec:main-result}, but 
only relies on the 
universal property of the free cc-bicategory $\freeCartClosedBicat{\allTypes\baseTypes}$ (defined 
in Construction~\ref{constr:free-cc-bicat}). 
Nonetheless, the development 
highlights the core of the normalisation-by-evaluation 
argument as just described.

Fix a set of base types $\baseTypes$ 
and an interpretation $h: \baseTypes \to \bicatX$ in 
a \mbox{cc-bicategory} $\ccBicat{\bicatX}$. This extends 
to an 
interpretation $\allTypes\baseTypes \to \bicatX$ we also denote by $h$. Now let 
$\ccBicat{\altCat}$ be a \mbox{2-category} with strict products and exponentials and
$(F, \prodPres, \expPres) : 
	\ccBicat{\bicatX} \to \ccBicat{\altCat}$
be any cc-pseudofunctor. Writing $F_0$ for the underlying set map
$ob(\bicatX) \to ob(\altCat)$, one obtains
an interpretation $F_0 \circ h : \baseTypes \to \altCat$. One thereby obtains a 
\emph{weak} interpretation in 
$\bicatX$ and a \emph{strict} interpretation 
in $\altCat$. The situation is described by the following commutative diagram:
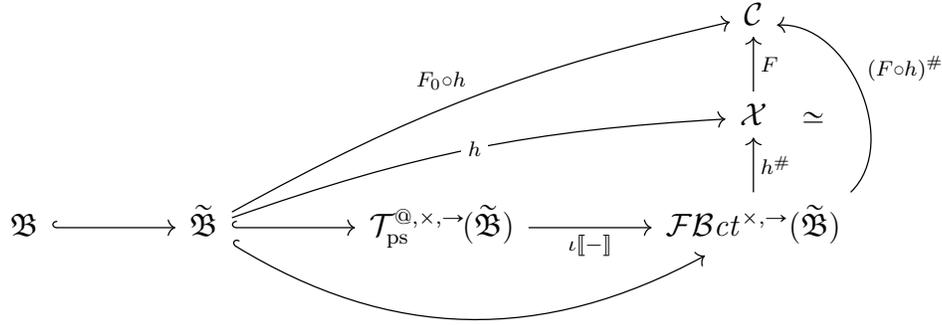
\begin{td}[column sep = 4em]
\: &
\: &
\: &
\altCat \\

\: &
\: &
\: &
\bicatX
\arrow[swap]{u}{F} \\

\baseTypes 
\arrow[hookrightarrow]{r} &
\allTypes{\baseTypes}
\arrow[bend right, hookrightarrow]{rr}
\arrow[bend left = 8]{uurr}{F_0 \circ h}
\arrow[bend left = 8]{urr}[description]{h}
\arrow[hookrightarrow]{r} &
\syncloneAtClosed{\allTypes\baseTypes}
\arrow{r}[swap]{\inc\sem{-}} &
\freeCartClosedBicat{\allTypes\baseTypes} 
\arrow[swap]{u}{\ext{h}}
\arrow[swap, bend right = 70]{uu}{\ext{(F \circ h)}}
\arrow[swap, bend right = 45, phantom, yshift=-1mm]{uu}{\simeq}
\end{td}
Now, the composite $F \circ \ext{h}$ is a cc-pseudofunctor, so by 
Lemma~\ref{lem:syntactic-model-uniqueness-up-to-equivalence} there exists 
an equivalence 
$\ext{(F_0 \circ h)} \simeq F \circ \ext{h} : \freeCartClosedBicat{\allTypes\baseTypes} \to 
\altCat$.
Denote this by 
$(\natTrans, \natCell) : 
	F \circ \ext{h} \To \ext{(F_0 \circ h)}$.
For any 1-cell $t : \Gamma \to A$
in $\freeCartClosedBicat{\allTypes\baseTypes}$, one therefore obtains an iso-commuting square
\begin{td}[column sep = 6em]
(F \circ \ext{h}){\Gamma} 
\arrow[phantom]{dr}[description]{\twocellIso{\natCell_{t}}}
\arrow[swap]{d}{\natTrans_{\Gamma}}
\arrow{r}{(F \circ \ext{h}){t}} &
(F \circ \ext{h}){A} 
\arrow{d}{\natTrans_A} \\
\ext{(F_0 \circ h)}{\Gamma}
\arrow[swap]{r}{\ext{(F_0 \circ h)}{t}} &
\ext{(F_0 \circ h)}{A}
\end{td}
Moreover, the naturality condition on $\natCell_t$ requires that, for any 
2-cell $\tau : t \To t' : \Gamma \to A$ in 
$\freeCartClosedBicat{\allTypes\baseTypes}$, the following commutes:
\begin{equation} \label{eq:yoneda-style-reduction-part-one}
\begin{tikzcd}[column sep = 9em]
\natTrans_A \circ (F \circ \ext{h}){(t)} 
\arrow[swap]{d}{\natCell_{t}}
\arrow{r}{\natTrans_A \circ (F \circ \ext{h}){(\tau)} } 
&
\natTrans_A \circ (F \circ \ext{h}){(t')} 
\arrow{d}{\natCell_{t'}} \\
\ext{(F_0 \circ h)}{(t)} \circ \natTrans_\Gamma
\arrow[swap]{r}
	{\ext{(F_0 \circ h)}{(\tau)}
		\circ \natTrans_\Gamma} &
\ext{(F_0 \circ h)}{(t')} \circ \natTrans_\Gamma
\end{tikzcd}
\end{equation}
But the cartesian closed structure of $\altCat$ is strict and the definition 
of the pseudofunctor $\ext{(F_0 \circ h)}$ only employs the canonical 
2-cells of the cc-bicategory structure, so 
$\ext{(F_0 \circ h)}{(\tau)}$
is the identity for every 2-cell $\tau$. To see this, one argues by induction 
on the definition of the cc-pseudofunctor $\ext{k}$ extending a map $k$ 
interpreting base types (Lemma~\ref{lem:free-cc-bicat-on-sig}).
It follows 
that~(\ref{eq:yoneda-style-reduction-part-one}) degenerates to the following:
\begin{equation} \label{eq:yoneda-style-reduction-part-two}
\begin{tikzcd}[column sep = 9em]
\natTrans_A \circ (F \circ \ext{h}){(t)} 
\arrow[swap]{d}{\natCell_{t}}
\arrow{r}{\natTrans_A \circ (F \circ \ext{h}){(\tau)} } 
&
\natTrans_A \circ (F \circ \ext{h}){(t')} 
\arrow{d}{\natCell_{t'}} \\
\ext{(F_0 \circ h)}{(t)} \circ \natTrans_\Gamma
\arrow[swap, equals]{r} &
\ext{(F_0 \circ h)}{(t')} \circ \natTrans_\Gamma
\end{tikzcd}
\end{equation}
Now, since $(\natTrans, \natCell)$ is an equivalence, every component
$\natTrans_X$ has a pseudoinverse. Let us denote this by 
$\psinv{\natTrans}_X$. 
From~(\ref{eq:yoneda-style-reduction-part-two}), one sees that the following 
commutes:
\begin{td}[column sep = 8em]
(F \circ \ext{h})(t)
\arrow[swap]{d}{\iso} 
\arrow{r}{(F \circ \ext{h}){(\tau)}} &
(F \circ \ext{h})(t') 
\arrow{d}{\iso} \\

(\psinv{\natTrans}_A \circ \natTrans_A) \circ
(F \circ \ext{h})(t)
\arrow[swap]{d}{\iso} 
\arrow[swap]{r}
	{(\psinv{\natTrans}_A 
		\circ \natTrans_A) 
		\circ (F \circ \ext{h}){(\tau)}} &
(\psinv{\natTrans}_A \circ \natTrans_A) \circ
(F \circ \ext{h})(t) 
\arrow{d}{\iso} \\

\psinv{\natTrans}_A 
	\circ 	\left(\natTrans_A 
	\circ (F \circ \ext{h})(t)\right) 
\arrow[swap]{d}{\psinv{\natTrans}_A \circ \natCell_t}
\arrow[swap]{r}
	{\psinv{\natTrans}_A 
		\circ \left(\natTrans_A 
		\circ (F \circ \ext{h}){(\tau)}\right)} &
\psinv{\natTrans}_A 
	\circ \left(\natTrans_A 
	\circ (F \circ \ext{h})(t')\right) 
\arrow{d}{\psinv{\natTrans}_A \circ \natCell_{t'}} \\

\psinv{\natTrans}_A \circ 
	\left( \ext{(F_0 \circ h)}(t) \circ \natTrans_\Gamma 
	\right) 
\arrow[equals]{r} &
\psinv{\natTrans}_A \circ 
	\left( \ext{(F_0 \circ h)}(t') \circ \natTrans_\Gamma 
	\right) 
\end{td}
One thereby sees that $(F \circ \ext{h}){\tau}$
is completely determined by a composite of 2-cells, none of which depend on 
$\tau$. 

\begin{prooflesspropn}
Let $\ccBicat{\bicatX}$ be a cc-bicategory , 
$\ccBicat\altCat$ be a 2-category with strict products and exponentials, 
and $(F, \prodPres, \expPres) : 
	\ccBicat{\bicatX} \to \ccBicat{\altCat}$
be any cc-pseudofunctor. Then
if $h : \allTypes{\baseTypes} \to \bicatX$ is the canonical extension
of an interpretation $\baseTypes \to \bicatX$ and 
$\tau : t \To t'$ is any 2-cell in $\freeCartClosedBicat{\allTypes\baseTypes}$, the 2-cell
$(F \circ \ext{h})(\tau)$ in $\altCat$
is completely determined by $t$ and $t'$. Hence, for any parallel pair 
of 2-cells $\tau, \sigma : t \To t'$ in $\freeCartClosedBicat{\allTypes\baseTypes}$, one has 
the equality 
$(F \circ \ext{h})(\tau) = (F \circ \ext{h})(\sigma)$.
\end{prooflesspropn}

Together with Proposition~\ref{prop:yoneda-coherence-for-cc-bicats}, one 
obtains the local coherence of $\freeCartClosedBicat{\allTypes\baseTypes}$, which completes our 
alternative proof of Theorem~\ref{thm:main-result-on-free-bicat}.

\begin{mythm} \label{thm:main-result-via-strictifying-first}
For any set of base types $\baseTypes$
and any pair of parallel 2-cells $\tau, \sigma : t \To t'$ in
$\freeCartClosedBicat{\allTypes\baseTypes}$, the equality
$\tau \equiv \sigma$ holds. 
\begin{proof}
Instantiate the preceding proposition with $h := \inc : \allTypes\baseTypes 
\hookrightarrow \freeCartClosedBicat{\allTypes\baseTypes}$ the inclusion and $F$ the biequivalence 
between a cc-bicategory and a 2-category with strict products and exponentials 
arising from 
Proposition~\ref{prop:yoneda-coherence-for-cc-bicats}. 
Note that $\ext{\inc} \simeq \id_{\freeCartClosedBicat{\allTypes\baseTypes}}$ by 
Lemma~\ref{lem:syntactic-model-uniqueness-up-to-equivalence}, so that
$F \circ \ext{\inc}$ is a biequivalence.
Then $F \circ \ext{\inc}$ is locally 
fully faithful, so 
$(F \circ \ext{\inc})(\tau) = 
(F \circ \ext{\inc})(\sigma)$
if and only if $\tau \equiv \sigma$. The result then follows from the preceding
proposition.
\end{proof}
\end{mythm}

Since $\freeCartClosedBicat{\allTypes\baseTypes} \simeq \syncloneAtClosed{\allTypes\baseTypes}$, this 
entails the 
local coherence of $\syncloneAtClosed{\sig}$. One therefore recovers 
Theorem~\ref{thm:main-theorem}. 

We end with some comments on the argument just presented.
First, as it stands it is not constructive. We make use of the coherence 
theorem for fp-bicategories (Proposition~\ref{prop:power-coherence}), for which 
one chooses a pseudoinverse to the inclusion of a bicategory into its 
image under the Yoneda embedding. This choice is only determined up to 
equivalence, so one does not obtain an explicit witness for the product 
structure.
Second, the argument relies crucially on the interplay between weak and strict 
structure. We use the strictness of $\Hom(\baseCat, \Cat)$ to obtain a 
strict cc-bicategory biequivalent to our original one, and then we use the 
strictness of this bicategory to 
degenerate~(\ref{eq:yoneda-style-reduction-part-one}) 
into~(\ref{eq:yoneda-style-reduction-part-two}). It is, therefore, a strategy 
that is only available in the higher-categorical setting.

\chapter{Conclusions} 
\label{chap:conclusions}

We leave a full investigation of the applications of the 
development in this thesis for future work. We do note, however, that the 
problem we posed in the introduction now disappears.

Consider a structure definable in any cartesian closed category. Examples 
include 
the canonical comonoid structure on any object, or the monoid structure on any 
endo-exponential. This definition is witnessed by a $\stlc$-term up to 
$\beta\eta$-equality, and hence---by Proposition~\ref{prop:stlc-up-to-iso}---by 
a $\langCartClosed$-term over the 
same signature, with $\beta\eta$-equalities replaced by rewrites. (Since we 
explicitly construct the 
correspondence between $\stlc$-terms and $\langCartClosed$-terms, this 
construction can be done via a terminating decision procedure.)
 These rewrites will provide the data required to define a bicategorical 
 version of the structure under consideration. Theorem~\ref{thm:main-theorem} 
 then entails that the required coherence axioms must hold. One thereby obtains 
 the following principle.

\begin{myprinciplelocal}
To show that a \emph{pseudo} structure may be constructed in any cartesian 
closed bicategory, it suffices to show that its strict version---that is, the 
image 
of the corresponding $\langCartClosed$-term in $\stlc$---may be constructed in 
any cartesian closed category.
\end{myprinciplelocal}

Applying this principle immediately entails the following results.

\begin{prooflesspropnlocal}
For any cc-bicategory,
\begin{enumerate}
\item Every object has a canonical commutative pseudo-comonoid structure, and
\item Every endo-exponential has a canonical pseudomonoid structure. \qedhere
\end{enumerate}
\end{prooflesspropnlocal}

\subsection*{Further work}

There are many interesting avenues for further work; we mention a few here.

\paragraph*{Extensions to $\langCartClosed$.}
It is natural to consider incorporating further type-theoretic constructions 
into $\langCartClosed$. One example would be sum types, 
corresponding to bicategorical coproducts. Extending the local coherence proof 
to this type theory would likely require a bicategorical development of 
\emph{Groethendieck logical relations}~\cite{Fiore1999sums}, with possible 
connections to the theory of stacks. 
A more ambitious development would be 
the inclusion of 
Martin-L{\"o}f style dependent types~\cite{MartinLof}. This would be 
particularly intriguing as the interpretation of 
these constructions in locally cartesian closed categories is, properly 
speaking, bicategorical~\cite{Clairambault2014}.

From a different perspective, Pitts has suggested considering the theory 
of fixpoints. In an unpublished manuscript~\cite{Pitts1987}, Pitts considered a 
calculus for \Def{initial fixpoint categories} (\Def{IFP-categories}): 
2-categories equipped with finite products 
and a notion of `initial algebra' on every endomorphism of the form
$A \xra{\seq{\id_A, a}} A \times B \xra{f} B$, representing a formal fixpoint 
construction. Other important examples in a similar vein 
include~\Def{algebraically complete} categories~\cite{Freyd1991}, or 
\Def{iteration (2-)theories}~\cite{Esik1999, Bloom2001}. 
The fact that 
bicategories represent a natural 
setting for `formal category theory' suggests considering constructions 
of type-theoretic interest (such as fixpoints) as well as constructions of 
category-theoretic interest (such as monads) as particular constructions within 
$\langBiclone$. 

An orthogonal line of development would be towards higher levels of 
categorical structure. One might, for example, extend to tricategories; 
restricting to unary contexts would recover a type theory for monoidal 
bicategories. (An alternative approach to the same result would be to introduce 
a \emph{linear} version of $\langBiclone$). It may even be possible to 
inductively 
generate higher levels of structure to recover some form of 
$\infty$-category. For these developments to be principled, the first 
consideration ought to be the appropriate correlate of biclones.

\paragraph*{Applications to higher category theory.}
Each extension to the type theory raises the question of its 
coherence. As 
outlined in the introduction to Chapter~\ref{chap:nbe}, there 
is a wealth of literature studying various forms of 
normalisation-by-evaluation for extensions to the simply-typed lambda calculus. 
It is plausible that their bicategorical correlates would lift to 
extensions of $\langCartClosed$. 
%
More speculatively, one might hope that by constructing higher-dimensional type 
theories and examining their relationship to well-understood classical type 
theories (in the style of Section~\ref{sec:STLC-vs-pseudoSTLC}, for 
instance), one may gain a better understanding of where 
coherence can be expected and---in the cases it cannot---why it fails.

This thesis also lays the groundwork for 
bicategorifying further category theoretic results. For instance, 
the conservative extension result of~\cite[\S 3]{Fiore2002sums} shares many 
tools with the normalisation-by-evaluation argument of~\cite{Fiore2002}, 
such as glueing and the relative hom-functor. It should be possible, therefore, 
to extend the bicategorical theory presented here to show that 
cc-bicategories are a conservative extension of 
fp-bicategories. 

\paragraph*{Higher-dimensional universal algebra.}
Moving away from type-theoretic concerns, there remains the question of the 
universal algebra associated to (mono-sorted) biclones. In the classical 
setting, it is well-known that the three components of the
monad--Lawvere  theory--clone triad are all equivalent. Biclones 
appear to represent one corner of the bicategorical version of this triad: 
whether pseudomonads and some bicategorical notion of Lawvere theory complete 
the picture remains to be seen.

\appendix

\part{~~Appendices~~}

\renewcommand{\thechapter}{\Alph{chapter}}
\renewcommand{\thesection}{\Alph{chapter}.\arabic{section}}
\renewcommand{\thefigure}{\thechapter.\arabic{figure}}

\chapter{An index of free structures and syntactic models}
\label{chap:summary:free-structures}

We summarise the various bicategorical free constructions and syntactic models 
employed throughout 
this thesis. As a rule of thumb, we use $\synclonesymbol$ to 
denote biclones (and their nuclei,~\ie restrictions to unary contexts) and $\termCatSymbol$ to denote bicategories.

\begin{table}[h]
{
\renewcommand{\arraystretch}{1.4}
\centering
\begin{tabular}{cccc}
\multicolumn{1}{l}{\textbf{Chapter~\ref{chap:biclone-lang}}} & 
\multicolumn{1}{l}{}                                                            
                                                        & 
\multicolumn{1}{l}{}                                         & 
\multicolumn{1}{l}{}                                \\ 
$\freeBiclone{\graph}$                                       & free biclone on 
a 2-multigraph                                                                  
                                        & 
Construction~\ref{constr:free-biclone}                       & 
\:\: \:\: p. \pageref{constr:free-biclone}                       \\
$\freeBicat{\graph}$                                         & free bicategory 
on a 
2-graph                                                                         
                                   &
 Lemma~\ref{lem:free-bicategory-on-a-multigraph}              & 
\:\: \:\: p. \pageref{lem:free-bicategory-on-a-multigraph}       \\
$\synclone{\graph}$                                          & syntactic 
biclone of $\langBiclone$ on a 
2-multigraph                                                                    
               &
 Construction~\ref{constr:langbiclone-syntactic-model}        & 
\:\: \:\: p. \pageref{constr:langbiclone-syntactic-model}        \\
$\urestrict{\synclone{\graph}}$                              & syntactic 
bicategory of $\langBicat$ on a 
2-graph                                                                         
              &
 Construction~\ref{constr:bicat-termcat}                      & 
\:\: \:\: p. \pageref{constr:bicat-termcat}                      \\
$\hirclone(\graph)$                                          & syntactic 
biclone of $\hir$ on a 
2-multigraph                                                                    
                       &
 Construction~\ref{constr:free-2-clone}                       & 
\:\: \:\: p. \pageref{constr:free-2-clone}                       \\
\multicolumn{1}{l}{}                                         & 
\multicolumn{1}{l}{}                                                            
                                                        & 
\multicolumn{1}{l}{}                                         &                  
                                   \\
\multicolumn{1}{l}{\textbf{Chapter~\ref{chap:fp-lang}}}      & 
\multicolumn{1}{l}{}                                                            
                                                        & 
\multicolumn{1}{l}{}                                         &                  
                                   \\
$\freeCartBiclone{\sig}$                                         & free cartesian 
biclone on a 
$\langCart$-signature                                                           
                            &
 Construction~\ref{constr:free-cart-biclone}                  & 
\:\: \:\: p. \pageref{constr:free-cart-biclone}                  \\
$\freeCartBicat{\sig}$                                           & free 
fp-bicategory on a unary 
$\langCart$-signature                                                           
                          &
 Lemma~\ref{lem:free-fp-bicat-from-free-biclone}              & 
\:\: \:\: p. \pageref{lem:free-fp-bicat-from-free-biclone}       \\
$\syncloneCart{\sig}$                                        & syntactic 
biclone of $\langCart$ on a 
$\langCart$-signature                                                           
                  &
 Construction~\ref{constr:langcart-syntactic-model}           & 
\:\: \:\: p. \pageref{constr:langcart-syntactic-model}           \\
$\urestrict{\syncloneCart{\sig}}$                            & 
\begin{tabular}[c]{@{}c@{}}syntactic model of type theory obtained by 
\\[-4pt]
restricting 
$\langCart$ to unary contexts \end{tabular} 
                                                  
                                  &
 Theorem~\ref{thm:unary-contexts-fp-bicat}                    & 
\:\: \:\: p. \pageref{thm:unary-contexts-fp-bicat}               \\
$\termCatContextExt(\sig)$                                   & 
\begin{tabular}[c]{@{}c@{}}extension of $\urestrict{\syncloneCart{\sig}}$ 
with\\[-4pt] context extension product structure\end{tabular}       & 
Construction~\ref{constr:context-cart-termcat}               & 
\:\: \:\: p. \pageref{constr:context-cart-termcat}               \\
\multicolumn{1}{l}{}                                         & 
\multicolumn{1}{l}{}                                                            
                                                        & 
\multicolumn{1}{l}{}                                         &                  
                                   \\
\multicolumn{1}{l}{\textbf{Chapter~\ref{chap:ccc-lang}}}      & 
\multicolumn{1}{l}{}                                                            
                                                        & 
\multicolumn{1}{l}{}                                         &                  
                                   \\
$\freeCartClosedBiclone{\sig}$                                         & free cartesian 
closed biclone on a 
$\langCartClosed$-signature                                                     
                     &
 Construction~\ref{constr:free-cc-biclone}                    & 
\:\: \:\: p. \pageref{constr:free-cc-biclone}                    \\
$\freeCartClosedBicat{\sig}$                                           & free 
cc-bicategory on a 
$\langCartClosed$-signature                                                     
                                &
 Construction~\ref{constr:free-cc-bicat}                      & 
\:\: \:\: p. \pageref{constr:free-cc-bicat}                      \\
$\syncloneCartClosed{\sig}$                                  & syntactic 
biclone of $\langCartClosed$ on a 
$\langCartClosed$-signature                                                     
            &
 Construction~\ref{constr:langcartclosed-syntactic-model}     & 
\:\: \:\: p. \pageref{constr:langcartclosed-syntactic-model}     \\
$\increstrict{\syncloneCartClosed{\sig}}$                      & 
nucleus of 
$\syncloneCartClosed{\sig}$ 
                                   &
 Construction~\ref{constr:restrict-of-langcartclosed-biclone} & 
\:\: \:\: p. \pageref{constr:restrict-of-langcartclosed-biclone} \\
$\syncloneAtClosed{\sig}$                                    & 
\begin{tabular}[c]{@{}c@{}}extension of 
$\increstrict{\syncloneCartClosed{\sig}}$ 
with\\[-4pt] context extension product structure\end{tabular} & 
Construction~\ref{constr:ccc-termcat}                        & 
\:\: \:\: p. \pageref{constr:ccc-termcat}                        \\ \hline
\end{tabular}
\caption{An index of free constructions and syntactic models 
\label{table:all-free-constructions}}
}
\end{table}

\let\oldcleartorecto\cleartorecto 
\let\cleartorecto\newpage         
\chapter{Cartesian closed structures} \label{chap:cartesian-closed-structures}
\let\cleartorecto\oldcleartorecto

We summarise the cartesian closed structures of $\Hom(\baseCat, \Cat)$ and 
$\gl{F}$.

\paragraph*{Cartesian closed structure on $\Hom(\baseCat, \Cat)$.}

Let $\baseCat$ be any 2-category. Then the 2-category $\Hom(\baseCat, \Cat)$ 
has finite products given pointwise and exponentials given as in the following 
table:

\begin{table}[h]
\centering
\renewcommand{\arraystretch}{2}
\begin{tabular}{c | c}
Exponential $\altexp{P}{Q}$ & 
$\lambda X^\baseCat \bind \Hom(\baseCat, \Cat)( \Yon X \times P, Q)$  \\ \hline

Evaluation 1-cell $\glued{\eval}_{P,Q}$ & 
$\lambda X^\baseCat \bind 
\lambda (\natTrans, \natCell)^{\Yon X \times P \To Q} \bind 
\lambda p^{PX} \bind 
\natTrans(X, \Id_X, p)$ 
\\ \hline

$\Lambda (\altNat, \altCell)^{R \times P \To Q}$ &
\begin{tabular}[c]{@{}l@{}}
$\lambda X^\baseCat \bind
\lambda r^{RX} \bind
\lambda A^\baseCat \bind
\lambda (h,p)^{\Yon (X, A) \times PA} \bind
\altNat\big(X, (Rh)(r), p\big)$ 
\\[-1em]
\quad with naturality witnessed by by Lemmas~\ref{lem:def-of-Lambda-part-one} 
and~\ref{lem:lem:def-of-Lambda-part-two}
\end{tabular} 
\\ \hline

Counit $\evalMod_{P,Q}(\altNat, \altCell)$ & 
$\lambda X^\baseCat \bind
\lambda (r,p)^{RX \times PX} \bind
\altNat\big(X, (\psi^R)^{-1}(r), p\big)$ \\ 
\hline

$\transExp{\modif}$ & 
defined by diagram~(\ref{eq:def-of-transModif})
\end{tabular}
\caption{Exponential structure in $\Hom(\baseCat, \Cat)$, from 
Section~\ref{sec:presheaves-cartesian-closed}
\label{table:ccc-structure-on-presheaves}}
\end{table}
Moreover, for a pseudofunctor $P : \op\baseCat \to \Cat$ and object 
$X \in \baseCat$ the exponential 
$\altexp{\Yon X}{P}$ in $\Hom(\op\baseCat, \Cat)$ is given by $P(- \times X)$, 
with structure summarised in Table~\ref{table:exponentiating-by-representable}.
\begin{table}[h]
\centering
\renewcommand{\arraystretch}{2}
\begin{tabular}{c | c}
Evaluation 1-cell $\glued{\eval}_{P,Q}$ & 
\begin{tabular}[c]{@{}l@{}}
$\lambda B^{\baseCat} \bind 
\lambda (p, h)^{P(B \times X) \times \baseCat(B,X)} \bind
P\big( \seq{\Id_B, h} \big)(p)$  \\[-1em]
\quad with naturality witnessed by Lemma~\ref{lem:YX-exponential-eval-map}
\end{tabular} 
\\ \hline

$\Lambda (\natTrans, \natCell)^{R \times \Yon X \To P}$ &
\begin{tabular}[c]{@{}l@{}}
$\lambda B^\baseCat \bind
\lambda r^{RB} \bind
\natTrans_{B \times X}\big( R(\pi_1)(r), \pi_2 \big)$ 
\\[-1em]
\quad with naturality witnessed by 
Corollary~\ref{cor:YX-exponentiating-lambda-mapping} 
\end{tabular} 
\\ \hline

Counit $\evalMod(\natTrans, \natCell)$ & 
defined by diagram~(\ref{eq:Yon-X-exponential-counit}) \\ 
\hline

$\transExp{\modif}$ & 
defined by diagram~(\ref{eq:Yon-X-exponential-trans-modif})
\end{tabular}
\caption{Exponential structure in $\Hom(\baseCat, \Cat)$, from 
Section~\ref{sec:exponentiating-by-Yon} 
\label{table:exponentiating-by-representable}}
\end{table}

\paragraph*{Cartesian closed structure on $\gl{\glueFun}$.}

Let $(\glueFun, \prodPres) : \fpBicat{\baseCat} \to \fpBicat{\altCat}$ be an 
fp-pseudofunctor between cc-bicategories and suppose that $\altCat$ 
has all pullbacks. Then $\gl{\glueFun}$ is cartesian closed, with structure 
given as 
in Tables~\ref{table-two}--\ref{table:glueing-ccc-structure}.

\begin{table}[h]
\centering
\renewcommand{\arraystretch}{2}
\begin{tabular}{c | c}
Product $\prod_i (C_i, c_i, B_i)_i$ & 
$\big( \prod_i C_i, \prodPres \circ \prod_i c_i, \prod_i B_i \big)$  \\ \hline

Projection 1-cells $\glued{\pi}_k$ & 
$(\pi_k, \mu_k, \pi_k)$ for $\mu_k$ defined in~(\ref{eq:def-of-mu-k}) \\ \hline

$n$-ary tupling $\seq{\glued{t}_1, \dots, \glued{t}_n}$ 
for $\glued{t}_i := (t_i, \alpha_i, s_i)$ &
$(\seq{\ind{t}}, 
\glTup{\ind{\alpha}},
\seq{\ind{s}})$ for 
$\glTup{\ind{\alpha}}$ defined in~(\ref{eq:def-of-glTup}) \\ \hline

Counit $\glued{\epsilonTimes{}}$ & 
$k$th component is
$(\epsilonTimesInd{k}{\ind{f}}, \epsilonTimesInd{k}{\ind{g}})$ \\  \hline

\renewcommand{\arraystretch}{1.1}
\begin{tabular}[c]{@{}l@{}}
$\transTimes{\glued{\tau}_1, \dots, \glued{\tau}_n}$ for 
$\glued{\tau}_i := (\tau_i, \sigma_i) : \glued{\pi}_k \circ \glued{u} \To 
\glued{t}_i$  \\ 
\quad $(i=1,\dots,n)$
\end{tabular}  & 
$\big(\transTimes{\tau_1, \dots, \tau_n}, 
\transTimes{\sigma_1, \dots, \sigma_n}\big)$
\end{tabular}
\caption{Product structure in $\gl{\glueFun}$, from Section~\ref{sec:glued-products-constructed} \label{table-two}}
\end{table}

\begin{table}[ht]
\centering
\renewcommand{\arraystretch}{2}
\begin{tabular}{c | c}
Exponential $\exp{(C,c,B)}{(C',c',B')}$ & 
$(\glexp{C}{C'}, p_{c,c'}, \expobj{B}{B'})$ defined by the 
pullback~(\ref{eq:exponential-in-glued-bicat})  \\ \hline

Evaluation 1-cell $\glued{\eval}_{\glued{C}, \glued{C'}}$ & 
{
\begin{tabular}[c]{@{}l@{}}
$(\eval_{C, C'} \circ (q_{c, c'} \times C), 
\glued{E}_{\glued{C}, \glued{C'}}, 
\eval_{B,B'})$  \\[-1em]
\quad for $\evalMod_{\glued{C}, \glued{C'}}$ defined 
in~(\ref{eq:def-of-eval-map-witnessing-2-cell}) and (\ref{eq:def-of-glued-E})
\end{tabular} } 
\\ \hline

$\glued{\lambda}(t, \alpha, s)$ &
{
\begin{tabular}[c]{@{}l@{}}
$(\glLam(t), \Gamma_{c,c'}, \lambda s)$ for 
$\glLam(t)$ and $\Gamma_{c,c'}$ defined by \\[-1em] 
\quad UMP of  
pullback applied to $L_\alpha$~(\ref{eq:L-alpha-ump-diagram})
\end{tabular} }
\\ \hline

Counit $\glued{\epsilonExp}$ & 
$(\gluedCo{}, \epsilonExp)$ for $\gluedCo{}$ defined 
in~(\ref{eq:def-of-counit}) \\ \hline

$\transExp{\glued{\tau}}$ for $\glued{\tau} := (\tau, \sigma)$ & 
\begin{tabular}[c]{@{}l@{}}
$\big(\altTrans{\tau}, \transExp{\sigma} \big)$ for $\altTrans{\tau}$ defined 
by UMP of pullback  
\\[-1em]
\quad applied to fill-in defined 
in~(\ref{eq:def-of-Sigma-1-and-Sigma-2})
\end{tabular} 
\end{tabular}
\caption{Exponential structure in $\gl{\glueFun}$, from Section~\ref{sec:exponentials-in-glueing-constructed}
\label{table:glueing-ccc-structure}}
\end{table}

\vfill

\chapter{The type theory and its semantic interpretation} \label{chap:full-language}

\section{The type theory \texorpdfstring{$\langCartClosed$}{}}

%
%
%

\newcommand{\skiplength}{0cm}

\setlength{\floatsep}{5pt plus 0pt minus 2.0pt} 
\setlength{\textfloatsep}{0pt}

\vspace{\skiplength}

Fix a $\langCartClosed$-signature $\sig = (\baseTypes, \graph)$~(Definition~\ref{def:lang-cart-closed-sig} on page~\pageref{def:lang-cart-closed-sig}). 
We give the rules for the full type theory $\langCartClosed$. The type theories 
$\langBiclone$ and $\langCart$ are fragments of $\langCartClosed$, and the type 
theories $\langBicat$ and $\urestrict{\langCart}$ are respectively obtained by 
restricting $\langBiclone$ and $\langCart$ to unary contexts.

\vspace{\skiplength}
{\small	
\begin{rules}
\unaryRule	{\faketext}
			{\diamond \mathrm{\: ctx}}
			{\qquad\qquad} 
\binaryRule	{\Gamma \mathrm{\: ctx}}
			{x \notin \dom(\Gamma)}
			{\Gamma, x : A \mathrm{\: ctx}} 
			{$\big(A \in \allTypes{\baseTypes}\big)$} \vspace{-\treeskip}
			
\caption{Rules for contexts}
\end{rules}
}
\vspace{\skiplength}

{\small

\begin{rules}
\unaryRule	{\faketext}
			{x_1 : A_1, \dots, x_n : A_n \vdash x_k : A_k}
			{var $(1 \leq k \leq n)$}
\unaryRule	{c \in \graph(A_1, \dots, A_n; B)}
			{x_1 : A_1, \dots, x_n : A_n \vdash c(x_1, \dots, x_n) : B}
			{const}
		
\binaryRule	{x_1 : A_1, \dots, x_n : A_n \vdash t : B}
			{(\Delta \vdash u_i : A_i)_{i= 1, \dots, n}}
			{\Delta \vdash \hcomp{t}{x_1 \mapsto u_1, \dots, x_n \mapsto u_n} : B}
			{horiz-comp}
\unaryRule	{\faketext}
			{p : \prodop_n(A_1, \dots, A_n) \vdash \pi_k(p) : A_k}
			{$k$-proj ($1 \leq k \leq n$)}\vspace{0.5\treeskip}

\unaryRule	{\Gamma \vdash t_1 : A_1 \quad \dots \quad \Gamma \vdash t_n : A_n}
			{\Gamma \vdash \pair{t_1, \dots, t_n} : \prodop_n (A_1, \dots, A_n)}
			{$n$-tuple}

\unaryRule 	{\Gamma, x : A \vdash t : B} 
			{\Gamma \vdash \lam{x}{t} : \exptype{A}{B}} 
			{lam} 
\qquad
\unaryRule 	{\faketext} 
			{f : \exptype{A}{B}, x : A \vdash \evalterm(f,x) : B} 
			{eval} 
\vspace{-\treeskip}

\caption{Introduction rules for terms}
\end{rules}

\begin{rules}

\begin{prooftree}
\AxiomC{$x_1 : A_1, \dots, x_n : A_n \vdash t : B$}
\RightLabel{\scriptsize $\subid{}$-intro}
\UnaryInfC{$x_1 : A_1, \dots, x_n : A_n \vdash \subid{t} : \birewrite{t}{\hcomp{t}{x_i \mapsto x_i}} : B$}
\noLine
\UnaryInfC{$x_1 : A_1, \dots, x_n : A_n \vdash \subid{t}^{-1} : \birewrite{\hcomp{t}{x_i \mapsto x_i}}{t} : B$}
\end{prooftree}

\begin{prooftree}
\AxiomC{$x_1 : A_1, \dots, x_n : A_n \vdash x_k : A_k$}
\AxiomC{$(\Delta \vdash u_i : A_i)_{i = 1, \dots, n}$}
\RightLabel{\scriptsize $\indproj{k}{}$-intro $(1 \leq k \leq n)$}
\BinaryInfC{$\Delta \vdash \indproj{k}{u_1, \dots, u_n} : \birewrite{\hcomp{x_k}{x_i \mapsto u_i}}{u_k} : A_k$}
\noLine
\UnaryInfC{$\Delta \vdash \indproj{-k}{u_1, \dots, u_n} :  \birewrite{u_k}{\hcomp{x_k}{x_i \mapsto u_i}} : A_k$}
\end{prooftree}

\begin{small}
\begin{prooftree}
\AxiomC{$(\Delta \vdash u_j : A_j)_{j = 1, \dots m}$}
\noLine
\UnaryInfC{$(x_1 : A_1, \dots, x_m : A_m \vdash v_i : B_i)_{i = 1, \dots, n}$}
\noLine
\UnaryInfC{$y_1 : B_1, \dots, y_n : B_n \vdash t : C$}
\RightLabel{\scriptsize $\assoc{}$-intro}
\UnaryInfC{$\Delta \vdash \assoc{t, \ind{v}, \ind{u}} : \rewrite{\hcomp{\hcomp{t}{y_i \mapsto v_i}}{x_j \mapsto u_j}}{\hcomp{t}{y_i \mapsto \hcomp{v_i}{x_j \mapsto u_j}}} : C$}
\noLine
\UnaryInfC{$\Delta \vdash \assoc{t, \ind{v}, \ind{u}}^{-1} : \rewrite{\hcomp{t}{y_i \mapsto \hcomp{v_i}{x_j \mapsto u_j}}}{\hcomp{\hcomp{t}{y_i \mapsto v_i}}{x_j \mapsto u_j}} : C$}
\end{prooftree}
\end{small}

\caption{Introduction rules for structural rewrites}
\end{rules}

}

{\small

\begin{rules}
\unaryRule	{\Gamma \vdash t : A}
			{\Gamma \vdash \id_t : \birewrite{t}{t} : A}
			{$\id$-intro}
\unaryRule		{\constrewr \in \graph(A_1, \dots, A_n; B)(c, c')}
				{x_1: A_1, \dots, x_n : A_n \vdash \constrewr(x_1, \dots, x_n) : \rewrite{c(x_1, \dots, x_n)}{c'(x_1, \dots, x_n)} : B}
				{2-const}
\unaryRule	{\Gamma \vdash t_1 : A_1 \qquad \dots \qquad \Gamma \vdash t_n : A_n}
			{\Gamma \vdash \epsilonTimesInd{k}{t_1, \dots, t_n} : \rewrite{\hcomp{\pi_k}{\pair{t_1, \dots, t_n}}}{t_k} : A_k} 
			{$\epsilonTimesInd{k}{}$-intro ($1 \leq k \leq n$)}

\begin{small}			
\binaryRule	{\Gamma \vdash u : \prodop_n(A_1, \dots, A_n)}
			{(\Gamma \vdash \alpha_i : \rewrite{\hcomp{\pi_i}{u}}{t_i} : A_i)_{i=1, \dots, n}}
			{\Gamma \vdash \transTimes{\alpha_1, \dots, \alpha_n} : \rewrite{u}{\pair{t_1, \dots, t_n}} : \prodop_n(A_1, \dots, A_n)}
			{$\transTimes{\alpha_1, \dots, \alpha_n}$-intro}
\end{small} \vspace{-\treeskip}
\begin{prooftree} 
\AxiomC{$\Gamma, x : A \vdash t : B$} 
\RightLabel{\scriptsize $\epsilonExpRewr{}$-intro} 
\UnaryInfC{$\Gamma, x : A \vdash \epsilonExpRewr{t} : \rewrite{\genevalterm{\wkn{(\lam{x}{t})}{x}, x}}{t} : B$} 
\end{prooftree} 

\begin{prooftree} 
\AxiomC{$\Gamma, x : A \vdash t : B$ \:\quad\: $\Gamma \vdash u : \exptype{A}{B}$} 
\noLine 
\UnaryInfC{$\Gamma, x : A \vdash \alpha : \rewrite{\genevalterm{\wkn{u}{x}, x}}{t} : B$} 
\RightLabel{{\scriptsize $\transExp{x \bind \alpha}$-intro}} 
\UnaryInfC{$\Gamma \vdash \transExp{x \bind \alpha} : \rewrite{u}{\lam{x}{t}} : \exptype{A}{B}$} 
\end{prooftree} 
\vspace{-0.5\treeskip}

\caption{Introduction rules for basic rewrites}
\end{rules}

\begin{rules}
\binaryRule	{\Gamma \vdash \tau : \rewrite{t}{t'} : A}
			{\Gamma \vdash \tau' : \rewrite{t'}{t''} : A}
			{\Gamma \vdash \tau' \vertsub{t, t',t''} \tau : \rewrite{t}{t''} : A}
			{vert-comp}
\binaryRule	{x_1 : A_1, \dots, x_n : A_n \vdash \tau : \rewrite{t}{t'} : B}
			{(\Delta \vdash \sigma_i : \rewrite{u_i}{u'_i} : A_i)_{i = 1, \dots, n}}
			{\Delta \vdash \horizComp{\tau}{x_i \mapsto \sigma_i} : \rewrite{\hcomp{t}{x_i \mapsto u_i}}{\hcomp{t'}{x_i \mapsto u_i'}} : B}
			{horiz-comp}
\caption{Composition operations for rewrites}
\end{rules}

\begin{rules}
	\unaryRule	{\Gamma \vdash t_1 : A_1 \qquad \dots \qquad \Gamma \vdash t_n : A_n}
			{\Gamma \vdash \epsilonTimesInd{-k}{t_1, \dots, t_n} : \rewrite{t_k}{\hcomp{\pi_k}{\pair{t_1, \dots, t_n}}} : A_k} 
			{$\epsilonTimesInd{-k}{}$-intro $(1 \leq k \leq n)$}
\unaryRule	{\Gamma \vdash t : \prodop_n(A_1, \dots, A_n)}
			{\Gamma \vdash \etaTimes{t}^{-1} : \rewrite{\pair{\hcomp{\pi_1}{t}, \dots, \hcomp{\pi_n}{t}}}{t} : \prodop_n(A_1, \dots, A_n)}
			{$\etaTimes{}^{-1}$-intro}
			
	\unaryRule	{\Gamma \vdash u : \exptype{A}{B}}
			{\Gamma \vdash \etaExp{u}^{-1} : \rewrite{\lam{x}{\genevalterm{\wkn{u}{x}, x}}}{u} : \exptype{A}{B}}
			{$\etaExp{}^{-1}$-intro}
\begin{bprooftree} 
\AxiomC{$\Gamma, x : A \vdash t : B$} 
\RightLabel{\scriptsize $\epsilonExpRewr{}^{-1}$-intro} 
\UnaryInfC{$\Gamma, x : A \vdash \epsilonExpRewr{t}^{-1} : \rewrite{t}{\genevalterm{\wkn{(\lam{x}{t})}{x}, x}} : B$} 
\end{bprooftree}\vspace{\treeskip} 
	
	\caption{Introduction rules for pseudo cartesian closed structure}
\end{rules}
}
\vspace{\skiplength}


{\small
\begin{rules}

\unaryRule	{\Gamma \vdash \tau : \rewrite{t}{t'} : A}
			{\Gamma \vdash \tau \vert \id_t \equiv \tau : \rewrite{t}{t'} : A}
			{$\vert$-right-unit}
\unaryRule	{\Gamma \vdash \tau : \rewrite{t}{t'} : A}
{\Gamma \vdash \tau \equiv \id_{t'} \vert \tau : \rewrite{t}{t'} : A}
{$\vert$-left-unit}

\trinaryRule	{\Gamma \vdash \tau'' : \rewrite{t''}{t'''} : A}
				{\Gamma \vdash \tau' : \rewrite{t'}{t''} : A}
				{\Gamma \vdash \tau : \rewrite{t}{t'} : A}
				{\Gamma \vdash (\tau'' \vert \tau') \vert \tau \equiv \tau'' \vert (\tau' \vert \tau) : \rewrite{t}{t'''} : A}
				{$\vert$-assoc}
\caption{Categorical structure of vertical composition}
\end{rules}

\begin{rules}

\binaryRule	{x_1 : A_1, \dots, x_n : A_n \vdash t : B}
			{(\Delta \vdash u_i : A_i)_{i = 1, \dots, n}}
			{\Delta \vdash \hcomp{\id_t}{x_i \mapsto u_i} \equiv \id_{\hcomp{t}{x_i \mapsto u_i}} : \rewrite{\hcomp{t}{x_i \mapsto u_i}}{\hcomp{t}{x_i \mapsto u_i}} : B}
			{$\id$-preservation} \vspace{-\treeskip}

\begin{footnotesize}
\begin{prooftree}
\AxiomC{$x_1 : A_1, \dots, x_n : A_n \vdash \tau : \rewrite{t}{t'} : B$}
\noLine
\UnaryInfC{$x_1 : A_1, \dots, x_n : A_n \vdash \tau' : \rewrite{t'}{t''} : B$}
\AxiomC{$(\Delta \vdash \sigma_i : \rewrite{u_i}{u_i'} : A_i)_{i=1,\dots, n}$}
\noLine
\UnaryInfC{$(\Delta \vdash \sigma_i' : \rewrite{u_i'}{u_i''} : A_i)_{i=1,\dots, n}$}

\RightLabel{{\scriptsize interchange}}
\BinaryInfC{$\Delta \vdash \horizComp{\tau'}{x_i \mapsto \sigma_i'} \vert \horizComp{\tau}{x_i \mapsto \sigma_i} \equiv \horizComp{(\tau' \vert \tau)}{x_i \mapsto \sigma_i' \vert \sigma_i} : \rewrite{\hcomp{t}{x_i \mapsto u_i}}{\hcomp{t''}{x_i \mapsto u_i''}} : B$}
\end{prooftree}
\end{footnotesize}

\caption{Preservation rules}
\end{rules}
}

{\small

\begin{rules}
\unaryRule	{(\Delta \vdash \sigma_i : \rewrite{u_i}{u_i'} : A_i)_{i = 1, \dots, n}}
			{\Delta \vdash \indproj{k}{u_1', \dots, u_n'} \vert \hcomp{x_k}{x_i \mapsto \sigma_i} \equiv \sigma_k \vert \indproj{k}{u_1, \dots, u_n} : \rewrite{\hcomp{x_k}{x_i \mapsto u_i}}{u_k'} : A_k}
			{$(1 \leq k \leq n)$}

\unaryRule	{x_1 : A_1, \dots, x_n : A_n \vdash \tau : \rewrite{t}{t'} : B}
			{x_1 : A_1, \dots, x_n : A_n \vdash \subid{t'} \vert \tau \equiv \hcomp{\tau}{x_i \mapsto x_i} \vert \subid{t} : \rewrite{t}{\hcomp{t'}{x_i \mapsto x_i}} : B}
			{} \vspace{-\treeskip}
			
\begin{footnotesize}
\begin{prooftree}
\alwaysNoLine
\AxiomC{$(\Delta \vdash \mu_j : \rewrite{u_j}{u_j'} : A_j)_{j = 1, \dots m}$}
\noLine
\UnaryInfC{$(x_1 : A_1, \dots, x_m : A_m \vdash \sigma_i : \rewrite{v_i}{v_i'} : B_i)_{i = 1, \dots, n}$}
\noLine
\UnaryInfC{$y_1 : B_1, \dots, y_n : B_n \vdash \tau : \rewrite{t}{t'} : C$}
\alwaysSingleLine
\UnaryInfC{$\Delta \vdash \assoc{t', \ind{v}, \ind{u}} \vert 
\hcomp{\hcomp{\tau}{y_i \mapsto \sigma_i}}{x_j \mapsto \mu_j}\equiv 
\hcomp{\tau}{y_i \mapsto \hcomp{\sigma_i}{x_j \mapsto \mu_j}} \vert  \assoc{t, 
\ind{v}, \ind{u}}$\hspace{29mm}}
\noLine
\UnaryInfC{\hspace{65mm}$: \rewrite{\hcomp{\hcomp{t}{y_i \mapsto v_i}}{x_j 
\mapsto u_j}}{\hcompsmall{t'}{y_i \mapsto \hcompsmall{v_i'}{x_j \mapsto u_j'}}} 
: C$}
\end{prooftree} 
\end{footnotesize}
\caption{Naturality rules for structural rewrites}
\end{rules}

\begin{rules}
\begin{prooftree}
\AxiomC{$x_1 : A_1, \dots, x_n : A_n \vdash t : B$}
\AxiomC{$(\Delta \vdash u_i : A_i)_{i = 1, \dots, n}$}
\BinaryInfC{$\Delta \vdash \hcompsmall{t}{x_i \mapsto \indproj{i}{\ind{u}}} 
\vert \assoc{t,\ind{x},\ind{u}} \vert \hcomp{\subid{t}}{x_i \mapsto u_i} \equiv 
\id_{\hcomp{t}{x_i \mapsto u_i}} : \rewrite{\hcomp{t}{x_i \mapsto 
u_i}}{\hcomp{t}{x_i \mapsto u_i}} : B$}
\end{prooftree} \vspace{0\treeskip}

\begin{footnotesize}
\begin{prooftree}
\alwaysNoLine
\AxiomC{$(\Delta \vdash u_j : A_j)_{j = 1, \dots m}$}
\UnaryInfC{$(x_1 : A_1, \dots, x_m : A_m \vdash v_i : B_i)_{i = 1, \dots, n}$}
\AxiomC{$(y_1 : B_1, \dots, y_n : B_n \vdash w_j : C_k)_{k = 1, \dots, l}$}
\UnaryInfC{$z_1 : C_1, \dots, z_l : C_l \vdash t : D$}
\alwaysSingleLine
\BinaryInfC{$\Delta \vdash \hcomp{t}{z_k \mapsto \assoc{w_k, \ind{v}, \ind{u}}} \vert \assoc{t, \hcomp{\ind{w}}{y_j \mapsto v_j}, \ind{u}} \vert \hcomp{\assoc{t, \ind{w}, \ind{v}}}{x_j \mapsto u_j}$\hspace{40mm}}
\noLine
\UnaryInfC{\hspace{-7mm}$\equiv \assoc{t, \ind{w}, \hcomp{\ind{v}}{x_j \mapsto u_i}} \vert \assoc{\hcomp{t}{z_k \mapsto w_k}, \ind{v}, \ind{u}}$}
\noLine
\UnaryInfC{\hspace{40mm}$: \rewrite{\hcomp{\hcomp{\hcomp{t}{z_k \mapsto w_k}}{y_i \mapsto v_i}}{x_j \mapsto u_j}}{\hcomp{t}{z_k \mapsto \hcomp{w_k}{y_i \mapsto \hcomp{v_i}{x_j \mapsto u_j}}}} : D$}
\end{prooftree}
\end{footnotesize}

\caption{Biclone laws}
\end{rules}

\begin{rules}
\unaryRule	{\Gamma \vdash \alpha_1 : \rewrite{\hcomp{\pi_1}{u}}{t_1} : A_1 \quad \dots \quad \Gamma \vdash \alpha_n : \rewrite{\hcomp{\pi_n}{u}}{t_n} : A_n}
			{\Gamma \vdash \alpha_k \equiv \epsilonTimesInd{k}{t_1, \dots, t_n} \vert \hcomp{\pi_k}{\transTimes{\alpha_1, \dots, \alpha_n}} :\rewrite{\hcomp{\pi_k}{u}}{t_k} : A_k}
			{U1 ($1 \leq k \leq n$)}
			
\begin{small}
\unaryRule	{\Gamma \vdash \gamma : \rewrite{u}{\pair{t_1, \dots, t_n}} : \prodop_n(A_1, \dots, A_n)}
			{\Gamma \vdash \gamma \equiv \transTimes{\epsilonTimesInd{1}{\ind{t}} \vert \hcomp{\pi_1}{\gamma}, \dots, \epsilonTimesInd{n}{\ind{t}} \vert \hcomp{\pi_n}{\gamma}} : \rewrite{u}{\pair{t_1, \dots, t_n}} : \prodop_n(A_1, \dots, A_n)}
			{U2}
\end{small}
\unaryRule{\big(\Gamma \vdash \alpha_i \equiv  \alpha'_i : \rewrite{\hcomp{\pi_i}{u}}{t_i} : A_i\big)_{i = 1, \dots, n} }
		{\Gamma \vdash \transTimes{\alpha_1, \dots, \alpha_n} \equiv \transTimes{\alpha'_1, \dots, \alpha'_n} : \rewrite{u}{\pair{t_1, \dots, t_n}} : \prodop_n(A_1, \dots, A_n)}
			{cong}
\caption{Universal property of $\transTimes{\alpha}$}
\end{rules}

\begin{rules}
\begin{bprooftree}
\AxiomC{$\Gamma, x : A \vdash \alpha : \rewrite{\genevalterm{\wkn{u}{x}, x}}{t} : B$}
\RightLabel{\scriptsize U1}
\UnaryInfC{$\Gamma, x : A \vdash \alpha \equiv \epsilonExpRewr{t} \vert 
\genevalterm{\wkn{\transExplr{x \bind \alpha}}{x}, x} : 
	\rewrite{\genevalterm{\wkn{u}{x}, x}}{t} : B$}
\end{bprooftree}\vspace{\treeskip}

\begin{bprooftree}
\AxiomC{$\Gamma \vdash \gamma : \rewrite{u}{\lam{x}{t}} : \exptype{A}{B}$}
\RightLabel{\scriptsize U2}
\UnaryInfC{$\Gamma \vdash \gamma \equiv 
	\transExplr{x \bind \epsilonExpRewr{t} \vert 
		\genevalterm{\wkn{\gamma}{x}, x}} 
	: 
	\rewrite{u}{\lam{x}{t}} 
	: 
	\exptype{A}{B} $}
\end{bprooftree}\vspace{\treeskip}
			
\unaryRule	{\Gamma, x : A \vdash \alpha \equiv \alpha' : \rewrite{\genevalterm{\wkn{u}{x}, x}}{t} : B}		
			{\Gamma \vdash \transExp{x \bind \alpha} \equiv \transExp{x \bind \alpha'} : \rewrite{u}{\lam{x}{t}} : \exptype{A}{B}}
			{cong} 
\vspace{-0.5\treeskip}
\caption{Universal property of $\transExp{\alpha}$}
\end{rules}

}


{\small

\begin{rules}
	\begin{small}
\unaryRule	{\Gamma \vdash t_1 : A_1 \qquad \dots \qquad \Gamma \vdash t_n : A_n}
			{\Gamma \vdash \epsilonTimesInd{-k}{t_1, \dots, t_n} \vert \epsilonTimesInd{k}{t_1, \dots, t_n} \equiv \id_{\hcomp{\pi_k}{\pair{t_1, \dots, t_n}}}  : \rewrite{\hcomp{\pi_k}{\pair{t_1, \dots, t_n}}}{\hcomp{\pi_k}{\pair{t_1, \dots, t_n}}} : A_k} 
			{}
\end{small}

\unaryRule	{\Gamma \vdash t_1 : A_1 \qquad \dots \qquad \Gamma \vdash t_n : A_n}
			{\Gamma \vdash \epsilonTimesInd{k}{t_1, \dots, t_n} \vert \epsilonTimesInd{-k}{t_1, \dots, t_n} \equiv \id_{t_k} : \rewrite{t_k}{t_k} : A_k} 
			{}		
						
\unaryRule	{\Gamma \vdash t : \prodop_n(A_1, \dots, A_n)}
			{\Gamma \vdash \etaTimes{t}^{-1} \vert \etaTimes{t} \equiv \id_{t} : \rewrite{t}{t} : \prodop_n(A_1, \dots, A_n)}
			{}

\unaryRule	{\Gamma \vdash t : \prodop_n(A_1, \dots, A_n)}
			{\Gamma \vdash \etaTimes{t} \vert \etaTimes{t}^{-1}  \equiv \id_{\pair{\hcomp{\pi_1}{t}, \dots, \hcomp{\pi_n}{t}}} : \rewrite{\pair{\hcomp{\ind{\pi}}{t}}}{\pair{\hcomp{\ind{\pi}}{t}}} : \prodop_n(A_1, \dots, A_n)}
			{}
	\small
\unaryRule	{\Gamma \vdash u : \exptype{A}{B}}
			{\Gamma \vdash \etaExp{u} \vert \etaExp{u}^{-1} \equiv \id_{\lam{x}{\genevalterm{\wkn{u}{x}, x}}} : \rewrite{\lam{x}{\genevalterm{\wkn{u}{x}, x}}}{\lam{x}{\genevalterm{\wkn{u}{x}, x}}} : \exptype{A}{B}}
			{}
\normalsize
\vspace{0.5\treeskip}

\begin{bprooftree} 
\AxiomC{$\Gamma \vdash u : \exptype{A}{B}$} 
\UnaryInfC{$\Gamma \vdash \etaExp{u}^{-1} \vert \etaExp{u} \equiv \id_{u} : \rewrite{u}{u} : \exptype{A}{B}$} 
\end{bprooftree}
\quad
\begin{bprooftree} 
\AxiomC{$\Gamma, x : A \vdash t : B$} 
\UnaryInfC{$\Gamma, x : A \vdash \epsilonExpRewr{t} \vert \epsilonExpRewr{t}^{-1} \equiv \id_t  : \rewrite{t}{t} : B$} 
\end{bprooftree}

\small
\begin{prooftree} 
\AxiomC{$\Gamma, x : A \vdash t : B$} 
\UnaryInfC{$\Gamma, x : A \vdash \epsilonExpRewr{t}^{-1} \vert \epsilonExpRewr{t} \equiv \id_{\genevalterm{\wkn{(\lam{x}{t})}{x}, x}} : \rewrite{\genevalterm{\wkn{(\lam{x}{t})}{x}, x}}{\genevalterm{\wkn{(\lam{x}{t})}{x}, x}} : B$} 
\end{prooftree}
\normalsize

	\caption{Invertibility rules for pseudo cartesian closed structure}
\end{rules}

\begin{rules}
	\unaryRule {\Gamma \vdash t : B} 
{\Gamma \vdash \subid{t}^{-1} \vert \subid{t} \equiv \id_t : \rewrite{t}{t} : B} 
{} 
\unaryRule {x_1 : A_1, \dots, x_n : A_n \vdash t : B} 
{x_1 : A_1, \dots, x_n : A_n \vdash \subid{t} \vert \subid{t}^{-1} \equiv \id_t : \rewrite{\hcomp{t}{x_i \mapsto x_i}}{\hcomp{t}{x_i \mapsto x_i}} : B} 
{} 

{\small
\begin{prooftree} 
\AxiomC{$x_1 : A_1, \dots, x_n : A_n \vdash u_1 : A_1 \quad \dots \quad x_1 : A_1, \dots, x_n : A_n \vdash u_n : A_n$} 
\RightLabel{{\scriptsize $(1 \leq k \leq n$)}} 
\UnaryInfC{$x_1 : A_1, \dots, x_n : A_n \vdash \indproj{-k}{\ind{u}} \vert \indproj{k}{\ind{u}} \equiv \id_{\hcomp{x_k}{x_i \mapsto u_i}} : \birewrite{\hcomp{x_k}{x_i \mapsto u_i}}{\hcomp{x_k}{x_i \mapsto u_i}} : A_k$} 
\end{prooftree}
}

\begin{prooftree} 
\AxiomC{$x_1 : A_1, \dots, x_n : A_n \vdash u : B$} 
\RightLabel{{\scriptsize $(1 \leq k \leq n$)}} 
\UnaryInfC{$x_1 : A_1, \dots, x_n : A_n \vdash \indproj{k}{\ind{u}} \vert \indproj{-k}{\ind{u}} \equiv \id_u : \birewrite{u}{u} : A$} 
\end{prooftree}

{ \small
\begin{prooftree} 
\alwaysNoLine 
\AxiomC{$(\Delta \vdash u_j : A_j)_{j = 1, \dots m}$} 
\UnaryInfC{$(x_1 : A_1, \dots, x_m : A_m \vdash v_i : B_i)_{i = 1, \dots, n}$} 
\AxiomC{$y_1 : B_1, \dots, y_n : B_n \vdash t : C$} 
\alwaysSingleLine 
\RightLabel{} 
\BinaryInfC
	{$\Delta \vdash 
		\assoc{t, \ind{v}, \ind{u}}^{-1} \vert 
			\assoc{t,\ind{v},\ind{u}} \equiv 
	\id_{\hcompsmall{\hcompsmall{t}{v_i}}{u_j}} 
			: 
	\rewrite
		{\hcomp{\hcomp{t}{y_i \mapsto v_i}}{x_j \mapsto u_j}}
		{\hcomp{\hcomp{t}{y_i \mapsto v_i}}{x_j \mapsto u_j}} : C$} 
\end{prooftree} 
}

{ \small
\begin{prooftree} 
\alwaysNoLine 
\AxiomC{$(\Delta \vdash u_j : A_j)_{j = 1, \dots m}$} 
\UnaryInfC{$(x_1 : A_1, \dots, x_m : A_m \vdash v_i : B_i)_{i = 1, \dots, n}$} 
\AxiomC{$y_1 : B_1, \dots, y_n : B_n \vdash t : C$} 
\alwaysSingleLine 
\RightLabel{} 
\BinaryInfC{$\Delta \vdash 
	\assoc{t, \ind{v}, \ind{u}} \vert \assoc{t,\ind{v},\ind{u}}^{-1} 
	\equiv 
	\id_{\hcompsmall{t}{\hcompsmall{v_i}{u_j}}} 
		: 
	\rewrite{\hcomp{t}{y_i \mapsto \hcomp{v_i}{x_j \mapsto u_j}}}
			{\hcomp{t}{y_i \mapsto \hcomp{v_i}{x_j \mapsto u_j}}} : C$} 
\end{prooftree} 
}
	\caption{Invertibility of structural rewrites \label{fig:full:invertibility-of-structurals}}
\end{rules}

\begin{rules}
	\unaryRule	{\Gamma \vdash \tau : \rewrite{t}{t'} : A}
			{\Gamma \vdash \tau \equiv \tau : \rewrite{t}{t'} : A}
			{refl}
\unaryRule	{\Gamma \vdash \tau \equiv \tau' : \rewrite{t}{t'} : A}
			{\Gamma \vdash \tau' \equiv \tau : \rewrite{t}{t'} : A}
			{symm}
\binaryRule	{\Gamma \vdash \tau' \equiv \tau'' : \rewrite{t}{t'} : A}
			{\Gamma \vdash \tau \equiv \tau' : \rewrite{t}{t'} : A}
			{\Gamma \vdash \tau \equiv \tau'' : \rewrite{t}{t'} : A}
			{trans}
\binaryRule	{\Gamma \vdash \tau' \equiv \sigma' : \rewrite{t'}{t''} : A}
			{\Gamma \vdash \tau \equiv \sigma : \rewrite{t}{t'} : A}
			{\Gamma \vdash (\tau' \vert \tau) \equiv (\sigma' \vert \sigma)  : \rewrite{t}{t''} : A}
			{}
\binaryRule	{x_1 : A_1, \dots, x_n : A_n \vdash \tau \equiv \tau' : \rewrite{t}{t'} : B}
			{(\Delta \vdash \sigma_i \equiv \sigma_i' : \rewrite{u_i}{u'_i} : A_i)_{i=1, \dots, n}}
			{\Delta \vdash \horizComp{\tau}{x_i \mapsto \sigma_i} \equiv \horizComp{\tau'}{x_i \mapsto \sigma_i'} : \rewrite{\hcomp{t}{x_i \mapsto u_i}}{\hcomp{t'}{x_i \mapsto u'_i}} : B}
			{}

	\caption{Congruence rules \label{fig:full:congruence-rules}}
\end{rules}

}


\vfill
\clearpage

\section{The semantic interpretation \texorpdfstring{of $\langCartClosed$}{}}
\label{sec:semantic-interpretation}

We employ the same notation as Example~\ref{ex:cc-bicategory-to-cc-biclone} 
(page~\pageref{ex:cc-bicategory-to-cc-biclone}).

\begin{mynotation}
For any 
$A_1, \,\dots\, , A_n, B \in \baseCat \:\: (n \in \Nat)$
in an fp-bicategory
$\fpBicat\baseCat$ there 
exists a canonical equivalence
\begin{equation*} 
e_{\ind{A}, B}
:
\prodop_{n+1}(A_1, \,\dots\, , A_n, B) 
	\leftrightarrows
\prodop_{2} \left( \prod_{n}(A_1, \,\dots\, , A_n), B\right) 
:
\psinv{e}_{\ind{A}, B} 
\end{equation*}
where 
$
e_{\ind{A}, B} := 
\seq{\seq{\pi_1, \,\dots\, , \pi_n}, \pi_{n+1}}
$
and 
$
\psinv{e}_{\ind{A}, B} := 
\seq{\pi_1 \circ \pi_1, \,\dots\, , \pi_n \circ \pi_1, \pi_2}
$.
We denote the witnessing 2-cells by
\begin{equation*}
\begin{aligned}
\un_{\ind{A}, B} &: \Id_{\prod_n(A_1, \,\dots\, , A_n) \times B} \To 
{e_{\ind{A}, B}  \circ \psinv{e}_{\ind{A}, B}} \\
\co_{\ind{A}, B} &: 
\psinv{e}_{\ind{A}, B} \circ e_{\ind{A}, B} 
\To 
\Id_{\prod_{n+1}(A_1, \,\dots\, , A_n, B)}
\end{aligned}
\end{equation*}
\end{mynotation}

\begin{myconstr}[Semantic interpretation of $\langCartClosed$]
\label{constr:interpretation-of-langCartClosed}
For any unary $\langCartClosed$-signature $\sig$, 
cc-bicategory $\ccBicat{\baseCat}$ and 
$\langCartClosed$-signature morphism
$h : \sig \to \baseCat$, the \Def{interpretation} 
$h\sem{-}$ 
of the syntax of $\langCartClosed(\sig)$ is defined by induction. 

\subproof{Types.}
\begin{align*}
h\sem{B} &:= hB &\text{ for } B \text{ a base type } \\
h\sem{\prodop_n(A_1, \,\dots\, , A_n)} &:= 
	\prodop_n\big( h\sem{A_1}, \,\dots\, , h\sem{A_n} \big) \\
h\sem{\exptype{A}{B}} &:= (\exp{h\sem{A}}{h\sem{B}})
\end{align*}
On contexts, we set
$h\sem{x_1 : A_1, \,\dots\, , x_n : A_n} := 
\prod_n{\big( h\sem{A_1}, \,\dots\, , h\sem{A_n} \big)}$. 

\subproof{Terms.} Let
$\Gamma := (x_i : A_i)_{i=1, \,\dots\, , n}$ be any context. 
\begin{align*}
h\semlr{\Gamma \vdash x_i : A_i} &:= \pi_i \\
h\semlr{\Gamma \vdash c(x_1, \,\dots\, , x_n) : B} &:= h(c)  \\
h\semlr{p : \prodop_m(B_1, \,\dots\, , B_m) \vdash \pi_i(p) : B_i} &:= \pi_i \\
h\semlr{\Gamma \vdash \pair{t_1, \,\dots\, , t_m} : \prodop_m(B_1, \,\dots\, , 
B_m)} &:= 
	\seqlr{ h\sem{\Gamma \vdash t_1 : B_1}, \,\dots\, , h\sem{\Gamma \vdash t_m 
	: 
	B_m}	
		} \\
h\sem{f : (\exptype{A}{B}), x : A \vdash \evalterm(f, x) : B} &:=
	\eval_{h\sem{A}, h\sem{B}} \\
h\sem{\Gamma \vdash \lam{x}{t} : \exptype{B}{C}} &:=
	\lambda{\left( 
		{h\sem{\Gamma, x : B \vdash t : C} \circ \psinv{e}_{\ind{A}, B}}	
	\right)} \\
h\sem{\Delta \vdash \hcomp{t}{x_i \mapsto u_i} : B} &:=
	h\sem{\Gamma \vdash t : B} \circ \seq{h\sem{\Delta \vdash u_i : A_i}}_i 
\end{align*}
We omit easily-recovered typing information for the purpose of readability.

\newpage
\subproof{Rewrites.} For composition, constants and products the 
definition is direct:
\begin{align*}
h\semlr{\Gamma \vdash \id_t : \rewrite{t}{t} : B} &:= \id_{h\semlr{t}} \\
h\semlr{\Gamma \vdash \tau' \vert \tau : \rewrite{t}{t''} : B} &:= 
	h\semlr{\tau'} \vert h\semlr{\tau} \\
h\semlr{\Delta \vdash 
		\hcomp{\tau}{x_i \mapsto \sigma_i} 
		:
		\rewrite{\hcomp{t}{x_i \mapsto u_i}}
				{\hcomp{t'}{x_i \mapsto u_i'}} 
		: 
		B} &:= 
	h\semlr{\tau} \circ \seq{h\semlr{\sigma_i}}_i \\
h\semlr{\Gamma \vdash \constrewr : \rewrite{c(\ind{x})}{c'(\ind{x})} :B} &:=
	h(\constrewr) \\ 
h\semlr{\Gamma \vdash \epsilonTimesInd{k}{t_1, \,\dots\, , t_m} 
		: \rewrite{\hpi{k}{\pair{t_1, \,\dots\, , t_m}}}{t_k} : B_k} &:=
	\epsilonTimesInd{k}{h\semlr{t_1}, \,\dots\, , h\semlr{t_m}} \\
h\semlr{\Gamma \vdash 
	\transTimes{\alpha_1, \,\dots\, , \alpha_m} 
	: 
	\rewrite{u}
		{\pair{t_1, \,\dots\, , t_m}}
	: 
	\prodop_m(B_1, \,\dots\, , B_m) 
	} &:=
	\transTimes{h\sem{\alpha_1}, \,\dots\, , h\sem{\alpha_m}}
\end{align*}
The structural rewrites are interpreted by composites of structural 
isomorphisms. For $\indproj{k}{}$ and $\subid{}$ 
one has:
\begingroup
\addtolength{\jot}{.5em}
\begin{align*}
h\sem{\indproj{k}{u_1, \,\dots\, , u_n}} &:=
	\pi_k \circ \seq{h\semlr{u_i}}_i \XRA{\epsilonTimesInd{k}{h\semlr{\ind{u}}}}
		h\semlr{u_k} \\
h\semlr{\subid{t}} &:= 
	h\semlr{t} \XRA{\iso}
	h\semlr{t} \circ \Id_{h\semlr{\Gamma}} \XRA{h\semlr{t} \circ \etaTimes{\Id}}
	h\semlr{t} \circ \seqlr{\ind{\pi} \circ h\semlr{\Gamma}} \XRA{\iso} 
	h\semlr{t} \circ \seqlr{\ind{\pi}} 
\end{align*}
\endgroup
For $\assoc{}$ one has
\begin{td}[column sep = 3em]
h\semlr{\hcompthree{t}{u_i}{v_j}}
\arrow{rr}{h\semlr{\assoc{t, \ind{u}, \ind{v}}}}
\arrow[equals]{d} &
\: &
h\semlr{\hcomp{t}{\hcomp{u_i}{\ind{v}}}} \\

\left(h\semlr{t} 	
	\circ \seqlr{h\semlr{u_i}}_i\right) 
	\circ \seqlr{h\semlr{v_j}}_j
\arrow[swap]{r}{\iso} &
h\semlr{t} 	
	\circ \big({\seqlr{h\semlr{u_i}}_i
	\circ \seqlr{h\semlr{v_j}}_j}\big)
\arrow[swap]{r}{h\semlr{t} \circ \postName} &
h\semlr{t} 
	\circ \seqlr{h\semlr{u_i}_i 
	\circ \seqlr{h\semlr{\ind{v}}}}_i
\arrow[equals]{u}
\end{td}

Finally we come to the exponential rewrites $\epsilonExp_{t}$ and 
$\transExp{x \bind \alpha}$.
Suppose that 
$\Gamma \vdash u : \exptype{B}{C}$. 
Then
\begingroup
\addtolength{\jot}{.5em}
\begin{align*}
h\semlr{\Gamma, x : B \vdash \heval{\wkn{u}{x}, x} : C} &=
	\eval_{h\sem{B}, h\sem{C}} \circ 
		\seqlr{h\sem{\Gamma, x : B \vdash \wkn{u}{x} : \exptype{B}{C}}, 
		\pi_{n+1}} \\
	&= \eval_{h\sem{B}, h\sem{C}} \circ 
			\seqlr{h\sem{\Gamma \vdash u : \exptype{B}{C}} \circ 
				\seq{\pi_1, \,\dots\, , \pi_n}, \pi_{n+1}}
\end{align*}
\endgroup
The interpretation
$h\semlr{\Gamma, x : B \vdash \epsilonExpRewr{t} : 
\rewrite{
\heval{\wkn{(\lam{x}{t})}{x}, x} }
{t} 
: C }$ 
is the following composite, in which we abbreviate 
$h\sem{\Gamma, x : B \vdash t : C}$ 
by 
$\hsem{t}^{\Gamma, x : B}$:
\label{diag:hsem-on-counit}
\begin{equation*}
\makebox[\textwidth]{
\begin{tikzcd}[column sep = -1em, ampersand replacement = \&]
\eval_{h\sem{B}, h\sem{C}} \circ 
			\seqlr{\lambda {(\hsem{t}^{\Gamma, x : B} \circ 
				\psinv{e}_{h\sem{\ind{A}},h\sem{B}})} 
			\circ \seq{\pi_1, \,\dots\, , \pi_n}, \pi_{n+1}} 
\arrow{r}
\arrow[swap]{dd}{\iso} \&
\hsem{t}^{\Gamma, x : B}  \\
\: \&
\hsem{t}^{\Gamma, x : B} \circ \Id_{\prod(h\sem{\ind{A}}) \times h\sem{B}} 
\arrow[swap]{u}{\iso} \\
\eval_{h\sem{B}, h\sem{C}} \circ 
			\seqlr{\lambda (\hsem{t}^{\Gamma, x : A} \circ 
			\psinv{e}_{h\sem{\ind{A}},h\sem{B}})
			\circ \seq{\pi_1, \,\dots\, , \pi_n}, \Id_{h\sem{B}} \circ 
			\pi_{n+1}} 
\arrow[swap]{dd}{\eval \circ \fuse^{-1}} \&
\: \\
\: \&
\hsem{t}^{\Gamma, x : B} 
	\circ \left(\psinv{e}_{h\sem{\ind{A}},h\sem{B}}
	\circ e_{h\sem{\ind{A}},h\sem{B}}\right) 
\arrow[swap]{uu}{\hsem{t}^{\Gamma, x : B} \circ \co_{h\sem{\ind{A}},h\sem{B}}} 
\\
\eval_{h\sem{B}, h\sem{C}} \circ 
		\left(	\big(\lambda {(\hsem{t}^{\Gamma, x : B} \circ 
					\psinv{e}_{h\sem{\ind{A}},h\sem{B}}) 
				\times h\sem{B}}\big) 
			\circ e_{h\sem{\ind{A}},h\sem{B}}\right) 
\arrow[swap]{dd}{\iso} \&
\: \\
\: \&
\left(\hsem{t}^{\Gamma, x : B} 
	\circ \psinv{e}_{h\sem{\ind{A}},h\sem{B}}\right) 
	\circ e_{h\sem{\ind{A}},h\sem{B}}
\arrow[swap]{uu}{\iso} \\
\left(\eval_{h\sem{B}, h\sem{C}} \circ 
			\big(\lambda {(\hsem{t}^{\Gamma, x : B} \circ 
					\psinv{e}_{h\sem{\ind{A}},h\sem{B}}) 
				\times h\sem{B}}\big)\right) 
			\circ e_{h\sem{\ind{A}},h\sem{B}}
\arrow[bend right = 10]{ur}[swap, yshift=-1mm]
	{\epsilonExp_{(\hsem{t} \circ \psinv{e})} \circ e} \&
\: 
\end{tikzcd}
}
\end{equation*}

On the other hand, for a judgement 
$(\Gamma, x : B \vdash \alpha : \rewrite{\heval{\wkn{u}{x}, x}}{t} : C)$, 
the interpretation of $\alpha$ has type
\begin{equation} \label{eq:interpretation-of-alpha}
\eval_{h\sem{B}, h\sem{C}} \circ 
	\seqlr{\hsem{\Gamma \vdash u : \exptype{B}{C}} \circ 
		\seq{\pi_1, \,\dots\, , \pi_n}, \pi_{n+1}} 
\To 
\hsem{\Gamma, x : B \vdash t : C}
\end{equation}
To interpret 
$(\Gamma \vdash \transExp{x \bind \alpha} 
	: \rewrite{u}{\lam{x}{t}} : \exptype{A}{B})$ 
using the universal property 
of exponentials, we distort~(\ref{eq:interpretation-of-alpha}) into a 
composite 
$h\sem{\alpha}^{\circ}$ as in the diagram below. We suppress the subscripts on 
$e_{\ind{A}, B}$ and
$\psinv{e}_{\ind{A}, B}$ to fit the diagram better onto the page.
\begin{equation*}
\makebox[\textwidth]{
\begin{tikzcd}[column sep = 2em, ampersand replacement = \&]
\eval_{h\sem{B}, h\sem{C}} \circ 
	(h\sem{u}^\Gamma \times h\sem{B}) 
\arrow{r}{h\sem{\alpha}^{\circ}}
\arrow[swap]{d}{\iso} \&
h\sem{t}^{\Gamma, x : B} \circ \psinv{e} \\
\left(\eval_{h\sem{B}, h\sem{C}} 
	\circ (h\sem{u}^\Gamma \times h\sem{B})\right) 
	\circ \Id_{\prod_2((\prod_n h\sem{\ind{A}}), h\sem{B})}
\arrow[swap]{d}
	{\eval 
		\circ (h\sem{u}^\Gamma \times h\sem{B}) 
		\circ \un_{\prod_2((\prod_n h\sem{\ind{A}}), h\sem{B}))}} \&
\: \\
\left(\eval_{h\sem{B}, h\sem{C}} 
	\circ (h\sem{u}^\Gamma \times h\sem{B})\right)
	\circ \left(e 
	\circ \psinv{e}\right)
\arrow[swap]{d}{\iso} \&
\: \\
\left(\eval_{h\sem{B}, h\sem{C}} 
	\circ \left((h\sem{u}^\Gamma \times h\sem{B})\right)
	\circ e \right)
	\circ \psinv{e}
\arrow[swap]{d}
	{\eval \circ 
		\fuse \circ \psinv{e}} \&
\: \\
\left(\eval_{h\sem{B}, h\sem{C}} 
	\circ \seqlr{h\sem{u}^\Gamma \circ \seq{\pi_1, \dots, \pi_n}, 
			\Id_{h\sem{B}} \circ \pi_{n+1}}\right) 
	\circ \psinv{e}
\arrow[swap]{r}{\iso} \&
\left(\eval_{h\sem{B}, h\sem{C}} 
	\circ \seqlr{h\sem{u}^\Gamma \circ \seq{\ind{\pi}}, \pi_{n+1}}\right) 
	\circ \psinv{e} 
\arrow[swap]{uuuu}{h\sem{\alpha}^{\Gamma, x : B} \circ \psinv{e}} 
\end{tikzcd}
}
\end{equation*}
The unlabelled arrow is 
$\eval_{h\sem{B}, h\sem{C}} \circ 
		\seq{h\sem{u}^\Gamma \circ \epsilonTimesInd{1}{}, 
				\Id_{h\sem{B}} \circ \epsilonTimesInd{2}{}} \circ 
		\psinv{e}_{h\sem{\ind{A}}, h\sem{B}}$.
Finally, then, one has
\[
h\sem{\Gamma \vdash \transExp{x \bind \alpha} : 
	\rewrite{u}{\lam{x}{t}} : 
	\exptype{B}{C}} := \transExplr{
			h\sem{\Gamma, x : B \vdash \alpha : 
				\rewrite{\heval{\wkn{u}{x}, x}}{t} : C}^{\circ}}
\] 
\end{myconstr}

%

\chapter{The universal property of a bicategorical pullback} 
\label{chap:extra-proofs}

\section*{\:\:}

\vspace{-2\baselineskip}
Recall the following definition of a pullback (Definition~\ref{def:pullback} on page~\pageref{def:pullback}).

\begin{defnlocal}
Let $\cospanCat$ (for `cospan') denote the category $\cospan$ and $\baseCat$ be 
any 
bicategory. A \Def{pullback} of the cospan
$(X_1 \xra{f_1} X_0 \xla{f_2} X_2)$ in $\baseCat$ is a bilimit for the strict
pseudofunctor $\cospanCat \to \baseCat$ determined by this cospan.
\end{defnlocal}

We translate this into a presentation closer to that for categorical 
pullbacks---namely, that given by Lemma~\ref{lem:pullback-ump} (page~\pageref{lem:pullback-ump})---by 
showing that, for any $F : \cospanCat \to \baseCat$, 
there exists an equivalence of categories
$\Hom(\cospanCat, \baseCat)(\Delta B, F) 
\simeq 
\squaresCat{B}{F}{\baseCat}$,
where each category $\squaresCat{B}{F}{\baseCat}$ consists of 
\mbox{iso-commuting} squares and fill-ins.

\begin{defnlocal}
Let $\baseCat$ be any bicategory, $B \in \baseCat$ and
$F : \cospanCat \to \baseCat$ be a pseudofunctor. The category 
$\squaresCat{B}{F}{\baseCat}$ has objects 
triples 
$(\gamma_1, \gamma_2, \overline\gamma)$, where
$\gamma_i : B \to Fi \:\: (i= 1, 2)$ and 
$\overline\gamma$ is an invertible 2-cell as in the diagram
\begin{equation*} 
\begin{tikzcd}[column sep = small, row sep = small]
\: &
B 
\arrow[phantom]{dd}[description]{\twocell{\overline{\gamma}}}
\arrow[swap]{dl}{\gamma_1}
\arrow{dr}{\gamma_2} &
\: \\
F1 
\arrow[swap]{dr}{Fh_1} &
\: &
F2 
\arrow{dl}{Fh_2} \\
\: &
F0 &
\:
\end{tikzcd}
\end{equation*}
Morphisms 
$(\gamma_1, \gamma_2, \overline{\gamma}) 
	\to (\delta_1, \delta_2, \overline{\delta})$
are pairs of 2-cells 
$\modif_i : \gamma_i \To \delta_i \:\: (i=1, 2)$ such that
\begin{equation*} \label{eq:fill-in-defn}
\begin{tikzcd}[column sep = 4.5em]
F(h_2) \circ \gamma_2 
\arrow{r}{F(h_2) \circ \modif_2}
\arrow[swap]{d}{\overline{\gamma}} &
F(h_2) \circ \delta_2
\arrow{d}{\overline{\delta}} \\
F(h_1) \circ \gamma_1 
\arrow[swap]{r}{F(h_1) \circ \modif_1} &
F(h_1) \circ \delta_1
\end{tikzcd}
\end{equation*}
The identity on $(\gamma_1, \gamma_2, \overline\gamma)$ is 
$(\id_{\gamma_1}, \id_{\gamma_2})$ and composition is as in $\baseCat$.
\end{defnlocal}
 
The next lemma provides the components of the required equivalence.

\begin{lemmalocal} \label{lem:trans-2-cell-equivalence}
Let $\baseCat$ be a bicategory, $\cospanCat$ be the category $\cospan$, and
$F  : \cospanCat \to \baseCat$ a 
pseudofunctor. Then, for any $B \in \baseCat$ there exists an equivalence of 
categories 
$\Hom(\cospanCat, \baseCat)(\Delta B, F) 
	\simeq \squaresCat{B}{F}{\baseCat}$, 
where 
$\Delta : \baseCat \to \Hom(\cospanCat, \baseCat)$ denotes the diagonal 
pseudofunctor.
\begin{proof}
We begin by defining functors
$
K : \Hom(\cospanCat, \baseCat)(\Delta B, F) 
	\leftrightarrows \squaresCat{B}{F}{\baseCat} : L 
$. Take $K$ first. For a pseudonatural transformation 
$(\natTrans, \natCell) : \Delta B \To F$ with components as in the square
\begin{td}
B 
\arrow[phantom]{dr}[description]{\twocell{\natCell_i}}
\arrow[swap]{d}{\natTrans_i}
\arrow{r}{\Id_B} &
B 
\arrow{d}{\natTrans_0} \\
Fi
\arrow[swap]{r}{Fh_i} &
F0
\end{td} 
we define
$K(\natTrans, \natCell) := 
	(\natTrans_{1}, \natTrans_{2}, \overline{\gamma}_{(\natTrans, \natCell)})$, 
where 
\begin{equation} \label{eq:image-of-K}
\gamma_{(\natTrans, \natCell)} :=  
	F(h_2) \circ \natTrans_2 
		\XRA{\natCell_2^{-1}} 
	\natTrans_0 \circ \Id_B 
		\XRA{\natCell_{1}} 
	F(h_1) \circ \natTrans_1
\end{equation}
For morphisms, suppose 
$\modif : (\natTrans, \natCell) \to (\altNat, \altCell)$ is a modification. One 
thereby obtains 2-cells $\modif_i : \natTrans_i \To \altNat_i \:\: (i=1,2)$, and
\begin{td}[column sep = 4em, row sep = 2.7em]
\: &
F(h_2) \circ \natTrans_2
\arrow[swap, bend right = 70]{dd}{\overline{\gamma}_{(\natTrans, \natCell)}}
\arrow[phantom]{dr}[description]{\equals{modif. law}}
\arrow{r}{F(h_2) \circ \modif_2}
\arrow[swap]{d}{\natCell_2^{-1}} &
F(h_2) \circ \altNat_2
\arrow{d}{\altCell_2^{-1}}
\arrow[bend left = 70]{dd}{\overline{\gamma}_{(\altNat, \altCell)}} \\

\: &
\natTrans_0 \circ \Id_B
\arrow[phantom]{dr}[description]{\equals{modif. law}}
\arrow[swap]{r}{\modif_0 \circ \Id_B}
\arrow[swap]{d}{\natCell_1} &
\altNat_0 \circ \Id_B
\arrow{d}{\altCell_1} 
&
\:  \\

\: &
F(h_1) \circ \natTrans_1
\arrow[swap]{r}{F(h_1) \circ \modif_1} &
F(h_1) \circ \altNat_1
\end{td}
So we may define $K(\modif) := (\modif_1, \modif_2)$. 

Going the other way, for a triple $(\gamma_1, \gamma_2, \overline\gamma)$ we 
define $L(\gamma_1, \gamma_2, \overline\gamma)$ to be the pseudonatural 
transformation with components
\begin{align*}
\altNat_i &:= B \xra{\gamma_i} Fi \qquad\qquad\qquad \text{ for } i = 1, 2 \\
\altNat_0 &:= B \xra{\gamma_2} F2 \xra{Fh_2} F0 
\end{align*}
and witnessing 2-cells
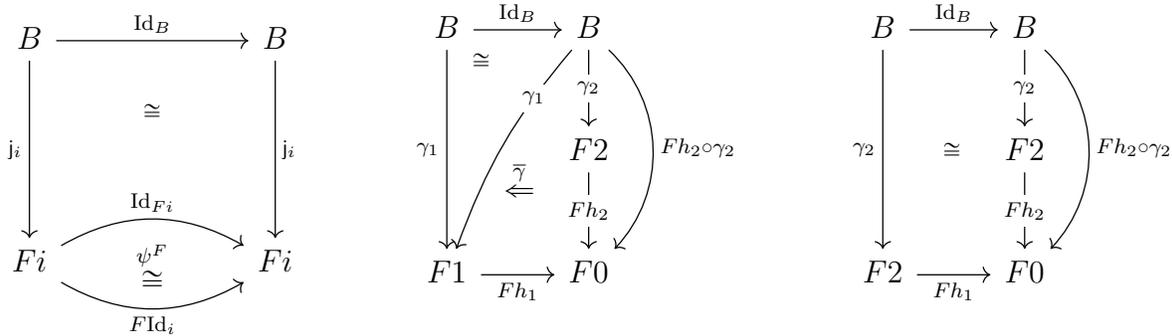
\begin{center}
\begin{tikzcd}[row sep = 2.5em]
B \arrow{rr}{\Id_B}
\arrow[swap]{dd}{\altNat_i} &
\arrow[phantom]{dd}[near start, description, font=\scriptsize]{\iso}
 &
B 
\arrow{dd}{\altNat_i} \\
\: &
\: &
\: \\
Fi 
\arrow[bend left]{rr}{\Id_{Fi}}
\arrow[phantom]{rr}[description]{\twocellIso{\psi^F}}
\arrow[swap, bend right]{rr}{F\Id_i} &
\: & 
Fi 
\end{tikzcd}
\qquad\quad
\begin{tikzcd}[row sep = 2.5em]
B 
\arrow[phantom]{dr}[description, font=\scriptsize, very near start]{\iso}
\arrow{r}{\Id_B}
\arrow[swap]{dd}{\gamma_1} &
B 
\arrow{d}[description]{\gamma_2}
\arrow[phantom]{ddl}[description, yshift=-4mm]{\twocell{\overline{\gamma}}} 
\arrow[bend right = 10]{ddl}[description, near start]{\gamma_1} 
\arrow[bend left = 45]{dd}{Fh_2 \circ \gamma_2} \\
\: &
F2 
\arrow{d}[description]{Fh_2} \\
F1 
\arrow[swap]{r}{Fh_1} & 
F0 
\end{tikzcd}
\qquad\quad
\begin{tikzcd}[row sep = 2.5em]
B 
\arrow[phantom]{ddr}[description, font=\scriptsize]{\iso}
\arrow{r}{\Id_B}
\arrow[swap]{dd}{\gamma_2} &
B 
\arrow{d}[description]{\gamma_2} 
\arrow[bend left = 45]{dd}{Fh_2 \circ \gamma_2} \\
\: &
F2 
\arrow{d}[description]{Fh_2} \\
F2 
\arrow[swap]{r}{Fh_1} & 
F0 
\end{tikzcd}
\end{center}
The naturality condition is trivial---there are no non-identity 2-cells in 
$\cospanCat$---and the unit law holds by definition, so the only thing to check 
is the associativity law. For this one must verify the axiom for each of the 
possible composites in $\cospanCat$, namely 
$\Id_i \circ \Id_i$, $\Id_0 \circ h_i$, and $h_i \circ \Id_i$. This is a 
long exercise.

On morphisms, for any $(\altModif_1, \altModif_2)$ in 
$\squaresCat{B}{F}{\baseCat}$,
we define $L(\altModif_1, \altModif_2)$ to be the modification with components
\begin{align*}
\altModif_i &:= \natTrans_i \XRA{\altModif_i} \altNat_i  \hspace{4cm} (i=1, 2) 
\\
\altModif_0 &:= F(h_2) \circ \natTrans_2 \XRA{F(h_2) \circ \altModif_2} 
F(h_2) \circ \altNat_2
\end{align*}
The only thing to check is the modification axiom, 
which we need to verify for the maps $h_1, h_2$ and $\Id_0, \Id_1, \Id_2$. Each 
of these is a simple calculation.

It remains to show that $K$ and $L$ form an equivalence. The composite 
$K \circ L$ is the identity. On the other hand, 
$LK(\natTrans, \natCell)$ has components $\natTrans_i$ for $i=1,2$
and $Fh_2 \circ \natTrans_2$ for $i=0$. One may then check that 
setting $\modif^{(\natTrans, \natCell)}_i := \id_{\natTrans_i}$ for $i = 1,2$ 
and
$\modif^{(\natTrans, \natCell)}_0 := \big( 
	Fh_2 \circ \natTrans_2 
	\XRA{\natCell_2^{-1}}
	\natTrans_0 \circ \Id_B
	\XRA\iso
	\natTrans_0 \big)
$
defines a modification $LK(\natTrans,\natCell) \to (\natTrans,\natCell)$. 
It remains to show that the modifications $\modif^{(\natTrans, \natCell)}$ 
are natural in $(\natTrans, \natCell)$. The $i=1$ and $i=2$ cases are trivial, 
and for $i=0$ one sees that, for any 
$\altModif : (\natTrans, \natCell) \to (\altNat, \altCell)$, 
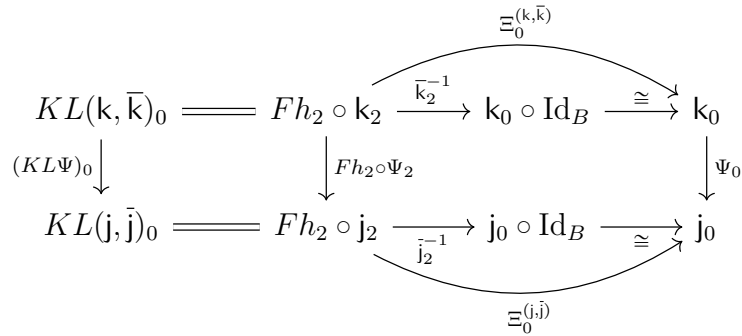
\begin{td}
KL(\natTrans, \natCell)_0
\arrow[equals]{r}
\arrow[swap]{d}{(KL\altModif)_0} &
Fh_2 \circ \natTrans_2
\arrow[bend left]{rr}{\modif^{(\natTrans, \natCell)}_0}
\arrow{r}{\natCell_2^{-1}}
\arrow{d}{Fh_2 \circ \altModif_2} &
\natTrans_0 \circ \Id_B
\arrow{r}{\iso} &
\natTrans_0 
\arrow{d}{\altModif_0} \\

KL(\altNat, \altCell)_0
\arrow[equals]{r} &
Fh_2 \circ \altNat_2
\arrow[swap]{r}{\altCell_2^{-1}}
\arrow[bend right, swap]{rr}{\modif^{(\altNat, \altCell)}_0} &
\altNat_0 \circ \Id_B
\arrow[swap]{r}{\iso} &
\altNat_0 
\end{td}
as required. It follows that 
$L \circ K \iso \id_{\Hom(\cospanCat, \baseCat)(\Delta B, F)}$, which completes 
the proof.
\end{proof}
\end{lemmalocal} 

The mapping $B \mapsto \squaresCat{B}{F}{\baseCat}$ extends to a 
pseudofunctor as follows. For $f : B' \to B$, we define 
$\squaresCat{f}{F}{\baseCat} : 
	\squaresCat{B}{F}{\baseCat} \to \squaresCat{B'}{F}{\baseCat}$
by setting 
$(\squaresCat{f}{F}{\baseCat})(\gamma_1, \gamma_2, \overline{\gamma})
	:= (\gamma_1 \circ f, \gamma_2 \circ f, \overline{\gamma} \circ f)$. 
Then for $\alpha : f \To f'$, the natural transformation 
$\squaresCat{\alpha}{F}{\baseCat}$ 
has components 
$\gamma_i \circ \alpha :
	\gamma_i \circ f \to \gamma_i \circ f'$.
This defines a pseudofunctor with unit and associativity witnessed by 
structural isomorphisms. In fact this pseudofunctor is equivalent to 
$\Hom(\cospanCat, \baseCat)(\Delta (-), F)$.

\begin{lemmalocal}
Let $\baseCat$ be a bicategory, $\cospanCat$ be the category $\cospan$, and
$F  : \cospanCat \to \baseCat$ a 
pseudofunctor. Then, writing 
$K_B : \Hom(\cospanCat, \baseCat)(\Delta B, F) \to  \squaresCat{B}{F}{\baseCat}$
for the functor constructed in Lemma~\ref{lem:trans-2-cell-equivalence},
the diagram below commutes for any $f : B' \to B$ in $\baseCat$:
\begin{td}[column sep = 3em]
\Hom(\cospanCat, \baseCat)(\Delta B, F)
\arrow[swap]{d}{K_B} 
\arrow{r}[yshift=2mm]{\Hom(\cospanCat, \baseCat)(\Delta f, F)} &
\Hom(\cospanCat, \baseCat)(\Delta B', F)
\arrow{d}{K_{B'}} \\
\squaresCat{B}{F}{\baseCat}
\arrow[swap]{r}{\squaresCat{f}{F}{\baseCat}} &
\squaresCat{B'}{F}{\baseCat}
\end{td}
\begin{proof}
For a pseudonatural transformation $(\natTrans,\natCell) : \Delta B \To F$, 
$(\squaresCat{f}{F}{\baseCat} \circ K_B)(\natTrans,\natCell)$
is the triple with 1-cells $\natTrans_1 \circ f$ and $\natTrans_2 \circ f$
and 2-cell 
\[
Fh_2 \circ (\natTrans_2 \circ f)
	\XRA{\iso} 
	(Fh_2 \circ \natTrans_2) \circ f
	\XRA{\gamma_{(\natTrans, \natCell)}}
	(Fh_1 \circ \natTrans_2) \circ f
	\XRA{\iso}
	Fh_1 \circ (\natTrans_2 \circ f)
\]
Here $\gamma_{(\natTrans, \natCell)}$ is the composite defined 
in~(\ref{eq:image-of-K}). 

On the other hand, writing $f_\ast := \Hom(\cospanCat, \baseCat)(\Delta f, F)$, 
one has that $f_\ast(\natTrans, \natCell)$ is the pseudonatural transformation 
with components $\natTrans_i \circ f$ and witnessing 2-cells given by composing 
$\natCell$ with the evident structural isomorphism:
\begin{td}
B'
\arrow[phantom]{dr}[description, yshift = 1mm]{\iso}
\arrow[swap]{d}{f} 
\arrow{r}{\Id_{B'}} &
B' 
\arrow{d}{f} \\
B 
\arrow[swap]{d}{\natTrans_i}
\arrow{r}{\Id_B}
\arrow[phantom]{dr}[description]{\twocell{\natCell_{i}}} &
B
\arrow{d}{\natTrans_i} \\
Fi \arrow[swap]{r}{Fh_i} &
F0
\end{td}
A short calculation shows that applying $K_{B'}$ to this pseudonatural 
transformation yields exactly 
$(\squaresCat{f}{F}{\baseCat} \circ K_B)(\natTrans,\natCell)$.
\end{proof}
\end{lemmalocal}

It follows that the functors $K_B$ are the components of a pseudonatural 
transformation. Since each $K_B$ is an equivalence, one obtains the following.

\begin{prooflesscorlocal}
Let $\baseCat$ be a bicategory, $\cospanCat$ be the category $\cospan$, and
$F  : \cospanCat \to \baseCat$ a 
pseudofunctor. Then
$\Hom(\cospanCat, \baseCat)(\Delta (-), F) \simeq  
\squaresCat{(-)}{F}{\baseCat}$
in $\Hom(\op{\baseCat}, \Cat)$. 
\end{prooflesscorlocal}

We can now use the fact that biequivalences preserve biuniversal arrows to
rephrase the universal property of a bicategorical pullback. For any bicategory 
$\baseCat$, let 
$({X_1 \xra{f_1} X_0 \xla{f_2} X_2})$ be any cospan 
and let $F$ be the strict pseudofunctor
$\cospanCat \to \baseCat$ it determines. The pullback of this cospan, when it 
exists, is a biuniversal arrow 
$(P, \lambda : \Delta P \To F)$
consisting of an object $P \in \baseCat$ and a pseudonatural transformation
$\lambda : \Delta P \To F$. The universal property then requires that, for any 
other pseudonatural transformation $\gamma : \Delta Q \To F$ there exists a 
1-cell $u : Q \to P$ and a universal modification 
$\epsilon : \lambda \circ \Delta u \To \gamma$, 
such that both the unit and the counit $\epsilon$ are invertible. 

We pass this data through the equivalence $K$. The pseudonatural transformations
$\lambda$ and $\gamma$ become iso-commuting squares:
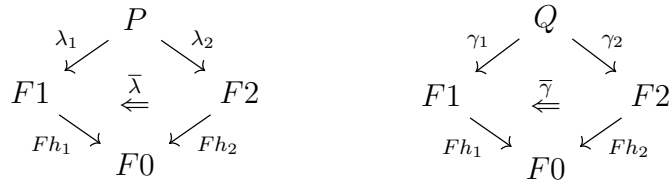
\begin{center}
\begin{tikzcd}[column sep = small, row sep = small]
\: &
P
\arrow[phantom]{dd}[description]{\twocell{\overline{\lambda}}}
\arrow[swap]{dl}{\lambda_1}
\arrow{dr}{\lambda_2} &
\: \\
F1 
\arrow[swap]{dr}{Fh_1} &
\: &
F2 
\arrow{dl}{Fh_2} \\
\: &
F0 &
\:
\end{tikzcd}
\qquad\qquad
\begin{tikzcd}[column sep = small, row sep = small]
\: &
Q 
\arrow[phantom]{dd}[description]{\twocell{\overline{\gamma}}}
\arrow[swap]{dl}{\gamma_1}
\arrow{dr}{\gamma_2} &
\: \\
F1 
\arrow[swap]{dr}{Fh_1} &
\: &
F2 
\arrow{dl}{Fh_2} \\
\: &
F0 &
\:
\end{tikzcd}
\end{center}
The pseudonatural transformation $\lambda \circ \Delta u$ then becomes 
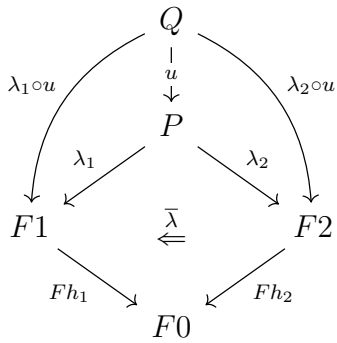
\begin{td}
\: &
Q 
\arrow[bend right]{ddl}[swap]{\lambda_1 \circ u}
\arrow[bend left]{ddr}{\lambda_2 \circ u}
\arrow{d}[description]{u} &
\: \\
\: &
P
\arrow[phantom]{dd}[description]{\twocell{\overline{\lambda}}}
\arrow[swap]{dl}{\lambda_1}
\arrow{dr}{\lambda_2} &
\: \\
F1 
\arrow[swap]{dr}{Fh_1} &
\: &
F2 
\arrow{dl}{Fh_2} \\
\: &
F0 &
\:
\end{td}
and the counit $\epsilon$ becomes a pair of 2-cells 
$\epsilon_i : \lambda_i \circ u \To \gamma_i$ which is universal among 2-cells 
satisfying the following:
\begin{td}[column sep = 4.5em]
Fh_2 \circ (\lambda_2 \circ u)
\arrow{r}{Fh_2 \circ \epsilon_2}
\arrow[swap]{d}{\iso} &
Fh_2 \circ \gamma_2
\arrow{ddd}{\overline{\gamma}} \\
(Fh_2 \circ \lambda_2) \circ u
\arrow[swap]{d}{\overline{\lambda} \circ u} &
\: \\
(Fh_1 \circ \lambda_1) \circ u
\arrow[swap]{d}{\iso} &
\: \\
Fh_1 \circ (\lambda_1 \circ u) 
\arrow[swap]{r}{Fh_1 \circ \epsilon_1} &
Fh_1 \circ \gamma_1
\end{td}
Starting this diagram from $(Fh_2 \circ \lambda_2) \circ u$ and inverting the 
isomorphisms, one obtains the fill-in requirement from 
Lemma~\ref{lem:pullback-ump}. One may now see that the remaining conditions of 
Lemma~\ref{lem:pullback-ump} are exactly those making $\epsilon$ universal.

\printnomenclature
\label{nomenclature}
\bibliographystyle{alpha}
\bibliography{bicats_db}

\end{document}